\documentclass[11pt]{amsart}
\usepackage[colorlinks,citecolor=red,pagebackref,hypertexnames=false]{hyperref}

\usepackage{verbatim} %for \begin{comment}\end{comment}
\usepackage{placeins}

\makeatletter
\def\@tocline#1#2#3#4#5#6#7{\relax
	\ifnum #1>\c@tocdepth % then omit
	\else
	\par \addpenalty\@secpenalty\addvspace{#2}%
	\begingroup \hyphenpenalty\@M
	\@ifempty{#4}{%
		\@tempdima\csname r@tocindent\number#1\endcsname\relax
	}{%
		\@tempdima#4\relax
	}%
	\parindent\z@ \leftskip#3\relax \advance\leftskip\@tempdima\relax
	\rightskip\@pnumwidth plus4em \parfillskip-\@pnumwidth
	#5\leavevmode\hskip-\@tempdima
	\ifcase #1
	\or\or \hskip 1em \or \hskip 2em \else \hskip 3em \fi%
	#6\nobreak\relax
	\dotfill\hbox to\@pnumwidth{\@tocpagenum{#7}}\par
	\nobreak
	\endgroup
	\fi}
\makeatother

\usepackage{tikz,amsthm,amsmath,amstext,amssymb,amscd,epsfig,euscript, mathrsfs,pspicture,multicol, mathtools,graphpap,graphics,graphicx,times,enumerate,subfig,sidecap,wrapfig,color}
\IfFileExists{dsfont.sty}{\usepackage{dsfont}}{}

\usepackage{xcolor,booktabs, multirow} %to create nice tables for notation
\RequirePackage{caption}% typesetting of captions
\usepackage{enumitem}

\usepackage{ hyperref}
\usepackage{enumerate}

\numberwithin{equation}{section}

\def\caD{{\mathcal{D}}}
\def\bbD{{\mathbb{D}}}
\def\Ch{{\mathrm{Ch}}}
\def\Tree{{\mathrm{Tree}}}
\def\Top{{\mathrm{Top}}}
\def\Next{{\mathrm{Next}}}
\def\Mdb{{\mathrm{MDb}}}
\def\Stop{{\mathbf{St}}}
\def\Stopp{{\mathrm{Stop}}}

\def\sd{\mathfrak{d}} %%% mathfrak d

\def\ve{\varepsilon}

\def\I{\mathop\textup{I}}

\def\Rn1{\mathbb{R}^{n+1}}
\def\el{\ell}

\def\ext{\mathop\mathrm{ext}} 	
 	
\DeclareMathOperator{\diam}{diam}
\def\dim{\mathop\mathrm{dim}} 					%dimension
\def\dist{\textup{dist}} 						%distance
\def\supp{\mathop\mathrm{supp}}					%support
\def\loc{\mathop\mathrm{loc}}						%locally
\def\tr{\mathop\mathrm{tr }}					%cofactor
\def\term{\mathop\mathrm{term }}					%cofactor
\def\op{\mathop\mathrm{op }}
\renewcommand{\div}{\mathop\mathrm{div }}			%divergence
\newcommand{\ceil}[1]{\left\lceil #1 \right\rceil}

\def\Xint#1{\mathchoice
	{\XXint\displaystyle\textstyle{#1}}%
	{\XXint\textstyle\scriptstyle{#1}}%
	{\XXint\scriptstyle\scriptscriptstyle{#1}}%
	{\XXint\scriptscriptstyle\scriptscriptstyle{#1}}%
	\!\int}
\def\XXint#1#2#3{{\setbox0=\hbox{$#1{#2#3}{\int}$ }
		\vcenter{\hbox{$#2#3$ }}\kern-.58\wd0}}

\def\avint{\Xint-}

\theoremstyle{plain}
\newtheorem{theorem}{Theorem}

\newtheorem{corollary}[theorem]{Corollary}

\newtheorem{lemma}[theorem]{Lemma}

\newtheorem{proposition}[theorem]{Proposition}

\theoremstyle{definition}

\newtheorem{definition}[theorem]{Definition}
\newtheorem{remark}[theorem]{Remark}

\numberwithin{equation}{section}
\numberwithin{theorem}{section}

\IfFileExists{mathb10.tfm}{
	\DeclareFontFamily{U}{mathb}{\hyphenchar\font45} 
	\DeclareFontShape{U}{mathb}{m}{n}{
		<5> <6> <7> <8> <9> <10> gen * mathb
		<10.95> mathb10 <12> <14.4> <17.28> <20.74> <24.88> mathb12
	}{}
	\DeclareSymbolFont{mathb}{U}{mathb}{m}{n}

	\DeclareMathSymbol{\toitself}      {3}{mathb}{"FD}  %*

}{\providecommand{\toitself}{\circlearrowleft}}

\newcommand{\vv}{\vspace{2mm}}
\newcommand{\vvv}{\vspace{4mm}}

\def\R{\mathbb{R}}

\begin{document}

	\title[Local $Tb$ theorems, Neumann problem, and caloric measure]{Local $Tb$ theorems with applications to the Neumann problem and caloric measure}
	
	\author[A. Merlo]{Andrea Merlo}
	\address{Andrea Merlo\\ Departamento de Matem\'aticas, Universidad del Pa\'is Vasco (UPV/EHU), Barrio Sarriena s/n 48940 Leioa, Spain and\\
		Ikerbasque, Basque Foundation for Science, Bilbao, Spain.}
	\email{andrea.merlo@ehu.eus}
	
	\author[M. Mourgoglou]{Mihalis Mourgoglou}
	\address{Mihalis Mourgoglou\\ Departamento de Matem\'aticas, Universidad del Pa\'is Vasco (UPV/EHU), Barrio Sarriena s/n 48940 Leioa, Spain and\\
		Ikerbasque, Basque Foundation for Science, Bilbao, Spain.}
	\email{michail.mourgoglou@ehu.eus}

	\author[C. Puliatti]{Carmelo Puliatti}
	\address{Carmelo Puliatti \\Departament de Matemàtiques, Universitat Autònoma de Barcelona, 08193 Bellaterra (Barcelona), Catalonia.}
	\email{carmelo.puliatti@uab.cat}

	\thanks{A.M. was  supported by IKERBASQUE, and partially supported by the project PID2024-157724NB-I00, funded by MICIU/AEI/10.13039/501100011033/FEDER, EU, by the Europa Excelencia grant EUR2025-165042, funded by the Ministerio de Ciencia, Innovación y Universidades and Agencia Estatal de Investigación, and by  IT1875-26 (Basque Government).
		M.M. was supported by IKERBASQUE and partially supported by the project PID2024-157724NB-I00, funded by MICIU/AEI/10.13039/501100011033/FEDER, EU, by IT-1615-22, and IT1875-26 (Basque Government).
		C.P. was supported by the Ram\'on y Cajal fellowship grant RYC2023-045598-I (MICIU, Spain), and partially supported by the grants PID2023-150984NB-I00 and PID2024-155320NB-I00 funded by MICIU/AEI/10.13039/501100011033/ FEDER, EU}
	
	\newcommand{\mih}[1]{\marginpar{\color{red} \scriptsize \textbf{Mi:} #1}}
	\newcommand{\car}[1]{\marginpar{\color{blue} \scriptsize \textbf{Carmelo:} #1}}
	\maketitle

\begin{abstract}
	We prove a local $Tb$ theorem with testing conditions on metric balls, for possibly non-doubling measures of polynomial growth on geometrically doubling metric spaces. It resolves the maximal-truncation version of Hofmann's local $Tb$ problem for all $1<p,q\leq\infty$, without assuming antisymmetry of the Calder\'on-Zygmund kernel.
	Our theorem allows accretive systems of complex testing measures satisfying a quantitative small-set condition; these measures need not be absolutely continuous
	with respect to the underlying measure. Testing requires only a	one-level $L^0$ bound for maximal truncations, with an alternative for balls on which the underlying measure concentrates in a set of small relative diameter. We also remove antisymmetry from the known buffered and unbuffered fixed-truncation criteria.
	The proof combines a mixed corona decomposition adapted to the two testing systems with two-measure $L^2$ estimates for suppressed operators.
	We obtain the following applications:
	\begin{itemize}
		\item If $\Omega\subset\mathbb R^{n+1}$, $n\geq2$, is a bounded domain with Ahlfors regular boundary satisfying the interior corkscrew condition, and the Neumann problem for the Laplacian in $\Omega$ is solvable with data in $L^p$ for some $1<p<\infty$, then $\partial\Omega$ is uniformly rectifiable. The proof applies the matrix-valued version of the theorem to the Cauchy-Clifford operator.
		\item In the parabolic space $\mathbb P^n$, local time-directed $A_\infty$-type conditions on caloric and adjoint caloric measure imply
		$L^2(\sigma)$-boundedness of the spatial gradient of the heat single layer potential for boundary measures $\sigma$ with upper parabolic $(n+1)$-growth. This gives a parabolic counterpart of the harmonic-measure theorem of Mourgoglou and Tolsa without Ahlfors regularity, corkscrew or Harnack chain assumptions.
	\end{itemize}
\end{abstract}

	\tableofcontents

	{
		\section{Introduction}
		A central line of research in harmonic analysis aims to determine sufficient conditions for the boundedness of singular integral operators (SIOs) in terms of their action on suitable testing functions, and has given rise to so-called $T1$- and $Tb$-type theorems. In this paper, we work in the framework of \textit{local} $Tb$ theorems, in which the $L^2$-boundedness of a SIO is deduced from testing conditions on functions $b_B$, and $b^*_B$ for the adjoint operator, adapted to a suitable family of balls $B$ in the ambient space. Motivated by applications to PDEs, Steve Hofmann posed in 2008, during his plenary short course at the 8th International Conference on Harmonic Analysis and Partial Differential Equations held in El Escorial (Madrid, Spain), the problem of whether suitable $L^p$-integrability conditions on the testing functions yield $L^2$-boundedness of a singular integral operator; see \cite[Section 3.3.1]{Ho10} or the survey by S. Hofmann and A. McIntosh \cite[Section 5.1.2]{HoMc11}. It turned out to be a delicate question, even in the case of SIOs associated with the Lebesgue measure on Euclidean spaces, and several advances have appeared in the recent literature in non-homogeneous harmonic analysis, as we discuss in Section \ref{section:hystory}.
		
		\vv
		Our main result, Theorem \ref{theorem:main_theorem_hofmann_problem}, is a local $Tb$ theorem with testing conditions on metric balls, for possibly non-doubling measures of polynomial growth on geometrically doubling metric spaces. It applies to general Calder\'on-Zygmund kernels, with two unrelated accretive systems of finite complex testing measures for the operator and for its adjoint.
		
		The testing assumptions are formulated on suitable families of doubling balls with small boundary, and they do not refer to any dyadic lattice. Passing from these assumptions to a dyadic testing system is delicate, since the accretivity on a ball need not persist on its dyadic subcells, and the testing balls associated with nested cells need not be nested. Our proof merges this geometry of balls with two independent random dyadic systems and a fixed David-Mattila lattice, by means of a mixed corona decomposition and two-measure estimates for suppressed operators.
		
		The testing measures satisfy a quantitative small-set condition, and they need not be absolutely continuous with respect to the underlying measure. In the operator-testing alternative, only a one-level $L^0$ bound for the maximal truncations is required on a contracted testing ball, outside a set of small relative measure. We also allow an alternative for the balls on which a fixed proportion of the underlying measure is concentrated on a set of small relative diameter. We remark that allowing this  dichotomy is not only a technical improvement but crucial to obtain the parabolic version of \cite{MT20} in full generality\footnote{Otherwise we need to assume that the poles are corkscrew points, which we do not.}.
		
		For testing functions with $L^p$ and $L^q$ size bounds, the theorem resolves the maximal-truncation version of Hofmann's local $Tb$ problem in the full range $1<p,q\leq\infty$, see Corollary \ref{corollary:hofmann_problem_maximal}. Combined with the Cotlar-type argument of \cite{MMT18}, it also gives buffered and unbuffered criteria for fixed truncations, in the corresponding ranges of exponents, for kernels that are not necessarily antisymmetric and with testing conditions on metric balls, see Corollary \ref{corollary:hofmann_problem}. 
	
	\vv
	We also obtain two applications to elliptic and parabolic PDEs. In Section \ref{section:intro_neumann}, we show that, if $\Omega\subset\mathbb R^{n+1}$, $n\geq2$, is a bounded domain with Ahlfors regular boundary satisfying the interior corkscrew condition, and the Neumann problem for the Laplacian in $\Omega$ is solvable with data in $L^p$ for some $1<p<\infty$, then $\partial\Omega$ is uniformly rectifiable, see Theorem \ref{theorem:neumann_UR}. The proof relies on a version of Theorem \ref{theorem:main_theorem_hofmann_problem} with matrix-valued kernels and testing functions, Theorem \ref{theorem:matrix_Tb}, that we apply to the Cauchy-Clifford operator, with testing functions constructed from solutions of the Neumann problem.
%	A priori, these testing functions are only in $L^p$, and the Neumann problem only controls their maximal truncations in $L^p$.
	In Section \ref{section_applications}, we show that local time-directed $A_\infty$-type conditions on caloric and adjoint caloric measure imply the $L^2$-boundedness of the spatial gradient of the heat single layer potential, see Theorem \ref{theorem:application_caloric_measure}.
	}
	
	\vvv
	
	Following the terminology in \cite{HM12}, we say that a metric space $(X,d)$ is \textit{geometrically doubling} if there exists $D\in \mathbb N$ such that, for every $x\in X$ and $r>0$, the open ball $B(x,r)\coloneqq \{y\in X: d(x,y)<r\}$ can be covered by at most $D$ balls $B(x_j, r/2)$, with $x_j \in X$.
	If $B\subset X$ is a ball, we also use the notation $r(B)$ for its radius and, if $B=B(x,r)$ and $\lambda>0$, we define its dilation $\lambda B\coloneqq B(x,\lambda r)$.
	
	\begin{definition}\label{def:CZ_kernel}
		We say that $K\colon (X\times X) \setminus \{(x,x):x\in X\}\to \mathbb C$ is an $n$-dimensional \textit{Calder\'on-Zygmund kernel} if there exist $C>0$ and $\gamma \in (0,1]$ such that:
		\begin{enumerate}
			\item For all $x,y\in X$ with $x\neq y$ we have
			\[
			|K(x,y)|\leq \frac{C}{d(x,y)^n}.
			\]
			\item For all $x,x',y\in X$ with $d(x,y)\geq 2d(x,x')$ we have
			\[
			|K(x,y)-K(x',y)|+|K(y,x)-K(y,x')|\leq C\, \frac{d(x,x')^\gamma}{d(x,y)^{n+\gamma}}.
			\] 
		\end{enumerate}
	\end{definition}
	
	Given a Calder\'on-Zygmund kernel $K$, a non-negative Radon measure $\mu$ on $X$, and $\varepsilon>0$ we define its associated \textit{singular integral operator} (SIO) as 
	\[
	T_\mu f(x)\coloneqq \int K(x,y)f(y)\, d\mu(y)\, \qquad \text{ for }f\in L^1(\mu) \text{ and }x\not \in \supp(f\mu)
	\]
	and its associated \textit{$\varepsilon$-truncations} 
	\[
	T_{\mu,\varepsilon}f(x)=\int_{d(x,y)>\varepsilon} K(x,y)f(y)\, d\mu(y).
	\]
	These integrals are absolutely convergent for $f\in L^1(\mu)$ and, if $\mu$ has polynomial growth as in \eqref{eq:def_growth} below, also for $f\in L^p(\mu)$, $1<p<\infty$.
	We say that $T_\mu$ is bounded on $L^p(\mu)$ for $1\leq p<\infty$ (resp. from $L^{1}(\mu)$ to $L^{1,\infty}(\mu)$) if $T_{\mu,\varepsilon}$ is bounded on $L^p(\mu)$ (resp. from $L^{1}(\mu)$ to $L^{1,\infty}(\mu)$) uniformly in $\varepsilon>0$. In this case, we write
	\[
	\|T_\mu\|_{L^p(\mu)\to L^p(\mu)}\coloneqq \sup_{\varepsilon>0}\|T_{\mu,\varepsilon}\|_{L^p(\mu)\to L^p(\mu)}.
	\]
	We also define the \textit{maximal singular integral} as
	\(
	T_{\mu,\sharp}f(x)=\sup_{\varepsilon>0}|T_{\mu, \varepsilon}f(x)|.
	\)
	
	We further understand that $T_\varepsilon\mu(x)\coloneqq T_{\mu,\varepsilon}1(x)$ if $\mu$ is finite. This definition also generalizes to complex measures, and we define $T_\sharp\mu(x)$ analogously.
	Finally, given a singular integral operator $T$ associated with a Calder\'on-Zygmund kernel $K(x,y)$, we denote by $T^*$ its adjoint, which corresponds to the SIO associated with the kernel $K(y,x)$.
	
	For $n>0$, we say that a Radon measure  $\mu\in M_+(X)$ has \textit{growth of degree $n$}, or, for brevity, \textit{$n$-growth} with constant $\mathfrak{C}_0>0$, if
	\begin{equation}\label{eq:def_growth}
		\mu(B(x,r))\leq \mathfrak{C}_0 \,r^n, \qquad \text{ for all }x\in X  \text{ and all }r>0.
	\end{equation}
	
	{Before stating the precise testing assumptions, we briefly explain their form. We allow finite complex measures in place of testing functions, with a total variation bound and a small-set condition requiring that sets of small relative $\mu$-measure carry only a small fraction of the total variation, in the spirit of \cite[Theorem 1.6]{MMV21}.
		On each testing ball, we allow two alternatives: either the maximal	truncations are bounded on a contracted ball outside a set of small relative $\mu$-measure, or a fixed proportion of the underlying measure is concentrated on a set of small relative diameter. The second alternative allows us to handle balls where the measure concentrates near the poles in our application to caloric measure.
		
%		{As mentioned above, instead of testing functions, we state our result under more general testing assumptions, and consider finite complex measures supported on balls. The size conditions on the testing functions are replaced by the following small-set condition, which is a quantitative absolute continuity condition with respect to $\mu$ at a fixed relative scale, in the spirit of \cite[Theorem 1.6]{MMV21}. The testing conditions on the operator are formulated as \textit{a dichotomy}: on each ball $B$, either a fixed proportion of $\mu(B)$ is concentrated on a set of diameter at most $\eta\,r(B)$, which forces a large increase of the density of $\mu$ at a smaller scale, or the maximal truncations of the testing measure are bounded outside a set of small measure, by a constant that may depend on $\eta$. We remark that in our application in Theorem \ref{theorem:application_caloric_measure}, as for harmonic measure in \cite{MT20}, the testing measures are normalized caloric measures with poles in $B$, whose maximal truncations can only be controlled outside a ball of radius $\eta\,r(B)/2$ centered at the pole, and the first alternative accounts for the balls on which the measure concentrates near the pole.
			
			\begin{definition}[Accretive measure testing system]
				\label{def:accretive_testing_system}
				Let $\mathscr B$ be a family of balls, and fix	$\mathfrak C_1\geq1$, $\varepsilon_0\in(0,1)$, and $\tau,\zeta\in(0,1)$. An accretive measure testing system for $T$ consists of finite complex Radon measures
				$\{\widehat\nu_B\}_{B\in\mathscr B}$ such that
				\begin{equation}\label{eq:measure_testing_size}
				|\widehat\nu_B|(X\setminus B)=0,\qquad\widehat\nu_B(B)=\mu(B),\qquad
					\|\widehat\nu_B\|\leq\mathfrak C_1\mu(B).
				\end{equation}
				For every Borel set $A\subset B$, the following condition holds:
				\begin{equation}\label{eq:measure_testing_small_sets}
					\text{If }\mu(A)\leq\varepsilon_0\, \mu(B),
					\qquad\text{ then }\qquad|\widehat\nu_B|(A)\leq\frac{\|\widehat\nu_B\|}{200\mathfrak C_1}.
				\end{equation}
				Moreover, for every $\eta\in(0,1/10)$ there is a finite constant $\mathfrak C_2(\eta)>0$ such that, for each $B\in\mathscr B$, there is a Borel set $U_B=U_B(\eta)\subset B$ satisfying one	of the following alternatives:
				\begin{enumerate}
					\item[(i)] There is a Borel set $V_B\subset U_B$ such that
					\begin{equation}\label{eq:testing_opt_1}
						\diam(V_B)\leq\eta\,  r(B)\qquad \text{ and }\qquad \mu(V_B)\geq\frac\zeta2\, \mu(B).
					\end{equation}
					\item[(ii)] We have
					\begin{equation}\label{eq:testing_opt_2}
						\mu(U_B)<\zeta\mu(B),\qquad	T_\sharp\widehat\nu_B(x)\leq\mathfrak C_2(\eta)
						\quad\text{for every }x\in(1-\tau)B\setminus U_B.
					\end{equation}
				\end{enumerate}
				An accretive measure testing system for $T^*$ is defined in the same way, with $T^*$ in place of $T$, and we always assume that the systems for $T$ and $T^*$ have the same constants $\mathfrak C_1$, $\varepsilon_0$, $\tau$, $\zeta$, and $\mathfrak C_2(\eta)$.
			\end{definition}
			
			Alternative \eqref{eq:testing_opt_2} is equivalent to the one-level $L^0$ condition
			\begin{equation}\label{eq:measure_testing_L0}
				\mu\bigl(\{x\in(1-\tau)B:	T_\sharp\widehat\nu_B(x)>\mathfrak C_2(\eta)\}\bigr)<\zeta\mu(B),
			\end{equation}
			since one can take $U_B$ to be this superlevel set. In particular, no integrability of $T_\sharp\widehat\nu_B$ is required, and the level $\mathfrak C_2(\eta)$ may be arbitrarily large.
			Similarly, condition \eqref{eq:measure_testing_small_sets} does not imply that $\widehat\nu_B$ is absolutely continuous with respect to $\mu$. %Moreover, the proof only uses that $\widehat\nu_B$ is concentrated on $B$, that is, $|\widehat\nu_B|(X\setminus B)=0$, rather than $\supp\widehat\nu_B\subset B$.
			\vv
			
			We now state the main result of the paper.
			\begin{theorem}\label{theorem:main_theorem_hofmann_problem}
				Let $(X,d)$ be geometrically doubling with constant $D>1$, let $\mu$ be a Radon measure with $n$-growth constant $\mathfrak C_0$, for some $n\geq 1$, and let $T$ have an $n$-dimensional Calder\'on-Zygmund kernel with size and regularity constants $C$ and $\gamma\in(0,1)$. Fix $\delta\in(0,1]$ and $\Lambda\geq100$. There is $\beta_0=\beta_0(D,n,\delta,\Lambda)>1$ such that the following holds for $\beta\geq\beta_0$ and $\mathfrak C_3\geq220$.
				
				Let $\mathscr B$ consist of all balls $B=B(x,r)$ with $x\in\supp\mu$, $0<r\leq\diam(\supp\mu)$, and	$\mathfrak C_3$-small boundary with respect to $\mu$, satisfying
				\begin{equation}\label{eq:improved_doubling}
					\mu(\Lambda B)\leq\beta\, \mu(\delta B).
				\end{equation}
				We fix $\mathfrak C_1\geq1$ and $\varepsilon_0\in(0,1)$, and define the parameters
				\begin{equation}\label{eq:measure_testing_tau_zeta}
					\tau_0=\min\Bigl\{\frac1{100},
					\frac{\varepsilon_0}{4\beta\mathfrak C_3}\Bigr\}
					\qquad\text{ and }\qquad 	\zeta_0=\min\Bigl\{\frac1{100},\frac{\varepsilon_0}{8}\Bigr\}.
				\end{equation}
				Suppose that, for $0<\tau<\tau_0$ and $0<\zeta<\zeta_0$, $T$ and $T^*$ admit accretive measure testing systems on $\mathscr B$ in the sense of Definition
				\ref{def:accretive_testing_system}.	Then $T_\mu$ is bounded on $L^2(\mu)$.
				
				The bound depends only on $D,n,\mathfrak C_0,C,\gamma,\delta,\Lambda,\beta,	\mathfrak C_3,\mathfrak C_1,\varepsilon_0,\zeta$ and
				$\mathfrak C_2(\eta_*)$, where $\eta_*$ is the value chosen	in \eqref{eq:choice_eta_Stop_1}. In particular, the thresholds
				\eqref{eq:measure_testing_tau_zeta} are independent of	$\mathfrak C_2(\eta_*)$.
			\end{theorem}
			
			Note that the thresholds \eqref{eq:measure_testing_tau_zeta} give $2\zeta+\tau\beta\mathfrak C_3\leq\varepsilon_0/2$. Hence, if both systems satisfy alternative \eqref{eq:testing_opt_2} on $B$, then, by \eqref{eq:improved_doubling} and the small boundary property \eqref{eq:def_small_boundary}, the union of the sets $U_B$ of the two systems and of the annulus $B\setminus(1-\tau)B$ has measure at most $\varepsilon_0\mu(B)/2$. Thus, by \eqref{eq:measure_testing_small_sets}, this union carries at most a fraction $(200\mathfrak C_1)^{-1}$ of the total variation of each of the two testing measures on $B$.
			
		}

	\vv
	
	For an explicit expression of $\beta_0$ in terms of the other structural parameters, we refer the reader to Remark \ref{remark:choice_of_beta}.
	
	Theorem \ref{theorem:main_theorem_hofmann_problem} applies, in particular, to testing functions with $L^p$ and $L^q$ size bounds. Let $1<p,q\leq\infty$, and, for a ball $B$, let $b_B\in L^p(\mu)$ and $b_B^*\in L^q(\mu)$ be supported on $B$ and such that \[\Bigl|\int_B b_B\,d\mu\Bigr|=\Bigl|\int_B b^*_B\,d\mu\Bigr|=\mu(B).\] 
	
	Multiplying the testing functions by constants of modulus one, we may assume that $\int_B b_B\,d\mu=\int_B b^*_B\,d\mu=\mu(B)$. Let $A\geq1$ be such that $\|b_B\|_{L^p(\mu)}\leq A\mu(B)^{1/p}$ and $\|b_B^*\|_{L^q(\mu)}\leq A\mu(B)^{1/q}$, with the usual interpretation when an exponent is infinite. By H\"older's inequality, every Borel set $E\subset B$ satisfies
	\[
\int_E	|b_B| \,d\mu \leq A\mu(B)\, \Bigl(\frac{\mu(E)}{\mu(B)}\Bigr)^{1/p'},
	\]
	and, moreover, $\mu(B)\leq\|b_B\mu\|\leq A\mu(B)$. The same estimates hold for $b_B^*$, with $q$ in place of $p$. Hence, the measures $\widehat\nu_B\coloneqq b_B\mu$ and $\widehat\nu_B^*\coloneqq b_B^*\mu$ satisfy \eqref{eq:measure_testing_size} and \eqref{eq:measure_testing_small_sets}, with the constants of Definition \ref{def:accretive_testing_system} given by
	\begin{equation}\label{eq:measure_testing_Lp_bridge}
		\mathfrak C_1=A\qquad \text{ and }\qquad 
		\varepsilon_0=(200A^2)^{-\max\{p',q'\}},
	\end{equation}
	where we understand that  $\infty'=1$.
	Finally, if for some $s,L>0$ we have
	\[
	\sup_{t>0}t^s\mu\bigl(\{x\in(1-\tau)B:T_\sharp(b_B\mu)(x)>t\}\bigr)\leq L\, \mu(B),
	\]
	then \eqref{eq:measure_testing_L0} holds, for every $\eta$, with level $(2L/\zeta)^{1/s}$, and the same holds for $b_B^*$ and $T^*$. 
	
	By Chebyshev's inequality, for any $p,q\in(1,\infty]$, Theorem \ref{theorem:main_theorem_hofmann_problem} gives the $L^2(\mu)$-boundedness of $T_\mu$ if, for every $B\in\mathscr B$, there are testing functions $b_B$ and $b^*_B$ as above such that
	\[
	\int_B (T_{\mu,\sharp} b_B)^{q'}\, d\mu+\int_B (T_{\mu,\sharp}^* b_B^*)^{p'}\, d\mu\lesssim\mu(B).
	\]
	
	More precisely, we obtain the following corollary.
	\begin{corollary}\label{corollary:hofmann_problem_maximal}
		Let $(X,d)$ be a geometrically doubling metric space with constant $D>1$, let $\mu$ be a Radon measure on $X$ with $n$-growth with constant $\mathfrak C_0>0$, for some $n\geq1$, and let $T_\mu$ be an $n$-dimensional singular integral operator associated with a Calder\'on-Zygmund kernel with constants $C>0$ and $\gamma\in(0,1)$. Let $1<p,q\leq\infty$, and let $\mathscr B$ be as in Theorem \ref{theorem:main_theorem_hofmann_problem}, for some $\delta\in(0,1]$, $\Lambda\geq100$, $\beta\geq\beta_0$, and $\mathfrak C_3\geq220$. Suppose that there are constants $\mathfrak C_1,\mathfrak C_2>0$ such that, for every $B\in\mathscr B$, there exist functions $b_B\in L^p(\mu)$ and $b^*_B\in L^q(\mu)$ supported on $B$ satisfying
		\[
		\begin{array}{@{}l@{\hspace{2.5cm}}l@{}}
			\displaystyle
			\Bigl|\int_B b_B \, d\mu\Bigr| = \mu(B), &
			\displaystyle
			\Bigl|\int_B b_B^* \, d\mu\Bigr| = \mu(B),
			\\[0.8em]
			\|b_B\|_{L^p(\mu)}\leq\mathfrak C_1\,\mu(B)^{1/p},
			&
			\|b^*_B\|_{L^q(\mu)}\leq\mathfrak C_1\,\mu(B)^{1/q},
			\\[0.8em]
			\displaystyle
			\int_B (T_{\mu,\sharp} b_B)^{q'} \, d\mu \leq \mathfrak C_2\, \mu(B),
			&
			\displaystyle
			\int_B (T_{\mu,\sharp}^* b_B^*)^{p'} \, d\mu \le \mathfrak C_2\, \mu(B),
		\end{array}
		\]
		with the convention $\infty'=1$. Then, the operator $T_\mu$ is bounded on $L^2(\mu)$.
	\end{corollary}
	
	\vv
	
	Combined with the Cotlar-type inequality of H. Martikainen, X. Tolsa, and the second named author \cite[Corollary 3.6]{MMT18}, Theorem \ref{theorem:main_theorem_hofmann_problem} also yields the following result on Hofmann's problem, see Section \ref{section:hystory}, for kernels that are not necessarily antisymmetric.
	
	\begin{corollary}\label{corollary:hofmann_problem}
		Let $(X,d)$ be a geometrically doubling metric space with constant $D>1$, let $\mu$ be a Radon measure on $X$ with $n$-growth with constant $\mathfrak C_0>0$, for some $n\geq1$, and let $T_\mu$ be an $n$-dimensional singular integral operator associated with a Calder\'on-Zygmund kernel with constants $C>0$ and $\gamma\in(0,1)$. Let $p,q\in(1,2]$, $\delta\in(0,1]$, and $\Lambda\geq100$. There exist $\beta\geq\beta_0$ and $\mathfrak C_3\geq220$, depending only on $D$, $n$, $\delta$, and $\Lambda$, such that the following holds, with $\mathscr B$ as in Theorem \ref{theorem:main_theorem_hofmann_problem}.
		
		Suppose that there are constants $\mathfrak C_1,\mathfrak C_2>0$ such that, for every $B\in\mathscr B$, there exist functions $b_B$ and $b^*_B$ supported on $B$ satisfying
		\[
		\begin{array}{@{}l@{\hspace{2.5cm}}l@{}}
			\displaystyle
			\int_B b_B \, d\mu = \mu(B), &
			\displaystyle
			\int_B b_B^* \, d\mu = \mu(B),
			\\[0.8em]
			\displaystyle
			\int_B |b_B|^p \, d\mu \le \mathfrak C_1\, \mu(B),
			&
			\displaystyle
			\int_B |b_B^*|^q \, d\mu \leq \mathfrak C_1\, \mu(B),
			\\[0.8em]
			\displaystyle
			\sup_{\varepsilon>0}\int_{B_1} |T_{\mu,\varepsilon} b_B|^{q'} \, d\mu \leq \mathfrak C_2\, \mu(B),
			&
			\displaystyle
			\sup_{\varepsilon>0}\int_{B_1} |T_{\mu,\varepsilon}^* b_B^*|^{p'} \, d\mu \le \mathfrak C_2\, \mu(B),
		\end{array}
		\]
		where
		\[
		B_1 \coloneqq
		\begin{cases}
			B &\qquad \qquad \text{ if } \frac{1}{p}+\frac{1}{q}<1+\frac{1}{n\max\{p,q\}},\\[4pt]
			2B &\qquad\qquad \text{otherwise}.
		\end{cases}
		\]
		Then, the operator $T_\mu$ is bounded on $L^2(\mu)$.
	\end{corollary}
	
	In particular, the unbuffered testing conditions, with $B_1=B$, suffice for all $p,q\in(1,2]$ if $n=1$, and for $p,q\in[2-\varepsilon(n),2]$, for some $\varepsilon(n)>0$, in general, while the buffered ones, with $B_1=2B$, suffice for all $p,q\in(1,2]$. These are the ranges obtained in \cite{MMT18} for measures with $n$-growth in $\R^d$ and antisymmetric kernels.
	
	As in \cite{MMT18}, Corollary \ref{corollary:hofmann_problem} follows from \cite[Corollary 3.6]{MMT18}, whose proof does not use the antisymmetry of the kernel. The argument there extends to balls in geometrically doubling metric spaces, with the dilated Hardy-Littlewood maximal function \eqref{eq:def_centered_maximal} in place of the Euclidean maximal functions, and gives, if $\beta$ and $\mathfrak C_3$ are large enough,
	\[
	\int_{(1-\tau)B}(T_{\mu,\sharp}b_B+T^*_{\mu,\sharp}b^*_B)\,d\mu\lesssim_\tau\mu(B)
	\]
	for every $B\in\mathscr B$ and every $\tau\in(0,1)$. Instead of the big pieces $Tb$ theorem, which requires antisymmetry, we then apply Theorem \ref{theorem:main_theorem_hofmann_problem}, together with \eqref{eq:measure_testing_Lp_bridge} and Chebyshev's inequality.
	
	We stress that neither the doubling property of $\mu$ nor the antisymmetry of $K$ is assumed in Theorem \ref{theorem:main_theorem_hofmann_problem} and Corollary \ref{corollary:hofmann_problem}. We discuss the main ideas of the proof in Section \ref{section:introduction_outline}.

	\begin{remark}\label{remark:range_n}
		We state Theorem \ref{theorem:main_theorem_hofmann_problem} and Corollary \ref{corollary:hofmann_problem} for $n\geq1$, which covers the application in Theorem \ref{theorem:application_caloric_measure}, and we assume $n\geq1$ throughout the paper. The same arguments apply to $0<n<1$, after replacing the exponent $n$ in the suppressed kernel \eqref{eq:suppressed_kernel} by $1$, so that Lemma \ref{lemma:kernel_suppressed} still follows from Appendix \ref{appendix:suppression_metric_spaces}, now with $M=1$.
	\end{remark}
	
	\begin{remark}\label{remark:upper_doubling}
		A natural extension of Theorem \ref{theorem:main_theorem_hofmann_problem} concerns measures satisfying the \textit{upper doubling} condition
		\begin{equation}\label{eq:def_upper_doubling}
			\mu(B(x,r))\leq \lambda(x,r), \qquad \text{ for all }x\in X  \text{ and all }r>0,
		\end{equation}
		where $\lambda(\cdot, \cdot): X  \times (0,\infty) \to (0,\infty)$ satisfies the following properties:
		\begin{enumerate}
			\item[(i)]  $\lambda (x, \cdot)$ is non-decreasing for any fixed $x \in X$,
			\item[(ii)]  $\lambda(x,2r)\le C_\lambda\,\lambda(x,r)$,  for all $x\in X, r>0$,
		\end{enumerate}
		and the associated \emph{Calder\'on-Zygmund kernel} $K\colon (X\times X)\setminus\{(x,x): x\in X\}\to\mathbb C$ is such that there exist $\gamma \in (0,1]$ and $C>0$ such that:
		\begin{itemize}
			\item For all $x,y\in X$, with $x\neq y$, we have
			\[
			|K(x,y)|\leq \frac{C}{\lambda(x,d(x,y))}.
			\]
			\item For all $y,x,x'\in X$ with $m\coloneqq \max\{d(y,x),d(y,x')\}$ and $2d(x,x')\leq m$, we have
			\[
			|K(x,y)-K(x',y)| + |K(y,x)-K(y,x')| \leq \frac{C}{\lambda(y,m)}\, \Bigl(\frac{d(x,x')}{m}\Bigr)^\gamma.
			\]
		\end{itemize}
		
		The proof in this paper is presented for $n$-growth measures on metric spaces. Although we do not extend the whole local $Tb$ argument to this setting, Appendix \ref{appendix:suppression_metric_spaces} contains the construction of suppressed kernels for upper doubling measures.
	\end{remark}

	\subsection{History and motivation}\label{section:hystory}
	The first version of the $T1$ theorem was proved by G. David and J.-L. Journé in \cite{DJ84}, and the classical reference for the $Tb$ theorem in homogeneous spaces is the subsequent work by G. David, J.-L. Journé, and S. Semmes in \cite{DJS85}. Their work was further extended in the very influential paper \cite{NTrV03}, in which F. Nazarov, S. Treil, and A. Volberg proved a (global) $Tb$ theorem in a non-homogeneous setting; see also \cite{NTrV99} and Volberg's monograph \cite{Vo03}.
	
	Nazarov-Treil-Volberg's articles are central to the literature also because of the techniques pioneered therein, such as, among others, the randomization of dyadic grids and the surgery arguments, the construction of suppressed kernels, and their associated SIOs, which have since become fundamental tools in a wide range of settings and now play a central role in non-homogeneous harmonic analysis and related areas.  Their $Tb$ theorem turned out to be fundamental in the solution of several deep problems in harmonic analysis such as X. Tolsa's proof of the semiadditivity of analytic capacity in \cite{To03} and T. Hyt\"onen's $A_2$ theorem \cite{Hy12}. For other versions of the $Tb$ theorem related to Vitushkin's conjecture, we also refer to the works of G. David  \cite{Da98} and G. David-P. Mattila \cite{DM00}.
	To highlight the role of $T1$/$Tb$ theorems in modern harmonic analysis, one can also see the bilinear methods by H. Martikainen and E. Vuorinen in \cite{MV21}, the characterization of Varopoulos' $\varepsilon$-approximability of bounded harmonic functions in terms of uniform rectifiability by J. Garnett, X. Tolsa, and the second named author in \cite{GMT18}, and the papers  \cite{Ma12a}, \cite{Ma12b}, and  \cite{HM14}. Further connections can be found in the interplay with PDEs and free-boundary problems such as the rectifiability of harmonic measure in \cite{AHMMMTV16} and its generalization to two-phase problems and to elliptic operators (\cite{AMT17a}, \cite{AMT20}, \cite{AMTV19}, \cite{MT20}, \cite{PPT21}, \cite{Pu19}, \cite{MMP25}).
	
	Moreover, there has been growing interest in extending several of the above results to more general settings, including parabolic and Kolmogorov spaces, sub-Riemannian geometry, and more general metric settings; see, for instance, \cite{BHHLN22}, \cite{BHHLN26}, \cite{He24}, \cite{HJ25}, \cite{LN21}, \cite{MPT22}, \cite{MP21}. For example, we note that, in \cite{HMP24}, the global $Tb$ theorem for a suitable non-antisymmetric SIO was the fundamental tool used by J. Hern\'andez, J. Mateu, and L. Prat to establish the semi-additivity of a symmetric version of the Lipschitz caloric capacity for the fractional heat operator $(-\Delta_x)^{1/2}+\partial_t$. Thus, establishing a local $Tb$ theorem in the generality of Theorem \ref{theorem:main_theorem_hofmann_problem} opens the way to a wide range of potential future applications and investigations, see also Section \ref{section_applications}.
	
	The first local $Tb$ theorem was proved by M. Christ in \cite[Main Theorem 10]{Ch90} in the doubling (homogeneous) setting. Namely, Christ established the $L^2(\mu)$-boundedness of a SIO by replacing the need for a globally defined testing condition with $L^\infty(\mu)$ testing assumptions on a family of functions $b_Q\in L^\infty(\mu)$ supported on cubes $Q\subset \R^d$ that satisfy the pseudo-accretivity condition
	\[	
	\biggl|\int_Q b_Q\, d\mu\biggr|\gtrsim \mu(Q),
	\]
	together with an analogous pseudo-accretive system for the adjoint operator.
	
	Christ's result was extended through a delicate use of random dyadic grids to the case of non-doubling measures by Nazarov, Treil, and Volberg in \cite{NTrV02}, assuming that $b_Q, b^*_Q\in L^\infty(\mu)$ and $Tb_Q, T^*b^*_Q$ satisfy suitable bounded $L^2$-mean oscillation-type conditions, uniformly in $Q$. The non-homogeneous $Tb$ theorem of \cite{NTrV03} has been further generalized to the framework of geometrically doubling metric measure spaces by T. Hyt\"onen and H. Martikainen in \cite{HM12}, using a metric variant of dyadic systems and an almost-covering of the space by balls. Both techniques are crucial for formulating Theorem \ref{theorem:main_theorem_hofmann_problem} in its full generality.
	
	The interest in lowering the integrability assumptions on the testing functions arose with the work of P. Auscher, S. Hofmann, C. Muscalu, T. Tao, and C. Thiele \cite{AHMTT02}, who proved a local $Tb$ theorem with local dyadic testing conditions on the pseudo-accretive system of the form
	\begin{equation}\label{eq:pq_testing}
		\int_Q |b_Q|^p+|b^*_Q|^q+ |Tb_Q|^{q'}+ |T^*b^*_Q|^{p'} \lesssim |Q|,
	\end{equation}
	in the model case of Calder\'on-Zygmund operators with perfect cancellation. The case of general SIOs is much more delicate. 
	
	S. Hofmann proved in the unpublished work \cite{Ho07} that, for the Lebesgue measure on $\R^d$, a sufficient condition for a SIO associated with a Calder\'on-Zygmund kernel as in Definition \ref{def:CZ_kernel} to be bounded on $L^2$ is the existence of pseudo-accretive systems of testing functions $\{b_Q\}$ and $\{b^*_Q\}$, with $b_Q$ and $b^*_Q$ supported on $Q$, such that, for every cube $Q\subset\R^d$,
	\begin{equation}\label{eq:pq_testing_hofmann}
		\int_Q |b_Q|^p+|b^*_Q|^p+ |Tb_Q|^{2}+ |T^*b^*_Q|^{2} \lesssim |Q|
	\end{equation}
	for some $p>2$.
	Hofmann obtained this result by building upon \cite{AHMTT02}, and his study was in turn motivated by applications to layer potentials in connection with $L^2$-solvability of boundary value problems, see \cite{AAAHK11}.
	
	The question of whether \eqref{eq:pq_testing} implies $L^2$ boundedness for arbitrary $p,q\in(1,\infty)$ is known as \textit{Hofmann's local $Tb$ problem}.
	
	P. Auscher and Q.X. Yang proved in \cite{AY09} the analogue of \cite{Ho07} for $p=2$ in the testing conditions \eqref{eq:pq_testing_hofmann}, by reducing the study to perfect dyadic operators, and hence invoking \cite{AHMTT02}. In spaces of homogeneous type, P. Auscher and E. Routin gave in \cite{AR13} a direct proof of the local $Tb$ theorem under the testing conditions \eqref{eq:pq_testing} in the case $1/p + 1/q\leq 1$, and they extended it to exponents with $1/p+1/q>1$, assuming an additional weak boundedness-type property related to Hardy inequalities. This article relied on the so-called \textit{twisted martingale transform}. Establishing a suitable variant of these transforms played a crucial role in the paper \cite{LM16}, where M. Lacey and H. Martikainen proved, in the non-homogeneous setting, a local $Tb$ theorem for $p=q=2$.

	Hofmann's problem has been solved by T. Hyt\"onen and F. Nazarov in \cite{HN20} for $p,q\in (1,\infty)$ under the ``buffered'' testing conditions on dilates of cubes $Q\subset \R^d$, and more specifically
	\[
	\int_Q |b_Q|^p + |b^*_Q|^q \lesssim |Q|\qquad\text{ and }\qquad \int_{2Q}|Tb_Q|^{q'}  +  |T^*b^*_Q|^{p'}\lesssim |Q|
	\]
	which, combined with a Cotlar-type lemma for ``buffered'' accretive systems (see \cite[Lemma 3.2]{HN20}), made it possible to control maximal truncations of the SIO. We remark that, by Hardy's inequalities, in the case of the Lebesgue measure and if $1/p+1/q\leq 1$, the unbuffered testing conditions \eqref{eq:pq_testing} imply the buffered ones. The buffered theorem in \cite{HN20} follows from their local $Tb$ theorem for maximal truncations \cite[Theorem 2.13]{HN20}, which allows uniform $L^0$ testing conditions and does not assume antisymmetry. Their extension to doubling metric measure spaces relies on Christ's dyadic cubes, see \cite[Section 6]{HN20}. 
	
	A different Cotlar-type approach was combined with the good lambda inequality by H. Martikainen, X. Tolsa, and the second named author in \cite{MMT18} to fully solve Hofmann's problem in $\mathbb R$ for $p,q\in (1,2]$ and, in $\mathbb R^d$ for $d>1$, for the range of exponents $p,q\in [2-\varepsilon (d), 2]$ with $\varepsilon(d)>0$.  Corollary \ref{corollary:hofmann_problem} extends the results of \cite{MMT18} for measures with $n$-growth in $\R^d$, which are proved for antisymmetric kernels, to kernels that are not necessarily antisymmetric and to geometrically doubling metric spaces: the unbuffered testing conditions give the $L^2$-boundedness for the range of exponents in \cite{MMT18}, and the buffered ones for all $p,q\in(1,2]$. 
	We do not know whether the unbuffered fixed-truncation theorem holds in this generality for the full range $1<p,q\leq\infty$. 
	Finally, we highlight that local $Tb$ theorems have also been studied for square functions; we refer to \cite{GH17}, \cite{LM17}, \cite{MM15}, \cite{MMV17}, \cite{MMV21}, and the discussion therein. For other related works, we also refer the reader to \cite{HV15}, \cite{LV16}, \cite{BKW21}, \cite{Vi24}, and the references therein. For further expositions on the methods in \cite{NTrV03}, we also refer to  \cite[Chapter 5]{To14}.
	
	\subsection{The Neumann problem}\label{section:intro_neumann}
	Let $n\geq2$. In this subsection and in Section \ref{section:neumann}, $\Omega\subset\mathbb R^{n+1}$ is a bounded domain, and we set  $\sigma\coloneqq\mathcal H^n|_{\partial \Omega}$, and $\delta_\Omega(X)\coloneqq\dist(X,\partial \Omega)$. We assume that there exist $C_A\geq1$ and $c_I\in(0,1)$ such that, for every $x\in \partial \Omega$ and $0<r\leq\diam \partial \Omega$,
	\begin{equation}\label{eq:neumann_ADR_corkscrew}
		C_A^{-1}\, r^n\leq\sigma(B(x,r))\leq C_A\, r^n\qquad\text{ and }\qquad B(X_{x,r},c_Ir)\subset\Omega\cap B(x,r)\ \text{ for some }X_{x,r}\in\Omega,
	\end{equation}
	that is, $\sigma$ is \textit{Ahlfors regular} and $\Omega$ satisfies the \textit{interior corkscrew condition}.
	
	For $a>1$ and a function $U$ on $\Omega$ with values in a normed space, we define the \textit{nontangential maximal function}
	\begin{equation}\label{eq:neumann_N}
		N_aU(x)\coloneqq\sup\bigl\{|U(X)|:X\in\Omega,\ |X-x|<a\,\delta_\Omega(X)\bigr\},\qquad \text{ for }x\in \partial\Omega,
	\end{equation}
	with the convention $\sup\varnothing\coloneqq0$. We say that $\partial \Omega$ is \textit{uniformly $n$-rectifiable} if there are $\theta>0$ and $M\geq1$ such that, for every $x\in \partial \Omega$ and $0<r\leq\diam \partial \Omega$, there is an $M$-Lipschitz map $\rho\colon\mathbb R^n\to\mathbb R^{n+1}$ with
	\[
	\mathcal H^n\bigl(\partial \Omega\cap B(x,r)\cap\rho(B_{\mathbb R^n}(0,r))\bigr)\geq\theta\, r^n.
	\]

	If $\Omega$ is a Lipschitz domain, with outer unit normal $\nu$, the \textit{$L^p$ Neumann problem} for the Laplacian asks, for every $h\in L^p(\sigma)$ with $\int_{\partial \Omega}h\,d\sigma=0$, for a harmonic function $u$ in $\Omega$ such that $\nabla u$ converges nontangentially $\sigma$-a.e.\ on $\partial \Omega$, with $\partial_\nu u=h$, and $\|N_a(\nabla u)\|_{L^p(\sigma)}\leq C\|h\|_{L^p(\sigma)}$. In rough domains, the outer normal and the boundary values of $\nabla u$ need not exist, and the boundary condition is formulated in the weak sense \eqref{eq:neumann_hypothesis} below which, by the divergence theorem, agrees with the classical one in Lipschitz domains, see \cite{MT24b, To26}.
	\begin{theorem}\label{theorem:neumann_UR}
		Let $\Omega\subset \mathbb R^{n+1}$, $n\geq 2$, be a bounded domain which satisfies the interior corkscrew condition and has Ahlfors regular boundary, see \eqref{eq:neumann_ADR_corkscrew}.
		Suppose that, for some $1<p<\infty$ and $a>1$, the $L^p$ Neumann problem is solvable, namely, for every bounded real Borel function $h$ on $\partial \Omega $ with $\int_{\partial \Omega}h\,d\sigma=0$ there exists a real harmonic function $u_h\in W^{1,2}(\Omega)$ such that
		\begin{equation}\label{eq:neumann_hypothesis}
			\int_\Omega\nabla u_h\cdot\nabla\varphi\,dX=\int_{\partial \Omega }h\,\varphi\,d\sigma,\qquad\text{ for every }\varphi\in C^\infty_c(\mathbb R^{n+1}),
		\end{equation}
		and
		\begin{equation}\label{eq:neumann_estimate}
			\|N_a(\nabla u_h)\|_{L^p(\sigma)}\leq C_N\,\|h\|_{L^p(\sigma)}.
		\end{equation}
		Then $\partial \Omega$ is uniformly $n$-rectifiable, with constants $\theta$ and $M$ depending only on $n$, $C_A$, $c_I$, $p$, $a$, and $C_N$.
	\end{theorem}

	We remark that the weak formulation \eqref{eq:neumann_hypothesis} is imposed on the whole topological boundary, and that it involves no normal, no trace, and no density or extension property of $W^{1,2}(\Omega)$. Moreover, we only use one solution for each datum, so that neither uniqueness nor a linear choice of the solutions is needed, and we do not assume any Harnack chain condition on $\Omega$. The hypotheses of Theorem \ref{theorem:neumann_UR} hold if the Neumann problem for the Laplacian is solvable in $L^p$ in the sense of \cite{MT24b}, where the solutions are the variational ones and the (modified) nontangential maximal function $\widetilde N_a$ is defined with $L^2$ averages on Whitney balls, since the subharmonicity of $|\nabla u_h|^2$ gives $N_a(\nabla u_h)\leq\widetilde N_a(\nabla u_h)$.

	The $L^p$ Neumann problem for the Laplacian in Lipschitz domains was solved by D. Jerison and C. Kenig \cite{JK81} for $p=2$, by means of a Rellich identity, which shows that the normal and the tangential derivatives of a harmonic function have comparable $L^2(\sigma)$ norms. G. Verchota \cite{Ve84} obtained the solution by the method of layer potentials, and B. Dahlberg and C. Kenig \cite{DK87} solved the problem for $1<p<2+\varepsilon$, with data in the Hardy space $H^1$ at the endpoint $p=1$. This range of exponents is sharp \cite{DK87}, and we refer to \cite{KS08} for necessary and sufficient conditions for $p>2$ and for convex domains, to \cite{Br95} for data in Hardy spaces $H^p$ with $p<1$, and to the monograph \cite{Ke94}. For divergence form elliptic operators with real coefficients in Lipschitz domains, we refer to the works of C. Kenig and J. Pipher \cite{KP93, KP95}, of C. Kenig and D. Rule \cite{KR09}, of M. Dindo\v{s}, J. Pipher, and D. Rule \cite{DPR17}, and of M. Dindo\v{s}, S. Hofmann, and J. Pipher \cite{DHP23}, and, for complex coefficients and elliptic systems, to \cite{AAH08, AAM10, AA11, AAAHK11, AMM13, HMM15, Sh07} and the references therein.

	Beyond Lipschitz domains, much less is known. In planar chord-arc domains, the solvability of the Neumann problem in $L^p$ is equivalent to that of the regularity problem in $L^p$ and of the Dirichlet problem in $L^{p'}$, see \cite{JK82}. In higher dimensions, S. Hofmann, M. Mitrea, and M. Taylor \cite{HMT10} showed that, for every $1<p<\infty$, the Neumann problem is solvable in $L^p$ in regular Semmes-Kenig-Toro domains whose outer unit normal has sufficiently small $\mathrm{BMO}$ norm, by means of the invertibility of boundary layer potentials. This approach has been developed into a rich Calder\'on-Zygmund theory for layer potentials on uniformly rectifiable sets, which gives the solvability of the Dirichlet, regularity, Neumann, and transmission problems for large classes of weakly elliptic systems in Ahlfors regular domains whose outer unit normal has small mean oscillation, with data in Lebesgue, Hardy, Morrey, Herz-type, and other function spaces, see \cite{MMMMM22, MMM23a, MMM23, LM24, MiTa26, Gu25} and the references therein. 
	
	In every bounded chord-arc domain, the Dirichlet problem is solvable in $L^{p'}$ for some $1<p'<\infty$ \cite{DJ90, Se90}, and the regularity problem in $L^p$ for some $1<p<\infty$ \cite{MT24}. Whether the same holds for the Neumann problem, if $n\geq2$, is an open problem posed by C. Kenig \cite[Problem 3.2.2]{Ke94}, see also \cite[Question 2.5]{Tor10} and the survey \cite{To26}. One of the difficulties is that, unlike the Dirichlet problem, the Neumann problem has no maximum principle, which prevents transferring solvability from subdomains. Recently, J. Feneuil and L. Li \cite{FL24} and S. Hofmann and D. Sparrius \cite{HS25} studied the Neumann function and the Poisson-Neumann problem in chord-arc domains, and extended the extrapolation of solvability of \cite[Theorem 6.3]{KP93}. Moreover, the second named author and X. Tolsa \cite{MT24b} proved the solvability of the Neumann problem in $L^p$, for $1<p<2$, in chord-arc domains with very big pieces of superdomains on which the problem is uniformly solvable, under additional assumptions on the regularity problem and on the boundary.

	For the Dirichlet and regularity problems, the role of the geometry of the boundary is much better understood. In domains with Ahlfors regular boundary satisfying the interior corkscrew condition, the solvability of the Dirichlet problem in $L^{p'}$, for some $1<p'<\infty$, or, equivalently, the weak-$A_\infty$ property of harmonic measure, implies the uniform rectifiability of the boundary \cite{HMU14, HLMN17, MT20}, and these properties are characterized geometrically in \cite{AHMMT20}. In the same class of domains, the solvability of the regularity problem in $L^p$ implies that of the Dirichlet problem in $L^{p'}$ \cite{MT24}, see also \cite[Theorem 5.4]{KP93} for Lipschitz domains, and hence the uniform rectifiability of the boundary. We refer to \cite{MT24, MPT25, GMT25} for the converse implication and further results, and to \cite{Ga25} for one-sided Rellich inequalities and uniform rectifiability. Moreover, the first named author recently proved the WALA conjecture in \cite{Me26}, that is, that Ahlfors regular sets satisfying the weak approximation of Lipschitz functions by affine functions are uniformly rectifiable, and, in a forthcoming part of this work, joint with the second named author, this is applied to show that, if $\partial\Omega$ is Ahlfors regular, the solvability of the regularity problem for uniformly elliptic operators in divergence form with bounded coefficients implies the uniform rectifiability of $\partial\Omega$. Theorem \ref{theorem:neumann_UR} is the counterpart of these results for the Neumann problem.
	To the best of our knowledge, this is the first result in which the solvability of the Neumann problem yields geometric information on the boundary.
	\vv

	We briefly describe the proof of Theorem \ref{theorem:neumann_UR}, which is given in Section \ref{section:neumann}. We work with the real Clifford algebra $\mathrm{Cl}_{n+1}$ generated by $e_1,\dots,e_{n+1}$, with $e_j^2=-1$, and with the \textit{Cauchy-Clifford operator} $\mathcal C_\sigma$, that is, the singular integral operator on $\partial \Omega$ with kernel $-\omega_{n}^{-1}(x-y)|x-y|^{-n-1}$, where $\omega_{n}\coloneqq\mathcal H^{n}(\mathbb S^{n})$, acting by left multiplication on $\mathrm{Cl}_{n+1}$-valued functions, whose components are the Riesz transforms, see Section \ref{subsection:clifford}. 
	
	Clifford analysis provides a natural framework for the higher-dimensional Cauchy integral, whose $L^2$-boundedness on Lipschitz surfaces follows from \cite{CMM82} by the method of rotations. This boundedness was also studied using Clifford analysis in \cite{Mur85, Mc89, LMS92}; see also \cite{Se89, GLQ93, Mi94}. S. Hofmann, M. Mitrea, and M. Taylor	\cite{HMT10} used the Cauchy-Clifford operator to study boundary value problems on uniformly rectifiable domains; see also \cite{MMMMM22, MMM23}. 
	
	For every ball $B$ centered on $\partial \Omega$, we use the hypotheses of Theorem \ref{theorem:neumann_UR} to construct a Neumann solution $u_B$, whose datum equals $1$ on $2B$ up to a compensating term, and a density $b_B=g_B+H_B$, where $g_B$ is a real function equal to $1$ on $2B$ and $H_B$ is bivector-valued, such that $\mathcal C_\sigma b_B=\nabla u_B$ in $\Omega$. In the smooth case, $b_B$ is the density of the Cauchy reproducing formula, whose scalar part is the Neumann datum and whose bivector part consists of tangential derivatives, see Section \ref{subsection:neumann_proof}. Since normals and traces are not available, we construct $H_B$ from the distributional curl of the gradient of a Neumann solution extended by zero outside $\Omega$. Then the Neumann estimate \eqref{eq:neumann_estimate} controls the size of $b_B$ and, with a Cotlar-type inequality, the maximal truncations of $\mathcal C_\sigma b_B$ in $L^p(\sigma)$. Since left multiplication by a bivector is skew-adjoint, the testing function $\mathfrak b_B\coloneqq\chi_BL_{b_B}$ is the identity plus a skew-adjoint operator on $B$, which is the accretivity condition of Theorem \ref{theorem:matrix_Tb}, a version of Theorem \ref{theorem:main_theorem_hofmann_problem} for singular integral operators with matrix-valued kernels and testing functions. This gives the $L^2(\sigma)$-boundedness of $\mathcal C_\sigma$, and hence of the Riesz transform, and Theorem \ref{theorem:neumann_UR} follows from the solution of the codimension one David-Semmes problem by F. Nazarov, X. Tolsa, and A. Volberg in \cite{NToV14}.
	
	Our proof is related to the first-order approach of P. Auscher, A. Axelsson, and A. McIntosh \cite{AAM10}, see also \cite{AAH08, AA11, AMM13} and the references therein, in which a divergence form elliptic equation is rewritten as a first-order system for the conormal gradient, and the Neumann and the regularity problems prescribe, respectively, the normal and the tangential components of its boundary trace. For the Laplacian, this system says that $\mathsf D\nabla u=0$ for every harmonic function $u$, where $\mathsf D$ is the Dirac operator, see Section \ref{subsection:clifford}. Then, the Cauchy reproducing formula recovers $\nabla u$ from its boundary values. In our proof, the normal and the tangential information is encoded in the scalar and the bivector parts of the density $b_B$, which plays the role of this boundary trace. However, we do not use this structure to solve a boundary value problem. Since neither normals nor traces are available, we construct $b_B$ in the sense of distributions, and we use the Neumann estimate to show that $\chi_BL_{b_B}$ is a local testing function for $\mathcal C_\sigma$, whose $L^2(\sigma)$-boundedness then follows from Theorem \ref{theorem:matrix_Tb}.
	
	Our construction of the testing functions is also related to the work of A. Grau de la Herr\'an and S. Hofmann \cite{GH17} on local $Tb$ theorems for square functions. In their application to layer potentials for elliptic operators with $t$-independent coefficients in the half-space, the testing functions combine the conormal derivative and the tangential gradient of a difference of fundamental solutions (with poles at interior and exterior corkscrew points associated with the testing cubes), and Green's formula gives the testing estimates. Their vector-valued $Tb$ theorem requires the accretivity of the conormal component and the smallness of a weighted average of the tangential one. Our testing functions also combine normal and tangential information from a single solution, but there are two crucial differences. First, we construct them from weak solutions of the Neumann problem on a boundary which is not known to be rectifiable, and we obtain the tangential part as a boundary density by a distributional argument, without using traces of the gradient. Second, the tangential part is the bivector $H_B$, and the symmetric part of $L_{b_B}$ is $g_BI$, since $L_{H_B}$ is skew-adjoint. Hence, the Neumann datum alone gives the accretivity of $\mathfrak b_B$, without any smallness condition on $H_B$. Moreover, we use these testing functions to control the maximal truncations of the Cauchy-Clifford operator, and hence to deduce the uniform rectifiability of $\partial\Omega$ from the solvability of the Neumann problem, rather than to obtain square function bounds.

	\vv
			We now make further remarks on the role of Theorem \ref{theorem:main_theorem_hofmann_problem}. The Neumann estimate controls $\mathcal C_\sigma b_B$, which combines the Riesz transforms of $g_B$ and of the coefficients of $H_B$, but not the Riesz transform of $g_B$ alone. Hence, scalar testing functions are not available, and neither the scalar local $Tb$ theorems of Section \ref{section:hystory}, including Theorem \ref{theorem:main_theorem_hofmann_problem}, nor the $Tb$ theorems for operator-valued kernels of \cite{Hy06, HV15}, whose testing functions are scalar, apply directly. Vector-valued and matrix-valued testing functions appear in $Tb$ arguments for square functions \cite{AHLMT02, AKM06, GH17} and in the global Clifford-valued $Tb$ theorem of \cite{GLQ93} but, to the best of our knowledge, no local $Tb$ theorem for SIOs with matrix-valued testing functions is available. Moreover, although the measure is Ahlfors regular and the kernel is antisymmetric, the testing functions are only in $L^p$, where $p<2$ is the natural range, see \cite{DK87, KP93, FL24, HS25}. 
			
			Thus, a matrix-valued version of the local $Tb$ theorems with $L^2$ testing conditions of \cite{AR13, LM16} would suffice if $p\geq2$. On the other hand, for $p<2$, Theorem \ref{theorem:neumann_UR} could also be obtained via a matrix-valued version of the results of \cite{MMV21}, or of those of \cite{HN20} for Ahlfors regular measures, which allow testing conditions on the maximal truncations below the duality exponents.
			The required local $Tb$ theorem, Theorem \ref{theorem:matrix_Tb}, follows from the proof of Theorem \ref{theorem:main_theorem_hofmann_problem} with only a few modifications, see Appendix \ref{appendix:matrix_Tb}, which also shows the robustness of our method.
			In particular, it may be possible that matrix-valued local $Tb$ theorems which would suffice to prove Theorem \ref{theorem:neumann_UR} can be obtained by adapting the proofs in \cite{HN20, MMV21}. However, we do not pursue this in the present paper and we believe that Theorem \ref{theorem:matrix_Tb} is of independent interest.
	
	\vv

	\subsection{Applications to the heat equation} \label{section_applications}

		Our second application concerns absolute continuity of caloric measure and gradients of single layer potentials for the heat equation.
		
		We denote by $\mathbb P^n\coloneqq\mathbb R^n\times\mathbb R$ the parabolic space, with points $\mathbf X=(X,t)$ and parabolic norm $\|\mathbf X\|\coloneqq\max\{|X|,|t|^{1/2}\}$. With the distance induced by  $\|\cdot\|$, $\mathbb P^n$ is a geometrically doubling metric space of homogeneous dimension $(n+2)$, whose open balls are \[B(\mathbf X,r)=\bigl\{(Y,s):|Y-X|<r,\ |s-t|<r^2\bigr\}, \qquad \mathbf X\in \mathbb P^n, r>0.\]
		For $\mathbf X, \mathbf Y\in \mathbb P^n$, $\mathbf X\neq \mathbf Y$, we denote by 
		\[
		\Gamma(\mathbf X,\mathbf Y)\coloneqq\frac{\chi_{\{t>s\}}}{(4\pi(t-s))^{n/2}}\exp\Bigl(-\frac{|X-Y|^2}{4(t-s)}\Bigr)
		\]
		the fundamental solution of the heat equation $\partial_t u- \Delta_Xu=0$, and let $K(\mathbf X,\mathbf Y)\coloneqq\nabla_X\Gamma(\mathbf X,\mathbf Y)$, where $\nabla_X$ denotes the spatial gradient. Each component of $K$ is an $(n+1)$-dimensional Calder\'on-Zygmund kernel for the parabolic distance, with regularity exponent $1$, see \cite[Lemma 5.4]{MPT22}. If $\sigma$ is a measure on $\mathbb P^n$ with $(n+1)$-growth, we consider the associated operator
		\[
		T_\sigma f(\mathbf X)=\int K(\mathbf X,\mathbf Y)f(\mathbf Y)\, d\sigma(\mathbf Y), \qquad f\in L^1(\sigma),
		\]
		to be interpreted in the usual sense of truncations, which corresponds to the \textit{spatial gradient of the single layer potential} for the heat equation.
		
		Let $\Omega\subset\mathbb P^n$ be open and connected, and let $\sigma$ be a non-negative Radon measure supported on $\partial\Omega$ with $(n+1)$-growth. For $\mathbf P\in\Omega$, we denote by $\omega^{\mathbf P}$ and $\omega_*^{\mathbf P}$ the caloric and adjoint caloric measures of $\Omega$ with pole at $\mathbf P$, see Section \ref{section:application_caloric_measure}. We assume the following \textit{$A_\infty$-type} hypotheses, which are quantitative and time-directed versions of the absolute continuity of $\omega^{\mathbf P}$ and $\omega_*^{\mathbf P}$ with respect to $\sigma$. As in \cite{MT20}, we formulate them in terms of measures, rather than densities.
		\begin{definition}[The classes $\mathcal A_\infty^+(\sigma;\kappa,\kappa')$ and $\mathcal A_\infty^-(\sigma;\kappa,\kappa')$]\label{definition:weak_A_infty}
			Fix $\delta,\vartheta,\kappa,\kappa'\in(0,1)$, $\Lambda\geq100$, $\mathfrak C_3\geq220$, and a sufficiently large structural constant $\beta$.
			\begin{itemize}
				\item We write $\omega\in\mathcal A_\infty^+(\sigma;\kappa,\kappa')$ if the following holds for every ball $B=B(\mathbf X_0,r)$, $\mathbf X_0=(X_0, t_0)$, centered on $\supp\sigma$, with $0<r\leq\diam(\supp\sigma)$, $\mathfrak C_3$-small boundary and $\sigma(\Lambda B)\leq\beta\sigma(\delta B)$: there exists $\mathbf P_B^+=(X_{P^+},t_{P^+})\in\Omega\cap\delta B$ satisfying
				\begin{equation}\label{eq:time_forward_point}
					t_0+\vartheta(\delta r)^2\leq t_{P^+}\leq t_0+(\delta r)^2,
				\end{equation}
				such that, for every Borel set $E\subset B$, the following property holds:
				\begin{equation}\label{eq:A_infty_type_condition}
					\text{If }\sigma(E)\leq\kappa\, \sigma(B),\qquad\text{ then }\qquad \omega^{\mathbf P_B^+}(E)\leq\kappa'\, \omega^{\mathbf P_B^+}(B).
				\end{equation}
				\item The property $\omega_*\in\mathcal A_\infty^-(\sigma;\kappa,\kappa')$ is defined by the analogous conditions, with $\omega_*^{\mathbf P_B^-}$ in place of $\omega^{\mathbf P_B^+}$ in \eqref{eq:A_infty_type_condition}, and poles $\mathbf P_B^-\in\Omega\cap\delta B$ satisfying
				\begin{equation}\label{eq:time_backward_point}
					t_0-(\delta r)^2\leq t_{P^-}\leq t_0-\vartheta(\delta r)^2.
				\end{equation}
			\end{itemize}
		\end{definition}
		
		We remark that, unlike reverse H\"older conditions, Definition \ref{definition:weak_A_infty} does not require $\omega^{\mathbf P_B^+}|_B$ and $\omega_*^{\mathbf P_B^-}|_B$ to be absolutely continuous with respect to $\sigma$. On the other hand, if $\omega^{\mathbf P_B^+}|_B\ll\sigma|_B$ and its density $k_B^+$ satisfies the reverse H\"older-type	 inequality
		\[
		\Bigl(\avint_B(k_B^+)^{s}\,d\sigma\Bigr)^{1/s}\leq C\, \frac{\omega^{\mathbf P_B^+}(B)}{\sigma(B)},
		\]
		for some $s>1$ and $C\geq1$, then H\"older's inequality gives 
		\[
		\omega^{\mathbf P_B^+}(E)\leq C\,(\sigma(E)/\sigma(B))^{1/s'}\,\omega^{\mathbf P_B^+}(B), \qquad \text{ for every Borel set }E\subset B.
		\]
		Hence, \eqref{eq:A_infty_type_condition} holds for every $\kappa'\in(0,1)$, with $\kappa\coloneqq(\kappa'/C)^{s'}$. The same applies to $\omega_*^{\mathbf P_B^-}$.
		\vv
		
	In corkscrew domains with Ahlfors-regular boundaries, suitable reverse	H\"older estimates and $A_\infty$-type conditions for elliptic measure	with respect to surface measure are equivalent to the solvability of the Dirichlet problem with data in $L^p$ for some $1<p<\infty$.
	This equivalence holds for second-order uniformly elliptic operators in divergence form with real $L^\infty$ coefficients; see \cite{MT24}.
	For the Laplacian, these conditions also imply uniform rectifiability of the boundary; see \cite{HLMN17,MT20} and the references therein.
	
	 In the parabolic setting, E. Fabes and S. Salsa showed in \cite{FS83} that, in Lipschitz cylinders, the caloric Poisson kernel satisfies a reverse H\"older inequality with exponent $2$, but for time-varying domains the picture is more delicate. R. Kaufman and J.-M. Wu showed in \cite{KW88} that caloric measure may fail to be absolutely continuous with respect to surface measure on the graph of a $\mathrm{Lip}(1,1/2)$ function, while J. L. Lewis and M. A. M. Murray proved in \cite{LM95} that it is an $A_\infty$ weight if the graph has, in addition, a half-order time derivative in parabolic $\mathrm{BMO}$, see also \cite{HL96}. These graphs are the model for the parabolic uniformly rectifiable sets introduced by S. Hofmann, J. L. Lewis, and K. Nystr\"om in \cite{HLN03, HLN04}. More recently, A. Genschaw and S. Hofmann proved in \cite{GeH20} a weak reverse H\"older inequality for caloric measure, equivalent to the $L^p$ solvability of the Dirichlet problem, assuming a time-backward Ahlfors regularity condition on the boundary, and S. Bortz, S. Hofmann, J. M. Martell, and K. Nystr\"om showed in \cite{BHMN23, BHMN25} that, for domains with time-symmetric Ahlfors regular boundary satisfying an interior corkscrew condition, this property implies parabolic uniform rectifiability.
		
		Definition \ref{definition:weak_A_infty} is a local version of these conditions, where the poles $\mathbf P_B^+$ and $\mathbf P_B^-$ lie in $\delta B$, in the future and in the past of the center of $B$, respectively, at time distance comparable to $(\delta r(B))^2$. We remark that we work under minimal geometric assumptions. In particular, we do not assume that $\sigma$ is Ahlfors regular, nor any corkscrew, Harnack chain, or capacity density condition on $\Omega$.
		
		\vv
		
		Time-directed reverse H\"older and $A_\infty$-type conditions also appear in the theory of parabolic Muckenhoupt weights, introduced by J. Kinnunen and O. Saari in \cite{KS16a, KS16b} in connection with the doubly nonlinear equation $\partial_t(|u|^{p-2}u)-\mathrm{div}(|\nabla u|^{p-2}\nabla u)=0$, which gives the heat equation for $p=2$, see \cite{KMYZ24, KM24, KM25} and also \cite{Saw86, Be11}
		and the references therein.
		There, $w$ is a weight on $\mathbb P^n$ with respect to Lebesgue measure, and the parabolic reverse H\"older inequality $(\avint_{R^-}w^{s})^{1/s}\leq C\avint_{R^+}w$ controls the integrability of $w$ on the lower part $R^-$ of a parabolic rectangle $R$ by its average on the upper part $R^+$, possibly separated from $R^-$ by a time lag. For instance, by the weak Harnack inequality, positive supersolutions of the heat equation in $\mathbb P^n$ are parabolic Muckenhoupt weights, see \cite{KS16a}.
		In Definition \ref{definition:weak_A_infty}, instead, both sides of \eqref{eq:A_infty_type_condition} involve the same ball $B$, and the time direction is encoded in the position of the pole $\mathbf P_B^+$, on which the measure $\omega^{\mathbf P_B^+}$ depends, so that $\mathcal A_\infty^+(\sigma;\kappa,\kappa')$ is not a condition on a single weight. Moreover, the theory of parabolic Muckenhoupt weights relies on the doubling property of Lebesgue measure, or of the underlying measure on a space of homogeneous type as in \cite{KMYZ24}, through Calder\'on-Zygmund decompositions and chaining arguments, which yield the self-improvement of parabolic reverse H\"older inequalities and several characterizations of the parabolic $A_\infty$ and reverse H\"older classes, see \cite{KM24, KM25}. Since $\sigma$ only has $(n+1)$-growth, we require \eqref{eq:A_infty_type_condition} only on balls with small boundary such that $\sigma(\Lambda B)\leq\beta\sigma(\delta B)$, with fixed parameters $\kappa$ and $\kappa'$, and we do not use any self-improvement property.
	
	\begin{theorem}\label{theorem:application_caloric_measure}
		Let $n\geq1$, let $\Omega\subset\mathbb P^n$ be open and connected, and let $\sigma$ be a non-negative Radon measure supported on $\partial\Omega$ such that $\sigma(B(\mathbf X,r))\leq\mathfrak C_0r^{n+1}$ for every $\mathbf X\in\mathbb P^n$ and $r>0$. Suppose that $\omega\in\mathcal A_\infty^+(\sigma;\kappa,\kappa')$ and $\omega_*\in\mathcal A_\infty^-(\sigma;\kappa,\kappa')$, for some $\kappa\in(0,1)$ and $\kappa'\in(0,1/200]$, where $\delta$ is sufficiently small depending only on $n$ and $\vartheta$, and $\beta$ is sufficiently large depending only on $n$, $\delta$, and $\Lambda$.
		
		Then $T_\sigma$ is bounded on $L^2(\sigma)$, that is,
		\[
		\sup_{\varepsilon>0}\|T_{\sigma,\varepsilon}\|_{L^2(\sigma)\to L^2(\sigma;\mathbb R^n)}<\infty.
		\] 
		The bound depends only on $n$, $\mathfrak C_0$, and the parameters in Definition \ref{definition:weak_A_infty}.
	\end{theorem}
	
	Theorem \ref{theorem:application_caloric_measure} is a parabolic analogue of the result \cite[Theorem 1.1]{MT20} by the second named author and X. Tolsa, where an $A_\infty$-type condition on harmonic measure, with poles depending on the ball, is shown to imply the $L^2$-boundedness of the Riesz transform with respect to measures with $n$-growth in $\mathbb R^{n+1}$, for arbitrary open sets. In \cite{MT20}, either the measure concentrates near the pole, or the $T1$ theorem applied to the normalized harmonic measure $\sigma(B)\,\omega^{x_B}|_B/\omega^{x_B}(B)$ gives a big piece of the ball on which the Riesz transform is bounded, and the global bound follows by a good-lambda argument. Instead, we test directly with the normalized caloric and adjoint caloric measures in Theorem \ref{theorem:main_theorem_hofmann_problem}.
	Since the testing measures are non-negative, \eqref{eq:measure_testing_size} holds with $\mathfrak C_1=1$, and the conditions in Definition \ref{definition:weak_A_infty} with $\kappa'\leq1/200$ give the small-set condition \eqref{eq:measure_testing_small_sets}. In \cite{MT20}, the analogous condition is assumed with an arbitrary $\kappa'\in(0,1)$. Moreover, a parabolic version of the estimates in \cite{MT20}, based on local boundedness of the gradient of a positive adjoint caloric function  and on a parabolic Bourgain-type lower bound for caloric measure, controls the maximal truncations by a maximal function of the testing measures away from the poles. The concentration alternative \eqref{eq:testing_opt_1} takes care of the balls on which $\sigma$ concentrates near the poles. Since caloric measure only sees the past of its pole, the measures with poles in the future of $B$ test $T^*$, those with poles in the past test $T$, and we work with balls whose time-forward and time-backward halves have comparable measure, see Lemma \ref{lemma:caloric_auxiliary_ball}.
	
	Finally, we remark that Theorem \ref{theorem:application_caloric_measure} is linked to \textit{parabolic uniform rectifiability}.
	Indeed, in the elliptic setting, the result in \cite{MT20}, combined with the solution of the codimension one David-Semmes problem by F. Nazarov, X. Tolsa, and A. Volberg \cite{NToV14}, shows that, if $\sigma$ is Ahlfors regular, the $A_\infty$-type condition on harmonic measure implies the uniform rectifiability of $\supp\sigma$, see also \cite{HLMN17}. A parabolic analogue of this argument would require that, for an $(n+1)$-dimensional Ahlfors regular measure $\sigma$ on $\mathbb P^n$, the $L^2(\sigma)$-boundedness of $T_\sigma$ implies the parabolic uniform rectifiability of $\supp\sigma$. The converse implication holds, see \cite{HJ25} and the references therein. However, J. Hoffman and B. Jaye showed in \cite{HJ25} that this direct analogue of the David-Semmes problem fails, by constructing measures that are not parabolic uniformly rectifiable, on which $T_\sigma$ and a large class of parabolic singular integrals are bounded. They also proved that the boundedness of all parabolic Calder\'on-Zygmund operators implies parabolic uniform rectifiability under an additional packing condition that rules out these examples, and they conjectured that the boundedness of $T_\sigma$ suffices under the same condition. A positive answer to their conjecture, combined with Theorem \ref{theorem:application_caloric_measure}, would show that, under this condition, the assumptions of Theorem \ref{theorem:application_caloric_measure} imply parabolic uniform rectifiability, see also Section \ref{sec:related_future}.

	\vv
	
	\subsection{Discussion of the proof, main ideas, and outline of the paper}\label{section:introduction_outline}
	
	In the case of a singular integral operator that is not necessarily antisymmetric, the natural assumptions for a local $Tb$ theorem involve two distinct testing systems, one associated with $T$ and one associated with its adjoint $T^*$, and this difficulty is substantial rather than technical. Indeed, in our proof of Theorem \ref{theorem:main_theorem_hofmann_problem}, each of the families $\{\widehat\nu_B\}_B$ and $\{\widehat\nu^*_B\}_B$ generates its own exceptional sets, suppressions, martingale difference decompositions, and corona decomposition. These objects act independently, so the one-measure techniques of F. Nazarov, S. Treil, and A. Volberg \cite{NTrV99, NTrV03} cannot be applied directly. The approach to non-homogeneous local $Tb$ theorems developed in \cite{MMV21} and \cite{MMT18}, which combines a non-homogeneous good-lambda method, the big pieces $Tb$ theorem of \cite{NTrV99}, and a change of measure argument, was implemented there for antisymmetric kernels, and it was pointed out in \cite{MMV21} that the general case involves genuine difficulties. The main achievement of our proof is to make these two a priori separate structures interact. To this end, we develop a framework that brings them into a single tree-wise analysis by means of suppressed singular integral operators and a \textit{mixed corona decomposition}. Its key ingredient is a tree-wise $L^2$-estimate with respect to two different measures.
	
	We first describe the overall strategy. After regularizing the kernel at small and large scales, which makes all the operators a priori bounded and allows us to absorb small multiples of their norms, see Section \ref{subsection:bounded_testing_reduction}, the proof consists of two steps. In the first one, which is the core of the paper, we prove that, for every ball $B_0\in\mathscr B$ and every $g\in L^\infty(\mu)$ supported in $B_0$,
	\[
	\int_{B_0}\bigl(T_{\mu,\sharp}g+(T^*)_{\mu,\sharp}g\bigr)\,d\mu\lesssim\|g\|_{L^\infty(\mu)}\,\mu(B_0),
	\]
	see \eqref{eq:corollary49_maximal_sec8}. This estimate resembles a local $T1$ theorem, but it only provides $L^1$ bounds for $T_{\mu,\sharp}\chi_{B_0}$ and $(T^*)_{\mu,\sharp}\chi_{B_0}$, whereas the classical local $T1$ theorems require $L^\infty$ or $\mathrm{BMO}$-type bounds.
	
	At this point, the generality of our testing assumptions plays a role in the proof itself. By Chebyshev's inequality, the measures $\mu|_B$ satisfy the $L^0$ condition \eqref{eq:measure_testing_L0}, and, in the second step, we run our construction once more with these testing measures. In this way, we find, in every doubling David-Mattila cell $Q$, a set $\mathcal G_Q$ with $\mu(\mathcal G_Q)\gtrsim\mu(Q)$ such that $T_{\mu|_{\mathcal G_Q}}$ is bounded on $L^2(\mu|_{\mathcal G_Q})$, and a good-lambda inequality in the spirit of \cite[Section 9]{MT20}, see Proposition \ref{theorem:MMT_theorem_A1}, upgrades these big pieces to the full $L^2(\mu)$-boundedness of $T$.
	
	\vv
	Our techniques involve three different geometric objects, and part of the difficulty is to make them compatible with one another.
	\begin{itemize}
		\item \textit{Metric balls.} The testing conditions are formulated on balls, and we work with the testing balls $\widetilde B(R)\in\mathscr B$ associated with doubling David-Mattila cells $R$, which have small boundary, see Appendix \ref{appendix:DM_testing_balls}. All the exceptional families associated with density, growth, and the operator are defined by conditions on balls, and the weak boundedness property is proved on random balls with small boundary. This allows us to carry out the whole argument in geometrically doubling metric spaces.
		\item \textit{Random dyadic systems} $\bbD(w)$ of T. Hyt\"onen and A. Kairema \cite{HK12}, which depend on a parameter $w$ in a probability space. We use two independent such systems, one for each testing system, to detect the cubes on which the testing measures are not accretive, and to implement the martingale difference decompositions and the probabilistic arguments of Sections \ref{section:suppressed_mixed_good} and \ref{sec:main_theorem_antisym}.
		\item \textit{A fixed David-Mattila lattice} $\caD_0$ associated with $\mu$, see Lemma \ref{lemma:DM_cubes_one_sided}, on which we construct the corona decompositions of both testing systems. Using the same lattice for $T$ and $T^*$ is crucial to make the two decompositions interact, and the same lattice is used in the good-lambda argument.
	\end{itemize}
	
	The first difficulty concerns the testing measures themselves. We write $\widehat\nu_B=\widehat b_B\,|\widehat\nu_B|$, with $|\widehat b_B|=1$, for the \textit{polar decomposition} of the testing measure. Thus, the testing measure is replaced by the bounded function $\widehat b_B$, at the expense of replacing, at the level of each ball $B$, the underlying measure $\mu$ by $|\widehat\nu_B|$. The latter need not have polynomial growth, need not be comparable with $\mu$, and need not be absolutely continuous with respect to $\mu$. Our small-set condition does not require a density in any $L^p(\mu)$ with $p>1$.
	
	There is a further difficulty in the operator testing. For each $\eta\in(0,1/10)$ and each testing ball $B$ on which the operator-testing alternative \eqref{eq:testing_opt_2} holds, there is an exceptional set $U_B=U_B(\eta)\subset B$ with $\mu(U_B)<\zeta\mu(B)$ such that $T_\sharp\widehat\nu_B\leq\mathfrak C_2(\eta)$ on $(1-\tau)B\setminus U_B$, where, for each $\eta$, the level $\mathfrak C_2(\eta)$ is uniform over the testing balls. No operator bound is assumed on $U_B$ or on $B\setminus(1-\tau)B$, and no integrability of the maximal truncations is imposed. If this alternative does not hold, the concentration alternative \eqref{eq:testing_opt_1} holds instead. Recovering the required growth, comparability, accretivity, and testing properties on a quantitatively large set is one of the main steps of our argument.
	
	In Section \ref{section:exceptional_set}, we construct, for each testing measure, an exceptional set outside of which $\widehat\nu_B$ is accretive on the cubes of a random dyadic system, $|\widehat\nu_B|$ and $\mu$ are comparable, $|\widehat\nu_B|$ has polynomial growth, and the maximal truncations of $\widehat\nu_B$ are bounded. The small-set condition \eqref{eq:measure_testing_small_sets} and the $L^0$ condition \eqref{eq:measure_testing_L0} imply that this set, together with the annulus $B\setminus(1-\tau)B$, carries a small portion of the total variation of $\widehat\nu_B$, so that its complement contains a set $F$ with $|\widehat\nu_B|(F)\gtrsim\|\widehat\nu_B\|$, see Lemma \ref{lemma:final_lemma_exceptional_set}.
%	We stress that no integrability of the maximal truncations is needed, which is also allowed by the case $u=0$ of \cite[Theorem 2.13]{HN20}, formulated through a uniform decay of the distribution tails. Our testing condition is formulated at one level, outside a set of small relative $\mu$-measure in the contracted ball, for complex testing measures satisfying the small-set condition, and it is combined with the concentration alternative.
	
	Then, as in \cite{NTrV99}, we suppress the kernel at the scale $\Xi(x)\coloneqq\dist(x,F)$ by setting
	\[
	K_\Xi(x,y)\coloneqq K(x,y)\,\frac{d(x,y)^{2n}}{d(x,y)^{2n}+\Xi(x)^n\,\Xi(y)^n}.
	\]
	This product-type suppression applies to general Calder\'on-Zygmund kernels on metric spaces, see Lemma \ref{lemma:kernel_suppressed} and Appendix \ref{appendix:suppression_metric_spaces}, and the suppressed operator satisfies the testing condition on the whole ball, see Lemma \ref{eq:lemm_bd_suppressed_2}.
	
	Next, we implement the coronization argument. The exceptional sets depend on the random parameter $w$, and so do the suppressions. By a probabilistic argument in the spirit of \cite{NTrV99}, we replace them by deterministic suppressions which still vanish on a set of measure comparable to the measure of the ball, see Lemmas \ref{lemma:lemma_probabilistic_sets} and \ref{lem:big piece of zero set}. For each testing system, we then construct a \textit{corona decomposition} on $\caD_0$:
	\begin{itemize}
		\item If the concentration alternative \eqref{eq:testing_opt_1} in the testing dichotomy holds on the testing ball of a cell $R$, we stop at a descendant of $R$ with much larger density, and the choice of $\eta=\eta_*$ in \eqref{eq:choice_eta_Stop_1} gives a packing condition for the densities of the top cells.
		\item  Otherwise, we stop at the maximal cells whose size is small compared with the deterministic suppression, which satisfy a Carleson packing condition.
		\item In Section \ref{section:mixed_corona}, we merge the two corona decompositions into a single \textit{mixed corona decomposition} by intersecting their trees, and we reduce the local estimate above to uniform $L^2$ bounds, on each tree of the mixed corona decomposition, for the operator suppressed by the maximum of the two suppressions, see Corollary \ref{corollary:49}. 
	\end{itemize}
	This is the step in which the two testing systems start to interact. In particular, each testing system is only used on its own good set, and no lower bound on the measure of the set where both systems are accretive is needed.
	The second main difficulty arises on the trees of the mixed corona decomposition. Each tree is associated with two nested David-Mattila cells $R_1$ and $R_2$, whose testing balls $B_1\coloneqq\widetilde B(R_1)$ and $B_2\coloneqq\widetilde B(R_2)$ carry a testing measure $\widehat\nu_1$ for $T$ and a testing measure $\widehat\nu_2$ for $T^*$, respectively. However, the partial order of the cells is not inherited by the balls, which need not be nested nor have comparable radii, and the natural $Tb$ argument bounds the suppressed operator from $L^2(|\widehat\nu_1|)$ to $L^2(|\widehat\nu_2|)$, with two different measures, neither of which is $\mu$.
	
	To transfer bounds between these spaces and $L^2(\mu)$, we prove in Section \ref{subsection:mixed_cotlar} a \textit{Cotlar inequality for different measures}, Lemma \ref{lemma:mixed_cotlar_transfer}, which extends \cite[Lemmas 5.25-5.27]{To14}: if a suppressed operator is bounded from $L^2(\mathfrak m_1)$ to $L^2(\mathfrak m_2)$, then its maximal truncations are bounded from $L^2(\mathfrak m_1)$ to $L^2(\mathfrak m_3)$, provided that, above the scale of the suppression, the three measures have polynomial growth and are mutually comparable on balls. By the construction of the exceptional sets, this is the case for $\mu$, $|\nu_1|$, and $|\nu_2|$, where $\nu_i\coloneqq\widehat\nu_i+\mu|_{X\setminus B_i}$. In Section \ref{section:L2bd}, we use this inequality to bound the norm on $L^2(\mu)$ of a localization of the suppressed operator to $B_2\times B_1$ by its norm from $L^2(|\widehat\nu_1|)$ to $L^2(|\widehat\nu_2|)$, see Lemma \ref{lemma:final_global_bound}.
	The latter norm is estimated in the Key Lemma \ref{lem:key_lemma_good_functions}, which is the heart of the paper. Its proof is inspired by the arguments of \cite{NTrV99}, \cite[Chapter 5]{To14}, and \cite{MV21}, with two independent random dyadic systems and \textit{martingale differences} adapted to $|\nu_i|$ and $\widehat b_i$.
	
	The possible imbalance between $B_1$ and $B_2$ creates an obstacle in the estimates of the terms involving the top cubes of the martingale decompositions, which does not appear in the one-measure case. We overcome it by a \textit{preliminary randomization} which, with large probability, allows us to work with good top cubes containing both balls, see \eqref{eq:large_probability_top_cubes} and Lemma \ref{lemma:unbalanced_lemma}.
	
	Moreover, since $K$ is not antisymmetric, the weak boundedness property is not automatic, and its proof is a delicate part of the argument. We prove it in Section \ref{section:wbp} for balls with small boundary with respect to both $|\nu_1|$ and $|\nu_2|$, by combining the small boundary property with the testing conditions on the maximal truncations, which control the operator acting on the restriction of $\widehat\nu_1$ to the complement of a ball, and we apply it on the random balls whose doubles meet both good sets. This produces an error term involving the operator norm, which we absorb at the end of the argument, after a careful choice of the (several) \textit{auxiliary parameters} that avoids circular arguments, see Table \ref{table:parameters}.
	
	Finally, in Section \ref{section:probabilisti_argument_subsection}, we combine the tree-wise estimate with a \textit{probabilistic averaging argument} over the random suppressions. The norms of the operators with random and deterministic suppressions are controlled by the norm of the averaged operator, again by means of Lemma \ref{lemma:mixed_cotlar_transfer}, see Lemma \ref{lemma:comparison_averaged_operators}, and the Key Lemma bounds the latter by a constant plus a small multiple of itself, see Lemma \ref{lemma:final_boundedness_suppressed_operators}. This gives the tree-wise bounds, uniformly over the trees and the truncations, and hence the local estimate above. The second step and the good-lambda argument of Section \ref{subsection:good_lambda} conclude the proof of Theorem \ref{theorem:main_theorem_hofmann_problem} in Section \ref{subsection:completion_main_theorem}. The preliminaries are collected in Section \ref{section:preliminaries_and_notation}, Theorem \ref{theorem:neumann_UR} is proved in Section \ref{section:neumann}, Theorem \ref{theorem:application_caloric_measure} is proved in Section \ref{section:application_caloric_measure}, and the appendices contain the constructions, in geometrically doubling metric spaces, of the David-Mattila testing balls, of random and Whitney-type cubes, and of the suppressed kernels, as well as the proof of Theorem \ref{theorem:matrix_Tb}.

	For the reader’s convenience, we collect the main parameters of the proof in Table \ref{table:parameters}.
	
	\subsection{Related results}\label{sec:related_future} While we were at an advanced stage in the preparation of this manuscript, X. Tolsa informed the second named author that M. Fernàndez-Vilaseca simultaneously and independently proved, in \cite{Fe26}, a big pieces \(Tb\) theorem with two testing functions for SIOs with non-antisymmetric kernels, assuming that the measure of the joint accretivity set is bounded below by a fixed positive proportion of the total measure.
%	
%	\vv
	We also announce an upcoming paper by John Hoffman. He proves \cite[Conjecture 1.3]{HJ25}, which is a parabolic analogue of the \textit{codimension one David-Semmes problem}, see \cite{DS91, NToV14}. More precisely, let $\mu$ be an $(n+1)$-dimensional Ahlfors regular measure on $\mathbb P^n$ with unbounded support, and assume that the parabolic Riesz transform is bounded on $L^2(\mu)$ and that $\mu$ satisfies the packing condition HSDC (\textit{higher spatial dimension condition}), see \cite[Section 3]{HJ25}. The latter means that the balls on which $\mu$ is close, in a weak sense, to a product measure $\mathcal L^n\otimes\gamma$, where $\gamma$ is a measure on the time axis, satisfy a Carleson packing condition. Then $\supp\mu$ is parabolic uniformly rectifiable, in the sense of \cite{HLN03, HLN04}, that is, the $\beta_2$-numbers of $\mu$ with respect to vertical hyperplanes satisfy a Carleson condition.
	He will also show that, if $\sigma$ is Ahlfors regular and has unbounded support, the assumptions of Theorem \ref{theorem:application_caloric_measure} imply that $\sigma$ satisfies HSDC. Hence, Hoffman's result and Theorem \ref{theorem:application_caloric_measure} will imply that, in this case, $\supp\sigma$ is parabolic uniformly rectifiable.

	\subsection*{Acknowledgments}
	The main ideas and key arguments of this paper were developed during the second-named author's visit to the Erwin Schr\"odinger International Institute for Mathematics and Physics (ESI), Vienna, for the thematic programme ``Free Boundary Problems'' in the fall of 2025.
	The proof of Theorem \ref{theorem:neumann_UR} on the free boundary Neumann problem was finalized during the second-named author's visit to the Centre de Recerca Matem\`atica (Bellaterra, Barcelona) for the program ``Analysis of Free Boundary Problems'' in the fall of 2026. He thanks both institutes for their hospitality.

	\subsection*{AI use} 
	
The proof strategy, the mathematical arguments, and an advanced draft of this paper were developed by the authors, with AI used only for editorial assistance. In the final stages of the preparation of the paper, the authors used LLM tools to check the proofs in detail, to proofread the text, to improve the presentation, and to generate Table \ref{table:parameters}. The proof checks helped the authors detect  some small gaps and inaccuracies, which were fixed by them. The authors take full responsibility for the content of the manuscript.

	\begin{table}[!ht]
		\centering
		\footnotesize
		\renewcommand{\arraystretch}{1.1}
		\setlength{\tabcolsep}{4pt}
		\begin{tabular}{@{}p{0.15\textwidth}p{0.44\textwidth}p{0.36\textwidth}@{}}
			\toprule
			Parameter &   & Reference and dependence\\
			\midrule
			\multicolumn{3}{@{}l}{\textit{Structural parameters}}\\
			$D$ & geometric doubling constant of $X$ & Theorem \ref{theorem:main_theorem_hofmann_problem}\\
			$n$, $\mathfrak C_0$ & growth exponent and growth constant of $\mu$ & \eqref{eq:def_growth}\\
			$C$, $\gamma$ & size and regularity constants of $K$ & Definition \ref{def:CZ_kernel}\\
			$\delta$, $\Lambda$, $\beta$ & doubling condition $\mu(\Lambda B)\leq\beta\mu(\delta B)$ for the balls of $\mathscr B$ & \eqref{eq:improved_doubling}; $\beta\geq\beta_0(D,n,\delta,\Lambda)$, see Remark \ref{remark:choice_of_beta}\\
			$\mathfrak C_3$ & small boundary constant of the balls of $\mathscr B$ & \eqref{eq:def_small_boundary}\\
			$K_0$, $A_0$ & parameters of the David-Mattila lattice & Lemma \ref{lemma:DM_cubes_one_sided}; $K_0$ large depending on $D$, $A_0>5000K_0$\\
			$M_*$, $\beta_*$ & testing balls: $R\subset M_*\widetilde B(R)$ and $\mu(R)\leq\beta_*\mu(\delta\widetilde B(R))$ & Lemma \ref{lemma:construction_testing_ball_DM}; depend on $D$, $\delta$, $\Lambda$, the lattice\\
			\midrule
			\multicolumn{3}{@{}l}{\textit{Testing parameters}}\\
			$\mathfrak C_1$ & bound for the total variation of the testing measures & \eqref{eq:measure_testing_size}\\
			$\varepsilon_0$ & relative scale of the small-set condition & \eqref{eq:measure_testing_small_sets}\\
			$\tau$, $\zeta$ & width of the annulus and bound for $\mu(U_B)$ in the $L^0$ alternative & \eqref{eq:testing_opt_2}; $\tau<\tau_0$, $\zeta<\zeta_0$, see \eqref{eq:measure_testing_tau_zeta}\\
			$\eta_*$ & concentration scale in \eqref{eq:testing_opt_1} & \eqref{eq:choice_eta_Stop_1}; depends on $n$, $\zeta$, the lattice, and the testing balls\\
			$\mathfrak C_2(\eta_*)$ & level of the $L^0$ condition & \eqref{eq:measure_testing_L0}; normalized to $1$ in Section \ref{subsection:bounded_testing_reduction}\\
			\midrule
			\multicolumn{3}{@{}l}{\textit{Exceptional sets and suppressions, Sections \ref{section:exceptional_set} and \ref{sec:corona}}}\\
			$\eta=(2\mathfrak C_1)^{-1}$ & accretivity threshold & Lemma \ref{lemma:estim_T_i}\\
			$m_1$, $m_2$ & thresholds for high and low density of $|\widehat\nu|$ with respect to $\mu$ & Lemma \ref{lemma:estim_nu_1}; depend on $\mathfrak C_1$, $\varepsilon_0$\\
			$\widehat{\mathbf m}$ & threshold for the dyadic density of $|\widehat\nu|$ with respect to $\mu$ & Lemma \ref{lemma:estim_S_m_hat}; depends on $\mathfrak C_1$, $\varepsilon_0$\\
			$\mathbf m$ & growth threshold for $|\widehat\nu|$ & Lemma \ref{lemma:estim_ex_H}; depends on $n$, $\varepsilon_0$, $\mathfrak C_0$, $\mathfrak C_1$\\
			$\alpha$ & threshold for the maximal truncations of $\widehat\nu$ & Lemma \ref{lemma:measures_exceptional_set_S}; $\alpha\geq2\max\{\alpha'_E,\mathfrak C_2\}$\\
			$\theta$ & width of the annulus removed from the good set & $\theta=\min\{1/4,\varepsilon_0/(2\beta\mathfrak C_3)\}$, see Lemma \ref{lemma:exceptional_set_new_contraction}\\
			$\delta_1$ & $1-\delta_1$ is a lower bound for the proportion of $|\widehat\nu|$ in the good set & $\delta_1=1-9/(20\mathfrak C_1)$, see Lemma \ref{lemma:final_lemma_exceptional_set}\\
			$\delta_{\mathrm{pr}}$ & probability threshold for the deterministic suppressions & $\delta_{\mathrm{pr}}=(1-\delta_1)/8$, see Lemma \ref{lemma:lemma_probabilistic_sets}\\
			$\mathbf c_T$ & packing constant of the top cells & Lemma \ref{lem:big piece of zero set}; depends on $\beta$, $\beta_*$, $D$, $\mathfrak C_1$, $\varepsilon_0$\\
			$c_\sharp$ & bound for the suppressed maximal truncations of $\widehat\nu$ & Lemma \ref{eq:lemm_bd_suppressed_2}\\
			\midrule
			\multicolumn{3}{@{}l}{\textit{Key Lemma and probabilistic argument, Sections \ref{section:suppressed_mixed_good} and \ref{sec:main_theorem_antisym}}}\\
			$\kappa$ & probability that the top cubes are not good & \eqref{eq:large_probability_top_cubes}; $\kappa<\delta_{\mathrm{pr}}^2/2$\\
			$m$, $\alpha=\frac{\gamma}{2(n+\gamma)}$ & goodness parameter and goodness exponent & Definition \ref{def:good_bad_cubes}\\
			$\lambda$ & width of the boundary layers of the dyadic cubes & Lemma \ref{lem:key_lemma_good_functions}\\
			$\upsilon$, $M_\Lambda$ & non-covering probability and small boundary constant of the random balls & Lemma \ref{lem:key_lemma_good_functions} and \eqref{eq:small_boundary_random_balls}; $M_\Lambda=M_\Lambda(m,\upsilon)$\\
			$a$ & enlargement of the balls in the weak boundedness property & Lemma \ref{lem:key_lemma_good_functions}\\
			\bottomrule
		\end{tabular}
		\caption{Main parameters of the proof. The parameters $m$, $\lambda$, $\upsilon$, and $a$ are fixed in this order in the proof of Lemma \ref{lemma:final_boundedness_suppressed_operators}, so that the coefficient of the operator norm in \eqref{eq:absorbing_terms} is at most $1/2$.}
		\label{table:parameters}
	\end{table}
	
	\FloatBarrier
	\section{Preliminaries and notation}\label{section:preliminaries_and_notation}

	\subsection{General notation and definitions}
	
	Given $a, b>0$, we write $a\lesssim b$ if there exists $C>0$ such that $a\leq Cb$, and $a\lesssim_t b$ to specify that the constant $C$ depends on
	the parameter $t$. We write $a\approx b$ to mean $a\lesssim b\lesssim a$, and define $a\approx_t b$ similarly.
	
	For $A\subset X$, we denote by $\chi_A$ its characteristic function.  We denote by $M(X;\mathbb C)$ the space of finite complex Radon measures on $X$ and, for $\nu \in M(X; \mathbb C)$, we denote by $\|\nu\|\coloneqq |\nu|(X)$ its total variation.
	We also write $M_+(X)$ for the family of non-negative Radon measures on $X$. All the Radon measures we consider are finite on bounded sets. For $\mathfrak m\in M_+(X)$, we denote by $L^1_{\loc}(\mathfrak m)$ the space of the functions that are $\mathfrak m$-integrable on every bounded set.
	
	We recall, see for instance \cite[Lemma 2.2]{HM12}, that, if $(X,d)$ is a geometrically doubling metric space with parameter $D$, then every ball $B(x,r)\subset X$ can contain the centers $x_i$ of at most $D\alpha^{-\log_2 D}$ disjoint balls $B(x_i,\alpha r)$, for every $\alpha\in (0,1]$. The parameter $\mathfrak d\coloneqq \log_2 D$ is often referred to as the \textit{doubling dimension} of $X$.

	A ball $B(x,r)$ is said to be $\mathfrak m$-$(a,\beta)$-doubling if
	\begin{equation}\label{eq:doubling}
		\mathfrak m(B(x,ar))\leq \beta\, \mathfrak m(B(x,r)),
	\end{equation}
	and we denote by $\mathfrak B(\mathfrak m; a, \beta)$ the family of $\mathfrak m$-$(a,\beta)$-doubling balls in $X$.

	\begin{remark}\label{remark:doubling}
		If $B\subset X$ satisfies the doubling-type assumption \eqref{eq:improved_doubling} for a measure $\mu$ as in the statement of the main theorem, then $B$ is both $\mu$-$(5, \beta)$- and $\mu$-$(100, \beta)$-doubling. Indeed, we have
		\[
		\mu(5\, B)\leq \mu(100\, B)\leq \mu (\Lambda\, B)\leq \beta\, \mu ( \delta\, B)\leq \beta\, \mu(B).
		\]
	\end{remark}

	Let $n$ be the growth exponent as in \eqref{eq:def_growth}. Given a ball $B\subset X$, we define its density
	\begin{equation}\label{eq:density}
		\Theta_{\mathfrak m}(B)\coloneqq \frac{\mathfrak m (B)}{r(B)^n}. 
	\end{equation}
	\vv
	
	{ Given $a>0$ and a non-negative Radon measure $\mathfrak m$ on $X$, we say that a ball $B\subset X$ centered at $x_B$ has \textit{$a$-small boundary} if, for every $\lambda>0$ we have
		\begin{equation}\label{eq:def_small_boundary}
			\mathfrak m\bigl(\bigl\{x\in 5B: |\dist (x,x_B)-r(B)|\leq \lambda\,  r(B)\bigr\}\bigr)\leq a \lambda\,  \mathfrak m(5B).
		\end{equation}
		We remark that, for $0<\lambda\leq 1$, both $B\setminus (1-\lambda) B$ and $(1+\lambda)B\setminus B$ are subsets of
		\[
		\bigl\{x\in 5B: |\dist (x,x_B)-r(B)|\leq \lambda\,  r(B)\bigr\}.
		\]
		
		\begin{lemma}\label{lemma:balls_small_boundary}
			Let $\mathfrak m$ be a non-negative Radon measure on $X$. For $x\in X$, $r>0$, and $0<\varepsilon\leq 1$, there exists $\rho \in [r, (1+\varepsilon)r]$ such that $B(x,\rho)$ has $\mathfrak C_\varepsilon$-small boundary with respect to $\mathfrak m$, where $\mathfrak C_\varepsilon\coloneqq 20(1+\varepsilon)/\varepsilon.$
		\end{lemma}
		\begin{proof}
			We denote $B\coloneqq B(x,r)$ and $\widetilde B\coloneqq B(x,5r)$. We also define $\phi_x\colon X\to \R$ by $\phi_x(y)\coloneqq d(x,y)$, and the push-forward measure $\sigma_x \coloneqq {\phi_x}_\sharp (\mathfrak m|_{\widetilde B})$ on $\R$.
			
			We denote by $\mathcal L^1$ the Lebesgue measure on $\mathbb R$.	
			By the $M(\R)\to L^{1,\infty}(\mathcal L^1)$ boundedness of the centered maximal operator
			\[
			M\sigma_x (s) \coloneqq \sup_{t>0}\frac{\sigma_x((s-t, s+t))}{2t}, \qquad \text{ for }s\in \R,
			\]
			we have
			\begin{equation}\label{eq:weak11_metric}
				\mathcal L^1\bigl(\bigl\{s\in \R: M\sigma_x(s)>\lambda\bigr\}\bigr)\leq \frac{5}{\lambda}\, \|\sigma_x\|.
			\end{equation}
			
			If $\mathfrak m(\widetilde B)=0$, every $\rho\in [r,(1+\varepsilon)r]$ proves the lemma. If $\mathfrak m (\widetilde B)>0$, for $\lambda = 10\mathfrak m (\widetilde B)/(\varepsilon r)$ the bound \eqref{eq:weak11_metric} yields
			\[
			\mathcal L^1\bigl(\bigl\{s\in[r,(1+\varepsilon)r] : M\sigma_x(s)>\lambda\bigr\}\bigr)\leq \frac{\varepsilon r}{2},
			\]
			so there exists $t\in [r,(1+\varepsilon)r]$ such that $M\sigma_x(t)\leq \lambda$, which concludes the proof of the lemma.
		\end{proof}
		
		We remark that the proof of Lemma \ref{lemma:balls_small_boundary} actually implies that the set of radii $\rho$ with the property in the statement has $\mathcal L^1$ measure at least $\varepsilon r/2$, so there are infinitely many possible choices. In particular, $\rho$ can be chosen in the open interval $(r,(1+\varepsilon)r)$. Similarly, given two non-negative Radon measures $\mathfrak m_1$ and $\mathfrak m_2$, running the same proof with $\lambda=20\,\mathfrak m_i(\widetilde B)/(\varepsilon r)$, $i=1,2$, shows that there exists $\rho\in(r,(1+\varepsilon)r)$ such that $B(x,\rho)$ has $2\mathfrak C_\varepsilon$-small boundary with respect to both $\mathfrak m_1$ and $\mathfrak m_2$.

	}
	
	For $G\subset X$, we denote by $\mathfrak m|_G$ the restriction of $\mathfrak m$ to $G$ and, for $p\in [1,\infty]$ we also write $L^p(G, \mathfrak m)\coloneqq L^p(\mathfrak m|_G).$
	Given two $\mathfrak m$-measurable functions $f,g\colon X\to \mathbb C$ with $fg\in L^1(\mathfrak m)$, we define the bilinear pairing
	\begin{equation}\label{eq:def_langle_rangle_m}
		\langle f,g\rangle_{\mathfrak m}\coloneqq \int f\, g\, \, d\mathfrak m
	\end{equation}
	and, if $E\subset X$ is a measurable set with $0<\mathfrak m(E)<\infty$, we denote
	\[
	\avint_E f\, d\mathfrak m\coloneqq \frac{1}{\mathfrak m(E)}\int_E f\, d\mathfrak m.
	\]
	In particular, $T^*$ is the transpose of $T$ with respect to this bilinear pairing.

	\subsection{Maximal functions}
	Let $\mathfrak m$ be a positive Radon measure on $X$ and $\mathfrak n \in M(X;\mathbb C)$.
	We define the \textit{non-centered dilated Hardy-Littlewood maximal function} as
	\begin{equation}\label{eq:def_centered_maximal}
		\widetilde{\mathcal M}_{\mathfrak m}{\mathfrak n}(x)\coloneqq \sup_{B\ni x} \frac{|\mathfrak n|(B)}{\mathfrak m(5B)}\qquad \text{ for }x\in X.
	\end{equation}
	For $f\in L^1_{\loc}(\mathfrak m)$, we denote $	\widetilde{\mathcal M}_{\mathfrak m}f\coloneqq 	\widetilde{\mathcal M}_{\mathfrak m}(f\, \mathfrak m)$, which defines a bounded sublinear operator from $M(X;\mathbb C)$ to $L^{1,\infty}(\mathfrak m)$ and from $L^p(\mathfrak m)$ to $L^p(\mathfrak m)$, for $p\in (1,\infty)$, with norms depending only on $D$ and the integrability exponent, see \cite[Proposition 3.5]{Hy10}, or  \cite[Theorem 9.32]{To14} for its Euclidean analogue. In particular, by the $5r$-covering lemma, applied to the balls $B'$ with bounded radii such that $|\mathfrak n|(B')>\lambda\,\mathfrak m(5B')$, for every $\lambda>0$ we have
	\begin{equation}\label{eq:measure_testing_maximal_constants}
		\mathfrak m\bigl(\{\widetilde{\mathcal M}_{\mathfrak m}\mathfrak n>\lambda\}\bigr)\leq\frac{\|\mathfrak n\|}{\lambda}
		\qquad\text{ and }\qquad
		\mathfrak m\bigl(\{M^r_{n}\mathfrak n>\lambda\}\bigr)\leq\frac{5^n\mathfrak C_0\,\|\mathfrak n\|}{\lambda},
	\end{equation}
	where, in the second estimate, $\mathfrak m$ has $n$-growth with constant $\mathfrak C_0$, so that $\mathfrak m(B(x,5\rho))\leq5^n\mathfrak C_0\,\rho^n$.

	\vv
	
	We write
	\[
	M^r_{n}{\mathfrak n}(x)\coloneqq \sup_{\rho>0} \frac{|\mathfrak n|(B(x,\rho))}{\rho^n}, \qquad \text{ for }x\in X,
	\]
	to denote the \textit{$n$-dimensional radial maximal function}. If $\mathfrak m\in M_+(X)$ has $n$-growth with constant $\mathfrak C_0>0$, $M^r_{n}$ defines a bounded sublinear operator from $M(X; \mathbb C)$ to $L^{1,\infty}(\mathfrak m)$ with norm depending on $\mathfrak C_0$ and $D$, see \cite[p. 53]{To14} for the proof in Euclidean spaces.
	For $f\in L^1_{\loc}(\mathfrak m)$, when $n$ is clear from the context, we also write $M^r_{\mathfrak m}f\coloneqq M^r_{n}(f\, \mathfrak m)$.

	\vv
	
	\subsection{David-Mattila-type cells}\label{sec:dyadic_grids}
	The following lemma summarizes the properties of the lattice of dyadic ``cubes'' adapted to a compactly supported Radon measure $\sigma$ on a geometrically doubling metric space, which we refer to as $\sigma$-\emph{cubes} or \emph{cells}. This dyadic system has been constructed by David and Mattila in \cite[Theorem 3.2]{DM00}, for measures on the Euclidean space, and further extended to the setting of geometrically doubling spaces by A. Volberg and P. Zorin-Kranich in \cite[Lemma 2.7]{VZ18}\footnote{In \cite{VZ18}, the authors construct the dyadic system of cells within the framework of \textit{upper doubling, geometrically doubling metric measure spaces}, assuming the existence of a suitable \textit{dominating function} $\lambda$, see also Appendix \ref{appendix:suppression_metric_spaces}.}.
	\begin{lemma}\label{lemma:DM_cubes}
		Let $(X,d)$ be a geometrically doubling metric space with geometrically doubling constant $D$, and denote $\sd\coloneqq \log_2(D)$.
		Let $\sigma$ be a Radon measure on $X$ whose support $W\coloneqq \supp(\sigma)$ is compact and has positive diameter. For $K_0 > 1$, one can choose $A_0 > 5000\,  K_0$ such that there exists a sequence of partitions $(\mathcal D_{\sigma,k})_{k\geq 0}$ of $W$ into Borel subsets $Q$, which we will refer to as $\sigma$-cubes, with the following properties:
		\begin{itemize}
			\item  For each integer $k\geq 0$, $W$ is the disjoint union of the $\sigma$-cubes $Q$, for $Q\in\mathcal D_{\sigma,k}$. If $k < \ell$,
			$Q \in \mathcal D_{\sigma,\ell}$, and $R\in\mathcal D_{\sigma,k}$, then either $Q\cap R=\varnothing$ or $Q\subset R$.
			\item 
			For each $k\geq 0$ and each $\sigma$-cube $Q \in \mathcal D_{\sigma,k}$, there is a ball $B(Q) = B(z_Q ,r(Q))$ such that
			\begin{equation*}
				\begin{split}
					z_Q\in W, \,\,\qquad  A_0^{-k}\,\diam(W) \leq r(Q)\leq K_0\,A_0^{-k} \,\diam(W),\\
					W\cap B(Q)\subset Q \subset W\cap 28 B(Q)= W\cap B(z_Q,28 r(Q)),
				\end{split}
			\end{equation*}
			and the balls $5B(Q),Q\in\mathcal D_{\sigma,k}$ are pairwise disjoint.
			\item The $\sigma$-cubes $Q\in\mathcal D_{\sigma,k}$ have small boundaries. Namely, for $Q\in\mathcal D_{\sigma,k}$ and $\ell\geq 0$, if we denote
			\begin{equation*}
				\begin{split}
					N_\ell^{\textup{int}}(Q)\coloneqq   \bigl\{x\in Q: \dist (x,W\setminus Q)<A_0^{-k-\ell}\,\diam(W)\bigr\},\\
					N_\ell^{\ext}(Q)\coloneqq   \bigl\{x\in W\setminus Q:\dist(x,Q)<A_0^{-k-\ell}\,\diam(W)\bigr\},
				\end{split}
			\end{equation*}
			and
			\(
			N_\ell(Q)\coloneqq   N_\ell^{\textup{int}}(Q)\cup N_\ell^{\ext}(Q),
			\)
			we have 
			\begin{equation*}
				\sigma(N_\ell(Q))\leq \vartheta_{DM}^{\ell}\sigma\big(90B(Q)\big), \qquad  \text{ for }\qquad \vartheta_{DM}\coloneqq  \frac{C_D K_0^{3\sd+1}}{A_0}.
			\end{equation*}
			Here $C_D$ depends only on the geometric doubling constant.
			\item Let us denote by $\mathcal{D}^{db}_{\sigma,k}$ the family of $\sigma$-cubes \( Q \in \mathcal{D}_{\sigma,k} \) for which  
			\begin{equation}\label{eq:doubling_david_mattila_balls}
				\sigma\big(100B(Q)\big) \leq K_0\,  \sigma(B(Q)).
			\end{equation}  
			
			For \( Q \in \mathcal{D}_{\sigma,k} \setminus \mathcal{D}^{db}_{\sigma,k} \), we have that \( r(Q) = A_0^{-k}\,\diam(W) \).  
			Moreover, for any \( \ell \geq 1 \) with \( 100^\ell \leq K_0 \) and \( Q \in \mathcal{D}_{\sigma,k} \setminus \mathcal{D}^{db}_{\sigma,k} \), we have  
			\begin{equation}\label{ndb}
				\sigma\big(100B(Q)\big) \leq K_0^{-\ell} \sigma\big(100^{\ell+1}B(Q)\big).
			\end{equation}
		\end{itemize}
	\end{lemma}
	
	\vvv
	Finally, we denote
	\[
	\mathcal{D}_\sigma \coloneqq \bigcup_{k\geq 0} \mathcal{D}_{\sigma, k},\qquad \qquad \mathcal{D}^{db}_\sigma \coloneqq \bigcup_{k\geq 0} \mathcal{D}^{db}_{\sigma, k},
	\]
	and choose \( A_0 \) large enough so that  
	\begin{equation}\label{assumption_a_0}
		A_0>5000 K_0\qquad \text{ and }\qquad \vartheta_{DM} <1.
	\end{equation}
	
	\begin{remark} No small-boundary property is guaranteed for balls $B(Q)$ associated with doubling David-Mattila cells $Q$. In particular, these balls need not belong to the family $\mathscr B$ of Theorem \ref{theorem:main_theorem_hofmann_problem}, so no testing measures are associated with them.
		This difficulty is overcome in Appendix \ref{appendix:DM_testing_balls} by introducing some balls $\widetilde B(Q)$, centered at points of $Q$, that are admissible for the testing hypotheses.
	\end{remark}
	\vv

	For measures with possibly unbounded support, we need a construction of the David-Mattila lattice with radii comparable to $sA_0^{-k}$, $k\geq0$, where $s>0$.
	
	\begin{lemma}[Auxiliary DM-lattice $\caD^{(s)}_\sigma$]\label{lemma:DM_cubes_one_sided}
		Let $(X,d)$ be a geometrically doubling metric space with geometrically doubling constant $D$, and denote $\sd\coloneqq \log_2(D)$. Let $s>0$.
		Let $\sigma$ be a Radon measure on $X$ and denote $W\coloneqq \supp(\sigma)$. For $K_0 > 1$, one can choose $A_0 > 5000\,  K_0$ such that there exists a sequence of partitions $(\mathcal D^{(s)}_{\sigma,k})_{k\geq 0}$ of $W$ into Borel subsets $Q$ with the following properties:
		\begin{itemize}
			\item  For each integer $k\geq 0$, $W$ is the disjoint union of $Q$, for $Q\in\mathcal D^{(s)}_{\sigma,k}$. If $k < \ell$,
			$Q \in \mathcal D^{(s)}_{\sigma,\ell}$, and $R\in\mathcal D^{(s)}_{\sigma,k}$, then either $Q\cap R=\varnothing$ or $Q\subset R$.
			\item 
			For each $k\geq 0$ and each $Q \in \mathcal D^{(s)}_{\sigma,k}$, there is a ball $B^{(s)}(Q) = B(z_Q ,r^{(s)}(Q))$ such that
			\begin{equation}\label{eq:new_radii_modified_DM}
				\begin{split}
					z_Q\in W, \,\,\qquad  s\, A_0^{-k} \leq r^{(s)}(Q)\leq s\, K_0\,A_0^{-k},\\
					W\cap B^{(s)}(Q)\subset Q \subset W\cap 28 B^{(s)}(Q)= W\cap B(z_Q,28 r^{(s)}(Q)),
				\end{split}
			\end{equation}
			and the balls $5B^{(s)}(Q),Q\in\mathcal D^{(s)}_{\sigma,k}$ are pairwise disjoint.
			\item For $Q\in\mathcal D^{(s)}_{\sigma,k}$ and $\ell\geq 0$, if we denote
			\begin{equation*}
				\begin{split}
					N_\ell^{\textup{int},(s)}(Q)\coloneqq   \bigl\{x\in Q: \dist (x,W\setminus Q)<s\, A_0^{-k-\ell}\bigr\},\\
					N_\ell^{\ext,(s)}(Q)\coloneqq   \bigl\{x\in W\setminus Q:\dist(x,Q)<s\, A_0^{-k-\ell}\bigr\},
				\end{split}
			\end{equation*}
			and
			\(
			N^{(s)}_\ell(Q)\coloneqq   N_\ell^{\textup{int},(s)}(Q)\cup N_\ell^{\ext,(s)}(Q),
			\)
			we have 
			\begin{equation*}
				\sigma(N^{(s)}_\ell(Q))\leq \vartheta_{DM}^{\ell}\sigma\big(90B^{(s)}(Q)\big), \qquad { \text{ for }\qquad \vartheta_{DM}\coloneqq  \frac{C_D K_0^{3\sd+1}}{A_0}.}
			\end{equation*}
			\item Let us denote by $\mathcal{D}^{db,(s)}_{\sigma,k}$ the family of the sets \( Q \in \mathcal{D}^{(s)}_{\sigma,k} \) for which  
			\begin{equation}\label{eq:doubling_david_mattila_balls_sec8}
				\sigma\big(100B^{(s)}(Q)\big) \leq K_0\,  \sigma(B^{(s)}(Q)).
			\end{equation}  
			
			For \( Q \in \mathcal{D}^{(s)}_{\sigma,k} \setminus \mathcal{D}^{db,(s)}_{\sigma,k} \), we have that \( r^{(s)}(Q) = s\, A_0^{-k}\).  
			Moreover, for any \( \ell \geq 1 \) with \( 100^\ell \leq K_0 \) and \( Q \in \mathcal{D}^{(s)}_{\sigma,k} \setminus \mathcal{D}^{db,(s)}_{\sigma,k} \), we have  
			\begin{equation}\label{ndb_sec8}
				\sigma\big(100B^{(s)}(Q)\big) \leq K_0^{-\ell} \sigma\big(100^{\ell+1}B^{(s)}(Q)\big).
			\end{equation}
		\end{itemize}
		Here $C_D$ depends only on $D$, and $A_0$ is chosen large enough that $\vartheta_{DM}<1$. Finally, we denote
		\[
		\mathcal{D}^{(s)}_\sigma \coloneqq \bigcup_{k\geq 0} \mathcal{D}^{(s)}_{\sigma, k},\qquad \qquad \mathcal{D}^{db,(s)}_\sigma \coloneqq \bigcup_{k\geq 0} \mathcal{D}^{db,(s)}_{\sigma, k}.
		\]
	\end{lemma}
	\begin{proof}
		The existence of a lattice with the properties listed above follows from \cite[Lemma 2.7]{VZ18}, applied to the metric space $(X,s^{-1}d)$, which is geometrically doubling and has geometrically doubling constant $D$, as well. In particular, the non-doubling estimate there reads $\sigma(cB^{(s)}(Q))\leq K_0^{-1}\sigma(100cB^{(s)}(Q))$ for $1\leq c\leq K_0$, which implies \eqref{ndb_sec8}.
	\end{proof}
	
	Cells are understood together with their generation indices. Distinct generations may give the same underlying subset, and ancestry and strict descent refer to these indexed cells.
	
	\vv
	
	Applying the construction in Appendix \ref{appendix:DM_testing_balls} to $(X,s^{-1}d)$ gives testing balls $\widetilde B^{(s)}(Q)$ for $Q\in\caD^{db,(s)}_\sigma$, with the same inclusions and measure estimates. If $Q\in \caD^{(s)}_{\sigma,k}$, we define the adapted side-length as
	\[
	\ell^{(s)}(Q)\coloneqq s\, 56\, K_0\, A^{-k}_0.
	\]

	\subsubsection{Auxiliary constants and remarks}\label{section:auxiliary_DM_constants}
	We introduce the technical constants used to construct the testing balls in Appendix \ref{appendix:DM_testing_balls}. They depend on $\delta$, $\Lambda$, and the parameters of the David-Mattila lattice. A short description of the role of the constants and their first appearance is added to guide the reader.
	
	We define $j_*$ as the minimum integer $j\geq 0$ such that $K_0 \vartheta^j_{DM}\leq 1/2$, and we set $c_*\coloneqq A_0^{-j_*}/K_0$. We also write $r_Q\coloneqq r(Q)$. 
	With this choice, we have the following estimates.
	\begin{lemma}\label{lemma:definition_F_Q_separation}
		For every $Q\in \caD^{db}_{\sigma, k}$ we have that 
		\[
		\sigma\bigl(Q\setminus N_{j_*}^{\textup{int}}(Q)\bigr)\geq \frac{1}{2}\, \sigma(Q)
		\]
		and
		\[
		\dist\bigl(x, \supp(\sigma)\setminus Q\bigr)\geq c_*\, r_Q, \qquad \text{ for every }\qquad x\in Q\setminus N_{j_*}^{\textup{int}}(Q).
		\]
	\end{lemma}
	\begin{proof}
		By Lemma \ref{lemma:DM_cubes} and the doubling of $Q$, we have
		\[
		\sigma\bigl(N_{j_*}^{\textup{int}}(Q)\bigr)\leq\vartheta_{DM}^{j_*}\, \sigma(90B(Q))\leq K_0\vartheta_{DM}^{j_*}\, \sigma(B(Q))\leq\frac12\,\sigma(Q).
		\]
		For $x\in Q\setminus N_{j_*}^{\textup{int}}(Q)$, the definition of $N_{j_*}^{\textup{int}}(Q)$ and $r_Q\leq K_0A_0^{-k}\diam(W)$ give
		\[
		\dist(x,\supp(\sigma)\setminus Q)\geq A_0^{-k-j_*}\diam(W)\geq c_*r_Q.\qedhere 
		\]
	\end{proof}
	
	For later use, we also set:
	\begin{itemize}
		\item $\displaystyle M_*\coloneqq 224\Lambda/c_*.$ It is the dilation for which the testing ball contains its cell, namely $Q\subset M_* \widetilde B(Q)$, see \eqref{eq:structure_testing_balls}. 
		\item $\displaystyle A_{\Lambda, *}\coloneqq 512 \Lambda/c_*.$ This constant appears in \eqref{eq:inclusions_84_100_B_B_tilde}, and is associated with enlarged testing balls.
		\item $\displaystyle m_*\coloneqq 4+\ceil{\log_2(M_*/\delta)}=4+\ceil{\log_2(224\Lambda/(c_*\delta))}$. It is the number of dyadic scales between the radius $28 r_Q$ of a ball containing $Q$ and the radius $\delta R_Q/2$, where $R_Q\coloneqq c_* r_Q/(4\Lambda)$ is the smallest admissible radius of $\widetilde B(Q)$, see \eqref{eq:property_def_m_*}.
		\item $\displaystyle \beta_*\coloneqq 2D^{m_*}.$ It is the doubling constant of testing balls, see Lemma \ref{lemma:construction_testing_ball_DM}, and it is defined by a standard covering argument.
	\end{itemize}
	\begin{remark}[Choice of $\beta$]\label{remark:choice_of_beta}
		In Theorem \ref{theorem:main_theorem_hofmann_problem}, it suffices to assume
		\[
		\beta>\max\Bigl\{\beta_*, K_0, \Bigl(\frac{1.1 \Lambda}{\delta}\Bigr)^n, D^{\ceil{\log_2(1.1\Lambda/\delta)}} \Bigr\}.
		\]
		
		The term $\beta_*$ implies that the testing balls of doubling David-Mattila cells, which satisfy \eqref{eq:thin_boundary_testing_ball_DM}, belong to the family of admissible balls $\mathscr B$ of Theorem \ref{theorem:main_theorem_hofmann_problem}. The last two terms on the right-hand side are used to guarantee the existence of balls $B$ centered at $x$ with $\mu(1.1\Lambda\, B)\leq \beta\, \mu(\delta\, B)$, see \cite[Lemma 3.3]{Hy10}, where the factor $1.1$ allows us to enlarge the radius of such a ball by a certain factor so that, by Lemma \ref{lemma:balls_small_boundary}, we can reduce to a ball with small boundary.
	\end{remark}
	\vv
	
	Here, we list some quantities and properties associated with $Q\in\mathcal D_\sigma$:  
	\begin{itemize}  
		\item \( J(Q) \coloneqq k \), which we refer to as the \textit{generation} of \( Q \).  
		\item \( \ell(Q) \coloneqq 56 K_0 A_0^{-k}\,\diam(W) \), that we call \textit{side length}. We note that  
		\begin{equation}\label{eq:comparab_r_l_0}
			\diam(Q)\leq\diam(28B(Q))\leq56r(Q)\leq\ell(Q),
		\end{equation}  
		and that $r(Q)\approx\ell(Q)$, with constants depending on $K_0$.   Moreover,
		\begin{equation}\label{eq:comparab_r_l}
			\frac{1}{1000}\, K_0^{-1}\ell(Q)\leq r(Q)\leq \ell(Q).
		\end{equation}
		\item Denoting by \( z_Q \) the \textit{center} of \( Q \), we define \( B_Q \coloneqq 28 B(Q) = B(z_Q, 28r(Q)) \), so that  
		\begin{equation*}
			W \cap \frac{1}{28} B_Q \subset Q \subset B_Q.
		\end{equation*}
		\item For $Q\in\mathcal D_\sigma,$ we define $\mathcal D_\sigma(Q)\coloneqq \{P\in\mathcal D_\sigma: P\subset Q\}$.
		\item If $Q\subset Q'$ are cells of $\caD_\sigma$ and $A_0$ is big enough, then $30 B(Q)\subset 30 B(Q')$, see also \cite[Section 3]{VZ18}.
		\item Every $Q\in \caD_{\sigma}$ can be covered up to a $\sigma$-null set by a union of cells in $\caD^{db}_\sigma$, see \cite[Lemma 5.28]{DM00}.
	\end{itemize}
	
	For a measure with $n$-growth we write
	\[
	\Theta_\mu(100B(Q))\coloneqq\frac{\mu(100B(Q))}{(100r(Q))^n}.
	\]
	More generally, in the upper-doubling setting of Remark \ref{remark:upper_doubling}, the denominator may be replaced by $\lambda(z_Q,100r(Q))$. We assume in this case that the dominating function is nondecreasing in its second variable and satisfies
	\[
	\lambda(x,2r)\leq\mathbf C_\lambda\, \lambda(x,r),\qquad
	\lambda(x,r)\leq\mathbf C_\lambda\, \lambda(y,r)\qquad\text{if }d(x,y)\leq r.
	\]
	Such a dominating function is supplied by Lemma \ref{lemma:aux_dominating_function_3}. For the $n$-dimensional density above we take $\lambda(x,r)=r^n$ and $\mathbf C_\lambda=2^n$.
	
	\begin{lemma}\label{lemma:lemma_density_nonboubling_chain}
		Let $K_0\geq10^4$, let $k_0$ be the largest nonnegative integer such that $100^{k_0+1}\leq K_0$, and assume that $K_0$ and $A_0$ satisfy
		\[
		K_0^{k_0/2}>\mathbf C_\lambda^{\lceil\log_2 A_0\rceil+1}.
		\]
		Let $Q\subset R$ be cells of the same David-Mattila lattice associated with $\mu$, and suppose that every cell $S$ with $Q\subsetneq S\subsetneq R$ is non-doubling. No doubling assumption is imposed on $Q$ or $R$. Then
		\[
		\Theta_\mu(100B(Q))\lesssim K_0^{-(J(Q)-J(R)-1)_+k_0/2}\, \Theta_\mu(100B(R))
		\]
		and
		\[
		\sum_{Q\subset S\subset R}\Theta_\mu(100B(S))\lesssim\Theta_\mu(100B(R)).
		\]
		The implicit constants depend only on the parameters of the lattice and of the dominating function.
	\end{lemma}
	\begin{proof}
		For the proof we refer, for instance, to \cite[Lemma 3.6]{VZ18} and \cite[Lemma 4.4]{To16}.
	\end{proof}
	These conditions hold if $A_0=C_D'K_0^{3\sd+2}$, with $C_D'$ sufficiently large for \eqref{assumption_a_0}, and $K_0$ sufficiently large depending on $D$ and $\mathbf C_\lambda$. Indeed,
	\[
	k_0\geq\frac{\log K_0}{\log100}-2,\qquad
	\lceil\log_2 A_0\rceil+1\leq\log_2 C_D'+(3\sd+2)\log_2K_0+2,
	\]
	which imply $K_0^{k_0/2}>\mathbf C_\lambda^{\lceil\log_2 A_0\rceil+1}$ for $K_0$ sufficiently large.
	\vv

	Finally, given two non-negative Radon measures $\mathfrak m, \mathfrak n$ on $X$ we define the \textit{dyadic maximal Hardy-Littlewood operator} with respect to $\caD_\sigma$ as
	\begin{equation}\label{eq:dyadic_maximal_function}
		M_{\mathfrak m, \caD}\mathfrak n(x)\coloneqq \sup_{\substack{Q\in \mathcal{D}_\sigma,\, Q\ni x \\ \mathfrak m(Q)>0}} \frac{\mathfrak n(Q)}{\mathfrak m(Q)}
	\end{equation}	
	which defines a bounded operator from the finite measures in $M_+(X)$ to $L^{1,\infty}(\mathfrak m)$. Moreover, for $1<p\leq \infty$, it defines a bounded operator from $L^p(\mathfrak m)$ to $L^p(\mathfrak m)$ whose norm depends only on $p$. See, for instance, \cite[Theorem 5.7]{To14} for a proof using a dyadic grid in $X$, which works unchanged for \eqref{eq:dyadic_maximal_function}.

	\subsection{Dyadic cubes in geometrically doubling metric spaces}\label{section:dyadic_cubes_gdms}

	Let $\el\in(0,1)$ be a fixed parameter such that $144\el\leq 1$. We recall the dyadic system of so-called cubes in the sense of Hytönen-Kairema \cite{HK12}, which we use as the suitable analogue of standard dyadic cubes on the whole space, and that in the original work was formulated for geometrically doubling \textit{quasi-metric} spaces.
	
	\begin{theorem}[\cite{HK12}, Theorem 2.2]\label{thm:HK22}
		Let $(X,d)$ be a geometrically doubling metric space with parameter $D$. Let $0<c_0\leq C_0<\infty$ and suppose that
		\[
		12C_0\el\leq c_0.
		\]
		Assume that, for every $k\in\mathbb Z$, we are given points $\{z^k_\alpha\}_{\alpha\in A_k}$ such that
		\[
		d(z^k_\alpha,z^k_\beta)\geq c_0\, \el^k\quad(\alpha\neq\beta),	\qquad \min_{\alpha\in A_k}d(x,z^k_\alpha)< C_0\, \el^k,\qquad\text{ for every }x\in X.
		\]
		Then, there exist families of sets
		\[
		\widetilde Q^k_\alpha\subset Q^k_\alpha\subset \overline Q^k_\alpha,	\qquad k\in\mathbb Z,
		\]
		called open, half-open and closed dyadic cubes, such that $\widetilde Q^k_\alpha$ and $\overline Q^k_\alpha$ are the interior and closure of $Q^k_\alpha$ respectively,
		\[
		\text{if }j\geq k,\text{ then either }Q^j_\beta\subset Q^k_\alpha\text{ or }Q^k_\alpha\cap Q^j_\beta=\varnothing,
		\]
		\[
		X=\bigcup_{\alpha\in A_k}Q^k_\alpha\qquad\text{disjointly, for every }k\in\mathbb Z,
		\]
		and
		\[
		B(z^k_\alpha,c_1\el^k)\subset Q^k_\alpha\subset B(z^k_\alpha,C_1\el^k)\eqqcolon B(Q^k_\alpha),
		\qquad c_1\coloneqq\frac{c_0}{3},\quad C_1\coloneqq2C_0.
		\]
		Moreover, if $k\leq j$ and $Q^j_\beta\subset Q^k_\alpha$, then $B(Q^j_\beta)\subset B(Q^k_\alpha).$
		In particular,
		\begin{equation}\label{eq:diam_vs_length_pre}
			\mathrm{diam} Q^k_\alpha\leq 2C_1\el^k.
		\end{equation}
	\end{theorem}
	\vv
	
	For $Q= Q^k_\alpha$ as in Theorem \ref{thm:HK22}, we also denote $\ell(Q)\coloneqq \el^k,$
	so that \eqref{eq:diam_vs_length_pre} reads
	\begin{equation}\label{eq:diam_vs_length}
		\diam(Q)\leq 2C_1\, \ell(Q).
	\end{equation}

	\begin{remark}
		For a geometrically doubling metric space $(X,d)$ with parameter $D$, we have that any ball $B(x,r)$ can be covered by at most $D\,t^{-\log_2D}$ balls of radii $t r$, for every $t \in (0,1]$. Moreover, any ball $B(x,r)$ contains at most $D\, t^{-\log_2D}$ centers $x_i$ of pairwise disjoint balls $B(x_i, t r)$. For a proof, we refer to Lemma 2.3 in \cite{Hy10}.
	\end{remark}
	
	We now define a randomization of the dyadic system of cubes.
	We say that a collection $\{x^k_\alpha\}_{k,\alpha}\subset X$ is a \textit{family of reference points} if for every $k\in\mathbb Z$,
	\[
	d(x^k_\alpha,x^k_\beta)\geq \el^k\quad(\alpha\neq\beta), \qquad \min_{\alpha\in A_k}d(x,x^k_\alpha)<\el^k\quad\text{for every }x\in X.
	\]
	
	\begin{theorem}[\cite{HK12}, Theorem 5.1]\label{thm:HK51}
		Given reference points $\{x^k_\alpha\}_{k,\alpha}\subset X$, suppose that $96\el\leq 1$. Then there exists a probability space $(\Omega,\mathbb P)$ such that every $w\in\Omega$ defines a dyadic system
		\[
		\bbD(w)=\{Q^k_\alpha(w)\}_{k,\alpha},
		\]
		related to new dyadic points $\{z^k_\alpha(w)\}_{k,\alpha}$, with the properties in Theorem \ref{thm:HK22}. Moreover,
		\[
		\Omega=\prod_{k\in\mathbb Z}\Omega_k  	\qquad w=(w_k)_{k\in\mathbb Z},
		\]
		where the coordinates $w_k$ are independent,
		\[
		z^k_\alpha(w)=z^k_\alpha(w_k),
		\]
		and there exists $\tau_0>0$ such that, whenever $(k+1,\beta)$ is a child of $(k,\alpha)$ in the construction,
		\[
		\mathbb P\bigl(\{w\in\Omega:z^k_\alpha(w)=x^{k+1}_\beta\}\bigr)\geq \tau_0.
		\]
	\end{theorem}
	
	\vv		
	For further properties of dyadic ``random'' cubes used in our manuscript, we refer to Appendix \ref{appendix:random_cubes}. We also remark that, in the Euclidean setting, random dyadic grids are obtained by means of random translations of the standard dyadic grid, see for instance \cite[Section 5]{To14}.
	For $w\in \Omega$, we set 
	\[
	\bbD_{k}(w)\coloneqq \{Q\in\bbD(w):\ell(Q)=\el^k\},\qquad k\in\mathbb Z,
	\]
	so that $\bbD(w)=\bigcup_{k\in \mathbb Z} \bbD_k(w).$

	\vv
	
	\section{The exceptional sets and the induced suppressed operator}\label{section:exceptional_set}
	
	{
		\subsection{Normalization of measure testing}\label{subsection:bounded_testing_reduction}
		The proof of Theorem \ref{theorem:main_theorem_hofmann_problem} can be reduced to the following normalized statement, which is proved in the remainder of the paper.
		
		\begin{proposition}\label{proposition:bounded_testing_normalized}
			Let $(X,d)$, $\mu$, $T$, and $\mathscr B$ be as in Theorem \ref{theorem:main_theorem_hofmann_problem}, let $\mathfrak C_1\geq1$ and $\varepsilon_0\in(0,1)$, and let $0<\tau<\tau_0$ and $0<\zeta<\zeta_0$, with $\tau_0$ and $\zeta_0$ as in \eqref{eq:measure_testing_tau_zeta}. Suppose that the measures $\{\widehat\nu_B\}_{B\in\mathscr B}$ and $\{\widehat\nu_B^*\}_{B\in\mathscr B}$ satisfy \eqref{eq:measure_testing_size} and \eqref{eq:measure_testing_small_sets}, and that, with $\eta_*$ as in \eqref{eq:choice_eta_Stop_1}, for every $B\in\mathscr B$ there exists a Borel set $U_B\subset B$ such that either \eqref{eq:testing_opt_1} holds with $\eta= \eta_*$ for some Borel set $V_B\subset U_B$ or
			\begin{equation}\label{eq:bounded_normalized_testing}
				\mu(U_B)<\zeta\mu(B)\qquad \text{ and }\qquad
				T_\sharp\widehat\nu_B\leq1\quad \text{ on }(1-\tau)B\setminus U_B,
			\end{equation}
			and that the analogous condition holds for $T^*$ and $\widehat\nu_B^*$. Then, we have
			\begin{equation}\label{eq:bounded_normalized_conclusion}
				\sup_{\varepsilon>0}
				\|T_{\mu,\varepsilon}\|_{L^2(\mu)\to L^2(\mu)}
				\leq C_{\rm norm},
			\end{equation}
			where $C_{\rm norm}$ depends only on $D$, $n$, $\mathfrak C_0$, $C$, $\gamma$, $\delta$, $\Lambda$, $\beta$, $\mathfrak C_3$, $\mathfrak C_1$, $\varepsilon_0$, and $\zeta$.
		\end{proposition}
		
		We first show that it suffices to prove Proposition \ref{proposition:bounded_testing_normalized} for auxiliary kernels with cut-offs. The associated operators are a priori bounded on $L^2(\mu)$, with bounds depending on the cut-off parameters, whereas their Calder\'on-Zygmund constants do not depend on them, so the main task is to bound their norms independently of the cut-offs. The a priori finiteness of the norms is crucial to perform absorption arguments, see \eqref{eq:absorbing_terms}.

		\vv
		
		To this end, we fix smooth functions $\varphi,\psi\colon[0,\infty)\to[0,1]$, with $\varphi$ non-decreasing and $\psi$ non-increasing, such that $\varphi=0$ on $[0,1]$, $\varphi=1$ on $[2,\infty)$, $\psi=1$ on $[0,1]$, and $\psi=0$ on $[2,\infty)$. For $0<h<L$, we denote by $T^{h,L}$ the singular integral operator associated with the kernel
		\begin{equation}\label{eq:bounded_reduction_regularization}
			K^{h,L}(x,y)\coloneqq K(x,y)\,\varphi\Bigl(\frac{d(x,y)}{h}\Bigr)\,\psi\Bigl(\frac{d(x,y)}{L}\Bigr).
		\end{equation}
		
		With $\eta_*$ as in the statement of Theorem \ref{theorem:main_theorem_hofmann_problem} (so, as in \eqref{eq:choice_eta_Stop_1}), we also set
		\begin{equation}\label{eq:bounded_testing_normalization}
			A\coloneqq1+3\,\mathfrak C_2(\eta_*),\qquad \qquad \mathcal T\coloneqq A^{-1}\,T^{h,L},
		\end{equation}
		and we prove the lemma that reduces the main theorem to the proposition above.
		\begin{lemma}\label{lemma:reduction_cutoffs}
			For every $0<h<L$, the operator $\mathcal T$ and the testing measures of Theorem \ref{theorem:main_theorem_hofmann_problem} satisfy the hypotheses of Proposition \ref{proposition:bounded_testing_normalized}, with Calder\'on-Zygmund constants independent of $h$ and $L$. Consequently, Proposition \ref{proposition:bounded_testing_normalized} gives
			\begin{equation}\label{eq:bounded_testing_uniform_regularization}
				\|T^{h,L}_\mu\|_{L^2(\mu)\to L^2(\mu)}\leq C_{\rm norm}\,A,
				\qquad \text{ for }0<h<L,
			\end{equation}
			and Theorem \ref{theorem:main_theorem_hofmann_problem} holds with 
			\begin{equation}\label{eq:bounded_testing_uniform_regularization2}
				\sup_{\varepsilon>0}\|T_{\mu,\varepsilon}\|_{L^2(\mu)\to L^2(\mu)}\leq C_{\rm norm}\,A+C\,2^n\mathfrak C_0.
			\end{equation}
		\end{lemma}
		
		\begin{proof}
			For $t\geq 0$, we set $a(t)\coloneqq\varphi(t/h)\,\psi(t/L)$, and note that $|a'(t)|\lesssim1/t$.
			Moreover, since the distance is $1$-Lipschitz in each variable, $K^{h,L}$ is a Calder\'on-Zygmund kernel with constants independent of $h$ and $L$. For every $\sigma\in M(X;\mathbb C)$, $x\in X$, and $\rho>0$, Fubini's theorem gives \[T^{h,L}_\rho\sigma(x)=a(\rho)\,T_\rho\sigma(x)+\int_\rho^\infty a'(t)\,T_t\sigma(x)\,dt,\] which, since $0\leq a\leq1$ and $\int_0^\infty|a'|\leq2$, yields $	(T^{h,L})_\sharp\sigma\leq3\,T_\sharp\sigma$ and the analogous estimate for the adjoint.
			
			Hence, if \eqref{eq:testing_opt_2} holds with $\eta=\eta_*$, then $\mathcal T_\sharp\widehat\nu_B\leq3\,\mathfrak C_2(\eta_*)/A<1$ on $(1-\tau)B\setminus U_B$, that is,
			\eqref{eq:bounded_normalized_testing} holds for $\mathcal T$, and the same holds for $\mathcal T^*$ and $\widehat\nu_B^*$. The remaining hypotheses of Proposition \ref{proposition:bounded_testing_normalized} do not involve the operator. Since $K^{h,L}(x,y)=0$ whenever $d(x,y)\leq h$, the truncations of $\mathcal T_\mu$ at radii smaller than $h$ coincide with $\mathcal T_\mu$, so \eqref{eq:bounded_normalized_conclusion} for $\mathcal T$ gives \eqref{eq:bounded_testing_uniform_regularization}.
			
			Let $f$ be a bounded function with bounded support, and let $\varepsilon>0$. Letting $L\to\infty$ in \eqref{eq:bounded_testing_uniform_regularization} with $h=\varepsilon$, Fatou's lemma gives $\|T^{\varepsilon,\infty}_\mu f\|_{L^2(\mu)}\leq C_{\rm norm}\,A\,\|f\|_{L^2(\mu)}$, where $T^{\varepsilon,\infty}$ is associated with the kernel $K(x,y)\,\varphi(d(x,y)/\varepsilon)$. The difference $T_{\mu,\varepsilon}-T^{\varepsilon,\infty}_\mu$ is an integral operator whose associated kernel is supported in $\{\varepsilon<d(x,y)<2\varepsilon\}$ and bounded by $C\varepsilon^{-n}$, so Schur's test and the $n$-growth of $\mu$ bound its norm by $C\,2^n\mathfrak C_0$. Finally, \eqref{eq:bounded_testing_uniform_regularization2} follows by density.
		\end{proof}
		
		By Lemma \ref{lemma:reduction_cutoffs}, it suffices to prove Proposition \ref{proposition:bounded_testing_normalized} for $T=\mathcal T$, with constants independent of $h$ and $L$, and from now on we write $T$ for $\mathcal T$. In particular, since the auxiliary kernel $K^{h,L}$ vanishes outside $\{h<d(x,y)<2L\}$, Schur's test and the $n$-growth of $\mu$ show that $T_\mu$ is a priori bounded on $L^2(\mu)$, as desired. Moreover, the testing measures keep their original normalization, total-variation bound, and small-set threshold.
		
	}
	
	\vv
	Throughout this section, we consider a geometrically doubling metric space $(X,d)$ with parameter $D$.

	{
		\subsection{The construction of the exceptional sets}
		Let $B\subset X$ be a ball. Suppose that $\mu$ is a non-negative Radon
		measure and $\widehat\nu$ is a complex Radon measure on $X$. Fix
		$\mathfrak C_0>0$, $\mathfrak C_1\geq1$, $\mathfrak C_2>0$,
		$\mathfrak C_3>1$, $\beta>5^{\log_2D}$, and
		$\varepsilon_0\in(0,1)$. We assume the following properties.
		\begin{enumerate}[label=(A\arabic*),ref=(A\arabic*)]
			\item[(A0)] The ball $B$ has $\mathfrak C_3$-small boundary with
			respect to $\mu$, and
			\begin{equation}\label{eq:improved_doubling_ex_set}
				\mu(5B)\leq\beta\mu(B).
			\end{equation}
			\item\label{ex_item_1} The measure $\mu$ has $n$-growth with
			constant $\mathfrak C_0$.
			\item\label{ex_item_2} We have $|\widehat\nu|(X\setminus B)=0$.
			\item\label{ex_item_3} We have $|\widehat\nu(B)|=\mu(B)>0$.
			\item\label{ex_item_4} We have
			\begin{equation}\label{eq:item_4_ref}
				\|\widehat\nu\|\leq\mathfrak C_1\mu(B).
			\end{equation}
			\item\label{ex_item_5} For every Borel set $A\subset B$ with
			$\mu(A)\leq\varepsilon_0\mu(B)$, we have
			\begin{equation}\label{eq:abs_cont}
				|\widehat\nu|(A)\leq
				\frac{\|\widehat\nu\|}{200\mathfrak C_1}.
			\end{equation}
			\item\label{ex_item_6} There is a Borel set $U_B\subset B$ such that
			\begin{equation}\label{eq:hp_6}
				\mu(U_B)\leq\varepsilon_0\,\mu(B)\qquad\text{ and }\qquad T_\sharp\widehat\nu\leq\mathfrak C_2\quad\text{ on }B\setminus U_B.
			\end{equation}
		\end{enumerate}
		Possibly by normalizing the measure $\widehat\nu$ by a constant of modulus one, we may assume that $\widehat\nu(B)=\mu(B)$. Set
		\[
		\nu\coloneqq\widehat\nu+\mu|_{X\setminus B},
		\qquad\text{ so that }\qquad |\nu|=|\widehat\nu|+\mu|_{X\setminus B}.
		\]
		Write $\widehat\nu=\widehat b\,|\widehat\nu|$ on $B$, and extend $\widehat b$ by $1$ on $X\setminus B$. Then
		\[
		\nu=\widehat b\,|\nu|,
		\qquad\text{ where }\qquad |\widehat b|=1\quad |\nu|\text{-a.e.}
		\]
		The assumptions concern $\widehat\nu$, whereas the stopping conditions for random cubes use $\nu$ and $|\nu|$.

		Assumption \ref{ex_item_6} has the same form as the normalized testing condition \eqref{eq:bounded_normalized_testing}. In the proof of Proposition \ref{proposition:bounded_testing_normalized}, it holds with $\mathfrak C_2=1$ once the annulus $B\setminus(1-\tau)B$ is added to the set in \eqref{eq:bounded_normalized_testing}, see Section \ref{sec:corona}.
		
	}	
	
	\subsubsection*{\underline{The set of non-accretive cubes}} 
	Next, fix $\eta \in (0,1)$ to be chosen later. 
	Let $(\Omega, \mathbb P)$ be the probability space as in Theorem \ref{thm:HK51} and, for $w\in \Omega$, let $\bbD(w)$ denote its associated grids of random cubes as defined in the previous section. { We denote by $\mathcal A(w)$ the collection of all dyadic cubes $S\in\bbD(w)$ such that $S\cap B\neq\varnothing$ and
		\[
		|\nu(S)|<\eta\,|\nu|(S).
		\]
	}
	We then set
	\[
	\mathscr T(w) \coloneqq \bigcup_{S \in \mathcal A(w)} S.
	\]

	\begin{lemma}\label{lemma:estim_T_i}
		For $\eta = (2\mathfrak C_1)^{-1}$, we have
		\[
		|{\widehat\nu}|(\mathscr T(w)) \leq (1-\eta)\,|{\widehat\nu}|(B), \qquad \text{ for all }w\in \Omega.
		\]
	\end{lemma}
	
	\begin{proof}
		{
			Let $w\in\Omega$ and $N\in\mathbb N$. Denote by $\mathcal A_N(w)$ the maximal cubes in $\mathcal A(w)$ whose generation is at least $-N$, and set $\mathscr T_N(w)\coloneqq\bigcup_{S\in\mathcal A_N(w)}S$. These cubes are pairwise disjoint. For every $S\in\mathcal A_N(w)$, the stopping condition gives
			\[
			\operatorname{Re}\widehat\nu(S)+\mu(S\setminus B)
			\leq\left|\widehat\nu(S)+\mu(S\setminus B)\right|
			<\eta\bigl(|\widehat\nu|(S)+\mu(S\setminus B)\bigr),
			\]
			and hence $\operatorname{Re}\widehat\nu(S)\leq\eta|\widehat\nu|(S)$. Therefore,
			\[
			\mu(B)=\operatorname{Re}\widehat\nu(B)
			\leq |\widehat\nu|(B\setminus\mathscr T_N(w))
			+\sum_{S\in\mathcal A_N(w)}\operatorname{Re}\widehat\nu(S)
			\leq |\widehat\nu|(B\setminus\mathscr T_N(w))+\eta|\widehat\nu|(B).
			\]
			Since $\mathscr T_N(w)\uparrow\mathscr T(w)$ and $|\widehat\nu|(B)<\infty$, letting $N\to\infty$ yields
			\[
			\mu(B)\leq |\widehat\nu|(B\setminus\mathscr T(w))+\eta|\widehat\nu|(B).
			\]
		}
		Combining the above inequalities, \eqref{eq:item_4_ref}, and using the definition of $\eta$, we obtain
		\[
		|{\widehat\nu}|\bigl(B \setminus \mathscr T(w)\bigr)\geq \mathfrak C_1^{-1}|{\widehat\nu}|(B) - \eta\,|{\widehat\nu}|(B)
		= \eta\,|{\widehat\nu}|(B),
		\]
		which immediately implies
		\[
		|{\widehat\nu}|\bigl(\mathscr T(w)\bigr)\leq (1-\eta)\,|{\widehat\nu}|(B).\qedhere
		\]
	\end{proof}
	
	\vv
	
	\subsubsection*{\underline{The sets of high and low density}}
	
	For  $\mathfrak{n} \in M(X; \mathbb C)$, $ \mathfrak m \in M_+(X)$, and $\lambda>0$, we define
	\begin{equation}\label{eq:def_e_lambda_mu_nu}
		E_{\lambda}(\mathfrak m,\mathfrak n) \coloneqq \{x\in B:\widetilde{\mathcal M}_{\mathfrak m}\mathfrak n(x)>  \lambda \}.
	\end{equation}
	
	Let $\mu$ and ${\widehat\nu}$ be as in assumptions \ref{ex_item_1}-\ref{ex_item_6}. 
	{ For $x\in B$, $m_1>\mathfrak C_1$, and $m_2>1$, we set}
	\begin{equation*}
		r_1(x) \coloneqq 
		\begin{cases}
			\sup \bigl\{r>0:\frac{|\widehat \nu|(B(x,r))}{\mu(B(x,5 r))}>m_1\bigr\}, \qquad &\text{ if }\widetilde{\mathcal M}_{\mu}\widehat \nu(x)>m_1,\\
			0, \qquad &\text{ otherwise,}
		\end{cases}
	\end{equation*}
	and 
	\begin{equation*}
		r_2(x) \coloneqq  
		\begin{cases}
			\sup \bigl\{r>0:\frac{{\mu(B(x,r)\cap B)}}{ |{\widehat\nu}|(B(x,5r))}>m_2\bigr\}, \qquad &\text{ if }\widetilde{\mathcal M}_{|{\widehat\nu}|}({\mu|_B})(x)>m_2,\\
			0, \qquad &\text{ otherwise.}
		\end{cases}
	\end{equation*}
	
	Finally, setting
	\[
	\widehat{\mathscr S}_1 \coloneqq E_{m_1}(\mu,\widehat \nu)\qquad\text{ and }\qquad \widehat{\mathscr S}_2 \coloneqq E_{m_2}(|{\widehat\nu}|,{\mu|_B}),
	\]
	we define the exceptional sets
	\[
	\mathscr S_1 \coloneqq \bigcup_{x \in \widehat{\mathscr S}_1} B(x,r_1(x))\qquad \text{ and }\qquad \mathscr S_2 \coloneqq \bigcup_{x \in \widehat{\mathscr S}_2} B(x,r_2(x)).
	\]
	{ Here the supremum of an empty set of admissible radii is defined to be zero, and $B(x,0)=\varnothing$. The radii are finite because the numerator measures are concentrated on $B$, while $m_1>\mathfrak C_1$ and $m_2>1$. We write $\mathscr S\coloneqq\mathscr S_1\cup\mathscr S_2$.}
	
	\vv
	
	\begin{lemma}\label{lemma:estim_nu_1}
		For $m_1\geq\mathfrak C_1/\varepsilon_0$ and $m_2\geq200\,\mathfrak C_1$, we have
		\begin{equation}\label{eq:meas_estim_H_1}
			|{\widehat\nu}|(\mathscr S_{1}) \leq \frac{\eta}{100}|{\widehat\nu}|(B)
			\qquad \text{and} \qquad
			|{\widehat\nu}|(\mathscr S_{2}) \leq \frac{\eta}{100}|{\widehat\nu}|(B).
		\end{equation}
	\end{lemma}
	\begin{proof}
		We denote
		\[
		E^1_{\lambda} \coloneqq E_\lambda(\mu, \widehat \nu)= \{x\in B : \widetilde{\mathcal M}_\mu\widehat \nu(x) > \lambda\},
		\]
		where $\widetilde{\mathcal M}_\mu$ is the maximal function defined in
		\eqref{eq:def_centered_maximal}.
		By \eqref{eq:measure_testing_maximal_constants} and \ref{ex_item_4}, we have
		\begin{equation}\label{eq:bound_S_1_1}
			\mu(E^1_{\lambda})\leq \frac{\|\widehat \nu\|}{\lambda}	\leq \frac{\mathfrak C_1\, \mu(B)}{\lambda}.
		\end{equation}
		
		{
			Let $y\in\mathscr S_1\cap B$. There exists $x\in B$ with $d(x,y)<r_1(x)$. By the definition of $r_1(x)$, we can choose $r>d(x,y)$ such that
			\begin{equation}\label{eq:bd_below_aux}
				|\widehat \nu|(B(x,r))>m_1\, \mu(B(x,5r)).
			\end{equation}
			Since $y\in B(x,r)$, the definition \eqref{eq:def_centered_maximal} gives $\widetilde{\mathcal M}_\mu\widehat \nu(y)>m_1$. Thus $\mathscr S_1\cap B\subset E^1_{m_1}$. Choosing $m_1\geq \mathfrak C_1/\varepsilon_0$, \eqref{eq:bound_S_1_1} gives $\mu(E^1_{m_1})\leq\varepsilon_0\mu(B)$. Since $\widehat\nu$ is concentrated on $B$, \eqref{eq:abs_cont} implies
			\begin{equation}\label{eq:aux_lemma_size_H_1_a}
				|\widehat\nu|(\mathscr S_1)=|\widehat\nu|(\mathscr S_1\cap B)\leq|\widehat\nu|(E^1_{m_1})\leq\frac{\eta}{100}\, |\widehat\nu|(B).
			\end{equation}
		}
		
		{
			To prove the second estimate in \eqref{eq:meas_estim_H_1}, define
			\[
			E^2_\lambda\coloneqq E_\lambda(|\widehat \nu|, \mu|_B)=\{x\in B:\widetilde{\mathcal M}_{|\widehat\nu|}(\mu|_B)(x)>\lambda\}.
			\]
			By \eqref{eq:measure_testing_maximal_constants}, we have
			\[
			|\widehat\nu|(E^2_\lambda)\leq\frac{\mu(B)}{\lambda}\leq\frac{ |\widehat\nu|(B)}{\lambda},
			\]
			where the last inequality follows from \ref{ex_item_3}. If $y\in\mathscr S_2\cap B$, there exists $x\in B$ with $d(x,y)<r_2(x)$. Choose $r>d(x,y)$ such that
			\[
			\mu(B(x,r)\cap B)>m_2\, |\widehat\nu|(B(x,5r)).
			\]
			Since $y\in B(x,r)$, we have $\widetilde{\mathcal M}_{|\widehat\nu|}(\mu|_B)(y)>m_2$. Thus $\mathscr S_2\cap B\subset E^2_{m_2}$ and, recalling that $\widehat\nu$ is concentrated on $B$, we conclude that
			\[
			|\widehat\nu|(\mathscr S_2)=|\widehat\nu|(\mathscr S_2\cap B)\leq|\widehat\nu|(E^2_{m_2})\leq\frac{|\widehat\nu|(B)}{m_2}.
			\]
			Taking $m_2\geq 200\mathfrak C_1=100/\eta$ proves the assertion.
		}
	\end{proof}
	
	\vv

	\begin{lemma}\label{lemma_RN_away_from_S}
		{ In the following formula, the Radon-Nikodym derivative is that of the restriction of $|\widehat\nu|$ to $B\setminus\mathscr S$ with respect to $\mu$.}
		Assume that  $\beta > 5^{\log_2 D}$ and that $m_1, m_2$ are as in Lemma \ref{lemma:estim_nu_1}. Then $|\widehat\nu||_{B\setminus\mathscr S}$ and $\mu|_{B\setminus\mathscr S}$ are mutually absolutely continuous, and
		\begin{equation}\label{eq:bounds_radon_nyk_derivative}
			{\frac{1}{\beta m_2}}\leq \frac{d|{\widehat\nu}|}{d\mu}(x)\leq \beta\,m_1,
			\qquad \text{for } |{\widehat\nu}|\text{-a.e. } x \in B \setminus \mathscr S.
		\end{equation}
	\end{lemma}
	
	\begin{proof}
		Put $G=B\setminus\mathscr S$. The definitions imply, for $x\in G$ and every sufficiently small $r>0$,
		\[
		|\widehat \nu|(B(x,r))\leq m_1\mu(B(x,5r)),\qquad\mu(B(x,r))\leq m_2|\widehat\nu|(B(x,5r)).
		\]
		We first justify absolute continuity on $G$. More generally, if positive Radon measures $\rho,\lambda$ satisfy $\rho(B(x,r))\leq M\lambda(B(x,5r))$ for all sufficiently small $r$ at each $x\in G$, then $\rho|_G\ll\lambda$. Indeed, let $A\subset G$ with $\lambda(A)=0$, and choose a bounded open set $O\supset A$ with $\lambda(O)<\varepsilon$. Fix $\beta'>25^{\log_2D}$. By \cite[Lemma 3.3]{Hy10}, for $\rho$-almost every $x\in A$ there are arbitrarily small $\rho$-$(25,\beta')$-doubling balls $B_x$ centered at $x$, with $25B_x\subset O$ and with the stated density inequality valid. Applying the $5r$-covering lemma to $5B_x$, choose pairwise disjoint $5B_i$ whose fivefold dilates cover $A$ up to a $\rho$-null set. Then
		\[
		\rho(A)\leq\sum_i\rho(25B_i)\leq\beta'\sum_i\rho(B_i)\leq\beta'M\sum_i\lambda(5B_i)\leq\beta'M\varepsilon.
		\]
		Letting $\varepsilon\downarrow0$ proves the claim. Applying it with $(\rho,\lambda,M)=(|\widehat\nu|,\mu,m_1)$ and with $(\rho,\lambda,M)=(\mu,|\widehat\nu|,m_2)$ gives $|\widehat\nu||_G\ll\mu$ and $\mu|_G\ll|\widehat\nu|$.
		
		Let $\varphi=d(|\widehat\nu||_G)/d\mu$, extended by zero off $G$. By \cite[Lemma 3.3 and Corollary 3.6]{Hy10}, at $\mu$-almost every $x\in G$ there are $(5,\beta)$-doubling balls $B_j=B(x,r_j)$ for $\mu$, with $r_j\to 0$ and $5B_j\subset B$, along which
		\[
		\begin{aligned}
			\varphi(x)
			&=\lim_{j\to\infty}\frac{|\widehat\nu|(B_j\cap G)}{\mu(B_j)}\leq\limsup_{j\to\infty}\frac{|\widehat\nu|(B_j)}{\mu(B_j)}\leq 
			m_1\, \limsup_{j\to \infty}\frac{\mu(5B_j)}{\mu(B_j)}\leq \beta\, m_1.
		\end{aligned}
		\]
		This inequality also holds $|\widehat \nu|$-almost everywhere on $G$.
		
		Finally, let $\psi=d(\mu|_G)/d|\widehat\nu|$. Differentiating $\psi$ along $(5,\beta)$-doubling balls $B_j=B(x,r_j)$ for $|\widehat\nu|$, with $5B_j\subset B$, gives for $|\widehat\nu|$-almost every $x\in G$,
		\[
		\psi(x)=\lim_{j\to\infty}
		\frac{\mu(B_j\cap G)}{|\widehat\nu|(B_j)}
		\leq\limsup_{j\to\infty}
		\frac{\mu(B_j)}{|\widehat\nu|(B_j)}
		\leq\beta m_2.
		\]
		Since the restricted measures are mutually absolutely continuous, $\psi=1/\varphi$ almost everywhere on $G$. This proves the lower bound in \eqref{eq:bounds_radon_nyk_derivative}.
	\end{proof}
	
	\vvv

	We next define the high- and low-density sets for the random grid $\bbD(w)$.
	Let $\widehat{\mathbf m}>1$ be a large parameter, to be chosen properly later. { We denote by $\mathcal S^{\widehat{\mathbf m},+}(w)$ the family of all cubes $Q\in\bbD(w)$ that satisfy
		\begin{equation}\label{eq:definition_ex_set_s_+}
			|\nu|(Q)>\widehat{\mathbf m}\,\mu(Q).
		\end{equation}
		Moreover, we define $\mathcal S^{\widehat{\mathbf m},-}(w)$ to be the collection of all cubes $Q\in\bbD(w)$ such that
		\begin{equation}\label{eq:definition_ex_set_s_-}
			|\nu|(Q)<\widehat{\mathbf m}^{-1}\mu(Q).
		\end{equation}
	}
	Finally, we define 
	\[
	\mathscr S^{\widehat{\mathbf m}, \pm}(w)\coloneqq \bigcup_{Q\in \mathcal S^{\widehat{\mathbf m},\pm}(w)}Q \qquad\text{ and }\qquad\mathscr S^{\widehat{\mathbf m}}(w)\coloneqq \mathscr S^{\widehat{\mathbf m}, +}(w)\cup \mathscr S^{\widehat{\mathbf m}, -}(w).
	\]
	We have the following estimate.
	\begin{lemma}\label{lemma:estim_S_m_hat}
		If $\widehat{\mathbf m}$ is chosen big enough depending on $\mathfrak C_1$ and $\varepsilon_0$, then
		\[
		|{\widehat\nu}|(\mathscr S^{\widehat{\mathbf m}}(w))\leq \frac{\eta}{50}\, |{\widehat\nu}|(B).
		\]
	\end{lemma}
	\begin{proof}
		For $N\in\mathbb N$, denote by $\mathcal S_N^{\widehat{\mathbf m},\pm}(w)$ the maximal cubes in $\mathcal S^{\widehat{\mathbf m},\pm}(w)$ whose generation is at least $-N$, and let $\mathscr S_N^{\widehat{\mathbf m},\pm}(w)$ be their unions. The cubes in each family are disjoint. If $Q\in\mathcal S_N^{\widehat{\mathbf m},+}(w)$, then \eqref{eq:definition_ex_set_s_+} gives
		\[
		|\widehat\nu|(Q)>\widehat{\mathbf m}\mu(Q\cap B)+(\widehat{\mathbf m}-1)\mu(Q\setminus B)\geq\widehat{\mathbf m}\, \mu(Q\cap B).
		\]
		Consequently,
		\[
		\mu(\mathscr S_N^{\widehat{\mathbf m},+}(w)\cap B)
		\leq\widehat{\mathbf m}^{-1}\sum_{Q\in\mathcal S_N^{\widehat{\mathbf m},+}(w)}|\widehat\nu|(Q)
		\leq\widehat{\mathbf m}^{-1}|\widehat\nu|(B)
		\leq\widehat{\mathbf m}^{-1}\mathfrak C_1\,\mu(B).
		\]
		Thus, for $\widehat{\mathbf m}\geq\mathfrak C_1\varepsilon_0^{-1}$, property \ref{ex_item_5} yields
		\[
		|\widehat\nu|\bigl(\mathscr S_N^{\widehat{\mathbf m},+}(w)\bigr)
		=|\widehat\nu|\bigl(\mathscr S_N^{\widehat{\mathbf m},+}(w)\cap B\bigr)
		\leq\frac{|\widehat\nu|(B)}{200\mathfrak C_1}
		\leq\frac{\eta}{100}|\widehat\nu|(B).
		\]
		For $Q\in\mathcal S_N^{\widehat{\mathbf m},-}(w)$, \eqref{eq:definition_ex_set_s_-} implies
		\[
		|\widehat\nu|(Q)<\widehat{\mathbf m}^{-1}\mu(Q\cap B)
		-(1-\widehat{\mathbf m}^{-1})\mu(Q\setminus B)
		\leq\widehat{\mathbf m}^{-1}\mu(Q\cap B).
		\]
		Therefore,
		\[
		|\widehat\nu|\bigl(\mathscr S_N^{\widehat{\mathbf m},-}(w)\bigr)
		\leq\widehat{\mathbf m}^{-1}\sum_{Q\in\mathcal S_N^{\widehat{\mathbf m},-}(w)}\mu(Q\cap B)
		\leq\widehat{\mathbf m}^{-1}\mu(B)
		\leq\frac{\eta}{100}|\widehat\nu|(B)
		\]
		whenever $\widehat{\mathbf m}\geq100/\eta$. The unions $\mathscr S_N^{\widehat{\mathbf m},\pm}(w)$ increase to $\mathscr S^{\widehat{\mathbf m},\pm}(w)$. Letting $N\to\infty$ proves the lemma for $\widehat{\mathbf m}\geq\max\{\mathfrak C_1\varepsilon_0^{-1},100/\eta\}$.
	\end{proof}
	\vv
	
	\subsubsection*{\underline{The set of non-Ahlfors regular balls}}
	
	Fix {$\mathbf m>0$}, to be chosen sufficiently large. For any $x\in X$ we define
	\begin{equation}\label{eq:definition_r_i}
		\mathcal{R}(x)\coloneqq 
		\begin{cases}
			\sup\bigl\{r>0: |{\widehat\nu}|(B(x,r))>\mathbf m\, r^n\bigr\}, & \text{if } M^r_{n}{\widehat\nu}(x)>\mathbf m,\\
			0, & \text{otherwise}.
		\end{cases}
	\end{equation}
	We remark that the radii $\mathcal R(x)$ are defined for every $x\in X$, and not on the support of ${\widehat\nu}$ only.
	We then define the sets
	\begin{equation}\label{eq:def_non_ahlfors}
		\widehat{\mathscr H}\coloneqq \{x\in X: M^r_{n}{\widehat\nu}(x)>\mathbf m\}
		\qquad\text{ and }\qquad
		\mathscr H\coloneqq \bigcup_{x\in \widehat{\mathscr H}} B(x,\mathcal{R}(x)).
	\end{equation}

	Thus, if $B(x,r)$ is a ball such that $|{\widehat\nu}|(B(x,r))>\mathbf m\, r^n$,
	then $B(x,r)\subset \mathscr H$.
	For $x\in X$ and $r>\mathcal R(x)$, the $n$-growth of $\mu$ also gives
	\[
	|\nu|(B(x,r))\leq|\widehat\nu|(B(x,r))+\mu(B(x,r))\leq(\mathbf m+\mathfrak C_0)r^n.
	\]
	We also have the following estimate for the $|{\widehat\nu}|$-measure of the
	exceptional set $\mathscr H$.
	
	\begin{lemma}\label{lemma:estim_ex_H}
		For {$\mathbf m$} chosen sufficiently large depending on
		$n, \varepsilon_0, \mathfrak C_0$, and $\mathfrak C_1$, we have
		\[
		|{\widehat\nu}|(\mathscr H)\leq \frac{\eta}{100}\, |{\widehat\nu}|(B).
		\]
	\end{lemma}
	\begin{proof}
		{
			For $\lambda>0$, set $E_\lambda=\{M^r_n\widehat\nu>\lambda\}$. By \eqref{eq:measure_testing_maximal_constants} and \ref{ex_item_4},
			\[
			\mu(E_\lambda)\leq\frac{5^n\, \mathfrak C_0}{\lambda}\|\widehat\nu\|
			\leq\frac{5^n\, \mathfrak C_0\, \mathfrak C_1}{\lambda}\, \mu(B).
			\]
			
			If $y\in\mathscr H$, we choose $x$ such that $d(x,y)<\mathcal R(x)$. By \eqref{eq:definition_r_i}, there exists $r>d(x,y)$ such that $|\widehat\nu|(B(x,r))>\mathbf m \, r^n$. Then
			\begin{equation}\label{eq:aux_lemma_size_H_2_i_2}
				|\widehat\nu|(B(y,2r))\geq|\widehat\nu|(B(x,r))	>\mathbf m \, r^n=2^{-n}\mathbf m\, (2r)^n.
			\end{equation}
			Thus $\mathscr H\subset E_{\mathbf m/2^n}$. For $\mathbf m\geq 10^n\, \mathfrak C_0\, \mathfrak C_1/\varepsilon_0$, we have $\mu(B\cap E_{\mathbf m/2^n})\leq\varepsilon_0\mu(B)$. Since $\widehat\nu$ is concentrated on $B$, the property \ref{ex_item_5} implies
			\[
			|\widehat\nu|(\mathscr H)
			\leq|\widehat\nu|(B\cap E_{\mathbf m/2^n})
			\leq\frac{\|\widehat\nu\|}{200\mathfrak C_1}
			\leq\frac{\eta}{100}|\widehat\nu|(B).\qedhere
			\]
		}
	\end{proof}
	
	\vv
	
	\subsubsection*{\underline{The set where $T{\widehat\nu}$ is large}}

	{
		{
			Fix $\alpha>1$, to be chosen sufficiently large later. For $x\in B$, we define
			\begin{equation}\label{eq:def_e_i}
				e(x)\coloneqq
				\begin{cases}
					\displaystyle\sup\{\varepsilon>0:|T_\varepsilon\widehat\nu(x)|>\alpha\},
					&\text{if }T_\sharp\widehat\nu(x)>\alpha,\\[1ex]
					0,    &\text{otherwise}.
				\end{cases}
			\end{equation}
			Note that these radii are finite, since $|T_\varepsilon\widehat\nu(x)|\leq C\varepsilon^{-n}\|\widehat\nu\|$. Next, we set
			\begin{equation}\label{eq:exceptional_set_S_i}
				\widehat{\mathscr E}\coloneqq\{x\in B:T_\sharp\widehat\nu(x)>\alpha\}
				\qquad\text{ and }\qquad
				\mathscr E\coloneqq\bigcup_{x\in\widehat{\mathscr E}}B(x,e(x)).
			\end{equation}
			The proof of the following lemma follows the scheme of \cite[Lemma 5.2]{To14}.
			
			\begin{lemma}\label{lemma:measures_exceptional_set_S}
				If $y\in\mathscr E\setminus\mathscr H$, then
				\[
				T_\sharp\widehat\nu(y)>\alpha-\alpha'_E,
				\qquad \text{ for }\quad
				\alpha'_E\coloneqq C\mathbf m\,2^n\frac{3\cdot2^\gamma-2}{2^\gamma-1}.
				\]
				Moreover, if $\alpha\geq2\max\{\alpha'_E,\mathfrak C_2\}$, then $B\cap\mathscr E\subset\mathscr H\cup U_B$.
			\end{lemma}
			\begin{proof}
				Let $y\in\mathscr E\setminus\mathscr H$ and choose $x\in\widehat{\mathscr E}$ with $d(x,y)<e(x)$. By \eqref{eq:def_e_i}, there exists $\varepsilon$ such that
				\[
				d(x,y)<\varepsilon\leq e(x)\qquad \text{ and }\qquad  |T_\varepsilon\widehat\nu(x)|>\alpha.
				\]
				The triangle inequality gives
				\begin{equation}\label{eq:estim_S_aux_0}
					\begin{split}
						|T_\varepsilon\widehat\nu(x)-T_\varepsilon\widehat\nu(y)|
						&\leq |T_\varepsilon(\widehat\nu|_{B(y,2\varepsilon)})(x)|
						+|T_\varepsilon(\widehat\nu|_{B(y,2\varepsilon)})(y)|\\
						&\quad+|T_\varepsilon(\widehat\nu|_{X\setminus B(y,2\varepsilon)})(x)
						-T_\varepsilon(\widehat\nu|_{X\setminus B(y,2\varepsilon)})(y)|\\
						&\eqqcolon\mathcal I_1+\mathcal I_2+\mathcal I_3.
					\end{split}
				\end{equation}
				Since $y\notin\mathscr H$, all balls centered at $y$ satisfy the $n$-growth bound with constant $\mathbf m$. Hence the standard kernel estimates yield
				\[
				\mathcal I_1+\mathcal I_2
				\leq2C\varepsilon^{-n}|\widehat\nu|(B(y,2\varepsilon))
				\leq2C\mathbf m\,2^n.
				\]
				For $z\notin B(y,2\varepsilon)$, we have $d(x,z)\geq d(y,z)-d(x,y)>\varepsilon$ and $d(y,z)>2d(x,y)$. Hence the regularity of $K$ gives
				\[
				\begin{aligned}
					\mathcal I_3
					&\leq C\varepsilon^\gamma
					\int_{X\setminus B(y,2\varepsilon)}
					\frac{d|\widehat\nu|(z)}{d(y,z)^{n+\gamma}}\leq C\mathbf m\,\frac{2^{n+\gamma}}{2^\gamma-1}.
				\end{aligned}
				\]
				
				So, we have
				\[
				|T_\varepsilon\widehat\nu(y)|
				\geq |T_\varepsilon\widehat\nu(x)|-\alpha'_E
				>\alpha-\alpha'_E,
				\]
				which proves the first assertion. If $\alpha\geq2\max\{\alpha'_E,\mathfrak C_2\}$ and $y\in B\cap\mathscr E\setminus\mathscr H$, the first assertion gives $T_\sharp\widehat\nu(y)>\alpha-\alpha'_E\geq\alpha/2\geq\mathfrak C_2$, so $y\in U_B$ by \ref{ex_item_6}.
			\end{proof}
		}
	}
	
	\vvv
	We remark that the thresholds $m_1$, $m_2$, $\widehat{\mathbf m}$, $\mathbf m$, and $\alpha$ depend only on $n$, $C$, $\gamma$, $\mathfrak C_0$, $\mathfrak C_1$, $\mathfrak C_2$, and $\varepsilon_0$. In particular, by Lemma \ref{lemma:reduction_cutoffs}, they do not depend on $h$ and $L$, and the same holds for the constant $c_\sharp$ in Lemma \ref{eq:lemm_bd_suppressed_2}.
	
	\subsubsection*{\underline{The exceptional set $\mathscr W^0(w)$}} We now collect all the exceptional sets constructed so far, and define
	\begin{equation}\label{eq:def_tot_exceptional_set}
		\mathscr W^0(w)\coloneqq
		\mathscr T(w)\cup \mathscr S \cup \mathscr S^{\widehat{\mathbf m}}(w)\cup \mathscr H \cup \mathscr E \cup U_B.
	\end{equation}
	\vv
	\begin{lemma}\label{lemma:final_lemma_exceptional_set}
		Let us assume that the constants in the constructions of $\mathscr S$, $\mathscr S^{\widehat{\mathbf m}}(w)$, $\mathscr H$, and  $\mathscr E$ are chosen as in
		Lemmas \ref{lemma:estim_T_i}, \ref{lemma:estim_nu_1}, \ref{lemma:estim_S_m_hat}, \ref{lemma:estim_ex_H}, and \ref{lemma:measures_exceptional_set_S}.
		Then, for $\delta_1\coloneqq 1-\tfrac{9}{10}\eta  =1-\tfrac{9}{10} (2\mathfrak C_1)^{-1}$, we have that 
		\begin{equation}\label{eq:estim_aux_prob1}
			|{\widehat\nu}|\bigl(B\setminus \mathscr W^0(w)\bigr)	\geq (1-\delta_1)\,|{\widehat\nu}|(B),\qquad \text{for all } w\in\Omega.
		\end{equation}
		Moreover, 
		\begin{equation}\label{eq:estim_aux_prob2}
			\mu\bigl(B\setminus \mathscr W^0(w)\bigr)\geq (\beta\,m_1)^{-1}\,(1-\delta_1)\,\mu(B),\qquad \text{for all } w\in\Omega.
		\end{equation}
	\end{lemma}
	\begin{proof}
		Since $\widehat\nu$ is concentrated on $B$, assumptions \ref{ex_item_5} and \ref{ex_item_6} give $|\widehat\nu|(U_B)\leq\|\widehat\nu\|/(200\mathfrak C_1)=\eta\,|\widehat\nu|(B)/100$, and Lemma \ref{lemma:measures_exceptional_set_S} gives $B\cap\mathscr E\subset\mathscr H\cup U_B$, so that $\mathscr E$ does not increase the $|\widehat\nu|$-measure of $\mathscr W^0(w)$. Hence, by Lemmas \ref{lemma:estim_T_i}, \ref{lemma:estim_nu_1}, \ref{lemma:estim_S_m_hat}, and \ref{lemma:estim_ex_H}, we have
		\[
		|{\widehat\nu}|(\mathscr W^0(w))\leq \Bigl(1-\eta+\frac{6}{100}\eta\Bigr)|{\widehat\nu}|(B).
		\]
		Thus, the choice $\delta_1= 1-\tfrac{9}{10}\eta$ concludes the proof of \eqref{eq:estim_aux_prob1}, which, in turn, by \eqref{eq:bounds_radon_nyk_derivative} and \ref{ex_item_3}, implies \eqref{eq:estim_aux_prob2}.
	\end{proof}
	
	\vv

	\vv

	\subsubsection*{\underline{The joint stopping family}}

	{
		{
			Suppose that $w\in\Omega$ and $Q_0\in\bbD(w)$ satisfy $B\subset Q_0$. With $\widehat{\mathbf m}$ chosen as above, denote by $\Stop(Q_0,w)$ the maximal cubes in
			\[
			\bigl\{S\in\mathcal A(w)\cup\mathcal S^{\widehat{\mathbf m},+}(w)\cup\mathcal S^{\widehat{\mathbf m},-}(w):S\subsetneq Q_0\bigr\}.
			\]
			Thus a cube belongs to one of the stopping families if it satisfies the non-accretivity condition or one of \eqref{eq:definition_ex_set_s_+} and \eqref{eq:definition_ex_set_s_-}.
			
			\begin{lemma}\label{lemma:lemma38}
				The cubes in $\Stop(Q_0,w)$ are pairwise disjoint and it holds
				\[
				\bigcup_{S\in\Stop(Q_0,w)}S	=\mathscr T(w)\cup\mathscr S^{\widehat{\mathbf m}}(w).
				\]
				In particular, this union does not depend on the choice of $Q_0$. If $S\in\Stop(Q_0,w)$ and $S^*$ is its dyadic parent, then
				\[
				|\nu(S^*)|\geq\eta\, |\nu|(S^*)\qquad \text{ and }\qquad \widehat{\mathbf m}^{-1}\mu(S^*)\leq|\nu|(S^*)\leq\widehat{\mathbf m}\mu(S^*).
				\]
				The same inequalities hold for every cube $Q\subset Q_0$ that is not contained in a cube of $\Stop(Q_0,w)$.
			\end{lemma}
			\begin{proof}
				Let $\mathcal F\coloneqq\mathcal A(w)\cup\mathcal S^{\widehat{\mathbf m},+}(w)\cup\mathcal S^{\widehat{\mathbf m},-}(w)$. If $Q\in\bbD(w)$ is disjoint from $B$, then $\nu=|\nu|=\mu$ on $Q$, so $Q\notin\mathcal F$.
				
				On the other hand, if $B\subset Q$, the normalization $\widehat\nu(B)=\mu(B)$ and \eqref{eq:item_4_ref} yield
				\[
				\nu(Q)=\mu(Q)\leq|\nu|(Q)\leq\mathfrak C_1\, \mu(Q),
				\]
				so $Q\notin\mathcal F$, because $\eta=(2\mathfrak C_1)^{-1}$ and $\widehat{\mathbf m}\geq\mathfrak C_1$. Hence, every $S\in\mathcal F$ meets $B\subset Q_0$ and does not contain $Q_0$, so that $S\subsetneq Q_0$. Since the cubes between $S$ and $Q_0$ form a finite chain, $S$ is contained in a maximal cube of $\mathcal F$, and the maximal cubes of $\mathcal F$, which are pairwise disjoint, form $\Stop(Q_0,w)$. This proves the identity in the statement.
				
				If $S\in\Stop(Q_0,w)$, its parent $S^*\subset Q_0$ is not in $\mathcal F$, either because $S^*=Q_0$ or by maximality. Similarly, a cube $Q\subset Q_0$ not contained in a cube of $\Stop(Q_0,w)$ is not in $\mathcal F$. For such cubes, the asserted inequalities follow from the definitions of the three families if they meet $B$, and from $\nu=|\nu|=\mu$ on them otherwise.
			\end{proof}
			By Lemma \ref{lemma:lemma38} and \eqref{eq:def_tot_exceptional_set}, we have the set inclusion
			\[
			\bigcup_{S\in\Stop(Q_0,w)}S\subset\mathscr W^0(w).
			\]
		}
	}
	
	\vv
	\subsubsection*{\underline{The total exceptional set}}
	Let $R\in\caD^{db}_0$, where $\caD^{db}_0$ is the family of doubling cells of the David-Mattila lattice $\caD_0$ introduced at the beginning of Section \ref{sec:corona}, and let {$B\coloneqq\widetilde B(R)$ be its associated testing ball, see Lemma \ref{lemma:construction_testing_ball_DM}.}
	Let $\mathscr W^0(w)$ be the set constructed in \eqref{eq:def_tot_exceptional_set}, let $\theta\coloneqq \min\{1/4,\varepsilon_0/(2\beta\mathfrak C_3)\}$, and let ${ B_\theta\coloneqq(1-\theta)B}$. We define the \emph{total exceptional set}
	\begin{equation}\label{eq:W_new}
		{\mathscr W_R(w)\coloneqq\mathscr W^0(w),}	
	\end{equation}
	and we denote $G_R(w)\coloneqq B_\theta\setminus {\mathscr W_R}(w).$
	\begin{lemma}\label{lemma:exceptional_set_new_contraction}
		We put $\theta'\coloneqq \mathfrak C_3\beta\theta$, so that $\theta'\leq \varepsilon_0/2$. Then $B_\theta$ is $\mu$-$(5,\beta/(1-\theta'))$-doubling, and
		\begin{equation}\label{eq:meas_good_set_nu}
			|\widehat\nu|(G_R(w))\geq\frac{1-\delta_1}{2}\, |\widehat\nu|(\widetilde B(R)),
			\qquad w\in\Omega,
		\end{equation}
		\begin{equation}\label{eq:meas_good_set_mu}
			\mu(G_R(w))\geq(\beta m_1)^{-1}\frac{1-\delta_1}{2}\, \mu(\widetilde B(R)),
			\qquad w\in\Omega.
		\end{equation}
	\end{lemma}
	\begin{proof}
		{ The small-boundary and doubling properties of $B$ give
			\begin{equation}\label{eq:small_bd_tilde_B0}
				\mu(B\setminus B_\theta)\leq\mathfrak C_3\theta\, \mu(5B)\leq\theta'\mu(B).
			\end{equation}
			Consequently,
			\[
			\mu(B_\theta)\geq(1-\theta')\, \mu(B)\geq\frac{1-\theta'}{\beta}\, \mu(5B)\geq\frac{1-\theta'}{\beta}\, \mu(5B_\theta),
			\]
			which proves the doubling assertion. Since $\theta'\leq\varepsilon_0$, assumption \ref{ex_item_5} yields
			\begin{equation}\label{eq:small_bd_nu_tilde_b}
				|\widehat\nu|(B\setminus B_\theta)
				\leq\frac{\|\widehat\nu\|}{200\mathfrak C_1}
				\leq\frac{\eta}{100}|\widehat\nu|(B).
			\end{equation}
			In particular,
			\begin{equation}\label{eq:nu_1_B_0_tilde_B_2}
				|\widehat\nu|(B_\theta)\geq \Bigl(1-\frac{\eta}{100}\Bigr)|\widehat\nu|(B).
			\end{equation}
			Note that $\mathscr W_R(w)=\mathscr W^0(w)$, so Lemma \ref{lemma:final_lemma_exceptional_set} gives
			\[
			\begin{aligned}
				|\widehat\nu|(G_R(w))
				&\geq|\widehat\nu|(B\setminus\mathscr W^0(w))
				-|\widehat\nu|(B\setminus B_\theta)\\
				&\geq\Bigl(1-\delta_1-\frac{\eta}{100}\Bigr)\, |\widehat\nu|(B)
				\geq\frac{1-\delta_1}{2}|\widehat\nu|(B),
			\end{aligned}
			\]
			since $1-\delta_1=9\eta/10$. This proves \eqref{eq:meas_good_set_nu}. The bound \eqref{eq:meas_good_set_mu} follows from \eqref{eq:bounds_radon_nyk_derivative} and \ref{ex_item_3}.}
	\end{proof}
	
	\vv

	\subsection{Suppressed maximal functions and singular integral operators}\label{sec:aux_operator}
	{
		Let $\Xi\colon X\to[0,\infty)$ be $1$-Lipschitz, that is,
		\[
		|\Xi(x)-\Xi(y)|\leq d(x,y),\qquad \text{ for } x,y\in X.
		\]
		For a complex Radon measure $\mathfrak n$ on $X$, we define the suppressed radial maximal function by
		\[
		M^r_\Xi(\mathfrak n)(x)\coloneqq\sup_{\substack{\rho>0, \rho\geq\Xi(x)}}\frac{|\mathfrak n|(B(x,\rho))}{\rho^n}.
		\]
		For $f\in L^1_{\loc}(\mathfrak m)$, we write $M^r_{\Xi,\mathfrak m}f\coloneqq M^r_\Xi(f\mathfrak m)$. The proof of the next lemma is standard, but we include it for the reader's convenience.		
		\begin{lemma}\label{lem:suppressed maximal}
			Let $\mathfrak m_1,\mathfrak m_2$ be positive Radon measures on $X$ with maximal non-Ahlfors radii $\mathcal R_1,\mathcal R_2$ and corresponding growth constants $c_1,c_2$, so that
			\[
			\mathfrak m_j(B(x,r))\leq c_jr^n\qquad\text{ if }r>\mathcal R_j(x),\quad j=1,2.
			\]
			If $\Xi\geq\max\{\mathcal R_1,\mathcal R_2\}$, then $M^r_{\Xi,\mathfrak m_2}$ maps $L^1(\mathfrak m_2)$ to $L^{1,\infty}(\mathfrak m_1)$ and maps $L^p(\mathfrak m_2)$ to $L^p(\mathfrak m_1)$ for $1<p<\infty$. The constants depend only on $D,n,c_1,c_2$ and, for the latter estimate, on $p$.
		\end{lemma}
		\begin{proof}
			The growth estimates hold also at a positive radius $r=\mathcal R_j(x)$, by taking larger radii decreasing to $r$. For a level set of $M^r_{\Xi,\mathfrak m_2}f$, the $5r$-covering lemma selects disjoint balls $B_i\coloneqq B(x_i,r_i)$ with $r_i\geq\Xi(x_i)$ and
			\[
			\lambda r_i^n<\int_{B_i}|f|\,d\mathfrak m_2.
			\]
			The corresponding dilated balls form a covering of the level set so, applying the argument first to bounded collections of radii gives
			\begin{equation}\label{eq:supp_max_f_weak}
				\mathfrak m_1\bigl(\{M^r_{\Xi,\mathfrak m_2}f>\lambda\}\bigr)\leq\sum_i\mathfrak m_1(5B_i)\leq5^nc_1\sum_i r_i^n
				\leq\frac{5^nc_1}{\lambda}\, \|f\|_{L^1(\mathfrak m_2)}.
			\end{equation}
			The $L^p(\mathfrak m_2)\to L^p(\mathfrak m_1)$ boundedness follows by interpolation of \eqref{eq:supp_max_f_weak}  and the trivial $L^\infty(\mathfrak m_2)\to L^\infty(\mathfrak m_1)$ boundedness.
		\end{proof}
		
		Let $K$ be an $n$-dimensional Calder\'on-Zygmund kernel with constants $C$ and $\gamma$.
		We define
		\begin{equation}\label{eq:suppressed_kernel}
			K_\Xi(x,y)\coloneqq K(x,y)\, \frac{d(x,y)^{2n}}{d(x,y)^{2n}+\Xi(x)^n\Xi(y)^n},
		\end{equation}
		and
		\begin{equation}\label{eq:SIO operator}
			T_{\Xi,\mathfrak m,\varepsilon}f(x)	\coloneqq\int_{d(x,y)>\varepsilon}K_\Xi(x,y)f(y)\,d\mathfrak m(y).
		\end{equation}
		Note that, if $K(x,y)=-K(y,x)$, then \eqref{eq:suppressed_kernel} gives $K_\Xi(x,y)=-K_\Xi(y,x)$. 
		We also define the maximal suppressed SIO and its truncated version as
		\[
		T_{\Xi,\mathfrak m,\sharp}f(x)\coloneqq\sup_{\varepsilon>0}|T_{\Xi,\mathfrak m,\varepsilon}f(x)|,	\qquad \text{ and }\qquad T_{\Xi,\mathfrak m,\sharp,\delta}f(x)\coloneqq\sup_{\varepsilon>\delta}|T_{\Xi,\mathfrak m,\varepsilon}f(x)|.
		\]
		For positive Radon measures $\mathfrak m,\mathfrak n$, we understand that $T_{\Xi,\mathfrak m}$ is bounded from $L^2(\mathfrak m)$ to $L^2(\mathfrak n)$  if truncations are  uniformly bounded, with
		\[
		\|T_{\Xi,\mathfrak m}\|_{L^2(\mathfrak m)\to L^2(\mathfrak n)}\coloneqq\sup_{\varepsilon>0}\|T_{\Xi,\mathfrak m,\varepsilon}\|_{L^2(\mathfrak m)\to L^2(\mathfrak n)}.
		\]
		In the case of a complex measure $\mathfrak n$ on $X$, we also write $T_{\Xi,\varepsilon}\mathfrak n\coloneqq T_{\Xi,\mathfrak n,\varepsilon}1$ and $T_{\Xi,\sharp}\mathfrak n\coloneqq T_{\Xi,\mathfrak n,\sharp}1$.

		\begin{lemma}\label{lemma:kernel_suppressed}
			The kernel $K_\Xi$ in \eqref{eq:suppressed_kernel} has the following properties, with constants depending only on $n$, $\gamma$, and the constants of $K$.
			\begin{enumerate}
				\item $K_\Xi$ is an $n$-dimensional Calder\'on-Zygmund kernel with exponent $\gamma$ and constants depending only on $n$, $\gamma$, and $C$.
				\item If $\max\{\Xi(x),\Xi(y)\}>0$, then
				\begin{equation}\label{eq:bound_suppressed_kernel_1}
					|K_\Xi(x,y)|\leq\frac{C}{\max\{\Xi(x)^n,\Xi(y)^n\}}.
				\end{equation}
				\item If $\Xi(x)=0$ or $\Xi(y)=0$, then $K_\Xi(x,y)=K(x,y)$. Moreover, for every $\mathfrak n\in M(X;\mathbb C)$, every $x\in X$, and every $\varepsilon\geq\Xi(x)>0$, we have
				\begin{equation}\label{eq:bound_suppressed_kernel_2}
					|T_{\Xi,\varepsilon}\mathfrak n(x)-T_\varepsilon\mathfrak n(x)|
					\leq C\sup_{r>\varepsilon}\frac{|\mathfrak n|(B(x,r))}{r^n}.
				\end{equation}
			\end{enumerate}
		\end{lemma}
		
		\begin{proof}
			Since $\mu$ has $n$-growth with constant $\mathfrak C_0$, Appendix \ref{appendix:suppression_metric_spaces} applies to $(X,d,\mathfrak C_0^{-1}\mu)$ with the dominating function $\lambda(x,r)\coloneqq r^n$, and, by Definition \ref{def:CZ_kernel}, $K$ satisfies \eqref{eq:size} and \eqref{eq:holder_alpha}. In this case, the function $\Lambda$ of Lemma \ref{lemma:aux_dominating_function_3} is $\Lambda(x,r)=c_nr^n$, for some $c_n>0$, so we can take $d_\lambda=n$ and, since $n\geq1$, $M=n$. Then $K_\Xi$ in \eqref{eq:suppressed_kernel} is the kernel of Definition \ref{def:B-product-suppression}, and (1) follows from \eqref{eq:B-usual-size} and \eqref{eq:B-product-regularity}, (2) from \eqref{eq:B-suppressed-size}, and (3) from \eqref{eq:B-kernel-one-zero} and Lemma \ref{lem:B-hard-comparison} with $t=\varepsilon$.
		\end{proof}
		
		\vv
		
		We next estimate $T_{\Xi,\varepsilon}\mathfrak n$ in terms of $T_\varepsilon\mathfrak n$.
		\begin{lemma}\label{lemma:lemm_bd_suppressed_1}
			Let $\Xi\colon X\to[0,\infty)$ be $1$-Lipschitz and let $\mathfrak n\in M(X;\mathbb C)$. Suppose that $\alpha,m>0$, $x\in X$, and $R_0\geq0$ satisfy 	for $\varepsilon>R_0$ and $r>R_0$ the conditions
			\[
			|\mathfrak n|(B(x,r))\leq m\, r^n\qquad \text{ and }\qquad |T_\varepsilon\mathfrak n(x)|\leq\alpha.
			\]
			
			If $\Xi(x)\geq R_0$, then
			\begin{equation}\label{eq:lemm_bd_suppressed_1}
				|T_{\Xi,\varepsilon}\mathfrak n(x)|\leq\alpha + C m,\qquad \text{ for }\varepsilon>0.
			\end{equation}
		\end{lemma}
		
		\begin{proof}
			Let $h\coloneqq\Xi(x)$. If $h=0$, then $R_0=0$ and $T_{\Xi,\varepsilon}\mathfrak n(x)=T_\varepsilon\mathfrak n(x)$, so the assertion follows. If $\varepsilon>h>0$, the bound \eqref{eq:bound_suppressed_kernel_2} gives
			\begin{equation}\label{eq:lemm_bd_suppressed_3}
				|T_{\Xi,\varepsilon}\mathfrak n(x)|\leq|T_\varepsilon\mathfrak n(x)|+C\, \sup_{r>\varepsilon}\frac{|\mathfrak n|(B(x,r))}{r^n}
				\leq\alpha+Cm.
			\end{equation}
			If $0<\varepsilon\leq h$, then
			\[
			\begin{aligned}
				|T_{\Xi,\varepsilon}\mathfrak n(x)|
				&\leq\int_{\varepsilon<d(x,y)\leq2h}|K_\Xi(x,y)|\,d|\mathfrak n|(y)+\bigl|T_{\Xi,2h}\mathfrak n(x)\bigr|\\
				&\leq Ch^{-n}|\mathfrak n|(B(x,3h))
				+\alpha+Cm
				\leq\alpha+Cm.
			\end{aligned}
			\]
			Here we used \eqref{eq:bound_suppressed_kernel_1}, $|\mathfrak n|(B(x,3h))\leq m(3h)^n$, and \eqref{eq:lemm_bd_suppressed_3} with $\varepsilon=2h$.
		\end{proof}
		
		Let $\mu$ and $\widehat \nu$ satisfy assumptions \ref{ex_item_1}-\ref{ex_item_6}, with constants in the exceptional-set construction chosen as in Lemma \ref{lemma:final_lemma_exceptional_set}. Recall that the radii $\mathcal R$ and $e$ are defined in \eqref{eq:definition_r_i} and \eqref{eq:def_e_i}, respectively. If the suppression is such that
		\begin{equation}\label{eq:condition_suppression}
			\Xi(x)\geq\dist\bigl(x,B\setminus(\mathscr H\cup\mathscr E)\bigr) \qquad \text{ for every }\qquad x\in X, 
		\end{equation}
		then the inclusions $B(x,\mathcal R(x))\subset\mathscr H$ and $B(x,e(x))\subset\mathscr E$ imply  the bound from below
		\[
		\Xi(x)\geq\max\{\mathcal R(x),e(x)\},\qquad \text{ for } x\in B.
		\]
		By \eqref{eq:definition_r_i} and \eqref{eq:def_e_i}, Lemma \ref{lemma:lemm_bd_suppressed_1} applies with $R_0=\max\{\mathcal R(x),e(x)\}$, $m=\mathbf m$, and $\alpha$ as in Lemma \ref{lemma:measures_exceptional_set_S}.
		
		\begin{lemma}\label{eq:lemm_bd_suppressed_2}
			Suppose that $\Xi\colon X\to[0,\infty)$ is $1$-Lipschitz and satisfies \eqref{eq:condition_suppression}. Then there is $c_\sharp>0$, depending only on $n,D,\gamma,\varepsilon_0,C,\mathfrak C_0,\mathfrak C_1,\mathfrak C_2$, such that
			\[
			T_{\Xi,\sharp}\widehat\nu(x)\leq c_\sharp,
			\qquad x\in B.
			\]
		\end{lemma}
	}
	
	\vv

	\subsection{A Cotlar inequality for different measures}\label{subsection:mixed_cotlar}
	The following SIO boundedness result is a variant of \cite[Lemma 5.25]{To14} for different measures, formulated in geometrically doubling metric measure spaces, and it will be applied to the measures $\mu$ and $|\nu_i|=|\widehat\nu_i|+\mu|_{X\setminus B_i}$, for $i=1,2$, in the proof of Lemma \ref{lemma:properties_regularized_cubes}.
	
	\begin{lemma}\label{lemma:mixed_cotlar_transfer}
		Let $n>0$ and $0<\gamma\leq1$, and let $(X,d)$ be a geometrically doubling metric space with constant $D$. Let $\Xi\colon X\to[0,\infty)$ be a $1$-Lipschitz function, and let $\mathfrak m_1$, $\mathfrak m_2$, and $\mathfrak m_3$ be non-negative Radon measures on $X$. Assume that there exist constants $C,A\geq1$ such that, for every $x\in X$, every $r>0$ with $r\geq\Xi(x)$, and all $i,j\in\{1,2,3\}$, we have
		\begin{equation}\label{eq:mct_measure_bounds}
			\mathfrak m_i(B(x,r))\leq C\,r^n\qquad \text{ and }\qquad	\mathfrak m_i(B(x,r))\leq C\,\mathfrak m_j(B(x,Ar)).
		\end{equation}
		Let $k$ be an $n$-dimensional Calder\'on-Zygmund kernel with constants $C$ and $\gamma$, with regularity in both variables, and suppose that
		\begin{equation}\label{eq:mct_kernel_bounds}
			|k(x,y)|\leq C\min\bigl\{d(x,y)^{-n},\Xi(x)^{-n},\Xi(y)^{-n}\bigr\},
		\end{equation}
		with the formal convention $0^{-n}=+\infty$. Assume also that there exists $\delta_0>0$ such that $k(x,y)=0$ whenever $d(x,y)\leq\delta_0$. For every bounded function $f$ with bounded support, we define
		\[
		S_{\mathfrak m_1}f(x)\coloneqq\int_X k(x,y)f(y)\,d\mathfrak m_1(y),
		\]
		\[
		S_{\mathfrak m_1, \delta} f(x)\coloneqq\int_{d(x,y)>\delta}k(x,y)f(y)\,d\mathfrak m_1(y),\qquad \text{ for }\delta>0,
		\]
		and $S_{\mathfrak m_1, \sharp} f(x)\coloneqq\sup_{\delta>0}|S_{\mathfrak m_1, \delta} f(x)|$. Suppose that $S_{\mathfrak m_1}$ extends to a bounded operator from $L^2(\mathfrak m_1)$ to $L^2(\mathfrak m_2)$, and denote $\|S_{\mathfrak m_1}\|_{\op}\coloneqq\|S_{\mathfrak m_1}\|_{L^2(\mathfrak m_1)\to L^2(\mathfrak m_2)}$. Then, we have
		\begin{equation}\label{eq:mct_weak_bound}
			\|S_{\mathfrak m_1}\|_{L^1(\mathfrak m_1)\to L^{1,\infty}(\mathfrak m_2)}\lesssim (1+\|S_{\mathfrak m_1}\|_{\op})^2
		\end{equation}
		and
		\begin{equation}\label{eq:mct_maximal_bound}
			\|S_{\mathfrak m_1,\sharp} f\|_{L^2(\mathfrak m_3)}
			\lesssim(1+\|S_{\mathfrak m_1}\|_{\op})\, \|f\|_{L^2(\mathfrak m_1)},
			\qquad \text{ for }f\in L^2(\mathfrak m_1).
		\end{equation}
		The implicit constants depend only on $n$, $D$, $\gamma$, $C$, and $A$, and they are independent of $\delta_0$.
	\end{lemma}
	
	{
		The proof of Lemma \ref{lemma:mixed_cotlar_transfer} follows the ideas in \cite[Lemmas 5.25-5.27]{To14}, and it combines the non-homogeneous Calder\'on-Zygmund decomposition and Cotlar inequality (see \cite[Chapter 2]{To14}) with the comparisons in \eqref{eq:mct_measure_bounds}. We give the details only where the suppression $\Xi$ and the use of three different measures play a role.
		
		Throughout the rest of this subsection, we assume that the hypotheses of Lemma \ref{lemma:mixed_cotlar_transfer} hold, and, to simplify the notation, we write $S\coloneqq S_{\mathfrak m_1}$, $S_\delta\coloneqq S_{\mathfrak m_1,\delta}$, $S_\sharp\coloneqq S_{\mathfrak m_1,\sharp}$, and $\|S\|_{\op}\coloneqq\|S_{\mathfrak m_1}\|_{\op}$. We denote by $c\geq1$ the constant associated with the finite family of adjacent dyadic lattices in \cite[Theorem 4.1 and Remark 4.13]{HK12}, so that every ball $B$ is contained in a cell $Q$ of one of these lattices such that
		\begin{equation}\label{eq:ball_associated_cube}
			Q\subset c\, B\qquad\text{ and } \qquad c^{-1}r(B)\leq \ell(Q)\leq c\, r(B).
		\end{equation}
		We also say that such a cell $Q$ is \textit{associated with} $B$.
		Increasing $A$ if necessary, we assume that $A>200c$, and we fix $\alpha>100A^2c^3$ and $\beta>\max\{\alpha^{n+1},\alpha^{\log_2D}\}$, which are not related to the parameters with the same names in the previous subsections. We recall that a ball $B$ is $\mathfrak m$-$(\alpha,\beta)$-doubling if $\mathfrak m(\alpha B)\leq\beta\,\mathfrak m(B)$.
		
		Following the notation in \eqref{eq:def_centered_maximal}, for a locally finite complex measure $\mathfrak n$, we define
		\[
		\widetilde{\mathcal M}_{\Xi,\mathfrak m_2}\mathfrak n(x)\coloneqq\sup_{\substack{r>0, r\geq\Xi(x)}}
		\frac{|\mathfrak n|(B(x,r))}{\mathfrak m_2(B(x,10Ar))},\qquad \text{ for }x\in X,
		\]
		with the convention that the quantity on the right-hand side equals $0$ if both the numerator and the denominator are null. We also use the suppressed radial maximal function $M^r_\Xi$ defined in Section \ref{sec:aux_operator}, for which
		\[
		M^r_\Xi(h\mathfrak m_i)(x)=\sup_{\substack{r>0, r\geq\Xi(x)}} \frac{1}{r^n}\int_{B(x,r)}|h|\,d\mathfrak m_i.
		\]
		These two maximal operators satisfy the following bounds.
		
		\begin{lemma}\label{lemma:mct_maximal_functions}
			Let $i\in\{1,2\}$ and $1<s<\infty$. Then, for every $h\in L^s(\mathfrak m_i)$, we have
			\begin{equation}\label{eq:mct_maximal_mapping}
				\|\widetilde{\mathcal M}_{\Xi,\mathfrak m_2}(h\mathfrak m_i)\|_{L^s(\mathfrak m_3)}
				+\|M^r_\Xi(h\mathfrak m_i)\|_{L^s(\mathfrak m_3)}
				\lesssim_s\|h\|_{L^s(\mathfrak m_i)}.
			\end{equation}
		\end{lemma}
		
		\begin{proof}
			For $M^r_\Xi$, the estimate follows from the proof of Lemma \ref{lem:suppressed maximal}, with $\mathfrak m_3$ and $\mathfrak m_i$ in place of $\mathfrak m_1$ and $\mathfrak m_2$, and with $\Xi$ in place of $\mathcal R_1$ and $\mathcal R_2$, since \eqref{eq:mct_measure_bounds} provides the required growth. For $\widetilde{\mathcal M}_{\Xi,\mathfrak m_2}$, the same argument applies. Indeed, by \eqref{eq:mct_measure_bounds}, every ball $B=B(x,r)$ with $r\geq\Xi(x)$ satisfies \[\mathfrak m_3(5B)\leq C\,\mathfrak m_2(5A\, B)\leq C\,\mathfrak m_2(10A\, B),\] so the $5r$-covering lemma gives the boundedness from $L^1(\mathfrak m_i)$ to $L^{1,\infty}(\mathfrak m_3)$. Moreover, \eqref{eq:mct_measure_bounds} also gives $\widetilde{\mathcal M}_{\Xi,\mathfrak m_2}(h\mathfrak m_i)\leq C\|h\|_{L^\infty(\mathfrak m_i)}$, and \eqref{eq:mct_maximal_mapping} follows by interpolation.
		\end{proof}
		
		\vv
		The proof of \eqref{eq:mct_maximal_bound} relies on the following pointwise Cotlar-type inequality.
		\begin{lemma}\label{lemma:mct_pointwise_cotlar}
			Suppose that there exists $C_w>0$ such that
			\begin{equation}\label{eq:mct_normalized_weak}
				\mathfrak m_2\bigl(\{|Sh|>\lambda\}\bigr)	\leq\frac{C_w}{\lambda}\, \|h\|_{L^1(\mathfrak m_1)},
			\end{equation}
			for every $\lambda>0$ and every bounded function $h$ with bounded support. Then, for every bounded function $f$ with bounded support and every $x\in\supp\mathfrak m_3$, we have
			\begin{equation}\label{eq:mct_pointwise_cotlar}
				S_\sharp f(x)\lesssim \widetilde{\mathcal M}_{\Xi,\mathfrak m_2}(Sf\,\mathfrak m_2)(x)
				+\widetilde{\mathcal M}_{\Xi,\mathfrak m_2}(f\mathfrak m_1)(x)+M^r_\Xi(f\mathfrak m_1)(x),
			\end{equation}
			where the implicit constant depends only on $n$, $D$, $\gamma$, $C$, $A$, and $C_w$.
		\end{lemma}
		
		\begin{proof}
			We follow the standard proof of Cotlar's inequality for non-doubling measures, see \cite[Chapter 2]{To14}. Let $x\in\supp\mathfrak m_3$ and $\delta>0$. If $\delta<\Xi(x)$, then \eqref{eq:mct_kernel_bounds} gives \[|S_\delta f(x)-S_{\Xi(x)}f(x)|\lesssim M^r_\Xi(f\mathfrak m_1)(x),\] so we can assume without loss of generality that $\delta\geq\Xi(x)$.
			
			We set $s_0\coloneqq A\delta$ and $B_k\coloneqq B(x,\alpha^ks_0)$, for $k\geq0$. By \eqref{eq:mct_measure_bounds}, we have \[\mathfrak m_2(B_0)\geq C^{-1}\mathfrak m_3(B(x,\delta))>0\] and $\mathfrak m_2(B_k)\leq C\alpha^{kn}s_0^n$. Since $\beta>\alpha^n$, some $B_k$ is $\mathfrak m_2$-$(\alpha,\beta)$-doubling, and we denote by $m$ the smallest such $k$. We set $R\coloneqq\alpha^ms_0$ and $B\coloneqq B_m$, and we note that $\mathfrak m_2(B_k)\leq\beta^{k-m}\mathfrak m_2(B)$, for $0\leq k\leq m$.
			
			Since $Sf(y)=S(f\chi_{4B})(y)+S(f\chi_{X\setminus4B})(y)$, for every $y\in B$ we split
			\begin{equation}\label{eq:mct_cotlar_splitting}
				\begin{split}
					|S_\delta f(x)|
					&\leq|S_\delta f(x)-S_{4R}f(x)|+|S_{4R}f(x)-S(f\chi_{X\setminus 4B})(y)|+|S(f\chi_{4B})(y)|+|Sf(y)|\\
					&\eqqcolon\mathfrak I_1+\mathfrak I_2(y)+\mathfrak I_3(y)+\mathfrak I_4(y).
				\end{split}
			\end{equation}
			
			To bound $\mathfrak I_1$, we split the region $\{z:\delta<d(x,z)\leq4R\}$ into $\{z:\delta<d(x,z)\leq s_0\}$, the annuli $B_{k+1}\setminus B_k$, for $0\leq k\leq m-2$, and $\{z:\alpha^{m'}s_0\leq d(x,z)\leq4R\}$, where $m'\coloneqq\max\{m-1,0\}$. Since $\alpha^{k+1}s_0\geq\Xi(x)$ and $10A\alpha^{k+1}s_0\leq\alpha^{k+2}s_0$, the size estimate gives
			\begin{equation}\label{eq:mct_truncation_comparison}
				\begin{split}
					\mathfrak I_1&\lesssim\frac{1}{\delta^n}\int_{B(x,2s_0)}|f|\,d\mathfrak m_1+\sum_{k=0}^{m-2}\frac{1}{(\alpha^ks_0)^n}\int_{B_{k+1}}|f|\,d\mathfrak m_1+\frac{1}{(\alpha^{m'}s_0)^n}\int_{8B}|f|\,d\mathfrak m_1\\
					&\lesssim M^r_\Xi(f\mathfrak m_1)(x)+\widetilde{\mathcal M}_{\Xi,\mathfrak m_2}(f\mathfrak m_1)(x)\sum_{k=0}^{m-2}\frac{\mathfrak m_2(B_{k+2})}{(\alpha^ks_0)^n}
					\lesssim M^r_\Xi(f\mathfrak m_1)(x)+\widetilde{\mathcal M}_{\Xi,\mathfrak m_2}(f\mathfrak m_1)(x),
				\end{split}
			\end{equation}
			where, in the last inequality, we used that $\beta>\alpha^n$ and \[\mathfrak m_2(B_{k+2})\leq\beta^{k+2-m}\mathfrak m_2(B)\leq C\beta^{k+2-m}R^n.\]
			
			To bound $\mathfrak I_2(y)$, for $y\in B$, we use the regularity of $k$ in the first variable on the annuli $2^{\ell+2}R\leq d(x,z)<2^{\ell+3}R$, for $\ell\geq0$, and the size estimate on the sphere $d(x,z)=4R$, and we get
			\begin{equation}\label{eq:mct_far_oscillation}
				\mathfrak I_2(y)\lesssim\frac{1}{R^n}\, \int_{8B}|f|\,d\mathfrak m_1+\sum_{\ell\geq0}\frac{R^\gamma}{(2^{\ell+2}R)^{n+\gamma}}\, \int_{B(x,2^{\ell+3}R)}|f|\,d\mathfrak m_1\lesssim M^r_\Xi(f\mathfrak m_1)(x).
			\end{equation}
			
			To bound $\mathfrak I_3(y)$ and $\mathfrak I_4(y)$, we choose $y$ suitably. By \eqref{eq:mct_normalized_weak} and by Chebyshev's inequality, respectively, each of the inequalities
			\[
			|S(f\chi_{4B})(y)|>\frac{4C_w}{\mathfrak m_2(B)}\,\int_{4B}|f|\,d\mathfrak m_1
			\qquad \text{ and }\qquad
			|Sf(y)|>\frac4{\mathfrak m_2(B)}\,\int_B|Sf|\,d\mathfrak m_2
			\]
			holds on a subset of $B$ of $\mathfrak m_2$-measure at most $\mathfrak m_2(B)/4$. Since $\mathfrak m_2(B)\geq\mathfrak m_2(B_0)>0$, we can choose $y\in B$ for which both inequalities fail. Since $\alpha>40A$ and $B$ is $\mathfrak m_2$-$(\alpha,\beta)$-doubling, we have $\mathfrak m_2(40AB)\leq\beta\,\mathfrak m_2(B)$, and hence
			\[
			\mathfrak I_3(y)\leq 4C_w\beta\,\widetilde{\mathcal M}_{\Xi,\mathfrak m_2}(f\mathfrak m_1)(x)
			\qquad \text{ and }\qquad
			\mathfrak I_4(y)\leq 4\beta\,\widetilde{\mathcal M}_{\Xi,\mathfrak m_2}(Sf\,\mathfrak m_2)(x),
			\]
			which, together with \eqref{eq:mct_cotlar_splitting}, \eqref{eq:mct_truncation_comparison}, and \eqref{eq:mct_far_oscillation}, proves \eqref{eq:mct_pointwise_cotlar}.
		\end{proof}
		\vv
		
		Finally, we apply Lemmas \ref{lemma:mct_maximal_functions} and \ref{lemma:mct_pointwise_cotlar} to prove the opening result of this subsection.
		\begin{proof}[Proof of Lemma \ref{lemma:mixed_cotlar_transfer}]
			We first prove \eqref{eq:mct_weak_bound} for a bounded function $f$ with bounded support. Let $\lambda>0$. If $\mathfrak m_1(X)<\infty$ and $\lambda\leq\|f\|_{L^1(\mathfrak m_1)}/\mathfrak m_1(X)$, then \eqref{eq:mct_measure_bounds}, followed by monotone convergence as $r\to\infty$, gives $\mathfrak m_2(X)\leq C\mathfrak m_1(X)$, and hence
			\[
			\mathfrak m_2(\{|Sf|>2\lambda\})\leq\mathfrak m_2(X)\leq C\, \mathfrak m_1(X)
			\leq\frac{C}{\lambda}\, \|f\|_{L^1(\mathfrak m_1)}.
			\]
			Thus, we can assume that $\lambda>\|f\|_{L^1(\mathfrak m_1)}/\mathfrak m_1(X)$ whenever $\mathfrak m_1(X)<\infty$.
			
			We perform the Calder\'on-Zygmund decomposition of \cite[Lemma 2.14]{To14} as in \cite[Theorem A.2]{Con12}, on the adjacent dyadic lattices of \cite[Theorem 4.1]{HK12}, with the following choices. We say that a cell $Q$ is \textit{bad} if it is associated with a ball $B$, in the sense of \eqref{eq:ball_associated_cube}, such that $\int_Q|f|\,d\mathfrak m_1>\lambda\mathfrak m_1(A^2B)$. Since $f$ has bounded support, the choice of $\lambda$ implies that the side lengths of the bad cells are uniformly bounded, so every bad cell is contained in a maximal bad cell $Q_j$, that is, a bad cell such that $\ell(Q)\leq \ell(Q_j)$ for every bad cell $Q\supset Q_j$, of any of the lattices. Since two maximal bad cells of the same lattice are disjoint, the family $\{Q_j\}_j$ of maximal bad cells has bounded overlap. For every $j$, we fix a ball $B_j=B(z_j,r_j)$ such that $Q_j$ is associated with $B_j$ and
			\begin{equation}\label{eq:mct_bad_averages}
				\int_{Q_j}|f|\,d\mathfrak m_1>\lambda\, \mathfrak m_1(A^2B_j).
			\end{equation}
			By maximality, we also have
			\begin{equation}\label{eq:mct_larger_averages}
				\int_Q|f|\,d\mathfrak m_1\leq\lambda\, \mathfrak m_1(A^2B),
				\qquad \text{ if }Q_j\subset Q,\ \ell(Q_j)<\ell(Q),\ \text{ and }Q \text{ is associated with }B.
			\end{equation}
			Moreover, by \cite[Lemma 3.3]{Hy10} and the choice of $\beta$, for $\mathfrak m_1$-almost every $x\notin\bigcup_jQ_j$ there exist $\mathfrak m_1$-$(\alpha,\beta)$-doubling balls $B$ centered at $x$ with arbitrarily small radii. A cell $Q$ associated with $B$ contains $x$, so it is not contained in any $Q_j$, and hence it is not bad. Since $A^2<\alpha$, we get 
			\begin{equation*}
				\int_B|f|\,d\mathfrak m_1\leq\int_Q|f|\,d\mathfrak m_1\leq\lambda\,\mathfrak m_1(A^2B)\leq\beta\lambda\,\mathfrak m_1(B),
			\end{equation*}
			and, by the differentiation theorem along doubling balls \cite[Corollary 3.6]{Hy10}, we obtain $|f(y)|\leq\beta\lambda$, for $\mathfrak m_1$-almost every  $y\in X\setminus \bigcup_jQ_j$.
			
			We denote
			\[
			w_j\coloneqq\frac{\chi_{Q_j}}{\sum_k\chi_{Q_k}},
			\qquad m_j\coloneqq\int_{Q_j}|f|\,d\mathfrak m_1,
			\qquad \text{ and }\qquad t_j\coloneqq\max\{10A r_j,10\Xi(z_j)\},
			\]
			with $w_j\equiv 0$ outside $\bigcup_kQ_k$. For $k\geq0$, we set $D_{j,k}\coloneqq B(z_j,\alpha^kt_j)$, and we let $k_j$ be the smallest $k\geq0$ such that $D_{j,k}$ is $\mathfrak m_1$-$(\alpha,\beta)$-doubling, which exists as in the proof of Lemma \ref{lemma:mct_pointwise_cotlar}, since $t_j\geq\Xi(z_j)$ and $\mathfrak m_1(D_{j,0})\geq\mathfrak m_1(Q_j)>0$. We set $D_j\coloneqq D_{j,k_j}$, so that $r(D_j)=\alpha^{k_j}t_j$, we choose a cell $R_j$ associated with $D_j$, and we note that
			\begin{equation}\label{eq:mct_nondoubling_decay}
				\mathfrak m_1(D_{j,k})\leq \beta^{k-k_j}\mathfrak m_1(D_j),\qquad \text{ for }0\leq k\leq k_j.
			\end{equation}
			
			If $R_k\cap R_j\neq\varnothing$ and $r(D_k)\leq r(D_j)$, then $R_k\subset3cD_j$. Moreover, if $Q$ is a cell associated with $3cD_j$, then $Q_j\subset Q$ and $\ell(Q_j)\leq c r_j<3r(D_j)\leq\ell(Q)$, because $r(D_j)\geq 10Ar_j$. Thus, \eqref{eq:mct_larger_averages}, the choice of $\alpha$, and the choice of $k_j$ give
			\begin{equation}\label{eq:mct_compensator_mass}
				\int_{3cD_j}|f|\,d\mathfrak m_1\leq\lambda\, \mathfrak m_1(3A^2cD_j)
				\leq \beta\lambda\, \mathfrak m_1(D_j)\leq \beta\lambda\, \mathfrak m_1(R_j).
			\end{equation}
			Arguing as in the proof of \cite[Lemma 2.14]{To14}, these two properties provide functions $\varphi_j$ such that
			\begin{equation}\label{eq:mct_compensator_properties}
				\begin{split}
					&\supp\varphi_j\subset R_j,\qquad
					\int\varphi_j\,d\mathfrak m_1=\int w_jf\,d\mathfrak m_1,\\
					&\|\varphi_j\|_{L^\infty(\mathfrak m_1)}
					\leq\frac{2m_j}{\mathfrak m_1(R_j)},\qquad \text{ and }\qquad
					\sum_j|\varphi_j|\leq 4\beta\lambda.
				\end{split}
			\end{equation}
			Indeed, if the cells are finitely many and ordered so that $r(D_j)$ is non-decreasing in $j$, one sets inductively $\varphi_j\coloneqq\mathfrak m_1(G_j)^{-1}\bigl(\int w_jf\,d\mathfrak m_1\bigr)\chi_{G_j}$, where $G_j\coloneqq R_j\cap\{\sum_{k<j}|\varphi_k|\leq2\beta\lambda\}$. Since \begin{equation*}
				\int_{R_j}\sum_{k<j}|\varphi_k|\,d\mathfrak m_1\leq\int_{3cD_j}|f|\,d\mathfrak m_1\leq\beta\lambda\mathfrak m_1(R_j),
			\end{equation*}
			by \eqref{eq:mct_compensator_mass}, Chebyshev's inequality gives $\mathfrak m_1(G_j)\geq\mathfrak m_1(R_j)/2$. The case of countably many cells follows by a weak-* limiting argument.
			
			We now set
			\[
			g\coloneqq f\chi_{X\setminus\bigcup_jQ_j}+\sum_j\varphi_j
			\qquad \text{ and }\qquad b_j\coloneqq w_jf-\varphi_j,
			\]
			so that $f=g+\sum_jb_j$. The bounded overlap of the cells $Q_j$ and \eqref{eq:mct_compensator_properties} imply
			\begin{equation}\label{eq:mct_decomposition_bounds}
				\begin{split}
					&\|g\|_{L^1(\mathfrak m_1)}\lesssim\|f\|_{L^1(\mathfrak m_1)},\qquad
					\|g\|_{L^\infty(\mathfrak m_1)}\lesssim\lambda,\qquad
					\sum_jm_j\lesssim\|f\|_{L^1(\mathfrak m_1)},\\
					&\int b_j\,d\mathfrak m_1=0,\qquad
					\supp b_j\subset cD_j,\qquad \text{ and }\qquad
					\|b_j\|_{L^1(\mathfrak m_1)}\leq2m_j.
				\end{split}
			\end{equation}
			
			We set $\mathcal J\coloneqq\{j:A r_j\geq\Xi(z_j)\}$ and $E\coloneqq\bigcup_{j\in\mathcal J}A B_j$. Since $|Sh(x)|\leq C\delta_0^{-n}\|h\|_{L^1(\mathfrak m_1)}$ for every $x\in X$ and every bounded function $h$ with bounded support, the integral defining $S$ extends to $L^1(\mathfrak m_1)$, and it agrees $\mathfrak m_2$-almost everywhere with the bounded extension of $S$ on $L^1(\mathfrak m_1)\cap L^2(\mathfrak m_1)$. In particular, $Sf=Sg+\sum_jSb_j$, and Chebyshev's inequality gives
			\begin{equation}\label{eq:mct_weak_estimate_proof}
				\mathfrak m_2(\{|Sf|>2\lambda\})
				\leq\mathfrak m_2(E)+\frac{\|Sg\|_{L^2(\mathfrak m_2)}^2}{\lambda^2}
				+\frac1\lambda\sum_j\int_{X\setminus E}|Sb_j|\,d\mathfrak m_2.
			\end{equation}
			By \eqref{eq:mct_measure_bounds} and the definition of $\mathcal J$, \eqref{eq:mct_bad_averages}, and \eqref{eq:mct_decomposition_bounds}, we have
			\[
			\mathfrak m_2(E)\leq\sum_{j\in\mathcal J}\mathfrak m_2(A B_j)
			\leq C\, \sum_{j\in\mathcal J}\mathfrak m_1(A^2B_j)
			\leq\frac{C}{\lambda}\, \sum_jm_j
			\lesssim\frac{\|f\|_{L^1(\mathfrak m_1)}}{\lambda},
			\]
			and, moreover, \eqref{eq:mct_decomposition_bounds} gives $\|Sg\|_{L^2(\mathfrak m_2)}^2\leq\|S\|_{\op}^2\|g\|_{L^2(\mathfrak m_1)}^2\lesssim \|S\|_{\op}^2\lambda\, \|f\|_{L^1(\mathfrak m_1)}$.
			To bound the last term in \eqref{eq:mct_weak_estimate_proof}, we fix $j$ and we denote $k'_j\coloneqq\max\{k_j-1,0\}$. Since $|Sb_j|\leq|S(w_jf)|+|S\varphi_j|$, and since $D_{j,k'_j}$ is the union of $D_{j,0}$ and of the annuli $D_{j,k+1}\setminus D_{j,k}$, for $0\leq k\leq k_j-2$, we split
			\begin{equation}\label{eq:mct_bad_splitting}
				\begin{split}
					\int_{X\setminus E}|Sb_j|\,d\mathfrak m_2
					&\leq\int_{X\setminus4cD_j}|Sb_j|\,d\mathfrak m_2
					+\int_{4cD_j}|S\varphi_j|\,d\mathfrak m_2
					+\int_{D_{j,0}\setminus E}|S(w_jf)|\,d\mathfrak m_2\\
					&\qquad+\sum_{k=0}^{k_j-2}\int_{D_{j,k+1}\setminus D_{j,k}}|S(w_jf)|\,d\mathfrak m_2
					+\int_{4cD_j\setminus D_{j,k'_j}}|S(w_jf)|\,d\mathfrak m_2\\
					&\eqqcolon\mathfrak I_1+\mathfrak I_2+\mathfrak I_3+\mathfrak I_4+\mathfrak I_5,
				\end{split}
			\end{equation}
			where the sum is understood as zero when $k_j<2$.
			
			The terms $\mathfrak I_1$, $\mathfrak I_2$, $\mathfrak I_4$, and $\mathfrak I_5$ are estimated as in the proof of the weak $(1,1)$ estimate in \cite[Chapter 2]{To14}, since \eqref{eq:mct_measure_bounds} can be applied to balls centered at $z_j$ with radius at least $t_j\geq10\,\Xi(z_j)$.
			
			More precisely, cancellation, kernel regularity, and the growth of $\mathfrak m_2$ on the annuli $2^{\ell+2}c\,r(D_j)\leq d(x,z_j)<2^{\ell+3}c\,r(D_j)$ give \[\mathfrak I_1\lesssim\|b_j\|_{L^1(\mathfrak m_1)}\lesssim m_j.\]
			
			By \eqref{eq:mct_measure_bounds}, the choice of $\alpha$, and the choice of $k_j$, we have $\mathfrak m_2(4cD_j)\leq C\mathfrak m_1(4AcD_j)\leq C\beta\,\mathfrak m_1(R_j)$, so the Cauchy-Schwarz inequality and \eqref{eq:mct_compensator_properties} give
			\begin{equation*}
				\begin{split}	
					\mathfrak I_2&\leq\mathfrak m_2(4cD_j)^{1/2}\,\|S\|_{\op}\,\|\varphi_j\|_{L^2(\mathfrak m_1)}\\
					&\lesssim \|S\|_{\op}\,\mathfrak m_1(R_j)^{1/2}\bigl(\|\varphi_j\|_{L^\infty(\mathfrak m_1)}\|\varphi_j\|_{L^1(\mathfrak m_1)}\bigr)^{1/2}
					\lesssim \|S\|_{\op}\,m_j.
				\end{split}
			\end{equation*}
			Moreover, we have $d(x,y)\geq\alpha^kt_j/2$ for $x\notin D_{j,k}$ and $y\in Q_j$, and $\mathfrak m_2(D_{j,k+1})\leq C\mathfrak m_1(D_{j,k+2})$, because $A<\alpha$. Hence, the size estimate, \eqref{eq:mct_nondoubling_decay}, \eqref{eq:mct_measure_bounds}, and the fact that $\beta>\alpha^n$ yield
			\begin{equation*}
				\begin{split}
					\mathfrak I_4+\mathfrak I_5&\lesssim m_j\sum_{k=0}^{k_j-2}\frac{\mathfrak m_1(D_{j,k+2})}{(\alpha^kt_j)^n}+m_j\,\frac{\mathfrak m_2(4cD_j)}{(\alpha^{k'_j}t_j)^n}
					\\
					&\lesssim m_j\,\frac{\mathfrak m_1(D_j)}{r(D_j)^n}\,\alpha^{2n}\sum_{k=0}^{k_j-2}\Bigl(\frac{\alpha^n}{\beta}\Bigr)^{k_j-k-2}+\alpha^n\,m_j
					\lesssim m_j.
				\end{split}
			\end{equation*}
			
			The term $\mathfrak I_3$ is where the suppression plays a role. If $j\notin\mathcal J$, then $t_j=10\,\Xi(z_j)$ and $\Xi(y)\geq\Xi(z_j)-cr_j\geq\tfrac12\Xi(z_j)$ for $y\in Q_j$, so \eqref{eq:mct_kernel_bounds} and \eqref{eq:mct_measure_bounds} give
			\[
			\mathfrak I_3\leq\int_{D_{j,0}}|S(w_jf)|\,d\mathfrak m_2\lesssim\frac{m_j}{\Xi(z_j)^n}\,\mathfrak m_2\bigl(B(z_j,10\,\Xi(z_j))\bigr)
			\lesssim m_j.
			\]
			If $j\in\mathcal J$, then $t_j=10A r_j$ and $AB_j\subset E$. Since $d(x,y)\geq(A-c)r_j\geq A r_j/2$ for $x\in D_{j,0}\setminus A B_j$ and $y\in Q_j$, the size estimate and \eqref{eq:mct_measure_bounds} yield
			\[
			\mathfrak I_3\leq\int_{D_{j,0}\setminus A B_j}|S(w_jf)|\,d\mathfrak m_2
			\lesssim\frac{m_j}{(A r_j)^n}\,\mathfrak m_2\bigl(B(z_j,10A r_j)\bigr)
			\lesssim m_j.
			\]
			Gathering the bounds for $\mathfrak I_1$, $\mathfrak I_2$, $\mathfrak I_3$, $\mathfrak I_4$, and $\mathfrak I_5$ in \eqref{eq:mct_bad_splitting}, we obtain 
			\begin{equation*}
				\int_{X\setminus E}|Sb_j|\,d\mathfrak m_2\lesssim(1+\|S\|_{\op})m_j,
			\end{equation*}
			which, together with the bounds for $\mathfrak m_2(E)$ and $\|Sg\|_{L^2(\mathfrak m_2)}$ and with \eqref{eq:mct_decomposition_bounds}, plugged into \eqref{eq:mct_weak_estimate_proof}, yields
			\[
			\mathfrak m_2(\{|Sf|>2\lambda\})
			\lesssim\frac{(1+\|S\|_{\op})^2}{\lambda}\, \|f\|_{L^1(\mathfrak m_1)}.
			\]
			This proves \eqref{eq:mct_weak_bound} for bounded functions with bounded support. The corresponding statement for functions in $L^1(\mathfrak m_1)$ follows by a standard density argument.
			
			We now prove \eqref{eq:mct_maximal_bound}. We assume without loss of generality that $\|S\|_{\op}\leq 1$; indeed, replacing $S$ and $k$ by $(1+\|S\|_{\op})^{-1}S$ and $(1+\|S\|_{\op})^{-1}k$, respectively, the kernel estimates remain valid with the same constants. Then, \eqref{eq:mct_weak_bound} shows that \eqref{eq:mct_normalized_weak} holds with a constant $C_w$ depending only on $n$, $D$, $\gamma$, $C$, and $A$.
			
			Let $f$ be a bounded function with bounded support. By Lemma \ref{lemma:mct_pointwise_cotlar}, the inequality \eqref{eq:mct_pointwise_cotlar} holds for every $x\in\supp\mathfrak m_3$. By \eqref{eq:mct_maximal_mapping} with $i=2$ and $s=2$, the first term on its right-hand side satisfies
			\[
			\bigl\|\widetilde{\mathcal M}_{\Xi,\mathfrak m_2}(Sf\,\mathfrak m_2)\bigr\|_{L^2(\mathfrak m_3)}
			\lesssim\|Sf\|_{L^2(\mathfrak m_2)}
			\leq\|f\|_{L^2(\mathfrak m_1)}.
			\]
			Moreover, \eqref{eq:mct_maximal_mapping} with $i=1$ and $s=2$ bounds the other two terms in \eqref{eq:mct_pointwise_cotlar}, and we obtain $\|S_\sharp f\|_{L^2(\mathfrak m_3)}\lesssim\|f\|_{L^2(\mathfrak m_1)}$.
			This proves \eqref{eq:mct_maximal_bound} for bounded functions with bounded support.
			The estimate extends to functions in $L^2(\mathfrak m_1)$ by a standard density argument.
		\end{proof}
	}

	\vvv
	
	\section{The corona construction and the first reduction}\label{sec:corona}
	
	Under the normalized measure-testing assumptions of Proposition \ref{proposition:bounded_testing_normalized}, we denote by  $B_0\subset X$ a ball with $\mathfrak C_3$-small boundary which satisfies the  doubling-type property
	\begin{equation}\label{eq:db_sec_4_start}
		\mu(\Lambda\, B_0)\leq \beta\, \mu(\delta\,	B_0).
	\end{equation}
	In particular, we assume that $B_0$ is a testing ball in the sense of Definition \ref{def:accretive_testing_system}, and $B_0\in \mathscr B$. It is useful to split the balls in the family $\mathscr B$ according to their role in the dichotomy of the testing, and treat cases separately.  Namely, in the corona construction, we put $B\in\mathscr B'$ if at least one of the two systems uses alternative \eqref{eq:testing_opt_1} on $B$. Choose the set $V_B$ from
	one such system. Thus
	\begin{equation}\label{eq:alternative_1_chapter_4}
		\diam(V_B)\leq \eta_*\, r(B) \qquad \text{ and }\qquad \mu(V_B)\geq\frac{\zeta}{2}\, \mu(B).
	\end{equation}
	We denote $\mathscr B''\coloneqq \mathscr B\setminus\mathscr B'$. Both the system for $T$ and the system for $T^*$ satisfy \eqref{eq:bounded_normalized_testing} on every
	$B\in\mathscr B''$. Here $\eta_*$ is the concentration parameter chosen in \eqref{eq:choice_eta_Stop_1}.

	Let $\caD_0\coloneqq\caD^{(s)}_\mu$, where $s\coloneqq 2 r(B_0)$ and $\caD^{(s)}_\mu$ is the David-Mattila lattice of Lemma \ref{lemma:DM_cubes_one_sided}. We denote its doubling cells by $\caD^{db}_0$. For $R\in\caD_0$, we write
	\[
	\caD_0(R)\coloneqq\{Q\in\mathcal D_0:Q\subset R\}\qquad\text{ and }\qquad \caD^{db}_0(R)\coloneqq\mathcal D_0(R)\cap\mathcal D^{db}_0.
	\]
	We denote $\mu_0\coloneqq\mu|_{B_0}$. The cells in the corona construction below are contained in $B_0$, and hence satisfy $\mu_0(Q)=\mu(Q)$. Their associated balls and doubling properties are taken with respect to $\mu$, as in Appendix \ref{appendix:DM_testing_balls}. The construction in that appendix also applies to $\caD_0$.

	\subsection*{The construction of the exceptional set and the suppression}
	Let $R_1,R_2\in\caD^{db}_0$ be contained in $B_0$,
	and let $B_i\coloneqq\widetilde B(R_i)$, for $i=1,2$, be their testing balls, see Appendix \ref{appendix:DM_testing_balls}. We assume that $B_1,B_2\in\mathscr B''$, we denote by $\widehat\nu_1$ the testing measure for $T$ on $B_1$ and by $\widehat\nu_2$ the testing measure for $T^*$ on $B_2$, and we set
	\[
	\nu_i\coloneqq\widehat\nu_i+\mu|_{X\setminus B_i},\qquad i=1,2.
	\]
	
	The measures $\mu$ and $\widehat\nu_i$ satisfy \eqref{eq:improved_doubling_ex_set} and \ref{ex_item_1}-\ref{ex_item_6} with $B=B_i$ and $\mathfrak C_2=1$, where $T$ is replaced by $T^*$ if $i=2$. Indeed, \eqref{eq:improved_doubling_ex_set} follows from $B_i\in\mathscr B$ and Remark \ref{remark:doubling}, while \ref{ex_item_1}-\ref{ex_item_5} follow from the $n$-growth of $\mu$, \eqref{eq:measure_testing_size}, and \eqref{eq:measure_testing_small_sets}. To check \ref{ex_item_6}, we let $U_{B_i}$ be the set in \eqref{eq:bounded_normalized_testing} for the corresponding operator, and we set
	\[
	\widetilde U_{B_i}\coloneqq U_{B_i}\cup\bigl(B_i\setminus(1-\tau)B_i\bigr).
	\]
	By the small-boundary and doubling properties of $B_i$ and by \eqref{eq:measure_testing_tau_zeta}, we have
	\[
	\mu(\widetilde U_{B_i})\leq(\zeta+\mathfrak C_3\beta\tau)\,\mu(B_i)\leq\varepsilon_0\,\mu(B_i),
	\]
	and \eqref{eq:bounded_normalized_testing} gives $T_\sharp\widehat\nu_1\leq1$ on $B_1\setminus\widetilde U_{B_1}$. The same holds for $T^*$ and $\widehat\nu_2$, and this proves \ref{ex_item_6} with $U_B\coloneqq\widetilde U_{B_i}$ and $\mathfrak C_2=1$.
	
	\vvv
	
	Hence, we can apply the construction of Section \ref{section:exceptional_set} to $\widehat\nu_1$ and $\widehat\nu_2$, with independent parameters $w_1,w_2\in\Omega$. We denote by $\mathscr W_{R_1}(w_1)$ and $\mathscr W_{R_2}(w_2)$ the associated total exceptional sets in \eqref{eq:W_new}, and by $\mathscr T_{i,R_i}(w_i)$, $\mathscr S_{i,R_i}$, $\mathscr S^{\widehat{\mathbf m}}_{i,R_i}(w_i)$, $\mathscr H_{i,R_i}$, and $\mathscr E_{i,R_i}$ the sets in their construction, for $i=1,2$. We also write $B_{i,\theta}\coloneqq(1-\theta)B_i$ and $G_i\coloneqq G_{R_i}(w_i)\subset B_{i,\theta}$, and we recall that $\delta_1=1-9/(20\mathfrak C_1)\in(1/2,1)$, as in Lemma \ref{lemma:final_lemma_exceptional_set}, is the same for $i=1,2$.
	
	By Lemma \ref{lemma_RN_away_from_S}, the restriction of $|\widehat\nu_i|$ to $B_i\setminus\mathscr S_{i,R_i}$ is absolutely continuous with respect to $\mu$. Hence, by inner regularity, there exists a $\sigma$-compact set $K_i\subset(B_i\setminus\mathscr S_{i,R_i})\cap\supp\mu$ such that $|\widehat\nu_i|\bigl((B_i\setminus\mathscr S_{i,R_i})\setminus K_i\bigr)=0$, and we set
	\[
	F_i(w_i)\coloneqq G_i\cap K_i.
	\]
	Since $G_i\subset B_i\setminus\mathscr S_{i,R_i}$ and, by Lemma \ref{lemma:construction_testing_ball_DM}, $B_i\cap\supp\mu\subset R_i$, we have $F_i(w_i)\subset B_{i,\theta}\cap R_i$ and, by Lemma \ref{lemma:exceptional_set_new_contraction},
	\[
	|\widehat\nu_i|(F_i(w_i))=|\widehat\nu_i|(G_i)\geq\frac{1-\delta_1}{2}\,|\widehat\nu_i|(B_i),\qquad\text{ for }w_i\in\Omega.
	\]
	In particular, $F_i(w_i)\neq\varnothing$. 
	
	\begin{remark}
		It is worth noting that the map $(w,x)\mapsto\dist(x,F_i(w))$ is measurable. 
		Indeed, since $\Omega$ is standard Borel and the completion $\overline X$ of $X$ is a Polish space, for every $t>0$ the set $\{(w,x):\dist(x,F_i(w))<t\}$ is the projection of a Borel subset of $\Omega\times\overline X\times K_i$, and hence it is universally measurable. Thus, $(w,x)\mapsto\dist(x,F_i(w))$ is measurable with respect to the completed measures, which we use from now on.
	\end{remark}

	\vvv
	
	For $i=1,2$, we finally define the suppressions
	\begin{equation}\label{eq:definition_suppression_Phi_R_i}
		\Phi_{i,R_i,w_i}(x)\coloneqq\dist\bigl(x,F_i(w_i)\bigr),\qquad x\in X.
	\end{equation}
	Since $F_i(w_i)\neq\varnothing$, the function $\Phi_{i,R_i,w_i}$ is finite and $1$-Lipschitz, and its zero set is $\overline{F_i(w_i)}$. The parameters $w_1$ and $w_2$ are independent, and so are $\Phi_{1,R_1,w_1}$ and $\Phi_{2,R_2,w_2}$. By Lemma \ref{lemma:lemma_probabilistic_sets} below, there is a subset of $R_i$ of $\mu$-measure comparable to $\mu(R_i)$ whose points belong to the zero set of $\Phi_{i,R_i,w_i}$ with probability larger than a fixed $\delta_{\mathrm{pr}}>0$.
	
	\vv
	
	\begin{lemma}\label{lemma:lemma_probabilistic_sets}
		Let $R_i\in \caD^{db}_0$, $i=1,2$. Set
		\[
		\delta_{\mathrm{pr}}\coloneqq \frac{1-\delta_1}{8}\qquad\text{and}\qquad \tilde{\delta}_1\coloneqq \frac{1-\delta_1}{7+\delta_1}.
		\]
		Then there exist sets $\widetilde G_{i,R_i}\subset B_i\cap R_i$, for $i=1,2$, such that
		\[
		\mathbb P\bigl(\{w\in\Omega: \Phi_{i,R_i,w}(x)=0\}\bigr)>\delta_{\mathrm{pr}},
		\qquad \text{for all } x\in \widetilde G_{i,R_i}.
		\]
		Moreover, these sets satisfy
		\begin{equation}\label{eq:big piece of good set}
			\mu(\widetilde G_{i,R_i}) \geq \beta_*^{-1}\,(m_1 \beta)^{-1}\,\tilde{\delta}_1\,\mu(R_i).
		\end{equation}
		
	\end{lemma}
	
	\begin{proof}
		For $x\in B_i$, we define
		\[
		\mathfrak p_{i}(x)\coloneqq \mathbb P\bigl(\bigl\{w\in \Omega: x\in  F_i(w)\bigr\}\bigr).
		\]
		{
			We  have that 
			\begin{equation*}
				\begin{split}
					\int_{B_i} \mathfrak p_i(x)\, d|\widehat\nu_i|(x)&=\int_{w\in \Omega}\int \chi_{F_i(w)}(x)\, d|\widehat\nu_i|(x)\, d\mathbb P(w)\\
					&=\int_{w\in \Omega}|\widehat\nu_i|\bigl(F_i(w)\bigr)\, d\mathbb P(w)\overset{\eqref{eq:meas_good_set_nu}}{\geq}\frac{1-\delta_1}{4}\, |\widehat\nu_i|(B_i).
				\end{split} 
			\end{equation*}
			Thus, if we define
			\[
			\widetilde G_{i,R_i}\coloneqq \bigl\{x\in B_i: \mathfrak p_i(x)>(1-\delta_1)/8	\bigr\},
			\]
			in light of the latter estimate,  we obtain that
			\begin{equation*}
				\begin{split}
					|\widehat\nu_i|(B_i\setminus \widetilde G_{i,R_i})&\leq \frac{8}{7+\delta_1}\, \int_{B_i} (1-\mathfrak p_i(x))\, d|\widehat\nu_i|(x)\\
					&=\frac{8}{7+\delta_1}\, \biggl(|\widehat\nu_i|(B_i)-\int_{B_i}\mathfrak p_i(x)\, d|\widehat\nu_i|(x)\biggr) \leq \frac{6+2\delta_1}{7+\delta_1}\, |\widehat\nu_i|(B_i).
				\end{split}
			\end{equation*}
			This implies that
			\begin{equation*}
				|\widehat\nu_i|(\widetilde G_{i,R_i})\geq \frac{1-\delta_1}{7+\delta_1}\, |\widehat\nu_i|(B_i).
			\end{equation*}
			Since $\mathfrak p_i(x)>0$ implies $x\in F_i(w)$ for some $w$, and every $F_i(w)$ avoids the fixed ball-density exceptional set $\mathscr S_{i,R_i}$, \eqref{eq:bounds_radon_nyk_derivative} gives 
			\begin{equation*}
				\mu(\widetilde G_{i,R_i})\geq (\beta m_1)^{-1}\, \frac{1-\delta_1}{7+\delta_1}\, \mu(B_i).
			\end{equation*}
			
			Moreover, if  $x\not \in F_i(w)$ for all $w\in \Omega$, then  $\mathfrak p_i(x)=0$  and  so
			\begin{equation}\label{eq:inclusion_G_1_tilde}
				\widetilde G_{i,R_i}\subset  \bigcup_{w\in \Omega} F_i(w) \subset B_{i,\theta} \setminus  \bigcap_{w\in \Omega} \mathscr W_{R_i}(w).
			\end{equation}
			
			Observe that if $x \in 	\widetilde G_{i,R_i}$, then 
			\[
			\mathbb  P\bigl(\bigl\{w\in \Omega: \Phi_{i,R_i, w} (x) = 0\bigr\}\bigr) = \mathbb  P\bigl(\bigl\{w\in \Omega: 	\dist(x, F_i(w) ) = 0\bigr\}\bigr) \geq \mathfrak p_i(x)> 
			\delta_{\mathrm{pr}}.
			\]
			The property \eqref{eq:doubling_property_testing_DM_balls}, then, implies that $\mu(B_i)\geq \beta_*^{-1}\, \mu(R_i)$ and concludes the proof of the lemma.
		}
	\end{proof}
	
	\vv
	
	For $\delta_{\mathrm{pr}}$ as in Lemma \ref{lemma:lemma_probabilistic_sets}, we now define
	\begin{equation}\label{eq:definition_suppressions_2a}
		\Phi_{i,R_i}(x)\coloneqq \inf_{A\subset \Omega, \mathbb P(A)> \delta_{\mathrm{pr}}}\sup_{w\in A}\Phi_{i,R_i, w}(x).
	\end{equation}
	The infimum is over measurable sets $A$. Since every $F_i(w)$ is nonempty and contained in $B_i$, the functions in this formula are finite. Taking suprema and then infima preserves their common $1$-Lipschitz bound, so $\Phi_{i,R_i}$ is also $1$-Lipschitz.

	\vv
	\begin{lemma}\label{lem:big piece of zero set}
		There exists $\mathbf c_T\in(0,1)$ depending on $\beta, \beta_*, D, \mathfrak C_1$, and $\varepsilon_0$ such that
		\[
		\mu\bigl(\left\{ x\in R_i :  \Phi_{i,R_i}(x)=0   \right\}\bigr) \geq \mathbf c_T\, \mu(R_i), \qquad\text{ for } i=1,2.
		\]
	\end{lemma}
	
	\begin{proof}
		Let $x \in \widetilde G_{i,R_i}$. By Lemma \ref{lemma:lemma_probabilistic_sets}, the set
		\[
		A^i_x \coloneqq \left\{ w \in \Omega : \Phi_{i,R_i,w}(x)=0 \right\}
		\]
		satisfies $\mathbb{P}(A^i_x)>\delta_{\mathrm{pr}}$. Hence, for every $x \in \widetilde G_{i,R_i}$,
		\[
		\Phi_{i,R_i}(x) \le \sup_{w \in A^i_x} \Phi_{i,R_i,w}(x)=0.
		\]
		Thus, $ \widetilde G_{i,R_i} \subset \left\{ x\in R_i :  \Phi_{i,R_i}(x)=0   \right\}	 $, which, by \eqref{eq:big piece of good set},  completes the proof.
	\end{proof}
	
	\vvv
	
	\subsection{The corona decomposition}\label{subsection:corona_decomposition}
	We define the stopping cells separately for the alternatives \eqref{eq:testing_opt_1} and \eqref{eq:testing_opt_2} in the definition of testing balls associated with doubling David-Mattila cells.
	
	{
		\subsubsection*{The families $\Stopp'$ and $\Stopp''$}
		Let $i=1,2$, let $R_i\in\caD^{db}_0$ be contained in $B_0$, and let $B_i\coloneqq\widetilde B(R_i)$. We first assume that $B_i\in\mathscr B'$, and we let $V_{B_i}$ be as in \eqref{eq:alternative_1_chapter_4}. For $k\geq0$, we denote by $\ell_k\coloneqq56K_0\,s\,A_0^{-k}$ the side length of the cells of generation $k$, and we choose $k$ such that
		\[
		\frac{\eta_*\,r(B_i)}{A_0}<\ell_k\leq\eta_*\,r(B_i).
		\]
		Since $r(B_i)\leq r(R_i)\leq\ell(R_i)$, by \eqref{eq:structure_testing_balls} and \eqref{eq:comparab_r_l_0}, we have $k>J(R_i)$. Moreover, $B_i\cap\supp\mu\subset R_i$ so, by Lemma \ref{lemma:construction_testing_ball_DM}, every cell $P$ of generation $k$ that meets $V_{B_i}\cap\supp\mu$ is contained in $R_i$. Furthermore, if $v\in V_{B_i}\cap\supp\mu$ is fixed, by \eqref{eq:comparab_r_l_0} and \eqref{eq:new_radii_modified_DM}, every such cell $P$ satisfies
		\[
		d(z_P,v)\leq\diam(P)+\diam(V_{B_i})\leq2\eta_*\,r(B_i)\qquad\text{ and }\qquad r(P)\geq s\,A_0^{-k}>\frac{\eta_*\,r(B_i)}{56K_0A_0}.
		\]
		
		Note that the balls $5B(P)$ are pairwise disjoint, thus the geometric doubling property implies that there are at most $N$ such cells, where $N$ depends only on $D$, $K_0$, and $A_0$. These cells cover $V_{B_i}\cap\supp\mu$, so \eqref{eq:alternative_1_chapter_4} and \eqref{eq:doubling_property_testing_DM_balls} give a cell $\widetilde P_{B_i}\subset R_i$ of generation $k$ such that
		\begin{equation}\label{eq:properties_tilde_P_B}
			\mu(\widetilde P_{B_i})\geq\frac{\zeta}{2N}\,\mu(B_i)\geq\alpha_*\,\mu(R_i),\qquad\text{ where }\quad\alpha_*\coloneqq\frac{\zeta}{2N\beta_*}.
		\end{equation}
		Moreover, $\ell(\widetilde P_{B_i})=\ell_k\leq\eta_*\,r(B_i)\leq\eta_*\,\ell(R_i)$ and, by \eqref{eq:structure_testing_balls} and \eqref{eq:comparab_r_l}, $\ell(\widetilde P_{B_i})\approx\eta_*\,\ell(R_i)$. We denote by $P_{B_i}$ the smallest doubling cell containing $\widetilde P_{B_i}$, possibly $P_{B_i}=\widetilde P_{B_i}$, which exists because $R_i$ is doubling. Thus, $\widetilde P_{B_i}\subset P_{B_i}\subset R_i$. For a cell $Q\in\caD_0$, we write $\Theta_\mu(Q)\coloneqq\mu(Q)/\ell(Q)^n$, and we have the following estimate.
		\begin{lemma}\label{lemma:density_control_corona_1}
			For $R_i\in \caD^{db}_0$ such that $B_i\in \mathscr B'$, we have that 
			\[
			\Theta_\mu(P_{B_i})\gtrsim \Theta_\mu (\widetilde P_{B_i})\geq \alpha_* \eta_*^{-n}\, \Theta_\mu(R_i).
			\]
		\end{lemma}
		\begin{proof}
			If $P_{B_i}\neq\widetilde P_{B_i}$, every cell $S$ with $\widetilde P_{B_i}\subset S\subsetneq P_{B_i}$ is non-doubling, so Lemma \ref{lemma:lemma_density_nonboubling_chain} gives $\Theta_\mu(100B(\widetilde P_{B_i}))\lesssim\Theta_\mu(100B(P_{B_i}))$, which is trivial if $P_{B_i}=\widetilde P_{B_i}$. Since $\widetilde P_{B_i}\subset100B(\widetilde P_{B_i})$, $r(\widetilde P_{B_i})\leq\ell(\widetilde P_{B_i})$, $P_{B_i}$ is doubling, and $\ell(P_{B_i})\lesssim r(P_{B_i})$, we obtain
			\[
			\Theta_\mu(\widetilde P_{B_i})\leq100^n\,\Theta_\mu(100B(\widetilde P_{B_i}))\lesssim\Theta_\mu(100B(P_{B_i}))\lesssim\Theta_\mu(P_{B_i}).
			\]
			Moreover, since $\ell(\widetilde P_{B_i})\leq\eta_*\,\ell(R_i)$, \eqref{eq:properties_tilde_P_B} gives
			\[
			\Theta_\mu(\widetilde P_{B_i})=\frac{\mu(\widetilde P_{B_i})}{\ell(\widetilde P_{B_i})^n}\geq 
			\alpha_*\,\eta_*^{-n}\,\Theta_\mu(R_i).\qedhere
			\]
		\end{proof}

		Since $\mu(P_{B_i})\geq\mu(\widetilde P_{B_i})\geq\alpha_*\,\mu(R_i)$, Lemma \ref{lemma:density_control_corona_1} gives \[\Theta_\mu(P_{B_i})\,\mu(P_{B_i})\geq c_0\,\alpha_*^2\,\eta_*^{-n}\,\Theta_\mu(R_i)\,\mu(R_i),\]
		where $c_0>0$ depends only on $n$ and the parameters of the David-Mattila lattice. We fix $\eta_*\in(0,1/10)$ such that $c_0\,\alpha_*^2\,\eta_*^{-n}\geq2$, which depends only on $n$, $\zeta$, and the parameters of the David-Mattila and testing-ball constructions, so that
		\begin{equation}\label{eq:choice_eta_Stop_1}
			\Theta_\mu(P_{B_i})\, \mu(P_{B_i}) \geq 2\, \Theta_\mu(R_i) \,\mu(R_i).
		\end{equation}
		In particular, $P_{B_i}\subsetneq R_i$. We define $\Stopp'_i(R_i)$ as the family of cells $S\subset R_i$ that belong to the same generation as $P_{B_i}$. These cells are pairwise disjoint proper subcells of $R_i$, cover $R_i$, and include $P_{B_i}$.
	}
	
	\vv
	
	{
		Next, let $i=1,2$, and suppose that $B_i=\widetilde B(R_i)\in \mathscr B''$. By \eqref{eq:definition_suppressions_2a}, $\Phi_{i,R_i}$ is $1$-Lipschitz. For $x\in R_i$ we define $Q_x\coloneqq\{x\}$ if $\Phi_{i,R_i}(x)=0$. Otherwise, we denote by $Q_x$ the cell of largest side length in $\mathcal D_0(R_i)$ containing $x$ such that
		\begin{equation}\label{eq:reg_cells_sec4}
			\ell(Q_x)<\frac{1}{480\, A_0} \inf_{y\in Q_x} \Phi_{i, R_i}(y),
		\end{equation}
		which we refer to as a \textit{regularized cell}. By Lemma \ref{lem:big piece of zero set}, there exists a point $z_i\in R_i$ with $\Phi_{i,R_i}(z_i)=0$, so $R_i$ does not satisfy \eqref{eq:reg_cells_sec4}. On the other hand, sufficiently small descendants $P\ni x$ satisfy
		\[
		\inf_P\Phi_{i,R_i}
		\geq\Phi_{i,R_i}(x)-\operatorname{diam}(P)
		\geq\Phi_{i,R_i}(x)-\ell(P)
		>480A_0\ell(P)
		\]
		whenever $\ell(P)<\Phi_{i,R_i}(x)/(480A_0+1)$. Thus $Q_x$ exists, is a proper subcell of $R_i$, and its parent lies in $\mathcal D_0(R_i)$ and fails \eqref{eq:reg_cells_sec4}.
	}
	
	{	
		We define $\Stopp''_i(R_i)$ as the collection of distinct regularized cells $Q_x$ with $x\in R_i$ and $\Phi_{i,R_i}(x)>0$. By maximality and nesting, these cells are pairwise disjoint, and each is contained in $R_i\cap\{\Phi_{i,R_i}>0\}$.
	}
	
	{
		\begin{lemma}\label{lemma:regularized_cubes_properties}
			The regularized cells satisfy the following properties:
			\begin{enumerate}
				\item If $Q\in \Stopp''_i(R_i)$ and $x\in B(z_Q, 50\,  \ell(Q))$, then $\Phi_{i,R_i}(x)\approx \ell(Q)$. More specifically,
				\begin{equation}\label{eq:comparability_reg_supp_prop_1}
					430\,A_0\, \ell(Q)< \Phi_{i,R_i}(x)\leq 481\, A_0^2\, \ell(Q).
				\end{equation}
				\item If $Q,Q'\in \Stopp''_i(R_i)$ are such that $B(z_Q, 50\, \ell(Q))\cap B(z_{Q'}, 50 \ell(Q'))\neq \varnothing$, then $\ell(Q)\approx \ell(Q')$. More specifically,
				\[
				\frac{1}{2 A_0}\, \ell(Q')\leq \ell(Q)\leq 2\, A_0\, \ell(Q').
				\]
				\item  It holds that 
				\begin{equation}\label{eq:packing_0_stop}
					\sum_{P\in \Stopp''_i(R_i)}\mu(P)\leq (1-  \mathbf c_T)\, \mu(R_i).
				\end{equation}
			\end{enumerate}
		\end{lemma}
		\begin{proof}
			Let $Q\in \Stopp''_i(R_i)$ and $x\in B(z_Q, 50\, \ell(Q))$. The definition \eqref{eq:reg_cells_sec4} yields $\Phi_{i,R_i}(z_Q)> 480\,A_0\,  \ell(Q)$, and by $1$-Lipschitz continuity we obtain
			\[
			\Phi_{i,R_i}(x)\geq\Phi_{i,R_i}(z_Q)-d(x,z_Q)>(480A_0-50)\ell(Q)\geq430A_0\ell(Q).
			\]
			By definition of regularized cell, the David-Mattila parent $\widehat Q$ of $Q$ does not satisfy  \eqref{eq:reg_cells_sec4}. Hence, we have
			\begin{equation}\label{eq:bd_below_parent_phi_reg}
				\inf_{y\in \widehat Q}\Phi_{i,R_i}(y)\leq 480\, A_0\, \ell(\widehat Q)= 480\, A_0^2\, \ell(Q).
			\end{equation}
			For $k\in\mathbb N$, let $y_k\in\widehat Q$ satisfy $\Phi_{i,R_i}(y_k)<\inf_{\widehat Q}\Phi_{i,R_i}+1/k$. Since $z_Q\in Q\subset\widehat Q$,
			\[
			d(x,y_k)\leq d(x,z_Q)+\operatorname{diam}(\widehat Q)
			<(50+A_0)\ell(Q)<2A_0\ell(Q).
			\]
			Thus Lipschitz continuity and \eqref{eq:bd_below_parent_phi_reg} yield
			\begin{equation}\label{eq:bd_below_parent_phi_reg2}
				\begin{split}
					\Phi_{i,R_i}(x)&\leq\Phi_{i,R_i}(y_k)+d(x,y_k)\\
					&<(480A_0^2+2A_0)\ell(Q)+1/k<481A_0^2\ell(Q)+1/k.
				\end{split}
			\end{equation}
			Thus, the upper bound in \eqref{eq:comparability_reg_supp_prop_1} follows from \eqref{eq:bd_below_parent_phi_reg} and taking the limit as $k\to \infty$ in \eqref{eq:bd_below_parent_phi_reg2}.
			
			To prove (2), choose $x\in B(z_Q,50\ell(Q))\cap B(z_{Q'},50\ell(Q'))$. By \eqref{eq:comparability_reg_supp_prop_1},
			\[
			430A_0\ell(Q)<\Phi_{i,R_i}(x)\leq481A_0^2\ell(Q')
			\]
			and
			\[
			\ell(Q)<\frac{481}{430}A_0\ell(Q')<2A_0\ell(Q').
			\]
			Interchanging $Q$ and $Q'$ gives the other inequality. To prove the third assertion it suffices to observe that, since the cells in $\Stopp''_i(R_i)$ are disjoint, we have
			\begin{equation}
				\sum_{P\in \Stopp''_{i}(R_i)}\mu(P)= \mu \Bigl(\bigcup_{S \in \Stopp''_{i}(R_i)} S \Bigr)\leq \mu\bigl(\{x\in R_i: \Phi_{i,R_i}(x)>0\}\bigr),
			\end{equation}
			and conclude by Lemma \ref{lem:big piece of zero set}.
		\end{proof}
	}
	
	\vv
	
	\subsubsection*{The coronization}
	We now iterate the construction of $\Stopp'_i(R_i)$ and $\Stopp''_i(R_i)$.
	We denote by $\Top^0$ the collection of maximal doubling cells $Q\in\caD_0$ contained in $B_0$.	Thus, $\Top^0$ is a pairwise disjoint family of doubling cells which covers $B_0$ modulo a $\mu$-null set.
	
	Next, for $Q\in \caD_0$ we denote by $\Mdb(Q)$ a family of maximal disjoint doubling cells in $\caD^{db}_0$ that are contained in $Q$ and cover $Q$ modulo a $\mu$-null set. We allow equality, so $\Mdb(Q)=\{Q\}$ when $Q$ is doubling.
	
	{
		{
			For $R_i\in \caD^{db}_0$, we write
			\[
			\Next_i(R_i)\coloneqq 
			\begin{cases}
				\bigcup_{Q\in \Stopp'_i(R_i)}\Mdb(Q), \qquad \text{ if }\widetilde B(R_i)\in \mathscr B',\\
				\bigcup_{Q\in \Stopp''_i(R_i)}\Mdb(Q), \qquad \text{ if }\widetilde B(R_i)\in \mathscr B''.
			\end{cases}
			\]
			
			Denoting $\Top^0_1\coloneqq \Top^0_2\coloneqq \Top^0$, for $k\geq 0$ we inductively define
			\[
			\Top^{k+1}_i\coloneqq \bigcup_{P \in \Top^k_i}\Next_i(P).
			\]
			
			Thus, we finally denote $\Top_i\coloneqq \bigcup_{k=0}^\infty \Top^k_i$ and we further split top cells as the union of 
			\[
			\Top'_i\coloneqq \bigl\{Q\in \Top_i: \widetilde B(Q)\in \mathscr B'\bigr\}
			\]
			and
			\[
			\Top''_i\coloneqq \bigl\{Q\in \Top_i: \widetilde B(Q)\in \mathscr B''\bigr\}= \Top_i\setminus \Top'_i.
			\]

			In the following lemma, we prove that $\Top''_i$ is a Carleson family.
			\begin{lemma}\label{lemma:packing_top_2}
				For every $S\in \caD_0$, we have
				\begin{equation}\label{eq:packing_top_2}
					\sum_{Q\in \Top''_i, Q\subset S}\mu(Q)\leq \mathbf c^{-1}_T\, \mu(S).
				\end{equation}
			\end{lemma}
			\begin{proof}
				Let $R\in \Top''_i$. Since $\Stopp''_i(R)$ is a family of pairwise disjoint cells and every David-Mattila cell is covered, up to a $\mu$-null set, by its maximal doubling subcells, by \eqref{eq:packing_0_stop} we have
				\begin{equation}\label{eq:sum_next_1}
					\begin{split}
						\sum_{Q\in \Next_i(R)}\mu(Q)&=\sum_{P\in \Stopp''_i(R)}\mu(P)\leq (1-\mathbf c_T)\, \mu(R).
					\end{split}
				\end{equation}
				Next, for $R\in \Top''_i$ we define
				\[
				\operatorname{Succ}(R)\coloneqq \bigl\{Q\in \Top''_i: Q\subsetneq R\text{ and there is no } P\in \Top''_i \text{ s.t. }Q\subsetneq P\subsetneq R\bigr\},
				\]
				which is a family of pairwise disjoint cells. Moreover, by construction, every cell in $\operatorname{Succ}(R)$ is contained in an element of $\Next_i(R)$, so
				\begin{equation}\label{eq:successor_measure}
					\sum_{Q\in \operatorname{Succ}(R)}\mu(Q)\leq \sum_{Q\in \Next_i(R)}\mu(Q)\overset{\eqref{eq:sum_next_1}}{\leq} (1-\mathbf c_T)\, \mu(R).
				\end{equation}
				
				To prove the lemma, we apply \eqref{eq:successor_measure} iteratively. Namely, for $S\in \caD_0$, we denote by $\mathcal F_0(S)$ the family of cells $Q\in \Top''_i$ that are contained in $S$ and are maximal with respect to set inclusion, and we recursively define
				\[
				\mathcal F_k(S)\coloneqq \bigcup_{Q\in \mathcal F_{k-1}(S)} \operatorname{Succ}(Q)\, \qquad \text{ for }k\geq 1.
				\]
				Hence, cells in $\mathcal F_k(S)$ are pairwise disjoint and
				\[
				\{Q\in \Top''_i: Q\subset S\}= \bigcup_{k\geq 0} \mathcal F_k(S),
				\]
				so by \eqref{eq:successor_measure} we obtain
				\begin{equation*}
					\begin{split}
						\sum_{Q\in \Top''_i, Q\subset S}\mu(Q)&= \sum_{k\geq 0}\sum_{Q\in \mathcal F_k(S)}\mu(Q)\leq \sum_{k\geq 0}(1-\mathbf c_T)^k\, \mu(S)= \mathbf c_T^{-1}\, \mu(S),
					\end{split}
				\end{equation*}
				where in the last equality we used that $\mathbf c_T\in (0,1)$.
			\end{proof}
			
			Then, we apply \eqref{eq:choice_eta_Stop_1} and Lemma \ref{lemma:packing_top_2} to obtain a packing condition for cells in $\Top_i$.
			
			\begin{lemma}\label{lemma:packing_top_density}
				With $\eta_*$ as in \eqref{eq:choice_eta_Stop_1}, there exists $\mathbf c'>0$ depending on $n$, $\mathfrak C_0$, and $\mathbf c_T$ such that, for every $S\in \caD_0$, we have
				\begin{equation}\label{eq:packing_top_density}
					\sum_{Q\in \Top_i, \, Q\subset S}\Theta_\mu(Q)\, \mu(Q)\leq \mathbf c'\, \mu(S).
				\end{equation}
			\end{lemma}
			\begin{proof}
				It suffices to prove \eqref{eq:packing_top_density} for $S\in\Top_i$ and, for a general cell $S$, to sum the resulting bounds over the pairwise disjoint maximal cells of $\Top_i$ contained in $S$. For $k\geq 0$ we define
				\begin{equation*}
					\Top_{i,(k)}\coloneqq \bigcup^k_{j=0} \Top^{j}_i \qquad \text{ and }\qquad	\Top'_{i,(k)}\coloneqq \bigcup^k_{j=0}\bigl(\Top'_i\cap \Top^{j}_i\bigr)
				\end{equation*}
				and $\Top_{i,(k)}(S)\coloneqq\{P\in\Top_{i,(k)}:P\subset S\}$, with $\Top'_{i,(k)}(S)\coloneqq\Top_{i,(k)}(S)\cap\Top'_i$.
				The families $\Top''_{i,(k)}$ and $\Top''_{i,(k)}(S)$ are defined analogously.
				
				Let $\eta_*$ be chosen as in \eqref{eq:choice_eta_Stop_1}. Every cell $P\in \Top_i$ is doubling, so the $n$-growth of $\mu$ yields
				\begin{equation}\label{eq:bound_density_doubling_cell}
					\begin{split}
						\Theta_\mu(P)&\leq\frac{\mu(28B(P))}{\ell(P)^n}\leq\mathfrak C_0\, \Bigl(\frac{28r(P)}{\ell(P)}\Bigr)^n\leq\frac{\mathfrak C_0}{2^n}.
					\end{split}
				\end{equation}
				
				Moreover, if $P\in \Top'_i$, then $P_{\widetilde B(P)}\in\Stopp'_i(P)$ is doubling, so that $P_{\widetilde B(P)}\in \Next_i(P)$, and we have
				\begin{equation}\label{eq:choice_eta_Stop_Cor}
					\sum_{Q\in \Next_i(P)}\Theta_\mu(Q)\, \mu(Q) \geq\Theta_\mu\bigl(P_{\widetilde B(P)}\bigr)\,\mu\bigl(P_{\widetilde B(P)}\bigr)\overset{\eqref{eq:choice_eta_Stop_1}}{\geq} 2\, \Theta_\mu(P)\, \mu(P).
				\end{equation}
				Hence, noting that $\Top''_{i,(k)}(S)\subseteq\Top''_i(S)$,
				\begin{equation}\label{eq:packing_top_densitya}
					\begin{split}
						&\sum_{P\in \Top_{i,(k)}(S)}\Theta_\mu(P)\mu(P)=\sum_{P\in \Top'_{i,(k)}(S)}\Theta_\mu(P)\mu(P) + \sum_{P\in \Top''_{i,(k)}(S)}\Theta_\mu(P)\mu(P)\\
						&\qquad\overset{\eqref{eq:choice_eta_Stop_Cor}}{\leq} \frac{1}{2}\, \sum_{P\in \Top'_{i,(k)}(S)}\sum_{Q\in \Next_i(P)}\Theta_\mu(Q)\, \mu(Q) + \sum_{P\in \Top''_{i,(k)}(S)}\Theta_\mu(P)\mu(P)\\
						&\qquad \overset{\eqref{eq:bound_density_doubling_cell}}{\leq}\frac{1}{2}\sum_{P\in  \Top_{i,(k)}(S)}\Theta_\mu(P)\, \mu(P) + \frac{\mathfrak C_0}{2^n}\, \mu(S) + \sum_{P\in \Top''_{i}(S)}\Theta_\mu(P)\mu(P),
					\end{split}
				\end{equation}
				where the new successors at level $(k+1)$ are pairwise disjoint subsets of $S$, so their density sum is at most $2^{-n}\mathfrak C_0\mu(S)$.
				
				By Lemma \ref{lemma:packing_top_2}, $\Top''_i$ is a Carleson family which, together with the density bound \eqref{eq:bound_density_doubling_cell} implies
				\begin{equation}\label{eq:packing_top_densityb}
					\sum_{P\in \Top''_{i}(S)}\Theta_\mu(P)\mu(P)\leq \frac{\mathfrak C_0}{2^n}\, \sum_{P\in \Top''_{i}(S)}\mu(P) \leq \frac{\mathfrak C_0}{2^n}\, \mathbf c_T^{-1}\, \mu(S).
				\end{equation}	
				Thus, we gather \eqref{eq:packing_top_densitya}, \eqref{eq:packing_top_densityb}, and obtain
				\[
				\sum_{P\in \Top_{i,(k)}(S)}\Theta_\mu(P)\mu(P)\leq \frac{2\, \mathfrak C_0}{2^n}\, (1+ \mathbf c_T^{-1})\,  \mu(S)\eqqcolon \mathbf c'\, \mu(S).
				\]
				Taking the limit as $k\to \infty$ in the latter bound concludes the proof of the lemma.
			\end{proof}
		}
	}

	To simplify notation, for $R_i\in \Top_i$ we write
	\[
	\Stopp_i(R_i)\coloneqq
	\begin{cases}
		\Stopp'_i(R_i), \qquad \text{ if }R_i \in \Top'_i,\\
		\Stopp''_i(R_i), \qquad \text{ if }R_i \in \Top''_i.
	\end{cases}
	\]
	For $R_i\in \Top_i$, we define
	\begin{equation}\label{eq:tree_case_2}
		\Tree_{i}(R_i)\coloneqq \bigl\{Q\in \caD_0(R_i): \nexists\,  S\in \Stopp_{i}(R_i) \text{ such that }Q\subseteq S \bigr\}.
	\end{equation}
	
	For $R_i\in \Top_i$, the tree defined in \eqref{eq:tree_case_2} satisfies the following properties:
	
	\begin{enumerate}
		\item The cell $R_i$ is the maximal element of $\Tree_{i}(R_i)$ with respect to inclusion. Namely, $Q\subset R_i$ for every $Q\in \Tree_{i}(R_i)$.\label{property_1_tree}
		\item If $Q,Q'\in \Tree_{i}(R_i)$, $Q\subset Q'$, and $P\in \mathcal D_0(R_i)$ is such that $Q\subset P\subset Q'$, then $P\in \Tree_{i}(R_i)$.\label{property_2_tree}
		\item If $Q\in\Tree_i(R_i)$, each child of $Q$ either belongs to $\Tree_i(R_i)$ or is a stopping cell in $\Stopp_i(R_i)$.\label{property_3_tree}
	\end{enumerate}
	
	Finally, for  $R_i\in \Top_i$ and  $Q\in \Stopp_{i}(R_i)$, we define
	\[
	\mathcal T_{i,Q}\coloneqq \bigl\{P\in \caD_0(Q): \nexists \,S\in \Mdb(Q) \text{ such that }P\subseteq S \bigr\}.
	\]
	If $Q$ is doubling, then $\mathcal T_{i,Q}=\varnothing$. Otherwise, every cell in $\mathcal T_{i,Q}$ is non-doubling. By \eqref{eq:tree_case_2} and the definition of $\mathcal T_{i,Q}$, the families $\Tree_i(R_i)$ and $\mathcal T_{i,Q}$ are pairwise disjoint as $R_i\in\Top_i$ and $Q\in\Stopp_i(R_i)$ vary.

	\vvv
	\subsection{The mixed corona construction}\label{section:mixed_corona}
	Since the cells in $\Top^0$ are pairwise disjoint and cover $B_0$ modulo a $\mu$-null set, Lemmas \ref{lemma:packing_top_2} and \ref{lemma:packing_top_density} give
	\[
	\sum_{Q\in\Top_i''}\mu(Q)\leq\mathbf c_T^{-1}\, \mu(B_0)
	\]
	and 
	\[
	\sum_{Q\in\Top_i}\Theta_\mu(Q)\mu(Q)\leq\mathbf c'\, \mu(B_0).
	\]

	For $R_1\in \Top''_1$ and $R_2\in \Top''_2$ we define
	\[
	\widetilde{\Tree}(R_1,R_2)\coloneqq \Tree_1(R_1)\cap\Tree_2(R_2),
	\]
	which we also refer to as a \textit{mixed tree}.
	If $\Tree_1(R_1)\cap \Tree_2(R_2)\neq \varnothing$, then  $\widetilde{\Tree}(R_1,R_2)$ has a maximal element $R$ with respect to inclusion, which coincides with the smaller of the cells $R_1$ and $R_2$. It also satisfies (\ref{property_2_tree}). Moreover, if $Q\in\widetilde{\Tree}(R_1,R_2)$ and a child $P$ of $Q$ does not belong to $\widetilde{\Tree}(R_1,R_2)$, then $P\in\Stopp_1(R_1)\cup\Stopp_2(R_2)$.
	Hence in this case, with a slight abuse of notation, we write
	\[
	\widetilde{\Tree}(R)\coloneqq \widetilde{\Tree}(R_1,R_2).
	\]

	Finally, we define the initial intermediate family
	\begin{equation}\label{eq:initial_intermediate_tree}
		\mathcal T_0\coloneqq \{Q\in\mathcal D_0:Q\not\subseteq R\text{ for every } R\in\Top^0\}.
	\end{equation}		
	We denote by $\widetilde{\Top}$ the roots of the nonempty mixed trees. For each $R\in\widetilde{\Top}$, there is a unique pair $(R_1,R_2)\in\Top_1''\times\Top_2''$ such that $R\in\Tree_1(R_1)\cap\Tree_2(R_2)$. Moreover,
	\[
	\widetilde{\Top}\subset\Top''_1\cup\Top''_2.
	\]
	Hence Lemma \ref{lemma:packing_top_2} implies the following estimate.
	\begin{lemma}\label{lemma:packing2_mixed}
		The collection $\widetilde{\Top}$ is a Carleson family. More precisely, for every $S\in \caD_0$, we have
		\[
		\sum_{Q\in \widetilde{\Top}, Q\subset S}\mu(Q)\leq 2 \,\mathbf c_T^{-1}\, \mu(S),
		\]
		where $\mathbf c_T$ is as in Lemma \ref{lem:big piece of zero set}.
	\end{lemma}
	\begin{proof}
		By the definition of  $\widetilde{\Top}$ and Lemma \ref{lemma:packing_top_2},  for $S\in \caD_0$, we get that
		\[
		\sum_{Q\in \widetilde{\Top}, Q\subset S}\mu(Q)\leq \sum_{Q\in \Top''_1, Q\subset S}\mu(Q) + \sum_{Q\in \Top''_2, Q\subset S}\mu(Q)\leq 2\, \mathbf c_T^{-1}\, \mu(S).\qedhere
		\]
	\end{proof}

	Finally, we gather all the families of cells defined in the mixed corona construction, and note that
	\begin{equation}\label{eq:decomposition_D_0}
		\begin{split}
			\caD_0 =\mathcal T_0 &\cup \Bigl[\bigcup_{R\in \widetilde{\Top}}\widetilde{\Tree}(R)\Bigr]\cup\Bigl[\bigcup_{i=1,2}\bigcup_{R_i\in \Top'_i}\Bigl(\Tree_i(R_i)\cup\bigcup_{Q\in \Stopp'_i(R_i)}\mathcal T_{i,Q}\Bigr)\Bigr]\\
			&\cup\Bigl[\bigcup_{i=1,2}\bigcup_{R_i\in \Top''_i}\bigcup_{Q\in \Stopp''_i(R_i)}\mathcal T_{i,Q}\Bigr].
		\end{split}
	\end{equation}
	\vv

	\subsection{SIOs associated with the mixed corona decomposition}
	The families $\widetilde{\Tree}(R)$, $R\in\widetilde{\Top}$, are pairwise disjoint and disjoint from the remaining families in \eqref{eq:decomposition_D_0}. Moreover, a cell may belong to the remaining families for both $i=1$ and $i=2$ so, in \eqref{eq:splitting_T_in_sigmas_tree} below, its absolute value is then counted twice.
	
	Now, we define the SIO associated with each family on the right-hand side of \eqref{eq:decomposition_D_0}. To do so, for $Q\in \caD_0$, we set
	\begin{equation}
		\rho(Q)\coloneqq sA_0^{-J(Q)}=\frac{\ell(Q)}{56K_0},
	\end{equation}
	and note that 
	\begin{equation}\label{eq:comparability_r_rho}
		\rho(Q)\leq r(Q)\leq K_0\, \rho(Q).
	\end{equation}
	If $Q\in\caD_0\setminus\caD^{db}_0$, then $r(Q)=\rho(Q)$.
	
	Then, for $\varepsilon>0$, $Q\in \caD_0$, and $f\in L^1_{\loc}(\mu)$ we define the truncated operator
	\begin{equation}\label{eq:SIO_associated_with_cell}
		T_{\varepsilon,Q}(f\mu)(x)\coloneqq \chi_Q(x)\, \int_{\max(\varepsilon, 2 \rho(Q)/A_0) < d(x,y) \leq 2\rho(Q)} K(x,y)\, f(y)\,d\mu(y), \qquad x\in X.
	\end{equation}
	We claim that, for every $f\in L^1_{\loc}(\mu)$ whose support is contained in $B_0$ and $\mu$-almost every $x\in B_0$, we have
	\begin{equation}\label{eq:decomposition_T_mu_ep}
		T_{\mu,\varepsilon}f(x)=\sum_{Q\in \caD_0}T_{\varepsilon,Q}(f\,\mu)(x).
	\end{equation}
	
	To prove \eqref{eq:decomposition_T_mu_ep}, for $\mu$-a.e. $x\in \supp(\mu|_{B_0})$, let $Q_k(x)$ denote the unique cell in $\caD_0$ of generation $k$ which contains $x$.
	Since $\rho(Q_{k+1}(x))=\rho(Q_k(x))/A_0$ and
	\[
	d(x,y)\leq2r(B_0)=s<2\rho(Q_0(x)), \qquad \text{ for }x,y\in \supp(\mu|_{B_0}),
	\]
	the half-open annuli in \eqref{eq:SIO_associated_with_cell} form a disjoint partition of $\{y\in \supp(\mu|_{B_0}):d(x,y)>\varepsilon\}$.
	
	Moreover,
	\begin{equation}\label{eq:annular_size_estimate}
		|T_{\varepsilon,Q}(f\mu)(x)|\lesssim C\, A_0^n\, \frac{\mu(30\, B(Q)\cap B_0)}{\rho(Q)^n}\, \|f\|_{L^\infty(\mu)}\chi_Q(x),
	\end{equation}
	since $x\in Q\subset 28\, B(Q)$ and $\rho(Q)\leq r(Q)$ imply $B(x,2\rho(Q))\cap \supp(\mu|_{B_0})\subset 30\, B(Q)$.
	
	For $R\in\widetilde{\Top}$, we define the main \textit{mixed tree operators}
	\begin{equation}\label{eq:mixed_tree_operator}
		\widetilde{\mathcal K}_{\varepsilon,R}(f\mu)(x)\coloneqq \sum_{Q\in\widetilde{\Tree}(R)}T_{\varepsilon,Q}(f\mu)(x)
	\end{equation}
	and 
	\begin{equation}\label{eq:mixed_sum_operator}
		\widetilde\Sigma_\varepsilon f(x)\coloneqq\sum_{R\in\widetilde{\Top}}|\widetilde{\mathcal K}_{\varepsilon,R}(f\mu)(x)|.
	\end{equation}
	For $R_i\in\Top_i$ and	$P\in\Stopp_i(R_i)$, we denote by
	\begin{equation}\label{eq:intermediate_tree_operator}
		\mathcal K_{\varepsilon,i,P}(f\mu)(x)\coloneqq \sum_{Q\in\mathcal T_{i,P}}T_{\varepsilon,Q}(f\mu)(x)
	\end{equation}
	the SIO associated with the intermediate trees between stopping cells and the following top cells.
	We also define
	\begin{align*}
		\Sigma_{0,\varepsilon}f(x)&\coloneqq \sum_{Q\in\mathcal T_0}|T_{\varepsilon,Q}(f\mu)(x)|,\\
		\Sigma_{i,\varepsilon}'f(x)	&\coloneqq\sum_{R_i\in\Top_i'}\Bigl(\sum_{Q\in\Tree_i(R_i)}|T_{\varepsilon,Q}(f\mu)(x)|	+\sum_{P\in\Stopp_i'(R_i)}	\sum_{Q\in\mathcal T_{i,P}}|T_{\varepsilon,Q}(f\mu)(x)|\Bigr),\\
		\Sigma_{i,\varepsilon}''f(x)&\coloneqq\sum_{R_i\in\Top_i''}\sum_{P\in\Stopp_i''(R_i)}\sum_{Q\in\mathcal T_{i,P}}	|T_{\varepsilon,Q}(f\mu)(x)|,\qquad i=1,2.
	\end{align*}

	Thus, combining \eqref{eq:decomposition_D_0} and \eqref{eq:decomposition_T_mu_ep}, we obtain the pointwise bound
	\begin{equation}\label{eq:splitting_T_in_sigmas_tree}
		|T_{\mu,\varepsilon}f(x)|\leq \widetilde\Sigma_\varepsilon f(x) + \Sigma_{0,\varepsilon}f(x) + \sum_{i=1,2}\bigl(\Sigma_{i,\varepsilon}'f(x) + \Sigma_{i,\varepsilon}''f(x)\bigr).
	\end{equation}
	
	In the following lemmas, we estimate the $L^1(\mu|_{B_0})$-norms of each term on the right-hand side of \eqref{eq:splitting_T_in_sigmas_tree}, assuming that the operators are applied to $L^\infty(\mu)$-functions supported in $B_0$.
	
	We recall that $\mu_0=\mu|_{B_0}$ and $\mu_0\leq\mu$, and remark that for every corona cell $Q\subset B_0$, we have $\mu_0(Q)=\mu(Q)$. Then, \eqref{eq:comparability_r_rho} and the $n$-growth of $\mu$ imply
	\begin{equation}\label{eq:annular_density_normalization}
		\frac{\mu_0(30\,B(Q))}{\rho(Q)^n}\lesssim_{K_0,n}\Theta_{\mu}(100\, B(Q))\lesssim \mathfrak C_0.
	\end{equation}
	If $Q\subset B_0$ is doubling, its doubling property improves \eqref{eq:annular_density_normalization} to
	\begin{equation}\label{eq:doubling_endpoint_density_comparison}
		\frac{\mu_0(30\,B(Q))}{\rho(Q)^n}\lesssim_{K_0,n}\Theta_{\mu}(100 \, B(Q))\lesssim_{K_0,n}\Theta_\mu(Q).
	\end{equation}
	\vv

	We first bound the operator associated with cells in $\mathcal T_0$.
	\begin{lemma}\label{lemma:estimate_sigma_0_mixed corona}
		For every $g\in L^\infty(\mu)$ supported in $B_0$, we have
		\[
		\int_{B_0}|\Sigma_{0,\varepsilon} g(x)|\, d\mu(x)\leq \mathbf c^0_\Sigma\,  \|g\|_{L^\infty(\mu|_{B_0})}\, \mu(B_0),
		\]
		where $\mathbf c^0_\Sigma>0$ depends only on $n$, $C$, $A_0$, $K_0$, $\mathfrak C_0$, $\beta$, and $\mathfrak C_3$.
	\end{lemma}
	\begin{proof}
		Let $x\in B_0\cap\supp\mu$ be contained in arbitrarily small doubling cells, which holds for $\mu$-almost every $x\in B_0$. We denote by $R_x$ the cell of $\Top^0$ that contains $x$, and by $P_x$ the largest cell that contains $x$ and is contained in $B_0$, which exists because $B_0$ is open. Since the cells in $\Top^0$ are pairwise disjoint, a cell that contains $x$ belongs to $\mathcal T_0$ if and only if it strictly contains $R_x$. Hence, we split
		\[
		\Sigma_{0,\varepsilon}g(x)=\sum_{\substack{Q\ni x, Q\not\subset B_0}}\bigl|T_{\varepsilon,Q}(g\mu)(x)\bigr|+\sum_{R_x\subsetneq Q\subset P_x}\bigl|T_{\varepsilon,Q}(g\mu)(x)\bigr|\eqqcolon\mathfrak I_1(x)+\mathfrak I_2(x).
		\]
		
		To estimate $\mathfrak I_1(x)$, we set $t(x)\coloneqq\dist(x,X\setminus B_0)$, and we note that $\mathfrak I_1(x)=0$ if $X=B_0$. Every cell $Q$ in the first sum contains a point outside $B_0$, so $\ell(Q)\geq\diam(Q)\geq t(x)$. Since $\ell(Q)\leq56K_0s=112K_0\,r(B_0)$ and the side lengths of consecutive generations differ by a factor $A_0$, the first sum has at most $C_{A_0,K_0}\bigl(1+\log^+(r(B_0)/t(x))\bigr)$ terms. Thus, \eqref{eq:annular_size_estimate} and \eqref{eq:annular_density_normalization} give
		\[
		\mathfrak I_1(x)\lesssim\Bigl(1+\log^+\frac{r(B_0)}{t(x)}\Bigr)\,\|g\|_{L^\infty(\mu|_{B_0})}.
		\]
		
		To estimate $\mathfrak I_2(x)$, we note that the cells $Q$ with $R_x\subsetneq Q\subset P_x$ are non-doubling, by the maximality of $R_x$. Hence, \eqref{eq:annular_size_estimate}, \eqref{eq:annular_density_normalization}, and Lemma \ref{lemma:lemma_density_nonboubling_chain}, applied to the parent of $R_x$ and to $P_x$, give
		\[
		\mathfrak I_2(x)\lesssim\sum_{R_x\subsetneq Q\subset P_x}\Theta_\mu(100B(Q))\,\|g\|_{L^\infty(\mu|_{B_0})}\lesssim\Theta_\mu(100B(P_x))\,\|g\|_{L^\infty(\mu|_{B_0})}\leq\mathfrak C_0\,\|g\|_{L^\infty(\mu|_{B_0})}.
		\]
		
		Finally, we write $B_0=B(z_0,r_0)$, so that $t(x)\geq r_0-d(x,z_0)$. For $u>0$, the $\mathfrak C_3$-small boundary of $B_0$ and \eqref{eq:db_sec_4_start} give
		\[
		\mu\bigl(\{x\in B_0:t(x)<r_0e^{-u}\}\bigr)\leq\mathfrak C_3\,e^{-u}\,\mu(5B_0)\leq\mathfrak C_3\,\beta\,e^{-u}\,\mu(B_0),
		\]
		and hence
		\[
		\int_{B_0}\log^+\Bigl(\frac{r_0}{t(x)}\Bigr)\,d\mu(x)=\int_0^\infty\mu\bigl(\{x\in B_0:t(x)<r_0e^{-u}\}\bigr)\,du\leq\mathfrak C_3\,\beta\,\mu(B_0).
		\]
		Integrating the estimates for $\mathfrak I_1(x)$ and $\mathfrak I_2(x)$ over $B_0$, we obtain
		\[
		\int_{B_0}\Sigma_{0,\varepsilon}g\,d\mu\lesssim(1+\mathfrak C_3\,\beta)\,\|g\|_{L^\infty(\mu|_{B_0})}\,\mu(B_0).\qedhere
		\]
	\end{proof}
	\vv

	We next estimate $\Sigma'_{i,\varepsilon}$ and $\Sigma''_{i,\varepsilon}$. We recall the constants $\mathbf c_T$ and $\mathbf c'$ in \eqref{eq:packing_0_stop} and \eqref{eq:packing_top_density}, respectively.
	\begin{lemma}\label{lemma:estimate_sigma_tree_prime}
		If $g\in L^\infty(\mu)$ is supported in $B_0$, there exists $\mathbf c'_\Sigma>0$ such that
		\begin{equation}\label{eq:bound_integral_sigma_1_a}
			\int_{B_0}|\Sigma'_{1,\varepsilon} g(x)|\, d\mu(x) + 	\int_{B_0}|\Sigma'_{2,\varepsilon} g(x)|\, d\mu(x) \leq \mathbf c'_\Sigma\,  \mathbf c'\, \|g\|_{L^\infty(\mu|_{B_0})}\, \mu(B_0).
		\end{equation}
		The constant $\mathbf c'_\Sigma$ depends on $n$, $C$, $A_0$, $K_0$, $\mathfrak C_0$, and the chosen $\eta_*$.
	\end{lemma}
	\begin{proof}
		Let $R_i\in \Top'_i$, $i=1,2$. The cell $R_i$ is doubling, so \eqref{eq:doubling_endpoint_density_comparison} yields
		\[
		\frac{\mu_0(30\, B(R_i))}{\rho(R_i)^n}\lesssim_{K_0, n} \Theta_\mu(R_i).
		\]
		
		Moreover, if $Q\in\Tree_i(R_i)$ is such that $Q\neq R_i$, the definition of $\Stopp'$-cells implies that $\rho(Q)\gtrsim\eta_*\,\rho(R_i)$ and $\rho(Q)\leq\rho(R_i)$, so 
		\[
		\frac{\mu_0(30\, B(Q))}{\rho(Q)^n}\lesssim_{\eta_*,n} \frac{\mu(100\, B(R_i))}{\rho(R_i)^n}\lesssim_{K_0,n}\Theta_\mu(R_i),
		\]
		Also, $\rho(Q)\gtrsim\eta_*\rho(R_i)$ implies that $J(Q)-J(R_i)\lesssim1+\log(1/\eta_*)$ for $Q\in\Tree_i(R_i)$. The bound \eqref{eq:annular_size_estimate}, then, yields
		\begin{equation}\label{eq:lem_sigma_top_1_a}
			\begin{split}
				\int_{B_0}\sum_{Q\in \Tree_i(R_i)}|T_{\varepsilon, Q}(g\, \mu)|\, d\mu&\lesssim \int_{B_0}\sum_{Q\in \Tree_i(R_i)}\frac{\mu_0(30B(Q))}{\rho(Q)^n}\, \|g\|_{L^\infty(\mu|_{B_0})}\chi_Q(x)\, d\mu(x)\\
				&\lesssim_{\eta_*} \Theta_{\mu}(R_i)\mu(R_i)\, \|g\|_{L^\infty(\mu|_{B_0})}.
			\end{split}
		\end{equation}
		
		If  $P\in \Stopp'_i(R_i)$, all intermediate cells in $\mathcal T_{i,P}$ are non-doubling. Hence, the same summation argument as in Lemma \ref{lemma:estimate_sigma_0_mixed corona} implies that, for $\mu$-almost every $x\in P$, we have
		\begin{equation}\label{eq:estimate_sigma_0_mixed corona_aux}
			\sum_{Q\in \mathcal T_{i,P}, x\in Q}\frac{\mu_0(30\, B(Q))}{\rho(Q)^n}\lesssim\Theta_\mu(100B(P))\lesssim_{\eta_*,n}\Theta_\mu(R_i),
		\end{equation}
		where the last bound uses $\rho(P)\gtrsim\eta_*\rho(R_i)$, $100B(P)\subset100B(R_i)$, and the doubling property of $R_i$. If $P$ is doubling, the intermediate family is empty and the estimate is immediate.
		
		So, since cells in $\Stopp'_i(R_i)$ are pairwise disjoint, we obtain
		\begin{equation}\label{eq:lem_sigma_top_1_b}
			\begin{split}
				\int_{B_0}\sum_{P\in \Stopp'_i(R_i)}\sum_{Q\in \mathcal T_{i,P}}|T_{\varepsilon, Q}(g\, \mu)|\, d\mu&\overset{\eqref{eq:annular_size_estimate}}{\lesssim} \int_{B_0}\sum_{P\in \Stopp'_i(R_i)}\Theta_\mu(R_i)\, \|g\|_{L^\infty(\mu|_{B_0})}\chi_{P}(x)\, d\mu(x)\\
				&\lesssim \Theta_{\mu}(R_i)\mu(R_i)\, \|g\|_{L^\infty(\mu|_{B_0})}.
			\end{split}
		\end{equation}
		
		Thus, we gather \eqref{eq:lem_sigma_top_1_a}, \eqref{eq:lem_sigma_top_1_b}, the packing condition \eqref{eq:packing_top_density}, and obtain		
		\[
		\int_{B_0}|\Sigma'_{i,\varepsilon} g(x)|\, d\mu(x)\lesssim \sum_{R_i\in \Top'_i}\Theta_\mu(R_i)\, \mu(R_i)\, \|g\|_{L^\infty(\mu|_{B_0})}\leq \mathbf c'\, \|g\|_{L^\infty(\mu|_{B_0})}\, \mu(B_0).\qedhere
		\]
	\end{proof}

	\vv
	\begin{lemma}\label{lemma:estimate_sigma_tree_non_doubling}
		If $g\in L^\infty(\mu)$ is supported in $B_0$, we have
		\begin{equation}\label{eq:bound_integral_sigma_1_a_sec4_double_prime}
			\int_{B_0}|\Sigma''_{1,\varepsilon} g(x)|\, d\mu(x) + 	\int_{B_0}|\Sigma''_{2,\varepsilon} g(x)|\, d\mu(x) \leq \mathbf c_\Sigma'' \mathbf c^{-1}_T\, \|g\|_{L^\infty(\mu|_{B_0})}\, \mu(B_0),
		\end{equation}
		where $\mathbf c''_\Sigma>0$ depends on $n$, $C$, $A_0$, $K_0$, and $\mathfrak C_0$.
	\end{lemma}
	\begin{proof}
		The bound \eqref{eq:bound_integral_sigma_1_a_sec4_double_prime} can be proved similarly to \eqref{eq:lem_sigma_top_1_b}. Indeed, let $i=1,2$ and $R_i\in \Top''_i$.
		Note that, if $P\in \Stopp''_i(R_i)$, then all intermediate cells in $\mathcal T_{i,P}$ are non-doubling. So, as in \eqref{eq:estimate_sigma_0_mixed corona_aux}, we can apply once again the summation argument in Lemma \ref{lemma:estimate_sigma_0_mixed corona} and, if $P\in \Stopp''_i(R_i)$, by \eqref{eq:annular_density_normalization}, for $\mu$-almost every $x\in P$ we obtain
		\begin{equation}\label{eq:estimate_sigma_0_mixed corona_aux2}
			\sum_{Q\in \mathcal T_{i,P}, x\in Q} \frac{\mu_0(30\, B(Q))}{\rho(Q)^n}\lesssim \frac{\mu(100\, B(P))}{\rho(P)^n}\lesssim \mathfrak C_0.
		\end{equation}
		Since cells in $\Stopp''_i(R_i)$ are pairwise disjoint subsets of $R_i$, by \eqref{eq:estimate_sigma_0_mixed corona_aux2} and \eqref{eq:annular_size_estimate} we have
		\begin{equation}\label{eq:lem_sigma_top_1_b_sec4_double_prime}
			\begin{split}
				\int_{B_0}\sum_{P\in \Stopp''_i(R_i)}\sum_{Q\in \mathcal T_{i,P}}|T_{\varepsilon, Q}(g\, \mu)|\, d\mu&\overset{}{\lesssim} \int_{B_0}\sum_{P\in \Stopp''_i(R_i)} \|g\|_{L^\infty(\mu|_{B_0})}\chi_{P}(x)\, d\mu(x)\\
				&\leq\mu(R_i) \|g\|_{L^\infty(\mu|_{B_0})}.
			\end{split}
		\end{equation}
		
		Recall that, by Lemma \ref{lemma:packing_top_2}, $\Top''_i$ is a Carleson family, so we conclude that
		\[
		\int_{B_0}|\Sigma''_{i,\varepsilon} g(x)|\, d\mu(x)\lesssim \sum_{R_i\in \Top''_i}\mu(R_i)\, \|g\|_{L^\infty(\mu|_{B_0})}\leq \mathbf c^{-1}_T\, \|g\|_{L^\infty(\mu|_{B_0})}\, \mu(B_0).\qedhere
		\]
	\end{proof}

	\vv
	We now estimate $\widetilde\Sigma_\varepsilon$ under assumption \eqref{eq:top_L2_bd_assumption}.
	Let $R_i \in \Top''_i$ for $i=1,2$ be such that $\widetilde{\Tree}(R_1, R_2) \neq \varnothing$, hence either $R_1\subset R_2$ or $R_2\subset R_1$. We define
	\begin{equation}\label{eq:Phi-i}
		\widetilde \Phi_{R_1, R_2}(x) \coloneqq \max\big(\Phi_{1,R_1}(x), \Phi_{2,R_2}(x)\big).
	\end{equation}
	With yet another minor abuse of notation, for $R \in \widetilde{\Top}$ we write
	\begin{equation}\label{eq:tilde Phi}
		\widetilde \Phi_R(x) \coloneqq \widetilde \Phi_{R_1, R_2}(x),
	\end{equation}
	where we recall that $R$ denotes the smaller of the two cells $R_1$ and $R_2$, and note that $\widetilde\Phi_R$ is a $1$-Lipschitz function.
	
	\begin{lemma}\label{lemma:estim_sigma_tree}
		Let $\varepsilon>0$ and assume that
		\begin{equation}\label{eq:top_L2_bd_assumption}
			\sup_{R \in \widetilde \Top} \| T_{\widetilde{\Phi}_R, \mu|_{84B(R)}, \varepsilon}\|_{L^2(\mu|_{84B(R)})  \to L^2(\mu|_{84B(R)}) } \leq \mathfrak{C_\mathcal{T} }.
		\end{equation}
		Then, for every $g\in L^\infty(\mu)$ supported in $B_0$, we have
		\begin{equation}\label{eq:bound_integral_sigma_1}
			\int_{B_0}| \widetilde \Sigma_{\varepsilon} g(x)|\, d\mu(x)  \leq \mathbf c'''_\Sigma \, \mathbf c_T^{-1}\, (1+ \mathfrak C_{\mathcal T})\|g\|_{L^\infty(\mu|_{B_0})} \,\mu(B_0),
		\end{equation}
		where $\mathbf c'''_\Sigma>0$ depends only on $n$, $C$, $\gamma$, $\mathfrak C_0$, $A_0$, and $K_0$.
	\end{lemma}
	
	{
		\begin{proof}
			We fix $R\in\widetilde{\Top}$, and assume that it is the root of $\widetilde{\Tree}(R_1,R_2)$. We first prove that, for every $x\in R$,
			\begin{equation}\label{eq:claim_bd:tree:SIo}
				|\widetilde{\mathcal K}_{\varepsilon,R}(g\mu)(x)|
				\lesssim\|g\|_{L^\infty(\mu|_{B_0})}
				+|T_{\widetilde\Phi_R,\varepsilon}(\chi_{84B(R)}g\mu)(x)|.
			\end{equation}
			We denote by $\widetilde{\Stopp}(R)$ the family of the maximal cells of \[\{Q\in\Stopp''_1(R_1)\cup\Stopp''_2(R_2):Q\subset R\},\] which are proper subcells of $R$ because $R\in\Tree_1(R_1)\cap\Tree_2(R_2)$. For $i=1,2$, every $x\in R$ with $\Phi_{i,R_i}(x)>0$ belongs to a cell of $\Stopp''_i(R_i)$, which is contained in $R$ because $R\in\Tree_i(R_i)$, and hence to a cell of $\widetilde{\Stopp}(R)$.
			
			Let $x\in R$. If $x$ belongs to a cell $Q\in\widetilde{\Stopp}(R)$, then, by the maximality of $Q$, the cells of $\widetilde{\Tree}(R)$ that contain $x$ are the cells $P$ with $Q\subsetneq P\subset R$. Thus, the definition of $T_{\varepsilon,P}$ gives
			\begin{equation}\label{eq:case_1_suppressed_tree}
				\widetilde{\mathcal K}_{\varepsilon,R}(g\mu)(x)=\int_{\max(\varepsilon,a(x))<d(x,y)\leq2\rho(R)}K(x,y)\,g(y)\,d\mu(y),\qquad\text{ where }\quad a(x)\coloneqq2\rho(Q).
			\end{equation}
			Moreover, if $Q\in\Stopp''_1(R_1)$, then \eqref{eq:reg_cells_sec4} gives $\Phi_{1,R_1}(x)>480A_0\,\ell(Q)$ and, if $\Phi_{2,R_2}(x)>0$, then $x$ belongs to a cell $S\in\Stopp''_2(R_2)$ such that $S\subset Q$, by the maximality of $Q$. Hence, \eqref{eq:comparability_reg_supp_prop_1}, applied to $Q$ and to $S$, and the analogous argument for $Q\in\Stopp''_2(R_2)$ give
			\begin{equation}\label{eq:comparability_tilde_Phi_R}
				480A_0\,\ell(Q)<\widetilde\Phi_R(x)\leq481A_0^2\,\ell(Q).
			\end{equation}
			If $x$ does not belong to any cell of $\widetilde{\Stopp}(R)$, then $\widetilde\Phi_R(x)=0$, every cell of $\caD_0(R)$ that contains $x$ belongs to $\widetilde{\Tree}(R)$, and \eqref{eq:case_1_suppressed_tree} holds with $a(x)\coloneqq0$.
			
			In both cases, we set $\varepsilon_x\coloneqq\max(\varepsilon,\widetilde\Phi_R(x))$, and we split
			\[
			\begin{split}
				\widetilde{\mathcal K}_{\varepsilon,R}(g\mu)(x)&=T_{\widetilde\Phi_R,\varepsilon_x}(\chi_{84B(R)}g\mu)(x)+\bigl(T_{\varepsilon_x}-T_{\widetilde\Phi_R,\varepsilon_x}\bigr)(\chi_{84B(R)}g\mu)(x)\\
				&\quad+\bigl(\widetilde{\mathcal K}_{\varepsilon,R}(g\mu)(x)-T_{\varepsilon_x}(\chi_{84B(R)}g\mu)(x)\bigr)\eqqcolon\mathfrak I_1+\mathfrak I_2+\mathfrak I_3.
			\end{split}
			\]
			If $\varepsilon_x=\varepsilon$, then $\mathfrak I_1=T_{\widetilde\Phi_R,\varepsilon}(\chi_{84B(R)}g\mu)(x)$. Otherwise, $\widetilde\Phi_R(x)>\varepsilon$, and \eqref{eq:bound_suppressed_kernel_1} and the $n$-growth of $\mu$ give
			\[
			\begin{split}
				|\mathfrak I_1-T_{\widetilde\Phi_R,\varepsilon}(\chi_{84B(R)}g\mu)(x)|&\leq\int_{\varepsilon<d(x,y)\leq\widetilde\Phi_R(x)}|K_{\widetilde\Phi_R}(x,y)|\,|g(y)|\,d\mu(y)\\
				&\leq C\,\frac{\mu(B(x,2\widetilde\Phi_R(x)))}{\widetilde\Phi_R(x)^n}\,\|g\|_{L^\infty(\mu|_{B_0})}\leq C\,\mathfrak C_0\,\|g\|_{L^\infty(\mu|_{B_0})}.
			\end{split}
			\] If $\widetilde\Phi_R(x)>0$, then \eqref{eq:bound_suppressed_kernel_2} and the $n$-growth of $\mu$ give $|\mathfrak I_2|\lesssim\mathfrak C_0\,\|g\|_{L^\infty(\mu|_{B_0})}$, while $\mathfrak I_2=0$ if $\widetilde\Phi_R(x)=0$, because $K_{\widetilde\Phi_R}(x,\cdot)=K(x,\cdot)$. To estimate $\mathfrak I_3$, we note that $x\in R\subset28B(R)$ gives $B(x,2\rho(R))\subset84B(R)\subset B(x,112\,r(R))$. Hence, by \eqref{eq:case_1_suppressed_tree}, $\mathfrak I_3$ is a difference of integrals of $K(x,\cdot)\,g$ over two sets whose symmetric difference is contained in
			\[
			\{y:a(x)<d(x,y)\leq\widetilde\Phi_R(x)\}\cup\{y:2\rho(R)<d(x,y)<112\,r(R)\}.
			\]
			By \eqref{eq:comparability_tilde_Phi_R} and \eqref{eq:comparability_r_rho}, the ratio between the outer and the inner radius of each of these annuli is bounded by a constant depending only on $A_0$ and $K_0$, so the size estimate of $K$ and the $n$-growth of $\mu$ give $|\mathfrak I_3|\lesssim\|g\|_{L^\infty(\mu|_{B_0})}$. This proves \eqref{eq:claim_bd:tree:SIo}.
			
			Since $\widetilde{\mathcal K}_{\varepsilon,R}(g\mu)$ vanishes outside $R$, \eqref{eq:claim_bd:tree:SIo} gives
			\begin{equation}\label{eq:splitting_0_supp_tree}
				\begin{split}
					\int_{B_0}|\widetilde\Sigma_\varepsilon g|\,d\mu&\leq\sum_{R\in\widetilde{\Top}}\int_R|\widetilde{\mathcal K}_{\varepsilon,R}(g\mu)|\,d\mu\\
					&\lesssim\sum_{R\in\widetilde{\Top}}\|g\|_{L^\infty(\mu|_{B_0})}\,\mu(R)+\sum_{R\in\widetilde{\Top}}\int_R \bigl|T_{\widetilde\Phi_R,\varepsilon}(\chi_{84B(R)}g\mu)\bigr|\,d\mu\eqqcolon\mathfrak J_1+\mathfrak J_2.
				\end{split}
			\end{equation}
			
			To bound $\mathfrak J_1$ we apply the Carleson packing of Lemma \ref{lemma:packing2_mixed}, which gives
			\begin{equation}\label{eq:splitting_1_supp_tree}
				\mathfrak J_1\leq2\,\mathbf c_T^{-1}\,\|g\|_{L^\infty(\mu|_{B_0})}\,\mu(B_0).
			\end{equation}
			
			To estimate $\mathfrak J_2$, we fix $R\in\widetilde{\Top}$ and we set $\sigma\coloneqq\mu|_{84B(R)}$. Since $R$ is doubling, we have \[\mu(84B(R))\leq\mu(100B(R))\leq K_0\,\mu(B(R))\leq K_0\,\mu(R),\] so the Cauchy-Schwarz inequality and \eqref{eq:top_L2_bd_assumption} give
			\[
			\begin{split}
				\int_R|T_{\widetilde\Phi_R,\varepsilon}(\chi_{84B(R)}g\mu)|\,d\mu&\leq\mu(R)^{1/2}\,\|T_{\widetilde\Phi_R,\sigma,\varepsilon}g\|_{L^2(\sigma)}\leq\mathfrak C_{\mathcal T}\,\mu(R)^{1/2}\,\|g\|_{L^2(\sigma)}\\
				&\leq K_0^{1/2}\,\mathfrak C_{\mathcal T}\,\|g\|_{L^\infty(\mu|_{B_0})}\,\mu(R).
			\end{split}
			\]
			Summing over $R\in\widetilde{\Top}$ and applying Lemma \ref{lemma:packing2_mixed} once again, we obtain
			\[
			\mathfrak J_2\leq2K_0^{1/2}\,\mathbf c_T^{-1}\,\mathfrak C_{\mathcal T}\,\|g\|_{L^\infty(\mu|_{B_0})}\,\mu(B_0),
			\]
			which, together with \eqref{eq:splitting_0_supp_tree} and \eqref{eq:splitting_1_supp_tree}, concludes the proof of the lemma.
		\end{proof}
	}
	
	\vv
	
	Finally, gathering all the estimates of this section, we obtain the following result.
	\begin{corollary}\label{corollary:49}
		Let $B_0$ be a ball with $\mathfrak C_3$-small boundary with respect to $\mu$ that satisfies the doubling-type assumption \eqref{eq:db_sec_4_start}. Assume that there exists  $ \mathfrak{C_\mathcal{T}} >0$, independent of $\varepsilon$,  such that
		\begin{equation}\label{eq:corollary:49b}
			\sup_{R \in \widetilde \Top(B_0)} \| T_{\widetilde{\Phi}_R,\mu|_{84B(R)}, \varepsilon}\|_{L^2(\mu|_{84B(R)})  \to L^2(\mu|_{84B(R)}) } \leq \mathfrak{C_\mathcal{T} }.
		\end{equation}
		Then for every $g \in L^\infty(\mu)$ supported in $B_0$ and every $\varepsilon>0$, it holds that
		\begin{equation}\label{eq:corollary:49}
			\|T_{\mu, \varepsilon} g \|_{L^1(\mu|_{B_0}) } + 	\|T^*_{\mu, \varepsilon} g \|_{L^1(\mu|_{B_0}) } \lesssim \|g\|_{L^\infty(\mu|_{B_0})} \,\mu(B_0),
		\end{equation}
		where the implicit constant depends on $A_0,K_0,n,C,\gamma,\mathfrak C_0,\mathfrak C_3,\beta,\eta_*, \mathbf c_T, \mathbf c', \mathfrak{C_\mathcal{T} }$, but not on $\varepsilon$.
	\end{corollary}
	\begin{proof}
		By \eqref{eq:splitting_T_in_sigmas_tree}, Lemma \ref{lemma:estimate_sigma_0_mixed corona}, Lemma \ref{lemma:estimate_sigma_tree_prime}, Lemma \ref{lemma:estimate_sigma_tree_non_doubling}, and Lemma \ref{lemma:estim_sigma_tree}, for every $g\in L^\infty(\mu)$ supported in $B_0$, it holds
		\[
		\int_{B_0} |T_{\mu, \varepsilon} g|\, d\mu \leq \bigl(\mathbf c^0_\Sigma + \mathbf c'_\Sigma\, \mathbf c'+ \mathbf c''_\Sigma\, \mathbf c^{-1}_T+ \mathbf c'''_\Sigma \, \mathbf c_T^{-1}\, (1+ \mathfrak C_{\mathcal T})\bigr)\|g\|_{L^\infty(\mu|_{B_0})}\, \mu(B_0).
		\]
		The same arguments apply to $T^*_\varepsilon$, since the suppression factor is symmetric and the localized $L^2$ bound also holds for the transpose. This proves the corollary.
	\end{proof}
	\vv

	\section{Bounds of different norms for the suppressed operator}\label{section:L2bd}
	
	{
		Let $R_1\in\Top''_1$ and $R_2\in\Top''_2$ be such that $\widetilde{\Tree}(R_1,R_2)\neq\varnothing$, so that either $R_1\subset R_2$ or $R_2\subset R_1$. Let $i=1,2$. We recall the following objects from Section \ref{sec:corona}:
		\begin{itemize}
			\item The testing ball $B_i=\widetilde B(R_i)$.
			\item The testing measures $\widehat\nu_1$ for $T$ on $B_1$ and $\widehat\nu_2$ for $T^*$ on $B_2$.
			\item The measure $\nu_i=\widehat\nu_i+\mu|_{X\setminus B_i}$.
			\item The sets $F_i(w_i)\subset(1-\theta)B_i$, for $w_i\in\Omega$.
			\item The suppressions $\Phi_{i,R_i,w_i}$ and $\Phi_{i,R_i}$ from \eqref{eq:definition_suppression_Phi_R_i} and \eqref{eq:definition_suppressions_2a}, respectively.
		\end{itemize} 
		
		We set $\mu_i\coloneqq\mu|_{B_i}$, and we note that $|\nu_i|=|\widehat\nu_i|+\mu|_{X\setminus B_i}$, since $|\widehat\nu_i|(X\setminus B_i)=0$.
		For $\mathtt w=(w_1,w_2)\in\Omega^2$ we define
		\begin{equation}\label{eq:definition_suppressions_2ab}
			\Xi_i(x)\coloneqq\Xi_{i,R_i,w_i}(x)
			\coloneqq\max\bigl\{\Phi_{i,R_i}(x),\Phi_{i,R_i,w_i}(x)\bigr\},\qquad i=1,2,
		\end{equation}
		and set
		\begin{equation}\label{eq:definition_suppressions_2b}
			\Xi_{\mathtt w}(x)
			\coloneqq\max\bigl\{\Xi_{1,R_1,w_1}(x),\Xi_{2,R_2,w_2}(x)\bigr\}.
		\end{equation}
		These functions are non-negative and $1$-Lipschitz. In this section, the implicit constants depend only on $n$, $D$, $\gamma$, $\varepsilon_0$, $C$, $\mathfrak C_0$, $\mathfrak C_1$, $\mathfrak C_3$, $\beta$, and the fixed David-Mattila parameters. In particular, they are independent of $R_1$, $R_2$, $\mathtt w$, and of the truncation parameter $\varepsilon>0$.
		\vv
		
		\underline{\it Growth and comparison of the measures.}
		Let $i=1,2$. We denote by $r_{1,i}$ and $r_{2,i}$ the high- and low-density radii defined after \eqref{eq:def_e_lambda_mu_nu}, and by $\mathcal R_i$ the radius in \eqref{eq:definition_r_i}, all associated with $\widehat\nu_i$ and $B_i$. By \eqref{eq:def_tot_exceptional_set} and \eqref{eq:W_new}, the set $F_i(w_i)\subset G_i$ does not intersect $\mathscr S_{i,R_i}\cup\mathscr H_{i,R_i}$. The latter union contains the balls $B(x,r_{j,i}(x))$, for $x\in B_i$ and $j=1,2$, and the balls $B(x,\mathcal R_i(x))$, for $x\in X$. Hence, by \eqref{eq:definition_suppression_Phi_R_i} and \eqref{eq:definition_suppressions_2ab}, we have
		\begin{equation}\label{eq:density_radii_sec6}
			r_{j,i}(x)\leq\dist(x,F_i(w_i))
			\leq\Xi_i(x)\leq\Xi_{\mathtt w}(x),
			\qquad x\in B_i,\quad j=1,2,
		\end{equation}
		and $\mathcal R_i(x)\leq\dist(x,F_i(w_i))\leq\Xi_{\mathtt w}(x)$ for $x\in X$. Therefore, by the definitions of the exceptional radii, for $x\in B_i$ we have
		\begin{align}
			|\widehat\nu_i|(B(x,t))&\leq m_1\,\mu(B(x,5t)),\qquad \text{ for } t>r_{1,i}(x),\label{eq:high_density_sec6}\\
			\mu_i(B(x,t))&\leq m_2\,|\widehat\nu_i|(B(x,5t)), \qquad\text{ for } t>r_{2,i}(x),\label{eq:low_density_sec6}
		\end{align}
		and, for $x\in X$ and $t>\Xi_{\mathtt w}(x)$, we have
		\begin{equation}
			|\widehat\nu_i|(B(x,t))\leq\mathbf m\,t^n.\label{eq:growth_nu_sec6}
		\end{equation}
		
		If $x\in X$ and $t>0$, then $t+\Xi_{\mathtt w}(x)>\Xi_{\mathtt w}(x)$. Hence, \eqref{eq:growth_nu_sec6} and the $n$-growth of $\mu$ give, for every $\mathfrak m\in\{\mu,\mu_i,|\widehat\nu_i|,|\nu_i|\}$, the bound
		\begin{equation}\label{eq:growth_auxiliary_measures_sec6}
			\mathfrak m(B(x,t))\leq(\mathbf m+\mathfrak C_0)\,(t+\Xi_{\mathtt w}(x))^n,
			\qquad x\in X, t>0.
		\end{equation}
		\vv
		
		The measures $\mu$, $|\nu_1|$, and $|\nu_2|$ can also be compared above the scale of the suppression. Namely, for all $\mathfrak m,\mathfrak m'\in\{\mu,|\nu_1|,|\nu_2|\}$, $x\in X$, and $t>0$ with $t\geq\Xi_{\mathtt w}(x)$, we have
		\begin{equation}\label{eq:comparison_auxiliary_measures_sec6}
			\mathfrak m(B(x,t))\leq(1+m_1)(1+m_2)\,\mathfrak m'\bigl(B(x,121\, t)\bigr).
		\end{equation}
		To prove \eqref{eq:comparison_auxiliary_measures_sec6}, let $x\in X$ and $t>0$ with $t\geq\Xi_{\mathtt w}(x)$. We first show that
		\begin{equation*}
			|\widehat\nu_i|(B(x,t))\leq m_1\,\mu(B(x,11t))
			\qquad\text{and}\qquad
			\mu_i(B(x,t))\leq m_2\,|\widehat\nu_i|(B(x,11t)).
		\end{equation*}
		If $B(x,t)\cap B_i=\varnothing$, both left-hand sides vanish. Otherwise, we choose $z\in B(x,t)\cap B_i$. Since $\Xi_{\mathtt w}$ is $1$-Lipschitz, \eqref{eq:density_radii_sec6} gives $r_{j,i}(z)\leq\Xi_{\mathtt w}(z)<2t$, for $j=1,2$. Moreover, $B(x,t)\subset B(z,2t)$ and $B(z,10t)\subset B(x,11t)$. Thus, the two inequalities follow from \eqref{eq:high_density_sec6} and \eqref{eq:low_density_sec6}, applied with $z$ and $2t$ in place of $x$ and $t$. Since \[\mu(B(x,t)\setminus B_i)\leq\min\bigl\{\mu\bigl(B(x,11t)\bigr),|\nu_i|\bigl(B(x,11t)\bigr)\bigr\},\] we deduce that
		\[
		|\nu_i|(B(x,t))\leq(1+m_1)\,\mu\bigl(B(x,11t)\bigr)
		\qquad\text{and}\qquad
		\mu(B(x,t))\leq(1+m_2)\,|\nu_i|\bigl(B(x,11t)\bigr).
		\]
		
		These inequalities yield \eqref{eq:comparison_auxiliary_measures_sec6} for $\mathfrak m=|\nu_i|$ and $\mathfrak m'\in\{|\nu_1|,|\nu_2|\}$. The remaining cases readily follow.
		
		\vv
		
		\subsection{The localized operator}
		Let $i=1,2$, and let $z_{B_i}$ be the center of $B_i$. We define
		\[
		\chi_i(x)\coloneqq
		\min\left\{1,\max\left\{0,
		\frac{2\bigl(r(B_i)-d(x,z_{B_i})\bigr)}{\theta r(B_i)}
		\right\}\right\}.
		\]
		Then $\chi_i=1$ on $(1-\theta/2)B_i$ and $\chi_i=0$ on $X\setminus B_i$. The triangle inequality also implies that
		\[
		|\chi_i(x)-\chi_i(x')|\leq 2\, \frac{d(x,x')}{\theta r(B_i)}.
		\]
		Moreover, if $\chi_i(x)\neq1$, then $d(x,z_{B_i})>(1-\theta/2)r(B_i)$. Hence, since $F_i(w_i)\subset(1-\theta)B_i$, we obtain
		\begin{equation}\label{eq:collar_bound_sec6}
			\Xi_{\mathtt w}(x)
			\geq\dist(x,F_i(w_i))
			\geq d(x,z_{B_i})-(1-\theta)r(B_i)
			\geq\frac\theta2\, r(B_i).
		\end{equation}
		
		Set $r_0\coloneqq\min\{r(B_1),r(B_2)\}$. Let $\varphi$ and $\psi$ be the functions fixed in Section \ref{subsection:bounded_testing_reduction}. For $\varepsilon>0$, we define the localized kernel
		\begin{equation}\label{eq:definition_auxiliary_kernel_sec6}
			L_\varepsilon(x,y)
			\coloneqq\chi_2(x)\chi_1(y)
			\psi\!\left(\frac{d(x,y)}{r_0}\right)
			\varphi\!\left(\frac{d(x,y)}\varepsilon\right)K_{\Xi_{\mathtt w}}(x,y)
		\end{equation}
		and we write $L_\varepsilon(x,y)= \varrho_\varepsilon(x,y)\,K_{\Xi_{\mathtt w}}(x,y)$, where
		\[
		\varrho_\varepsilon(x,y)\coloneqq\chi_2(x)\,\chi_1(y)\,a_\varepsilon(d(x,y))
		\qquad\text{and}\qquad
		a_\varepsilon(t)\coloneqq\varphi(t/\varepsilon)\,\psi(t/r_0),\quad t\geq0.
		\]
		Then, for a positive measure $\mathfrak m$, we define the ``localized'' operator 
		\[
		S_{\mathfrak m,\varepsilon}f(x)\coloneqq\int L_\varepsilon(x,y)f(y)\,d\mathfrak m(y), \qquad f\in L^1_{\loc}(\mathfrak m).
		\]
		
		We note that, since $\varphi'(t/\varepsilon)$ and $\psi'(t/r_0)$ vanish unless $\varepsilon\leq t\leq2\varepsilon$ and $r_0\leq t\leq2r_0$, respectively, we have $|a_\varepsilon'(t)|\lesssim1/t$. Moreover, $0\leq a_\varepsilon\leq1$, and $a_\varepsilon(t)=0$ for $t\geq2r_0$.
		
		Let $x,x',y\in X$ with $d(x,x')\leq d(x,y)/2$, so that $d(x',y)\geq d(x,y)/2$. 
		Since $0\leq\chi_1,\chi_2\leq1$ and $r(B_2)\geq r_0$, the Lipschitz bound for $\chi_2$ and the mean value theorem give
		\[
		\begin{split}
			|\varrho_\varepsilon(x,y)-\varrho_\varepsilon(x',y)|
			&\leq|\chi_2(x)-\chi_2(x')|\,a_\varepsilon(d(x,y))+ \bigl|a_\varepsilon(d(x,y))-a_\varepsilon(d(x',y))\bigr|\\
			&\lesssim\frac{d(x,x')}{r_0}\,\chi_{\{d(x,y)<2r_0\}}+\frac{d(x,x')}{d(x,y)}
			\lesssim\frac{d(x,x')}{d(x,y)}.
		\end{split}
		\]
		Hence, by Lemma \ref{lemma:kernel_suppressed}, $|\varrho_\varepsilon|\leq1$, and $0<\gamma\leq1$,
		\[
		\begin{split}
			\bigl|L_\varepsilon(x,y)-L_\varepsilon(x',y)\bigr|
			&\leq\bigl|K_{\Xi_{\mathtt w}}(x,y)-K_{\Xi_{\mathtt w}}(x',y)\bigr|
			+|\varrho_\varepsilon(x,y)-\varrho_\varepsilon(x',y)|\,|K_{\Xi_{\mathtt w}}(x',y)|\\
			&\lesssim\frac{d(x,x')^\gamma}{d(x,y)^{n+\gamma}}+\frac{d(x,x')}{d(x,y)^{n+1}}
			\lesssim\frac{d(x,x')^\gamma}{d(x,y)^{n+\gamma}}.
		\end{split}
		\]
		Since $r(B_1)\geq r_0$, the same argument applies to the second variable, and we obtain
		\begin{equation}\label{eq:regularity_auxiliary_kernel_sec6}
			\bigl|L_\varepsilon(x,y)-L_\varepsilon(x',y)\bigr|
			+\bigl|L_\varepsilon(y,x)-L_\varepsilon(y,x')\bigr| \lesssim\frac{d(x,x')^\gamma}{d(x,y)^{n+\gamma}},
			\qquad d(x,x')\leq\tfrac12d(x,y).
		\end{equation}
		
		Let $\Psi_\varepsilon\coloneqq\Xi_{\mathtt w}+\varepsilon$. Since $L_\varepsilon(x,y)=0$ if $d(x,y)\leq\varepsilon$, Lemma \ref{lemma:kernel_suppressed} gives
		\begin{equation}\label{eq:size_auxiliary_kernel_sec6}
			\begin{split}
				|L_\varepsilon(x,y)|
				&\lesssim\min\bigl\{d(x,y)^{-n},\varepsilon^{-n},
				\Xi_{\mathtt w}(x)^{-n},\Xi_{\mathtt w}(y)^{-n}\bigr\}\\
				&\lesssim\min\bigl\{d(x,y)^{-n},\Psi_\varepsilon(x)^{-n},
				\Psi_\varepsilon(y)^{-n}\bigr\}.
			\end{split}
		\end{equation}
		The function $\Psi_\varepsilon$ is positive and $1$-Lipschitz. Let $\mathfrak m\in\{\mu,|\nu_1|,|\nu_2|\}$, $x\in X$, and $t\geq\Psi_\varepsilon(x)$. Then $t\geq\Xi_{\mathtt w}(x)$, so that \eqref{eq:comparison_auxiliary_measures_sec6} applies, and \eqref{eq:growth_auxiliary_measures_sec6} gives \[\mathfrak m(B(x,t))\leq2^n(\mathbf m+\mathfrak C_0)\,t^n.\] Therefore, by \eqref{eq:regularity_auxiliary_kernel_sec6} and \eqref{eq:size_auxiliary_kernel_sec6}, the kernel $L_\varepsilon$ and its transpose, together with $\Psi_\varepsilon$ and any $\mathfrak m_1,\mathfrak m_2,\mathfrak m_3\in\{\mu,|\nu_1|,|\nu_2|\}$, satisfy the hypotheses of Lemma \ref{lemma:mixed_cotlar_transfer} other than the $L^2$-bound, with $\delta_0=\varepsilon$ and with constants independent of $\varepsilon$ and $\mathtt w$.
		
		\begin{lemma}\label{lemma:l51comparison}
			For every $\varepsilon>0$ and every choice $\mathfrak m_i\in\{\mu_i,|\widehat\nu_i|\}$, $i=1,2$,
			\begin{equation}\label{eq:bound_auxiliary_difference_sec6}
				\|T_{\Xi_{\mathtt w},\mathfrak m_1,\varepsilon}
				-S_{\mathfrak m_1,\varepsilon}\|_
				{L^2(\mathfrak m_1)\to L^2(\mathfrak m_2)}\lesssim1.
			\end{equation}
		\end{lemma}
		\begin{proof}
			Since $\mathfrak m_i$ is concentrated on $B_i$, the kernel of the operator in \eqref{eq:bound_auxiliary_difference_sec6} is
			\[
			K_{\Xi_{\mathtt w}}(x,y)\chi_{\{d(x,y)>\varepsilon\}}-L_\varepsilon(x,y)
			=\mathfrak D_1(x,y)+\mathfrak D_2(x,y),
			\qquad x\in B_2,\quad y\in B_1,
			\]
			where we split
			\[
			\begin{split}
				\mathfrak D_1(x,y)&\coloneqq\Bigl[\chi_{\{d(x,y)>\varepsilon\}}-\varphi\Bigl(\frac{d(x,y)}\varepsilon\Bigr)\Bigr]K_{\Xi_{\mathtt w}}(x,y),\\
				\mathfrak D_2(x,y)&\coloneqq\Bigl[1-\chi_2(x)\chi_1(y)\psi\Bigl(\frac{d(x,y)}{r_0}\Bigr)\Bigr]\varphi\Bigl(\frac{d(x,y)}\varepsilon\Bigr)K_{\Xi_{\mathtt w}}(x,y).
			\end{split}
			\]
			
			Since $\varphi=0$ on $[0,1]$ and $\varphi=1$ on $[2,\infty)$, the kernel $\mathfrak D_1(x,y)$ vanishes unless $\varepsilon<d(x,y)<2\varepsilon$, and $|\mathfrak D_1|\leq|K_{\Xi_{\mathtt w}}|$. Hence, by Lemma \ref{lemma:kernel_suppressed} and \eqref{eq:growth_auxiliary_measures_sec6}, for $\mathfrak m\in\{\mathfrak m_1,\mathfrak m_2\}$ and $z\in X$ we have
			\[
			\int|\mathfrak D_1(z,y)|\,d\mathfrak m(y)+\int|\mathfrak D_1(x,z)|\,d\mathfrak m(x)
			\lesssim\frac{\mathfrak m(B(z,2\varepsilon))}{\max\{\varepsilon,\Xi_{\mathtt w}(z)\}^n}
			\lesssim\frac{(2\varepsilon+\Xi_{\mathtt w}(z))^n}{\max\{\varepsilon,\Xi_{\mathtt w}(z)\}^n}\leq3^n.
			\]
			By Schur's test, the integral operator with kernel $\mathfrak D_1$ is bounded from $L^2(\mathfrak m_1)$ to $L^2(\mathfrak m_2)$, with norm $\lesssim1$.
			
			If $\mathfrak D_2(x,y)\neq0$, then $d(x,y)>r_0$, $\chi_2(x)\neq1$, or $\chi_1(y)\neq1$. Thus, by \eqref{eq:collar_bound_sec6} and $r(B_i)\geq r_0$, we have \[\max\bigl\{d(x,y),\Xi_{\mathtt w}(x),\Xi_{\mathtt w}(y)\bigr\}\geq\frac{\theta r_0}{2}.\] Since $|\mathfrak D_2|\leq|K_{\Xi_{\mathtt w}}|$, Lemma \ref{lemma:kernel_suppressed} implies
			\begin{equation}\label{eq:positive_bound_auxiliary_difference_sec6}
				|\mathfrak D_2(x,y)|
				\lesssim\frac{1}{\max\bigl\{d(x,y),\Xi_{\mathtt w}(x),\Xi_{\mathtt w}(y)\bigr\}^n}
				\lesssim\frac{1}{\bigl(r_0+d(x,y)+\Xi_{\mathtt w}(x)+\Xi_{\mathtt w}(y)\bigr)^n}
				\eqqcolon H(x,y).
			\end{equation}
			
			Let $z\in X$, $\mathfrak m\in\{\mathfrak m_1,\mathfrak m_2\}$, and $r_z\coloneqq r_0+\Xi_{\mathtt w}(z)$. Since $H(x,z)\leq(r_z+d(x,z))^{-n}$, we have $H(\cdot,z)\leq2^n(2^kr_z)^{-n}$ on $B(z,2^kr_z)\setminus B(z,2^{k-1}r_z)$, for $k\geq1$, and on $B(z,r_z)$. Hence, by \eqref{eq:growth_auxiliary_measures_sec6} and $\Xi_{\mathtt w}(z)\leq r_z$,
			\[
			\int H(x,z)^2\,d\mathfrak m(x)
			\lesssim\sum_{k\geq0}\frac{\mathfrak m(B(z,2^kr_z))}{(2^kr_z)^{2n}}
			\lesssim\sum_{k\geq0}\frac{(2^kr_z+\Xi_{\mathtt w}(z))^n}{(2^kr_z)^{2n}}
			\lesssim r_z^{-n}\leq r_0^{-n}.
			\]
			Since $H$ is symmetric, for $y\in X$ we obtain
			\begin{align}
				\int H(x,y)^2\,d\mathfrak m_2(x) +	\int H(y,x)^2\,d\mathfrak m_1(x)&\lesssim r_0^{-n}.\label{eq:integral_square_difference_sec6}
			\end{align}
			By \eqref{eq:item_4_ref} and the $n$-growth of $\mu$, for $i=1,2$ we have \[\mathfrak m_i(X)\leq\mathfrak C_1\mu(B_i)\leq\mathfrak C_1\mathfrak C_0\,r(B_i)^n.\] Integrating \eqref{eq:integral_square_difference_sec6}, we get
			\[
			\iint H(x,y)^2\,d\mathfrak m_1(y)\,d\mathfrak m_2(x)
			\lesssim r_0^{-n}\min\{\mathfrak m_1(X),\mathfrak m_2(X)\}\lesssim1.
			\]
			By \eqref{eq:positive_bound_auxiliary_difference_sec6} and the Cauchy-Schwarz inequality, the integral operator with kernel $\mathfrak D_2$ is bounded from $L^2(\mathfrak m_1)$ to $L^2(\mathfrak m_2)$, with norm $\lesssim1$. Together with the bound for $\mathfrak D_1$, this proves \eqref{eq:bound_auxiliary_difference_sec6}.
		\end{proof}
		\vv
		
		\subsection{Comparison of the norms} As an application of Lemma \ref{lemma:l51comparison} and the properties of the measures discussed at the beginning of the section, we obtain the following comparison result for the $L^2$-operator norms with different measures.
		\begin{lemma}\label{lemma:properties_regularized_cubes}
			There exists a constant $\hat{\mathbf c}_0>0$, depending only on $n,D,\gamma,\varepsilon_0,C,\mathfrak C_0,\mathfrak C_1, \mathfrak C_3,\beta$, and the fixed David-Mattila parameters, such that, for every $\varepsilon>0$,
			\begin{equation}\label{eq:main_bd_lemma_widetilde_mu_0}
				\|T_{\Xi_{\mathtt w},\mu_1,\varepsilon}\|_
				{L^2(\mu_1)\to L^2(|\widehat\nu_2|)}
				+\|T_{\Xi_{\mathtt w},|\widehat\nu_1|,\varepsilon}\|_
				{L^2(|\widehat\nu_1|)\to L^2(\mu_2)} \leq\hat{\mathbf c}_0\,
				\Bigl(1+\|T_{\Xi_{\mathtt w},|\widehat\nu_1|,\varepsilon}\|_
				{L^2(|\widehat\nu_1|)\to L^2(|\widehat\nu_2|)}\Bigr).
			\end{equation}
			Moreover, there exists a constant $\hat{\mathbf c}>0$, with the same dependence, such that
			\begin{equation}\label{eq:main_bd_lemma_widetilde_mu_ab}
				\|T_{\Xi_{\mathtt w},\mu_1,\varepsilon}\|_
				{L^2(\mu_1)\to L^2(\mu_2)}
				\leq\hat{\mathbf c}\,
				\Bigl(1+\|T_{\Xi_{\mathtt w},|\widehat\nu_1|,\varepsilon}\|_
				{L^2(|\widehat\nu_1|)\to L^2(|\widehat\nu_2|)}\Bigr).
			\end{equation}
		\end{lemma}
		\begin{proof}
			In this proof, we use the shorthand notation $\|\cdot\|_{\mathfrak m\to\mathfrak m'}\coloneqq\|\cdot\|_{L^2(\mathfrak m)\to L^2(\mathfrak m')}$ and
			note that $\|T_{\Xi_{\mathtt w},|\widehat\nu_1|,\varepsilon}\|_{|\widehat\nu_1|\to|\widehat\nu_2|}<\infty$, since the kernel of $T_{\Xi_{\mathtt w},|\widehat\nu_1|,\varepsilon}$ is bounded by $C\varepsilon^{-n}$ and the measures $|\widehat\nu_1|$ and $|\widehat\nu_2|$ are finite.

			Let $\mathfrak m_1\in\{\mu,|\nu_1|\}$ and $\mathfrak m_2\in\{\mu,|\nu_2|\}$, so that $\mathfrak m_i|_{B_i}\in\{\mu_i,|\widehat\nu_i|\}$. Since the kernel of the localized operator is such that $L_\varepsilon(x,y)=0$ unless $(x,y)\in B_2\times B_1$, we obtain  $\|S_{\mathfrak m_1,\varepsilon}\|_{\mathfrak m_1\to\mathfrak m_2}=\|S_{\mathfrak m_1|_{B_1},\varepsilon}\|_{\mathfrak m_1|_{B_1}\to\mathfrak m_2|_{B_2}}$. Hence, \eqref{eq:bound_auxiliary_difference_sec6} gives
			\begin{equation}\label{eq:bound_auxiliary_operator_sec6}
				\Bigl|\|T_{\Xi_{\mathtt w},\mathfrak m_1|_{B_1},\varepsilon}\|_{\mathfrak m_1|_{B_1}\to\mathfrak m_2|_{B_2}}-\|S_{\mathfrak m_1,\varepsilon}\|_{\mathfrak m_1\to\mathfrak m_2}\Bigr|\leq C.
			\end{equation}
			In particular, 
			\begin{equation}\label{eq:bound_l52}
				\|S_{|\nu_1|,\varepsilon}\|_{|\nu_1|\to|\nu_2|}\leq \|T_{\Xi_{\mathtt w},|\widehat\nu_1|,\varepsilon}\|_{|\widehat\nu_1|\to|\widehat\nu_2|}+C<\infty.
			\end{equation}
			Thus, by \eqref{eq:bound_auxiliary_operator_sec6} and \eqref{eq:bound_l52} with $(\mathfrak m_1,\mathfrak m_2)$ equal to $(|\nu_1|,\mu)$, $(\mu,|\nu_2|)$, and $(\mu,\mu)$, the estimates \eqref{eq:main_bd_lemma_widetilde_mu_0} and \eqref{eq:main_bd_lemma_widetilde_mu_ab} follow from
			\begin{equation}\label{eq:comparison_norms_S_sec5}
				\|S_{|\nu_1|,\varepsilon}\|_{|\nu_1|\to\mu}+\|S_{\mu,\varepsilon}\|_{\mu\to|\nu_2|}+\|S_{\mu,\varepsilon}\|_{\mu\to\mu}
				\lesssim1+\|S_{|\nu_1|,\varepsilon}\|_{|\nu_1|\to|\nu_2|}.
			\end{equation}
			
			As observed after \eqref{eq:size_auxiliary_kernel_sec6}, Lemma \ref{lemma:mixed_cotlar_transfer} applies to $L_\varepsilon$ and to its transpose $L^t_\varepsilon(x,y)\coloneqq L_\varepsilon(y,x)$, with $\Xi=\Psi_\varepsilon$, $\delta_0=\varepsilon$, and any three measures $\mathfrak m,\mathfrak m',\mathfrak m''\in\{\mu,|\nu_1|,|\nu_2|\}$. We denote by $S^t_{\mathfrak m,\varepsilon}$ the integral operator with kernel $L^t_\varepsilon$ with respect to $\mathfrak m$ so that, by duality, $\|S^t_{\mathfrak m',\varepsilon}\|_{\mathfrak m'\to\mathfrak m}=\|S_{\mathfrak m,\varepsilon}\|_{\mathfrak m\to\mathfrak m'}$. Since $L_\varepsilon(x,y)=0$ if $d(x,y)\leq\varepsilon$, the truncations of $S_{\mathfrak m,\varepsilon}$ at levels $\delta<\varepsilon$ coincide with $S_{\mathfrak m,\varepsilon}$. Hence, \eqref{eq:mct_maximal_bound} gives
			\begin{equation}\label{eq:transfer_output_measure_sec5}
				\|S_{\mathfrak m,\varepsilon}\|_{\mathfrak m\to\mathfrak m''}\lesssim1+\|S_{\mathfrak m,\varepsilon}\|_{\mathfrak m\to\mathfrak m'}
			\end{equation}
			whenever the right-hand side is finite, and the same holds for $S^t_{\mathfrak m,\varepsilon}$. Applying \eqref{eq:transfer_output_measure_sec5} three times, we obtain
			\[
			\begin{split}
				\|S_{|\nu_1|,\varepsilon}\|_{|\nu_1|\to\mu}
				&\lesssim1+\|S_{|\nu_1|,\varepsilon}\|_{|\nu_1|\to|\nu_2|},\\
				\|S_{\mu,\varepsilon}\|_{\mu\to|\nu_2|}=\|S^t_{|\nu_2|,\varepsilon}\|_{|\nu_2|\to\mu}
				&\lesssim1+\|S^t_{|\nu_2|,\varepsilon}\|_{|\nu_2|\to|\nu_1|}=1+\|S_{|\nu_1|,\varepsilon}\|_{|\nu_1|\to|\nu_2|},
			\end{split}
			\]
			and 
			\[
			\|S_{\mu,\varepsilon}\|_{\mu\to\mu}\lesssim1+\|S_{\mu,\varepsilon}\|_{\mu\to|\nu_2|}.
			\]
			By \eqref{eq:bound_l52} and the second inequality, all the right-hand sides are finite, and the inequalities above prove \eqref{eq:comparison_norms_S_sec5}.
		\end{proof}
		
		Let $R$ be the smaller of the cells $R_1$ and $R_2$, and set $M\coloneqq2M_*$. Since either $R_1\subset R_2$ or $R_2\subset R_1$, \eqref{eq:nesting_inclusion_testing} and \eqref{eq:inclusions_84_100_B_B_tilde} give
		\begin{equation}\label{eq:inclusions_enlarged_balls_sec5}
			84B(R)\subset M\widetilde B(R)\subset MB_1\cap MB_2.
		\end{equation}
		\begin{lemma}\label{lemma:final_global_bound}
			For every $\varepsilon>0$, we have
			\begin{equation}\label{eq:enalrged bound small}
				\|T_{\Xi_{\mathtt w},\mu|_{84B(R)},\varepsilon}\|_{L^2(\mu|_{84B(R)})\to L^2(\mu|_{84B(R)})}
				\lesssim1+\|T_{\Xi_{\mathtt w},|\widehat\nu_1|,\varepsilon}\|_{L^2(|\widehat\nu_1|)\to L^2(|\widehat\nu_2|)}.
			\end{equation}
			The same estimate holds with $\Xi_{\mathtt w}$ replaced, on both sides, by the suppression $\widetilde\Phi_R$ from \eqref{eq:tilde Phi}.
		\end{lemma}
		\begin{proof}
			We use the notation of the proof of Lemma \ref{lemma:properties_regularized_cubes}, and we set $\widehat\mu_i\coloneqq\mu|_{MB_i}$, for $i=1,2$, and $\widehat\mu\coloneqq\mu|_{M\widetilde B(R)}$. Since $\mu|_{84B(R)}\leq\widehat\mu$, by \eqref{eq:inclusions_enlarged_balls_sec5}, it suffices to prove \eqref{eq:enalrged bound small} with $\widehat\mu$ in place of $\mu|_{84B(R)}$.
			
			Let $i=1,2$. Since $\chi_i=0$ on $X\setminus B_i$, \eqref{eq:collar_bound_sec6} gives $\Xi_{\mathtt w}\geq\theta r(B_i)/2$ on $X\setminus B_i$. Thus, by the standard kernel bounds of Lemma \ref{lemma:kernel_suppressed}, $|K_{\Xi_{\mathtt w}}(x,y)|\lesssim r(B_i)^{-n}$ if $x\notin B_i$ or $y\notin B_i$. Since $\mu(MB_i)\leq\mathfrak C_0M^nr(B_i)^n$, it follows that
			\[
			\iint_{MB_i\times MB_i}\bigl(\chi_{X\setminus B_i}(x)+\chi_{X\setminus B_i}(y)\bigr)|K_{\Xi_{\mathtt w}}(x,y)|^2\,d\mu(y)\,d\mu(x)
			\lesssim r(B_i)^{-2n}\mu(MB_i)^2\lesssim1.
			\]
			Hence, by the Cauchy-Schwarz inequality, for every $\varepsilon>0$ we have
			\begin{equation}\label{eq:L2_1}
				\bigl\|\chi_{X\setminus B_i}T_{\Xi_{\mathtt w},\widehat\mu_i,\varepsilon}\bigr\|_{\widehat\mu_i\to\widehat\mu_i}
				+\bigl\|T_{\Xi_{\mathtt w},\widehat\mu_i,\varepsilon}(\chi_{X\setminus B_i}\,\cdot\,)\bigr\|_{\widehat\mu_i\to\widehat\mu_i}\lesssim1.
			\end{equation}
			
			Let $f\in L^2(\widehat\mu)$, extended by zero outside $M\widetilde B(R)$. Since $M\widetilde B(R)\subset MB_1\cap MB_2$, we can split
			\[
			T_{\Xi_{\mathtt w},\widehat\mu,\varepsilon}f
			=\chi_{B_2}T_{\Xi_{\mathtt w},\mu_1,\varepsilon}f
			+\chi_{B_2}T_{\Xi_{\mathtt w},\widehat\mu_1,\varepsilon}(\chi_{X\setminus B_1}f)
			+\chi_{X\setminus B_2}T_{\Xi_{\mathtt w},\widehat\mu_2,\varepsilon}f
			\eqqcolon\mathfrak I_1+\mathfrak I_2+\mathfrak I_3.
			\]
			Since $\widehat\mu|_{B_2}\leq\mu_2$ and $\|f\|_{L^2(\mu_1)}\leq\|f\|_{L^2(\widehat\mu)}$, \eqref{eq:main_bd_lemma_widetilde_mu_ab} gives
			\[
			\|\mathfrak I_1\|_{L^2(\widehat\mu)}\leq\|T_{\Xi_{\mathtt w},\mu_1,\varepsilon}\|_{\mu_1\to\mu_2}\|f\|_{L^2(\mu_1)}\lesssim \bigl(1+\|T_{\Xi_{\mathtt w},|\widehat\nu_1|,\varepsilon}\|_{L^2(|\widehat\nu_1|)\to L^2(|\widehat\nu_2|)}\bigr)\,\|f\|_{L^2(\widehat\mu)}.
			\]
			Moreover $\widehat\mu\leq\widehat\mu_i$ and $\|f\|_{L^2(\widehat\mu_i)}=\|f\|_{L^2(\widehat\mu)}$, for $i=1,2$, so \eqref{eq:L2_1} gives
			\[
			\|\mathfrak I_2\|_{L^2(\widehat\mu)}+\|\mathfrak I_3\|_{L^2(\widehat\mu)}
			\leq\bigl\|T_{\Xi_{\mathtt w},\widehat\mu_1,\varepsilon}(\chi_{X\setminus B_1}f)\bigr\|_{L^2(\widehat\mu_1)}
			+\bigl\|\chi_{X\setminus B_2}T_{\Xi_{\mathtt w},\widehat\mu_2,\varepsilon}f\bigr\|_{L^2(\widehat\mu_2)}
			\lesssim\|f\|_{L^2(\widehat\mu)},
			\]
			which proves \eqref{eq:enalrged bound small} for the SIO associated with the suppression $\Xi_{\mathtt w}$.
			
			Finally, let $i=1,2$ and $w_i\in\Omega$. The set $F_i(w_i)$ is contained in $(1-\theta)B_i$ and does not intersect the balls $B(x,r_{j,i}(x))$, for $x\in B_i$ and $j=1,2$, and $B(x,\mathcal R_i(x))$, for $x\in X$. Hence, $\Phi_{i,R_i,w_i}$ is bounded from below by $r_{1,i}$ and $r_{2,i}$ on $B_i$, and by $\mathcal R_i$ and $\dist(\cdot,(1-\theta)B_i)$ on $X$. These lower bounds do not depend on $w_i$, so, by \eqref{eq:definition_suppressions_2a}, they also hold for $\Phi_{i,R_i}$, and thus for $\widetilde\Phi_R=\max\{\Phi_{1,R_1},\Phi_{2,R_2}\}$. Since $\widetilde\Phi_R$ is also $1$-Lipschitz, and these are the only properties of $\Xi_{\mathtt w}$ that we used in this section, all the estimates of this section hold with $\widetilde\Phi_R$ in place of $\Xi_{\mathtt w}$.
		\end{proof}

	}

	\section{The Key Lemma for the suppressed operator}\label{section:suppressed_mixed_good}
	
	\subsubsection*{The setup}
	We keep the setting and the notation of Section \ref{section:L2bd}. In particular, $R_1\in\Top''_1$ and $R_2\in\Top''_2$ are such that $\widetilde{\Tree}(R_1,R_2)\neq\varnothing$, $B_i=\widetilde B(R_i)$, and $\widehat\nu_1$ and $\widehat\nu_2$ are the testing measures for $T$ on $B_1$ and for $T^*$ on $B_2$, respectively. As observed in Section \ref{sec:corona}, for $i=1,2$ the measures $\mu$ and $\widehat\nu_i$ satisfy assumptions \ref{ex_item_1}-\ref{ex_item_6} with $B=B_i$ and $\mathfrak C_2=1$, where $T$ is replaced by $T^*$ if $i=2$. We refer to \eqref{eq:hp_6} for $T$ on $B_1$ and for $T^*$ on $B_2$ as the maximal truncation assumptions.
	
	We use $\Omega$ as in Section \ref{section:preliminaries_and_notation}, $\Omega^2=\Omega\times\Omega$, and the notation $\mathtt w=(w_1, w_2)\in \Omega^2$. For $w_i\in\Omega$, we consider the two independent dyadic grids $\bbD(w_i)$ and the total exceptional sets $\mathscr W_i(w_i)$. 
	
	Let $i=1,2$. By \eqref{eq:measure_testing_size}, we have $\widehat\nu_i(B_i)=\mu(B_i)$. We write $\widehat\nu_i=\widehat b_i\,|\widehat\nu_i|$ on $B_i$, and we extend $\widehat b_i$ by $1$ on $X\setminus B_i$. Note that, since $\nu_i=\widehat\nu_i+\mu|_{X\setminus B_i}$, we have $\nu_i=\widehat b_i\,|\nu_i|$ and $|\widehat b_i|=1$ $|\nu_i|$-almost everywhere.
	\vv
	
	\textbf{By symmetry, we may assume that $R_1\subset R_2$.} Otherwise, it suffices to interchange the indices and to replace $T$ by $T^*$. Then \eqref{eq:nesting_inclusion_testing_0} and \eqref{eq:nesting_inclusion_testing} give $r(B_1)\leq r(B_2)$ and $2M_*B_1\subset2M_*B_2$.
	
	\vv
	
	Let $\kappa\in(0,1)$, to be fixed in Section \ref{section:probabilisti_argument_subsection}. We fix the generations $k_1$ and $k_2$ given by Lemma \ref{lem:random_capture_goodness}, applied to $B_1$, $B_2$, and $\kappa$. For $i,j=1,2$ and $w_j\in\Omega$, we denote by $Q_i(w_j)$ the cube of $\bbD_{k_i}(w_j)$ that contains the center of $B_i$. Note that $Q_i(w_j)$ depends only on $w_j$.
	
	\vv

	As is customary in $Tb$ theorems, see \cite{NTrV03}, we adopt the notion of good and bad cubes, depending on the exponent $\alpha\coloneqq\gamma/(2(n+\gamma))\in(0,1)$. We denote by $\eta>0$ the boundary exponent in Theorem \ref{thm:HK56}, and set $\mathfrak d=\log_2(D)$, where $D$ is the geometrically doubling constant of $X$. We also recall that $\el\in(0,1)$ is the scaling factor in Theorem \ref{thm:HK22}.

	\begin{definition}[Good and bad cubes]\label{def:good_bad_cubes}
		Let $\mathtt w=(w_1,w_2)$ and $m\in\mathbb N$. A cube $Q\in\bbD(w_1)$ is \textit{bad} if there exists $R\in\bbD(w_2)$ with $\ell(R)\geq\el^{-m}\ell(Q)$ such that
		\[
		\operatorname{dist}(Q,R)\leq\ell(Q)^\alpha\ell(R)^{1-\alpha}
		\quad\text{and}\quad
		\operatorname{dist}(Q,X\setminus R)	\leq\ell(Q)^\alpha\ell(R)^{1-\alpha}.
		\]
		Otherwise $Q$ is \textit{good}. We denote these families by	$\bbD_1^{\mathcal B}(\mathtt w,m)$ and $\bbD_1^{\mathcal G}(\mathtt w,m)$, respectively,  omitting $m$ when it is fixed. The corresponding families $\bbD_2^{\mathcal B}(\mathtt w,m)$ and $\bbD_2^{\mathcal G}(\mathtt w,m)$ in the second grid are defined by interchanging the indices.
	\end{definition}

	By Lemma \ref{lem:random_capture_goodness}, for $0<\kappa\ll 1$ and $m$ sufficiently large, to be fixed in Section \ref{section:probabilisti_argument_subsection}, there exists $L=L(\kappa)>0$, depending also on the other structural constants in Lemma \ref{lem:random_capture_goodness}, such that, denoting $\mathbb P^{\Omega^2}\coloneqq \mathbb P\otimes \mathbb P$, good cubes occur with large probability. Namely, we have
	\begin{equation}\label{eq:large_probability_top_cubes}
		\begin{split}
			\mathbb P^{\Omega^2}\Bigl(&\Bigl\{\mathtt w=(w_1,w_2)\in \Omega^2:\text{ for every }i=1,2,\,j=1,2\,
			\text{there exists } Q_i(w_j)\in\bbD(w_j),\\
			&\qquad\qquad\text{such that }6M_*B_i\subset Q_i(w_j),\, Q_1(w_j) \subset Q_2(w_j) \text{ for }j=1,2,\\
			&\qquad\qquad\mathrm{diam} Q_i(w_j)\leq L\, r(B_i),\text{ and }\,
			Q_i(w_j)\in\bbD_j^{\mathcal G}(\mathtt w,m)
			\Bigr\}\Bigr)>1-\kappa.
		\end{split}
	\end{equation}

	Hence, we denote by $\mathscr G\subset \Omega^2$ the event in parentheses in \eqref{eq:large_probability_top_cubes}.
	For the rest of the argument, we hence fix $(w_1,w_2)\in \mathscr G$ and consider its associated cubes $Q_i(w_j)$, $i,j=1,2$. In the following subsections, all probabilistic arguments are carried out by considering the conditional probability of the event belonging to $\mathscr G$, namely
	\[
	\mathbb P^{\Omega^2}_{\mathscr G}(A)\coloneqq \frac{\mathbb P^{\Omega^2}(A\cap \mathscr G)}{\mathbb P^{\Omega^2}(\mathscr G)}, \qquad \text{ for }A\subset \Omega^2,
	\]
	and its associated expectations.

	\vv
	
	{
		{
			\subsubsection*{Adapted martingale differences}
			Let $\mathtt w\in\mathscr G$ and $i=1,2$. We use $Q_2(w_i)$ as the top cube of the martingale decomposition associated with $\bbD(w_i)$. By the definition of $\mathscr G$, we have $B_1\subset Q_1(w_i)\subset Q_2(w_i)$ and $B_2\subset Q_2(w_i)$. We extend functions in $L^2(|\widehat\nu_i|)$ by zero on $X\setminus B_i$, so that their norms in $L^2(|\widehat\nu_i|)$ and $L^2(|\nu_i|)$ coincide. We remark that the martingale sums below run over the cubes of $\bbD(w_i)$ contained in $Q_2(w_i)$, whereas goodness is defined in terms of the full grids.
			\begin{definition}[``Terminal'' and ``transit'' cubes]
				Let $i=1,2$. A cube $Q\in\bbD(w_i)$ is \textit{terminal} if $Q\subset\mathscr W_i(w_i)$, and \textit{transit} otherwise. We denote these families by $\bbD_i^{\term}$ and $\bbD_i^{\tr}$, respectively.
			\end{definition}
			Since $B_i\subset Q_2(w_i)$, Lemma \ref{lemma:lemma38} applies with $Q_0=Q_2(w_i)$, $B=B_i$, and $\nu=\nu_i$. The cubes of $\Stop(Q_2(w_i),w_i)$ are contained in $\mathscr W_i(w_i)$, so a transit cube $Q\subset Q_2(w_i)$ is not contained in any of them. Moreover, $Q_2(w_i)$ is itself a transit cube, since $\varnothing\neq F_i(w_i)\subset B_i\setminus\mathscr W_i(w_i)\subset Q_2(w_i)\setminus\mathscr W_i(w_i)$. Hence, Lemma \ref{lemma:lemma38} gives
			\[
			\widehat{\mathbf m}^{-1}\mu(Q)\leq|\nu_i|(Q)\leq\widehat{\mathbf m}\,\mu(Q)
			\qquad\text{and}\qquad
			|\nu_i(Q)|\geq\frac{|\nu_i|(Q)}{2\mathfrak C_1},
			\qquad Q\in\bbD_i^{\tr},\quad Q\subset Q_2(w_i).
			\]
			In particular, since $\nu_i=\widehat b_i\,|\nu_i|$, we have $|\langle\widehat b_i\rangle_{|\nu_i|,Q}|\geq(2\mathfrak C_1)^{-1}$ whenever $|\nu_i|(Q)>0$. 
			
			\begin{remark}\label{remark:growth_children_transit_cubes}
				Let $Q\in\bbD_i^{\tr}$, let $P$ be a dyadic child of $Q$, and set $C\coloneqq2C_1\el^{-1}$. Then
				\[
				|\nu_i|(B(z,t))\lesssim t^n,\qquad \text{ for }  z\in P,\quad t\geq\ell(P).
				\]
				Indeed, let $x\in Q\setminus\mathscr W_i(w_i)$. Since $d(x,z)\leq\diam Q\leq2C_1\ell(Q)=C\ell(P)\leq Ct$, we have $B(z,t)\subset B(x,(C+1)t)$. Moreover, $x\notin\mathscr H_{i,R_i}$, because $\mathscr H_{i,R_i}\subset\mathscr W_i(w_i)$. Thus, \eqref{eq:definition_r_i} and \eqref{eq:def_non_ahlfors} give $|\widehat\nu_i|(B(x,r))\leq\mathbf m\,r^n$ for every $r>0$. Hence,
				\[
				|\nu_i|(B(z,t))\leq|\widehat\nu_i|\bigl(B(x,(C+1)t)\bigr)+\mu(B(z,t))\leq\bigl(\mathbf m\,(C+1)^n+\mathfrak C_0\bigr)\,t^n.
				\]
			\end{remark}
			\vv
			
			Let $\mathfrak m\in M_+(X)$.
			For $f\in L^1_{\loc}(\mathfrak m)$ and a cube $Q$ in $X$ with $\mathfrak m(Q)\neq 0$, we write
			\begin{equation}\label{eq:def_Xi_f}
				\langle f\rangle_{\mathfrak m, Q} \coloneqq \frac{1}{\mathfrak m(Q)}\int_Q f\, d\mathfrak m.
			\end{equation}
			
			We also note that 
			\[
			\int_{Q_2(w_i)}\widehat b_i\,d|\nu_i|=\mu(Q_2(w_i)),\qquad
			\mu(Q_2(w_i))\leq|\nu_i|(Q_2(w_i))
			\leq\mathfrak C_1\mu(Q_2(w_i)),
			\]
			and define the \textit{top averages} by
			\begin{equation}\label{eq:top_averages_def}
				E_i f\coloneqq\chi_{Q_2(w_i)}\widehat b_i\,
				\frac{\int_{Q_2(w_i)}f\,d|\nu_i|}{\mu(Q_2(w_i))}\qquad \text{ and }\qquad 
				E_i^*f\coloneqq\chi_{Q_2(w_i)}\,
				\frac{\int_{Q_2(w_i)}f\widehat b_i\,d|\nu_i|}{\mu(Q_2(w_i))}.
			\end{equation}
			For $i=1,2$ and $A\subset X$ with $|\nu_i|(A)>0$, we often use the shorthand notation $\langle f\rangle_{i,A}\coloneqq \langle f\rangle_{|\nu_i|,A}$.
			
			For $Q\in \bbD_i$, we set $\Ch(Q)=\varnothing$ if $|\nu_i|(Q)=0$ and, if $|\nu_i|(Q)>0$, we denote by $\Ch(Q)$ the family of dyadic children $P$ of $Q$ such that $|\nu_i|(P)>0$. For $Q\in \bbD^{\tr}_i$ and $f\in L^1_{\loc}(|\nu_i|)$ we define
			\begin{equation}\label{eq:def_martingales}
				\Delta_{i,Q}f(x)\coloneqq 
				\begin{cases}
					0, &\qquad\text{ for } x\in X\setminus \bigcup_{P\in \Ch(Q)} P,  \\
					\displaystyle\biggl(\frac{\langle f\rangle_{i,P}}{\langle \widehat b_i\rangle_{i,P}}-\frac{\langle f\rangle_{i, Q}}{\langle \widehat b_i\rangle_{i, Q}}\biggr)\,\widehat b_i(x), &\qquad\text{ for } x\in P, \text{ for some } P \in \Ch(Q)\cap \bbD^{\tr}_i,\\
					\displaystyle f(x)-\frac{\langle f\rangle_{i, Q}}{\langle \widehat b_i\rangle_{i,Q}}\, \widehat b_i(x),&\qquad  \text{ for } x\in P, \text{ for some } P \in \Ch(Q)\cap \bbD^{\term}_i.
				\end{cases}
			\end{equation}
			We also set $\Delta_{i,Q}\coloneqq0$ if $Q$ is terminal or $Q\not\subset Q_2(w_i)$.
			
			A standard computation, see \cite[Lemma 5.10]{To14} for a close analogue, shows that
			\begin{equation}\label{eq:property_id_mart_diff}
				\Delta_{i,Q}\, \Delta_{i,R}=
				\begin{cases}
					\Delta_{i,Q}, \qquad &\text{ if }R=Q,\\
					0, \qquad &\text{ if }R\in \bbD^{\tr}_i, R\neq Q,
				\end{cases}
			\end{equation}
			and every $f\in L^2(|\nu_i|)$ supported in $Q_2(w_i)$ can be decomposed as	
			\begin{equation}\label{eq:dec_f_mart}
				f= E_i f + \sum_{Q\in \bbD^{\tr}_i} \Delta_{i,Q} f,
			\end{equation}
			where the sum is unconditionally convergent in $L^2(|\nu_i|)$. Furthermore, see \cite[Lemma 5.11]{To14}, there exists a constant $\tilde c=\tilde c(\mathfrak C_1)>0$ such that we have the martingale square-function bound
			\begin{equation}\label{eq:dec_f_mart_2}
				\tilde c^{-1}\, \|f\|^2_{L^2(|\nu_i|)}\leq \|E_{i} f\|^2_{L^2(|\nu_i|)} + \sum_{Q\in \bbD^{\tr}_i}\|\Delta_{i,Q} f\|^2_{L^2(|\nu_i|)}\leq \tilde c \, \|f\|^2_{L^2(|\nu_i|)}.
			\end{equation}
			
			Here we use $|\widehat b_i|=1$ $|\nu_i|$-almost everywhere and the uniform lower bound for $|\langle\widehat b_i\rangle_{i,Q}|$ on transit cubes. Thus $\tilde c$ depends only on the common accretivity constant
			$(2\mathfrak C_1)^{-1}$.
			It is also possible to prove that the adjoint of $\Delta_{i,Q}$ is given by
			\[
			\Delta^*_{i,Q}f(x)\coloneqq 
			\begin{cases}
				0, &\qquad\text{ for } x\in X\setminus \bigcup_{P\in \Ch(Q)} P,  \\
				\displaystyle\frac{\langle f\,\widehat b_i\rangle_{i,P}}{\langle \widehat b_i\rangle_{i, P}}-\frac{\langle f\, \widehat b_i\rangle_{i,Q}}{\langle \widehat b_i\rangle_{i,Q}}, &\qquad\text{ for } x\in P, \text{ for some } P \in \Ch(Q)\cap \bbD^{\tr}_i,\\
				\displaystyle f(x)-\frac{\langle f\,\widehat  b_i\rangle_{i,Q}}{\langle \widehat b_i\rangle_{i,Q}},&\qquad  \text{ for } x\in P, \text{ for some } P \in \Ch(Q)\cap \bbD^{\term}_i,
			\end{cases}
			\]
			which is such that
			\[
			\Delta^*_{i,Q}\, \Delta^*_{i,R}=
			\begin{cases}
				\Delta^*_{i,Q}, \qquad &\text{ if }R=Q,\\
				0, \qquad &\text{ if }R\in \bbD^{\tr}_i, R\neq Q.
			\end{cases}
			\]
			
			Similarly to \eqref{eq:dec_f_mart}, a function $g\in L^2(|\nu_i|)$ supported in $Q_2(w_i)$ can be decomposed as
			\begin{equation}\label{eq:dec_f_mart_adj}
				g= E^*_i g + \sum_{Q\in \bbD^{\tr}_i} \Delta^*_{i,Q} g,
			\end{equation}
			where the sum is unconditionally convergent in $L^2(|\nu_i|)$ and a bound analogous to \eqref{eq:dec_f_mart_2} holds.
			
			\vv
		}
	}
	
	\subsubsection*{Projection operators}
	Next, we define the projection operators associated with good and bad cubes.
	
	\begin{definition}\label{def:p_good_p_bad_not_truncated}
		Let $\mathtt w = (w_1, w_2)\in \Omega^2$, $f\in L^2(|\nu_1|)$, and $g\in L^2(|\nu_2|)$. We define the operators
		\[
		\mathcal P^1_{\mathcal G}(\mathtt w)f\coloneqq \sum_{Q\in \bbD^{\mathcal G}_{1}(\mathtt w)\cap\bbD_1^{\tr} }\Delta_{1,Q}f, \qquad \qquad \mathcal P^1_{\mathcal B}(\mathtt w)f\coloneqq \sum_{Q\in \bbD^{\mathcal B}_{1}(\mathtt w)\cap\bbD_1^{\tr} } \Delta_{1,Q}f,
		\]
		and 
		\[
		\mathcal P^2_{\mathcal G}(\mathtt w)g\coloneqq \sum_{R\in \bbD^{\mathcal G}_{2}(\mathtt w)\cap\bbD_2^{\tr}}\Delta_{2,R}g, \qquad \qquad \mathcal P^2_{\mathcal B}(\mathtt w)g\coloneqq \sum_{R\in \bbD^{\mathcal B}_{2}(\mathtt w)\cap\bbD_2^{\tr} }\Delta_{2,R}g.
		\]
	\end{definition}

	Given the probability space $(\Omega, \mathbb P^\Omega)$, we denote by $\mathbb E_w\coloneqq \mathbb E^\Omega_w$ the expectation over $w\in \Omega$.
	We also denote by $\widetilde{\mathbb E}_{\mathtt w}$ the conditional expectation over $(w_1,w_2)\in \Omega\times\Omega$ with respect to $\mathbb P^{\Omega^2}_{\mathscr G}$.
	
	{
		\begin{lemma}\label{lemma:probability_bad_1}
			There is a quantity $\varepsilon_{\mathcal B}(m)\to0$ as $m\to\infty$
			such that, for fixed $w_1$,
			\[
			\mathbb E_{w_2}\biggl\|
			\sum_{Q\in\bbD^{\mathcal B}_1(w_1,w_2)\cap\bbD_1^{\tr}}
			\Delta_{1,Q}f\biggr\|_{L^2(|\nu_1|)}
			\leq\varepsilon_{\mathcal B}(m)\|f\|_{L^2(|\nu_1|)}.
			\]
			The analogous estimate holds with the indices interchanged.
		\end{lemma}
		\begin{proof}
			Fix $w_1\in\Omega$. The cube $Q_2(w_1)$, the families $\bbD_1^{\tr}$ and $\bbD_1^{\term}$, and the martingale differences $\Delta_{1,Q}$ depend only on $w_1$. If $B_1\cup B_2\not\subset Q_2(w_1)$, we set $\Delta_{1,Q}\coloneqq0$ for every $Q$. In this way, the operators $\Delta_{1,Q}$ are defined for every $w_1\in\Omega$, agree with \eqref{eq:def_martingales} on $\mathscr G$, and satisfy the upper bound for the martingale differences in \eqref{eq:dec_f_mart_2} uniformly in $w_1$. Only the family of bad cubes depends on $w_2$.
			
			The argument leading to \eqref{eq:probability_estimate_aux} applies to every $Q\in\bbD(w_1)$, and gives
			\[
			\mathbb P_{w_2}(\{Q\text{ is bad}\})
			\leq\frac{C_2(1+2C_1)^\eta}{1-\ell^{\alpha\eta}}\,\ell^{\alpha\eta m}
			\eqqcolon p_{\mathcal B}(m).
			\]
			Moreover, $\int\Delta_{1,Q}f\,d|\nu_1|=0$, so that $E_1\Delta_{1,Q}f=0$. Hence, by \eqref{eq:property_id_mart_diff} and the lower bound in \eqref{eq:dec_f_mart_2}, we have 
			\[
			\bigl\|\sum_{Q\in\mathcal Q}\Delta_{1,Q}f\bigr\|^2_{L^2(|\nu_1|)}\leq\tilde c\sum_{Q\in\mathcal Q}\|\Delta_{1,Q}f\|^2_{L^2(|\nu_1|)},\qquad \text{ for every } \mathcal Q\subset\bbD_1^{\tr}. 
			\]
			If $B_1\cup B_2\not\subset Q_2(w_1)$, both sides vanish.
			Applying this with $\mathcal Q=\bbD^{\mathcal B}_1(w_1,w_2)\cap\bbD_1^{\tr}$, Tonelli's theorem and the upper bound in \eqref{eq:dec_f_mart_2} give
			\begin{equation*}
				\begin{split}
					\mathbb E_{w_2}\Bigl\|\sum_{Q\in\bbD^{\mathcal B}_1(w_1,w_2)\cap\bbD_1^{\tr}}\Delta_{1,Q}f\Bigr\|^2_{L^2(|\nu_1|)}
					&\leq\tilde c\sum_{Q\in\bbD_1^{\tr}}\mathbb P_{w_2}(\{Q\text{ is bad}\})\,\|\Delta_{1,Q}f\|^2_{L^2(|\nu_1|)}\\
					&\leq\tilde c^{2}\,p_{\mathcal B}(m)\,\|f\|^2_{L^2(|\nu_1|)}.
				\end{split}
			\end{equation*}
			The lemma follows from the Cauchy-Schwarz inequality, with $\varepsilon_{\mathcal B}(m)\coloneqq\tilde c\,p_{\mathcal B}(m)^{1/2}$.
		\end{proof}
	}
	
	\vv
	As an immediate corollary, we get the following bounds on projections and their associated conditional expectations.
	\begin{corollary}\label{lemma:probability_bad_2}
		Let $\mathtt w=(w_1,w_2)\in \Omega^2$, $f\in L^2(|\nu_1|)$, and $g\in L^2(|\nu_2|)$. Then 
		\begin{equation}\label{eq:bound_exp_bad}
			\widetilde {\mathbb E}_{\mathtt w}\|\mathcal P^1_{\mathcal B}(\mathtt w)f\|_{L^2(|\nu_1|)}\leq \, \frac{\varepsilon_{\mathcal B}(m)}{\mathbb P^{\Omega^2}(\mathscr G)} \, \|f\|_{L^2(|\nu_1|)}\quad \text{ and }\,  \widetilde {\mathbb E}_{\mathtt w}\|\mathcal P^2_{\mathcal B}(\mathtt w)g\|_{L^2(|\nu_2|)}\leq \, \frac{\varepsilon_{\mathcal B}(m)}{\mathbb P^{\Omega^2}(\mathscr G)} \, \|g\|_{L^2(|\nu_2|)}.
		\end{equation}
		
		We also have
		\begin{equation}\label{eq:bound_exp_good}
			\|\mathcal P^1_{\mathcal G}f\|_{L^2(|\nu_1|)}\lesssim \|f\|_{L^2(|\nu_1|)}\qquad \text{ and }\qquad \|\mathcal P^2_{\mathcal G}g\|_{L^2(|\nu_2|)}\lesssim \|g\|_{L^2(|\nu_2|)},
		\end{equation}
		where the implicit constants do not depend on $(w_1,w_2)$.
	\end{corollary}
	\begin{proof}
		The bounds \eqref{eq:bound_exp_bad} are direct consequences of Lemma \ref{lemma:probability_bad_1} and the definition of $\mathbb P^{\Omega^2}(\mathscr G)$, recalling that the grids associated with $w_1$ and $w_2$ are independent. The inequalities \eqref{eq:bound_exp_good} follow from \eqref{eq:dec_f_mart_2}.
	\end{proof}
	
	\vv
	
	\subsection{The Key Lemma}
	
	Our goal is to bound the norm of $T_{\Xi_{\mathtt w},|\widehat\nu_1|,\varepsilon}$ from $L^2(|\widehat\nu_1|)$ to $L^2(|\widehat\nu_2|)$, uniformly in $\varepsilon$ and $\mathtt w\in\mathscr G$. As in \cite{NTrV03}, in Section \ref{sec:main_theorem_antisym} we split $f$ and $g$ into their good and bad parts. The bad parts are controlled, on average over $\mathtt w$, by Corollary \ref{lemma:probability_bad_2}. The good parts are controlled by the Key Lemma below, which bounds their averaged pairing by a constant plus a small multiple of the norm of the operator, which can then be absorbed.
	
	We write
	\[
	\|T_{\Xi_{\mathtt w},\varepsilon}\|_{\op}\coloneqq\|T_{\Xi_{\mathtt w},|\widehat\nu_1|,\varepsilon}\|_{L^2(|\widehat\nu_1|)\to L^2(|\widehat\nu_2|)}.
	\]
	Since $|K_{\Xi_{\mathtt w}}(x,y)|\lesssim d(x,y)^{-n}$, by Lemma \ref{lemma:kernel_suppressed}, the Cauchy-Schwarz inequality gives 
	\[
	\|T_{\Xi_{\mathtt w},\varepsilon}\|_{\op}\lesssim\varepsilon^{-n}|\widehat\nu_1|(B_1)^{1/2}|\widehat\nu_2|(B_2)^{1/2},
	\]
	uniformly in $\mathtt w$. Hence,
	\begin{equation}\label{eq:definition_N_epsilon_sec6}
		N_\varepsilon\coloneqq\sup_{\mathtt w\in\mathscr G}\|T_{\Xi_{\mathtt w},\varepsilon}\|_{\op}<\infty.
	\end{equation}
	The main result of this section is the following.
	
	\begin{lemma}[Key Lemma]\label{lem:key_lemma_good_functions}
		Let $0<a<1/4$, $0<\lambda,\upsilon<1$, and fix the goodness parameter
		$m$. Let $N_\varepsilon$ be as in \eqref{eq:definition_N_epsilon_sec6}. For $f\in L^2(|\widehat\nu_1|)$ and
		$g\in L^2(|\widehat\nu_2|)$, extended by zero outside their testing
		balls, we have
		\begin{equation}\label{eq:bound_supp_good_functions}
			\begin{split}
				&\widetilde{\mathbb E}_{\mathtt w}
				\Bigl|\bigl\langle T_{\Xi_{\mathtt w},|\widehat\nu_1|,\varepsilon}
				(\mathcal P^1_{\mathcal G}f),
				\mathcal P^2_{\mathcal G}g\bigr\rangle_{|\widehat\nu_2|}\Bigr|\\
				&\quad\leq\Bigl[C_{m,\lambda,\upsilon}(1+a^{-n})+C_{m,\upsilon}(M_\Lambda a)^{1/2}(1+N_\varepsilon)\\
				&\qquad\qquad\qquad+C_m\Bigl(\frac{(2m+1)\lambda^\eta+\upsilon}{\mathbb P^{\Omega^2}(\mathscr G)}\Bigr)^{1/2}(1+N_\varepsilon)\Bigr]\, \|f\|_{L^2(|\widehat\nu_1|)}\|g\|_{L^2(|\widehat\nu_2|)}.
			\end{split}
		\end{equation}
		Here, $M_\Lambda\geq1$ is the small-boundary constant in \eqref{eq:small_boundary_random_balls}, which is fixed once $m$ and $\upsilon$ are fixed. The constants $C_{m,\lambda,\upsilon}$, $C_{m,\upsilon}$, and $C_m$ depend only on the indicated parameters and the structural constants, and are independent of $a$, $\varepsilon$, and the random parameters.
	\end{lemma}
	
	For fixed $m$, the coefficient of $N_\varepsilon$ in \eqref{eq:bound_supp_good_functions} can be made arbitrarily small by choosing $\lambda,\upsilon$ sufficiently small, and then choosing $a$ sufficiently small. We prove the weak boundedness estimate for $S_{|\nu_1|,\varepsilon}$ in Lemma \ref{lemma:key_supp_WBP_in_text}.
	
	We remark that the value of the constants in the bound \eqref{eq:bound_supp_good_functions} does not depend on the truncation $\varepsilon$.
	The proof of the Key Lemma \ref{lem:key_lemma_good_functions} follows the scheme of the Nazarov-Treil-Volberg proof.
	However, due to the a priori lack of antisymmetry of $T$, the proof requires several adaptations. In particular, before applying the martingale difference decomposition, we prove that the Key Lemma can be equivalently stated for the localized operator $S_{|\nu_1|, \varepsilon}$, for which the testing is also satisfied.
	
	\subsubsection*{Reduction to the localized operator}
	We recall that the kernel $L_\varepsilon$ of the localized operator $S_{|\nu_1|,\varepsilon}$ of Section \ref{section:L2bd}, defined in \eqref{eq:definition_auxiliary_kernel_sec6}, vanishes unless $(x,y)\in B_2\times B_1$ and $d(x,y)>\varepsilon$, and that $L_\varepsilon$ and its transpose satisfy \eqref{eq:regularity_auxiliary_kernel_sec6} and \eqref{eq:size_auxiliary_kernel_sec6}, uniformly in $\varepsilon$ and $\mathtt w$. For a positive measure $\mathfrak m$, we denote by $S^*_{\mathfrak m,\varepsilon}$ the integral operator with kernel $L_\varepsilon(y,x)$ with respect to $\mathfrak m$, so that $\langle S_{|\nu_1|,\varepsilon}f,g\rangle_{|\nu_2|}=\langle f,S^*_{|\nu_2|,\varepsilon}g\rangle_{|\nu_1|}$. We also use the shorthand notation
	\[
	\langle\cdot,\cdot\rangle\coloneqq\langle\cdot,\cdot\rangle_{|\nu_1|},
	\qquad
	\langle\cdot,\cdot\rangle_*\coloneqq\langle\cdot,\cdot\rangle_{|\nu_2|},
	\qquad
	\|S_\varepsilon\|_{\op}\coloneqq\|S_{|\nu_1|,\varepsilon}\|_{L^2(|\nu_1|)\to L^2(|\nu_2|)}.
	\]
	By \eqref{eq:bound_auxiliary_operator_sec6} with $\mathfrak m_i=|\nu_i|$, we have
	\begin{equation}\label{eq:localized_norm_comparison_sec5}
		\bigl|\|S_\varepsilon\|_{\op}-\|T_{\Xi_{\mathtt w},\varepsilon}\|_{\op}\bigr|\lesssim1,
	\end{equation}
	and, in particular, $\|S_\varepsilon\|_{\op}\lesssim1+N_\varepsilon$ for every $\mathtt w\in\mathscr G$. Moreover, since $|\nu_i|=|\widehat\nu_i|$ on $B_i$, Lemma \ref{lemma:l51comparison} with $\mathfrak m_i=|\widehat\nu_i|$ gives, for $f\in L^2(|\nu_1|)$ and $g\in L^2(|\nu_2|)$,
	\begin{equation}\label{eq:localized_pairing_comparison_sec5}
		\Bigl|\bigl\langle T_{\Xi_{\mathtt w},|\widehat\nu_1|,\varepsilon}f,g\bigr\rangle_{|\widehat\nu_2|}
		-\bigl\langle S_{|\nu_1|,\varepsilon}f,g\bigr\rangle_*\Bigr|
		\lesssim\|f\|_{L^2(|\nu_1|)}\|g\|_{L^2(|\nu_2|)}.
	\end{equation}
	By \eqref{eq:bound_exp_good}, it therefore suffices to prove \eqref{eq:bound_supp_good_functions} with $\langle S_{|\nu_1|,\varepsilon}(\mathcal P^1_{\mathcal G}f),\mathcal P^2_{\mathcal G}g\rangle_*$ in place of the pairing on its left-hand side.
	
	\subsubsection*{Testing for the localized operator}
	In the proof of the Key Lemma, $\widehat b_1$ and $\widehat b_2$ play the role of the accretive functions of a $Tb$ theorem. The corresponding testing conditions are the bounds \eqref{eq:localized_maximal_testing_sec5} below for the truncations of $S_{|\nu_1|,\varepsilon}\widehat b_1$ and $S^*_{|\nu_2|,\varepsilon}\widehat b_2$, which we now deduce from the maximal truncation assumptions.
	
	Let $i=1,2$. Since $F_i(w_i)\subset G_i\subset B_i\setminus(\mathscr H_{i,R_i}\cup\mathscr E_{i,R_i})$, the construction of the suppression yields
	\[
	\Xi_{\mathtt w}(x)\geq\dist(x,F_i(w_i))\geq\dist\bigl(x,B_i\setminus(\mathscr H_{i,R_i}\cup\mathscr E_{i,R_i})\bigr),\qquad x\in X,
	\]
	that is, $\Xi_{\mathtt w}$ satisfies \eqref{eq:condition_suppression} with $B=B_i$. Hence, by the maximal truncation assumptions and Lemma \ref{eq:lemm_bd_suppressed_2},
	\[
	T_{\Xi_{\mathtt w},\sharp}\widehat\nu_1\leq c_\sharp\quad\text{on }B_1
	\qquad\text{and}\qquad
	(T^*)_{\Xi_{\mathtt w},\sharp}\widehat\nu_2\leq c_\sharp\quad\text{on }B_2.
	\]
	Moreover, since $\chi_i=0$ on $X\setminus B_i$, \eqref{eq:collar_bound_sec6} gives $\Xi_{\mathtt w}\geq\theta r(B_i)/2$ on $X\setminus B_i$. Thus, \eqref{eq:bound_suppressed_kernel_1}, \eqref{eq:item_4_ref}, and the $n$-growth of $\mu$ imply
	\[
	T_{\Xi_{\mathtt w},\sharp}\widehat\nu_1(x)\leq\int|K_{\Xi_{\mathtt w}}(x,y)|\,d|\widehat\nu_1|(y)
	\lesssim\frac{|\widehat\nu_1|(B_1)}{(\theta r(B_1))^n}\lesssim1,\qquad x\in X\setminus B_1,
	\]
	and, in the same way, $(T^*)_{\Xi_{\mathtt w},\sharp}\widehat\nu_2\lesssim1$ on $X\setminus B_2$. Therefore, we obtain
	\begin{equation}\label{eq:global_testing_before_localization_sec5}
		\|T_{\Xi_{\mathtt w},\sharp}\widehat\nu_1\|_{L^\infty(X)}
		+\|(T^*)_{\Xi_{\mathtt w},\sharp}\widehat\nu_2\|_{L^\infty(X)}\lesssim1.
	\end{equation}
	
	For $\rho>0$, we denote by $S_{|\nu_1|,\varepsilon,\rho}$ and $S^*_{|\nu_2|,\varepsilon,\rho}$ the truncations of $S_{|\nu_1|,\varepsilon}$ and $S^*_{|\nu_2|,\varepsilon}$ at level $\rho$, respectively. 
	
	\begin{lemma}\label{lemma:localized_maximal_testing_sec5}
		For every $x\in X$, we have that
		\begin{equation}\label{eq:localized_maximal_testing_sec5}
			\sup_{\rho>0}\bigl|S_{|\nu_1|,\varepsilon,\rho}\widehat b_1(x)\bigr|
			+\sup_{\rho>0}\bigl|S^*_{|\nu_2|,\varepsilon,\rho}\widehat b_2(x)\bigr|\lesssim 1,
		\end{equation}
		where the implicit constant depends only on the structural parameters.
	\end{lemma}
	\begin{proof}
		Let $x\in X$ and $\rho>0$, and let $a_\varepsilon(t)=\varphi(t/\varepsilon)\psi(t/r_0)$ be as in Section \ref{section:L2bd}. Since $\widehat b_1|\nu_1|=\widehat\nu_1$ on $B_1$ and $\chi_1=0$ on $X\setminus B_1$, we can split
		\[
		\begin{split}
			S_{|\nu_1|,\varepsilon,\rho}\widehat b_1(x)
			&=\chi_2(x)\int_{d(x,y)>\rho}a_\varepsilon(d(x,y))\,K_{\Xi_{\mathtt w}}(x,y)\,d\widehat\nu_1(y)\\
			&\qquad-\chi_2(x)\int_{d(x,y)>\rho}\bigl(1-\chi_1(y)\bigr)\,a_\varepsilon(d(x,y))\,K_{\Xi_{\mathtt w}}(x,y)\,d\widehat\nu_1(y)
			\eqqcolon\mathfrak I_1-\mathfrak I_2.
		\end{split}
		\]
		Since $0\leq a_\varepsilon\leq1$ and $\int_0^\infty|a_\varepsilon'|\leq2$, Fubini's theorem gives, as in the proof of Lemma \ref{lemma:reduction_cutoffs},
		\[
		\mathfrak I_1=\chi_2(x)\Bigl(a_\varepsilon(\rho)\,T_{\Xi_{\mathtt w},\rho}\widehat\nu_1(x)+\int_\rho^\infty a_\varepsilon'(t)\,T_{\Xi_{\mathtt w},t}\widehat\nu_1(x)\,dt\Bigr),
		\]
		so that $|\mathfrak I_1|\leq3\,T_{\Xi_{\mathtt w},\sharp}\widehat\nu_1(x)\lesssim1$, by \eqref{eq:global_testing_before_localization_sec5}. If $\chi_1(y)\neq1$, then $\Xi_{\mathtt w}(y)\geq\theta r(B_1)/2$, by \eqref{eq:collar_bound_sec6}. Hence, as above,
		\[
		|\mathfrak I_2|\leq\int_{B_1}\bigl(1-\chi_1(y)\bigr)\,|K_{\Xi_{\mathtt w}}(x,y)|\,d|\widehat\nu_1|(y)
		\lesssim\frac{|\widehat\nu_1|(B_1)}{(\theta r(B_1))^n}\lesssim1.
		\]
		This proves the bound for the first term in \eqref{eq:localized_maximal_testing_sec5}. The second one follows in the same way, by interchanging the variables and the roles of $\chi_1$ and $\chi_2$, and by using $(T^*)_{\Xi_{\mathtt w},\sharp}\widehat\nu_2$ in place of $T_{\Xi_{\mathtt w},\sharp}\widehat\nu_1$.
	\end{proof}
	
	Increasing $c_\sharp$ if necessary, we may assume that the left-hand sides of \eqref{eq:global_testing_before_localization_sec5} and \eqref{eq:localized_maximal_testing_sec5} are bounded by $c_\sharp$, which depends only on the structural parameters.

	\vvv
	Before proving the Key Lemma for the operator $S_{|\nu_1|,\varepsilon}$, we introduce a truncated version of the projection operators associated with the martingale difference decompositions.
	First, note that if $f$ and $g$ are as in the statement of Lemma \ref{lem:key_lemma_good_functions} and we define
	\begin{equation}\label{eq:md_decomposition}
		f_{K_1}=E_1f + \sum_{Q\in \bbD^{\tr}_1, \ell(Q) > \el^{K_1} }\Delta_{1,Q} f \qquad \text{ and }\qquad g_{K_2}=E_2g + \sum_{P\in \bbD^{\tr}_2, \ell(P) > \el^{K_2}}\Delta_{2,P}g, 
	\end{equation}
	then, by \eqref{eq:dec_f_mart}, 
	\begin{equation}\label{eq:md_decomposition_2}
		f_{K_1} \to f \,\, \text{ in } \,\, L^2(|\nu_1|) \quad  \text{ and } \quad g_{K_2} \to g  \,\,\text{ in } \,\, L^2(|\nu_2|) \quad  \text{ unconditionally}.
	\end{equation}
	
	Hence, we provide a truncated analogue of Definition \ref{def:p_good_p_bad_not_truncated}.
	\begin{definition}\label{def:p_good_p_bad_truncated}
		Let $K_1,K_2>0$, $\mathtt w = (w_1, w_2)\in \Omega^2$, $f\in L^2(|\nu_1|)$ and $g\in L^2(|\nu_2|)$. We define
		\[
		\mathcal P^{K_1}_{\mathcal G,1}(\mathtt w)f\coloneqq \sum_{Q\in \bbD^{\mathcal G}_{1}(\mathtt w)\cap\bbD_1^{\tr},\,  \ell(Q)\geq \el^{K_1} }\Delta_{1,Q}f, \qquad \qquad \mathcal P^{K_1}_{\mathcal B,1}(\mathtt w)f\coloneqq \sum_{Q\in \bbD^{\mathcal B}_{1}(\mathtt w)\cap\bbD_1^{\tr},\,  \ell(Q)\geq \el^{K_1} } \Delta_{1,Q}f,
		\]
		and 
		\[
		\mathcal P^{K_2}_{\mathcal G,2}(\mathtt w)g\coloneqq \sum_{R\in \bbD^{\mathcal G}_{2}(\mathtt w)\cap\bbD_2^{\tr},\,  \ell(R)\geq \el^{K_2}}\Delta_{2,R}g, \qquad \qquad \mathcal P^{K_2}_{\mathcal B,2}(\mathtt w)g\coloneqq \sum_{R\in \bbD^{\mathcal B}_{2}(\mathtt w)\cap\bbD_2^{\tr},\,  \ell(R)\geq \el^{K_2} }\Delta_{2,R}g.
		\]
	\end{definition}
	\vv
	
	To simplify notation, for $\mathtt w=(w_1,w_2)\in \Omega^2$ and $K_1,K_2>1$, we write $f_{\mathtt w}\coloneqq f_{\mathtt w}^{K_1}\coloneqq \mathcal P^{K_1}_{\mathcal G,1}(\mathtt w)f$ and $g_{\mathtt w}\coloneqq g_{\mathtt w}^{K_2}\coloneqq \mathcal P^{K_2}_{\mathcal G,2}(\mathtt w)g$.
	The truncated singular integrals satisfy
	\[
	\|S_{|\nu_1|,\varepsilon}\|_{L^2(|\nu_1|)\to L^2(|\nu_2|)}<\infty \qquad \text{ and }\qquad \|S^*_{|\nu_2|,\varepsilon}\|_{L^2(|\nu_2|)\to L^2(|\nu_1|)}<\infty.
	\]
	So, \eqref{eq:md_decomposition_2} yields
	\[
	S_{|\nu_1|,\varepsilon}f^{K_1}_{\mathtt w} = S_{|\nu_1|,\varepsilon}\bigl(E_1 f^{K_1}_{\mathtt w}\bigr)+ \sum_{Q\in \bbD^{\tr}_1} S_{|\nu_1|,\varepsilon}\bigl(\Delta_{1,Q}f^{K_1}_{\mathtt w}\bigr)\overset{{K_1}\to \infty}{\longrightarrow}S_{|\nu_1|,\varepsilon}\bigl(\mathcal P^1_{\mathcal G}(\mathtt w)f\bigr)
	\]
	and 
	\[
	S^*_{|\nu_2|,\varepsilon}g^{K_2}_{\mathtt w}= S^*_{|\nu_2|,\varepsilon}(E_2g^{K_2}_{\mathtt w}) + \sum_{P\in \bbD^{\tr}_2} S^*_{|\nu_2|,\varepsilon}(\Delta_{2,P}g^{K_2}_{\mathtt w})\overset{{K_2}\to \infty}{\longrightarrow}S^*_{|\nu_2|,\varepsilon}\bigl(\mathcal P^2_{\mathcal G}(\mathtt w)g\bigr),
	\]
	where the sums converge in $L^2(|\nu_2|)$ and $L^2(|\nu_1|)$, respectively.
	
	\begin{remark}
		The sums appearing in the expressions of $S_{|\nu_1|,\varepsilon} f^{K_1}_{\mathtt w}$ and $S^*_{|\nu_2|,\varepsilon}g^{K_2}_{\mathtt w}$ contain only finitely many nonzero terms; therefore, we may interchange the order in summations. Moreover, all the bounds that we prove in the remainder of this section are independent of $K_1$ and $K_2$, so we may safely pass to the limit. For this reason, we simply write $f_{\mathtt w}$ and $g_{\mathtt w}$. Finally, we remark that, in the terminology of \cite[Chapter 5]{To14}, studying the contributions of $S_{|\nu_1|,\varepsilon} f^{K_1}_{\mathtt w}$ and $S^*_{|\nu_2|,\varepsilon}g^{K_2}_{\mathtt w}$ corresponds to considering \emph{very good functions}.
	\end{remark}

	\vv

	\subsection{Proof of Key Lemma \ref{lem:key_lemma_good_functions}}\label{section:proof:key_lemma}
	By \eqref{eq:localized_pairing_comparison_sec5}, we estimate the bilinear form of $S_{|\nu_1|,\varepsilon}$. Recall that $f_{\mathtt w}$ and $g_{\mathtt w}$ are supported in $Q_2(w_1)$ and $Q_2(w_2)$, respectively, and both top cubes contain $B_1\cup B_2$. The martingale decompositions give
	\begin{equation}\label{eq:martingale_difference_decomposition_extended}
		\begin{split}
			\langle S_{|\nu_1|,\varepsilon}f_{\mathtt w},g_{\mathtt w}\rangle_*
			&=-\langle S_{|\nu_1|,\varepsilon}(E_1f_{\mathtt w}),E_2g_{\mathtt w}\rangle_*
			+\langle S_{|\nu_1|,\varepsilon}(E_1f_{\mathtt w}),g_{\mathtt w}\rangle_*\\
			&\quad+\langle S_{|\nu_1|,\varepsilon}f_{\mathtt w},E_2g_{\mathtt w}\rangle_*
			+\sum_{Q\in\bbD_1^{\tr},\,R\in\bbD_2^{\tr}}
			\langle S_{|\nu_1|,\varepsilon}(\Delta_{1,Q}f_{\mathtt w}),\Delta_{2,R}g_{\mathtt w}\rangle_*\\
			&\eqqcolon\mathcal I_1(\mathtt w)+\mathcal I_2(\mathtt w)
			+\mathcal I_3(\mathtt w)+\mathcal I_4(\mathtt w).
		\end{split}
	\end{equation}
	
	Since $\int\Delta_{i,Q}h\,d|\nu_i|=0$ for every $Q\in\bbD_i^{\tr}$ and $h\in L^2(|\nu_i|)$, see the proof of Lemma \ref{lemma:probability_bad_1}, we have $E_1f_{\mathtt w}=0$ and $E_2g_{\mathtt w}=0$, so that $\mathcal I_1(\mathtt w)$, $\mathcal I_2(\mathtt w)$, and $\mathcal I_3(\mathtt w)$ vanish. However, the following bounds for these terms hold with arbitrary $f\in L^2(|\nu_1|)$ and $g\in L^2(|\nu_2|)$ in place of $f_{\mathtt w}$ and $g_{\mathtt w}$, and they are used in this generality in Section \ref{sec:main_theorem_antisym}.
	We recall that, by the testing assumptions, we have
	\begin{equation}\label{eq:term_I_1_1}
		|\widehat\nu_i|(B_i)\leq\mathfrak C_1\mu(B_i)\qquad\text{ and }\qquad 
		\mu(Q_2(w_i))\leq|\nu_i|(Q_2(w_i))\leq\mathfrak C_1\mu(Q_2(w_i)).
	\end{equation}
	The Cauchy-Schwarz inequality and \eqref{eq:term_I_1_1} yield
	\[
	\biggl|\frac{\int_{Q_2(w_i)}u\,d|\nu_i|}{\mu(Q_2(w_i))}\biggr|
	\lesssim\mu(Q_2(w_i))^{-1/2}\|u\|_{L^2(|\nu_i|)},
	\]
	so $\|E_i u\|_{L^2(|\nu_i|)}\lesssim\|u\|_{L^2(|\nu_i|)}$.
	By the choice of the cutoffs in $L_\varepsilon$, we have the identities 
	\[
	S_{|\nu_1|,\varepsilon}\bigl(\chi_{Q_2(w_1)}\widehat b_1\bigr)
	=S_{|\nu_1|,\varepsilon}\widehat b_1
	\qquad\text{ and }\qquad 
	S^*_{|\nu_2|,\varepsilon}\bigl(\chi_{Q_2(w_2)}\widehat b_2\bigr)
	=S^*_{|\nu_2|,\varepsilon}\widehat b_2.
	\]
	Moreover, the two outputs are supported in $B_2$ and $B_1$, respectively. Consequently, \eqref{eq:localized_maximal_testing_sec5} and \eqref{eq:term_I_1_1} imply
	\[
	\bigl\|S_{|\nu_1|,\varepsilon}\widehat b_1\bigr\|_{L^2(|\nu_2|)}^2
	\lesssim\mu(B_2)\leq\mu(Q_2(w_1)),\qquad
	\bigl\|S^*_{|\nu_2|,\varepsilon}\widehat b_2\bigr\|_{L^2(|\nu_1|)}^2
	\lesssim\mu(B_1)\leq\mu(Q_2(w_2)),
	\]
	and combining these bounds with the controls on the top averages, we obtain
	\begin{equation}\label{eq:final_bound_I_1}
		|\mathcal I_1(\mathtt w)|
		\leq\bigl\|S_{|\nu_1|,\varepsilon}(E_1f_{\mathtt w})\bigr\|_{L^2(|\nu_2|)}
		\|E_2g_{\mathtt w}\|_{L^2(|\nu_2|)}
		\lesssim\|f_{\mathtt w}\|_{L^2(|\nu_1|)}\|g_{\mathtt w}\|_{L^2(|\nu_2|)}.
	\end{equation}
	We also have the following estimates.
	\begin{lemma}\label{lemma:unbalanced_lemma}
		We have
		\begin{equation}\label{eq:inbalanced_estimate_1}
			\|S_{|\nu_1|,\varepsilon}\widehat b_1\|_{L^2(|\nu_2|)}
			\lesssim\mu(Q_2(w_1))^{1/2}
		\end{equation}
		and
		\begin{equation}\label{eq:inbalanced_estimate_2}
			|\mathcal I_2(\mathtt w)|
			\lesssim\|f_{\mathtt w}\|_{L^2(|\nu_1|)}\|g_{\mathtt w}\|_{L^2(|\nu_2|)}.
		\end{equation}
	\end{lemma}
	\begin{proof}
		The first estimate was proved above. For the second, by the definition of $E_1$,
		\[
		|\mathcal I_2(\mathtt w)|
		\leq\frac{\Bigl|\int_{Q_2(w_1)}f_{\mathtt w}\,d|\nu_1|\Bigr|}{\mu(Q_2(w_1))}\, 
		\|S_{|\nu_1|,\varepsilon}\widehat b_1\|_{L^2(|\nu_2|)}
		\|g_{\mathtt w}\|_{L^2(|\nu_2|)}
		\lesssim\|f_{\mathtt w}\|_{L^2(|\nu_1|)}\|g_{\mathtt w}\|_{L^2(|\nu_2|)}.\qedhere
		\]
	\end{proof}
	Finally, by \eqref{eq:localized_maximal_testing_sec5} and \eqref{eq:term_I_1_1}, we obtain
	\begin{equation}\label{eq:final_bound_I_3_alternative}
		\begin{split}
			|\mathcal I_3(\mathtt w)|
			&\leq\frac{\Bigl|\int_{Q_2(w_2)}g_{\mathtt w}\,d|\nu_2|\Bigr|}{\mu(Q_2(w_2))}\, 
			\|S^*_{|\nu_2|,\varepsilon}\widehat b_2\|_{L^2(|\nu_1|)}
			\|f_{\mathtt w}\|_{L^2(|\nu_1|)}\\
			&\lesssim\|f_{\mathtt w}\|_{L^2(|\nu_1|)}\|g_{\mathtt w}\|_{L^2(|\nu_2|)}.
		\end{split}
	\end{equation}
	
	\vv
	
	It remains to estimate $\mathcal I_4(\mathtt w)$, which requires a more delicate analysis.
	We recall that $\alpha=\gamma/(2(n+\gamma))$, which belongs to $(0,1)$.
	We say that $Q\in \bbD^{\tr}_1$ and $R\in \bbD^{\tr}_2$ are \textit{distant} if 
	\begin{equation}\label{eq:def_distant_cubes}
		\dist(Q,R)\geq \min(\ell(Q),\ell(R))^\alpha\cdot  \max(\ell(Q), \ell(R))^{1-\alpha}.
	\end{equation}
	Note that, if $Q$ and $R$ are both good, disjoint, and not distant, by Definition \ref{def:good_bad_cubes} it holds that
	\[
	\el^{m}\ell(R)<\ell(Q)<\el^{-m}\ell(R).
	\]
	
	Hence, we split $\mathcal I_4(\mathtt w)$ as
	\begin{equation}
		\begin{split}
			\mathcal I_4(\mathtt w) &= \sum_{Q\in \bbD^{\tr}_1, R \in \bbD^{\tr}_2} \bigl\langle S_{|\nu_1|,\varepsilon}(\Delta_{1,Q} f_{\mathtt w}), \Delta_{2,R} g_{\mathtt w}\bigr\rangle_*\\
			&=\sum_{\substack{Q\in \bbD^{\tr}_1, R \in \bbD^{\tr}_2\\ Q, R\text{ distant }}} \quad + \sum_{\substack{Q\in \bbD^{\tr}_1, R \in \bbD^{\tr}_2\\ Q\cap R=\varnothing\\ Q, R\text{ not distant }}} \quad + \sum_{\substack{Q\in \bbD^{\tr}_1, R \in \bbD^{\tr}_2\\Q\cap R\neq \varnothing \\ \el^{m}\ell(R)\leq \ell(Q)\leq \el^{-m}\ell(R)}} \quad + \sum_{\substack{Q\in \bbD^{\tr}_1, R \in \bbD^{\tr}_2\\Q\cap R\neq \varnothing \\ \ell(Q)<\el^{m}\ell(R) \text{ or }\ell(Q)>\el^{-m}\ell(R)}} \quad \\
			&\eqqcolon \mathcal S_1(\mathtt w)+ \mathcal S_2(\mathtt w)+ \mathcal S_3(\mathtt w)+ \mathcal S_4(\mathtt w).
		\end{split}
	\end{equation}
	
	We estimate these four terms separately.
	
	\vv
	
	\subsection*{{Estimate of $\mathcal S_1$}}
	We first write
	\begin{equation*}
		\begin{split}
			|\mathcal S_1(\mathtt w)|&\leq \sum_{\substack{Q\in \bbD^{\tr}_1, R \in \bbD^{\tr}_2\\ Q, R\text{ distant }}}\bigl|\bigl\langle S_{|\nu_1|,\varepsilon}(\Delta_{1,Q} f_{\mathtt w}), \Delta_{2,R} g_{\mathtt w}\bigr\rangle_*\bigr|= \sum_{\substack{Q\in \bbD^{\tr}_1, R \in \bbD^{\tr}_2\\ Q, R\text{ distant }\\ \ell(Q)\leq \ell(R)}} \quad + \sum_{\substack{Q\in \bbD^{\tr}_1, R \in \bbD^{\tr}_2\\ Q, R\text{ distant }\\\ell(Q)> \ell(R)}}\quad \eqqcolon \mathcal S'_{1}+ \mathcal S_1''.
		\end{split}
	\end{equation*}
	
	The \textit{long distance} is defined as
	\[
	D(Q,R)\coloneqq \ell(Q)+ \ell(R)+ \dist(Q,R), \qquad \text{ for }Q\in\bbD_1, R\in\bbD_2,
	\]
	and it is used in the bound of $\mathcal  S_1'$ via the following lemma.
	
	\begin{lemma}\label{lemma:lem_514}
		Let $Q_1, Q_2$ be cubes in $X$. For $i=1,2$, assume that $\varphi_i\in L^2(|\nu_i|)$ are such that $\supp(\varphi_i)\subset Q_i$, $\dist(Q_1, \supp(\varphi_2))\geq \ell(Q_1)$ and $\int \varphi_1\, d|\nu_1|=0.$
		Then, we have that 
		\begin{equation}
			\bigl|\bigl\langle S_{|\nu_1|,\varepsilon}\varphi_1, \varphi_2\bigr\rangle_*\bigr|\leq C\,  \frac{\ell(Q_1)^\gamma}{\dist\bigl(Q_1, \supp(\varphi_2)\bigr)^{n+\gamma}}\, \|\varphi_1\|_{L^1(|\nu_1|)}\|\varphi_2\|_{L^1(|\nu_2|)}.
		\end{equation}
		Moreover, if $\ell(Q_1)\leq \ell(Q_2)$ and there exists $\mathfrak c'>0$ such that 
		\[
		\dist\bigl(Q_1, \supp(\varphi_2)\bigr)\geq \mathfrak c'\, \ell(Q_1)^\alpha\, \ell(Q_2)^{1-\alpha},
		\]
		then, denoting by $C_1$ the constant in \eqref{eq:diam_vs_length},
		\begin{equation}\label{eq:bound_zero_mean}
			\bigl|\bigl\langle S_{|\nu_1|,\varepsilon}\varphi_1, \varphi_2\bigr\rangle_*\bigr|\lesssim_{C, \mathfrak c', C_1} \frac{\ell(Q_1)^{\gamma/2}\, \ell(Q_2)^{\gamma/2}}{D(Q_1,Q_2)^{n+\gamma}}\, \|\varphi_1\|_{L^1(|\nu_1|)}\|\varphi_2\|_{L^1(|\nu_2|)}.
		\end{equation}
	\end{lemma}
	\begin{proof}
		The proof is a minor variation of the proof of \cite[Lemma 5.14]{To14}, using Calder\'on-Zygmund kernel estimates and the choice $\alpha=\gamma/(2(n+\gamma))$ in the definition of distant cubes.
	\end{proof}
	\vv
	
	Since the cubes $Q$ and $R$ in $\mathcal S_1'$ are distant and $\ell(Q)\leq \ell(R)$, Lemma \ref{lemma:lem_514} yields
	\begin{equation}\label{eq:s_1_prime_aux}
		\mathcal S_1'\lesssim_{C} \sum_{\substack{Q\in \bbD^{\tr}_1, R \in \bbD^{\tr}_2\\ \ell(Q)\leq \ell(R)}} \frac{\ell(Q)^{\gamma/2}\, \ell(R)^{\gamma/2}}{D(Q,R)^{n+\gamma}}\, |\nu_1|(Q)^{1/2}|\nu_2|(R)^{1/2}\, \|\Delta_{1,Q}f_{\mathtt w}\|_{L^2(|\nu_1|)}\, \|\Delta_{2,R}g_{\mathtt w}\|_{L^2(|\nu_2|)}.
	\end{equation}
	
	{	
		For transit cubes $Q\in \bbD^{\tr}_1$ and $R\in \bbD^{\tr}_2$ with positive measure, we define
		\[
		T_{Q,R}\coloneqq \frac{\ell(Q)^{\gamma/2}\ell(R)^{\gamma/2}}
		{D(Q,R)^{n+\gamma}}\,|\nu_1|(Q)^{1/2}|\nu_2|(R)^{1/2},\qquad \text{ if }\ell(Q)\leq\ell(R)
		\]
		and set $T_{Q,R}\coloneqq 0$ otherwise. We claim that
		\begin{equation}\label{eq:sec5_two_measure_schur}
			\begin{split}
				\sum_{\ell(Q)\leq\ell(R)}T_{Q,R}|\nu_1|(Q)^{1/2}
				&\lesssim |\nu_2|(R)^{1/2},\\
				\sum_{\ell(R)\geq\ell(Q)}T_{Q,R}|\nu_2|(R)^{1/2}
				&\lesssim |\nu_1|(Q)^{1/2}.
			\end{split}
		\end{equation}
		Note that, by Schur's test, \eqref{eq:sec5_two_measure_schur} implies that the matrix $(T_{Q,R})$ defines a bounded operator on the corresponding $\ell^2$ spaces, with norm bounded independently of the random grids.
		
		We now prove \eqref{eq:sec5_two_measure_schur}. Let $i=1,2$, and let $\mathcal A\subset\bbD_i^{\tr}$ be a nonempty family of cubes with the same side length, all contained in a ball $B(z,t)$. Since the cubes in $\mathcal A$ are transit, there exists $x\in\bigcup_{P\in\mathcal A}P\setminus\mathscr W_i(w_i)$. As in Remark \ref{remark:growth_children_transit_cubes}, we have $x\notin\mathscr H_{i,R_i}$, so \eqref{eq:definition_r_i}, \eqref{eq:def_non_ahlfors}, and the $n$-growth of $\mu$ give $|\nu_i|(B(x,\rho))\leq(\mathbf m+\mathfrak C_0)\,\rho^n$ for every $\rho>0$. Moreover, the cubes in $\mathcal A$ are pairwise disjoint and contained in $B(z,t)\subset B(x,2t)$, thus we have
		\begin{equation}\label{eq:packing_transit_cubes}
			\sum_{P\in\mathcal A}|\nu_i|(P)=|\nu_i|\Bigl(\bigcup_{P\in\mathcal A}P\Bigr)			\leq|\nu_i|(B(x,2t))
			\leq2^n(\mathbf m+\mathfrak C_0)\, t^n.
		\end{equation}
		
		Let $R\in\bbD_2^{\tr}$ and denote $r\coloneqq\ell(R)$. Since $D(Q,R)\geq r$, we can group the cubes $Q\in\bbD_1^{\tr}$ with $\ell(Q)\leq r$ according to the values $k,j\geq0$ such that $\ell(Q)=\el^kr$ and $2^jr\leq D(Q,R)<2^{j+1}r$. Moreover, by \eqref{eq:diam_vs_length}, each group is contained in a ball of radius $(2C_1+1)\,2^{j+1}r$, so \eqref{eq:packing_transit_cubes} yields
		\[
		\sum_{\substack{Q\in\bbD_1^{\tr}\\ \ell(Q)\leq r}}\frac{\ell(Q)^{\gamma/2}\,r^{\gamma/2}}{D(Q,R)^{n+\gamma}}\,|\nu_1|(Q)
		\lesssim\sum_{k,j\geq0}\frac{(\el^kr)^{\gamma/2}\,r^{\gamma/2}}{(2^jr)^{n+\gamma}}\,(2^jr)^n
		=\sum_{k,j\geq0}\el^{k\gamma/2}\,2^{-j\gamma}\lesssim1,
		\]
		which implies the first estimate in \eqref{eq:sec5_two_measure_schur}. The second estimate follows analogously.
		\vv
	}
	
	Thus, by \eqref{eq:s_1_prime_aux}, \eqref{eq:sec5_two_measure_schur}, \eqref{eq:dec_f_mart_2}, and Schur's test we obtain
	\[
	\mathcal S_1'\lesssim \|f_{\mathtt w}\|_{L^2(|\nu_1|)}\|g_{\mathtt w}\|_{L^2(|\nu_2|)}.
	\]
	
	We observe that, for $Q\in \bbD_1^{\tr}$ and $R\in \bbD_2^{\tr}$ we have
	\begin{equation}\label{eq:transpose_martingale_difference}
		\bigl\langle S_{|\nu_1|,\varepsilon}\Delta_{1,Q}f_{\mathtt w}, \Delta_{2,R}g_{\mathtt w} \bigr\rangle_* = \bigl \langle \Delta_{1,Q}f_{\mathtt w}, S^*_{|\nu_2|,\varepsilon}\Delta_{2,R}g_{\mathtt w}\bigr\rangle,
	\end{equation}
	so the same argument can be applied to estimate $\mathcal S_1''$, and we finally get
	\[
	|\mathcal S_1(\mathtt w)|\lesssim \|f_{\mathtt w}\|_{L^2(|\nu_1|)}\, \|g_{\mathtt w}\|_{L^2(|\nu_2|)}\overset{\eqref{eq:dec_f_mart_2}}{\lesssim}  \|f\|_{L^2(|\nu_1|)}\|g\|_{L^2(|\nu_2|)},
	\]
	which implies that there exists $\mathbf c_{\mathcal S1}>0$ such that
	\begin{equation}\label{eq:bound_s_1_final}
		\widetilde{\mathbb E}_{\mathtt w}|\mathcal S_1(\mathtt w)|\leq \mathbf c_{\mathcal S1} \|f\|_{L^2(|\nu_1|)}\|g\|_{L^2(|\nu_2|)}.
	\end{equation}
	
	\vv
	
	\subsection*{{Estimates of $\mathcal S_2$ and $\mathcal S_3$: the diagonal terms}}
	
	{
		
		We recall that $\Delta_{1,Q}f_{\mathtt w}=0$ unless $Q$ is good, and $\Delta_{2,R}g_{\mathtt w}=0$ unless $R$ is good. Moreover, as observed above, if $Q\in \bbD^{\tr}_1$ and $R\in \bbD^{\tr}_2$ are good, disjoint, and not distant, then $\el^{m}\ell(R)<\ell(Q)<\el^{-m}\ell(R)$. Given two cubes $Q, R\subset X$, we also write $\ell(Q)\sim_m \ell(R)$ if $\el^{m}\, \ell(Q)\leq \ell(R)\leq \el^{-m}\, \ell(Q)$. 
		
		Moreover, for $\alpha= \gamma/(2(n+\gamma))$ as in \eqref{eq:def_distant_cubes},
		\[
		\dist(Q,R)\leq \min(\ell(Q), \ell(R))^\alpha\, \max(\ell(Q), \ell(R))^{1-\alpha}\leq \max(\ell(Q), \ell(R)).
		\]
		
		Thus we have 
		\begin{equation}\label{eq:bound_s_2_s_3}
			|\mathcal S_2 (\mathtt w) + \mathcal S_3(\mathtt w)|
			\leq \sum_{\substack{Q\in \bbD^{\tr}_1, R \in \bbD^{\tr}_2\\\dist(Q,R)\leq \max(\ell(Q), \ell(R)) \\ \ell(Q)\sim_m\ell(R)}}\bigl|\bigl\langle S_{|\nu_1|,\varepsilon}\Delta_{1,Q}f_{\mathtt w}, \Delta_{2,R}g_{\mathtt w}\bigr\rangle_*\bigr|.
		\end{equation}
		
		Let $\mathfrak D^*(\mathtt w)\subset \bbD^{\tr}_1\times \bbD^{\tr}_2 $  be the collection of pairs $(Q,R)$ such that $\dist(Q,R)\leq \max(\ell(Q), \ell(R))$ and $ \ell(Q)\sim_m\ell(R)$. We also denote by $\mathfrak D(\mathtt w)$ the collection of pairs $(I,J)\in \bbD(w_1)\times \bbD(w_2) $ such that $I\in \Ch(Q)$ and $J\in \Ch(R)$ for some $(Q,R)\in \mathfrak D^*(\mathtt w)$.
		
		Thus, the right-hand side of \eqref{eq:bound_s_2_s_3} reads
		\begin{equation}\label{eq:key_sum_diagonal}
			\sum_{(Q,R)\in \mathfrak D^*(\mathtt w)}\bigl|\bigl\langle S_{|\nu_1|,\varepsilon}\Delta_{1,Q}f_{\mathtt w}, \Delta_{2,R}g_{\mathtt w}\bigr\rangle_*\bigr|\leq \sum_{(I,J)\in \mathfrak D(\mathtt w)}\bigl|\bigl\langle S_{|\nu_1|,\varepsilon}\bigl(\chi_I \Delta_{1,I^*}f_{\mathtt w}\bigr), \chi_J\Delta_{2,J^*}g_{\mathtt w}\bigr\rangle_*\bigr|,
		\end{equation}
		where $I^*\in \bbD^{\tr}_1$ and $J^*\in \bbD^{\tr}_2$ are the dyadic parents of $I$ and $J$, respectively.
		If $(I,J)\in\mathfrak D(\mathtt w)$, then $\ell(J)\sim_m\ell(I)$ and $\dist(I^*,J^*)\leq\max(\ell(I^*),\ell(J^*))\leq\el^{-m-1}\ell(I)$. Hence, by \eqref{eq:diam_vs_length}, for fixed $I$, the cubes $J$ such that $(I,J)\in\mathfrak D(\mathtt w)$ belong to $2m+1$ generations and are contained in a ball of radius $C_m\ell(I)$. Since $X$ is geometrically doubling and each cube $J$ contains a ball of radius $c_1\ell(J)$, we obtain the first bound in
		\begin{equation}\label{eq:counting_pairs}
			\#\{J:(I,J)\in\mathfrak D(\mathtt w)\}\leq C_m
			\qquad\text{ and }\qquad
			\#\{I:(I,J)\in\mathfrak D(\mathtt w)\}\leq C_m,
		\end{equation}
		and the second one, for fixed $J$, follows analogously.
		
		We define $\mathfrak D^{\tr}_{\tr}(\mathtt w)\coloneqq \bigl\{(I,J)\in \mathfrak D(\mathtt w): \text{both }I \text{ and } J \text{ are transit}\bigr\}$.
		If $(I_1, I_2)\in  \mathfrak D^{\tr}_{\tr}(\mathtt w)$, we also denote
		\[
		D^j_{I_j}\phi \coloneqq \frac{\langle \phi\rangle_{j,I_j}}{\langle \widehat b_j\rangle_{j,I_j}}-\frac{\langle \phi\rangle_{j,I^*_j}}{\langle \widehat b_j\rangle_{j,I^*_j}}, \qquad \text{ for }j=1,2\text{ and }\phi \in L^2(|\nu_j|),
		\]
		where $I_j^*\in \bbD(w_j)$ is the dyadic parent of $I_j$. Thus, note that by the definition \eqref{eq:def_martingales}, if $(I,J)\in \mathfrak D^{\tr}_{\tr}(\mathtt w)$, we have
		\begin{equation}
			\chi_I\, D^1_{I}f_{\mathtt w}\, \widehat b_1= \chi_I\, \Delta_{1,I^*}f_{\mathtt w}\qquad \text{ and }\qquad \chi_J\, D^2_{J}g_{\mathtt w}\, \widehat b_2= \chi_J\, \Delta_{2,J^*}g_{\mathtt w},
		\end{equation}
		so
		\begin{equation}\label{eq:sum_surgery}
			\begin{split}
				\bigl|\bigl\langle S_{|\nu_1|,\varepsilon}\bigl(\chi_I \Delta_{1,I^*}f_{\mathtt w}\bigr), \chi_J\Delta_{2,J^*}g_{\mathtt w}\bigr\rangle_*\bigr|&=	\bigl|\bigl\langle S_{|\nu_1|,\varepsilon} \bigl(D^1_If_{\mathtt w}\,\chi_{I}\,\widehat b_1\bigr), D^2_Jg_{\mathtt w}\, \chi_J\, \widehat b_2\bigr\rangle_*\bigr|.
			\end{split}
		\end{equation}
		
		\vv
		\subsubsection*{The random almost-covering balls}

		In the next proposition, we state the properties of a ``random'' family of balls constructed by Hyt\"onen and Martikainen in the quasi-metric setting. 
		The balls have radii comparable to $\vartheta^k$, where $k\in \mathbb Z$ and $0<\vartheta<1/32$ replaces factor $1/2$ in the classical Euclidean dyadic system.
		The construction depends on a random parameter $\xi$, which encodes the random choices of the centers and of the dilation factors of the balls, in an auxiliary probability space $(\widehat\Omega,\widehat{\mathbb P})$, independent of $\mathtt w$. We denote by $\widehat{\mathbb E}$ the associated expectation. Taking a product over the scales $k\in\mathbb Z$, we use the same parameter $\xi$ at every scale.
		
		\begin{proposition}[see Proposition 8.2 in \cite{HM12}]\label{proposition:random_balls}
			Let $0<\vartheta<1/32$, $k\in \mathbb Z$, and $\upsilon\in (0,1)$.
			Then, there exist $\widetilde \Lambda\in (0,1)$ and integers $s,M\geq1$, depending only on $D$, $\vartheta$, and $\upsilon$, and, for every $\xi\in\widehat\Omega$, a countable family $\mathfrak B_\xi$ of balls in $X$ such that:
			\begin{itemize}
				\item If $B, B'\in \mathfrak B_{\xi},$ then 
				\begin{equation}\label{eq:radii_random_balls}
					\frac{1}{32}\, \vartheta^{k+n_1\, s}\leq	r(B)\leq \frac{1}{16}\, \vartheta^{k+n_1\, s}\qquad \text{ and }\qquad \frac{1}{32}\, \vartheta^{k+n_2\, s}\leq	r(B')\leq \frac{1}{16}\, \vartheta^{k+n_2\, s},
				\end{equation}
				for some $n_1, n_2= 0, 1, \ldots, M-1$. Moreover, if $B\neq B'$, then
				\[
				d(x_B,x_{B'})\geq r(B)+r(B')+\widetilde\Lambda\max(r(B),r(B')),
				\]
				where $x_B$ and $x_{B'}$ are the centers of $B$ and $B'$, respectively.
				\item We have
				\begin{equation}\label{eq:high_probability_balls}
					\widehat{\mathbb P}\Bigl(x\in \bigcup_{B\in \mathfrak B_{\xi}}B\Bigr)>1-\upsilon, \qquad \text{ for every } x\in X.
				\end{equation}
			\end{itemize}
			Moreover, the set $\{(\xi,x):x\in\bigcup_{B\in\mathfrak B_\xi}B\}$ is measurable.
		\end{proposition}
		
		In \cite{HM12}, the separation is stated for the balls themselves, but the separation of their centers follows from the construction. Indeed, if $n_1=n_2$ in \eqref{eq:radii_random_balls}, then 
		\[
		d(x_B,x_{B'})\geq\vartheta^{k+n_1s}/4\geq2(r(B)+r(B')).
		\]
		If $n_1<n_2$, then $x_{B'}\notin(1+\omega)B$ and $r(B')\leq2\vartheta^sr(B)<\omega\,r(B)/3$, where $\omega=\omega(\upsilon)\in(0,1)$ is the buffer parameter in \cite[Section 8]{HM12}, and $s$ is chosen so that $2\vartheta^s<\omega/3$. Hence, Proposition \ref{proposition:random_balls} holds with $\widetilde\Lambda=2\omega/3$. The construction can be carried out with centers in fixed countable nets and finitely many random choices at each level, which gives the measurability in Proposition \ref{proposition:random_balls}.
		
		Note that the balls in $\mathfrak B_{\xi}$ have comparable radii. Indeed, \eqref{eq:radii_random_balls} implies
		\[
		\frac{1}{2}\, \vartheta^{Ms}\leq \frac{r(B)}{r(B')}\leq 2\, \vartheta^{-Ms}.
		\]
		
		\vv
		
		For $\vartheta,$ $k$, and $\upsilon$ as in Proposition \ref{proposition:random_balls} and $\xi\in\widehat\Omega$, we denote by $\mathfrak B_\xi$ the associated family of balls. Moreover, if $M_\Lambda$ is chosen big enough depending on $\widetilde \Lambda$, by Lemma \ref{lemma:balls_small_boundary} and the remark following it, for every $B\in \mathfrak B_\xi$ there exists a ball $\widehat B$ centered at $x_B$ with $M_\Lambda$-$|\nu_1|$- and $M_\Lambda$-$|\nu_2|$- small boundary, whose radius can be chosen to depend measurably on $x_B$ and $r(B)$, such that
		\begin{equation}\label{eq:small_boundary_random_balls}
			B\subset \widehat B\subset \Bigl(1+\frac{\widetilde \Lambda}{8}\Bigr) B.
		\end{equation}
		
		{Let $\lambda\in (0,1)$ be a parameter to be fixed later. For $(I,J)\in \mathfrak D(\mathtt w)$, we choose $k\in \mathbb Z$ such that 
			\begin{equation}\label{eq:choice_k_random_balls}
				\vartheta^k\leq \lambda \min(\ell(I), \ell(J))<\vartheta^{k-1}.
			\end{equation}
			Thus $k$ depends on $\lambda$ and $\min(\ell(I),\ell(J))$.}
		\vv
		
		We denote $\widehat{\mathfrak B}_{\xi} \coloneqq \bigl\{\widehat B: B \in \mathfrak B_{\xi}\bigr\}.$ For distinct $\widehat B,\widehat B'\in\widehat{\mathfrak B}_{\xi}$, \eqref{eq:small_boundary_random_balls}, the triangle inequality, and Proposition \ref{proposition:random_balls} give 
		\[
		\dist(\widehat B,\widehat B')\geq d(x_B,x_{B'})-\Bigl(1+\frac{\widetilde\Lambda}{8}\Bigr)(r(B)+r(B'))\geq\frac{3}{4}\widetilde\Lambda\, \max(r(B),r(B')).
		\]
		Since $r(\widehat B)\leq2r(B)$ and $r(\widehat B')\leq2r(B')$, we obtain
		\begin{equation}\label{eq:separation_random_balls}
			\dist(\widehat B,\widehat B')\geq \frac{\widetilde \Lambda}{8}\, \max(r(\widehat B), r(\widehat B')).
		\end{equation}
		
		The bound \eqref{eq:high_probability_balls} readily implies 
		\[
		\widehat{\mathbb P}\Bigl(x\in \bigcup_{\widehat B\in \widehat{\mathfrak B}_{\xi}}\widehat B\Bigr)>1-\upsilon, \qquad \text{ for every }x\in X,
		\]
		so for every non-negative Radon measure $\mathfrak m$ on $X$ and any $A\subset X$ by Fubini's theorem we have
		\begin{equation}\label{eq:expectation_small_measure_general}
			\widehat{\mathbb E}\Bigl(\mathfrak m\Bigl(A\setminus \bigcup_{\widehat B\in \widehat{\mathfrak B}_{\xi}}\widehat B\Bigr)\Bigr)\leq \upsilon\, \mathfrak m(A).
		\end{equation}

		\vv
		\subsubsection*{The surgery argument}
		
		{			
			Let $(I,J)\in\mathfrak D(\mathtt w)$ and put $h\coloneqq \lambda\min(\ell(I),\ell(J))$. We use the family $\widehat{\mathfrak B}_{\xi}$ at the scale chosen in \eqref{eq:choice_k_random_balls} to decompose both $I$ and $J$. First, we remark that the balls in $\widehat{\mathfrak B}_{\xi}$ satisfy $6r(B)\leq h$.
			Then, we set
			\[
			\mathcal C_J\coloneqq\{B\in\widehat{\mathfrak B}_{\xi}:5B\subset J\},
			\qquad	I_\Delta=I\cap\bigcup_{B\in\mathcal C_J}B,
			\]
			\[
			I_{\mathrm{sep}}=\{x\in I:\operatorname{dist}(x,J)>h\},
			\qquad	I_\partial=I\setminus(I_{\mathrm{sep}}\cup I_\Delta).
			\]
			The sets $J_\Delta,J_{\mathrm{sep}},J_\partial$ are defined analogously.
			These are disjoint decompositions of $I$ and $J$, respectively, and subscript $\mathrm{sep}$ stands for ``separated''. Indeed, we have
			\begin{equation}\label{eq:separation_condition_cubes}
				\operatorname{dist}(I_{\mathrm{sep}},J)\geq h\qquad \text{ and }\qquad  \operatorname{dist}(J_{\mathrm{sep}},I)\geq h.
			\end{equation}
			
			For $s>0$, we denote 
			\[
			\partial_s J\coloneqq\{x\in X:\operatorname{dist}(x,J)\leq s,\operatorname{dist}(x,X\setminus J)\leq s\},
			\]
			and note that \eqref{eq:separation_condition_cubes} implies the set inclusion
			\[
			I_\partial\subset (I\cap\partial_h J)\cup\Bigl(I\setminus\bigcup_{B\in\widehat{\mathfrak B}_{\xi}}B\Bigr).
			\]
			Indeed, a covered point $x\in I\cap J$ lying in a ball $B$ for which $5B\not\subset J$ has distance at most $6r(B)\leq h$ from $X\setminus J$. 
			
			By \eqref{eq:choice_k_random_balls}, the scale of the family $\widehat{\mathfrak B}_{\xi}$ depends on the generation of $J$. For fixed $I$, we denote by $\mathcal K(I)$ the set of the scales $k$ given by \eqref{eq:choice_k_random_balls} for the $2m+1$ generations of cubes $J$ with $\ell(J)\sim_m\ell(I)$, so that $\mathcal K(I)$ contains at most $2m+1$ elements. For each $k\in\mathcal K(I)$, we apply Proposition \ref{proposition:random_balls} with $\upsilon/(2m+1)$ in place of $\upsilon$, so the constants $\widetilde\Lambda$ and $M$ of the corresponding families depend on $m$ and $\upsilon$. We write $\widehat{\mathfrak B}_{\xi,k}$ when the scale needs to be specified. Since $h\leq\lambda\ell(J)$, $\ell(J)\sim_m\ell(I)$, and the scale of $\widehat{\mathfrak B}_{\xi}$ belongs to $\mathcal K(I)$, the previous inclusion gives
			\begin{equation}\label{eq:inclusion_I_partial}
				\begin{split}
					I_\partial\subset{}&
					I\cap\bigcup_{\substack{J'\in\bbD(w_2)\\
							\ell(J')\sim_m\ell(I)}}
					\partial_{\lambda\ell(J')}J'\cup\bigcup_{k\in\mathcal K(I)}
					\biggl(I\setminus\bigcup_{B\in\widehat{\mathfrak B}_{\xi,k}}B\biggr)\\
					&\eqqcolon\widetilde I_{\partial,1}^{w_2}
					\cup\widetilde I_{\partial,2}^{\xi}
					\eqqcolon\widetilde I_\partial^{w_2,\xi}.
				\end{split}
			\end{equation}
			Here, the union is over all the cubes $J'\in\bbD(w_2)$ of the indicated generations, including the terminal ones, and the sets on the right-hand side do not depend on $J$.
			
			We also observe that, if $s>0$, $J'\in\bbD(w_2)$, $x\in\partial_sJ'$, and $J_x\in\bbD(w_2)$ is the cube of the same generation as $J'$ that contains $x$, then $\dist(x,X\setminus J_x)\leq s$. Indeed, if $x\in J'$, then $J_x=J'$. Otherwise, $J'\subset X\setminus J_x$, so $\dist(x,X\setminus J_x)\leq\dist(x,J')\leq s$.
			
			Finally, we write
			\[
			I_\Delta= \bigcup_{i\in\mathscr I_\Delta}B_{I,i}\qquad \text{ and }\qquad J_\Delta=\bigcup_{j\in\mathscr J_\Delta}B_{J,j},
			\]
			where $\{B_{I,i}\}_{i\in\mathscr I_\Delta}$ are the nonempty sets $B\cap I$ with $B\in\mathcal C_J$, and $\{B_{J,j}\}_{j\in\mathscr J_\Delta}$ are the nonempty sets $B'\cap J$ with $B'\in\widehat{\mathfrak B}_{\xi}$ and $5B'\subset I$. If $B_{I,i}=B\cap I$ and $B_{J,j}=B'\cap J$, then either $B=B'$, in which case $5B\subset I\cap J$ and $B_{I,i}=B_{J,j}=B$, or $B\neq B'$, in which case \eqref{eq:separation_random_balls} gives \[\dist(B_{I,i},B_{J,j})\geq\dist(B,B')\geq\frac{\widetilde\Lambda}{8}\max(r(B),r(B')).\] 
		}

		To bound \eqref{eq:key_sum_diagonal} we first consider
		\[
		\mathcal S_{\mathrm{sep}}(\mathtt w)\coloneqq \sum_{(I,J)\in\mathfrak D(\mathtt w)}\bigl|\bigl\langle S_{|\nu_1|,\varepsilon}\bigl(\chi_{I_{\mathrm{sep}}} \Delta_{1,I^*}f_{\mathtt w}\bigr), \chi_J\Delta_{2,J^*}g_{\mathtt w}\bigr\rangle_*\bigr|.
		\]
		
		\begin{lemma}\label{lemma:sep_transit_transit}
			We have
			\[
			\mathcal S_{\mathrm{sep}}(\mathtt w)\lesssim_{\lambda, m} \|f_{\mathtt w}\|_{L^2(|\nu_1|)}\, \|g_{\mathtt w}\|_{L^2(|\nu_2|)}.
			\]
		\end{lemma}
		\begin{proof}
			If $(I,J)\in\mathfrak D(\mathtt w)$, $x\in I_{\mathrm{sep}}$, we note that the separation yields $\dist(x,J)\gtrsim_{\lambda, m}\ell(I)$.
			
			Thus, \eqref{eq:size_auxiliary_kernel_sec6} and the Cauchy-Schwarz inequality yield
			\begin{equation*}
				\begin{split}
					&\bigl|\bigl\langle S_{|\nu_1|,\varepsilon}\bigl(\chi_{I_{\mathrm{sep}}} \Delta_{1,I^*}f_{\mathtt w}\bigr), \chi_J\Delta_{2,J^*}g_{\mathtt w}\bigr\rangle_*\bigr|\lesssim_{\lambda, m}\frac{1}{\ell(I)^n}\, \|\chi_{I_{\mathrm{sep}}}\Delta_{1,I^*}f_{\mathtt w}\|_{L^1(|\nu_1|)}\, \|\chi_J \Delta_{2,J^*}g_{\mathtt w}\|_{L^1(|\nu_2|)}\\
					&\qquad \qquad\leq \frac{|\nu_1|(I_{\mathrm{sep}})^{1/2}\, |\nu_2|(J)^{1/2}}{\ell(I)^n}\|\chi_{I_{\mathrm{sep}}}\Delta_{1,I^*}f_{\mathtt w}\|_{L^2(|\nu_1|)}\, \|\chi_J \Delta_{2,J^*}g_{\mathtt w}\|_{L^2(|\nu_2|)}\\
					&\qquad\qquad \lesssim_m  \|\chi_{I_{\mathrm{sep}}}\Delta_{1,I^*}f_{\mathtt w}\|_{L^2(|\nu_1|)}\, \|\chi_J \Delta_{2,J^*}g_{\mathtt w}\|_{L^2(|\nu_2|)}\leq \|\chi_{I}\Delta_{1,I^*}f_{\mathtt w}\|_{L^2(|\nu_1|)}\, \|\chi_J \Delta_{2,J^*}g_{\mathtt w}\|_{L^2(|\nu_2|)},
				\end{split}
			\end{equation*}
			where we recall that $|\nu_1|(I)\lesssim \ell(I)^n$ and $|\nu_2|(J)\lesssim \ell(J)^n$ because $I$ and $J$ are children of transit cubes, see Remark \ref{remark:growth_children_transit_cubes}.
			
			Hence, the Cauchy-Schwarz inequality, \eqref{eq:counting_pairs}, and \eqref{eq:dec_f_mart_2} give 
			\begin{equation}
				\begin{split}
					&\sum_{(I,J)\in\mathfrak D(\mathtt w)}\bigl|\bigl\langle S_{|\nu_1|,\varepsilon}\bigl(\chi_{I_{\mathrm{sep}}} \Delta_{1,I^*}f_{\mathtt w}\bigr), \chi_J\Delta_{2,J^*}g_{\mathtt w}\bigr\rangle_*\bigr|\\
					&\qquad\qquad{\lesssim}_{\lambda, m} \sum_{(I,J)\in\mathfrak D(\mathtt w)}\|\chi_{I}\Delta_{1,I^*}f_{\mathtt w}\|_{L^2(|\nu_1|)}\, \|\chi_J \Delta_{2,J^*}g_{\mathtt w}\|_{L^2(|\nu_2|)}\lesssim \|f_{\mathtt w}\|_{L^2(|\nu_1|)}\, \|g_{\mathtt w}\|_{L^2(|\nu_2|)}.\qedhere
				\end{split}
			\end{equation}
		\end{proof}
		
		\vv
		Next, we analyze the contribution in \eqref{eq:key_sum_diagonal} corresponding to $I_{\partial}$, namely
		\[
		\mathcal S_{\partial}(\mathtt w) \coloneqq \sum_{(I,J)\in\mathfrak D(\mathtt w)}\bigl|\bigl\langle S_{|\nu_1|,\varepsilon}\bigl(\chi_{I_{\partial}} \Delta_{1,I^*}f_{\mathtt w}\bigr), \chi_J\Delta_{2,J^*}g_{\mathtt w}\bigr\rangle_*\bigr|.
		\]
	}
	
	For $\mathfrak m\in M_+(X)$ and a sequence $(\phi_j)_{j\in \mathbb Z}$ of $\mathfrak m$-measurable functions, we define
	\[
	\bigl\|(\phi_j)_{j\in \mathbb Z} \bigr\|_{L^2(\mathfrak m; \ell^2)}\coloneqq \Bigl\|\Bigl(\sum_{j\in \mathbb Z}|\phi_j|^2\Bigr)^{1/2}\Bigr\|_{L^2(\mathfrak m)}.
	\]
	
	The following estimate follows from the Cauchy-Schwarz inequality.
	\begin{lemma}\label{lemma:random_signs}
		If $\|(f_j)_{j\in \mathbb Z}\|_{L^2(|\nu_1|; \ell^2)}<\infty$ and $\|(g_j)_{j\in \mathbb Z}\|_{L^2(|\nu_2|; \ell^2)}<\infty$, then the following series converges absolutely, and
		\[
		\sum_{j\in \mathbb Z}\bigl|\bigl\langle S_{|\nu_1|,\varepsilon}f_j, g_j\bigr\rangle_*\bigr|\leq \|S_{|\nu_1|,\varepsilon}\|_{L^2(|\nu_1|)\to L^2(|\nu_2|)}\, \bigl\|(f_j)_{j\in \mathbb Z}\bigr\|_{L^2(|\nu_1|; \ell^2)}\, \bigl\|(g_j)_{j\in \mathbb Z}\bigr\|_{L^2(|\nu_2|; \ell^2)}.
		\]
	\end{lemma}
	\begin{proof}
		By the Cauchy-Schwarz inequality, for every $j\in\mathbb Z$ we have
		\[
		|\langle S_{|\nu_1|,\varepsilon}f_j, g_j\rangle_*|\leq\|S_{|\nu_1|,\varepsilon}\|_{L^2(|\nu_1|)\to L^2(|\nu_2|)}\,\|f_j\|_{L^2(|\nu_1|)}\,\|g_j\|_{L^2(|\nu_2|)}.
		\]
		The lemma follows by summing over $j$, by the Cauchy-Schwarz inequality in $\ell^2$, and by Tonelli's theorem, which gives $\sum_{j\in\mathbb Z}\|f_j\|^2_{L^2(|\nu_1|)}=\|(f_j)_{j\in\mathbb Z}\|^2_{L^2(|\nu_1|;\ell^2)}$, and the analogous identity for $(g_j)_{j\in\mathbb Z}$.
	\end{proof}
	
	The bound of $\mathcal  S_{\partial}(\mathtt w)$ also requires the control of the expectation of the $L^2(|\nu_1|)$-norm of a suitable square function.
	\begin{lemma}\label{lemma:lemma48}
		We have
		\[
		\begin{split}
			&\widetilde{\mathbb E}_{\mathtt w}\widehat{\mathbb E}
			\sum_{\substack{I\in\bbD(w_1)\\I^*\in\bbD_1^{\tr}}}
			\left\|\chi_{\widetilde I_\partial^{w_2,\xi}}
			\Delta_{1,I^*}f_{\mathtt w}\right\|_{L^2(|\nu_1|)}^2\leq
			\frac{C[(2m+1)\lambda^\eta+\upsilon]}
			{\mathbb P^{\Omega^2}(\mathscr G)}\,
			\|f\|_{L^2(|\nu_1|)}^2.
		\end{split}
		\]
		The constant $C$ includes the random-boundary and martingale constants.
	\end{lemma}
	\begin{proof}
		Let $w_1\in\Omega$, $I\in\bbD(w_1)$, and $x\in I$. If $x\in\partial_{\lambda\ell(J')}J'$ for some $J'\in\bbD(w_2)$, we denote by $J_x\in\bbD(w_2)$ the cube of the same generation as $J'$ that contains $x$. As observed after \eqref{eq:inclusion_I_partial}, we have $\dist(x,X\setminus J_x)\leq\lambda\ell(J')$, so $x$ belongs to the boundary zone $\partial(J_x,\lambda\ell(J'))$ of Theorem \ref{thm:HK56}. Since there are at most $(2m+1)$ generations of cubes $J'$ with $\ell(J')\sim_m\ell(I)$, Theorem \ref{thm:HK56} gives
		\begin{equation}\label{eq:small_probability_2}
			\mathbb P_{w_2}\bigl(\{x\in\widetilde I_{\partial,1}^{w_2}\}\bigr)\leq C_2\,(2m+1)\,\lambda^\eta.
		\end{equation}
		Moreover, $\mathcal K(I)$ contains at most $(2m+1)$ scales, and at each of them we use the failure probability $\upsilon/(2m+1)$ in Proposition \ref{proposition:random_balls}. Thus, by \eqref{eq:high_probability_balls} and the inclusion $B\subset\widehat B$, we have
		\begin{equation}\label{eq:expectation_first_part}
			\widehat{\mathbb P}\bigl(\{x\in\widetilde I_{\partial,2}^{\xi}\}\bigr)\leq\upsilon.
		\end{equation}
		
		We extend the martingale differences $\Delta_{1,Q}$ to every $w_1\in\Omega$ as in the proof of Lemma \ref{lemma:probability_bad_1}, so that the upper bound for the martingale differences in \eqref{eq:dec_f_mart_2} holds uniformly in $w_1$. By \eqref{eq:property_id_mart_diff} and Definition \ref{def:p_good_p_bad_truncated}, we have
		\[
		\Delta_{1,I^*}f_{\mathtt w}
		=\chi_{\{I^*\in\bbD_1^{\mathcal G}(\mathtt w),\ \ell(I^*)\geq\el^{K_1}\}}\Delta_{1,I^*}f.
		\]
		In particular, $|\Delta_{1,I^*}f_{\mathtt w}|\leq|\Delta_{1,I^*}f|$, and the function $\Delta_{1,I^*}f$ and the condition $I^*\in\bbD_1^{\tr}$ depend only on $w_1$.
		Therefore, by Tonelli's theorem, \eqref{eq:small_probability_2}, and \eqref{eq:expectation_first_part},
		\[
		\begin{split}
			&\widetilde{\mathbb E}_{\mathtt w}\widehat{\mathbb E}
			\sum_{I:\,I^*\in\bbD_1^{\tr}}\|\chi_{\widetilde I_\partial^{w_2,\xi}}\Delta_{1,I^*}f_{\mathtt w}\|_2^2\\
			&\quad\leq\frac{\mathbb E_{w_1}}{\mathbb P^{\Omega^2}(\mathscr G)}
			\sum_{I:\,I^*\in\bbD_1^{\tr}}\int_I
			\mathbb E_{w_2}\widehat{\mathbb E}\chi_{\widetilde I_\partial^{w_2,\xi}}(x)
			|\Delta_{1,I^*}f(x)|^2\,d|\nu_1|(x)\\
			&\quad\leq\frac{C[(2m+1)\lambda^\eta+\upsilon]}{\mathbb P^{\Omega^2}(\mathscr G)}
			\mathbb E_{w_1}\sum_{I:\,I^*\in\bbD_1^{\tr}}\|\chi_I\Delta_{1,I^*}f\|_2^2
			\lesssim\frac{(2m+1)\lambda^\eta+\upsilon}{\mathbb P^{\Omega^2}(\mathscr G)}\|f\|_2^2,
		\end{split}
		\]
		where the last inequality follows from \eqref{eq:dec_f_mart_2}.
	\end{proof}
	
	\vv
	We now prove the main estimate for $\mathcal  S_{\partial}(\mathtt w)$.
	\begin{lemma}\label{lemma:s_partial}
		We have
		\begin{equation}\label{eq:s_partial_1}
			\mathcal  S_{\partial}(\mathtt w)\lesssim_m (1+N_\varepsilon)\, \Bigl\|\Bigl(\sum_{I\in \bbD(w_1), I^*\in \bbD^{\tr}_1}|\chi_{\widetilde I^{w_2,\xi}_\partial}\,\Delta_{1,I^*}f_{\mathtt w}|^2\Bigr)^{1/2}\Bigr\|_{L^2(|\nu_1|)} \, \|g_{\mathtt w}\|_{L^2(|\nu_2|)}.
		\end{equation}
		
		Moreover, we have
		
		\begin{equation}\label{eq:s_partial_2}
			\widetilde{\mathbb E}_{\mathtt w}\widehat{\mathbb E}\mathcal S_{\partial}(\mathtt w)\lesssim_{C_2,m} \Bigl(\frac{(2m+1)\lambda^\eta+\upsilon}{\mathbb P^{\Omega^2}(\mathscr G)}\Bigr)^{1/2}\Bigl(1+N_\varepsilon\Bigr)\, \|f\|_{L^2(|\nu_1|)}\, \|g\|_{L^2(|\nu_2|)}.
		\end{equation}		
	\end{lemma}
	\begin{proof}
		Let $I\in\bbD(w_1)$. By \eqref{eq:inclusion_I_partial}, we have $I_\partial\subset\widetilde I_\partial^{w_2,\xi}$ for every $J$ such that $(I,J)\in\mathfrak D(\mathtt w)$.
		Hence, by \eqref{eq:localized_norm_comparison_sec5}, \eqref{eq:counting_pairs}, and the Cauchy-Schwarz inequality, we obtain
		\[
		\begin{split}
			\mathcal S_\partial(\mathtt w)
			&\leq\|S_\varepsilon\|_{\op}\sum_{(I,J)\in\mathfrak D(\mathtt w)}
			\bigl\|\chi_{\widetilde I_\partial^{w_2,\xi}}\Delta_{1,I^*}f_{\mathtt w}\bigr\|_{L^2(|\nu_1|)}\,
			\|\chi_J\Delta_{2,J^*}g_{\mathtt w}\|_{L^2(|\nu_2|)}\\
			&\lesssim_m(1+N_\varepsilon)\,
			\Bigl(\sum_{\substack{I\in\bbD(w_1)\\I^*\in\bbD_1^{\tr}}}
			\bigl\|\chi_{\widetilde I_\partial^{w_2,\xi}}\Delta_{1,I^*}f_{\mathtt w}\bigr\|_{L^2(|\nu_1|)}^2\Bigr)^{1/2}
			\Bigl(\sum_{\substack{J\in\bbD(w_2)\\J^*\in\bbD_2^{\tr}}}
			\|\chi_J\Delta_{2,J^*}g_{\mathtt w}\|_{L^2(|\nu_2|)}^2\Bigr)^{1/2}.
		\end{split}
		\]
		Then, the bound \eqref{eq:s_partial_1} follows by Tonelli's theorem and \eqref{eq:dec_f_mart_2}.
		
		Finally, \eqref{eq:s_partial_1}, \eqref{eq:dec_f_mart_2}, the Cauchy-Schwarz inequality, and Lemma \ref{lemma:lemma48} yield
		\[
		\begin{split}
			\widetilde{\mathbb E}_{\mathtt w}\widehat{\mathbb E}\mathcal S_\partial(\mathtt w)
			&\lesssim_m(1+N_\varepsilon)\,
			\Bigl(\widetilde{\mathbb E}_{\mathtt w}\widehat{\mathbb E}
			\sum_{\substack{I\in\bbD(w_1)\\I^*\in\bbD_1^{\tr}}}
			\bigl\|\chi_{\widetilde I_\partial^{w_2,\xi}}\Delta_{1,I^*}f_{\mathtt w}\bigr\|_{L^2(|\nu_1|)}^2\Bigr)^{1/2}
			\|g\|_{L^2(|\nu_2|)}\\
			&\lesssim_{C_2,m}\Bigl(\frac{(2m+1)\lambda^\eta+\upsilon}{\mathbb P^{\Omega^2}(\mathscr G)}\Bigr)^{1/2}
			(1+N_\varepsilon)\,\|f\|_{L^2(|\nu_1|)}\,\|g\|_{L^2(|\nu_2|)}.\qedhere
		\end{split}
		\]
		
	\end{proof}

	\vvv
	
	We are left with the study of the purely diagonal term
	\begin{equation}\label{eq:purely_diagonal_term}
		\mathcal S_{\Delta}(\mathtt w) \coloneqq \sum_{(I,J)\in\mathfrak D(\mathtt w)}\bigl|\bigl\langle S_{|\nu_1|,\varepsilon}\bigl(\chi_{I_{\Delta}} \Delta_{1,I^*}f_{\mathtt w}\bigr), \chi_J\Delta_{2,J^*}g_{\mathtt w}\bigr\rangle_*\bigr|.
	\end{equation}
	We further write
	\begin{equation*}
		\begin{split}
			&\mathcal S_{\Delta}(\mathtt w) \leq \sum_{(I,J)\in\mathfrak D(\mathtt w)}\bigl|\bigl\langle S_{|\nu_1|,\varepsilon}\bigl(\chi_{I_{\Delta}} \Delta_{1,I^*}f_{\mathtt w}\bigr), \chi_{J_\Delta}\Delta_{2,J^*}g_{\mathtt w}\bigr\rangle_*\bigr|\\
			&\qquad \qquad\qquad  + \sum_{(I,J)\in\mathfrak D(\mathtt w)}\bigl|\bigl\langle S_{|\nu_1|,\varepsilon}\bigl(\chi_{I_{\Delta}} \Delta_{1,I^*}f_{\mathtt w}\bigr), \chi_{J\setminus J_\Delta}\Delta_{2,J^*}g_{\mathtt w}\bigr\rangle_*\bigr|\eqqcolon \mathcal S_{\Delta,1}(\mathtt w)+ \mathcal S_{\Delta, 2}(\mathtt w).
		\end{split}
	\end{equation*}
	
	Since
	\[
	\mathcal S_{\Delta,2}= \sum_{(I,J)\in\mathfrak D(\mathtt w)}\bigl|\bigl\langle \chi_{I_{\Delta}} \Delta_{1,I^*}f_{\mathtt w}, S^*_{|\nu_2|,\varepsilon}\bigl( \chi_{J\setminus J_\Delta}\Delta_{2,J^*}g_{\mathtt w} \bigr) \bigr\rangle\bigr|, 
	\]
	we can estimate $\mathcal S_{\Delta,2}(\mathtt w)$ as $\mathcal S_{\mathrm{sep}}(\mathtt w)$ and $\mathcal S_{\partial}(\mathtt w)$, with the roles of the two grids interchanged, because $J\setminus J_\Delta=J_{\mathrm{sep}}\cup J_\partial$. More precisely, we obtain
	\begin{equation}\label{eq:s_delta_2}
		\begin{split}
			\widetilde{\mathbb E}_{\mathtt w}\widehat{\mathbb E}
			\mathcal S_{\Delta,2}(\mathtt w)
			\leq\biggl[C_{m,\lambda}
			+C_m\left(\frac{(2m+1)\lambda^\eta+\upsilon}
			{\mathbb P^{\Omega^2}(\mathscr G)}\right)^{1/2}
			(1+N_\varepsilon)\biggr]\|f\|_2\|g\|_2.
		\end{split}
	\end{equation}
	
	{
		We will use the following definition.
		\begin{definition}[Nice balls]\label{definition:nice_balls}
			A ball $B\in\widehat{\mathfrak B}_{\xi}$ is nice if $2B\cap F_1(w_1)\neq\varnothing$ and $2B\cap F_2(w_2)\neq\varnothing$.
		\end{definition}
		If $B$ is not nice, then $\Xi_{\mathtt w}\geq r(B)$ on $B$.
		If $B$ is nice and $5B\subset I\cap J$, then both $I$ and $J$ are transit, since $F_i(w_i)\subset X\setminus\mathscr W_i(w_i)$. If, in addition,
		\begin{equation}\label{eq:choice_of_constants}
			\lambda\min(\ell(I),\ell(J))
			\leq\frac{\theta}{5000A_0}\, \min(r(B_1),r(B_2)),
		\end{equation}
		then $2B\cap F_i(w_i)\neq\varnothing$ and $F_i(w_i)\subset(1-\theta)B_i$ imply
		\begin{equation}\label{eq:inclusion_choice_nice_balls}
			1000A_0B\subset(1-\theta/2)B_1\cap(1-\theta/2)B_2.
		\end{equation}
		
		{
			\subsubsection{Weak boundedness}\label{section:wbp}
			We first prove the weak boundedness property for the original testing measures $\widehat\nu_i$, and then apply the same argument to the localized operator. In the next two lemmas the pairings are taken with respect to $|\widehat\nu_2|$, and $\|T_{\Xi_{\mathtt w},\varepsilon}\|_{\op}$ denotes the norm from $L^2(|\widehat\nu_1|)$ to $L^2(|\widehat\nu_2|)$. By \eqref{eq:global_testing_before_localization_sec5}, increasing $c_\sharp$ by a structural constant if necessary, we have
			\[
			T_{\Xi_{\mathtt w},\sharp}\widehat\nu_1\leq c_\sharp\qquad \text{ and }\qquad 
			T^*_{\Xi_{\mathtt w},\sharp}\widehat\nu_2\leq c_\sharp\qquad\text{on }X.
			\]
			
			Before proving the weak boundedness lemmas, we state the following auxiliary result.
			\begin{lemma}\label{lemma:main_technical_lemma_wbp}
				There exists $\mathbf c>0$, depending only on the structural constants, such that, for every ball $B$ and $a\in(0,1/4)$,
				\begin{equation}\label{eq:wbp_I_2_1}
					\sup_{y\in B}\big|T_{\Xi_{\mathtt w},|\widehat\nu_1|,\varepsilon}\bigl(\chi_{((1+a)B)^c}\widehat b_1\bigr)(y)\big|\leq\mathbf c\,a^{-n}+c_\sharp
				\end{equation}
				and
				\begin{equation}\label{eq:wbp_I_2_2}
					\sup_{y\in B}\big|T^*_{\Xi_{\mathtt w},|\widehat\nu_2|,\varepsilon}\bigl(\chi_{((1+a)B)^c}\widehat b_2\bigr)(y)\big|\leq\mathbf c\,a^{-n}+c_\sharp.
				\end{equation}
			\end{lemma}
			\begin{proof}
				Write $r=r(B)$ and fix $y\in B$. By \eqref{eq:growth_nu_sec6},
				\begin{equation}\label{eq:growht_nu_1_1ballssec7}
					|\widehat\nu_i|(B(y,R))\lesssim R^n,\qquad \text{ for }R>\Xi_{\mathtt w}(y), i=1,2.
				\end{equation}
				Since $(1+a)B\subset B(y,3r)$, we have
				\begin{equation}\label{eq:splitting_xy_lemma72}
					\begin{split}
						T_{\Xi_{\mathtt w},|\widehat\nu_1|,\varepsilon}(\chi_{((1+a)B)^c}\widehat b_1)(y)
						&=T_{\Xi_{\mathtt w},\max\{\varepsilon,3r\}}\widehat\nu_1(y)\\
						&\quad+\int_{((1+a)B)^c\cap\{z:\,\varepsilon<d(y,z)\leq3r\}}K_{\Xi_{\mathtt w}}(y,z)\,d\widehat\nu_1(z).
					\end{split}
				\end{equation}
				The first term is bounded by $c_\sharp$, by \eqref{eq:global_testing_before_localization_sec5}. To bound the second term, note that $d(y,z)\geq ar$ on $((1+a)B)^c$. Therefore, the suppressed kernel estimate \eqref{eq:bound_suppressed_kernel_1} gives
				\begin{equation}\label{eq:holder_continuity_ball_enlarged}
					\begin{split}
						\int_{((1+a)B)^c\cap B(y,4r)}|K_{\Xi_{\mathtt w}}(y,z)|\,d|\widehat\nu_1|(z)
						&\lesssim\frac{|\widehat\nu_1|(B(y,4\max\{r,\Xi_{\mathtt w}(y)\}))}{\max\{ar,\Xi_{\mathtt w}(y)\}^n}\\
						&\overset{\eqref{eq:growht_nu_1_1ballssec7}}{\lesssim}\frac{\max\{r,\Xi_{\mathtt w}(y)\}^n}{\max\{ar,\Xi_{\mathtt w}(y)\}^n}\leq a^{-n}.
					\end{split}
				\end{equation}
				This proves \eqref{eq:wbp_I_2_1}. The estimate \eqref{eq:wbp_I_2_2} follows by the same argument applied to $T^*$, with the roles of $\widehat\nu_1$ and $\widehat\nu_2$ interchanged.
			\end{proof}
			
			\begin{lemma}[Key suppressed WBP]\label{lemma:key_supp_WBP_2}
				Let $a\in(0,1/4)$, and let $B\subset X$ be a ball with $M$-small boundary with respect to both $|\widehat\nu_1|$ and $|\widehat\nu_2|$. Then, there exists $\mathfrak C_{wb}>0$ such that
				\begin{equation}\label{eq:weak_boundedness_property_claim}
					\begin{split}
						\big|\langle T_{\Xi_{\mathtt w},|\widehat\nu_1|,\varepsilon}
						(\chi_B\widehat b_1),\chi_B\widehat b_2
						\rangle_{|\widehat\nu_2|}\big|
						&\leq a^{-n}\mathfrak C_{wb}
						|\widehat\nu_1|(B)^{1/2}|\widehat\nu_2|(B)^{1/2}\\
						&\quad+\mathfrak C_{wb}(Ma)^{1/2}
						\|T_{\Xi_{\mathtt w},\varepsilon}\|_{\op}
						|\widehat\nu_1|(5B)^{1/2}|\widehat\nu_2|(5B)^{1/2},
					\end{split}
				\end{equation}
				where $\mathfrak C_{wb}$ depends on the structural constants from the statement of Theorem \ref{theorem:main_theorem_hofmann_problem}, but not on $a$, $M$, $\mathtt w$, or $\varepsilon$.
			\end{lemma}
			
			\begin{proof}
				
				For $B$ and $a$ as in the statement, we write
				\begin{equation}\label{eq:split_wbp_I_1_I_2}
					\begin{split}
						\big|\langle T_{\Xi_{\mathtt w},|\widehat\nu_1|,\varepsilon}(\chi_B\widehat b_1),\chi_B\widehat b_2\rangle_{|\widehat\nu_2|}\big|
						&\leq\big|\langle T_{\Xi_{\mathtt w},|\widehat\nu_1|,\varepsilon}\widehat b_1,\chi_B\widehat b_2\rangle_{|\widehat\nu_2|}\big|\\
						&\quad+\big|\langle T_{\Xi_{\mathtt w},|\widehat\nu_1|,\varepsilon}(\chi_{B^c}\widehat b_1),\chi_B\widehat b_2\rangle_{|\widehat\nu_2|}\big|\eqqcolon I_1+I_2.
					\end{split}
				\end{equation}
				The term $I_1$ is controlled by \eqref{eq:global_testing_before_localization_sec5} and, since $T_{\Xi_{\mathtt w},|\widehat\nu_1|,\varepsilon}\widehat b_1=T_{\Xi_{\mathtt w},\varepsilon}\widehat\nu_1$, we obtain
				\begin{equation}\label{eq:wbp_1}
					I_1\leq\int_B|T_{\Xi_{\mathtt w},\varepsilon}\widehat\nu_1|\,d|\widehat\nu_2|\leq c_\sharp|\widehat\nu_2|(B).
				\end{equation}
				To deal with the second term, we write
				\[
				\begin{split}
					I_2&\leq\big|\langle T_{\Xi_{\mathtt w},|\widehat\nu_1|,\varepsilon}(\chi_{(1+a)B\setminus B}\widehat b_1),\chi_B\widehat b_2\rangle_{|\widehat\nu_2|}\big|\\
					&\qquad+\big|\langle T_{\Xi_{\mathtt w},|\widehat\nu_1|,\varepsilon}(\chi_{((1+a)B)^c}\widehat b_1),\chi_B\widehat b_2\rangle_{|\widehat\nu_2|}\big|\eqqcolon I'_2+I''_2.
				\end{split}
				\]
				Since $B$ has $M$-$|\widehat\nu_1|$-small boundary, we have
				\begin{equation}\label{eq:wbp_21}
					\begin{split}
						I'_2&\leq\|T_{\Xi_{\mathtt w},\varepsilon}\|_{\op}\|\chi_{(1+a)B\setminus B}\widehat b_1\|_{L^2(|\widehat\nu_1|)}\|\chi_B\widehat b_2\|_{L^2(|\widehat\nu_2|)}\\
						&\leq a^{1/2}M^{1/2}\|T_{\Xi_{\mathtt w},\varepsilon}\|_{\op}|\widehat\nu_1|(5B)^{1/2}|\widehat\nu_2|(5B)^{1/2}.
					\end{split}
				\end{equation}
				On the other hand, we bound the second summand by
				\begin{equation}\label{eq:wbp_22}
					I''_2\leq\sup_{y\in B}\big|T_{\Xi_{\mathtt w},|\widehat\nu_1|,\varepsilon}(\chi_{((1+a)B)^c}\widehat b_1)(y)\big|\,|\widehat\nu_2|(B).
				\end{equation}
				
				Gathering \eqref{eq:split_wbp_I_1_I_2}, \eqref{eq:wbp_1}, \eqref{eq:wbp_21}, \eqref{eq:wbp_22}, and Lemma \ref{lemma:main_technical_lemma_wbp}, we obtain, for some $\mathbf c'>0$ depending only on the structural constants,
				\begin{equation}\label{eq:wbp_av1}
					\begin{split}
						\big|\langle T_{\Xi_{\mathtt w},|\widehat\nu_1|,\varepsilon}(\chi_B\widehat b_1),\chi_B\widehat b_2\rangle_{|\widehat\nu_2|}\big|
						&\leq\mathbf c' a^{-n}|\widehat\nu_2|(B)\\
						&\quad+a^{1/2}M^{1/2}\|T_{\Xi_{\mathtt w},\varepsilon}\|_{\op}|\widehat\nu_1|(5B)^{1/2}|\widehat\nu_2|(5B)^{1/2}.
					\end{split}
				\end{equation}
				Since the operator norms of $T_{\Xi_{\mathtt w},|\widehat\nu_1|,\varepsilon}$ and $T^*_{\Xi_{\mathtt w},|\widehat\nu_2|,\varepsilon}$ coincide, we can repeat the argument for $T^*$ and obtain
				\begin{equation}\label{eq:wbp_av2}
					\begin{split}
						\big|\langle T_{\Xi_{\mathtt w},|\widehat\nu_1|,\varepsilon}(\chi_B\widehat b_1),\chi_B\widehat b_2\rangle_{|\widehat\nu_2|}\big|
						&\leq\mathbf c'a^{-n}|\widehat\nu_1|(B)\\
						&\quad+a^{1/2}M^{1/2}\|T_{\Xi_{\mathtt w},\varepsilon}\|_{\op}|\widehat\nu_1|(5B)^{1/2}|\widehat\nu_2|(5B)^{1/2}.
					\end{split}
				\end{equation}
				Taking the smaller of \eqref{eq:wbp_av1} and
				\eqref{eq:wbp_av2}, and using
				\[
				\min\{|\widehat\nu_1|(B),|\widehat\nu_2|(B)\}
				\leq |\widehat\nu_1|(B)^{1/2}|\widehat\nu_2|(B)^{1/2},
				\]
				proves \eqref{eq:weak_boundedness_property_claim}.
			\end{proof}
			
			We now prove the analogue of Lemma \ref{lemma:key_supp_WBP_2} for the localized operator $S_{|\nu_1|,\varepsilon}$.
			\begin{lemma}[Weak boundedness for the localized operator]
				\label{lemma:key_supp_WBP_in_text}
				Let $a\in(0,1/4)$, and let $B\subset X$ be a ball with $M_\Lambda$-small boundary with respect to both $|\nu_1|$ and $|\nu_2|$. Then, we have
				\begin{equation}\label{eq:localized_weak_boundedness_property_sec5}
					\begin{split}
						&\bigl|\bigl\langle S_{|\nu_1|,\varepsilon}(\chi_B\widehat b_1),\chi_B\widehat b_2\bigr\rangle_*\bigr|\\
						&\quad\leq C a^{-n}|\nu_1|(B)^{1/2}|\nu_2|(B)^{1/2}
						+C(M_\Lambda a)^{1/2}(1+N_\varepsilon)
						|\nu_1|(5B)^{1/2}|\nu_2|(5B)^{1/2},
					\end{split}
				\end{equation}
				where $C>0$ depends only on the structural constants, and not on $a$, $M_\Lambda$, $B$, $\mathtt w$, or $\varepsilon$. 
			\end{lemma}
			\begin{proof}
				To simplify the notation, we write $S\coloneqq S_{|\nu_1|,\varepsilon}$ and $r\coloneqq r(B)$, and we denote by $x_B$ the center of $B$. We split
				\begin{equation}\label{eq:splitting_localized_wbp}
					\begin{split}
						\bigl\langle S(\chi_B\widehat b_1),\chi_B\widehat b_2\bigr\rangle_*
						&=\bigl\langle S\widehat b_1,\chi_B\widehat b_2\bigr\rangle_*
						-\bigl\langle S(\chi_{(1+a)B\setminus B}\widehat b_1),\chi_B\widehat b_2\bigr\rangle_*\\
						&\quad-\bigl\langle S(\chi_{((1+a)B)^c}\widehat b_1),\chi_B\widehat b_2\bigr\rangle_*
						\eqqcolon\mathfrak I_1-\mathfrak I_2-\mathfrak I_3.
					\end{split}
				\end{equation}
				Since $L_\varepsilon(x,y)=0$ if $d(x,y)\leq\varepsilon$, we have $S\widehat b_1=S_{|\nu_1|,\varepsilon,\rho}\widehat b_1$ for $\rho<\varepsilon$. Hence, since $|\widehat b_2|=1$ $|\nu_2|$-almost everywhere, \eqref{eq:localized_maximal_testing_sec5} gives
				\[
				|\mathfrak I_1|\leq\int_B|S\widehat b_1|\,d|\nu_2|\leq c_\sharp\,|\nu_2|(B).
				\]
				By \eqref{eq:def_small_boundary} with $\lambda=a$, we have $|\nu_1|((1+a)B\setminus B)\leq M_\Lambda a\,|\nu_1|(5B)$. Thus, \eqref{eq:localized_norm_comparison_sec5} yields
				\[
				|\mathfrak I_2|\leq\|S_\varepsilon\|_{\op}\,|\nu_1|((1+a)B\setminus B)^{1/2}\,|\nu_2|(B)^{1/2}
				\lesssim(M_\Lambda a)^{1/2}(1+N_\varepsilon)\,|\nu_1|(5B)^{1/2}|\nu_2|(5B)^{1/2}.
				\]
				To bound $\mathfrak I_3$, we fix $y\in B$. Since $(1+a)B\subset B(y,3r)$ and $d(y,z)>ar$ for $z\notin(1+a)B$, we have
				\[
				S(\chi_{((1+a)B)^c}\widehat b_1)(y)
				=S_{|\nu_1|,\varepsilon,3r}\widehat b_1(y)
				+\int_{((1+a)B)^c\cap\{z:\,d(y,z)\leq3r\}}L_\varepsilon(y,z)\,\widehat b_1(z)\,d|\nu_1|(z),
				\]
				and, by \eqref{eq:size_auxiliary_kernel_sec6} and \eqref{eq:growth_auxiliary_measures_sec6}, the absolute value of the last integral is bounded by a constant times
				\[
				\frac{|\nu_1|(B(y,4r))}{\max\{ar,\Xi_{\mathtt w}(y)\}^n}
				\leq(\mathbf m+\mathfrak C_0)\,\frac{(4r+\Xi_{\mathtt w}(y))^n}{\max\{ar,\Xi_{\mathtt w}(y)\}^n}
				\lesssim a^{-n}.
				\]
				Therefore, by \eqref{eq:localized_maximal_testing_sec5}, we obtain
				\[
				|\mathfrak I_3|\leq\sup_{y\in B}\bigl|S(\chi_{((1+a)B)^c}\widehat b_1)(y)\bigr|\,|\nu_2|(B)
				\lesssim(c_\sharp+a^{-n})\,|\nu_2|(B).
				\]
				Since $a<1$, we gather the bounds for the terms in \eqref{eq:splitting_localized_wbp} and obtain
				\begin{equation}\label{eq:localized_wbp_one_side}
					\bigl|\bigl\langle S(\chi_B\widehat b_1),\chi_B\widehat b_2\bigr\rangle_*\bigr|
					\leq C a^{-n}|\nu_2|(B)
					+C(M_\Lambda a)^{1/2}(1+N_\varepsilon)\,|\nu_1|(5B)^{1/2}|\nu_2|(5B)^{1/2}.
				\end{equation}
				Moreover, $\langle S(\chi_B\widehat b_1),\chi_B\widehat b_2\rangle_*=\langle\chi_B\widehat b_1,S^*_{|\nu_2|,\varepsilon}(\chi_B\widehat b_2)\rangle$, and the operators $S$ and $S^*_{|\nu_2|,\varepsilon}$ have the same norm. Thus, the same argument applied to $S^*_{|\nu_2|,\varepsilon}$, with the second bound in \eqref{eq:localized_maximal_testing_sec5} and the $M_\Lambda$-small boundary of $B$ with respect to $|\nu_2|$, gives \eqref{eq:localized_wbp_one_side} with $|\nu_1|(B)$ in place of $|\nu_2|(B)$. Since $\min\{|\nu_1|(B),|\nu_2|(B)\}\leq|\nu_1|(B)^{1/2}|\nu_2|(B)^{1/2}$, this proves \eqref{eq:localized_weak_boundedness_property_sec5}.
			\end{proof}
		}
		
	}
	\vv
	
	We recall that, in the following lemma, the finite depth and separation constants of the almost-covering
	families depend on $m$ and $\upsilon$.
	\begin{lemma}\label{lemma:S_delta_1}
		There exist $C_{m,\lambda,\upsilon},C_{m,\upsilon}>0$, independent of $a$, $\varepsilon$, $\mathtt w$, and $\xi$, such that
		\[
		|\mathcal S_{\Delta,1}(\mathtt w)|
		\leq\left[C_{m,\lambda,\upsilon}(1+a^{-n})
		+C_{m,\upsilon}(M_\Lambda a)^{1/2}(1+N_\varepsilon)\right]
		\|f\|_{L^2(|\nu_1|)}\|g\|_{L^2(|\nu_2|)}.
		\]
		Here, $M_\Lambda$ is the small-boundary constant of the balls in $\widehat{\mathfrak B}_{\xi}$.
	\end{lemma}
	\begin{proof}
		Using $I_\Delta=\bigcup_iB_{I,i}$ and $J_\Delta=\bigcup_jB_{J,j}$, we split
		\[
		\mathcal S_{\Delta,1}(\mathtt w)\leq\sum_{(I,J)\in\mathfrak D(\mathtt w)}\sum_{i,j}
		\bigl|\bigl\langle S_{|\nu_1|,\varepsilon}(\chi_{B_{I,i}}\Delta_{1,I^*}f_{\mathtt w}),\chi_{B_{J,j}}\Delta_{2,J^*}g_{\mathtt w}\bigr\rangle_*\bigr|
		\eqqcolon\mathfrak S'_\Delta+\mathfrak S''_\Delta+\mathfrak S'''_\Delta,
		\]
		where $\mathfrak S'_\Delta$, $\mathfrak S''_\Delta$, and $\mathfrak S'''_\Delta$ are the sums over the pairs of pieces coming from distinct balls $B,B'\in\widehat{\mathfrak B}_{\xi}$, from the same ball $B$ that is not nice, and from the same nice ball $B$, respectively. By \eqref{eq:radii_random_balls} and \eqref{eq:choice_k_random_balls}, every $B\in\widehat{\mathfrak B}_{\xi}$ satisfies
		\begin{equation}\label{eq:bd_from_below_radius_random_ball}
			r(B)\geq c_{m,\upsilon}\vartheta^k
			\geq c_{m,\upsilon}\vartheta\lambda\, \min(\ell(I),\ell(J)).
		\end{equation}
		
		If $B\neq B'$, then \eqref{eq:separation_random_balls} and \eqref{eq:bd_from_below_radius_random_ball} give
		\begin{equation}\label{eq:S_delta_1_3}
			\operatorname{dist}(B_{I,i},B_{J,j})
			\geq c_{m,\lambda,\upsilon}\, \min(\ell(I),\ell(J)).
		\end{equation}
		Since the balls in $\widehat{\mathfrak B}_{\xi}$ are pairwise disjoint and $\ell(I)\sim_m\ell(J)$, \eqref{eq:S_delta_1_3} and \eqref{eq:size_auxiliary_kernel_sec6} imply that, for fixed $(I,J)$, the sum of these terms is at most $C_{m,\lambda,\upsilon}\,\ell(I)^{-n}\|\chi_I\Delta_{1,I^*}f_{\mathtt w}\|_{L^1(|\nu_1|)}\|\chi_J\Delta_{2,J^*}g_{\mathtt w}\|_{L^1(|\nu_2|)}$. Hence, the argument in the proof of Lemma \ref{lemma:sep_transit_transit} and \eqref{eq:dec_f_mart_2} give
		\begin{equation}\label{eq:boundSdelta1}
			\mathfrak S'_\Delta\leq C_{m,\lambda,\upsilon}\|f\|_{L^2(|\nu_1|)}\|g\|_{L^2(|\nu_2|)}.
		\end{equation}
		
		If $B=B'$, then $5B\subset I\cap J$ and both pieces coincide with $B$. If, in addition, $B$ is not nice, then $2B\cap F_i(w_i)=\varnothing$ for some $i\in\{1,2\}$, so $\Xi_{\mathtt w}\geq\operatorname{dist}(\cdot,F_i(w_i))\geq r(B)$ on $B$, and \eqref{eq:size_auxiliary_kernel_sec6} gives $|L_\varepsilon|\lesssim r(B)^{-n}$ on $B\times B$. Moreover, $|\nu_1|(I)\lesssim\ell(I)^n$ and $|\nu_2|(J)\lesssim\ell(J)^n$ by Remark \ref{remark:growth_children_transit_cubes}. Thus, the Cauchy-Schwarz inequality and \eqref{eq:bd_from_below_radius_random_ball} yield
		\begin{equation}\label{eq:S_delta_1_4}
			\bigl|\bigl\langle S_{|\nu_1|,\varepsilon}(\chi_B\Delta_{1,I^*}f_{\mathtt w}),\chi_B\Delta_{2,J^*}g_{\mathtt w}\bigr\rangle_*\bigr|
			\leq C_{m,\lambda,\upsilon}\,
			\|\chi_B\Delta_{1,I^*}f_{\mathtt w}\|_{L^2(|\nu_1|)}
			\|\chi_B\Delta_{2,J^*}g_{\mathtt w}\|_{L^2(|\nu_2|)}.
		\end{equation}
		Since the balls are disjoint and each cube occurs in at most $C_m$ pairs in $\mathfrak D(\mathtt w)$, the Cauchy-Schwarz inequality and \eqref{eq:dec_f_mart_2} give
		\begin{equation}\label{eq:boundSdelta2}
			\mathfrak S''_\Delta\leq C_{m,\lambda,\upsilon}\|f\|_{L^2(|\nu_1|)}\|g\|_{L^2(|\nu_2|)}.
		\end{equation}
		
		If $B$ is nice, then $2B$ meets $F_1(w_1)$ and $F_2(w_2)$. Since $2B\subset I\cap J$ and $F_i(w_i)\cap\mathscr W_i(w_i)=\varnothing$, both $I$ and $J$ are transit, so, by \eqref{eq:def_martingales}, $\Delta_{1,I^*}f_{\mathtt w}=D^1_If_{\mathtt w}\,\widehat b_1$ on $I$ and $\Delta_{2,J^*}g_{\mathtt w}=D^2_Jg_{\mathtt w}\,\widehat b_2$ on $J$. Hence, we have
		\begin{equation}\label{eq:S_delta_1_2}
			\bigl\langle S_{|\nu_1|,\varepsilon}(\chi_B\Delta_{1,I^*}f_{\mathtt w}),\chi_B\Delta_{2,J^*}g_{\mathtt w}\bigr\rangle_*
			=D^1_If_{\mathtt w}\,D^2_Jg_{\mathtt w}\,\bigl\langle S_{|\nu_1|,\varepsilon}(\chi_B\widehat b_1),\chi_B\widehat b_2\bigr\rangle_*,
		\end{equation}
		and Lemma \ref{lemma:key_supp_WBP_in_text} applies to the last pairing.
		The balls in $\widehat{\mathfrak B}_{\xi}$ have comparable radii and satisfy \eqref{eq:separation_random_balls}, so, by geometric doubling, the dilates $5B$ have overlap bounded by $C_{m,\upsilon}$. Since $5B\subset I\cap J$, the Cauchy-Schwarz inequality over the common nice balls gives
		\[
		\sum_{B}\bigl(|\nu_1|(B)^{1/2}|\nu_2|(B)^{1/2}+|\nu_1|(5B)^{1/2}|\nu_2|(5B)^{1/2}\bigr)
		\leq C_{m,\upsilon}|\nu_1|(I\cap J)^{1/2}|\nu_2|(I\cap J)^{1/2}.
		\]
		For fixed $I$, the sets $I\cap J$ are pairwise disjoint within each of the $2m+1$ generations of the cubes $J$, and the same holds with $I$ and $J$ interchanged. Moreover, $|D^1_If_{\mathtt w}|^2|\nu_1|(I)=\|\chi_I\Delta_{1,I^*}f_{\mathtt w}\|_{L^2(|\nu_1|)}^2$, and analogously for $J$. Therefore, the Cauchy-Schwarz inequality and \eqref{eq:dec_f_mart_2} give
		\[
		\begin{split}
			&\sum_{(I,J)\in\mathfrak D^{\tr}_{\tr}(\mathtt w)}
			|D^1_I f_{\mathtt w}D^2_Jg_{\mathtt w}|\,
			|\nu_1|(I\cap J)^{1/2}|\nu_2|(I\cap J)^{1/2}\\
			&\qquad\leq (2m+1)
			\Bigl(\sum_I|D^1_I f_{\mathtt w}|^2|\nu_1|(I)\Bigr)^{1/2}
			\Bigl(\sum_J|D^2_J g_{\mathtt w}|^2|\nu_2|(J)\Bigr)^{1/2}
			\lesssim_m\|f\|_{L^2(|\nu_1|)}\|g\|_{L^2(|\nu_2|)},
		\end{split}
		\]
		which, combined with \eqref{eq:S_delta_1_2}, Lemma \ref{lemma:key_supp_WBP_in_text}, and the previous bound, yields
		\begin{equation}\label{eq:boundSdelta3}
			\mathfrak S'''_\Delta
			\leq C_{m,\upsilon}
			\bigl[a^{-n}+(M_\Lambda a)^{1/2}(1+N_\varepsilon)\bigr]\,\|f\|_{L^2(|\nu_1|)}\|g\|_{L^2(|\nu_2|)}.
		\end{equation}
		Gathering \eqref{eq:boundSdelta1}, \eqref{eq:boundSdelta2}, and \eqref{eq:boundSdelta3}, we conclude the proof.
	\end{proof}
	
	\vv
	Now, we write
	\[
	\rho_m(\lambda,\upsilon)\coloneqq C_m\, \Bigl(\frac{(2m+1)\lambda^\eta+\upsilon}{\mathbb P^{\Omega^2}(\mathscr G)}\Bigr)^{1/2}
	\]
	and, combining the preceding lemma with \eqref{eq:s_delta_2}, we obtain
	\begin{equation}\label{eq:final_estimate_S_delta}
		\begin{split}
			\widetilde{\mathbb E}_{\mathtt w}\widehat{\mathbb E}
			\mathcal S_\Delta(\mathtt w)
			\leq\bigl[&C_{m,\lambda,\upsilon}(1+a^{-n})\\
			&\qquad\qquad+\bigl(C_{m,\upsilon}(M_\Lambda a)^{1/2}
			+\rho_m(\lambda,\upsilon)\bigr)(1+N_\varepsilon)\bigr]
			\|f\|_2\|g\|_2.
		\end{split}
	\end{equation}
	
	\begin{proposition}\label{prop:diagonal_terms}
		There are constants $C_{m,\lambda,\upsilon}$ and $C_{m,\upsilon}$,
		independent of $a$ and $\varepsilon$, such that
		\[
		\begin{split}
			\widetilde{\mathbb E}_{\mathtt w}
			|\mathcal S_2(\mathtt w)+\mathcal S_3(\mathtt w)|
			\leq\bigl[&C_{m,\lambda,\upsilon}(1+a^{-n})\\
			&\qquad+\bigl(C_{m,\upsilon}(M_\Lambda a)^{1/2}
			+\rho_m(\lambda,\upsilon)\bigr)(1+N_\varepsilon)\bigr]
			\|f\|_2\|g\|_2.
		\end{split}
		\]
	\end{proposition}
	\begin{proof}
		By \eqref{eq:bound_s_2_s_3}, \eqref{eq:key_sum_diagonal}, and the decomposition $I=I_{\mathrm{sep}}\cup I_\partial\cup I_\Delta$, for every value of $\xi$ we have $|\mathcal S_2(\mathtt w)+\mathcal S_3(\mathtt w)|\leq\mathcal S_{\mathrm{sep}}(\mathtt w)+\mathcal S_\partial(\mathtt w)+\mathcal S_\Delta(\mathtt w)$. Since the left-hand side does not depend on $\xi$, the proposition follows by taking expectations and applying Lemma \ref{lemma:sep_transit_transit}, \eqref{eq:dec_f_mart_2}, \eqref{eq:s_partial_2}, and \eqref{eq:final_estimate_S_delta}.
	\end{proof}

	\vv
	
	\subsection*{{Estimate of $\mathcal S_4$}}
	We split
	\begin{equation}\label{eq:splitting_s_4}
		\begin{split}
			\mathcal S_4(\mathtt w)&= \sum_{\substack{Q\in \bbD^{\tr}_1, R \in \bbD^{\tr}_2\\Q\cap R\neq \varnothing \\ \ell(Q)<\el^{m}\ell(R)}}\bigl\langle S_{|\nu_1|,\varepsilon}(\Delta_{1,Q} f_{\mathtt w}), \Delta_{2,R} g_{\mathtt w}\bigr\rangle_* \\
			&\qquad\qquad+ \sum_{\substack{Q\in \bbD^{\tr}_1, R \in \bbD^{\tr}_2\\Q\cap R\neq \varnothing\\ \ell(Q)>\el^{-m}\ell(R)}}\bigl\langle S_{|\nu_1|,\varepsilon}(\Delta_{1,Q} f_{\mathtt w}), \Delta_{2,R} g_{\mathtt w}\bigr\rangle_* \eqqcolon \mathcal S'_4+ \mathcal S''_4.
		\end{split}
	\end{equation}
	Note that it suffices to bound the first sum since, by  \eqref{eq:transpose_martingale_difference}, the second one can be estimated by the same arguments for $S^*_{|\nu_2|, \varepsilon}$ instead of $S_{|\nu_1|, \varepsilon}$.
	
	Let $Q\in \bbD^{\tr}_1$ and $R\in\bbD^{\tr}_2$ be such that $Q\cap R\neq \varnothing$ and $\ell(Q)<\el^{m}\ell(R)$.
	We also recall that, by construction of $f_{\mathtt w}$ and $g_{\mathtt w}$, the cubes $Q$ and $R$ in \eqref{eq:splitting_s_4} are good.
	Since the generations are discrete,
	$\ell(Q)\leq\el^{m+1}\ell(R)$. For each child $U$ of $R$, goodness
	therefore gives either
	$\operatorname{dist}(Q,U)>\ell(Q)^\alpha\ell(U)^{1-\alpha}$ or
	$\operatorname{dist}(Q,X\setminus U)>\ell(Q)^\alpha\ell(U)^{1-\alpha}$.
	As $Q$ meets $R$, exactly one child, denoted by $R_Q$, meets $Q$; it
	contains all of $Q$ and
	\[
	\operatorname{dist}(Q,X\setminus R_Q)>\ell(Q)^\alpha\ell(R_Q)^{1-\alpha}.
	\]
	We use this notation and define the families
	\[
	\mathscr S'_{4, \tr}\coloneqq \bigl\{(Q, R)\in \bbD^{\tr}_1\times \bbD^{\tr}_2: Q\subset R, \ell(Q)<\el^{m}\ell(R), R_Q \text{ transit}\bigr\}	
	\]
	and 
	\[
	\mathscr S'_{4,\term}\coloneqq \bigl\{(Q, R)\in \bbD^{\tr}_1\times \bbD^{\tr}_2: Q\subset R, \ell(Q)<\el^{m}\ell(R), R_Q \text{ terminal}\bigr\}.
	\]
	
	Hence, we further split the sum $\mathcal S'_4$ as 
	\[
	\mathcal S'_4= \sum_{(Q,R)\in \mathscr S'_{4, \tr}} \qquad + \sum_{(Q,R)\in \mathscr S'_{4, \term}}\eqqcolon \mathcal S'_{4,\tr} + \mathcal S'_{4,\term}.
	\]
	
	\vv
	\subsubsection*{\underline{The sum over $\mathscr S'_{4, \tr}$}}
	To bound the contribution of transit cubes, we apply \eqref{eq:transpose_martingale_difference} and write 
	\begin{equation}\label{eq:s_4_1_tr}
		\begin{split}
			\mathcal S'_{4,\tr}&= 	\sum_{(Q,R)\in \mathscr S'_{4, \tr}}\bigl\langle\Delta_{1,Q} f_{\mathtt w},S^*_{|\nu_2|,\varepsilon}\bigl(\chi_{R_Q}\Delta_{2,R} g_{\mathtt w}\bigr)\bigr\rangle \\
			&\qquad\qquad\qquad + \sum_{(Q,R)\in \mathscr S'_{4, \tr}}\bigl\langle\Delta_{1,Q} f_{\mathtt w},S^*_{|\nu_2|,\varepsilon}\bigl(\chi_{X\setminus R_Q}\Delta_{2,R} g_{\mathtt w}\bigr)\bigr\rangle .
		\end{split}
	\end{equation}
	
	To bound the far away term, we note that \[\dist(Q,X\setminus R_Q)>\ell(Q)^\alpha\ell(R_Q)^{1-\alpha}=\el^{1-\alpha}\ell(Q)^\alpha\ell(R)^{1-\alpha}\] and that $D(Q,R)\approx\ell(R)$, since $Q\subset R$. Hence, as in the estimate of $\mathcal S_1'$, see also \cite[Lemma 5.17]{To14}, Lemma \ref{lemma:lem_514} with $\mathfrak c'=\el^{1-\alpha}$, \eqref{eq:sec5_two_measure_schur} combined with Schur's test, and \eqref{eq:dec_f_mart_2} give
	\begin{equation}\label{eq:S_4_1_far_away_part}
		\sum_{(Q,R)\in \mathscr S'_{4, \tr}}\bigl|\bigl\langle\Delta_{1,Q} f_{\mathtt w},S^*_{|\nu_2|,\varepsilon}\bigl(\chi_{X\setminus R_Q}\Delta_{2,R} g_{\mathtt w}\bigr)\bigr\rangle\bigr|\lesssim \|f_{\mathtt w}\|_{L^2(|\nu_1|)}\, \|g_{\mathtt w}\|_{L^2(|\nu_2|)}.
	\end{equation}

	{
		To estimate the first sum in \eqref{eq:s_4_1_tr}, we define
		\[
		c_{R,Q}(g_{\mathtt w})\coloneqq \frac{\langle g_{\mathtt w}\rangle_{2,R_Q}}{\langle \widehat b_2\rangle_{2,R_Q}}- \frac{\langle g_{\mathtt w}\rangle_{2,R}}{\langle \widehat b_2\rangle_{2,R}}, \qquad \text{ for }g_{\mathtt w}\in L^2(|\nu_2|).
		\]
		So, for $g_{\mathtt w}\in L^2(|\nu_2|)$ and $(Q,R)\in \mathscr S'_{4,\tr}$ we have
		\begin{equation}\label{eq:property_coef_c_Q_R}
			\begin{split}
				\|\Delta_{2,R}g_{\mathtt w}\|^2_{L^2(|\nu_2|)}&\geq \int_{R_Q} |\Delta_{2,R}g_{\mathtt w}|^2\, d|\nu_2|\\
				&=|c_{R,Q}(g_{\mathtt w})|^2\int_{R_Q}|\widehat b_2|^2\, d|\nu_2|=|c_{R,Q}(g_{\mathtt w})|^2\, |\nu_2|(R_Q)
			\end{split}
		\end{equation}
		and
		\begin{equation}\label{eq:identity_coefficients_c_Q_R}
			\chi_{R_Q}\Delta_{2,R}g_{\mathtt w}= c_{R,Q}(g_{\mathtt w})\, \chi_{R_Q}\widehat b_2.
		\end{equation}
		{
			We recall that the parameter $\widehat{\mathbf m}$ is introduced in the construction of the exceptional set, see \eqref{eq:definition_ex_set_s_+} and \eqref{eq:definition_ex_set_s_-}. In the present argument, the measures are $|\nu_i|=|\widehat\nu_i|+\mu|_{X\setminus B_i}$ and the operator is $S_{|\nu_1|,\varepsilon}$ associated with \eqref{eq:definition_auxiliary_kernel_sec6}. Thus, the high- and low-density stopping conditions give
			\[
			\widehat{\mathbf m}^{-1}\, \mu(Q)\leq|\nu_i|(Q)\leq\widehat{\mathbf m}\, \mu(Q),\qquad Q\in\bbD_i^{\tr},\ Q\subset Q_2(w_i),\  i=1,2.
			\]
			We first record the consequence of these inequalities that is needed below.
			\begin{lemma}\label{lemma:pack_nu_12_parap}
				Let $\mathscr F\subset\bbD_1^{\tr}$ be a family of pairwise disjoint cubes contained in $Q_2(w_1)$. Then, for every $P\in\bbD_2^{\tr}$ with $P\subset Q_2(w_2)$ we have
				\[
				\sum_{Q\in\mathscr F,\,Q\subset P}|\nu_1|(Q)\leq\widehat{\mathbf m}^{2}\, |\nu_2|(P).
				\]
				The same assertion holds with the indices $1$ and $2$ interchanged.
			\end{lemma}
			\begin{proof}
				Since all the cubes in the sum are transit and $\mathscr F$ is a family of pairwise disjoint cubes, we have
				\[
				\sum_{Q\in\mathscr F,\,Q\subset P}|\nu_1|(Q)\leq\widehat{\mathbf m}\sum_{Q\in\mathscr F,\,Q\subset P}\mu(Q)\leq\widehat{\mathbf m}\, \mu(P)\leq\widehat{\mathbf m}^{2}\, |\nu_2|(P).\qedhere
				\]
			\end{proof}
			
			\vv
			In the following two lemmas, we use the coefficients $c_{R,Q}(g_{\mathtt w})$ from \eqref{eq:property_coef_c_Q_R}, with the averages taken with respect to $|\nu_2|$. By \eqref{eq:identity_coefficients_c_Q_R}, the local term in \eqref{eq:s_4_1_tr}, for the auxiliary operator, is
			\[
			\begin{split}
				&\sum_{(Q,R)\in\mathscr S'_{4,\tr}}\big\langle\Delta_{1,Q}f_{\mathtt w},S^*_{|\nu_2|,\varepsilon}(\chi_{R_Q}\Delta_{2,R}g_{\mathtt w})\big\rangle\\
				&\quad=\sum_{(Q,R)\in\mathscr S'_{4,\tr}}c_{R,Q}(g_{\mathtt w})\big\langle\Delta_{1,Q}f_{\mathtt w},S^*_{|\nu_2|,\varepsilon}\widehat b_2\big\rangle
				-\sum_{(Q,R)\in\mathscr S'_{4,\tr}}c_{R,Q}(g_{\mathtt w})\big\langle\Delta_{1,Q}f_{\mathtt w},S^*_{|\nu_2|,\varepsilon}(\chi_{X\setminus R_Q}\widehat b_2)\big\rangle\\
				&\quad\eqqcolon\widehat S_1-\widehat S_2.
			\end{split}
			\]
			\begin{lemma}\label{lemma:s_41_tr_s_2}
				We have that
				\[
				|\widehat S_2|\lesssim\el^{m\gamma/2}\widehat{\mathbf m}\,\|f_{\mathtt w}\|_{L^2(|\nu_1|)}\|g_{\mathtt w}\|_{L^2(|\nu_2|)},
				\]
				where the implicit constant depends only on the structural parameters.
			\end{lemma}
			\begin{proof}
				We denote by $z_Q$ the center of $Q$. As in \cite[display (5.49)]{To14}, the kernel estimate \eqref{eq:regularity_auxiliary_kernel_sec6}, the identity $\int\Delta_{1,Q}f_{\mathtt w}\,d|\nu_1|=0$, and the fact that $|\widehat b_2|\equiv1$ yield
				\begin{equation}\label{eq:s_41_tr_s_2_aux}
					\big|\big\langle\Delta_{1,Q}f_{\mathtt w},S^*_{|\nu_2|,\varepsilon}(\chi_{X\setminus R_Q}\widehat b_2)\big\rangle\big|
					\lesssim\int_{X\setminus R_Q}\frac{\ell(Q)^\gamma}{d(x,z_Q)^{n+\gamma}}\,d|\nu_2|(x)\,\|\Delta_{1,Q}f_{\mathtt w}\|_{L^1(|\nu_1|)}.
				\end{equation}
				Here the goodness parameter $m$ is chosen sufficiently large so that \eqref{eq:diam_vs_length} and the separation from $X\setminus R_Q$ allow us to apply the kernel regularity estimate.
				Since $R_Q$ is a child of the transit cube $R$, Remark \ref{remark:growth_children_transit_cubes} gives $|\nu_2|(B(z_Q,t))\lesssim t^n$ for $t\geq\ell(R_Q)$.
				Moreover, the goodness of $Q$ gives
				\[
				\dist(Q,X\setminus R_Q)\geq\ell(Q)^\alpha\ell(R_Q)^{1-\alpha}.
				\]
				Thus, since $\alpha(n+\gamma)=\gamma/2$, a standard decomposition into annuli gives
				\[
				\begin{split}
					\int_{X\setminus R_Q}\frac{\ell(Q)^\gamma}{d(x,z_Q)^{n+\gamma}}\,d|\nu_2|(x)
					&\lesssim\frac{\ell(Q)^\gamma\ell(R_Q)^n}{(\ell(Q)^\alpha\ell(R_Q)^{1-\alpha})^{n+\gamma}}+\sum_{k=0}^{\infty}\frac{\ell(Q)^\gamma(2^{k+1}\ell(R_Q))^n}{(2^k\ell(R_Q))^{n+\gamma}}\\
					&\lesssim\frac{\ell(Q)^{\gamma/2}}{\ell(R_Q)^{\gamma/2}}.
				\end{split}
				\]
				We combine this bound with \eqref{eq:s_41_tr_s_2_aux} and apply the Cauchy-Schwarz inequality to obtain
				\begin{equation}\label{eq:l518a}
					\begin{split}
						&\big|c_{R,Q}(g_{\mathtt w})\big\langle\Delta_{1,Q}f_{\mathtt w},S^*_{|\nu_2|,\varepsilon}(\chi_{X\setminus R_Q}\widehat b_2)\big\rangle\big|\\
						&\qquad\qquad\lesssim |c_{R,Q}(g_{\mathtt w})|\frac{\ell(Q)^{\gamma/2}}{\ell(R_Q)^{\gamma/2}}|\nu_1|(Q)^{1/2}\|\Delta_{1,Q}f_{\mathtt w}\|_{L^2(|\nu_1|)}.
					\end{split}
				\end{equation}
				Fix $j\geq m+1$, a cube $R$, and a transit child $S$ of $R$. For all $Q$ such that $R_Q=S$, the coefficient $c_{R,Q}(g_{\mathtt w})$ has the same value, and \eqref{eq:identity_coefficients_c_Q_R} gives
				\[
				|c_{R,Q}(g_{\mathtt w})||\nu_2|(S)^{1/2}=\|\chi_S\Delta_{2,R}g_{\mathtt w}\|_{L^2(|\nu_2|)}.
				\]
				The cubes $Q$ with $R_Q=S$ and $\ell(Q)=\el^j\ell(R)$ are pairwise disjoint. Hence, Lemma \ref{lemma:pack_nu_12_parap} and the Cauchy-Schwarz inequality imply
				\begin{equation}\label{eq:l518b}
					\begin{split}
						&\sum_{\substack{Q\in\bbD_1^{\tr}:\,R_Q=S\\\ell(Q)=\el^j\ell(R)}}|c_{R,Q}(g_{\mathtt w})|\frac{\ell(Q)^{\gamma/2}}{\ell(S)^{\gamma/2}}|\nu_1|(Q)^{1/2}\|\Delta_{1,Q}f_{\mathtt w}\|_{L^2(|\nu_1|)}\\
						&\quad\leq\el^{(j-1)\gamma/2}\widehat{\mathbf m}\,\|\chi_S\Delta_{2,R}g_{\mathtt w}\|_{L^2(|\nu_2|)}\biggl(\sum_{\substack{Q\in\bbD_1^{\tr}:\,R_Q=S\\\ell(Q)=\el^j\ell(R)}}\|\Delta_{1,Q}f_{\mathtt w}\|_{L^2(|\nu_1|)}^2\biggr)^{1/2}.
					\end{split}
				\end{equation}
				For fixed $j$, each $Q$ occurs for at most one pair $(R,S)$, and the children of $R$ are disjoint. Therefore, summing \eqref{eq:l518b} over $R$ and its transit children, using the Cauchy-Schwarz inequality, and applying \eqref{eq:dec_f_mart_2}, we obtain
				\[
				|\widehat S_2|\lesssim\widehat{\mathbf m}\sum_{j\geq m+1}\el^{(j-1)\gamma/2}\|f_{\mathtt w}\|_{L^2(|\nu_1|)}\|g_{\mathtt w}\|_{L^2(|\nu_2|)}\lesssim\el^{m\gamma/2}\widehat{\mathbf m}\,\|f_{\mathtt w}\|_{L^2(|\nu_1|)}\|g_{\mathtt w}\|_{L^2(|\nu_2|)}.\qedhere
				\]
			\end{proof}
			
			It remains to estimate $\widehat S_1$.
			\begin{lemma}\label{lemma:s_41_tr_s_1}
				We have that
				\[
				|\widehat S_1|=\biggl|\sum_{(Q,R)\in\mathscr S'_{4,\tr}}c_{R,Q}(g_{\mathtt w})\big\langle\Delta_{1,Q}f_{\mathtt w},S^*_{|\nu_2|,\varepsilon}\widehat b_2\big\rangle\biggr|\lesssim c_\sharp\widehat{\mathbf m}\,\|f_{\mathtt w}\|_{L^2(|\nu_1|)}\|g_{\mathtt w}\|_{L^2(|\nu_2|)}.
				\]
			\end{lemma}
			\begin{proof}
				Let $\widetilde g\coloneqq g_{\mathtt w}-E_2g_{\mathtt w}$. Then, for $(Q,R)\in\mathscr S'_{4,\tr}$,
				\begin{equation}\label{eq:par_aux_1}
					c_{R,Q}(g_{\mathtt w})=\frac{\langle g_{\mathtt w}\rangle_{2,R_Q}}{\langle\widehat b_2\rangle_{2,R_Q}}-\frac{\langle g_{\mathtt w}\rangle_{2,R}}{\langle\widehat b_2\rangle_{2,R}}=\frac{\langle\widetilde g\rangle_{2,R_Q}}{\langle\widehat b_2\rangle_{2,R_Q}}-\frac{\langle\widetilde g\rangle_{2,R}}{\langle\widehat b_2\rangle_{2,R}}.
				\end{equation}
				We denote by $\mathcal F_\Pi$ the family of good cubes $Q\in\bbD_1^{\tr}$ for which there exists $R$ with $(Q,R)\in\mathscr S'_{4,\tr}$. For $Q\in\mathcal F_\Pi$, let $R(Q)$ be the smallest cube among the transit children $R_Q$ occurring in this sum\footnote{To prove the existence of such $R(Q)$, it suffices to argue as in \cite[p. 161]{To14}.}. These children and their ancestors contained in $Q_2(w_2)$ form a chain up to the top cube $Q_2(w_2)$. Since $\widetilde g$ has zero average on that top cube, \eqref{eq:par_aux_1} gives
				\[
				\sum_{R:\,(Q,R)\in\mathscr S'_{4,\tr}}c_{R,Q}(g_{\mathtt w})=\frac{\langle\widetilde g\rangle_{2,R(Q)}}{\langle\widehat b_2\rangle_{2,R(Q)}}.
				\]
				Thus, we define the paraproduct
				\[
				\Pi_{\Xi_{\mathtt w}}(\phi)\coloneqq\sum_{Q\in\mathcal F_\Pi}\frac{\langle\phi\rangle_{2,R(Q)}}{\langle\widehat b_2\rangle_{2,R(Q)}}\,\Delta^*_{1,Q}S^*_{|\nu_2|,\varepsilon}\widehat b_2,\qquad\text{ for }\phi\in L^2(|\nu_2|),
				\]
				and, arguing as in \cite[p. 166]{To14}, we obtain the identity
				\begin{equation}\label{eq:writing_paraproduct_equivalence}
					\widehat S_1=\sum_{(Q,R)\in\mathscr S'_{4,\tr}}c_{R,Q}(g_{\mathtt w})\big\langle\Delta_{1,Q}f_{\mathtt w},S^*_{|\nu_2|,\varepsilon}\widehat b_2\big\rangle=\langle f_{\mathtt w},\Pi_{\Xi_{\mathtt w}}(\widetilde g)\rangle.
				\end{equation}
				We first consider finite sums and obtain the corresponding bounded operator by passage to the limit.
				
				To prove that $\Pi_{\Xi_{\mathtt w}}$ is bounded from $L^2(|\nu_2|)$ to $L^2(|\nu_1|)$, we argue by duality. For $\varphi\in L^2(|\nu_1|)$, we define
				\[
				a_Q\coloneqq\begin{cases}
					\displaystyle\frac{|\langle\Delta_{1,Q}\varphi,S^*_{|\nu_2|,\varepsilon}\widehat b_2\rangle|^2}{\|\Delta_{1,Q}\varphi\|_{L^2(|\nu_1|)}^2},&\text{if }\Delta_{1,Q}\varphi\neq0,\\
					0,&\text{if }\Delta_{1,Q}\varphi=0,
				\end{cases}
				\qquad A_R\coloneqq\sum_{\substack{Q\in\mathcal F_\Pi\\R(Q)=R}}a_Q.
				\]
				By \eqref{eq:property_id_mart_diff} and the Cauchy-Schwarz inequality,
				\[
				a_Q\leq\|\Delta^*_{1,Q}S^*_{|\nu_2|,\varepsilon}\widehat b_2\|_{L^2(|\nu_1|)}^2.
				\]
				For $S\in\bbD_2^{\tr}$, let $\mathcal F^+_\Pi(S)$ be the family of maximal cubes in $\mathcal F_\Pi$ contained in $S$. The locality and the square-function bound for the adjoint martingale differences, together with \eqref{eq:localized_maximal_testing_sec5}, give
				\[
				\begin{split}
					\sum_{\substack{R\in\bbD_2^{\tr}\\R\subset S}}A_R
					&\leq\sum_{P\in\mathcal F^+_\Pi(S)}\sum_{\substack{Q\in\mathcal F_\Pi\\Q\subset P}}\|\Delta^*_{1,Q}(\chi_P S^*_{|\nu_2|,\varepsilon}\widehat b_2)\|_{L^2(|\nu_1|)}^2\\
					&\lesssim\sum_{P\in\mathcal F^+_\Pi(S)}\|\chi_P S^*_{|\nu_2|,\varepsilon}\widehat b_2\|_{L^2(|\nu_1|)}^2
					\lesssim c_\sharp^2\sum_{P\in\mathcal F^+_\Pi(S)}|\nu_1|(P)
					\leq c_\sharp^2\widehat{\mathbf m}^{\,2}|\nu_2|(S),
				\end{split}
				\]
				where in the last inequality we used Lemma \ref{lemma:pack_nu_12_parap}. If $S$ is terminal, there is no transit cube $R\subset S$, so the corresponding Carleson sum is zero. Moreover, since $R$ is transit, accretivity gives $|\langle\widehat b_2\rangle_{2,R}|\gtrsim1$. Hence, the Cauchy-Schwarz inequality, \eqref{eq:dec_f_mart_2}, and the Carleson embedding theorem yield, for $\psi\in L^2(|\nu_2|)$,
				\begin{equation}\label{eq:L2bound_paraproduct_claim}
					\begin{split}
						|\langle\varphi,\Pi_{\Xi_{\mathtt w}}\psi\rangle|
						&\leq\biggl(\sum_{Q\in\mathcal F_\Pi}\|\Delta_{1,Q}\varphi\|_{L^2(|\nu_1|)}^2\biggr)^{1/2}\biggl(\sum_{R\in\bbD_2^{\tr}}\biggl|\frac{\langle\psi\rangle_{2,R}}{\langle\widehat b_2\rangle_{2,R}}\biggr|^2A_R\biggr)^{1/2}\\
						&\lesssim c_\sharp\widehat{\mathbf m}\,\|\varphi\|_{L^2(|\nu_1|)}\|\psi\|_{L^2(|\nu_2|)}.
					\end{split}
				\end{equation}
				Thus, by \eqref{eq:writing_paraproduct_equivalence}, \eqref{eq:L2bound_paraproduct_claim}, and the bound for $E_2$, we have
				\[
				|\widehat S_1|\lesssim c_\sharp\widehat{\mathbf m}\,\|f_{\mathtt w}\|_{L^2(|\nu_1|)}\|g_{\mathtt w}-E_2g_{\mathtt w}\|_{L^2(|\nu_2|)}
				\lesssim c_\sharp\widehat{\mathbf m}\,\|f_{\mathtt w}\|_{L^2(|\nu_1|)}\|g_{\mathtt w}\|_{L^2(|\nu_2|)}.\qedhere
				\]
			\end{proof}
		}
	}
	\vv
	
	{
		{
			We gather \eqref{eq:s_4_1_tr}, \eqref{eq:S_4_1_far_away_part}, and Lemmas \ref{lemma:s_41_tr_s_2} and \ref{lemma:s_41_tr_s_1}, and obtain
			\begin{equation}\label{eq:S_{4,transit}}
				|\mathcal S'_{4,\tr}|\lesssim\bigl(1+(1+c_\sharp)\widehat{\mathbf m}\bigr)\, \|f_{\mathtt w}\|_{L^2(|\nu_1|)}\|g_{\mathtt w}\|_{L^2(|\nu_2|)}.
			\end{equation}
			
			\subsubsection*{\underline{The sum over $\mathscr S'_{4,\term}$}}
			Recall that
			\[
			\mathscr S'_{4,\term}\coloneqq\bigl\{(Q,R)\in\bbD_1^{\tr}\times\bbD_2^{\tr}:Q\subset R,\ \ell(Q)<\el^m\ell(R),\ R_Q\text{ terminal}\bigr\}.
			\]
			We reindex this sum by the terminal child, which we now denote by $R$, and write $\widehat R$ for its dyadic parent. Thus $R\in\bbD_2^{\term}$, $\widehat R\in\bbD_2^{\tr}$, and $Q\subset R$. Since the side lengths belong to the same discrete sequence, the condition $\ell(Q)<\el^m\ell(\widehat R)$ is equivalent to $\ell(Q)\leq\el^m\ell(R)$. The terminal cubes $R$ with transit parent are pairwise disjoint.
			
			We use the following Whitney decomposition. We prove it in Appendix \ref{appendix:random_cubes}.
			\begin{lemma}\label{prop:whitney_cubes_metric}
				Let $(X,d)$ be a geometrically doubling metric space with parameter $D$, and let $\bbD$ be the system of dyadic cubes in Theorem \ref{thm:HK22}. For every $Q\in\bbD$ with $X\setminus Q\neq\varnothing$ there exists a family $\mathcal W(Q)$ of cubes $S\in\bbD$ contained in $Q$ such that:
				\begin{enumerate}
					\item Cubes in $\mathcal W(Q)$ are pairwise disjoint and cover the interior of $Q$.
					\item There exist $\mathbf C_W>\mathbf c_W>1$ such that
					\[
					\mathbf c_W\ell(S)\leq\dist(S,X\setminus Q)\leq\mathbf C_W\ell(S),\qquad S\in\mathcal W(Q).
					\]
					\item There exists $\mathbf c>0$ such that
					\begin{equation}\label{eq:bounded_overlap}
						\sum_{S\in\mathcal W(Q)}\chi_{2S}(x)\leq\mathbf c,\qquad x\in X,
					\end{equation}
					where $2S=\{x\in X:\dist(x,S)<\ell(S)\}$.
				\end{enumerate}
				The constants depend only on $D$ and the parameters in Theorem \ref{thm:HK22}.
			\end{lemma}
			
			For each good cube $Q$ occurring in the sum, let $S\in\mathcal W(R)$ be the unique cube containing $z_Q$. We remark that goodness implies $\dist(Q,X\setminus R)>0$, so this cube exists. Moreover,
			\[
			\ell(Q)^\alpha\ell(R)^{1-\alpha}\leq\dist(Q,X\setminus R)
			\leq\dist(z_Q,X\setminus R)\leq(\mathbf C_W+2C_1)\ell(S),
			\]
			where $\diam(S)\leq2C_1\ell(S)$. Since $\ell(Q)\leq\el^m\ell(R)$, choosing $m$ sufficiently large gives
			\begin{equation}\label{eq:ell_Q_ell_S}
				\diam(Q)\leq 2C_1\,\ell(Q)\leq\frac14\, \ell(S)\qquad\text{ and }\qquad\ell(Q)^\alpha\ell(R)^{1-\alpha}\lesssim\ell(S).
			\end{equation}
			In particular, $Q\subset2S$, although $Q$ need not be contained in $S$. Also $2S\subset R$, by the second property of the Whitney decomposition.
			
			We can therefore write
			\[
			\begin{split}
				\mathcal S'_{4,\term}
				&=\sum_{R,S,Q}\bigl\langle S_{|\nu_1|,\varepsilon}(\Delta_{1,Q}f_{\mathtt w}),\chi_{2S}\Delta_{2,\widehat R}g_{\mathtt w}\bigr\rangle_*\\
				&\quad+\sum_{R,S,Q}\bigl\langle S_{|\nu_1|,\varepsilon}(\Delta_{1,Q}f_{\mathtt w}),\chi_{\widehat R\setminus2S}\Delta_{2,\widehat R}g_{\mathtt w}\bigr\rangle_*
				\eqqcolon\widehat S^{\term}_1+\widehat S^{\term}_2.
			\end{split}
			\]
			Here $R$ ranges over terminal cubes with transit parent, $S\in\mathcal W(R)$, and $Q$ ranges over the good transit cubes such that $Q\subset R$, $\ell(Q)\leq\el^m\ell(R)$, and $z_Q\in S$.
			
			\begin{lemma}\label{lem:bound_s_term}
				If $m$ is chosen sufficiently large, depending only on the structural constants and $\alpha$, then
				\[
				|\widehat S^{\term}_1|+|\widehat S^{\term}_2|
				\lesssim\|f_{\mathtt w}\|_{L^2(|\nu_1|)}\|g_{\mathtt w}\|_{L^2(|\nu_2|)},
				\]
				uniformly in $\mathtt w$ and $\varepsilon$.
			\end{lemma}
			\begin{proof}
				To bound $\widehat S^{\term}_1$, we fix $R$ and $S\in\mathcal W(R)$, and we set $g'\coloneqq\chi_{2S}\Delta_{2,\widehat R}g_{\mathtt w}$. Since $R\subset\mathscr W_2(w_2)$ and $F_2(w_2)\cap\mathscr W_2(w_2)=\varnothing$, by \eqref{eq:definition_suppressions_2b} and the second property in Lemma \ref{prop:whitney_cubes_metric}, we have
				\[
				\Xi_{\mathtt w}(x)\geq\dist(x,F_2(w_2))\geq\dist(x,X\setminus R)\geq(\mathbf c_W-1)\,\ell(S),\qquad\text{ for }x\in2S.
				\]
				We set $t_S\coloneqq\ell(S)+\inf_{2S}\Xi_{\mathtt w}$, so that the previous bound gives \[t_S\leq\mathbf c_W(\mathbf c_W-1)^{-1}\inf_{2S}\Xi_{\mathtt w}.\] Hence, \eqref{eq:size_auxiliary_kernel_sec6} yields
				\begin{equation}\label{eq:bound_s_term_2}
					|L_\varepsilon(x,y)|\lesssim t_S^{-n},\qquad\text{ for }x,y\in2S.
				\end{equation}
				Moreover, $\diam(2S)\leq(2C_1+2)\,\ell(S)$ and $\Xi_{\mathtt w}$ is $1$-Lipschitz, so, for every $z\in2S$, we have $2S\subset B(z,(2C_1+3)\,\ell(S))$ and $\Xi_{\mathtt w}(z)\leq t_S+(2C_1+2)\,\ell(S)$.
				
				Thus, \eqref{eq:growth_auxiliary_measures_sec6} gives $|\nu_i|(2S)\lesssim t_S^n$ for $i=1,2$, and, by \eqref{eq:bound_s_term_2} and the Cauchy-Schwarz inequality, we obtain
				\[
				\|\chi_{2S}S^*_{|\nu_2|,\varepsilon}g'\|_{L^2(|\nu_1|)}\lesssim t_S^{-n}\,|\nu_1|(2S)^{1/2}\,|\nu_2|(2S)^{1/2}\,\|g'\|_{L^2(|\nu_2|)}\lesssim\|g'\|_{L^2(|\nu_2|)}.
				\]
				The cubes $Q$ assigned to $(R,S)$ are contained in $2S$, as observed before the statement. Hence, by \eqref{eq:property_id_mart_diff}, the Cauchy-Schwarz inequality, and the square-function bound for the adjoint martingale differences stated after \eqref{eq:dec_f_mart_adj}, we have
				\[
				\begin{split}
					\sum_Q\bigl|\bigl\langle S_{|\nu_1|,\varepsilon}(\Delta_{1,Q}f_{\mathtt w}),g'\bigr\rangle_*\bigr|
					&=\sum_Q\bigl|\bigl\langle\Delta_{1,Q}f_{\mathtt w},\Delta_{1,Q}^*(\chi_{2S}S^*_{|\nu_2|,\varepsilon}g')\bigr\rangle\bigr|\\
					&\leq\Bigl(\sum_Q\|\Delta_{1,Q}f_{\mathtt w}\|_{L^2(|\nu_1|)}^2\Bigr)^{1/2}
					\Bigl(\sum_Q\|\Delta_{1,Q}^*(\chi_{2S}S^*_{|\nu_2|,\varepsilon}g')\|_{L^2(|\nu_1|)}^2\Bigr)^{1/2}\\
					&\lesssim\Bigl(\sum_Q\|\Delta_{1,Q}f_{\mathtt w}\|_{L^2(|\nu_1|)}^2\Bigr)^{1/2}
					\|\chi_{2S}\Delta_{2,\widehat R}g_{\mathtt w}\|_{L^2(|\nu_2|)},
				\end{split}
				\]
				where the sums run over the cubes $Q$ assigned to $(R,S)$. Each $Q$ is assigned to at most one pair $(R,S)$, because the terminal cubes with transit parent are pairwise disjoint and so are the cubes in $\mathcal W(R)$. Moreover, by \eqref{eq:bounded_overlap} and the inclusion $2S\subset R$, grouping the terminal cubes $R$ with the same parent, we have
				\[
				\sum_{R,S}\|\chi_{2S}\Delta_{2,\widehat R}g_{\mathtt w}\|_{L^2(|\nu_2|)}^2
				\lesssim\sum_R\|\chi_R\Delta_{2,\widehat R}g_{\mathtt w}\|_{L^2(|\nu_2|)}^2
				\leq\sum_{P\in\bbD_2^{\tr}}\|\Delta_{2,P}g_{\mathtt w}\|_{L^2(|\nu_2|)}^2.
				\]
				Finally, summing the previous bound over $(R,S)$, the Cauchy-Schwarz inequality and \eqref{eq:dec_f_mart_2} give \[|\widehat S^{\term}_1|\lesssim\|f_{\mathtt w}\|_{L^2(|\nu_1|)}\|g_{\mathtt w}\|_{L^2(|\nu_2|)}.\]
				
				To bound $\widehat S^{\term}_2$, for each triple $(R,S,Q)$ in the sum we set $g''\coloneqq\chi_{\widehat R\setminus2S}\Delta_{2,\widehat R}g_{\mathtt w}$. Since $z_Q\in S$, \eqref{eq:ell_Q_ell_S} gives
				\[
				\dist(Q,\supp g'')\geq\ell(S)-\diam(Q)\geq\tfrac34\,\ell(S)\gtrsim\ell(Q)^\alpha\ell(R)^{1-\alpha}.
				\]
				In particular, $\diam(Q)\leq\dist(Q,\supp g'')/3$, so \eqref{eq:regularity_auxiliary_kernel_sec6} and the identity $\int\Delta_{1,Q}f_{\mathtt w}\,d|\nu_1|=0$ give
				\[
				\bigl|\bigl\langle S_{|\nu_1|,\varepsilon}(\Delta_{1,Q}f_{\mathtt w}),g''\bigr\rangle_*\bigr|
				\lesssim\frac{\ell(Q)^\gamma}{\dist(Q,\supp g'')^{n+\gamma}}\,
				\|\Delta_{1,Q}f_{\mathtt w}\|_{L^1(|\nu_1|)}\,\|g''\|_{L^1(|\nu_2|)}.
				\]
				Since $\alpha(n+\gamma)=\gamma/2$ and $D(Q,\widehat R)\approx\ell(\widehat R)\approx\ell(R)$, we have
				\[
				\frac{\ell(Q)^\gamma}{\dist(Q,\supp g'')^{n+\gamma}}
				\lesssim\frac{\ell(Q)^{\gamma/2}}{\ell(R)^{n+\gamma/2}}
				\approx\frac{\ell(Q)^{\gamma/2}\ell(\widehat R)^{\gamma/2}}{D(Q,\widehat R)^{n+\gamma}},
				\]
				so the Cauchy-Schwarz inequality yields
				\begin{equation}\label{eq:bound_s_term_1}
					\begin{split}
						&\bigl|\bigl\langle S_{|\nu_1|,\varepsilon}(\Delta_{1,Q}f_{\mathtt w}),g''\bigr\rangle_*\bigr|\\
						&\qquad\lesssim
						\frac{\ell(Q)^{\gamma/2}\ell(\widehat R)^{\gamma/2}}{D(Q,\widehat R)^{n+\gamma}}\,
						|\nu_1|(Q)^{1/2}|\nu_2|(\widehat R)^{1/2}\,
						\|\Delta_{1,Q}f_{\mathtt w}\|_{L^2(|\nu_1|)}\,
						\|\Delta_{2,\widehat R}g_{\mathtt w}\|_{L^2(|\nu_2|)}.
					\end{split}
				\end{equation}
				Each pair $(Q,\widehat R)$ occurs at most once in $\widehat S^{\term}_2$, and both cubes are transit with $\ell(Q)\leq\ell(\widehat R)$. Therefore, \eqref{eq:bound_s_term_1}, \eqref{eq:sec5_two_measure_schur} combined with Schur's test, and \eqref{eq:dec_f_mart_2} give $|\widehat S^{\term}_2|\lesssim\|f_{\mathtt w}\|_{L^2(|\nu_1|)}\|g_{\mathtt w}\|_{L^2(|\nu_2|)}$, which concludes the proof.
			\end{proof}
			
		}
	}
	
	As a consequence of Lemma \ref{lem:bound_s_term}, we obtain
	\begin{equation}\label{eq:s_1_hat_term}
		|\mathcal S'_{4,\term}|\lesssim\|f_{\mathtt w}\|_{L^2(|\nu_1|)}\|g_{\mathtt w}\|_{L^2(|\nu_2|)}
	\end{equation}
	and applying the same estimates to $S^*_{|\nu_2|,\varepsilon}$, as in \eqref{eq:transpose_martingale_difference}, we obtain
	\begin{equation}\label{eq:{eq_s_4_pp}}
		|\mathcal S_4''|\lesssim \|f_{\mathtt w}\|_{L^2(|\nu_1|)}\, \|g_{\mathtt w}\|_{L^2(|\nu_2|)}.
	\end{equation}
	
	Then, we gather \eqref{eq:S_{4,transit}}, \eqref{eq:s_1_hat_term}, and \eqref{eq:{eq_s_4_pp}}, and obtain a constant $\mathbf c_{\mathcal S4}>0$ such that
	\begin{equation}\label{eq:final_estimate_S_4}
		\begin{split}
			|\mathcal S_4(\mathtt w)|&= |\mathcal S'_4 + \mathcal S''_4|\leq  |\mathcal S'_{4,\tr}|+ |\mathcal S'_{4,\term}| + |\mathcal S_4''|\leq \mathbf c_{\mathcal S4} \, \|f_{\mathtt w}\|_{L^2(|\nu_1|)}\,  \|g_{\mathtt w}\|_{L^2(|\nu_2|)},
		\end{split}
	\end{equation}
	where $\mathbf c_{\mathcal S4}$ depends only on the structural parameters and $m$, and is independent of $\mathtt w$, $\lambda$, $\upsilon$, $a$, and $\varepsilon$.
	
	{
		\vv
		Finally, \eqref{eq:bound_s_1_final}, Proposition \ref{prop:diagonal_terms},
		and \eqref{eq:final_estimate_S_4} give
		\[
		\begin{split}
			\widetilde{\mathbb E}_{\mathtt w}|\mathcal I_4(\mathtt w)|
			\leq\Bigl[&C_{m,\lambda,\upsilon}(1+a^{-n})+\bigl(C_{m,\upsilon}(M_\Lambda a)^{1/2}+\rho_m(\lambda,\upsilon)\bigr)(1+N_\varepsilon)\Bigr]\,
			\|f\|_{L^2(|\nu_1|)}\|g\|_{L^2(|\nu_2|)}.
		\end{split}
		\]
		By \eqref{eq:localized_pairing_comparison_sec5}, \eqref{eq:bound_exp_good}, \eqref{eq:martingale_difference_decomposition_extended}, \eqref{eq:final_bound_I_1}, \eqref{eq:inbalanced_estimate_2}, and \eqref{eq:final_bound_I_3_alternative}, we obtain
		\[
		\begin{split}
			&\widetilde{\mathbb E}_{\mathtt w}\left|\left\langle T_{\Xi_{\mathtt w},|\widehat\nu_1|,\varepsilon}f_{\mathtt w},g_{\mathtt w}\right\rangle_{|\widehat\nu_2|}\right|\\
			&\quad\leq\widetilde{\mathbb E}_{\mathtt w}\left|\left\langle S_{|\nu_1|,\varepsilon}f_{\mathtt w},g_{\mathtt w}\right\rangle_*\right|+C\|f\|_{L^2(|\nu_1|)}\|g\|_{L^2(|\nu_2|)}\\
			&\quad\leq\widetilde{\mathbb E}_{\mathtt w}|\mathcal I_4(\mathtt w)|+C\|f\|_{L^2(|\nu_1|)}\|g\|_{L^2(|\nu_2|)}\\
			&\quad\leq\Bigl[C_{m,\lambda,\upsilon}(1+a^{-n})+\bigl(C_{m,\upsilon}(M_\Lambda a)^{1/2}+\rho_m(\lambda,\upsilon)\bigr)(1+N_\varepsilon)\Bigr]\|f\|_{L^2(|\nu_1|)}\|g\|_{L^2(|\nu_2|)}.
		\end{split}
		\]
		Letting $K_1,K_2\to\infty$ and applying Fatou's lemma gives \eqref{eq:bound_supp_good_functions}, since the zero extensions of $f$ and $g$ have the same norms for $|\nu_i|$ and $|\widehat\nu_i|$. This proves Lemma \ref{lem:key_lemma_good_functions}.
		
		For fixed $m$, we have $\rho_m(\lambda,\upsilon)\to0$ as $\lambda,\upsilon\to0$. Once $\lambda,\upsilon$ are fixed, $M_\Lambda$ and $C_{m,\upsilon}$ are fixed, and $C_{m,\upsilon}(M_\Lambda a)^{1/2}\to0$ as $a\to0$. We conclude the section with the important remark that the term $C_{m,\lambda,\upsilon}(1+a^{-n})$ is independent of $N_\varepsilon$.
	}
	
	\vv
	\section{Proof of the main theorem}\label{sec:main_theorem_antisym}

	\subsection{Averaged suppressed operators and the probabilistic argument}\label{section:probabilisti_argument_subsection}
	Under the normalized measure-testing assumptions of Proposition \ref{proposition:bounded_testing_normalized} and within the framework of the mixed corona decomposition of Section \ref{sec:corona}, let $R_i \in \Top''_i$, for $i=1,2$, be such that $\widetilde{\Tree}(R_1, R_2) \neq \varnothing$. In this case, either $R_1 \subset R_2$ or $R_2 \subset R_1$. Without loss of generality we assume that $R_1\subset R_2$, the other case being analogous. We also recall that, in \eqref{eq:Phi-i} and \eqref{eq:tilde Phi}, we defined the suppression
	\[
	\widetilde \Phi_{R_1, R_2}(x) \coloneqq \max\big(\Phi_{1,R_1}(x), \Phi_{2,R_2}(x)\big),
	\]
	and, for $R \in \widetilde{\Top}$,
	\(
	\widetilde \Phi_R(x) \coloneqq \widetilde \Phi_{R_1, R_2}(x).
	\)
	For the sake of simplicity, we write $\Phi_{w_i}, \Phi_i$, and $ \Phi$ instead of $\Phi_{i,R_i,w_i}, \Phi_{i,R_i},	\widetilde \Phi_R$ respectively, for $i=1,2$. Finally, recall that in \eqref{eq:definition_suppressions_2ab} and \eqref{eq:definition_suppressions_2b} we introduced
	\begin{equation*}
		\Xi_{i}	\coloneqq \Xi_{i,R_i, w_i}
		\coloneqq 
		\max\bigl\{\Phi_{i,R_i}, \Phi_{i,R_i, w_i}\bigr\} \qquad \text{ and }\qquad \Xi_{\mathtt w}
		\coloneqq 
		\max\bigl\{\Xi_{1,R_1, w_1}, \Xi_{2,R_2, w_2}\bigr\}.
	\end{equation*}

	\vv
	
	Now we define the \textit{averaged suppressed SIO}. Namely, given $\varepsilon>0$ and a non-negative Radon measure $\mathfrak m$ on $X$, we set
	\begin{equation}\label{eq:definition_averaged_T}
		\widetilde T_{\mathfrak m, \varepsilon} f(x) \coloneqq \int T_{\Xi_{\mathtt w },\mathfrak m, \ve} f(x)\, d\mathbb P^{\Omega^2}_{\mathscr G}(\mathtt w), \qquad x\in X.
	\end{equation}
	Initially, we take $f$ bounded and with bounded support. For the finite original testing measures $|\widehat\nu_i|$, the definitions also apply to $f\in L^2(|\widehat\nu_i|)$. By Fubini's theorem, $\widetilde T$ is associated with the kernel
	\begin{equation}\label{eq:averaged_kernel_def}
		\widetilde K(x,y) \coloneqq \int K_{\Xi_{\mathtt w }}(x,y) \, d\mathbb P^{\Omega^2}_{\mathscr G}(\mathtt w).
	\end{equation}

	Given a non-negative $1$-Lipschitz function $\Psi\colon X\to \mathbb R$, we define 
	\[
	{\mathcal T}_{\Psi,\mathfrak m, \varepsilon} f(x)\coloneqq \int_{d(x,y)\geq \max(\Psi(x),\,\ve)} K(x,y) \, f(y)\,d\mathfrak m(y) = \int \mathcal K_{\Psi,\ve}(x,y)\, f(y)\, d\mathfrak m(y),
	\]
	where
	\[
	\mathcal K_{\Psi,\ve}(x,y) \coloneqq K(x,y)\, \chi_{\{d(x,y)\geq\max(\Psi(x),\ve)\}}.
	\]
	Then, in the case $\Psi= \Xi_{\mathtt w}$, we define the averaged operator
	\[
	\widetilde{\mathcal T}_{\mathfrak m, \ve}f(x) \coloneqq \int {\mathcal T}_{\Xi_{\mathtt w},\mathfrak m, \ve} f(x)\, d\mathbb P^{\Omega^2}_{\mathscr G}(\mathtt w), \qquad x\in X.
	\]
	
	We also set
	\[
	\widetilde{\mathcal T}_{\mathfrak m,\sharp,\ve}f(x)\coloneqq\sup_{\rho\geq\ve}|\widetilde{\mathcal T}_{\mathfrak m,\rho}f(x)|
	\qquad \text{ and }\qquad  \widetilde T_{\mathfrak m,\sharp,\ve}f(x)\coloneqq\sup_{\rho\geq\ve}|\widetilde T_{\mathfrak m,\rho}f(x)|.
	\]

	We define
	\begin{equation}\label{eq:frak v}
		\begin{split}
			\mathfrak v_{x,\ve}(t) &\coloneqq \mathbb P^{\Omega^2}_{\mathscr G}\bigl( \{ \mathtt w \in \Omega^2:  \max(\Xi_{\mathtt w}(x),\varepsilon)\leq t\}\bigr), \qquad \text{ for }t\in \mathbb R,
		\end{split}
	\end{equation}
	which corresponds to the distribution function associated with the random variable $\max(\Xi_{\mathtt w}(x),\varepsilon)$.
	We remark that the original grids associated with $w_1$ and $w_2$ are independent. However, they may not be independent after conditioning on $\mathscr G$. We define
	\[
	A_i(t)\coloneqq \bigl\{w_i\in\Omega: t\geq \max(\Phi_{w_i}(x), \Phi_i(x), \varepsilon)\bigr\}, \qquad \text{ for } i=1,2,\, t\in \mathbb R,
	\]
	so \eqref{eq:frak v} yields
	\[
	\mathfrak v_{x,\ve}(t)= \mathbb P^{\Omega^2}_{\mathscr G}(A_1(t)\times A_2(t)).
	\]
	
	By Fubini's theorem and \eqref{eq:frak v}, $\widetilde{\mathcal T}$ is associated with the kernel
	\begin{equation}\label{eq:averaged kernel}
		\begin{split}
			\widetilde {\mathcal K}(x,y) &= \int K(x,y)\, \chi_{\{d(x,y)\geq \max(\Xi_{\mathtt w}(x),\varepsilon)\}} \,d\mathbb P^{\Omega^2}_{\mathscr G}(\mathtt w)= K(x,y)\, \mathfrak v_{x,\ve}(d(x,y)).
		\end{split}
	\end{equation}
	
	By \eqref{eq:averaged_kernel_def} and the kernel estimates for $K_{\Xi_{\mathtt w}}$, the function $\widetilde K$ is a Calder\'on-Zygmund kernel with the constants of Lemma \ref{lemma:kernel_suppressed}, which depend only on $n$, $\gamma$, and the constants of $K$.
	However, this is not necessarily the case for $\widetilde {\mathcal K}(x,y)$ in \eqref{eq:averaged kernel}.
	
	\vvv
	
	We recall that $\delta_{\mathrm{pr}}$ is the parameter defined in Lemma \ref{lemma:lemma_probabilistic_sets}.
	Then, we choose the error $\kappa$ in \eqref{eq:large_probability_top_cubes} so that $0<\kappa<\delta_{\mathrm{pr}}^2/2$, and define
	\[
	c_{\mathscr G}\coloneqq \frac{\delta_{\mathrm{pr}}^2-\mathbb P^{\Omega^2}(\Omega^2\setminus \mathscr G)}{\mathbb P^{\Omega^2}(\mathscr G)}\geq \frac{\delta_{\mathrm{pr}}^2-\kappa}{\mathbb P^{\Omega^2}(\mathscr G)}\geq \frac{\delta_{\mathrm{pr}}^2}{2}.
	\]
	Here we used that $\mathbb P^{\Omega^2}(\mathscr G)\leq1$ and that $\mathbb P^{\Omega^2}(\Omega^2\setminus\mathscr G)<\kappa$, which holds by \eqref{eq:large_probability_top_cubes} for every goodness parameter $m\geq m_0$, with $m_0=m_0(\kappa)$ as in Lemma \ref{lem:random_capture_goodness}. In the rest of this section, we assume that $m\geq m_0$.

	\begin{lemma}\label{lem:final probab 1}
		For $\varepsilon>0$, we have
		\begin{equation}\label{eq:Tphi<tildeT}
			\|\mathcal T_{\Phi,|\widehat\nu_1|,\ve}\|_{L^2(|\widehat\nu_1|)\to L^2(|\widehat\nu_2|)}\leq 3c_{\mathscr G}^{-1}\|\widetilde{\mathcal T}_{|\widehat\nu_1|,\sharp,\ve}\|_{L^2(|\widehat\nu_1|)\to L^2(|\widehat\nu_2|)}
		\end{equation}
		and, for every $\mathtt w\in\Omega^2$,
		\begin{equation}\label{eq:TXi<tildeT}
			\|\mathcal T_{\Xi_{\mathtt w},|\widehat\nu_1|,\ve}\|_{L^2(|\widehat\nu_1|)\to L^2(|\widehat\nu_2|)}\leq 3c_{\mathscr G}^{-1}\|\widetilde{\mathcal T}_{|\widehat\nu_1|,\sharp,\ve}\|_{L^2(|\widehat\nu_1|)\to L^2(|\widehat\nu_2|)}.
		\end{equation}
	\end{lemma}
	\begin{proof}
		The function $\mathfrak v_{x,\ve}$ defined in \eqref{eq:frak v} satisfies the following properties:
		\begin{enumerate}
			\item If $0\leq t<\max(\Phi(x),\ve)$, then $\mathfrak v_{x,\ve}(t)=0$.
			\item If $t\geq\max(\Phi(x),\ve)$, then $\mathfrak v_{x,\ve}(t)\geq c_{\mathscr G}$.
			\item The function $\mathfrak v_{x,\ve}(\cdot)$ is non-decreasing and right-continuous.
		\end{enumerate}
		Properties (1) and (3) follow from the definition of $\mathfrak v_{x,\ve}$ and the inequality $\Xi_{\mathtt w}\geq\Phi$. To prove the second, let $t\geq\max(\Phi(x),\ve)$. By the definition of $\Phi_i(x)$ in \eqref{eq:definition_suppressions_2a}, for every integer $j\geq1$ there exist $\Gamma_{i,j}\subset\Omega$, $i=1,2$, such that
		\[
		\mathbb P^\Omega(\Gamma_{i,j})>\delta_{\mathrm{pr}},\qquad \sup_{w_i\in\Gamma_{i,j}}\Phi_{w_i}(x)<t+1/j.
		\]
		Thus $\Gamma_{i,j}\subset A_i(t+1/j)$ and, since $A_i(t+1/j)\downarrow A_i(t)$, we obtain $\mathbb P^\Omega(A_i(t))\geq\delta_{\mathrm{pr}}$. By independence with respect to the unconditioned probability, we have
		\[
		\begin{split}
			\mathfrak v_{x,\ve}(t)&=\frac{\mathbb P^{\Omega^2}((A_1(t)\times A_2(t))\cap\mathscr G)}{\mathbb P^{\Omega^2}(\mathscr G)}\\
			&\geq\frac{\mathbb P^\Omega(A_1(t))\mathbb P^\Omega(A_2(t))-\mathbb P^{\Omega^2}(\Omega^2\setminus\mathscr G)}{\mathbb P^{\Omega^2}(\mathscr G)}\geq c_{\mathscr G}.
		\end{split}
		\]
		
		Let us turn to the proof of \eqref{eq:Tphi<tildeT}. We denote $\Phi^\varepsilon(x)\coloneqq\max(\Phi(x),\ve)$, which is also $1$-Lipschitz. Since $\mathfrak v_{x,\ve}(t)$ is identically zero for $t<\Phi^\varepsilon(x)$, we infer that
		\[
		\widetilde{\mathcal K}(x,y)=K(x,y)\mathfrak v_{x,\ve}(d(x,y))\chi_{\{d(x,y)\geq\Phi^\varepsilon(x)\}}=\mathfrak v_{x,\ve}(d(x,y))\mathcal K_{\Phi,\ve}(x,y).
		\]
		Thus,
		\begin{equation}\label{eq:K_phi-ve-bis}
			\mathcal K_{\Phi,\ve}(x,y)=\frac{1}{\mathfrak v_{x,\ve}(d(x,y))}\widetilde{\mathcal K}(x,y)\chi_{\{d(x,y)\geq\Phi^\varepsilon(x)\}},
		\end{equation}
		where the reciprocal is only used on $\{d(x,y)\geq\Phi^\varepsilon(x)\}$.
		
		Since $1/\mathfrak v_{x,\varepsilon}$ is right-continuous and non-increasing on $[\Phi^\varepsilon(x),\infty)$, we can define the Lebesgue-Stieltjes measure $\sigma$ on $\mathbb R$ by
		\[
		\sigma((a,b])=\frac{1}{\mathfrak v_{x,\ve}(a)}-\frac{1}{\mathfrak v_{x,\ve}(b)},\quad \Phi^\varepsilon(x)\leq a<b,\qquad \sigma(( -\infty,\Phi^\varepsilon(x)])=0.
		\]
		Observe that $\sigma$ is a finite positive Borel measure. Indeed, since $\mathfrak v_{x,\ve}(t)\to1$ as $t\to\infty$, property (2) gives
		\begin{equation}\label{eq:sigma-bound}
			\sigma(\mathbb R)=\frac{1}{\mathfrak v_{x,\ve}(\Phi^\varepsilon(x))}-1\leq c_{\mathscr G}^{-1}.
		\end{equation}
		Now, for $d(x,y)\geq\Phi^\varepsilon(x)$, we write
		\[
		\begin{split}
			\int_{t>\Phi^\varepsilon(x)}\widetilde{\mathcal K}(x,y)\chi_{\{d(x,y)\geq t\}}\,d\sigma(t)
			&=\widetilde{\mathcal K}(x,y)\sigma((\Phi^\varepsilon(x),d(x,y)])\\
			&=\widetilde{\mathcal K}(x,y)\left(\frac{1}{\mathfrak v_{x,\ve}(\Phi^\varepsilon(x))}-\frac{1}{\mathfrak v_{x,\ve}(d(x,y))}\right).
		\end{split}
		\]
		Therefore, by \eqref{eq:K_phi-ve-bis},
		\[
		\mathcal K_{\Phi,\ve}(x,y)=\frac{\widetilde{\mathcal K}(x,y)}{\mathfrak v_{x,\ve}(\Phi^\varepsilon(x))}-\int_{t>\Phi^\varepsilon(x)}\widetilde{\mathcal K}(x,y)\chi_{\{d(x,y)\geq t\}}\,d\sigma(t).
		\]
		For $t\geq\ve$, we have $\mathfrak v_{x,\ve}(u)\chi_{\{u\geq t\}}=\mathfrak v_{x,t}(u)$ for $u\in\mathbb R$, so $\widetilde{\mathcal K}(x,y)\chi_{\{d(x,y)\geq t\}}$ is the kernel of $\widetilde{\mathcal T}_{|\widehat\nu_1|,t}$. Hence, multiplying this identity by $f(y)$, integrating with respect to $|\widehat\nu_1|$ in $y$, and using Fubini's theorem, we get
		\[
		\mathcal T_{\Phi,|\widehat\nu_1|,\ve}f(x)=\frac{\widetilde{\mathcal T}_{|\widehat\nu_1|,\ve}f(x)}{\mathfrak v_{x,\ve}(\Phi^\varepsilon(x))}-\int_{t>\Phi^\varepsilon(x)}\widetilde{\mathcal T}_{|\widehat\nu_1|,t}f(x)\,d\sigma(t).
		\]
		From the preceding identity and \eqref{eq:sigma-bound}, we infer that
		\[
		|\mathcal T_{\Phi,|\widehat\nu_1|,\ve}f(x)|\leq c_{\mathscr G}^{-1}|\widetilde{\mathcal T}_{|\widehat\nu_1|,\ve}f(x)|+c_{\mathscr G}^{-1}\sup_{\rho\geq\ve}|\widetilde{\mathcal T}_{|\widehat\nu_1|,\rho}f(x)|\leq3c_{\mathscr G}^{-1}\widetilde{\mathcal T}_{|\widehat\nu_1|,\sharp,\ve}f(x).
		\]
		Taking the $L^2(|\widehat\nu_2|)$ norm proves \eqref{eq:Tphi<tildeT}.
		
		Let us now turn our attention to \eqref{eq:TXi<tildeT}. By the definition of $\mathcal K_{\Xi_{\mathtt w},\ve}$ and the inequality $\Xi_{\mathtt w}\geq\Phi$, we have
		\[
		\mathcal K_{\Xi_{\mathtt w},\ve}(x,y)=\mathcal K_{\Phi,\ve}(x,y)\chi_{\{d(x,y)\geq\max(\Xi_{\mathtt w}(x),\ve)\}}=\frac{\widetilde{\mathcal K}(x,y)}{\mathfrak v_{x,\ve}(d(x,y))}\chi_{\{d(x,y)\geq\max(\Xi_{\mathtt w}(x),\ve)\}}.
		\]
		This is the same identity as \eqref{eq:K_phi-ve-bis}, with $\Phi$ replaced by $\Xi_{\mathtt w}$. We define $\sigma$ on $(\max(\Xi_{\mathtt w}(x),\ve),\infty)$ as above and obtain
		\[
		|\mathcal T_{\Xi_{\mathtt w},|\widehat\nu_1|,\ve}f(x)|\leq3c_{\mathscr G}^{-1}\widetilde{\mathcal T}_{|\widehat\nu_1|,\sharp,\ve}f(x).
		\]
		Taking the $L^2(|\widehat\nu_2|)$ norm proves \eqref{eq:TXi<tildeT} and completes the proof of the lemma.
	\end{proof}
	\vvv

	{
		Using maximal function estimates, the bounds in Lemma \ref{lem:final probab 1} imply analogous estimates in terms of $T$ and $\widetilde T$.
		\begin{lemma}\label{lemma:comparison_averaged_operators}
			If $\ve>0$, we have
			\begin{equation}\label{eq:T-phi<tilde T bis}
				\|T_{\Phi,|\widehat\nu_1|,\ve}\|_{L^2(|\widehat\nu_1|)\to L^2(|\widehat\nu_2|)}\lesssim 1+\sup_{\rho\geq\ve}\|\widetilde T_{|\widehat\nu_1|,\rho}\|_{L^2(|\widehat\nu_1|)\to L^2(|\widehat\nu_2|)},
			\end{equation}
			and, for every $\mathtt w\in\Omega^2$,
			\begin{equation}\label{eq:T-xi<tilde T bis}
				\|T_{\Xi_{\mathtt w},|\widehat\nu_1|,\ve}\|_{L^2(|\widehat\nu_1|)\to L^2(|\widehat\nu_2|)}\lesssim 1+\sup_{\rho\geq\ve}\|\widetilde T_{|\widehat\nu_1|,\rho}\|_{L^2(|\widehat\nu_1|)\to L^2(|\widehat\nu_2|)},
			\end{equation}
			where the implicit constants depend on the structural parameters and $c_{\mathscr G}^{-1}$, but are independent of $\mathtt w$, $\ve$, and the parameters $\lambda$, $\upsilon$, and $a$ of Lemma \ref{lem:key_lemma_good_functions}. With $\kappa$ fixed as above, these constants can be chosen independently of the goodness parameter $m$, since $c_{\mathscr G}^{-1}\leq2\delta_{\mathrm{pr}}^{-2}$.
		\end{lemma}
		\begin{proof}
			In this proof, we write $\|\cdot\|\coloneqq\|\cdot\|_{L^2(|\widehat\nu_1|)\to L^2(|\widehat\nu_2|)}$ for brevity. Let $\Psi\in\{\Phi,\Xi_{\mathtt w}\}$. By the end of the proof of Lemma \ref{lemma:final_global_bound}, we have $\max\{\mathcal R_1,\mathcal R_2\}\leq\Phi\leq\Psi$, where $\mathcal R_i$ is the non-Ahlfors radius in \eqref{eq:definition_r_i} associated with $\widehat\nu_i$. Hence, Lemma \ref{lem:suppressed maximal} gives
			\begin{equation}\label{eq:maximal_bound_sec7}
				\|M^r_{\Psi,|\widehat\nu_1|}f\|_{L^2(|\widehat\nu_2|)}\lesssim\|f\|_{L^2(|\widehat\nu_1|)}.
			\end{equation}
			
			Let $x\in X$, $\rho\geq\ve$, and $\rho_x\coloneqq\max(\Psi(x),\rho)$. By Lemma \ref{lemma:kernel_suppressed}, applied on $\{y:\rho<d(x,y)<\rho_x\}$, on the sphere $\{y:d(x,y)=\rho_x\}$, and on $\{y:d(x,y)>\rho_x\}$, we have
			\begin{equation}\label{eq:lemmma82_gen_in_psi}
				\bigl|T_{\Psi,|\widehat\nu_1|,\rho}f(x)-\mathcal T_{\Psi,|\widehat\nu_1|,\rho}f(x)\bigr|\lesssim M^r_{\Psi,|\widehat\nu_1|}f(x).
			\end{equation}
			Applying \eqref{eq:lemmma82_gen_in_psi} with $\rho=\ve$, \eqref{eq:maximal_bound_sec7}, \eqref{eq:Tphi<tildeT}, and \eqref{eq:TXi<tildeT}, we obtain
			\begin{equation}\label{eq:aux_probab_argument_norms}
				\|T_{\Psi,|\widehat\nu_1|,\ve}\|\lesssim1+\|\mathcal T_{\Psi,|\widehat\nu_1|,\ve}\|\lesssim1+c_{\mathscr G}^{-1}\|\widetilde{\mathcal T}_{|\widehat\nu_1|,\sharp,\ve}\|.
			\end{equation}
			Moreover, since $\Xi_{\mathtt w}\geq\Phi$, we have $M^r_{\Xi_{\mathtt w},|\widehat\nu_1|}f\leq M^r_{\Phi,|\widehat\nu_1|}f$. Thus, averaging \eqref{eq:lemmma82_gen_in_psi} with $\Psi=\Xi_{\mathtt w}$ over $\mathtt w\in\mathscr G$, we get 
			\[
			\bigl|\widetilde T_{|\widehat\nu_1|,\rho}f-\widetilde{\mathcal T}_{|\widehat\nu_1|,\rho}f\bigr|\lesssim M^r_{\Phi,|\widehat\nu_1|}f, \qquad \text{ for }\rho\geq\ve,
			\]
			and \eqref{eq:maximal_bound_sec7} yields $\|\widetilde{\mathcal T}_{|\widehat\nu_1|,\sharp,\ve}\|\lesssim1+\|\widetilde T_{|\widehat\nu_1|,\sharp,\ve}\|$. Therefore, by \eqref{eq:aux_probab_argument_norms}, the lemma follows from
			\begin{equation}\label{eq:averaged_maximal_cotlar}
				\|\widetilde T_{|\widehat\nu_1|,\sharp,\ve}\|\lesssim 1+\sup_{\rho\geq\ve}\|\widetilde T_{|\widehat\nu_1|,\rho}\|.
			\end{equation}
			
			To prove \eqref{eq:averaged_maximal_cotlar}, we argue as in Section \ref{section:L2bd}, with $\Phi$ in place of $\Xi_{\mathtt w}$ and $\widetilde K$ in place of $K_{\Xi_{\mathtt w}}$. By the end of the proof of Lemma \ref{lemma:final_global_bound}, the function $\Phi$ satisfies the lower bounds of $\Xi_{\mathtt w}$ used in that section, so \eqref{eq:growth_auxiliary_measures_sec6}, \eqref{eq:comparison_auxiliary_measures_sec6}, and \eqref{eq:collar_bound_sec6} hold with $\Phi$ in place of $\Xi_{\mathtt w}$. Moreover, since $\Xi_{\mathtt w}\geq\Phi$, \eqref{eq:averaged_kernel_def} and Lemma \ref{lemma:kernel_suppressed} imply that $\widetilde K$ is a Calder\'on-Zygmund kernel such that
			\[
			|\widetilde K(x,y)|\lesssim\min\bigl\{d(x,y)^{-n},\Phi(x)^{-n},\Phi(y)^{-n}\bigr\}.
			\]
			With the measures $|\nu_i|$, the radius $r_0$, and the cutoffs $\chi_i$, $\psi$, and $\varphi$ of Section \ref{section:L2bd}, we define
			\[
			\widetilde L_\ve(x,y)\coloneqq\chi_2(x)\chi_1(y)\psi\Bigl(\frac{d(x,y)}{r_0}\Bigr)\varphi\Bigl(\frac{d(x,y)}\ve\Bigr)\widetilde K(x,y),
			\]
			and we denote by $\widetilde S_\ve$ and $\widetilde S_{\ve,\rho}$ the integral operators with respect to $|\nu_1|$ with kernels $\widetilde L_\ve$ and $\widetilde L_\ve\chi_{\{d(x,y)>\rho\}}$, respectively. Then \eqref{eq:regularity_auxiliary_kernel_sec6} and \eqref{eq:size_auxiliary_kernel_sec6} hold for $\widetilde L_\ve$, with $\Phi$ in place of $\Xi_{\mathtt w}$ and constants independent of $\ve$.
			
			Let $f\in L^2(|\widehat\nu_1|)$, extended by zero outside $B_1$, let $x\in B_2$, and let $\rho\geq\ve$. We split
			\begin{equation}\label{eq:splitting_averaged_cotlar}
				\begin{split}
					\widetilde T_{|\widehat\nu_1|,\rho}f(x)
					&=\widetilde S_{\ve,\rho}f(x)+\int_{d(x,y)>\rho}\Bigl(1-\varphi\Bigl(\frac{d(x,y)}\ve\Bigr)\Bigr)\widetilde K(x,y)f(y)\,d|\widehat\nu_1|(y)\\
					&\qquad\qquad+\int_{d(x,y)>\rho}\mathfrak D(x,y)f(y)\,d|\widehat\nu_1|(y)
					\eqqcolon\widetilde S_{\ve,\rho}f(x)+\mathfrak J_1(x)+\mathfrak J_2(x),
				\end{split}
			\end{equation}
			where 
			\[
			\mathfrak D(x,y)\coloneqq\Bigl(1-\chi_2(x)\chi_1(y)\psi\Bigl(\frac{d(x,y)}{r_0}\Bigr)\Bigr)\varphi\Bigl(\frac{d(x,y)}{\varepsilon}\Bigr)\, \widetilde K(x,y).
			\]
			
			Since $\varphi(t/\ve)=1$ for $t\geq2\ve$, the bound for $\widetilde K$ gives
			\[
			|\mathfrak J_1(x)|\lesssim\max(\ve,\Phi(x))^{-n}\int_{B(x,2\ve)}|f|\,d|\widehat\nu_1|\lesssim M^r_{\Phi,|\widehat\nu_1|}f(x).
			\]
			As in \eqref{eq:positive_bound_auxiliary_difference_sec6} and \eqref{eq:integral_square_difference_sec6}, we have $|\mathfrak D(x,y)|\lesssim H_\Phi(x,y)\coloneqq(r_0+d(x,y)+\Phi(x)+\Phi(y))^{-n}$ and
			\[
			\iint_{B_2\times B_1}H_\Phi(x,y)^2\,d|\widehat\nu_1|(y)\,d|\widehat\nu_2|(x)\lesssim1.
			\]
			Hence, $|\mathfrak J_2(x)|\lesssim\int H_\Phi(x,y)|f(y)|\,d|\widehat\nu_1|(y)$, and the Cauchy-Schwarz inequality bounds the $L^2(|\widehat\nu_2|)$ norm of the right-hand side by a constant times $\|f\|_{L^2(|\widehat\nu_1|)}$.
			
			To estimate $\widetilde S_{\ve,\rho}$, we note that $\varphi(t)=\int_1^2\varphi'(u)\chi_{\{t>u\}}\,du$ for $t\geq0$, so that
			\[
			\int\varphi\Bigl(\frac{d(x,y)}\ve\Bigr)\widetilde K(x,y)f(y)\,d|\widehat\nu_1|(y)=\int_1^2\varphi'(u)\,\widetilde T_{|\widehat\nu_1|,u\ve}f(x)\,du.
			\]
			Since $\varphi$ is non-decreasing and $\int_1^2\varphi'(u)\,du=1$, the norm of the operator on the left-hand side is at most $\sup_{\rho\geq\ve}\|\widetilde T_{|\widehat\nu_1|,\rho}\|$. Its kernel differs from $\widetilde L_\ve$ by $\mathfrak D$ on $B_2\times B_1$, while $\widetilde L_\ve$ vanishes outside $B_2\times B_1$. Thus, the bound for $H_\Phi$ gives
			\[
			\|\widetilde S_\ve\|_{L^2(|\nu_1|)\to L^2(|\nu_2|)}\lesssim1+\sup_{\rho\geq\ve}\|\widetilde T_{|\widehat\nu_1|,\rho}\|.
			\]
			Since $\widetilde L_\ve(x,y)=0$ for $d(x,y)\leq\ve$, Lemma \ref{lemma:mixed_cotlar_transfer}, applied as after \eqref{eq:size_auxiliary_kernel_sec6} with $\Phi+\ve$ in place of $\Xi$, $\mathfrak m_1=|\nu_1|$, and $\mathfrak m_2=\mathfrak m_3=|\nu_2|$, yields
			\[
			\Bigl\|\sup_{\rho>0}|\widetilde S_{\ve,\rho}f|\Bigr\|_{L^2(|\nu_2|)}
			\lesssim\bigl(1+\|\widetilde S_\ve\|_{L^2(|\nu_1|)\to L^2(|\nu_2|)}\bigr)\|f\|_{L^2(|\widehat\nu_1|)}
			\lesssim\Bigl(1+\sup_{\rho\geq\ve}\|\widetilde T_{|\widehat\nu_1|,\rho}\|\Bigr)\|f\|_{L^2(|\widehat\nu_1|)}.
			\]
			
			Gathering the bounds for the terms in \eqref{eq:splitting_averaged_cotlar}, taking the supremum over $\rho\geq\ve$, and using \eqref{eq:maximal_bound_sec7} with $\Psi=\Phi$, we obtain \eqref{eq:averaged_maximal_cotlar}.
		\end{proof}
		\vvv
	}

	We now prove the uniform $L^2(|\widehat\nu_1|)\to L^2(|\widehat\nu_2|)$ bounds for $\widetilde T$, $T_\Phi$, and $T_{\Xi_{\mathtt w}}$.
	\begin{lemma}\label{lemma:final_boundedness_suppressed_operators}
		There exists $c>0$ depending on the structural parameters of the theorem such that, for every $\delta>0$, we have
		\[
		\begin{split}
			&\|\widetilde T_{|\widehat\nu_1|,\delta}\|_{L^2(|\widehat\nu_1|)\to L^2(|\widehat\nu_2|)}
			+\|T_{\Phi,|\widehat\nu_1|,\delta}\|_{L^2(|\widehat\nu_1|)\to L^2(|\widehat\nu_2|)}+\sup_{\mathtt w\in\mathscr G}\|T_{\Xi_{\mathtt w},|\widehat\nu_1|,\delta}\|_{L^2(|\widehat\nu_1|)\to L^2(|\widehat\nu_2|)}\leq c,
		\end{split}
		\]
		where $c$ is independent of $\delta$, $R_1$, and $R_2$.
	\end{lemma}
	\begin{proof}
		For $\delta>0$, we set
		\[
		A_\delta\coloneqq\sup_{\rho\geq\delta}\|\widetilde T_{|\widehat\nu_1|,\rho}\|_{L^2(|\widehat\nu_1|)\to L^2(|\widehat\nu_2|)}.
		\]
		Since $|\widetilde K(x,y)|\lesssim d(x,y)^{-n}$, the Cauchy-Schwarz inequality gives 
		\[
		A_\delta\lesssim\delta^{-n}|\widehat\nu_1|(B_1)^{1/2}|\widehat\nu_2|(B_2)^{1/2}<\infty.
		\]
		
		Let $C'>0$ be the implicit constant in \eqref{eq:T-xi<tilde T bis}, which, by Lemma \ref{lemma:comparison_averaged_operators}, is independent of $m$. Then, for $\varepsilon\geq\delta$, the quantity $N_\varepsilon$ in \eqref{eq:definition_N_epsilon_sec6} satisfies
		\begin{equation}\label{eq:sec8_uniform_norm_comparison}
			N_\varepsilon\leq C'\, (1+A_\delta).
		\end{equation}
		
		Let $\varepsilon\geq\delta$, and let $f\in L^2(|\widehat\nu_1|)$ and $g\in L^2(|\widehat\nu_2|)$ have norm one. We extend them by zero outside $B_1$ and $B_2$, respectively, so that $\|f\|_{L^2(|\nu_1|)}=\|g\|_{L^2(|\nu_2|)}=1$.
		
		Let $\mathtt w\in\mathscr G$, and let $T_{\mathtt w}\coloneqq T_{\Xi_{\mathtt w},|\widehat\nu_1|,\varepsilon}$. Since $f$ and $g$ are supported in $Q_2(w_1)$ and $Q_2(w_2)$, respectively, \eqref{eq:dec_f_mart} gives $f=F_{\mathrm{good}}+F_{\mathrm{bad}}$ and $g=G_{\mathrm{good}}+G_{\mathrm{bad}}$, where
		\[
		\begin{aligned}
			F_{\mathrm{good}}&\coloneqq E_1f+\mathcal P^1_{\mathcal G}(\mathtt w)f,
			&\qquad F_{\mathrm{bad}}&\coloneqq\mathcal P^1_{\mathcal B}(\mathtt w)f,\\
			G_{\mathrm{good}}&\coloneqq E_2g+\mathcal P^2_{\mathcal G}(\mathtt w)g,
			&\qquad G_{\mathrm{bad}}&\coloneqq\mathcal P^2_{\mathcal B}(\mathtt w)g,
		\end{aligned}
		\]
		where we recall that $E_1$ and $E_2$ denote the top-averaging operators associated with the martingale difference decomposition construction, see \eqref{eq:top_averages_def}, and projections are introduced in Definition \ref{def:p_good_p_bad_not_truncated}.
		
		Note that the functions above need not be supported in the testing balls, so we set $f_{\mathrm{good}}\coloneqq\chi_{B_1}F_{\mathrm{good}}$ and $f_{\mathrm{bad}}\coloneqq\chi_{B_1}F_{\mathrm{bad}}$, and we define $g_{\mathrm{good}}$ and $g_{\mathrm{bad}}$ analogously. Then $f=f_{\mathrm{good}}+f_{\mathrm{bad}}$ and $g=g_{\mathrm{good}}+g_{\mathrm{bad}}$, and we split
		\begin{equation}\label{eq:probabilistic_splitting_1}
			\begin{split}
				\langle T_{\mathtt w}f,g\rangle_{|\widehat\nu_2|}
				&=\langle T_{\mathtt w}f_{\mathrm{good}},g_{\mathrm{good}}\rangle_{|\widehat\nu_2|}
				+\langle T_{\mathtt w}f_{\mathrm{good}},g_{\mathrm{bad}}\rangle_{|\widehat\nu_2|}
				+\langle T_{\mathtt w}f_{\mathrm{bad}},g\rangle_{|\widehat\nu_2|}\\
				&\eqqcolon\mathfrak J_1(\mathtt w)+\mathfrak J_2(\mathtt w)+\mathfrak J_3(\mathtt w).
			\end{split}
		\end{equation}
		
		Since $\|\chi_{B_i}h\|_{L^2(|\widehat\nu_i|)}\leq\|h\|_{L^2(|\nu_i|)}$, we have
		\[
		|\mathfrak J_2(\mathtt w)|+|\mathfrak J_3(\mathtt w)|
		\leq N_\varepsilon\, \bigl(\|F_{\mathrm{good}}\|_{L^2(|\nu_1|)}\|\mathcal P^2_{\mathcal B}(\mathtt w)g\|_{L^2(|\nu_2|)}
		+\|\mathcal P^1_{\mathcal B}(\mathtt w)f\|_{L^2(|\nu_1|)}\bigr).
		\]
		Hence, by \eqref{eq:dec_f_mart_2}, \eqref{eq:bound_exp_good}, \eqref{eq:bound_exp_bad}, and $\mathbb P^{\Omega^2}(\mathscr G)>1-\kappa>1/2$, we obtain
		\begin{equation}\label{eq:sec8_bad_pairing_bound}
			\widetilde{\mathbb E}_{\mathtt w}\bigl(|\mathfrak J_2(\mathtt w)|+|\mathfrak J_3(\mathtt w)|\bigr)
			\leq C_{\mathrm{bad}}\,\varepsilon_{\mathcal B}(m)N_\varepsilon,
		\end{equation}
		where $C_{\mathrm{bad}}>0$ is independent of $m$.
		
		To bound $\mathfrak J_1(\mathtt w)$, we note that 
		\[
		F_{\mathrm{good}}-E_1F_{\mathrm{good}}=\mathcal P^1_{\mathcal G}(\mathtt w)f
		\qquad \text{ and }\qquad G_{\mathrm{good}}-E_2G_{\mathrm{good}}=\mathcal P^2_{\mathcal G}(\mathtt w)g.
		\]
		
		Hence, \eqref{eq:localized_pairing_comparison_sec5}, applied twice, the splitting \eqref{eq:martingale_difference_decomposition_extended} with $F_{\mathrm{good}}$ and $G_{\mathrm{good}}$ in place of $f_{\mathtt w}$ and $g_{\mathtt w}$, and the bounds \eqref{eq:final_bound_I_1}, \eqref{eq:inbalanced_estimate_2}, and \eqref{eq:final_bound_I_3_alternative} for its first three terms give
		\[
		\bigl|\mathfrak J_1(\mathtt w)-\langle T_{\mathtt w}(\mathcal P^1_{\mathcal G}(\mathtt w)f),\mathcal P^2_{\mathcal G}(\mathtt w)g\rangle_{|\widehat\nu_2|}\bigr|\lesssim1,
		\]
		where we also used \eqref{eq:dec_f_mart_2} and \eqref{eq:bound_exp_good}. Thus, the Key Lemma \ref{lem:key_lemma_good_functions} gives
		\begin{equation}\label{eq:sec8_good_pairing_bound}
			\widetilde{\mathbb E}_{\mathtt w}|\mathfrak J_1(\mathtt w)|
			\leq C_{m,\lambda,\upsilon}\, (1+a^{-n})+\beta_m(\lambda,\upsilon,a)\, (1+N_\varepsilon),
		\end{equation}
		where, increasing $C_{m,\lambda,\upsilon}$ if necessary, we set
		\[
		\beta_m(\lambda,\upsilon,a)\coloneqq C_m\Bigl(\frac{(2m+1)\lambda^\eta+\upsilon}{\mathbb P^{\Omega^2}(\mathscr G)}\Bigr)^{1/2}+C_{m,\upsilon}(M_\Lambda a)^{1/2}.
		\]
		
		We now fix the parameters, independently of $\delta$ and $\varepsilon$. Since $C'$ and $C_{\mathrm{bad}}$ are independent of $m$, we first choose $m$ such that \eqref{eq:large_probability_top_cubes} holds and $C'C_{\mathrm{bad}}\varepsilon_{\mathcal B}(m)\leq1/4$. Moreover, $C_m$ is independent of $\lambda$ and $\upsilon$, so we can choose $\lambda$ and $\upsilon$ such that
		\[
		C'C_m\Bigl(\frac{(2m+1)\lambda^\eta+\upsilon}{\mathbb P^{\Omega^2}(\mathscr G)}\Bigr)^{1/2}\leq\frac18.
		\]
		Finally, since $M_\Lambda$ and $C_{m,\upsilon}$ are now fixed, we choose $a\in(0,1/4)$ such that $C'C_{m,\upsilon}(M_\Lambda a)^{1/2}\leq1/8$.
		
		By \eqref{eq:definition_averaged_T} and Fubini's theorem, we have
		\[
		\langle\widetilde T_{|\widehat\nu_1|,\varepsilon}f,g\rangle_{|\widehat\nu_2|}=\widetilde{\mathbb E}_{\mathtt w}\langle T_{\mathtt w}f,g\rangle_{|\widehat\nu_2|}.
		\]
		Hence, \eqref{eq:probabilistic_splitting_1}, \eqref{eq:sec8_bad_pairing_bound}, \eqref{eq:sec8_good_pairing_bound}, and \eqref{eq:sec8_uniform_norm_comparison} give, after taking the supremum over $f$, $g$, and $\varepsilon\geq\delta$,
		\begin{equation}\label{eq:absorbing_terms}
			A_\delta\leq C_*+C'\bigl(\beta_m(\lambda,\upsilon,a)+C_{\mathrm{bad}}\varepsilon_{\mathcal B}(m)\bigr)\, A_\delta
			\leq C_*+\frac12A_\delta,
		\end{equation}
		where $C_*$ depends on $m$, $\lambda$, $\upsilon$, and $a$, but not on $\delta$. Since $A_\delta<\infty$, we obtain $A_\delta\leq2C_*$, and the lemma follows from \eqref{eq:T-phi<tilde T bis} and \eqref{eq:T-xi<tilde T bis}.
	\end{proof}
	\vv

	\subsection{The good-lambda argument}\label{subsection:good_lambda}
	We use the dyadic system of Lemma \ref{lemma:DM_cubes_one_sided} and a good-lambda inequality to prove the following result.
	\begin{proposition}\label{theorem:MMT_theorem_A1}
		Let $\mu$ be a Radon measure on $X$ that has $n$-growth with constant $\mathfrak C_0$, let $\theta_0\in (0,1]$, and let $T$ be a singular integral operator associated with an $n$-dimensional Calder\'on-Zygmund kernel with parameters $\gamma$ and $C$. Suppose that for every $s>0$ and every cell $Q\in \caD^{db, (s)}_\mu$	there exists a Borel subset $\mathcal G_Q\subset Q$ with $\mu(\mathcal G_Q)\geq \theta_0 \, \mu(Q)$ such that $T_{\mu|_{\mathcal G_Q}}$ is a bounded operator from $L^2(\mu|_{\mathcal G_Q})$ to $L^2(\mu|_{\mathcal G_Q})$, with norm bounded uniformly in $Q$ and $s$.
		
		Then, if the parameter $K_0$ in the construction of the David-Mattila lattice is chosen big enough depending on $D$, for every $f\in L^2(\mu)$ we have
		\begin{equation}\label{eq:L2_boundedness_good_lambda}
			\sup_{\varepsilon>0}\|T_{\mu,\varepsilon}f\|_{L^2(\mu)}\lesssim \|f\|_{L^2(\mu)},
		\end{equation}
		where the implicit constant depends only on $D$, $n$, $\mathfrak C_0$, $C$, $\gamma$, $\theta_0$, the lattice parameters $K_0$ and $A_0$, and the uniform bound for the operators $T_{\mu|_{\mathcal G_Q}}$.
	\end{proposition}

	\vv
	By density and Fatou's lemma, it suffices to prove \eqref{eq:L2_boundedness_good_lambda} for $f\in L^2(\mu)$ with bounded support. Fix $h>0$ and write
	\[
	T_{\mu,\sharp,h}f(x)\coloneqq\sup_{\rho\geq h}|T_{\mu,\rho}f(x)|.
	\]
	If $\supp f\subset B(x_0,R)$, standard kernel estimates give $T_{\mu,\sharp,h}f\leq Ch^{-n}\|f\|_{L^1(\mu)}$ and
	\[
	\int_{X\setminus B(x_0,2R)}(T_{\mu,\sharp,h}f)^2\,d\mu
	\lesssim\|f\|_{L^1(\mu)}^2\sum_{j\geq1}(2^jR)^{-2n}\mu(B(x_0,2^{j+1}R))<\infty.
	\]
	Since $\mu(B(x_0,2R))<\infty$, we deduce that $T_{\mu,\sharp,h}f\in L^2(\mu)$. We prove, uniformly in $h$, that for every $0<\varepsilon<1$ there is $\varepsilon'>0$ such that, for every $\lambda>0$,
	\begin{equation}\label{eq:reduction_1_big_lambda}
		\begin{split}
			&\mu\bigl(\{x:T_{\mu,\sharp,h}f(x)>(1+\varepsilon)\lambda,
			\ \widetilde{\mathcal M}_\mu f(x)\leq\varepsilon'\lambda\}\bigr)\\
			&\qquad\qquad\qquad\leq\Bigl(1-\frac{\theta_0}{4}\Bigr)\,\mu\bigl(\{x:T_{\mu,\sharp,h}f(x)>\lambda\}\bigr).
		\end{split}
	\end{equation}
	We follow the scheme of the proof of \cite[Section 9]{MT20}, using the cells of $\caD^{(s)}_\mu$.

	For $s>0$ and $Q\in\caD^{db,(s)}_\mu$, we set $\sigma_Q\coloneqq\mu|_{\mathcal G_Q}$. The measure $\sigma_Q$ has $n$-growth and is nonatomic. Thus the metric weak-type theorem \cite[Corollary 1.7(ii)]{HLYY12}, applied to the $L^2(\sigma_Q)$-bounded operator $T$, gives
	\begin{equation}\label{eq:big_piece_measure_weak11_sec8}
		\mu\bigl(\{x\in\mathcal G_Q:T_\sharp\omega(x)>t\}\bigr)
		\leq\frac{C}{t}\, \|\omega\|,
		\qquad t>0,
	\end{equation}
	for every finite complex Borel measure $\omega$ on $X$. The constant is uniform in $Q$ and $s$. We remark that $\omega$ need not be supported on $\mathcal G_Q$.

	We first prove the following two lemmas.
	\begin{lemma}\label{lemma:level_sets_good_lambda}
		Let $(X,d)$ be a geometrically doubling metric space, and let $\mu$ be a Radon measure on $X$ that has $n$-growth with parameter $\mathfrak C_0$. Let $T$ be a singular integral operator associated with a Calder\'on-Zygmund kernel with parameters $C>0$ and $\gamma\in(0,1)$. For a function $f\in L^2(\mu)$ with bounded support and $h,\lambda>0$, we define
		\[
		\Omega_\lambda\coloneqq \bigl\{x\in \supp(\mu):T_{\mu,\sharp,h}f(x)>\lambda\bigr\}.
		\]
		
		For every $\lambda>0$, the set $\Omega_\lambda$ is bounded, $\mu(\Omega_\lambda)<\infty$, and $\Omega_\lambda$ is relatively open in $\supp(\mu)$.
	\end{lemma}
	\begin{proof}
		By the assumptions on $f$, there exist $x_0\in X$ and $R>0$ such that $f\equiv 0$ outside $B(x_0,R)$. Moreover, $f\in L^1(\mu)$ so, for every $x\in X\setminus \overline{B(x_0,R)}$ and $\rho>0$, we have
		\[
		|T_{\mu,\rho}f(x)|\leq \frac{C}{\dist(x,\overline{B(x_0,R)})^n}\, \int |f|\, d\mu,
		\]
		which implies that $\Omega_\lambda\subset B(x_0,R_\lambda)$ for a sufficiently large $R_\lambda>R$, depending on $f$ and $\lambda$. Thus $\Omega_\lambda$ is bounded and
		\[
		\mu(\Omega_\lambda)\leq \mu(B(x_0,R_\lambda))\leq \mathfrak C_0\, R_\lambda^n<\infty.
		\]
		
		It remains to prove that $\Omega_\lambda$ is relatively open. Let $x\in \Omega_\lambda$. By definition, there exists $\rho\geq h$ such that $|T_{\mu,\rho}f(x)|>\lambda$. Since $T_{\mu,r}f(x)\to T_{\mu,\rho}f(x)$ as $r\downarrow\rho$, and since the sphere $\{y:d(x,y)=r\}$ has null $|f|\mu$-measure for all but countably many $r>0$, there exists $r>\rho$ such that $|T_{\mu,r}f(x)|>\lambda$ and
		\[
		\int_{d(x,y)=r} |f(y)|\, d\mu(y)=0.
		\]
		
		If $z\to x$, then $K(z,y)\chi_{\{d(z,y)>r\}}f(y)\to K(x,y)\chi_{\{d(x,y)>r\}}f(y)$ for $(|f|\mu)$-almost every $y$, and the integrands are bounded by $Cr^{-n}|f(y)|$. By dominated convergence, $T_{\mu,r}f(z)\to T_{\mu,r}f(x)$. Hence there is a neighborhood $U$ of $x$ such that $|T_{\mu,r}f(z)|>\lambda$ for $z\in U$. It follows that $U\cap\supp(\mu)\subset\Omega_\lambda$, which completes the proof of the lemma.
	\end{proof}

	If $\Omega_\lambda=\varnothing$, then \eqref{eq:reduction_1_big_lambda} is immediate. Suppose first that $\varnothing\neq\Omega_\lambda\subsetneq\supp(\mu)$. This includes the unbounded-support case, because $\Omega_\lambda$ is bounded. Choose $\bar x\in\supp(\mu)\setminus\Omega_\lambda$ and set
	\begin{equation}\label{eq:choice_of_s}
		s\coloneqq2\sup_{x\in\Omega_\lambda}d(x,\bar x)>0.
	\end{equation}
	
	If a top-generation cell $Q\in\caD^{(s)}_{\mu,0}$ were contained in $\Omega_\lambda$, its center $z_Q$ would belong to $\Omega_\lambda$. Hence \[d(z_Q,\bar x)\leq s/2<s\leq r^{(s)}(Q),\] and \eqref{eq:new_radii_modified_DM} would give $\bar x\in\supp(\mu)\cap B^{(s)}(Q)\subset Q\subset\Omega_\lambda$, which is a contradiction. Thus no top-generation cell is contained in $\Omega_\lambda$, and the following decomposition applies. The remaining case $\Omega_\lambda=\supp(\mu)$ will be treated separately.
	\begin{lemma}\label{lemma:whitney_DM_s}
		Let $\Omega\subset \supp(\mu)$ be a non-empty bounded set, relatively open in $\supp(\mu)$, and $s>0$ be such that no cell of $\caD^{(s)}_{\mu,0}$ is contained in $\Omega$.
		Then, there exists a family $\mathcal W^{(s)}\subset \caD^{(s)}_\mu$ of pairwise disjoint cells such that
		\[
		\Omega =\bigcup_{Q\in\mathcal W^{(s)}}Q,
		\]
		and such that, for some constants $T_0>10^4$, depending on $A_0$ and $K_0$, and $N_0\geq1$, depending on $D$, $A_0$, and $K_0$, the following properties hold:
		\begin{enumerate}
			\item Every  $Q\in \mathcal W^{(s)}$ has a dyadic parent $\widehat Q$ in $\caD^{(s)}_\mu$ such that $\widehat Q\subset \Omega$.
			\item For every $Q\in \mathcal W^{(s)}$, we have $\supp(\mu)\cap 10^4\, B^{(s)}(Q)\subset\Omega$.
			\item For every $Q\in \mathcal W^{(s)}$, there exists $y_Q\in \supp(\mu)\setminus \Omega$ such that 
			\[
			d(x,y_Q)\leq T_0\, r^{(s)}(Q), \qquad \text{ for every }x\in Q.
			\]
			\item For every $Q\in \mathcal W^{(s)}$, there are at most $N_0$ cells $P\in \mathcal W^{(s)}$	such that
			\(
			10^4B^{(s)}(Q)\cap10^4B^{(s)}(P)\neq\varnothing.
			\)
			Further, for such cells $Q$, $P$, we have
			\begin{equation}\label{eq:comparison_l_whitney_DM_s}
				(4A_0)^{-1}\, \ell^{(s)}(P)\leq \ell^{(s)}(Q)\leq 4A_0\,  \ell^{(s)}(P).
			\end{equation}
			\item The family
			\[
			\mathcal W^{(s)}_{db}\coloneqq\mathcal W^{(s)}\cap \caD_\mu^{db,(s)}
			\]
			satisfies
			\begin{equation}\label{eq:big_piece_db_whitney_DM_aux}
				\mu\biggl(\bigcup_{Q\in \mathcal W^{(s)}_{db}}Q\biggr)\geq \frac{1}{2}\,\mu(\Omega),
			\end{equation}
			assuming the parameter $K_0$ in the construction of the David-Mattila lattice is chosen big enough depending only on $D$.
		\end{enumerate}
	\end{lemma}
	\begin{proof}
		Fix a constant $0<\mathfrak t<1$, depending only on the fixed lattice parameters, such that $\mathfrak t\leq\min\{(2A_0)^{-1},14\cdot10^{-4}\}$. We say that $Q\in \caD^{(s)}_\mu$ is $\mathfrak t$-admissible if
		\[
		\ell^{(s)}(Q)\leq \mathfrak t\, \dist (Q,\supp(\mu)\setminus \Omega).
		\]
		In particular, every $\mathfrak t$-admissible cell is contained in $\Omega$. Since $Q\subset28B^{(s)}(Q)$, we have $\diam(Q)\leq56\,r^{(s)}(Q)\leq\ell^{(s)}(Q)$ for every $Q\in\caD^{(s)}_\mu$.
		
		For $k\geq 0$ and $x\in \Omega$, we denote by $Q_k(x)$ the cell from $\caD^{(s)}_{\mu,k}$ that contains $x$. Since $\diam(Q_k(x))\leq\ell^{(s)}(Q_k(x))$, the triangle inequality gives
		\[
		\dist(Q_k(x), \supp(\mu)\setminus \Omega)\geq \dist(x,\supp(\mu)\setminus \Omega)-\ell^{(s)}(Q_k(x)).
		\]
		Since $\Omega$ is relatively open in $\supp(\mu)$, we have $\dist(x,\supp(\mu)\setminus\Omega)>0$. Hence, since $\ell^{(s)}(Q_k(x))\to0$ as $k\to\infty$, for $k$ sufficiently large we also have \[\ell^{(s)}(Q_k(x))\leq\mathfrak t(1+\mathfrak t)^{-1}\dist(x,\supp(\mu)\setminus\Omega).\] Combining the two inequalities gives
		\[
		\ell^{(s)}(Q_k(x))\leq\mathfrak t\,\dist(Q_k(x),\supp(\mu)\setminus\Omega).
		\]
		
		Thus, for $k$ large enough, the cell $Q_k(x)$ is $\mathfrak t$-admissible. We define $\mathcal W^{(s)}$ as the family of $\mathfrak t$-admissible cells having no $\mathfrak t$-admissible ancestor of smaller generation. This definition also applies when a cell and one of its ancestors coincide as subsets of $X$. The cells in $\mathcal W^{(s)}$ are pairwise disjoint and, since every $\mathfrak t$-admissible cell is contained in $\Omega$, their union is $\Omega$. We now prove (1)-(5). 
		
		By assumption, no cell from $\caD^{(s)}_{\mu,0}$ is admissible, so each $Q\in\mathcal W^{(s)}$ has a parent $\widehat Q$. Moreover,
		\[
		\diam(\widehat Q)\leq A_0\ell^{(s)}(Q)\leq A_0\mathfrak t\,\dist(Q,\supp(\mu)\setminus\Omega)<\dist(Q,\supp(\mu)\setminus\Omega).
		\]
		Since $Q\subset\widehat Q$, this implies $\widehat Q\subset\Omega$ and proves (1).
		
		To prove (2), let $Q\in \mathcal W^{(s)}$ and $z\in \supp(\mu)\cap 10^4\, B^{(s)}(Q)$. Since $r^{(s)}(Q)\leq\ell^{(s)}(Q)/56$ and $\mathfrak t\leq14\cdot10^{-4}$, the admissibility of $Q$ gives
		\begin{equation*}
			\begin{split}
				d(z,z_Q)< 10^4\, r^{(s)}(Q)&\leq \frac{10^4}{56}\, \ell^{(s)}(Q)\leq \frac{10^4}{56}\, \mathfrak t\, \dist(Q,\supp(\mu)\setminus \Omega)\\
				&\leq \frac{1}{4}\, \dist(Q,\supp(\mu)\setminus \Omega)\leq \frac{1}{4}\, \dist(z_Q,\supp(\mu)\setminus \Omega),
			\end{split}
		\end{equation*}
		which proves that $z\in \Omega$.
		
		Now, we prove (3). Let $Q\in \mathcal W^{(s)}$ and denote by $\widehat Q$ its dyadic parent in $\caD^{(s)}_\mu$. Since $\widehat Q$ is not admissible, there exist $y'\in \widehat Q$ and $y_Q\in \supp(\mu)\setminus \Omega$ such that 
		\[
		d(y_Q,y')<\frac{\ell^{(s)}(\widehat Q)}{\mathfrak t}= A_0\, \frac{\ell^{(s)}( Q)}{\mathfrak t}.
		\]
		Thus, by the triangle inequality, for every $x\in Q$ we have
		\begin{equation*}
			\begin{split}
				d(x,y_Q)&\leq \diam(\widehat Q) + d(y_Q,y')< A_0\, \ell^{(s)}(Q) + \frac{A_0\, \ell^{(s)}(Q)}{\mathfrak t}\\
				&\leq 56 A_0\, K_0\, (1+\mathfrak t^{-1})\, r^{(s)}(Q)\eqqcolon T_0\, r^{(s)}(Q).
			\end{split}
		\end{equation*}

		To prove (4), let $Q,P\in \mathcal W^{(s)}$ be such that $10^4B^{(s)}(Q)\cap 10^4 B^{(s)}(P)\neq \varnothing$. By the maximality of $Q,P$ in the construction of $\mathcal W^{(s)}$, both cells are admissible but their dyadic parents are not. Hence, by the triangle inequality, the proof of (3), and the choice of $\mathfrak t$, we obtain
		\begin{equation*}
			\begin{split}
				\ell^{(s)}(Q)&\leq\mathfrak t\, \dist(z_Q, \supp(\mu)\setminus \Omega)\\
				&\leq \mathfrak t\, \bigl(d(z_Q,z_P) + \dist(z_P, \supp(\mu)\setminus \Omega)\bigr)\\
				&\leq \mathfrak t\, \Bigl(\frac{10^4}{56}\bigl(\ell^{(s)}(Q)+ \ell^{(s)}(P)\bigr) + A_0(1+\mathfrak t^{-1})\ell^{(s)}(P)\Bigr)\\
				&\leq \frac{1}{4}\, \ell^{(s)}(Q) +  \frac{1}{4}\, \ell^{(s)}(P) + 2 A_0\, \ell^{(s)}(P),
			\end{split}
		\end{equation*}
		which implies that $\ell^{(s)}(Q)\leq 4A_0\,  \ell^{(s)}(P)$. Reversing the roles of $Q$ and $P$, we obtain $\ell^{(s)}(P)\leq 4A_0\,  \ell^{(s)}(Q)$, which proves \eqref{eq:comparison_l_whitney_DM_s}. In particular, \eqref{eq:comparison_l_whitney_DM_s} implies that the generations of $Q$ and $P$ can differ by at most $1$. Then, since centers of cells in $\caD^{(s)}_\mu$ are well-separated at each fixed generation, the bounded overlap property follows from a standard argument and the geometrically doubling property of $X$.
		
		We now prove (5). We define $\mathcal W^{(s)}_{nd}\coloneqq \mathcal W^{(s)}\setminus \mathcal W^{(s)}_{db}$. Since $K_0\geq100$, \eqref{ndb_sec8} with $\ell=1$ gives, for $Q\in\mathcal W^{(s)}_{nd}$, 
		\begin{equation}\label{eq:l87_a}
			\mu(Q)\leq \mu(100\, B^{(s)}(Q))\leq K_0^{-1}\, \mu(10^4\, B^{(s)}(Q)).
		\end{equation}
		We claim that there exists $C_D>0$ depending only on the geometrically doubling constant $D$ such that we have the bounded overlap property
		\begin{equation}\label{eq:bounded_overlap_whitney8}
			\sum_{Q\in \mathcal W^{(s)}_{nd}}\chi_{10^4\, B^{(s)}(Q)}\leq C_D.
		\end{equation}
		Note that item (2), \eqref{eq:l87_a} and \eqref{eq:bounded_overlap_whitney8} readily imply \eqref{eq:big_piece_db_whitney_DM_aux}. Indeed, 
		\begin{equation*}
			\begin{split}
				\mu\biggl(\bigcup_{Q\in \mathcal W^{(s)}_{db}}Q\biggr)&= \mu(\Omega) - \sum_{Q\in \mathcal W^{(s)}_{nd}}\mu(Q)\geq \mu(\Omega) -\frac{1}{K_0} \sum_{Q\in \mathcal W^{(s)}_{nd}}\mu(10^4\, B^{(s)}(Q))\\
				&\geq \mu(\Omega) - \frac{1}{K_0}\int_{\Omega} 	\sum_{Q\in \mathcal W^{(s)}_{nd}}\chi_{10^4\, B^{(s)}(Q)}\, d\mu\geq \Bigl(1-\frac{C_D}{K_0}\Bigr)\, \mu(\Omega)>\frac{1}{2}\, \mu(\Omega),
			\end{split}
		\end{equation*}
		provided that $C_D/K_0<1/2$.
		
		Hence, we are left only with proving \eqref{eq:bounded_overlap_whitney8}. Fix $x\in X$. 
		If $Q\in\mathcal W^{(s)}_{nd}$ is a cell of generation $k$ such that $x\in 10^4\, B^{(s)}(Q)$, then $r^{(s)}(Q)=s\, A^{-k}_0$ and its center belongs to $B(x,10^4\, s\, A_0^{-k})$. Furthermore, if $Q'$ is another cell  in $\mathcal W^{(s)}_{nd}$ of generation $k$ such that $x\in 10^4\, B^{(s)}(Q')$, then by Lemma \ref{lemma:DM_cubes_one_sided} we have $5\, B^{(s)}(Q)\cap 5\, B^{(s)}(Q')=\varnothing$. Thus, since $X$ is geometrically doubling, at each fixed generation there are at most $D^{\ceil{\log_2(10^4)}}$ balls whose contribution to the sum in \eqref{eq:bounded_overlap_whitney8} is nonzero.
		
		The claim \eqref{eq:bounded_overlap_whitney8} follows with $C_D\coloneqq 2\, D^{\ceil{\log_2(10^4)}}$ noting that if two cells $P,P'\in \mathcal W^{(s)}_{nd}$ are such that $x\in 10^4\, B^{(s)}(P)\cap 10^4\, B^{(s)}(P')$, then their generations differ by at most one.
	\end{proof}
	
	\vvv
	
	We now prove \eqref{eq:reduction_1_big_lambda}, which completes the proof of Proposition \ref{theorem:MMT_theorem_A1}. We denote by $E_\lambda$ the set on its left-hand side, and we note that its right-hand side equals $(1-\theta_0/4)\mu(\Omega_\lambda)$.
	
	We first assume that $\varnothing\neq\Omega_\lambda\subsetneq\supp(\mu)$, we let $s$ be as in \eqref{eq:choice_of_s}, and we denote by $\mathcal W^{(s)}$ the family given by Lemma \ref{lemma:whitney_DM_s} with $\Omega=\Omega_\lambda$. Let $Q\in\mathcal W^{(s)}_{db}$ and $r\coloneqq r^{(s)}(Q)$. Since $Q\subset28B^{(s)}(Q)$, there exist $N_D$ balls $V_j\coloneqq B(c_j,r)$, with $c_j\in Q$ and $N_D$ depending only on $D$, that cover $Q$.
	
	Let $j$ be such that $V_j\cap E_\lambda\cap\mathcal G_Q\neq\varnothing$, let $x_j\in V_j\cap E_\lambda\cap\mathcal G_Q$, and let $f_j\coloneqq f\chi_{B(c_j,4r)}$. Since $x_j\in B(c_j,4r)$, $\widetilde{\mathcal M}_\mu f(x_j)\leq\varepsilon'\lambda$, $B(c_j,20r)\subset100B^{(s)}(Q)$, and $Q$ is doubling, we have
	\begin{equation}\label{eq:goodlambda_local_mass_sec8}
		\|f_j\|_{L^1(\mu)}
		\leq\varepsilon'\lambda\,\mu(B(c_j,20r))
		\leq K_0\varepsilon'\lambda\,\mu(Q).
	\end{equation}
	
	Let $y_Q$ be as in item (3) of Lemma \ref{lemma:whitney_DM_s}, let $x\in Q\cap V_j$ and $\rho\geq h$, and let $u\coloneqq\max\{\rho,4(T_0+2)r\}$ and $B_j\coloneqq B(x_j,4(T_0+2)r)$. Since $B(c_j,4r)\subset B_j$, we split
	\begin{equation}\label{eq:splitting_far_comparison_sec8}
		\begin{split}
			T_{\mu,\rho}(f-f_j)(x)-T_{\mu,u}f(y_Q)
			&=\int_{B_j}\bigl(K(x,z)\chi_{A}(z)-K(y_Q,z)\chi_{\{d(y_Q,z)>u\}}(z)\bigr)f(z)\,d\mu(z)\\
			&\qquad+\int_{X\setminus B_j}\bigl(K(x,z)-K(y_Q,z)\bigr)\chi_{\{d(y_Q,z)>u\}}(z)f(z)\,d\mu(z)\\
			&\qquad+\int_{X\setminus B_j}K(x,z)\bigl(\chi_{\{d(x,z)>\rho\}}(z)-\chi_{\{d(y_Q,z)>u\}}(z)\bigr)f(z)\,d\mu(z)\\
			&\eqqcolon\mathfrak J_1+\mathfrak J_2+\mathfrak J_3,
		\end{split}
	\end{equation}
	where $A\coloneqq\{z:d(x,z)>\rho\}\setminus B(c_j,4r)$. We have $d(x,x_j)<2r$, $d(x,z)\geq3r$ for $z\notin B(c_j,4r)$, and $d(x,y_Q)\leq T_0r$. Moreover, the definition of $\widetilde{\mathcal M}_\mu$ and the $n$-growth of $\mu$ give
	\begin{equation}\label{eq:local_mass_maximal_sec8}
		\int_{B(x_j,t)}|f|\,d\mu\leq\widetilde{\mathcal M}_\mu f(x_j)\,\mu(B(x_j,5t))\lesssim t^n\,\widetilde{\mathcal M}_\mu f(x_j),\qquad t>0.
	\end{equation}
	
	Since the kernels in $\mathfrak J_1$ are bounded by $Cr^{-n}$, \eqref{eq:local_mass_maximal_sec8} gives $|\mathfrak J_1|\lesssim\widetilde{\mathcal M}_\mu f(x_j)$.
	
	If $z\notin B_j$, then $d(x,z)\geq d(x_j,z)/2\geq2T_0r\geq2d(x,y_Q)$. Hence, the regularity estimate and \eqref{eq:local_mass_maximal_sec8}, applied on the annuli $B(x_j,2^{k+1}t_j)\setminus B(x_j,2^kt_j)$, with $t_j\coloneqq4(T_0+2)r$ and $k\geq0$, give
	\[
	|\mathfrak J_2|\lesssim\int_{X\setminus B_j}\frac{r^\gamma|f(z)|}{d(x_j,z)^{n+\gamma}}\,d\mu(z)
	\lesssim\sum_{k\geq0}2^{-k\gamma}\,\widetilde{\mathcal M}_\mu f(x_j)\lesssim\widetilde{\mathcal M}_\mu f(x_j).
	\]
	Moreover $|d(x,z)-d(y_Q,z)|\leq T_0r\leq u/4$, and $d(x,z)\geq2(T_0+2)r$ for $z\notin B_j$, so the integrand in $\mathfrak J_3$ vanishes unless $z\in B(x_j,2u)$ and $d(x,z)\geq u/2$. Thus, \eqref{eq:local_mass_maximal_sec8} gives 
	\[
	|\mathfrak J_3|\lesssim u^{-n}\int_{B(x_j,2u)}|f|\,d\mu\lesssim\widetilde{\mathcal M}_\mu f(x_j).
	\]
	
	Gathering the bounds for the terms in \eqref{eq:splitting_far_comparison_sec8}, we obtain
	\begin{equation}\label{eq:goodlambda_far_comparison_sec8}
		\bigl|T_{\mu,\rho}(f-f_j)(x)-T_{\mu,u}f(y_Q)\bigr|
		\leq C\,\widetilde{\mathcal M}_\mu f(x_j),
	\end{equation}
	where $C$ is independent of $r$, $h$, and $\rho$.
	
	Since $u\geq h$, $y_Q\in\supp(\mu)\setminus\Omega_\lambda$, and $\widetilde{\mathcal M}_\mu f(x_j)\leq\varepsilon'\lambda$, \eqref{eq:goodlambda_far_comparison_sec8} implies 
	\[
	T_{\mu,\sharp,h}(f-f_j)(x)\leq(1+C\varepsilon')\lambda.
	\]
	We choose $\varepsilon'$ such that $C\varepsilon'\leq\varepsilon/2$. Then, for $x\in E_\lambda\cap\mathcal G_Q\cap V_j$, we have 
	\[
	T_\sharp(f_j\mu)(x)\geq T_{\mu,\sharp,h}f(x)-T_{\mu,\sharp,h}(f-f_j)(x)>\frac{\varepsilon}{2}\,\lambda.
	\]
	
	Hence, \eqref{eq:big_piece_measure_weak11_sec8} with $\omega=f_j\mu$ and \eqref{eq:goodlambda_local_mass_sec8} give
	\[
	\mu\bigl(E_\lambda\cap\mathcal G_Q\cap V_j\bigr)\leq\frac{2C}{\varepsilon\lambda}\,\|f_j\|_{L^1(\mu)}\leq2CK_0\,\frac{\varepsilon'}{\varepsilon}\,\mu(Q).
	\]
	Summing over $j$ and decreasing $\varepsilon'$ if necessary, we get $\mu(E_\lambda\cap\mathcal G_Q)\leq\theta_0\mu(Q)/2$, so that, since $\mu(\mathcal G_Q)\geq\theta_0\mu(Q)$, we have $\mu(E_\lambda\cap Q)\leq(1-\theta_0/2)\mu(Q)$. Since $E_\lambda\cap\supp(\mu)\subset\Omega_\lambda$, \eqref{eq:big_piece_db_whitney_DM_aux} gives
	\[
	\mu(E_\lambda)
	\leq\mu(\Omega_\lambda)
	-\frac{\theta_0}{2}\sum_{Q\in\mathcal W^{(s)}_{db}}\mu(Q)
	\leq\Bigl(1-\frac{\theta_0}{4}\Bigr)\mu(\Omega_\lambda).
	\]
	
	We now assume that $\Omega_\lambda=\supp(\mu)$ and, since otherwise \eqref{eq:reduction_1_big_lambda} is trivial, that $E_\lambda\neq\varnothing$. By Lemma \ref{lemma:level_sets_good_lambda}, $\supp(\mu)$ is bounded and $\mu(X)<\infty$. Let $s>\diam(\supp(\mu))$. Since the balls $5B^{(s)}(Q)$, $Q\in\caD^{(s)}_{\mu,0}$, are pairwise disjoint, centered in $\supp(\mu)$, and have radius larger than $\diam(\supp(\mu))$, there is a single top cell $Q=\supp(\mu)$. Moreover, $\supp(\mu)\subset B^{(s)}(Q)$, so $Q$ is doubling. Let $x\in E_\lambda$. Since $\widetilde{\mathcal M}_\mu f(x)\leq\varepsilon'\lambda$, a ball containing $x$ and $\supp(\mu)$ gives $\|f\|_{L^1(\mu)}\leq\varepsilon'\lambda\,\mu(X)$. Hence, \eqref{eq:big_piece_measure_weak11_sec8} with $\omega=f\mu$ gives $\mu(E_\lambda\cap\mathcal G_Q)\leq C\varepsilon'\mu(X)$, and, decreasing $\varepsilon'$ if necessary, we obtain
	\[
	\mu(E_\lambda)\leq\mu(X\setminus\mathcal G_Q)+\frac{\theta_0}{2}\,\mu(X)\leq\Bigl(1-\frac{\theta_0}{2}\Bigr)\,\mu(\Omega_\lambda),
	\]
	which proves \eqref{eq:reduction_1_big_lambda} also in this case.
	
	Finally, we choose $\varepsilon\in(0,1)$ such that $(1+\varepsilon)^2(1-\theta_0/4)<1$, and $\varepsilon'$ as above. Since 
	\[
	\bigl\{T_{\mu,\sharp,h}f>(1+\varepsilon)\lambda\bigr\}\subset E_\lambda\cup\bigl\{\widetilde{\mathcal M}_\mu f>\varepsilon'\lambda\bigr\},
	\]
	integrating \eqref{eq:reduction_1_big_lambda} with respect to $2\lambda\,d\lambda$ gives
	\[
	\|T_{\mu,\sharp,h}f\|_{L^2(\mu)}^2
	\leq(1+\varepsilon)^2\Bigl(1-\frac{\theta_0}4\Bigr)\,\|T_{\mu,\sharp,h}f\|_{L^2(\mu)}^2
	+\frac{(1+\varepsilon)^2}{(\varepsilon')^2}\,\|\widetilde{\mathcal M}_\mu f\|_{L^2(\mu)}^2.
	\]
	
	Since $T_{\mu,\sharp,h}f\in L^2(\mu)$ and $\widetilde{\mathcal M}_\mu$ is bounded on $L^2(\mu)$, this gives $\|T_{\mu,\sharp,h}f\|_{L^2(\mu)}\lesssim\|f\|_{L^2(\mu)}$, uniformly in $h$. Since, moreover, $|T_{\mu,h}f|\leq T_{\mu,\sharp,h}f$ for every $h>0$, we obtain \eqref{eq:L2_boundedness_good_lambda}, which completes the proof of Proposition \ref{theorem:MMT_theorem_A1}.
	
	{
		\subsection{Completion of the proof}\label{subsection:completion_main_theorem}
		We now prove Proposition \ref{proposition:bounded_testing_normalized} by verifying the hypotheses of Proposition \ref{theorem:MMT_theorem_A1}, and we obtain the required big pieces by applying the construction of Sections \ref{section:exceptional_set} and \ref{sec:corona}, together with the bounds of Section \ref{section:probabilisti_argument_subsection}, to the testing measures $\mu|_B$.
		
		\begin{proof}[Proof of Proposition \ref{proposition:bounded_testing_normalized}]
			Let $B_0\in\mathscr B$. By Lemmas \ref{lemma:final_boundedness_suppressed_operators} and \ref{lemma:final_global_bound}, applied to the corona decomposition of Section \ref{sec:corona} associated with $B_0$, we have
			\begin{equation}\label{eq:uniform_tree_bound_sec7}
				\sup_{\varepsilon>0}\sup_{R\in\widetilde\Top(B_0)}\|T_{\widetilde\Phi_R,\mu|_{84B(R)},\varepsilon}\|_{L^2(\mu|_{84B(R)})\to L^2(\mu|_{84B(R)})}\lesssim1,
			\end{equation}
			so the assumption \eqref{eq:corollary:49b} of Corollary \ref{corollary:49} holds with $\mathfrak C_{\mathcal T}$ independent of $\varepsilon$, and, since the kernel of $T$ vanishes if $d(x,y)\leq h$, by \eqref{eq:bounded_reduction_regularization}, the bound \eqref{eq:uniform_tree_bound_sec7} also holds for the whole operator $T_{\widetilde\Phi_R,\mu|_{84B(R)}}$.
			
			Let $g\in L^\infty(\mu)$ be supported in $B_0$. By \eqref{eq:annular_size_estimate}, \eqref{eq:claim_bd:tree:SIo}, Lemma \ref{lemma:packing2_mixed}, and the proofs of Lemmas \ref{lemma:estimate_sigma_0_mixed corona}, \ref{lemma:estimate_sigma_tree_prime}, and \ref{lemma:estimate_sigma_tree_non_doubling}, there exists a non-negative function $A_g$, independent of the truncation parameter, such that
			\[
			\int_{B_0}A_g\,d\mu\lesssim\|g\|_{L^\infty(\mu)}\, \mu(B_0)
			\]
			and, after taking the supremum over the truncation parameter in \eqref{eq:splitting_T_in_sigmas_tree}, we get
			\[
			T_{\mu,\sharp}g(x)\lesssim A_g(x)
			+\sum_{R\in\widetilde\Top(B_0)}\chi_R(x)\,
			T_{\widetilde\Phi_R,\sharp}(\chi_{84B(R)}g\,\mu)(x),\qquad\text{for $\mu$-a.e. }x\in B_0.
			\]
			
			For $R\in\widetilde\Top(B_0)$, Lemma \ref{lemma:mixed_cotlar_transfer}, applied with $\mathfrak m_1=\mathfrak m_2=\mathfrak m_3=\mu|_{84B(R)}$, $\Xi\equiv0$, and $\delta_0=h$, the bound \eqref{eq:uniform_tree_bound_sec7} for the whole operator, the Cauchy-Schwarz inequality, and $\mu(84B(R))\leq K_0\,\mu(R)$ give
			\[
			\int_R T_{\widetilde\Phi_R,\sharp}(\chi_{84B(R)}g\,\mu)\,d\mu
			\lesssim\mu(R)^{1/2}\|g\|_{L^2(\mu|_{84B(R)})}
			\lesssim\|g\|_{L^\infty(\mu)}\,\mu(R).
			\]
			Hence, summing over $R\in\widetilde\Top(B_0)$, applying Lemma \ref{lemma:packing2_mixed}, and arguing in the same way for $T^*$, we obtain that there exists $\mathtt C>0$, independent of $B_0$, $g$, $h$, and $L$, such that
			\begin{equation}\label{eq:corollary49_maximal_sec8}
				\int_{B_0}\bigl(T_{\mu,\sharp}g+(T^*)_{\mu,\sharp}g\bigr)\,d\mu
				\leq \mathtt C\,\|g\|_{L^\infty(\mu)}\,\mu(B_0).
			\end{equation}
			
			We may assume that $\mu\neq0$. We set $W\coloneqq\supp(\mu)$ and, for $s>0$, we let $\caD^{(s)}_\mu$ be the lattice of Lemma \ref{lemma:DM_cubes_one_sided} associated with $\mu$, with $K_0$ large enough as required in Proposition \ref{theorem:MMT_theorem_A1}. Let $R\in\caD^{db,(s)}_\mu$ and $B\coloneqq\widetilde B^{(s)}(R)$. By Lemma \ref{lemma:construction_testing_ball_DM} and Remark \ref{remark:choice_of_beta}, the ball $B$ is centered in $W$, has $40$-small boundary, and satisfies $B\cap W\subset R$, $\mu(R)\leq\beta\,\mu(B)$, and $\mu(\Lambda B)\leq\beta\,\mu(\delta B)$.
			\vvv
			
			\underline{\textit{Case I.}}
			We first assume that $r(B)\leq\diam W$, so that, since $\mathfrak C_3\geq40$, the ball $B$ belongs to $\mathscr B$. We apply the construction of Sections \ref{section:exceptional_set} and \ref{sec:corona} with $B_1=B_2=B$ and with the testing measures $\widehat\nu_B=\widehat\nu_B^*=\mu|_B$, so that the total-variation constant equals $1$ and the augmented measures equal $\mu|_B+\mu|_{X\setminus B}=\mu$, and we take the small-set threshold equal to $1/200$.
			
			By \eqref{eq:corollary49_maximal_sec8}, with $B$ in place of $B_0$ and $g=\chi_B$, and by Chebyshev's inequality, the set 
			\[
			U_B\coloneqq\{x\in B:T_\sharp(\mu|_B)(x)>200\,\mathtt C\}
			\]
			and the analogous set $U^*_B$ for $T^*$ have $\mu$-measure at most $\mu(B)/200$, so \ref{ex_item_6} holds for both operators, with $\mathfrak C_2=200\,\mathtt C$. 
			
			Since the augmented measures equal $\mu$, the non-accretive and dyadic-density stopping families are empty, and so are the sets $\mathscr S_1$ and $\mathscr S_2$ and the non-Ahlfors sets, if their thresholds exceed $1$ and the growth constant of $\mu$, respectively. Hence, by Lemma \ref{lemma:measures_exceptional_set_S}, we may choose the thresholds for $T$ and $T^*$ so that $B\cap\mathscr E_T\subset U_B$ and $B\cap\mathscr E_{T^*}\subset U^*_B$.
			
			Let $\theta$ be as in Section \ref{section:exceptional_set}, so that \eqref{eq:small_bd_tilde_B0} gives $\mu(B\setminus(1-\theta)B)\leq\mu(B)/400$. Since the stopping families, the sets $\mathscr S_1$ and $\mathscr S_2$, and the non-Ahlfors sets are empty, the total exceptional sets \eqref{eq:W_new} of the two systems are $\mathscr E_T\cup U_B$ and $\mathscr E_{T^*}\cup U^*_B$, and their intersections with $B$ are $U_B$ and $U^*_B$. Thus, since $B\cap W\subset R$ and $\mu(U_B\cup U^*_B)\leq\mu(B)/100$, inner regularity gives a compact set 
			\[
			F\subset\bigl((1-\theta)B\cap R\bigr)\setminus(U_B\cup U^*_B)\qquad \text{ with }\qquad \mu(F)\geq\mu(B)/2.
			\]
			
			Since $F$ has these properties and the stopping families are empty, we may take $F_1(w_1)=F_2(w_2)=F$ for all $w_1,w_2\in\Omega$, in place of the sets of Section \ref{sec:corona}, so that the random suppressions and their deterministic quantiles equal $\Psi\coloneqq\dist(\cdot,F)$, and $\Xi_{\mathtt w}=\Psi$ for every $\mathtt w\in\Omega^2$.
			
			Moreover, Lemma \ref{eq:lemm_bd_suppressed_2} gives the suppressed testing estimates for $\mu|_B$ and both operators, since $F$ avoids $\mathscr E_T\cup\mathscr E_{T^*}$, and $|\nu_i|(Q)=\mu(Q)$ for every transit cube $Q$ and $i=1,2$, since $|\nu_i|=\mu$. Lemma \ref{lemma:final_boundedness_suppressed_operators} applies with $\mathfrak C_2=200\,\mathtt C$, since its proof only uses the suppressed testing estimates of Lemma \ref{eq:lemm_bd_suppressed_2}, whose constants depend on this fixed $\mathfrak C_2$. Hence, by Lemma \ref{lemma:final_boundedness_suppressed_operators}, with $B_1=B_2=B$ and $\Xi_{\mathtt w}=\Psi$, the operators $T_{\Psi,\mu|_B,\varepsilon}$ are bounded on $L^2(\mu|_B)$ uniformly in $\varepsilon>0$. Since $K_\Psi=K$ on $F\times F$, by Lemma \ref{lemma:kernel_suppressed}, the set $\mathcal G_R\coloneqq F\subset R$ satisfies
			\begin{equation}\label{eq:big_piece_completion_sec7}
				\sup_{\varepsilon>0}\|T_{\mu|_{\mathcal G_R},\varepsilon}\|_{L^2(\mu|_{\mathcal G_R})\to L^2(\mu|_{\mathcal G_R})}\lesssim1
				\qquad\text{and}\qquad
				\mu(\mathcal G_R)\geq\frac{1}{2}\mu(B)\geq(2\beta)^{-1}\mu(R),
			\end{equation}
			with constants independent of $R$ and $s$.
			\vvv
			
			\underline{\textit{Case II.}}
			It remains to consider the case $r(B)>\diam W$, in which $W$ is bounded. Then, $a\coloneqq\diam W$ is positive, because $\mu$ has no atoms, and \eqref{eq:structure_testing_balls} gives 
			\[
			r^{(s)}(R)\geq(2\Lambda/c_*)\,r(B)>a.
			\]
			
			Since the center of $R$ belongs to $W$, we have $W\subset B^{(s)}(R)$, and hence $R=W$, by \eqref{eq:new_radii_modified_DM}. We set $s_0\coloneqq3a/2$. As in the proof of Proposition \ref{theorem:MMT_theorem_A1}, the lattice $\caD^{(s_0)}_\mu$ has the single top cell $W$, which is doubling, and, since $s_0\leq2a$, the testing ball of $W$ belongs to $\mathscr B$, by the observation following Lemma \ref{lemma:construction_testing_ball_DM}. Hence, the construction above, with $s_0$ in place of $s$, gives a set $\mathcal G_W\subset W$ that satisfies \eqref{eq:big_piece_completion_sec7} with $R=W$, and we set $\mathcal G_R\coloneqq\mathcal G_W$.
			
			Therefore, the hypotheses of Proposition \ref{theorem:MMT_theorem_A1} hold for every $s>0$, with $\theta_0=(2\beta)^{-1}$, and Proposition \ref{theorem:MMT_theorem_A1} gives the $L^2(\mu)$-boundedness of $T$, with a bound depending only on the testing and structural parameters. In particular, the bound is independent of the regularization parameters $h$ and $L$.
		\end{proof}
	}

	Since the bound in Proposition \ref{proposition:bounded_testing_normalized} is independent of $h$ and $L$, we can now remove the auxiliary cut-offs introduced in Section \ref{subsection:bounded_testing_reduction} and prove the main result of the paper.
	
	\begin{proof}[Proof of Theorem \ref{theorem:main_theorem_hofmann_problem}]
		The theorem follows from Proposition \ref{proposition:bounded_testing_normalized} and the reduction to the auxiliary kernel in Lemma \ref{lemma:reduction_cutoffs}.
	\end{proof}

\section{The Neumann problem and uniform rectifiability}\label{section:neumann}

In this section, we prove Theorem \ref{theorem:neumann_UR}. To this end, we apply a version of Theorem \ref{theorem:main_theorem_hofmann_problem} with matrix-valued testing functions, see Theorem \ref{theorem:matrix_Tb}, to the Cauchy-Clifford operator, whose components are the Riesz transforms. We use the notation of Section \ref{section:intro_neumann}. In particular, $n\geq2$ and \eqref{eq:neumann_ADR_corkscrew} holds, so that
$\sigma$ has $n$-growth with a constant depending only on $n$ and $C_A$.

Let $V$ be a real Hilbert space of finite dimension, with inner product $\langle\cdot,\cdot\rangle_V$, and let $\mathcal L(V)$ be the space of linear operators on $V$, endowed with the operator norm $\|\cdot\|$. We denote by $M^*$ the adjoint of $M\in\mathcal L(V)$.

An \textit{$\mathcal L(V)$-valued Calder\'on-Zygmund kernel} is a function $K\colon(X\times X)\setminus\{(x,x):x\in X\}\to\mathcal L(V)$ which satisfies the conditions of Definition \ref{def:CZ_kernel} with $\|\cdot\|$ in place of $|\cdot|$. Given such a kernel, a measure $\mu$ as in Theorem \ref{theorem:main_theorem_hofmann_problem}, and $\varepsilon>0$, for $F\in L^1(\mu;V)$ we set
\[
T_{\mu,\varepsilon}F(x)\coloneqq\int_{d(x,y)>\varepsilon}K(x,y)F(y)\,d\mu(y),
\]
and we define $T_{\mu,\sharp}F$ as the supremum over $\varepsilon>0$ of $|T_{\mu,\varepsilon}F|$. The operator $T_{\mu,\varepsilon}$ can be defined analogously if $F\in L^1(\mu;\mathcal L(V))$.

We denote by $T^*$ the SIO associated with the kernel $K(y,x)^*$, and we say that $T_\mu$ is \textit{bounded on} $L^2(\mu;V)$ if its truncations are uniformly bounded on $L^2(\mu;V)$.

\begin{theorem}\label{theorem:matrix_Tb}
Let $(X,d)$, $\mu$, and $\mathscr B$ be as in Theorem \ref{theorem:main_theorem_hofmann_problem}, and let $T$ be associated with an $\mathcal L(V)$-valued Calder\'on-Zygmund kernel with constants $C$ and $\gamma\in(0,1)$. Let $1< p,q\leq \infty$ and $C_b\geq1$. Suppose that, for every $B\in\mathscr B$, there are $\mathfrak b_B\in L^p(\mu|_B;\mathcal L(V))$ and $\mathfrak b^*_B\in L^q(\mu|_B;\mathcal L(V))$, extended by zero outside $B$, such that
\begin{equation}\label{eq:matrix_Tb_accretive}
	\langle\mathfrak b_B(x)v,v\rangle_V\geq|v|^2\quad\text{ and }\quad\langle\mathfrak b^*_B(x)v,v\rangle_V\geq|v|^2,\qquad\text{ for }v\in V\text{ and $\mu$-a.e. }x\in B,
\end{equation}
we have
\begin{equation}\label{eq:matrix_Tb_size}
	\int_B\|\mathfrak b_B\|^p\,d\mu\leq C_b^p\,\mu(B)\qquad\text{ and }\qquad\int_B\|\mathfrak b^*_B\|^q\,d\mu\leq C_b^q\,\mu(B),
\end{equation}
and that, for every $\tau\in(0,1)$, there is $C_\tau\geq1$, independent of $B$, such that
\begin{equation}\label{eq:matrix_Tb_testing}
	\int_{(1-\tau)B}\bigl(T_{\mu,\sharp}\mathfrak b_B\bigr)^p\,d\mu\leq C_\tau\,\mu(B)\qquad\text{ and }\qquad\int_{(1-\tau)B}\bigl(T^*_{\mu,\sharp}\mathfrak b^*_B\bigr)^q\,d\mu\leq C_\tau\,\mu(B).
\end{equation}
Then $T_\mu$ is bounded on $L^2(\mu;V)$. The bound depends only on $\dim V$, $p$, $q$, $C_b$, the constants $D,n,\mathfrak C_0,C,\gamma,\delta,\Lambda,\beta,\mathfrak C_3$ of Theorem \ref{theorem:main_theorem_hofmann_problem}, and $C_\tau$, for a value of $\tau$ depending only on these parameters.
\end{theorem}

We note that condition \eqref{eq:matrix_Tb_accretive} only involves the symmetric part of the testing functions. For instance, it holds if $\mathfrak b_B(x)$ is the identity plus a skew-adjoint operator, which is the case in the proof of Theorem \ref{theorem:neumann_UR}. We state Theorem \ref{theorem:matrix_Tb} in the generality of Theorem \ref{theorem:main_theorem_hofmann_problem}, since the proof is the same, although, in the proof of Theorem \ref{theorem:neumann_UR}, the measure is Ahlfors regular and the kernel is antisymmetric, see Section \ref{section:intro_neumann}, and we prove it in Appendix \ref{appendix:matrix_Tb}.

\subsection{Clifford algebra and the Cauchy-Clifford operator}\label{subsection:clifford}
Before proving the main application to the Neumann problem, we review the elements of Clifford analysis used in the proof. For more details, we refer the reader, for instance, to \cite[Sections 3-4]{HMT10} or \cite{BDS82}.

We denote by $\mathrm{Cl}_{n+1}$ the \textit{real Clifford algebra} generated by $e_1,\dots,e_{n+1}$, with the relations 
\[
e_j^2=-1\qquad \text{ and }\qquad e_je_k=-e_ke_j \quad \text{ for }j\neq k.
\]

Its elements are $b=\sum_Ib_Ie_I$, with $b_I\in\mathbb R$, where $e_\varnothing\coloneqq1$ and $e_I\coloneqq e_{i_1}\cdots e_{i_k}$ for $I=\{i_1<\dots<i_k\}\subset\{1,\dots,n+1\}$. We call $b_\varnothing$ the \textit{scalar part} of $b$, and \textit{bivector} an element of the form $H=\sum_{j<k}H_{jk}e_je_k$. We identify $\mathrm{Cl}_{n+1}$ with the real Hilbert space $V\coloneqq\mathbb R^{2^{n+1}}$, for which $\{e_I\}_I$ is an orthonormal basis, so that $|b|=(\sum_Ib_I^2)^{1/2}$.

We use the following conventions, which are standard in Clifford analysis.
\begin{enumerate}[label=(\roman*)]
	\item \textit{Products.} We denote the Clifford product by juxtaposition, and we write $X\cdot Y$ for the Euclidean inner product of $X,Y\in\mathbb R^{n+1}$.
	\item \textit{Scalars and vectors.} We identify $a\in\mathbb R$ with $a\,e_\varnothing$, and $X=(x_1,\dots,x_{n+1})\in\mathbb R^{n+1}$ with the \textit{vector} $\sum_jx_je_j$. Accordingly, an $\mathbb R^{n+1}$-valued function $F=(F_1,\dots,F_{n+1})$, such as the outer unit normal $\nu$ or the gradient $\nabla u$ of a real function $u$, is identified with the $\mathrm{Cl}_{n+1}$-valued \textit{vector field} $\sum_jF_je_j$. For $X,Y\in\mathbb R^{n+1}$, the relations above give
	\begin{equation}\label{eq:clifford_vectors}
		XY=-X\cdot Y+X\wedge Y,\qquad\text{ where }\qquad X\wedge Y\coloneqq\sum_{j<k}(x_jy_k-x_ky_j)\,e_je_k,
	\end{equation}
	so that, in particular, $XX=-|X|^2$.
	\item \textit{Distributions and measures.}
	A $\mathrm{Cl}_{n+1}$-valued distribution on $\mathbb R^{n+1}$ is a sum $T=\sum_IT_Ie_I$, where each $T_I$ is a real distribution; $\mathrm{Cl}_{n+1}$-valued functions and measures are defined analogously. Derivatives, integrals, and identities in $\mathcal D'(\mathbb R^{n+1})$ are understood coefficientwise, so that, for instance, $\partial_jT=\sum_I(\partial_jT_I)\,e_I$ and $\int g\,d\mu=\sum_I\bigl(\int g_I\,d\mu\bigr)\,e_I$. The Clifford product is not commutative,  so products and convolutions keep the order of the factors. That is, for $a\in\mathrm{Cl}_{n+1}$, a $\mathrm{Cl}_{n+1}$-valued function $S\in L^1_{\loc}(\mathbb R^{n+1})$, and a compactly supported $\mathrm{Cl}_{n+1}$-valued distribution $T$, we set
	\[
		aT\coloneqq\sum_{I,J}a_I\,T_J\,e_Ie_J\qquad\text{ and }\qquad S*T\coloneqq\sum_{I,J}(S_I*T_J)\,e_Ie_J,
	\]
	and $Ta$ is defined analogously. For a $\mathrm{Cl}_{n+1}$-valued function $g$ and a positive Radon measure $\mu$ such that $|g|\in L^1_{\loc}(\mu)$, we denote by $g\,\mu\coloneqq\sum_Ig_I\,\mu$ the $\mathrm{Cl}_{n+1}$-valued measure with density $g$ with respect to $\mu$. We write $\mathcal L^{n+1}$ for the Lebesgue measure on $\mathbb R^{n+1}$, and $dX$ in the integrals with respect to it.
\end{enumerate}

\vv
For $b\in\mathrm{Cl}_{n+1}$, we denote by $L_b\in\mathcal L(V)$ the \textit{left multiplication} by $b$.

\begin{remark}[Properties of left multiplication]
	For $b,c\in \mathrm{Cl}_{n+1}$, we have $L_bL_c=L_{bc}$ and, since $L_b1=b$, it also holds that $|b|\leq\|L_b\|\leq C_n|b|$.
	
	The operator $L_{e_j}$ is orthogonal, and $L_{e_j}^*=L_{e_j}^{-1}=-L_{e_j}$. Hence, 
	\[
	L_{e_je_k}^*=L_{e_k}L_{e_j}=-L_{e_je_k}, \qquad \text{ for }\qquad j\neq k
	\]
	and, for $x\in\mathbb R^{n+1}$ and every bivector $H$, we have
	\begin{equation}\label{eq:clifford_adjoints}
		L_x^*=-L_x,\qquad L_x^*L_x=|x|^2I,\qquad\text{ and }\qquad L_H^*=-L_H.
	\end{equation}
	In particular, $\|L_x\|=|x|$ and, if $b=g+H$, with $g\in\mathbb R$ and $H$ a bivector, then
	\begin{equation}\label{eq:clifford_accretive}
		\langle L_bv,v\rangle_V=g\,|v|^2,\qquad\text{ for }v\in V.
	\end{equation}
	The identities in \eqref{eq:clifford_adjoints} are special cases of the formula $L_b^*=L_{\bar b}$, where $\bar b$ is the Clifford conjugate of $b$, for which $\bar x=-x$ and $\bar H=-H$, see \cite[(4.6.12) and Lemma 4.44]{HMT10}.
\end{remark}
\vv

We denote by $\mathsf D\coloneqq\sum_je_j\partial_j$ the \textit{Dirac operator}, which acts from the left, that is, $\mathsf DF=\sum_je_j\,\partial_jF$ for a $\mathrm{Cl}_{n+1}$-valued distribution $F$. For clarity, we write $\mathsf D_RF\coloneqq\sum_j(\partial_jF)\,e_j$ for its action from the right, which only appears in \eqref{eq:cauchy_reproducing}. We say that $F$ is \textit{monogenic} if $\mathsf DF=0$. Then $\mathsf D^2=\mathsf D_R^2=-\Delta$ and, for a vector field $F=\sum_kF_ke_k$, we have
\begin{equation}\label{eq:dirac_div_curl}
	\mathsf DF=-\div F+\operatorname{curl}F,\qquad\text{ where }\qquad\operatorname{curl}F\coloneqq\sum_{j<k}(\partial_jF_k-\partial_kF_j)\,e_je_k,
\end{equation}
in the sense of distributions.
\begin{remark}
	On a real function $u$, both $\mathsf D$ and $\mathsf D_R$ act as the gradient, that is, $\mathsf Du=\mathsf D_Ru=\nabla u$, by the identification in (ii). This is the only overlap between the two notations: we write $\nabla u$ only for the gradient of a real function $u$, which is a vector field, and $\mathsf D$ for the Dirac operator acting on $\mathrm{Cl}_{n+1}$-valued functions and distributions. For instance, \eqref{eq:dirac_div_curl} gives $\mathsf D\nabla u=-\Delta u$, since $\operatorname{curl}\nabla u=0$, so that $\nabla u$ is monogenic if $u$ is harmonic.
\end{remark}

Let $\omega_n\coloneqq\mathcal H^n(\mathbb S^{n})$ and $\Gamma_\Delta(X)\coloneqq|X|^{1-n}/((n-1)\,\omega_n)$, so that $-\Delta\Gamma_\Delta=\delta_0$, where $\delta_0$ denotes the Dirac mass at the origin, and let
\[
\mathfrak E(X)\coloneqq\nabla\Gamma_\Delta(X)=-\frac{X}{\omega_n|X|^{n+1}}
\]
be the \textit{Cauchy-Clifford kernel}.

Since $\mathfrak E=\mathsf D\Gamma_\Delta=\mathsf D_R\Gamma_\Delta$, we have $\mathsf D\mathfrak E=\mathsf D_R\mathfrak E=-\Delta\Gamma_\Delta=\delta_0$. Hence, for every compactly supported $\mathrm{Cl}_{n+1}$-valued distribution $F$, we have
\begin{equation}\label{eq:cauchy_reproducing}
	\mathfrak E*(\mathsf DF)=\sum_j(\mathfrak Ee_j)*\partial_jF=(\mathsf D_R\mathfrak E)*F=F.
\end{equation}

\begin{definition}
	For $b\in L^1(\sigma;V)$, we define the \textit{Cauchy-Clifford potential}
	\[
	\mathcal C_\sigma b(X)\coloneqq\int_{\partial\Omega}\mathfrak E(X-y)\,b(y)\,d\sigma(y),\qquad\text{ for } X\in\mathbb R^{n+1}\setminus \partial\Omega.
	\]
\end{definition}

 We remark that $\mathcal C_\sigma b$ is smooth on $\mathbb R^{n+1}\setminus \partial\Omega$ and,
since $\mathfrak E$ is locally integrable, Fubini's theorem shows that the same integral defines a function in $L^1_{\loc}(\mathbb R^{n+1})$, which equals $\mathfrak E*(b\sigma)$.

We also denote by $\mathcal C_{\sigma,\varepsilon}$ and $\mathcal C_{\sigma,\sharp}$ the truncations and the maximal operator, with respect to $\sigma$, of the SIO associated with the kernel $K(x,y)\coloneqq L_{\mathfrak E(x-y)}$, as before Theorem \ref{theorem:matrix_Tb}.

By \eqref{eq:clifford_adjoints}, $K$ is an $\mathcal L(V)$-valued Calder\'on-Zygmund kernel, with $\|K(x,y)\|=|x-y|^{-n}/\omega_n$ and regularity exponent $1$, and the oddness of $\mathfrak E$ gives $K(y,x)^*=K(x,y)$. Moreover, $\mathcal C_{\sigma,\varepsilon}(L_b)=L_{\mathcal C_{\sigma,\varepsilon}b}$ for $b\in L^1(\sigma;V)$, where $L_b$ denotes the function $x\mapsto L_{b(x)}$.
\vv

Finally, for $1\leq j\leq n+1$, let 
\[
R_{j,\sigma,\varepsilon}f(x)\coloneqq\int_{|x-y|>\varepsilon}\frac{x_j-y_j}{|x-y|^{n+1}}f(y)\,d\sigma(y),
\]
and let $R_{\sigma,\varepsilon}\coloneqq(R_{1,\sigma,\varepsilon},\dots,R_{n+1,\sigma,\varepsilon})$ be the \textit{truncated $n$-dimensional Riesz transform}.

We remark that, for every real function $f\in L^2(\sigma)$, the Riesz and the Cauchy-Clifford transforms are related by the identity  
$
\mathcal C_{\sigma,\varepsilon}f=-\omega_{n}^{-1}\, \sum_j(R_{j,\sigma,\varepsilon}f)\,e_j,
$
which readily implies
\begin{equation}\label{eq:cauchy_riesz}
	\|R_{\sigma,\varepsilon}f\|_{L^2(\sigma;\mathbb R^{n+1})}=\omega_{n}\,\|\mathcal C_{\sigma,\varepsilon}f\|_{L^2(\sigma;V)}.
\end{equation}

\vvv

\subsection{The testing functions and the proof of Theorem \ref{theorem:neumann_UR}}\label{subsection:neumann_proof}

We start by presenting an overview of the construction of the testing functions. For every ball $B$ centered on $\partial \Omega$, Theorem \ref{theorem:matrix_Tb} requires a testing function for the Cauchy-Clifford operator which is accretive on $B$, has controlled size, and whose maximal truncations are controlled in $L^p(\sigma)$. The only information we have is the Neumann estimate \eqref{eq:neumann_estimate}, so we look for $\mathfrak b_B=\chi_BL_{b_B}$, where $b_B$ is a density whose Cauchy-Clifford potential is the gradient of a solution of a Neumann problem. The smooth case, which we sketch first, serves as a model for the construction of such a density.
\vv

\underline{\textit{The smooth case.}} Suppose that $\Omega\subset \mathbb R^{n+1}$ is smooth, with outer unit normal $\nu$, and let $w$ be smooth up to the boundary, with $\Delta w=f$ compactly supported in $\Omega$. 
For a $\mathrm{Cl}_{n+1}$-valued function $F$ smooth up to the boundary, extended by zero outside $\Omega$, the divergence theorem gives $\partial_j(\chi_\Omega F)=\chi_\Omega\,\partial_jF-\nu_jF\,\sigma$, and hence
\begin{equation}\label{eq:clifford_divergence}
	\mathsf D(\chi_\Omega F)=\chi_\Omega\,\mathsf DF-(\nu F)\,\sigma\qquad\text{ in }\mathcal D'(\mathbb R^{n+1}),
\end{equation}
which is the Clifford version of the divergence theorem, see \cite[(3.4.6)]{HMT10}. For $F=\nabla w$, we have $\mathsf D\nabla w=-\Delta w=-f$ in $\Omega$ and, by \eqref{eq:clifford_vectors}, we obtain the identity of $\mathrm{Cl}_{n+1}$-valued measures
\begin{equation}\label{eq:neumann_model}
	\mathsf D(\chi_\Omega\nabla w)=(-\nu\nabla w)\,\sigma-f\,\mathcal L^{n+1},\qquad\text{ where }\qquad-\nu\nabla w=\partial_\nu w-\nu\wedge\nabla w.
\end{equation}

Thus, the boundary term is the sum of the Neumann datum of $w$ and of a bivector whose coefficients are tangential derivatives of $w$. By \eqref{eq:dirac_div_curl}, the scalar and the bivector parts of \eqref{eq:neumann_model} are, respectively,
\[
-\div(\chi_\Omega\nabla w)=(\partial_\nu w)\,\sigma-f\,\mathcal L^{n+1}\qquad\text{ and }\qquad\operatorname{curl}(\chi_\Omega\nabla w)=-(\nu\wedge\nabla w)\,\sigma,
\]
and, in the general case, they will be given by the weak formulation of the Neumann problem and by Lemma \ref{lemma:neumann_curl}, respectively. Let $v\coloneqq-\Gamma_\Delta*f$, so that $\Delta v=f$ and $\nabla v=-\mathfrak E*f$. Since $\chi_\Omega\nabla w$ is compactly supported, \eqref{eq:cauchy_reproducing} and \eqref{eq:neumann_model} give
\[
\chi_\Omega\nabla w=\mathfrak E*\mathsf D(\chi_\Omega\nabla w)=\mathfrak E*\bigl((-\nu\nabla w)\,\sigma\bigr)-\mathfrak E*f=\mathcal C_\sigma(-\nu\nabla w)+\nabla v\qquad\text{ in }\mathbb R^{n+1}\setminus\partial\Omega.
\]
Since all the terms are continuous in $\Omega$, we obtain
\begin{equation}\label{eq:neumann_model_reproducing}
	\mathcal C_\sigma(-\nu\nabla w)=\nabla(w-v)\quad\text{ in }\Omega,
\end{equation}
and $w-v$ is harmonic in $\Omega$. If $f=0$, this is the Cauchy reproducing formula for the monogenic function $\nabla w$, see \cite[(4.7.13)]{HMT10}.

\vv

\underline{\textit{The general case.}} In our setting, neither the normal nor the traces of $\nabla w$ are available, and we proceed as follows:
\begin{itemize}
	\item We first fix a datum $g_B$ equal to $1$ on $2B$, and a smooth source $f_B$ with the same mass at a corkscrew point.
	\item Next, we obtain the weak solution $w_B$ of $\Delta w_B=f_B$ in $\Omega$, with Neumann datum $g_B$, from the hypotheses of Theorem \ref{theorem:neumann_UR}.
	\item We replace the tangential derivatives by the distributional curl of $\chi_\Omega\nabla w_B$, which is a bivector-valued measure $H_B\,\sigma$ given by Lemma \ref{lemma:neumann_curl}, and we set $b_B\coloneqq g_B+H_B$. Then \eqref{eq:neumann_model} holds with $w_B$ in place of $w$ and $b_B$ in place of $-\nu\nabla w$, so that \eqref{eq:neumann_model_reproducing} holds as well, see \eqref{eq:neumann_cauchy_identity}. Finally,we obtain proper bounds for $b_B$ and use this information to prove that $\mathfrak b_B\coloneqq \chi_B L_{b_B}$ satisfies the desired testing properties.
\end{itemize} 

From now on, the constants denoted by $C$ depend only on $n$, $C_A$, $c_I$, $p$, $a$, and $C_N$, unless otherwise specified. We also fix the aperture $a_0\coloneqq\max\{a,2/c_I\}$, and we use two estimates for nontangential maximal functions.

First, we note that $\|N_{a_0}U\|_{L^p(\sigma)}\leq C\|N_aU\|_{L^p(\sigma)}$, see \cite[Proposition 2.2]{HMT10}. We record the following standard auxiliary lemma.

\begin{lemma}
	If $0<t\leq\diam \partial\Omega$ and $U$ is measurable on $\Omega$, then
	\begin{equation}\label{eq:neumann_layer}
		\int_{\Omega\cap\{\delta_\Omega<t\}}|U|\,dX\leq C\,t\int_{\partial\Omega}N_{a_0}U\,d\sigma.
	\end{equation}
\end{lemma}
\begin{proof}
	Let $X\in\Omega$ be such that $\delta_\Omega(X)<t$, and let $x_X\in \partial\Omega$ be such that $|X-x_X|=\delta_\Omega(X)$. Every $y\in \partial\Omega\cap B(x_X,(a_0-1)\delta_\Omega(X))$ satisfies $|X-y|<a_0\,\delta_\Omega(X)$ and, by \eqref{eq:neumann_ADR_corkscrew} and since $\delta_\Omega(X)<\diam \partial\Omega$, this set has $\sigma$-measure at least $c\,\delta_\Omega(X)^n$, with $c>0$ depending only on $a_0$ and $C_A$. Hence, Tonelli's theorem gives
	\[
	\begin{split}
		\int_{\Omega\cap\{\delta_\Omega<t\}}|U|\,dX&\leq\frac1c\, \int_{\partial\Omega}\int_{\{X\in\Omega:\,|X-y|<a_0\delta_\Omega(X)<a_0t\}}\frac{|U(X)|}{\delta_\Omega(X)^n}\,dX\,d\sigma(y)\\
		&\leq\frac{a_0^n}c\, \int_{\partial\Omega}N_{a_0}U(y)\, \Bigl(\int_{B(y,a_0t)}\frac{dX}{|X-y|^n}\Bigr)\,d\sigma(y),
	\end{split}
	\]
	which proves \eqref{eq:neumann_layer}.
\end{proof}

In the next lemma, which provides the bivector part of the testing functions, we show that the distributional curl of a gradient extended by zero outside $\Omega$ is absolutely continuous with respect to $\sigma$, with density $H_{jk}$ controlled by the nontangential maximal function of the gradient.
\begin{lemma}\label{lemma:neumann_curl}
	Let $w\in W^{1,2}(\Omega)\cap C^\infty(\Omega)$ be such that $N_{a_0}(\nabla w)\in L^p(\sigma)$, and let $F\coloneqq\chi_\Omega\nabla w$. Then there are real functions $H_{jk}$ on $\partial\Omega$, $1\leq j<k\leq n+1$, such that
	\begin{equation}\label{eq:neumann_curl}
		\partial_jF_k-\partial_kF_j=H_{jk}\,\sigma\quad\text{ in }\mathcal D'(\mathbb R^{n+1})\qquad\text{ and }\qquad\Bigl(\sum_{j<k}|H_{jk}|^2\Bigr)^{1/2}\leq C\,N_{a_0}(\nabla w)\quad\sigma\text{-a.e.}
	\end{equation}
\end{lemma}

\begin{proof}
	Note that the field $F$ is integrable and compactly supported, since $\Omega$ is bounded and $\nabla w\in L^2(\Omega)$. For $0<t<\diam \partial\Omega/2$, we define
	\[
	\psi_t(X)\coloneqq\min\bigl\{1,\max\{0,t^{-1}\, \delta_\Omega(X)-1\}\bigr\},
	\]
	so that $\psi_t$ is Lipschitz, $\psi_t=0$ on $\{\delta_\Omega\leq t\}$, $\psi_t=1$ on $\{\delta_\Omega\geq2t\}$, and $|\nabla\psi_t|\leq1/t$.
	
	The field $F_t\coloneqq\psi_tF$ is Lipschitz and supported in a compact subset of $\Omega$, $F_t\to F$ in $L^1(\mathbb R^{n+1})$ as $t\to0$, and, since the mixed derivatives of $w$ commute, we have
	\[
	\partial_j(F_t)_k-\partial_k(F_t)_j=(\partial_j\psi_t)\,\partial_kw-(\partial_k\psi_t)\,\partial_jw.
	\]
	
	Let $\Psi=(\Psi_{jk})_{j<k}$ be a family of smooth compactly supported real-valued functions, and let $|\Psi|$ be its Euclidean norm. Since $F_t\to F$ in $L^1(\mathbb R^{n+1})$, for $j<k$, the pairing $\langle\partial_jF_k-\partial_kF_j,\Psi_{jk}\rangle$ is the limit, as $t\to0$, of the corresponding pairing for $F_t$. Moreover, 
	$$\langle\partial_j(F_t)_{k}-\partial_k(F_t)_{j},\Psi_{jk}\rangle=\int_\Omega\bigl((\partial_j\psi_t)\,\partial_kw-(\partial_k\psi_t)\,\partial_jw\bigr)\Psi_{jk}\,dX,$$
	since $F_t$ is Lipschitz and compactly supported in $\Omega$. Thus, $\nabla\psi_t$ vanishes almost everywhere outside $\{t<\delta_\Omega<2t\}$ and, by the Cauchy-Schwarz inequality and the identity $\sum_{j<k}(a_jb_k-a_kb_j)^2=|a|^2|b|^2-(a\cdot b)^2$, applied with $a=\nabla\psi_t$ and $b=\nabla w$, we obtain
	\[
	\begin{split}
		\Bigl|\sum_{j<k}\langle\partial_jF_k-\partial_kF_j,\Psi_{jk}\rangle\Bigr|&\leq\limsup_{t\to0}\frac1t\int_{\Omega\cap\{t<\delta_\Omega<2t\}}|\nabla w|\,|\Psi|\,dX\\
		&\leq C\limsup_{t\to0}\int_{\partial\Omega}N_{a_0}(\nabla w)(x)\sup_{|Y-x|<2a_0t}|\Psi(Y)|\,d\sigma(x),
	\end{split}
	\]
	where, in the last inequality, we applied \eqref{eq:neumann_layer}, with $2t$ in place of $t$, to $\chi_{\{t<\delta_\Omega<2t\}}|\nabla w||\Psi|$, and we used that $|Y-x|<2a_0t$ if $|Y-x|<a_0\,\delta_\Omega(Y)$ and $\delta_\Omega(Y)<2t$. Since $\Psi$ is uniformly continuous and $N_{a_0}(\nabla w)\in L^1(\sigma)$, dominated convergence yields
	\[
	\Bigl|\sum_{j<k}\langle\partial_jF_k-\partial_kF_j,\Psi_{jk}\rangle\Bigr|\leq C\int_{\partial\Omega}|\Psi|\,d\mathfrak m,\qquad\text{ where }\quad d\mathfrak m\coloneqq N_{a_0}(\nabla w)\,d\sigma.
	\]
	
	Hence, the sum on the left-hand side defines a linear functional of $\Psi|_{\partial\Omega}$ which, since the restrictions to $\partial\Omega$ of smooth compactly supported functions are dense in $L^1(\mathfrak m)$, extends to a bounded functional on $L^1(\mathfrak m;\mathbb R^{n(n+1)/2})$. By duality, there is $\kappa=(\kappa_{jk})_{j<k}$ such that $|\kappa|\leq C$ $\mathfrak m$-a.e. and the sum equals $\sum_{j<k}\int_{\partial\Omega}\kappa_{jk}\Psi_{jk}\,d\mathfrak m$. The functions $H_{jk}\coloneqq\kappa_{jk}N_{a_0}(\nabla w)$ satisfy \eqref{eq:neumann_curl}.
\end{proof}

In the smooth case, we recall that $H_{jk}=\nu_k\partial_jw-\nu_j\partial_kw$.
\vv

We now construct the Neumann datum and the source.
Let $B\coloneqq B(x_0,r)$, with $x_0\in \partial\Omega$ and $0<r\leq\diam \partial\Omega$. We fix $g_B\in C^\infty_c(3B)$ such that $0\leq g_B\leq1$ and $g_B=1$ on $2B$, and we set $m_B\coloneqq\int_{\partial\Omega}g_B\,d\sigma$. By \eqref{eq:neumann_ADR_corkscrew}, we have
\begin{equation}\label{eq:neumann_gB}
	\sigma(B)\leq m_B\leq C\sigma(B)\leq Cr^n\qquad\text{ and }\qquad\|g_B\|_{L^p(\sigma)}\leq C\sigma(B)^{1/p}.
\end{equation}

Since $m_B>0$, the datum $g_B$ is not admissible in \eqref{eq:neumann_hypothesis}, and we compensate its mass with a smooth source term inside $\Omega$, exploiting the geometric assumptions on the domain. Indeed, by the interior corkscrew condition \eqref{eq:neumann_ADR_corkscrew}, there is $X_B\in\Omega$ such that $B(X_B,4\varrho_B)\subset\Omega\cap B(x_0,r/2)$, where $\varrho_B\coloneqq c_Ir/8$. We fix a non-negative function $\phi_0\in C^\infty_c(B(0,1))$ with $\int\phi_0\,dX=1$, and we set
\begin{equation}\label{eq:neumann_source}
	f_B(X)\coloneqq m_B\,\varrho_B^{-n-1}\, \phi_0\Bigl(\frac{X-X_B}{\varrho_B}\Bigr)\qquad\text{ and }\qquad v_B\coloneqq-\Gamma_\Delta*f_B.
\end{equation}

We note that $f_B$ is supported in $B(X_B,\varrho_B)\subset\Omega$, $\int f_B\,dX=m_B$, $v_B\in C^\infty(\mathbb R^{n+1})$, $\Delta v_B=f_B$, and $\nabla v_B=-\mathfrak E*f_B$. Moreover, $f_B$ approximates $m_B$ times the Dirac mass at $X_B$, at the scale $\varrho_B\approx r$, so $\nabla v_B$ behaves like $-m_B\,\mathfrak E(\cdot-X_B)$, and we have the following estimates.

\begin{lemma}\label{lemma:neumann_vB}
	For every $X\in\mathbb R^{n+1}$ and $x\in \partial\Omega$, we have
	\begin{equation}\label{eq:neumann_vB_decay}
		|\nabla v_B(X)|\leq\frac{C\,m_B}{(r+|X-X_B|)^n}\qquad\text{ and }\qquad N_{a_0}(\nabla v_B)(x)\leq\frac{C\,m_B}{(r+|x-X_B|)^n}.
	\end{equation}
	In particular, $\|\nabla v_B\|_{L^p(\sigma)}+\|N_{a_0}(\nabla v_B)\|_{L^p(\sigma)}\leq C\sigma(B)^{1/p}$.
\end{lemma}

\begin{proof}
	If $|X-X_B|\leq2\varrho_B$, the local integrability of $\mathfrak E$ gives
	\[
	|\nabla v_B(X)|\leq C\,m_B\,\varrho_B^{-n-1}\int_{B(X,3\varrho_B)}\frac{dY}{|X-Y|^n}\leq C\,m_B\,\varrho_B^{-n}.
	\]
	
	On the other hand, in case $|X-X_B|>2\varrho_B$, we have $|X-Y|\geq|X-X_B|/2$ for $Y\in\supp f_B$, so that $|\nabla v_B(X)|\leq Cm_B|X-X_B|^{-n}$. Since $\varrho_B\approx r$, this proves the first estimate in \eqref{eq:neumann_vB_decay}. 
	
	We now prove the nontangential maximal function bound for $\nabla v_B$. If $x\in \partial\Omega$ and $Y\in\Omega$ are such that $|Y-x|<a_0\,\delta_\Omega(Y)$, then
	\(
	\delta_\Omega(Y)\leq|Y-X_B|+\delta_\Omega(X_B)\leq|Y-X_B|+\tfrac{r}{2},
	\)
	thus
	\[
	r+|x-X_B|\leq r+a_0\,\delta_\Omega(Y)+|Y-X_B|\leq(a_0+1)(r+|Y-X_B|),
	\]
	and the second estimate in \eqref{eq:neumann_vB_decay} follows from the first one. Finally, since $\sigma$ has $n$-growth, splitting $\partial\Omega$ into $\partial\Omega\cap B(X_B,r)$ and the annuli $\{x\in \partial\Omega:2^kr\leq|x-X_B|<2^{k+1}r\}$, $k\geq0$, we obtain
	\[
	\int_{\partial\Omega}\frac{m_B^p}{(r+|x-X_B|)^{np}}\,d\sigma(x)\leq C\,m_B^p\,r^{-n(p-1)}\Bigl(1+\sum_{k\geq0}2^{-kn(p-1)}\Bigr)\leq C\, \sigma(B),
	\]
	where we used $p>1$ and \eqref{eq:neumann_gB}. By \eqref{eq:neumann_vB_decay}, this concludes the proof.
\end{proof}

Next, we construct and study the Neumann solution with right-hand side $f_B$.
Before doing so, we note that since $\mathcal H^n(\partial\Omega)<\infty$, the domain $\Omega$ has finite perimeter, its reduced boundary $\partial^*\Omega$ is contained in $\partial\Omega$ and, denoting by $\nu_\Omega$ the measure-theoretic outer unit normal, the Gauss-Green formula
\[
\int_\Omega\div W\,dX=\int_{\partial^*\Omega}W\cdot\nu_\Omega\,d\mathcal H^n
\]
holds for every $W\in C^1_c(\mathbb R^{n+1};\mathbb R^{n+1})$, see for instance \cite[Section 2.2]{HMT10}. 

We set $q_B\coloneqq\nu_\Omega\cdot\nabla v_B$ on $\partial^*\Omega$ and $q_B\coloneqq0$ on $\partial\Omega\setminus\partial^*\Omega$. This is a bounded Borel function and, by the Gauss-Green formula applied to $W=\varphi\nabla v_B$, with $\varphi\in C^\infty_c(\mathbb R^{n+1})$, we have
\begin{equation}\label{eq:neumann_qB}
	\int_{\partial\Omega}q_B\,\varphi\,d\sigma=\int_\Omega\nabla v_B\cdot\nabla\varphi\,dX+\int_\Omega f_B\,\varphi\,dX.
\end{equation}
Choosing $\varphi=1$ on a neighborhood of $\overline\Omega$, we get $\int_{\partial\Omega}q_B\,d\sigma=m_B$, and Lemma \ref{lemma:neumann_vB} gives $\|q_B\|_{L^p(\sigma)}\leq C\sigma(B)^{1/p}$.

Therefore, $h_B\coloneqq g_B-q_B$ is a bounded real Borel function such that $\int_{\partial\Omega}h_B\,d\sigma=0$ and, by \eqref{eq:neumann_gB}, $\|h_B\|_{L^p(\sigma)}\leq C\sigma(B)^{1/p}$.

Let $u_B\coloneqq u_{h_B}$ be as in Theorem \ref{theorem:neumann_UR}, and let $w_B\coloneqq u_B+v_B$. By \eqref{eq:neumann_estimate}, the aperture comparison, and Lemma \ref{lemma:neumann_vB}, we then have
\begin{equation}\label{eq:neumann_N_wB}
	\|N_{a_0}(\nabla u_B)\|_{L^p(\sigma)}+\|N_{a_0}(\nabla v_B)\|_{L^p(\sigma)}\leq C\sigma(B)^{1/p},
\end{equation}
and, adding \eqref{eq:neumann_qB} to \eqref{eq:neumann_hypothesis} with $h=h_B$, we obtain
\begin{equation}\label{eq:neumann_weak_wB}
	\int_\Omega\nabla w_B\cdot\nabla\varphi\,dX=\int_{\partial\Omega}g_B\,\varphi\,d\sigma-\int_\Omega f_B\,\varphi\,dX,\qquad\text{ for }\varphi\in C^\infty_c(\mathbb R^{n+1}).
\end{equation}
Thus, $w_B$ is a weak solution of $\Delta w_B=f_B$ in $\Omega$ with boundary Neumann datum $g_B$, and $u_B=w_B-v_B$ plays the role of the harmonic function $w-v$ in the smooth case.
We remark that the normal $\nu_\Omega$ is only paired with the smooth field $\nabla v_B$, and that we do not need the extra assumption $\mathcal H^n(\partial\Omega\setminus\partial^*\Omega)=0$.

The function $w_B$ belongs to $W^{1,2}(\Omega)\cap C^\infty(\Omega)$, since $u_B\in W^{1,2}(\Omega)$ is harmonic and $v_B$ is smooth. Let $H_{jk,B}$ be the functions given by Lemma \ref{lemma:neumann_curl}, applied to $w_B$, and set
\[
H_B\coloneqq\sum_{j<k}H_{jk,B}\,e_je_k,\qquad b_B\coloneqq g_B+H_B,\qquad\text{ and }\qquad F_B\coloneqq\chi_\Omega\nabla w_B.
\]

By \eqref{eq:neumann_gB}, Lemma \ref{lemma:neumann_curl}, and \eqref{eq:neumann_N_wB}, we have
\begin{equation}\label{eq:neumann_bB_size}
	\|b_B\|_{L^p(\sigma;V)}\leq C\sigma(B)^{1/p}.
\end{equation}
Moreover, \eqref{eq:neumann_weak_wB} says that $-\div F_B=g_B\,\sigma-f_B\,\mathcal L^{n+1}$ in $\mathcal D'(\mathbb R^{n+1})$, so that \eqref{eq:dirac_div_curl} and \eqref{eq:neumann_curl} give
\[
\mathsf DF_B=b_B\,\sigma-f_B\,\mathcal L^{n+1}.
\]

Since $F_B$ is integrable and compactly supported, \eqref{eq:cauchy_reproducing} yields
\[
F_B=\mathfrak E*(b_B\sigma)-\mathfrak E*f_B=\mathcal C_\sigma b_B+\nabla v_B.
\]
In $\Omega$, we have $F_B=\nabla u_B+\nabla v_B$, and all the terms are continuous, so that
\begin{equation}\label{eq:neumann_cauchy_identity}
	\mathcal C_\sigma b_B(Y)=\nabla u_B(Y),\qquad\text{ for }Y\in\Omega.
\end{equation}
We remark that if $\Omega$ and $w_B$ are smooth up to the boundary, then $b_B=-\nu\nabla w_B$, and \eqref{eq:neumann_cauchy_identity} is \eqref{eq:neumann_model_reproducing} with $w=w_B$ and $f=f_B$.

Hence, $b_B$ has scalar part $1$ on $2B$ and satisfies the size bound \eqref{eq:neumann_bB_size}, and it remains to control its maximal truncations.

Since, by \eqref{eq:neumann_cauchy_identity}, the Cauchy-Clifford potential of $b_B$ is the gradient of a Neumann solution, whose nontangential maximal function is controlled by \eqref{eq:neumann_estimate}, it suffices to compare the truncation of $\mathcal C_\sigma b_B$ at $x\in \partial\Omega$ and at scale $\varepsilon$ with the value of $\mathcal C_\sigma b_B$ at a corkscrew point at distance comparable to $\varepsilon$ from $x$.

\begin{remark}
	The same principle is used for harmonic measure in \cite{AHMMMTV16, MT20}, and for caloric measure in Lemma \ref{lemma:caloric_doubling_truncation}, where the truncated singular integrals of harmonic and caloric measure are expressed, up to the kernel at the pole, by the gradient of a localized Green function, see \eqref{eq:bound_T_caloric0}. In our setting, the source $f_B$ plays the role of the pole.
\end{remark}

The next lemma, which is a variant of Cotlar's inequality, controls the maximal truncations of the Cauchy-Clifford operator by the nontangential maximal function of the Cauchy-Clifford potential.
\begin{lemma}\label{lemma:neumann_cotlar}
	For every $b\in L^p(\sigma;V)$, we have
	\[
	\mathcal C_{\sigma,\sharp}b(x)\leq N_{a_0}(\mathcal C_\sigma b)(x)+C\,\widetilde{\mathcal M}_\sigma(|b|)(x),\qquad\text{ for }x\in \partial\Omega.
	\]
\end{lemma}

\begin{proof}
	Let $x\in \partial\Omega$ and $\varepsilon>0$. Since $|x-z|\leq\diam \partial\Omega$ for $z\in \partial\Omega$, we may assume that $\varepsilon<\diam \partial\Omega$. By \eqref{eq:neumann_ADR_corkscrew}, there is $Y\in\Omega$ such that $|Y-x|<\varepsilon/4$ and $\delta_\Omega(Y)\geq c_I\varepsilon/4$, so that $|Y-x|<a_0\,\delta_\Omega(Y)$. We split
	\[
	\begin{split}
		\mathcal C_{\sigma,\varepsilon}b(x)-\mathcal C_\sigma b(Y)&=\int_{|x-z|>\varepsilon}\bigl(\mathfrak E(x-z)-\mathfrak E(Y-z)\bigr)\,b(z)\,d\sigma(z)\\
		&\qquad\qquad-\int_{|x-z|\leq\varepsilon}\mathfrak E(Y-z)\,b(z)\,d\sigma(z)\eqqcolon\mathcal I'-\mathcal I''.
	\end{split}
	\]
	If $|x-z|>\varepsilon$, then $|x-Y|<|x-z|/4$, and the smoothness of $\mathfrak E$ gives $|\mathfrak E(x-z)-\mathfrak E(Y-z)|\leq C\varepsilon|x-z|^{-n-1}$. Hence, by a standard decomposition of the domain of integration into dyadic annuli and using the $n$-growth of $\sigma$, we obtain
	\[
	|\mathcal I'|\leq C\sum_{k\geq0}2^{-k}\,\widetilde{\mathcal M}_\sigma(|b|)(x)\leq C\,\widetilde{\mathcal M}_\sigma(|b|)(x).
	\]
	On the other hand, $|Y-z|\geq\delta_\Omega(Y)\geq c_I\varepsilon/4$ for $z\in \partial\Omega$, so that the size estimates for $\mathfrak E$ yield
	\[
	|\mathcal I''|\leq C\varepsilon^{-n}\int_{B(x,2\varepsilon)}|b|\,d\sigma\leq C\,\widetilde{\mathcal M}_\sigma(|b|)(x).
	\]
	Since $|\mathcal C_\sigma b(Y)|\leq N_{a_0}(\mathcal C_\sigma b)(x)$, the lemma follows by taking the supremum over $\varepsilon$.
\end{proof}

We finally prove Theorem \ref{theorem:neumann_UR} by combining Theorem \ref{theorem:matrix_Tb} with the solution of the codimension one David-Semmes problem.
\begin{proof}[Proof of Theorem \ref{theorem:neumann_UR}]
	Let $B=B(x_0,r)$, with $x_0\in \partial\Omega$ and $0<r\leq\diam \partial\Omega$, and let $b_B$ be as above.
	By Lemma \ref{lemma:neumann_cotlar}, \eqref{eq:neumann_cauchy_identity}, \eqref{eq:neumann_N_wB}, \eqref{eq:neumann_bB_size}, and the $L^p(\sigma)$-boundedness of $\widetilde{\mathcal M}_\sigma$, we have
	\begin{equation}\label{eq:neumann_maximal_bB}
		\|\mathcal C_{\sigma,\sharp}b_B\|_{L^p(\sigma)}\leq\|N_{a_0}(\nabla u_B)\|_{L^p(\sigma)}+C\|b_B\|_{L^p(\sigma;V)}\leq C\sigma(B)^{1/p}.
	\end{equation}
	We set $\mathfrak b_B\coloneqq\chi_BL_{b_B}$. Since $g_B=1$ on $B$, \eqref{eq:clifford_accretive} gives 
	$\langle\mathfrak b_B(x)v,v\rangle_V=|v|^2$ for $x\in B$ and $v\in V$
	and, since $\|L_b\|\leq C_n|b|$, \eqref{eq:neumann_bB_size} gives $\int_B\|\mathfrak b_B\|^p\,d\sigma\leq C\sigma(B)$.
	
	Let $\tau\in(0,1)$ and $x\in(1-\tau)B$. Since $|x-z|\geq\tau r$ for $z\in \partial\Omega\setminus B$, H\"older's inequality, \eqref{eq:neumann_bB_size}, and the $n$-growth of $\sigma$ give
	\[
	\begin{split}
		\int_{\partial\Omega\setminus B}|\mathfrak E(x-z)|\,|b_B(z)|\,d\sigma(z)&\leq C\|b_B\|_{L^p(\sigma;V)}\Bigl(\int_{\{z\in \partial\Omega:\,|x-z|\geq\tau r\}}\frac{d\sigma(z)}{|x-z|^{np'}}\Bigr)^{1/p'}\\
		&\leq C\sigma(B)^{1/p}(\tau r)^{-n/p}\leq C\tau^{-n/p}.
	\end{split}
	\]
	Thus, since $\mathcal C_{\sigma,\varepsilon}\mathfrak b_B=L_{\mathcal C_{\sigma,\varepsilon}(\chi_Bb_B)}$, we have 
	\[
	\mathcal C_{\sigma,\sharp}\mathfrak b_B(x)\leq C_{n}\,\mathcal C_{\sigma,\sharp}b_B(x)+C\tau^{-n/p} \qquad \text{ for }x\in(1-\tau)B,
	\]
	and \eqref{eq:neumann_maximal_bB} yields
	\[
	\int_{(1-\tau)B}\bigl(\mathcal C_{\sigma,\sharp}\mathfrak b_B\bigr)^p\,d\sigma\leq C\tau^{-n}\,\sigma(B).
	\]
	Note that $K(y,x)^*=K(x,y)$, so the operator $\mathcal C^*_\sigma$ coincides with $\mathcal C_\sigma$, and we set $\mathfrak b^*_B\coloneqq\mathfrak b_B$ and $q\coloneqq p$.
	
	The testing functions $\mathfrak b_B$ are defined on every ball centered on $\partial\Omega$ with radius at most $\diam \partial\Omega$, hence the hypotheses of Theorem \ref{theorem:matrix_Tb} hold for $X=\mathbb R^{n+1}$, $\mu=\sigma$, the kernel $K$ of Section \ref{subsection:clifford}, with $\gamma=1/2$, and, for instance, $\delta=1$, $\Lambda=100$, $\beta=\beta_0$, and $\mathfrak C_3=220$. Thus, Theorem \ref{theorem:matrix_Tb} and \eqref{eq:cauchy_riesz} show that the truncated Riesz transforms $R_{\sigma,\varepsilon}$ are uniformly bounded on $L^2(\sigma)$, with a bound depending only on $n$, $C_A$, $c_I$, $p$, $a$, and $C_N$. Since $\sigma$ is Ahlfors regular, the result of Nazarov, Tolsa, and Volberg \cite{NToV14} shows that $\partial\Omega$ is uniformly rectifiable, with constants depending only on $n$, $C_A$, and this bound.
\end{proof}

\vv

	\section{Applications to caloric measure}\label{section:application_caloric_measure}

	{
		
		In this section we prove Theorem \ref{theorem:application_caloric_measure}, and we use the notation of Section \ref{section_applications}. We set
		\begin{equation}\label{eq:def_parabolic_ball_tf}
			B^+(\mathbf X,r)\coloneqq B(\mathbf X,r)\cap\{(Y,s):s>t\},
		\end{equation}
		and
		\begin{equation}\label{eq:def_parabolic_ball_tb}
			B^-(\mathbf X,r) \coloneqq B(\mathbf X,r)\cap\{(Y,s):s<t\}
		\end{equation}
		to denote the time-forward and time-backward half-balls, respectively. For a ball $B=B(\mathbf X,r)$, we also write $B^\pm\coloneqq B^\pm(\mathbf X,r)$.
		
		As in the introduction, for a finite measure $\nu$ on $\mathbb P^n$, $T_\varepsilon\nu$ and $T^*_\varepsilon\nu$ denote the truncated operators associated with the kernels $K(\mathbf X,\mathbf Y)= \nabla_X\Gamma(\mathbf X,\mathbf Y)$ and $K(\mathbf Y,\mathbf X)= \nabla_Y\Gamma(\mathbf Y,\mathbf X)$, respectively.
		
		For $f\in L^1(\sigma)$, we write $T_{\sigma,\varepsilon}f\coloneqq T_\varepsilon(f\sigma)$ and $T^*_{\sigma,\varepsilon}f\coloneqq T^*_\varepsilon(f\sigma)$, and $T_{\sigma,\sharp}f$ and $T^*_{\sigma,\sharp}f$ for the associated maximal singular integrals. Analogously, $T_\sharp\nu$ and $T^*_\sharp\nu$ denote the maximal truncations of a finite measure $\nu$.
	}
	
	\vv
	
	Throughout this section, $\Omega\subset\mathbb P^n$ is open and connected, $n\geq1$, and $\sigma$ is a non-negative Radon measure supported on $\partial\Omega$, with 
	\[
	\sigma(B(\mathbf X,r))\leq\mathfrak C_0r^{n+1}\qquad \text{ for every }\mathbf X\in \mathbb P^n, r>0.
	\]
	
	For a pole $\mathbf P\in \Omega$, caloric and adjoint caloric measures are defined by the Perron-Wiener-Brelot construction on their respective essential boundaries, with the point at infinity included when $\Omega$ is unbounded. We denote by $\omega^{\mathbf P}$ and $\omega_*^{\mathbf P}$ their restrictions to the finite essential boundaries $\partial_e\Omega, \partial^*_e\Omega\subset \partial \Omega$, extended by zero to the rest of $\partial\Omega$, so that $\omega^{\mathbf P}(\partial\Omega)\leq1$ and $\omega_*^{\mathbf P}(\partial\Omega)\leq1$. For more details on the construction, and precise definitions of those boundaries, we refer to \cite[Sections 2 and 3]{MP21} and \cite[Chapter 8]{Wa12}.  
	
	In Section \ref{subsection:caloric_proof}, we associate with each ball of the local $Tb$ theorem a smaller ball whose time-forward and time-backward halves have comparable measure, see Lemma \ref{lemma:caloric_auxiliary_ball}. The construction relies on the following lemma.
	
	{
		
		\begin{lemma}\label{lemma:halving_parabolic_balls}
			Let $\sigma$ be a non-negative Radon measure on $\mathbb P^n$ with $(n+1)$-growth. Let $E\subset \mathbb P^n$ be a compact set such that $0<\sigma(E)<\infty$. Then, there exists $\mathbf X_E=(X_E, t_E)\in E$ such that
			\begin{equation}\label{eq:halving_parabolic_balls}
				\sigma\bigl(\bigl\{(Y,s)\in E: s>t_E\bigr\}\bigr)=\sigma\bigl(\bigl\{(Y,s)\in E: s<t_E\bigr\}\bigr)=\frac{1}{2}\sigma(E).
			\end{equation}
		\end{lemma}
		\begin{proof}
			First, we observe that, for every $t_0\in \mathbb R$, the set $E\cap \{(Y,s):s=t_0\}$ has parabolic Hausdorff dimension at most $n$. By the $(n+1)$-growth assumption on $\sigma$, a standard density argument yields
			\[
			\sigma\bigl(E\cap \{(Y,s):s=t_0\}\bigr)=0, \qquad \text{ for every }t_0\in \mathbb R.
			\]
			Hence, if we denote by $\pi_t(X,t)\coloneqq t$ the projection onto the time variable and we define $\nu\coloneqq (\pi_t)_\sharp (\sigma|_E)$, we have that $\nu$ is nonatomic.
			
			We define the function $F_\nu(t)\coloneqq \nu((-\infty,t])$, $t\in \mathbb R$, and note that it is a continuous non-decreasing function such that 
			\[
			\lim_{t\to -\infty}F_\nu(t)=0\qquad \text{ and }\qquad \lim_{t\to\infty}F_\nu(t)=\sigma(E).
			\]
			Thus, there exists $t_E\in \mathbb R$ such that $t_E=\min\{t\in\mathbb R: F_\nu(t)\geq \sigma(E)/2\}$. By minimality, we also have that $\nu((t_E-\varepsilon, t_E+ \varepsilon))>0$ for every $\varepsilon>0$, so $t_E\in \supp(\nu)$. Since $E$ is compact, we have $\supp(\nu)\subset \pi_t(E)$, and we finally conclude that there exists a point $(X_E,t_E)\in E$ for which \eqref{eq:halving_parabolic_balls} holds, because $\nu$ is nonatomic.
		\end{proof}
	}

	{
		\subsection{Elements of heat potential theory}
		In this subsection, we collect some auxiliary potential-theoretic results.
		\begin{lemma}\label{lemma:caloric_boundary_potential}
			Let $\Omega$ be as in Theorem \ref{theorem:application_caloric_measure}, let $\mathbf Z\in\mathbb P^n$ and $r>0$, and let $E\subset\partial\Omega\cap B(\mathbf Z,r)$ be compact. The function 
			\[
			v(\mathbf X)\coloneqq \int_E\Gamma(\mathbf X,\mathbf Y)\,d\sigma(\mathbf Y)
			\]
			is continuous on $\mathbb P^n$, caloric in $\Omega$, vanishes at infinity, and $0\leq v\lesssim \mathfrak C_0\, r$.
		\end{lemma}
		\begin{proof}
			Standard kernel estimates and the upper growth of $\sigma$ give, for every $\mathbf X\in\mathbb P^n$ and $a>0$,
			\begin{equation}\label{eq:caloric_potential_near}
				\int_{E\cap B(\mathbf X,a)}\Gamma(\mathbf X,\mathbf Y)\,d\sigma(\mathbf Y)
				\lesssim\sum_{j\geq0}(2^{-j}a)^{-n}\sigma(B(\mathbf X,2^{-j}a))
				\lesssim\mathfrak C_0a.
			\end{equation}
			
			If $\mathbf X\in B(\mathbf Z, 2r)$, then $E\subset B(\mathbf X, 3r)$, so \eqref{eq:caloric_potential_near} with $a=3r$ gives $v(\mathbf X)\lesssim \mathfrak C_0\, r.$
			For $\mathbf X\in \mathbb P^n\setminus B(\mathbf Z,2r)$, it is enough to note that
			\[
			v(\mathbf X)\lesssim r^{-n}\sigma(E)\lesssim\mathfrak C_0r.
			\]
			The estimate \eqref{eq:caloric_potential_near} with $a\to 0$, which is uniform in $\mathbf X$, together with the continuity of the kernel off the diagonal, yields continuity of $v$. Since $E\subset\partial\Omega$, the function $v$ is caloric in $\Omega$, and finiteness of $\sigma(E)$ gives decay at infinity. 
		\end{proof}
		
		Extended by zero at infinity, the function $v$ belongs to both the upper and the lower Perron-Wiener-Brelot classes for its values on the essential boundary, since it is continuous on $\mathbb P^n$, see \cite{Wa12}, so it is their common solution. Since the term at infinity in the representation formula vanishes, as an immediate consequence of Lemma \ref{lemma:caloric_boundary_potential}, we obtain the identity
		\begin{equation}\label{eq:caloric_potential_representation}
			v(\mathbf P)=\int v(\mathbf W)\,d\omega^{\mathbf P}(\mathbf W),\qquad\text{ for }\mathbf P\in\Omega,
		\end{equation}
		and the corresponding statements hold for the adjoint caloric potential $v^*(\mathbf X)\coloneqq \int_E\Gamma(\mathbf Y,\mathbf X)\,d\sigma(\mathbf Y)$ and $\omega_*^{\mathbf P}$ replacing $v(\mathbf X)$ and  $\omega^{\mathbf P}$, respectively.
	}
	
	{
		The next lemma is a parabolic version of Bourgain's lemma, see also for instance \cite{GeH20}. Here $\mathcal H^{n+s}_{p,\infty}$ denotes the $(n+s)$-dimensional Hausdorff content for the parabolic metric, see \cite[Section 2]{MP21}. 
		\begin{lemma}\label{lemma:Bourgain_lemma_parabolic}
			Let $\mathbf X_0=(X_0,t_0)\in\partial\Omega$, $r>0$, $A\geq1$, and $\vartheta\in(0,1)$. Suppose that $\mathbf P^+\in\Omega$ satisfies $|X_{P^+}-X_0|\leq Ar$ and
			\[
			t_0+\vartheta r^2\leq t_{P^+}\leq t_0+\vartheta^{-1}r^2.
			\]
			Then there are constants $M_0=M_0(n,A,\vartheta)>1$ and $c=c(n,A,\vartheta)>0$ such that, for $M\geq M_0$,
			\begin{equation}\label{eq:caloric_bourgain_density}
				\omega^{\mathbf P^+}(B(\mathbf X_0,Mr))\geq	c\,\frac{\sigma(B^-(\mathbf X_0,r))}{\mathfrak C_0r^{n+1}}.
			\end{equation}
			More generally, for $0<s\leq2$, we have
			\begin{equation}\label{eq:caloric_bourgain_density2}
				\omega^{\mathbf P^+}(B(\mathbf X_0,Mr))\gtrsim_{n,s,A,\vartheta}
				\frac{\mathcal H^{n+s}_{p,\infty}\bigl(B^-(\mathbf X_0,r)\cap\partial\Omega \bigr)}{r^{n+s}}.
			\end{equation}
		\end{lemma}
		\begin{proof}
			Take a compact set $E\subset B^-(\mathbf X_0,r)\cap\partial\Omega$ and let $v$ be its potential as in Lemma \ref{lemma:caloric_boundary_potential}. The time and spatial separation give $v(\mathbf P^+)\geq c_1r^{-n}\sigma(E)$. On the other hand, if $M\geq M_0$,
			\[
			\sup_{\mathbf W\notin B(\mathbf X_0,Mr)}v(\mathbf W)\leq C(M-1)^{-n}r^{-n}\sigma(E)\leq\frac12\, c_1r^{-n}\sigma(E).
			\]
			
			Since $\omega^{\mathbf P^+}(\partial\Omega)\leq1$, \eqref{eq:caloric_potential_representation} and Lemma \ref{lemma:caloric_boundary_potential} imply
			\[
			c_1r^{-n}\sigma(E)\leq C\mathfrak C_0r\,\omega^{\mathbf P^+}(B(\mathbf X_0,Mr))	+\frac12\, c_1r^{-n}\sigma(E),
			\]
			which readily yields \eqref{eq:caloric_bourgain_density} by the inner regularity of $\sigma$. Finally, the bound \eqref{eq:caloric_bourgain_density2} follows from Frostman's lemma.
		\end{proof}
		
		A time reflection of Lemma \ref{lemma:Bourgain_lemma_parabolic} implies that, if $\mathbf P^-\in\Omega$, $|X_{P^-}-X_0|\leq Ar$, and $t_0-\vartheta^{-1}r^2\leq t_{P^-}\leq t_0-\vartheta r^2$, then
		\[
		\omega_*^{\mathbf P^-}(B(\mathbf X_0,Mr))\geq
		c\,\frac{\sigma(B^+(\mathbf X_0,r))}{\mathfrak C_0r^{n+1}},\qquad \text{ for }M\geq M_0.
		\]
	}
	
	{
		In particular, in our application, we take $\delta\leq M_0(n,1,\vartheta)^{-1}$. Let $B$ be a ball as in Definition \ref{definition:weak_A_infty} such that $\sigma((\delta B)^\pm)\gtrsim \sigma(B)$. Lemma \ref{lemma:Bourgain_lemma_parabolic} and its time reflection with $\delta r(B)$ in place of $r$, $A=1$, and $M=\delta^{-1}\geq M_0$, yield
		\begin{equation}\label{eq:nondegeneracy}
			\omega^{\mathbf P_B^+}(B)\gtrsim\frac{\sigma((\delta B)^-)}{\mathfrak C_0(\delta r(B))^{n+1}}
			\gtrsim\frac{\sigma(B)}{\mathfrak C_0r(B)^{n+1}},\qquad\text{ and }\qquad
			\omega_*^{\mathbf P_B^-}(B)\gtrsim\frac{\sigma(B)}{\mathfrak C_0r(B)^{n+1}},
		\end{equation}
		where the implicit constants depend on $n$, $\vartheta$, $\delta$, and the implicit constants in the assumption on $B$. This applies, in particular, to the balls of Lemma \ref{lemma:caloric_auxiliary_ball}.
	}

	\subsection{A local estimate for the cutoff potential}
	We first recall a consequence of the interior gradient estimate for solutions to the adjoint heat equation and Harnack's inequality.

	{
		\begin{lemma}\label{lemma:Cauchy_estimates_heat}
			Suppose that $u\geq0$ is adjoint caloric in $B(\mathbf Z,10r)$, where $\mathbf Z=(Z,t_Z)$, and put $\mathbf Z^+=(Z,t_Z+4r^2)$. Then
			\[
			|\nabla_Xu(\mathbf Z^+)|\leq\frac Cr\,u(\mathbf Y),\qquad \text{ for every }\mathbf Y\in B(\mathbf Z,r).
			\]
			For a non-negative caloric function, the analogous estimate holds at $\mathbf Z^-=(Z,t_Z-4r^2)$.
		\end{lemma}
		\begin{proof}
			The interior gradient estimate gives $|\nabla_Xu(\mathbf Z^+)|\lesssim r^{-1}\sup_{B(\mathbf Z^+,r)}u$. The time-coordinates of points in $B(\mathbf Z^+,r)$ belong to the interval $(t_Z+3r^2,t_Z+5r^2)$, whereas the times of $\mathbf Y\in B(\mathbf Z,r)$ belong to $(t_Z-r^2,t_Z+r^2)$. The adjoint Harnack inequality, applied inside $B(\mathbf Z,10r)$, bounds the supremum by $Cu(\mathbf Y)$, uniformly in $\mathbf Y$. The caloric estimate follows by time reversal.
		\end{proof}
	}
	
	{
		\begin{lemma}\label{lemma:caloric_cutoff_estimate}
			Let $\mathbf Z\in\supp\sigma$, $r>0$ and $\mathbf P^+\in\Omega\setminus B(\mathbf Z,100r)$. Choose $0\leq\varphi_{\mathbf Z,r}\leq1$, smooth, vanishing on $B(\mathbf Z,10r)$ and equal to $1$ outside $B(\mathbf Z,20r)$, and define
			\begin{equation}\label{eq:definition_u_+}
				u^+(\mathbf Y)\coloneqq \Gamma(\mathbf P^+,\mathbf Y)-\int\Gamma(\mathbf W,\mathbf Y)\varphi_{\mathbf Z,r}(\mathbf W)\,d\omega^{\mathbf P^+}(\mathbf W).
			\end{equation}
			Then $u^+\geq0$ is adjoint caloric in $B(\mathbf Z,10r)$ and, for $\mathbf Z^+\coloneqq (Z,t_Z+4r^2)$, we have
			\begin{equation}\label{eq:bound_T_caloric2}
				|\nabla_Xu^+(\mathbf Z^+)|\lesssim\mathfrak C_0\, 	\frac{\omega^{\mathbf P^+}(B(\mathbf Z,20r))}{\sigma(B(\mathbf Z,r))}.
			\end{equation}
		\end{lemma}
		\begin{proof}
			We first claim that, for $\mathbf P\in \Omega$ and  $\mathbf P\neq\mathbf Y$, we have
			\begin{equation}\label{eq:caloric_potential_domination}
				\int\Gamma(\mathbf W,\mathbf Y)\,d\omega^{\mathbf P}(\mathbf W)\leq\Gamma(\mathbf P,\mathbf Y).
			\end{equation}
			Note that, since $0\leq\varphi_{\mathbf Z,r}\leq1$, \eqref{eq:caloric_potential_domination} proves $u^+\geq0$.
			
			To prove the claim, it suffices to approximate a unit point mass at $\mathbf Y$ by non-negative smooth compactly supported densities. Their heat potentials are bounded non-negative supercaloric functions vanishing at infinity, and therefore dominate the caloric extensions of their boundary values. Passing to the limit by Fatou proves \eqref{eq:caloric_potential_domination}; at the diagonal we assign the kernel the value zero.
			
			The adjoint caloricity of $u^+$ in the whole $B(\mathbf Z,10r)$  follows because both the pole and the support of $\varphi_{\mathbf Z,r}\omega^{\mathbf P^+}$ lie outside that ball.
			
			We are left with proving \eqref{eq:bound_T_caloric2}. Let $E\subset\supp\sigma\cap B(\mathbf Z,r)$ be a compact set with $\sigma(E)\geq\sigma(B(\mathbf Z,r))/2$, and let $v$ be the potential in Lemma \ref{lemma:caloric_boundary_potential}. By Fubini's theorem and \eqref{eq:caloric_potential_representation}, we have
			\begin{equation}\label{eq:bound_T_caloric3}
				\int_Eu^+\,d\sigma
				=v(\mathbf P^+)-\int v\, \varphi_{\mathbf Z,r}\,d\omega^{\mathbf P^+}
				=\int v\, (1-\varphi_{\mathbf Z,r})\,d\omega^{\mathbf P^+}
				\lesssim\mathfrak C_0r\,\omega^{\mathbf P^+}(B(\mathbf Z,20r)).
			\end{equation}
			On the other hand, Lemma \ref{lemma:Cauchy_estimates_heat} and the inequality $u^+\geq0$ imply
			\begin{equation}\label{eq:bound_T_caloric1}
				|\nabla_Xu^+(\mathbf Z^+)|\lesssim\frac1{r\sigma(E)}\int_Eu^+\,d\sigma.
			\end{equation}
			Thus, gathering the last two estimates proves \eqref{eq:bound_T_caloric2}. 	
		\end{proof}
		
		A time reflection yields an estimate analogous to \eqref{eq:bound_T_caloric2} for the function 
		\[
		u^-(\mathbf Y)\coloneqq \Gamma(\mathbf Y,\mathbf P^-)-\int\Gamma(\mathbf Y,\mathbf W)\varphi_{\mathbf Z,r}(\mathbf W)\,d\omega_*^{\mathbf P^-}(\mathbf W),
		\]
		at the time-lagged point $\mathbf Z^-\coloneqq (Z,t_Z-4r^2)$.
	}
	
	{
		\subsection{The maximal truncations on a testing ball}
		Let $B$ be a ball as in Definition \ref{definition:weak_A_infty} such that
		\begin{equation}\label{eq:caloric_balanced_ball_condition}
			\sigma((\delta B)^-)\gtrsim\sigma(B)\qquad \text{ and }\qquad \sigma((\delta B)^+)\gtrsim\sigma(B),
		\end{equation}
		with fixed constants. In particular, this holds for the balls of Lemma \ref{lemma:caloric_auxiliary_ball}.  To apply the local $Tb$ theorem, we first define the testing measures
		
		\begin{equation}\label{eq:testing_functions_caloric}
			\nu_B^*\coloneqq\frac{\sigma(B)}{\omega^{\mathbf P_B^+}(B)}\,\omega^{\mathbf P_B^+}|_B, \qquad \nu_B\coloneqq\frac{\sigma(B)}{\omega_*^{\mathbf P_B^-}(B)}\,\omega_*^{\mathbf P_B^-}|_B,
		\end{equation}
		and we remark that the denominators in \eqref{eq:testing_functions_caloric} are positive by \eqref{eq:nondegeneracy}. Both measures are non-negative and concentrated on $B$, and by definition we have
		
		\begin{equation}\label{eq:caloric_meas_accretivity}
			\nu_B(B)=\nu_B^*(B)=\sigma(B).
		\end{equation}
		Moreover, by Definition \ref{definition:weak_A_infty}, for every Borel set $E\subset B$ we have that:
		\begin{equation}\label{eq:size-tests}
			\text{If }\sigma(E)\leq \kappa\, \sigma(B), \qquad\text{ then }\qquad \nu_B(E)\leq \kappa'\, \sigma(B) \text{ and }\nu^*_B(E)\leq \kappa'\, \sigma(B).
		\end{equation}
	}
	
	{
		To verify the alternative \eqref{eq:testing_opt_2} of Definition \ref{def:accretive_testing_system}, we need pointwise bounds for the maximal truncations of $T$ and $T^*$ applied to these measures. We recall that the non-centered dilated maximal function in \eqref{eq:def_centered_maximal} of a non-negative measure $\nu$ is given by
		\[
		\widetilde{\mathcal M}_\sigma \nu(\mathbf Z)=\sup_{B'\ni\mathbf Z,\,\sigma(5B')>0}\frac{\nu(B')}{\sigma(5B')}, \qquad \text{ for }\mathbf Z\in \mathbb P^n.
		\]
	}
	
	{
		\begin{proposition}\label{proposition:testing_balls}
			For $\tau,\eta\in(0,1/10)$, the measures in \eqref{eq:testing_functions_caloric} satisfy
			\begin{equation}\label{eq:bound_ball_caloric_testing_1}
				T_{\sharp}^*\nu_B^*(\mathbf Z)\leq C_{\tau,\eta}\,\bigl(1+\widetilde{\mathcal M}_\sigma \nu_B^*(\mathbf Z)\bigr),\quad
				\mathbf Z\in \bigl((1-\tau)B\cap \supp\sigma\bigr) \setminus B\bigl(\mathbf P_B^+,\eta r(B)/2\bigr),
			\end{equation}
			and
			\begin{equation}\label{eq:bound_ball_caloric_testing_2}
				T_{\sharp}\nu_B(\mathbf Z)\leq C_{\tau,\eta}\,\bigl(1+\widetilde{\mathcal M}_\sigma \nu_B(\mathbf Z)\bigr),\quad
				\mathbf Z\in\bigl((1-\tau)B\cap \supp\sigma\bigr)\setminus B\bigl(\mathbf P_B^-,\eta r(B)/2\bigr),
			\end{equation}
			where $C_{\tau,\eta}\coloneqq C\min\{\tau,\eta\}^{-n-1}$, and $C>0$ depends only on $n$, $\mathfrak C_0$, $\vartheta$, $\delta$, and the constants in \eqref{eq:caloric_balanced_ball_condition}.
		\end{proposition}
		
		The proof of the proposition above relies on the following lemma.
		\begin{lemma}\label{lemma:caloric_doubling_truncation}
			We write $\mathbf P^+\coloneqq \mathbf P_B^+$, $\omega^+\coloneqq \omega^{\mathbf P^+}$, $a_B\coloneqq \sigma(B)/\omega^+(B)$, and fix $\mathbf Z$ as in the statement of Proposition \ref{proposition:testing_balls}. We also assume that $\widetilde C\coloneqq \widetilde{\mathcal M}_\sigma \nu_B^*(\mathbf Z)<\infty$. Let $r_*\coloneqq 10^{-4}\min\{\tau,\eta\}r(B)$. If $0<R\leq r_*$ is such that
			\begin{equation}\label{eq:caloric_intermediate_doubling}
				\sigma(B(\mathbf Z,100R))\leq 100^{n+2}\, \sigma(B(\mathbf Z,R)),
			\end{equation}
			then
			\begin{equation}\label{eq:caloric_doubling_truncation}
				|T_{R}^*\nu_B^*(\mathbf Z)|\leq C_{\tau,\eta}(1+\widetilde C).
			\end{equation}
		\end{lemma}
		\begin{proof}
			By the definition of $\widetilde C$ and the upper growth of $\sigma$, for $r>0$ we have
			\begin{equation}\label{eq:caloric_one_point_growth}
				\nu_B^*(B(\mathbf Z,r))\leq \widetilde C\, \sigma(B(\mathbf Z,5r))\leq 5^{n+1} \mathfrak C_0\, \widetilde C\, r^{n+1}.
			\end{equation}
			
			Let $\mathbf Z_R\coloneqq (Z,t_Z+4R^2)$. Since $R\leq10^{-4}\min\{\tau,\eta\}r(B)$, we have $\mathbf P^+\notin B(\mathbf Z,100R)$, $B(\mathbf Z,100R)\subset B$, and $\|\mathbf P^+-\mathbf Z_R\|\geq\eta r(B)/4$. Hence, by Lemma \ref{lemma:caloric_cutoff_estimate} with $r=R$ and differentiating \eqref{eq:definition_u_+} with respect to the space variable, we obtain
			\begin{equation}\label{eq:bound_T_caloric0}
				\int K(\mathbf W,\mathbf Z_R)\varphi_{\mathbf Z,R}(\mathbf W)\,d\omega^+(\mathbf W)=\nabla_Xu^+(\mathbf Z_R)+K(\mathbf P^+,\mathbf Z_R).
			\end{equation}
			
			To prove the desired bound, we denote $\nu\coloneqq \nu_B^*=a_B\,\omega^+|_B$, and we split
			\begin{equation}\label{eq:splitting_l99}
				\begin{split}
					&T_{R}^*\nu(\mathbf Z)=a_B\, \nabla_X u^+(\mathbf Z_R) + a_B K(\mathbf P^+, \mathbf Z_R)  + \int_{R< \|\mathbf W-\mathbf Z\|< 20 R}K(\mathbf W, \mathbf Z)\, d\nu(\mathbf W)\\
					&\qquad + \int_{\|\mathbf W-\mathbf Z\|\geq 20 R}\bigl(K(\mathbf W, \mathbf Z)-K(\mathbf W, \mathbf Z_R)\bigr)\, d\nu(\mathbf W)\\
					&\qquad - \int_{B(\mathbf Z, 20 R)}K(\mathbf W, \mathbf Z_R)\, \varphi_{\mathbf Z, R}(\mathbf W)\, d\nu(\mathbf W) - a_B\, \int_{\mathbb P^n\setminus B}K(\mathbf W, \mathbf Z_R)\varphi_{\mathbf Z, R}(\mathbf W)\, d\omega^+(\mathbf W)\\
					&\qquad\eqqcolon a_B\, \nabla_X u^+(\mathbf Z_R) + a_B K(\mathbf P^+, \mathbf Z_R) + \mathfrak I_1 + \mathfrak I_2 + \mathfrak I_3 + \mathfrak I_4.
				\end{split}
			\end{equation}
			By \eqref{eq:bound_T_caloric2}, the definition of $\widetilde C$,  and \eqref{eq:caloric_intermediate_doubling}, we have	
			\begin{equation}\label{eq:caloric_doubling_gradient}
				\begin{split}
					a_B\,\bigl|\nabla_Xu^+(\mathbf Z_R)\bigr|
					&\lesssim \frac{\mathfrak C_0}{\sigma(B(\mathbf Z,R))}\,\nu(B(\mathbf Z,20R))\\
					&\leq\mathfrak C_0\, \widetilde C\, \frac{\sigma(B(\mathbf Z,100R))}{\sigma(B(\mathbf Z,R))}\lesssim\mathfrak C_0\, \widetilde C.
				\end{split}
			\end{equation}
			Moreover, \eqref{eq:nondegeneracy} and standard kernel estimates yield
			\begin{equation}\label{eq:bound_final_remark}
				a_B\bigl|K(\mathbf P^+,\mathbf Z_R)\bigr|\lesssim \frac{\mathfrak C_0r(B)^{n+1}}{(\eta r(B))^{n+1}}=\mathfrak C_0\,\eta^{-n-1}.
			\end{equation}
			
			By the size estimate and \eqref{eq:caloric_one_point_growth}, we obtain
			\[
			\begin{split}
				|\mathfrak I_1| + |\mathfrak I_3| &\lesssim R^{-n-1}\nu(B(\mathbf Z,20R))\lesssim\mathfrak C_0\, \widetilde C.
			\end{split}
			\]
			To bound the second integral, note that kernel regularity and \eqref{eq:caloric_one_point_growth} give
			\[
			|\mathfrak I_2|\lesssim \int_{\|\mathbf W-\mathbf Z\|\geq20R}\frac{R}{\|\mathbf W-\mathbf Z\|^{n+2}}\,d\nu(\mathbf W)\lesssim\mathfrak C_0\, \widetilde C\, \sum_{j\geq 0}2^{-j}\lesssim\mathfrak C_0\, \widetilde C.
			\]
			Since $\|\mathbf Z_R-\mathbf Z\|=2R$, the distance from $\mathbf Z_R$ to $\mathbb P^n\setminus B$ is at least $\tau r(B)/2$. Moreover, $\omega^+(\partial\Omega)\leq1$, so 
			\[
			|\mathfrak I_4|\lesssim a_B\int_{\mathbb P^n\setminus B}|K(\mathbf W,\mathbf Z_R)|\varphi_{\mathbf Z,R}(\mathbf W)\,d\omega^+(\mathbf W)\lesssim\frac{a_B}{(\tau r(B))^{n+1}}\lesssim\mathfrak C_0\,\tau^{-n-1},
			\] 
			which, together with the previous bounds for the other terms in \eqref{eq:splitting_l99}, allows us to conclude \eqref{eq:caloric_doubling_truncation}.
		\end{proof}
		
		\vvv
		The proof of Proposition \ref{proposition:testing_balls} combines the lemma above with the standard kernel estimates. 
		\begin{proof}[Proof of Proposition \ref{proposition:testing_balls}]
			We prove \eqref{eq:bound_ball_caloric_testing_1} and note that \eqref{eq:bound_ball_caloric_testing_2} readily follows by time reversal. Moreover, we fix $\mathbf Z$ as in the statement and we can assume without loss of generality that $\widetilde C\coloneqq \widetilde{\mathcal M}_\sigma \nu_B^*(\mathbf Z)<\infty$. As in the proof of the lemma above, we denote $r_*\coloneqq 10^{-4}\min\{\tau,\eta\}r(B)$ and $\nu\coloneqq \nu^*_B$.
			
			Let $\varepsilon>0$. If $\varepsilon\geq r_*/100$, kernel size and \eqref{eq:caloric_meas_accretivity} give
			\begin{equation}\label{eq:bound_for_T_z}
				|T_{\varepsilon}^*\nu(\mathbf Z)|\lesssim\varepsilon^{-n-1}\sigma(B)\lesssim\mathfrak C_0\min\{\tau,\eta\}^{-n-1}.
			\end{equation}
			
			Let us now assume that $\varepsilon< r_*/100$ and, for $j\geq0$, define $R_j \coloneqq 100^j\, \varepsilon$. We denote by $k\geq0$ the smallest integer for which either \eqref{eq:caloric_intermediate_doubling} holds with $R=R_k$ or we have $R_k\geq r_*/100$. By minimality, $R_k<r_*$, and
			\begin{equation}\label{eq:caloric_nondoubling_growth}
				\sigma(B(\mathbf Z,R_{j+1}))>100^{n+2}\sigma(B(\mathbf Z,R_j)),\qquad\text{ for } 0\leq j<k.
			\end{equation}
			In either case, \eqref{eq:caloric_doubling_truncation} or the preceding size estimate gives $|T_{R_k}^*\nu(\mathbf Z)|\leq C_{\tau,\eta}(1+\widetilde C)$.
			
			Moreover,
			\begin{equation}\label{eq:caloric_nondoubling_sum}
				\bigl|T_{\varepsilon}^*\nu(\mathbf Z)-T_{R_k}^*\nu(\mathbf Z)\bigr|\lesssim\sum_{j=0}^{k-1}\frac{\nu(B(\mathbf Z,2R_{j+1}))}{R_j^{n+1}}\leq \widetilde C\, \sum_{j=0}^{k-1}\frac{\sigma(B(\mathbf Z,R_{j+2}))}{R_j^{n+1}}.
			\end{equation}
			We assume without loss of generality that $k\geq 1$, since otherwise the sum in \eqref{eq:caloric_nondoubling_sum} is empty. Note that upper growth implies
			\begin{equation}\label{eq:p98a}
				\frac{\sigma(B(\mathbf Z,R_{k+1}))}{R_{k-1}^{n+1}}\lesssim 100^{2(n+1)}\mathfrak C_0.
			\end{equation}
			For $j\leq k-2$, iterating \eqref{eq:caloric_nondoubling_growth} and again using upper growth yields
			\begin{equation}\label{eq:p98b}
				\frac{\sigma(B(\mathbf Z,R_{j+2}))}{R_j^{n+1}} \leq 100^{-(n+2)(k-j-2)}\, \frac{\sigma(B(\mathbf Z,R_k))}{R_j^{n+1}}\leq 100^{2(n+1)}\, \mathfrak C_0\,100^{-(k-j-2)}.
			\end{equation}
			Hence, we gather \eqref{eq:caloric_nondoubling_sum}, \eqref{eq:p98a}, \eqref{eq:p98b}, and obtain
			\[
			|T_{\varepsilon}^*\nu(\mathbf Z)-T_{R_k}^*\nu(\mathbf Z)|\lesssim\mathfrak C_0\,\widetilde C,
			\]
			which, combined with \eqref{eq:caloric_doubling_truncation} and \eqref{eq:bound_for_T_z}, proves \eqref{eq:bound_ball_caloric_testing_1}. 
		\end{proof}
	}

	{
		\subsection{The exceptional sets on the selected balls and the proof of Theorem \ref{theorem:application_caloric_measure}}\label{subsection:caloric_proof}
		Fix $\tau,\eta\in(0,1/10)$ and a ball $B$ as in Definition \ref{definition:weak_A_infty} satisfying \eqref{eq:caloric_balanced_ball_condition}. We set
		\[
		V_B\coloneqq B\cap B\bigl(\mathbf P_B^-,\eta r(B)/2\bigr),\qquad V_B^*\coloneqq B\cap B\bigl(\mathbf P_B^+,\eta r(B)/2\bigr),
		\]
		and, for $A>1$, we define
		\begin{equation}\label{eq:caloric_exceptional_sets}
			U_B\coloneqq V_B\cup\bigl\{\mathbf Z\in B:\widetilde{\mathcal M}_\sigma \nu_B(\mathbf Z)>A\bigr\},\qquad
			U_B^*\coloneqq V_B^*\cup \bigl\{\mathbf Z\in B:\widetilde{\mathcal M}_\sigma \nu_B^*(\mathbf Z)>A\bigr\}.
		\end{equation}
		\begin{lemma}\label{lemma:lemma_alternative_application}
			For every $A>1$, we have
			\begin{equation}\label{eq:caloric_maximal_exceptional_measure}
				\sigma(U_B)\leq\sigma(V_B)+A^{-1}\sigma(B)\qquad \text{ and }\qquad 
				\sigma(U_B^*)\leq\sigma(V_B^*)+A^{-1}\sigma(B).
			\end{equation}
			Moreover,
			\begin{equation}\label{eq:caloric_maximal_testing}
				T_{\sharp}\nu_B(\mathbf Z)\leq C_{\tau,\eta}(1+A),\qquad \text{ for }
				\mathbf Z\in\bigl((1-\tau)B\cap\supp\sigma\bigr)\setminus U_B,
			\end{equation}
			and
			\begin{equation}\label{eq:caloric_adjoint_maximal_testing}
				T_{\sharp}^*\nu_B^*(\mathbf Z)\leq C_{\tau,\eta}(1+A),\qquad \text{ for } \mathbf Z\in\bigl((1-\tau)B\cap\supp\sigma\bigr)\setminus U_B^*,
			\end{equation}
			where $C_{\tau,\eta}$ is the constant in Proposition \ref{proposition:testing_balls}.
		\end{lemma}
		\begin{proof}
			The estimates in \eqref{eq:caloric_maximal_exceptional_measure} follow from \eqref{eq:measure_testing_maximal_constants} and \eqref{eq:caloric_meas_accretivity}. If $\mathbf Z\in((1-\tau)B\cap\supp\sigma)\setminus U_B$, then $\mathbf Z\notin B(\mathbf P_B^-,\eta r(B)/2)$ and $\widetilde{\mathcal M}_\sigma \nu_B(\mathbf Z)\leq A$, so \eqref{eq:bound_ball_caloric_testing_2} gives \eqref{eq:caloric_maximal_testing}. The proof of \eqref{eq:caloric_adjoint_maximal_testing} is analogous, with \eqref{eq:bound_ball_caloric_testing_1} in place of \eqref{eq:bound_ball_caloric_testing_2}.
		\end{proof}
		\vvv
		
		We prove Theorem \ref{theorem:application_caloric_measure} by applying Theorem \ref{theorem:main_theorem_hofmann_problem} with $(n+1)$ in place of $n$ and $\mu\coloneqq\sigma$, with $\delta_0\coloneqq1/2$ and $\Lambda_0\coloneqq100$ in place of $\delta$ and $\Lambda$, and with $\beta_0\coloneqq\beta_0(D,n+1,\delta_0,\Lambda_0)$ in place of $\beta$, where $D$ is the geometric doubling constant of $\mathbb P^n$. In what follows, $B=B(\mathbf X_0,r)$ is a ball of the corresponding family $\mathscr B$, so that, in particular, $\mathbf X_0\in\supp\sigma$ and $\sigma(\Lambda_0B)\leq\beta_0\,\sigma(\delta_0B)$, while $\delta$ and $\Lambda$ denote the parameters in Definition \ref{definition:weak_A_infty}. 
		\begin{lemma}\label{lemma:caloric_auxiliary_ball}
			There exist $N=N(n,\delta,\Lambda)\geq 1$ and a ball $B'=B(\mathbf X_E,\rho)$ with $40$-small boundary, where $\mathbf X_E\in B/2\cap\supp\sigma$ and $r/(8\Lambda)\leq\rho\leq r/(4\Lambda)$, such that
			\begin{equation}\label{eq:caloric_auxiliary_ball}
				\Lambda B'\subset B,\qquad
				\sigma\bigl((\delta B')^\pm\bigr)\geq\frac{\sigma(B)}{8\beta_0N},\qquad \text{ and }\qquad
				\sigma(\Lambda B')\leq4\beta_0N\,\sigma(\delta B').
			\end{equation}
		\end{lemma}
		\begin{proof}
			We put $R\coloneqq r/(8\Lambda)$. Since $\sigma(B/2)\geq\beta_0^{-1}\sigma(B)$, there exists a compact set $F\subset B/2\cap\supp\sigma$ with $\sigma(F)\geq\sigma(B)/(2\beta_0)$. We cover $F$ by $N$ balls of radius $\delta R/4$, where $N$ depends only on $n$, $\delta$, and $\Lambda$, and we choose a compact subset $E$ of the intersection of $F$ with one member of this cover such that $\sigma(E)\geq\sigma(B)/(4\beta_0N)$. By Lemma \ref{lemma:halving_parabolic_balls}, there exists $\mathbf X_E=(X_E,t_E)\in E$ which bisects $\sigma|_E$ in time. Moreover, since $E$ is contained in a ball of radius $\delta R/4$, we have $E\subset B(\mathbf X_E,\delta R/2)$. Finally, by Lemma \ref{lemma:balls_small_boundary} with $\varepsilon=1$, we choose $R\leq\rho\leq2R$ so that $B'\coloneqq B(\mathbf X_E,\rho)$ has $40$-small boundary.
			
			The inclusion in \eqref{eq:caloric_auxiliary_ball} holds because $\mathbf X_E\in B/2$ and $\Lambda\rho\leq r/4$. Since $E\subset\delta B'$, the choice of $\mathbf X_E$ gives
			\[
			\sigma\bigl((\delta B')^\pm\bigr)\geq\frac12\,\sigma(E)\geq\frac{\sigma(B)}{8\beta_0N}
			\]
			and
			\[
			\sigma(\Lambda B')\leq\sigma(B)\leq 4\beta_0N\,\sigma(E)\leq4\beta_0N\,\sigma(\delta B'),
			\]
			which concludes the proof.
		\end{proof}
		
		Let $B'$ be the ball of Lemma \ref{lemma:caloric_auxiliary_ball}. By construction and \eqref{eq:caloric_auxiliary_ball}, we have $r(B')\approx r(B)$ and $\sigma(B')\approx\sigma(B)$, and $B'$ satisfies \eqref{eq:caloric_balanced_ball_condition}. Moreover, $B'$ satisfies the hypotheses in Definition \ref{definition:weak_A_infty}, provided that the doubling constant $\beta$ in Definition \ref{definition:weak_A_infty} is at least $4\beta_0N$, which we assume from now on. Indeed, $B'$ is centered on $\supp\sigma$, its radius is at most $r\leq \diam(\supp\sigma)$, it has $\mathfrak C_3$-small boundary because $\mathfrak C_3\geq 40$, and $\sigma(\Lambda B')\leq \beta \sigma(\delta B')$. Hence, the measures $\nu_{B'}$ and $\nu_{B'}^*$ in \eqref{eq:testing_functions_caloric} are defined, and we set
		\[
		\widehat\nu_B\coloneqq\frac{\sigma(B)}{\sigma(B')}\,\nu_{B'}\qquad \text{ and }\qquad \widehat\nu_B^*\coloneqq\frac{\sigma(B)}{\sigma(B')}\,\nu_{B'}^*.
		\]
		By \eqref{eq:caloric_meas_accretivity}, these measures are non-negative, supported in $\overline{B'}\subset B$, and satisfy
		\begin{equation}\label{eq:caloric_g_accretivity}
			\widehat\nu_B(B)=\|\widehat\nu_B\|=\sigma(B)\qquad \text{ and }\qquad \widehat\nu_B^*(B)=\|\widehat\nu_B^*\|=\sigma(B).
		\end{equation}
		\vv
		
		Moreover, since $\sigma(B)\leq4\beta_0N\,\sigma(B')$ by \eqref{eq:caloric_auxiliary_ball}, the condition \eqref{eq:size-tests} applied to $B'$ yields, for every Borel set $A\subset B$,
		\begin{equation}\label{eq:caloric_small_sets_final}
			\text{if }\sigma(A)\leq\frac{\kappa}{4\beta_0N}\, \sigma(B),\qquad\text{ then }\qquad \widehat\nu_B(A)\leq\kappa'\sigma(B)\quad\text{ and }\quad\widehat\nu_B^*(A)\leq\kappa'\sigma(B).
		\end{equation}
		Indeed, if $\sigma(A)\leq\kappa\,\sigma(B)/(4\beta_0N)$, then $\sigma(A\cap B')\leq\kappa\,\sigma(B')$, so 
		\[
		\widehat\nu_B(A)=\frac{\sigma(B)}{\sigma(B')}\,\nu_{B'}(A\cap B')\leq\kappa'\sigma(B),
		\]
		and the same holds for $\widehat\nu_B^*$.
		\vv
		
		Let $\zeta\in(0,1)$ and $A\coloneqq 8/\zeta$.  Since $B'$ has $40$-small boundary and $5B'\subset B$, we have
		\[
		\sigma\bigl((1+h)B'\setminus(1-h)B'\bigr)\leq40\,h\,\sigma(B)\qquad \text{ for } h\in(0,1).
		\]
		
		Hence, we fix $h\in(0,1/10)$ such that
		\[
		W_B\coloneqq(1+h)B'\setminus(1-h)B'\qquad \text{ satisfies }\qquad \sigma(W_B)\leq\frac{\zeta}{8}\,\sigma(B).
		\]
		For $\eta\in(0,1/10)$, let $V_{B'}$, $U_{B'}$, $V_{B'}^*$, and $U_{B'}^*$ be the sets associated with $B'$ and $A$ in \eqref{eq:caloric_exceptional_sets}. We define
		\[
		V_B\coloneqq V_{B'},\qquad U_B\coloneqq U_{B'}\cup W_B\cup(B\setminus\supp\sigma),\qquad V_B^*\coloneqq V_{B'}^*,\qquad U_B^*\coloneqq U_{B'}^*\cup W_B\cup(B\setminus\supp\sigma),
		\]
		and we prove that the alternative in the following lemma holds.
		\begin{lemma}\label{lemma:caloric_final_alternative}
			There exists $\mathfrak C_2=\mathfrak C_2(\eta,\zeta)>0$, depending only on $\eta$, $\zeta$, and the parameters in Theorem \ref{theorem:application_caloric_measure}, such that the following holds.
			\begin{itemize}
				\item If $\sigma(V_B)\geq\zeta\sigma(B)/2$, then $V_B\subset U_B$ satisfies \eqref{eq:testing_opt_1}.
				\item If $\sigma(V_B)<\zeta\sigma(B)/2$, then $\sigma(U_B)<3\zeta\sigma(B)/4$ and
				\begin{equation}\label{eq:caloric_final_testing}
					T_{\sharp}\widehat\nu_B(\mathbf Z)\leq\mathfrak C_2,\qquad \mathbf Z\in B\setminus U_B.
				\end{equation}
			\end{itemize}
			The analogous statement holds for $V_B^*$, $U_B^*$, and $T_{\sharp}^*\widehat\nu_B^*$.
		\end{lemma}
		\begin{proof}
			The first alternative is immediate because $\diam(V_B)\leq\eta r(B')\leq\eta r(B)$, so we assume that $\sigma(V_B)<\zeta\sigma(B)/2$. By \eqref{eq:caloric_maximal_exceptional_measure} applied to $B'$, the choices of $A$ and $h$, and $\sigma(B\setminus\supp\sigma)=0$, we have
			\begin{equation*}
				\begin{split}
					\sigma(U_B)&\leq\sigma(V_{B'})+A^{-1}\sigma(B')+\sigma(W_B)\\
					&<\frac\zeta2\,\sigma(B)+\frac\zeta8\,\sigma(B)+\frac\zeta8\,\sigma(B)=\frac34\,\zeta\sigma(B).
				\end{split}
			\end{equation*}
			
			To prove \eqref{eq:caloric_final_testing}, we fix $\mathbf Z\in B\setminus U_B$. Since $W_B\subset U_B$, either $\mathbf Z\in(1-h)B'$ or $\mathbf Z\notin(1+h)B'$.
			
			In the first case, we have $\mathbf Z\in((1-h)B'\cap\supp\sigma)\setminus U_{B'}$, so \eqref{eq:caloric_maximal_testing} applied to $B'$ with $\tau=h$ and the bound $\sigma(B)\leq4\beta_0N\sigma(B')$, which follows from \eqref{eq:caloric_auxiliary_ball}, readily give
			\[
			T_{\sharp}\widehat\nu_B(\mathbf Z)=\frac{\sigma(B)}{\sigma(B')}\,T_{\sharp}\nu_{B'}(\mathbf Z)\leq 4\beta_0N\,C_{h,\eta}\, (1+A).
			\]
			
			In the second case, the support of $\widehat \nu_B$ lies at distance at least $h\,r(B')$ from $\mathbf Z$. Thus, the size estimate, \eqref{eq:caloric_g_accretivity}, the $(n+1)$-growth of $\sigma$, and $r(B')\geq r/(8\Lambda)$ give
			\[
			T_{\sharp}\widehat\nu_B(\mathbf Z)\lesssim\frac{\sigma(B)}{(h\, r(B'))^{n+1}}\lesssim\mathfrak C_0\, \Bigl(\frac{8\Lambda}{h}\Bigr)^{n+1}.
			\]
			Since $h$ and $A$ depend only on $\zeta$, this proves \eqref{eq:caloric_final_testing}. The proof for $\widehat \nu_B^*$ is analogous, with \eqref{eq:caloric_adjoint_maximal_testing} in place of \eqref{eq:caloric_maximal_testing}.
		\end{proof}
		
		We finally conclude the section by gathering all the results above and applying our local $Tb$ theorem.
		\begin{proof}[Proof of Theorem \ref{theorem:application_caloric_measure}]
			We apply Theorem \ref{theorem:main_theorem_hofmann_problem}, as described before Lemma \ref{lemma:caloric_auxiliary_ball}, to each component of $T_\sigma$, with the testing measures $\widehat\nu_B$ and $\widehat\nu_B^*$ defined above, for $B\in\mathscr B$.
			
			By \eqref{eq:caloric_g_accretivity}, the measures $\widehat\nu_B$ and $\widehat\nu_B^*$ satisfy \eqref{eq:measure_testing_size} with $\mathfrak C_1\coloneqq1$. Moreover, since $\kappa'\leq1/200$, \eqref{eq:caloric_small_sets_final} gives \eqref{eq:measure_testing_small_sets} with $\mathfrak C_1=1$ and $\varepsilon_0\coloneqq\kappa/(4\beta_0N)$, which depends only on the parameters in Theorem \ref{theorem:application_caloric_measure}.
			
			We fix $\tau\in(0,\tau_0)$ and $\zeta\in(0,\zeta_0)$, with $\tau_0$ and $\zeta_0$ as in \eqref{eq:measure_testing_tau_zeta}. Since \eqref{eq:caloric_final_testing} holds on the whole $B\setminus U_B$, Lemma \ref{lemma:caloric_final_alternative} shows that, for every $\eta\in(0,1/10)$, one of the alternatives \eqref{eq:testing_opt_1} and \eqref{eq:testing_opt_2} holds with $\mathfrak C_2(\eta)\coloneqq\mathfrak C_2(\eta,\zeta)$, both for $\widehat\nu_B$ and $T_\sigma$ and for $\widehat\nu_B^*$ and $T_\sigma^*$. The maximal truncations of each component of $T_\sigma$ and $T_\sigma^*$ are bounded by those of $T_\sigma$ and $T_\sigma^*$, respectively. Therefore, the families $\{\widehat\nu_B\}_{B\in\mathscr B}$ and $\{\widehat\nu_B^*\}_{B\in\mathscr B}$ are accretive measure testing systems for each component of $T_\sigma$ and for its adjoint, in the sense of Definition \ref{def:accretive_testing_system}. Finally, Theorem \ref{theorem:main_theorem_hofmann_problem} applies to each component of $T_\sigma$, and summing over the components concludes the proof of Theorem \ref{theorem:application_caloric_measure}.
		\end{proof}
		
		\vv
		\begin{remark}\label{remark:caloric_dependence_eta}
			The level $\mathfrak C_2(\eta,\zeta)$ in Lemma \ref{lemma:caloric_final_alternative} associated with the $L^0$ condition depends on $\eta$ through the constant $C_{h,\eta}$ of Proposition \ref{proposition:testing_balls}, which comes from the bound \eqref{eq:bound_final_remark}.
			The poles $\mathbf P_B^\pm$ of caloric measures may be arbitrarily close to $\partial\Omega$, so our argument does not avoid this dependence, and the alternative \eqref{eq:testing_opt_1} in the testing dichotomy takes care of the balls on which $\sigma$ concentrates near the poles.
			
			However, we note that the bound in Theorem \ref{theorem:main_theorem_hofmann_problem} depends on the levels only through $\mathfrak C_2(\eta_*)$, and $\eta_*$ depends only on $n$, $\zeta$, and the structural parameters, see \eqref{eq:choice_eta_Stop_1}. Thus, in Theorem \ref{theorem:application_caloric_measure}, the level is ultimately determined by $\zeta$, as in the function-testing setting of Corollary \ref{corollary:hofmann_problem}, where it does not depend on $\eta$ at all.
		\end{remark}
		
	}

	\appendix
	
	\section{David-Mattila, random, and Whitney-type cubes in metric spaces}\label{appendix:random_cubes}
	
	Throughout this appendix, $(X,d)$ is a geometrically doubling metric space with parameter $D$.

	{
		\subsection{Construction of the David-Mattila testing balls}\label{appendix:DM_testing_balls}
		{
			Let $\sigma$ be a non-negative Radon measure on $X$. If $\supp\sigma$ is compact and has positive diameter, let $\mathcal D_\sigma$ be the lattice of Lemma \ref{lemma:DM_cubes} and, otherwise, use the lattice $\mathcal D^{(s)}_\sigma$ of Lemma \ref{lemma:DM_cubes_one_sided}, for a fixed $s>0$. The arguments below also apply to the latter lattice when $\supp\sigma$ is compact. We write $\mathcal D_\sigma$ in both cases, and use the notation $r_Q$ and $N_j^{\mathrm{int}}(Q)$ for the radii and boundary strips. The auxiliary constants are chosen as in Section \ref{section:auxiliary_DM_constants}, with the same boundary-decay constant for both lattices.
			\begin{lemma}\label{lemma:construction_testing_ball_DM}
				For every $Q\in\caD^{db}_\sigma$, there exist $x_Q\in Q\setminus N_{j_*}^{\textup{int}}(Q)$ and $\rho_Q>0$ such that the ball $\widetilde B(Q)\coloneqq B(x_Q, \rho_Q)$ satisfies the following properties:
				\begin{enumerate}
					\item We have that $\rho_Q\approx r_Q$. More precisely, it holds
					\begin{equation}\label{eq:structure_testing_balls}
						\frac{c_*}{4\Lambda}\,  r_Q\leq \rho_Q\leq \frac{c_*}{2\Lambda}\, r_Q.
					\end{equation}
					\item We have the inclusion
					\[
					\Lambda\,  \widetilde B(Q)\cap \supp(\sigma)\subset Q\subset M_*\, \widetilde B(Q).
					\]
					\item $\widetilde B(Q)$ has the doubling property
					\begin{equation}\label{eq:doubling_property_testing_DM_balls}
						\sigma(\Lambda\, \widetilde B(Q))\leq \sigma(Q)\leq  \beta_*\, \sigma(\delta\, \widetilde B(Q)).
					\end{equation}
					\item $\widetilde B(Q)$ has $40$-small boundary with respect to $\sigma$. Namely, for every $0<\lambda\leq 1$ we have
					\begin{equation}\label{eq:thin_boundary_testing_ball_DM}
						\sigma\bigl(\bigl\{x\in 5\widetilde B(Q): |\dist (x,x_Q)-\rho_Q|\leq \lambda\,  \rho_Q\bigr\}\bigr)\leq 40	 \lambda\,  \sigma (5\widetilde B(Q)).
					\end{equation}
				\end{enumerate}
			\end{lemma}
			\begin{proof}
				First, we define the center $x_Q$ of the ball $\widetilde B(Q)$ using a covering argument applied to the set $F_Q\coloneqq Q\setminus N_{j_*}^{\textup{int}}(Q)$. 
				By the properties of David-Mattila cells, we have that $F_Q\subset Q\subset 28 B(Q)$. We define $R_Q\coloneqq \tfrac{c_*}{4\Lambda}r_Q$, and note that the definition of $m_*$ yields
				\begin{equation}\label{eq:property_def_m_*}
					\frac{28}{2^{m_*}}\, r_Q\leq \frac{\delta}{2}\, R_Q.
				\end{equation}
				
				By the geometrically doubling property, $28B(Q)$ can be covered by at most $D^{m_*}$ balls of radius $28r_Q/2^{m_*}$. Discard the balls that do not meet $F_Q$, and recenter each remaining ball at a point of its intersection with $F_Q$, doubling its radius. 
				By \eqref{eq:property_def_m_*}, the resulting balls have radius at most $\delta R_Q$. Since they are at most $D^{m_*}$ and $\sigma(F_Q)\geq \sigma(Q)/2$, by Lemma \ref{lemma:definition_F_Q_separation}, one of them, whose center we denote by $x_Q$, is such that $E_Q\coloneqq F_Q\cap B(x_Q, \delta R_Q)$ satisfies
				\begin{equation}\label{eq:lemma_testing_DM_ph}
					\sigma(E_Q)\geq \frac{\sigma(Q)}{2D^{m_*}}.
				\end{equation}
				To define the radius of $\widetilde B(Q)$, we apply Lemma \ref{lemma:balls_small_boundary}, which implies that there exists $\rho_Q\in [R_Q, 2R_Q]$ such that the ball $\widetilde B(Q)\coloneqq B(x_Q, \rho_Q)$ has $40$-small boundary with respect to the measure $\sigma$. Hence, the property (4) holds by construction.
				Property (1) holds by the choice of the radii.
				
				To prove (2), it is enough to observe that
				\[
				\Lambda\rho_Q\leq2\Lambda R_Q = \frac{c_*}{2}\, r_Q<\dist(x_Q, \supp(\sigma)\setminus Q)
				\]
				by Lemma \ref{lemma:definition_F_Q_separation}, since $x_Q\in F_Q$. Moreover, if $y\in Q$, then, since $x_Q, y\in Q\subset 28\, B(Q)$, we have
				\[
				d(y,x_Q)<56\,r_Q\leq\frac{224\Lambda}{c_*}\,\rho_Q=M_*\rho_Q.
				\]
				
				Finally, the first inequality in \eqref{eq:doubling_property_testing_DM_balls} follows from (2). To prove the second one, note that the inequality $\rho_Q\geq R_Q$ implies $B(x_Q, \delta R_Q)\subset \delta \widetilde B(Q)$. Thus, \eqref{eq:lemma_testing_DM_ph} yields
				\[
				\sigma(Q)\leq 2D^{m_*}\, \sigma(\delta \widetilde B(Q)) = \beta_*\, \sigma(\delta\, \widetilde B(Q)).\qedhere
				\]
			\end{proof}
			\vv
			
			Since, see Remark \ref{remark:choice_of_beta}, we assume $\beta>\beta_*$, the balls $\widetilde B(R)$ associated with doubling David-Mattila cells $R\in \caD^{db}_\sigma$ satisfy
			\[
			\sigma(\Lambda\, \widetilde B(R)) \leq \beta\, \sigma(\delta\, \widetilde B(R)),
			\]
			so they satisfy the doubling and small-boundary requirements of Theorem \ref{theorem:main_theorem_hofmann_problem}. The radius requirement is automatic when $\supp\sigma$ is unbounded. If $0<\diam W<\infty$, where $W=\supp\sigma$, it also holds for the lattice of Lemma \ref{lemma:DM_cubes} and for the lattice of Lemma \ref{lemma:DM_cubes_one_sided} with $s\leq2\diam W$, since \eqref{eq:structure_testing_balls} gives $\rho_Q\leq c_*K_0s/(2\Lambda)\leq\diam W$. This includes the scale used in Section \ref{sec:corona}.
			
			By a standard application of Lemma \ref{lemma:balls_small_boundary}, for any given non-negative Radon measure $\mathfrak m$ on $X$, the radius of $\widetilde B(Q)$ constructed in Lemma \ref{lemma:construction_testing_ball_DM} can be selected so that the ball satisfies all the properties above, and an additional analogue of \eqref{eq:thin_boundary_testing_ball_DM} with $\mathfrak m$ in place of $\sigma$ and $80$ in place of $40$. Indeed, in the proof of Lemma \ref{lemma:balls_small_boundary}, the excluded radii occupy at most half of the interval for $\sigma$ and, with the threshold doubled, at most one quarter for $\mathfrak m$, so a common admissible radius exists.
			
			Testing balls associated with doubling David-Mattila cells satisfy the following, additional properties.
			\begin{lemma}\label{lemma:inclusion_testing_balls}
				For $Q,R\in \caD^{db}_\sigma$ such that $Q\subset R$ and $J(Q)>J(R)$, the following properties hold:
				\begin{itemize}
					\item The radii $\rho_Q$ and $\rho_R$ are such that
					\begin{equation}\label{eq:nesting_inclusion_testing_0}
						\frac{\rho_Q}{\rho_R}\leq 2K_0 A^{-(J(Q)-J(R))}_0\leq \frac{2K_0}{A_0}<\frac{1}{2}.
					\end{equation}
					\item We have the inclusion
					\begin{equation}\label{eq:nesting_inclusion_testing}
						2M_*\, \widetilde B(Q)\subset 2M_*\, \widetilde B(R	).
					\end{equation}
					\item If $B(Q)$ denotes the ball associated with $Q$ as defined in Lemma \ref{lemma:DM_cubes} (or Lemma \ref{lemma:DM_cubes_one_sided}), we have
					\begin{equation}\label{eq:inclusions_84_100_B_B_tilde}
						84\, B(Q)\subset 2M_*\, \widetilde B(Q)\qquad \text{ and }\qquad 100 B(Q)\subset 2 A_{\Lambda, *}\widetilde B(Q).
					\end{equation}
				\end{itemize}
			\end{lemma}
			\begin{proof}
				
				The estimate \eqref{eq:structure_testing_balls} and the bounds $sA_0^{-J(P)}\leq r_P\leq K_0sA_0^{-J(P)}$, for $P\in\{Q,R\}$, where $s=\diam(\supp\sigma)$ for the lattice of Lemma \ref{lemma:DM_cubes}, give $\rho_Q/\rho_R\leq2r_Q/r_R\leq2K_0A_0^{-(J(Q)-J(R))}$, and \eqref{eq:nesting_inclusion_testing_0} follows from $A_0>5000K_0$.
				To prove \eqref{eq:nesting_inclusion_testing}, we first note that, since $x_Q\in Q$ and $Q\subset R\subset 28 B(R)$, we have
				\begin{equation}\label{eq:nesting_aux_1}
					d(x_Q, x_R)<56\, r_R\leq M_*\, \rho_R.
				\end{equation}
				For $y\in2M_*\widetilde B(Q)$, the triangle inequality and \eqref{eq:nesting_inclusion_testing_0} therefore give $d(y,x_R)<2M_*\rho_Q+M_*\rho_R<2M_*\rho_R$, which proves \eqref{eq:nesting_inclusion_testing}.	
				Finally, since $x_Q\in Q\subset28B(Q)$, \eqref{eq:structure_testing_balls} and the definitions of $M_*$ and $A_{\Lambda,*}$ give $d(y,x_Q)<112\,r_Q\leq2M_*\rho_Q$ for $y\in84B(Q)$, and $d(y,x_Q)<128\,r_Q\leq2A_{\Lambda,*}\rho_Q$ for $y\in100B(Q)$, which proves \eqref{eq:inclusions_84_100_B_B_tilde}.
			\end{proof}
			
		}
	}
	
	\subsection{Properties of random cubes} First, we consider the dyadic cubes associated with $(X,d)$ as in Theorem \ref{thm:HK22}.

	\begin{theorem}[\cite{HK12}, Theorem 5.6]\label{thm:HK56}
		Let $(\Omega,\mathbb P)$ be a probability space as in Theorem \ref{thm:HK51} with parameter $D$, and suppose that $144\el\leq 1$. For a dyadic cube $Q$, define the Hytönen-Kairema boundary zone by
		\[
		\partial(Q,\varepsilon)\coloneqq \{x\in\overline Q:\mathrm{dist}(x,X\setminus \widetilde Q)\leq \varepsilon\},
		\qquad \varepsilon>0.
		\]
		Then there exist constants $C_2,\eta>0$, with $C_2$ depending only on $c_1$ and $\el$ and $\eta$ depending on $\el$ and $\tau_0$, the constant constructed in Theorem \ref{thm:HK51}, such that, for every $x\in X$, every $\tau>0$ and every $k\in\mathbb Z$,
		\[
		\mathbb P\Bigl(\Bigl\{w\in\Omega:
		x\in\bigcup_{\alpha\in A_k}\partial(Q^k_\alpha(w),\tau\el^k)
		\Bigr\}\Bigr)\leq C_2\tau^\eta.
		\]
		Moreover, for every $x\in X$,
		\[
		\mathbb P\Bigl(\Bigl\{w\in\Omega:x\in\bigcup_{k,\alpha}\partial Q^k_\alpha(w)\Bigr\}\Bigr)=0,
		\]
		and, if $\mu$ is a positive $\sigma$-finite measure on $X$, then
		\[
		\mu\Bigl(\bigcup_{k,\alpha}\partial Q^k_\alpha(w)\Bigr)=0,	\qquad\text{ for } \mathbb P\text{-a.e. }w	\in\Omega.
		\]
	\end{theorem}

	For $i=1,2$, let $\bbD_i(w_i)$ be two independent copies of the random dyadic system constructed above. If $Q\in\bbD_i(w_i)$ and $Q=Q^k_\alpha(w_i)$ for some $k\in\mathbb Z$, we recall that $\ell(Q) = \el^k.$

	We also set
	\[
	\bbD_{i,k}(w_i)\coloneqq \{Q\in\bbD_i(w_i):\ell(Q)=\el^k\},\qquad k\in\mathbb Z,
	\]
	and we refer to \textit{good} and \textit{bad} cubes as in Definition \ref{def:good_bad_cubes}.
	Recall that a cube $Q\in\bbD_1(w_1)$ is bad if there exists $R\in\bbD_2(w_2)$ such that $\ell(R)\geq\el^{-m}\ell(Q)$ and
	\[
	\dist(Q,R)\leq\ell(Q)^\alpha\ell(R)^{1-\alpha},\qquad
	\dist(Q,X\setminus R)\leq\ell(Q)^\alpha\ell(R)^{1-\alpha}.
	\]
	The cube $Q$ is good if no such $R$ exists.

	{
		\begin{lemma}\label{lem:random_capture_goodness}
			Let $0<\alpha<1$ and $\kappa\in(0,1)$. There exist $L<\infty$ and
			$m_0\in\mathbb N$, depending only on $\el,C_1,C_2,\eta,\alpha,\kappa$
			and the structural constant $M_*$, with the following property.
			Let $B_i=B(x_i,r_i)$, $i=1,2$, where $r_1\leq r_2$ and
			$2M_*B_1\subset2M_*B_2$. One can fix generations $k_i$ before
			choosing either grid so that, if $Q_i(w_j)$ is the unique cube of
			$\bbD_{j,k_i}(w_j)$ containing $x_i$, then for every $m\geq m_0$,
			\[
			\begin{split}
				(\mathbb P\otimes\mathbb P)\Bigl(\bigl\{\mathtt w\in\Omega^2: &\ 6M_*B_i\subset Q_i(w_j),\quad\diam Q_i(w_j)\leq Lr_i,\\
				&\ Q_1(w_j)\subset Q_2(w_j),\quad Q_i(w_j)\in\bbD_j^{\mathcal G}(\mathtt w,m),\quad i,j=1,2\bigr\}\Bigr)
				>1-\kappa.
			\end{split}
			\]
			In particular, $L$ and the chosen generations are independent of $m\geq m_0$.
		\end{lemma}
		\begin{proof}
			We choose $\varepsilon\in(0,1)$ such that
			\[
			2(C_2+1)\varepsilon^\eta=\frac\kappa4,
			\]
			and let $k_i$ be the largest integer such that
			\[
			\el^{k_i}\geq\frac{6M_*r_i}{\varepsilon},
			\]
			which yields $\el^{k_i}<6M_*r_i/(\varepsilon\el)$. Then, we set
			\(
			E_i^j\coloneqq\{w_j:6M_*B_i\subset Q_i(w_j)\},
			\)
			and note that, if $w_j\notin E_i^j$, then
			\[
			\dist(x_i,X\setminus Q_i(w_j))\leq6M_*r_i\leq\varepsilon\el^{k_i},
			\]
			so $x_i\in\partial(Q_i(w_j),\varepsilon\el^{k_i})$. By Theorem \ref{thm:HK56}, we deduce that
			\begin{equation}\label{stimaintersezione}
				\mathbb P(E_1^j\cap E_2^j)
				\geq1-2C_2\varepsilon^\eta>1-\frac\kappa4,
				\qquad j=1,2.
			\end{equation}
			We also note that, if $w_j\in E_1^j\cap E_2^j$, then $x_1\in2M_*B_1\subset2M_*B_2\subset Q_2(w_j)$. Since $k_1\geq k_2$, the nesting property gives $Q_1(w_j)\subset Q_2(w_j)$.
			Moreover, we have
			\[
			\diam Q_i(w_j)\leq2C_1\el^{k_i}
			<\frac{12M_*C_1}{\varepsilon\el}r_i.
			\]
			
			Fix $w_1$ and put $Q=Q_i(w_1)$ and $k=k_i$.
			If $Q$ is bad, there are $h\leq k-m$ and $R\in\bbD_{2,h}(w_2)$
			such that the distances from $Q$ to both $R$ and $X\setminus R$
			are at most $\el^{\alpha k}\el^{(1-\alpha)h}$. Since $x_i\in Q$,
			the distances from $x_i$ to these two sets are at most
			\[
			t_h\coloneqq\el^{\alpha k}\el^{(1-\alpha)h}+2C_1\el^k.
			\]
			Let $S\in\bbD_{2,h}(w_2)$ contain $x_i$. If $x_i\in R$, then $S=R$ and $\dist(x_i,X\setminus S)\leq t_h$. Otherwise $R\subset X\setminus S$, so again $\dist(x_i,X\setminus S)\leq t_h$. Therefore $x_i\in\partial(S,t_h)$, and Theorem \ref{thm:HK56} gives
			\begin{equation}\label{eq:probability_estimate_aux}
				\begin{split}
					\mathbb P_{w_2}\bigl(\{Q_i(w_1)\text{ is bad}\}\bigr)
					&\leq C_2\sum_{h\leq k-m}
					\bigl(\el^{\alpha(k-h)}+2C_1\el^{k-h}\bigr)^\eta\\
					&\leq C_2(1+2C_1)^\eta\sum_{j=m}^\infty\el^{\alpha\eta j}
					=\frac{C_2(1+2C_1)^\eta}{1-\el^{\alpha\eta}}
					\el^{\alpha\eta m}.
				\end{split}
			\end{equation}
			Choose $m_0$ so that the last expression is less than $\kappa/16$
			for every $m\geq m_0$. Interchanging the two grids yields
			\begin{equation}\label{stimesxz}
				\mathbb P_{w_2}\bigl(\{Q_i(w_1)\text{ is bad}\}\bigr)<\frac\kappa{16}
				\qquad \text{ and }\qquad 
				\mathbb P_{w_1}\bigl(\{Q_i(w_2)\text{ is bad}\}\bigr)<\frac\kappa{16},
				\qquad i=1,2.
			\end{equation}
			By Fubini's theorem, \eqref{stimaintersezione}, and \eqref{stimesxz}, we finally obtain
			\[
			(\mathbb P\otimes\mathbb P)\Bigl(\bigcap_{j=1}^2(E_1^j\cap E_2^j)\cap\bigcap_{i,j=1}^2\{Q_i(w_j)\in\bbD_j^{\mathcal G}(\mathtt w,m)\}\Bigr)>1-\frac\kappa2-\frac{4\kappa}{16}>1-\kappa.
			\]
			This proves the lemma with $L=12M_*C_1/(\varepsilon\el)$.
		\end{proof}
	}
	
	\vv
	
	We conclude this appendix with the proof of Lemma \ref{prop:whitney_cubes_metric}, which provides the \emph{Whitney-type decomposition} of dyadic cubes used for the sum over $\mathscr S'_{4,\term}$ in Section \ref{section:suppressed_mixed_good}, see Lemma \ref{lem:bound_s_term}.
	
	\subsection{Proof of Lemma \ref{prop:whitney_cubes_metric}}
	{
		Fix $\mathbf c_W=2$. We consider all cubes $S\in\bbD$, with their generations, such that $S\subset Q$ and $\mathbf c_W\ell(S)\leq\dist(S,X\setminus Q)$, and define $\mathcal W(Q)$ to be the family of those cubes which are maximal with respect to dyadic ancestry. To see that this family covers the interior of $Q$, let $x\in\operatorname{int}Q$ and denote by $S_k(x)$ the cube of generation $k$ containing $x$. For all sufficiently large $k$, the diameter bound in Theorem \ref{thm:HK22} gives
		\[
		\dist(S_k(x),X\setminus Q)\geq\dist(x,X\setminus Q)-2C_1\el^k\geq\mathbf c_W\el^k.
		\]
		Thus $S_k(x)\subset Q$ and $\mathbf c_W\ell(S_k(x))\leq\dist(S_k(x),X\setminus Q)$. If an ancestor $S$ of $S_k(x)$ satisfies the same conditions, then, for any fixed $y\in X\setminus Q$,
		\[
		\mathbf c_W\ell(S)\leq\dist(S,X\setminus Q)\leq d(x,y).
		\]
		There are only finitely many such ancestors, since their side lengths tend to infinity. Consequently, $x$ belongs to a cube of $\mathcal W(Q)$. Each cube in $\mathcal W(Q)$ lies in $\operatorname{int}Q$, and maximality and nesting imply that these cubes are pairwise disjoint. This proves (1).
		
		The lower bound in (2) follows from the definition of $\mathcal W(Q)$. If $S\in\mathcal W(Q)$ and $\widehat S$ is its dyadic parent, then either $\widehat S\not\subset Q$, so $\dist(\widehat S,X\setminus Q)=0$, or maximality gives $\dist(\widehat S,X\setminus Q)<\mathbf c_W\ell(\widehat S)$. Hence,
		\[
		\dist(S,X\setminus Q)\leq\diam(\widehat S)+\dist(\widehat S,X\setminus Q)\leq(2C_1+\mathbf c_W)\ell(\widehat S)=\frac{2C_1+\mathbf c_W}{\el}\ell(S).
		\]
		This proves (2) with $\mathbf C_W=(2C_1+\mathbf c_W)/\el$.
	}
	
	To prove \eqref{eq:bounded_overlap}, we fix $x\in X$, we set
	\[
	\mathcal W_x\coloneqq\bigl\{S\in\mathcal W(Q):x\in2S\bigr\},
	\]
	and we may assume that $\mathcal W_x\neq\varnothing$. Let $S\in\mathcal W_x$. Since $\dist(x,S)<\ell(S)$, (2), the triangle inequality, and the diameter bound in Theorem \ref{thm:HK22} give
	\[
	\begin{split}
		(\mathbf c_W-1)\,\ell(S)&\leq\dist(S,X\setminus Q)-\dist(x,S)\leq\dist(x,X\setminus Q)\\
		&\leq\dist(x,S)+\diam S+\dist(S,X\setminus Q)\leq(\mathbf C_W+1+2C_1)\,\ell(S).
	\end{split}
	\]
	Hence, if we fix $S_0\in\mathcal W_x$ and set $\rho_W\coloneqq(\mathbf c_W-1)/(\mathbf C_W+1+2C_1)$, we have 
	\[
	\rho_W\ell(S_0)\leq\ell(S)\leq\rho_W^{-1}\ell(S_0)\qquad \text{ for every } S\in\mathcal W_x.
	\]
	
	Thus, since $B(z_S,c_1\ell(S))\subset S$, by Theorem \ref{thm:HK22}, and the cubes in $\mathcal W_x$ are pairwise disjoint, by (1), we obtain
	\begin{equation}\label{eq:separation_centers_whitney}
		d(z_S,z_{S'})\geq c_1\rho_W\ell(S_0),\qquad S,S'\in\mathcal W_x,\quad S\neq S'.
	\end{equation}
	Moreover, since $S\subset B(z_S,C_1\ell(S))$, we have
	\[
	d(x,z_S)\leq\dist(x,S)+C_1\ell(S)<(1+C_1)\rho_W^{-1}\ell(S_0),\qquad S\in\mathcal W_x.
	\]
	
	Now, we denote
	\[
	N\coloneqq\Bigl\lceil\log_2\frac{3(2+C_1)}{c_1\rho_W^2}\Bigr\rceil+1.
	\]
	By the geometric doubling property, the ball $B\bigl(x,(2+C_1)\rho_W^{-1}\ell(S_0)\bigr)$, which contains the centers $z_S$, for $S\in\mathcal W_x$, can be covered by $D^N$ balls of radius $(2+C_1)\rho_W^{-1}\ell(S_0)/2^N<c_1\rho_W\ell(S_0)/3$. By \eqref{eq:separation_centers_whitney}, each of these balls contains at most one center $z_S$, so that $\#\mathcal W_x\leq D^N$. Since $x\in X$ is arbitrary, this proves \eqref{eq:bounded_overlap} with $\mathbf c\coloneqq D^N$, and concludes the proof of Lemma \ref{prop:whitney_cubes_metric}.
	
	\section{Suppressed kernels associated with a dominating function}\label{appendix:suppression_metric_spaces}
	
	\subsection{Dominating functions}
	Throughout this appendix, as in \cite{Hy10}, we assume $(X,d,\mu)$ to be a geometrically doubling metric measure space with $\mu\neq 0$ and $\lambda\colon X\times(0,\infty)\to(0,\infty)$ a so-called \textit{dominating function} such that, for some $C_\lambda >0$, the following properties hold:
	\begin{enumerate}
		\item[(i)] The function $\lambda (x, \cdot)$ is non-decreasing and $\mu(B(x,r))\le \lambda(x,r)$, for all $x\in X$, $r>0$.
		\item[(ii)] We have that
		\begin{equation}\label{eq:UD2_property}
			\lambda(x,2r)\le C_\lambda\,\lambda(x,r), \qquad \text{ for all }x\in X, r>0.
		\end{equation}
	\end{enumerate}
	
	\vv
	
	We first construct a dominating function comparable to $\lambda$ and locally Lipschitz in the radius.
	
	\begin{lemma}\label{lemma:aux_lambda_2}
		Let $(X,d, \mu)$ be a geometrically doubling metric measure space, let $\lambda$ be a dominating function,	and define
		\[
		\widetilde\lambda(x,r)\coloneqq \frac{1}{\log 2}\, \int_r^{2r}\lambda(x,\rho)\,\frac{d\rho}{\rho}, \qquad \text{ for }x\in X, r>0.
		\]
		Then, for every $x\in X$ and $r>0$ we have
		\begin{equation}\label{eq:comparability_lambda_lambda_tilde}
			\lambda(x,r)\leq \widetilde\lambda(x,r)\leq C_\lambda \lambda(x,r)\qquad \text{ and }\qquad \widetilde \lambda(x,2r)\leq C_\lambda\, \widetilde \lambda(x,r).
		\end{equation}
		Moreover, there exists a constant $\mathbf C>0$ depending on $C_\lambda$ such that, for every $x\in X$, every $m>0$ and every $s,t\in[m/2,m]$, we have
		\begin{equation}\label{eq:regularity_estimate_lambda_tilde}
			\frac{|\widetilde\lambda(x,t)-\widetilde\lambda(x,s)|}{\widetilde\lambda(x,m)}
			\leq 	\mathbf C\, \frac{|t-s|}{m} .
		\end{equation}
	\end{lemma}
	\begin{proof}
		Let 	$\rho\in[r,2r]$. Properties (i) and (ii) imply
		\[
		\lambda(x,r)\leq \lambda(x,\rho)\le \lambda(x,2r) \leq C_\lambda \, \lambda(x,r),
		\]
		which readily implies the first estimate in \eqref{eq:comparability_lambda_lambda_tilde}. Moreover, a change of variables and \eqref{eq:UD2_property} give the second one. 
		
		We now prove the regularity estimate \eqref{eq:regularity_estimate_lambda_tilde}. By symmetry, we assume that $s<t$. Then, the definition of $\widetilde \lambda$ implies
		\[
		\widetilde\lambda(x,t)-\widetilde\lambda(x,s)	=	\frac{1}{\log 2} \, \Bigl(\int_t^{2t}\lambda(x,\rho)\,\frac{d\rho}{\rho}-\int_s^{2s}\lambda(x,\rho)\,\frac{d\rho}{\rho}\Bigr).
		\]
		Hence, denoting by $\triangle$ the symmetric difference, we have
		\begin{equation*}\label{eq:tilde_lambda_1}
			|\widetilde\lambda(x,t)-\widetilde\lambda(x,s)|	\leq	\frac{1}{\log 2}\, 	\int_{[s,2s]\triangle[t,2t]}\lambda(x,\rho)\,\frac{d\rho}{\rho}.
		\end{equation*}
		
		Since $s,t\in[m/2,m]$, we have the inclusion $[s,2s]\triangle[t,2t]\subset [m/2,2m]$. Thus, for
		$\rho\in[m/2,2m]$, monotonicity and doubling give
		\begin{equation}\label{eq:tilde_lambda_2}
			\lambda(x,\rho)\le \lambda(x,2m)\le C_\lambda \lambda(x,m)\overset{\eqref{eq:comparability_lambda_lambda_tilde}}{\leq} C_\lambda \widetilde\lambda(x,m).
		\end{equation}
		Therefore, we obtain	
		\begin{equation}\label{eq:tilde_lambda_3}
			|\widetilde\lambda(x,t)-\widetilde\lambda(x,s)|\leq \frac{C_\lambda }{\log 2}\, 	\widetilde\lambda(x,m)\int_{[s,2s]\triangle[t,2t]}\frac{d\rho}{\rho}\leq 2\,  \frac{C_\lambda }{\log 2}\, \log\Bigl(\frac{t}{s}\Bigr)\, 	\widetilde\lambda(x,m),
		\end{equation}
		where the last inequality follows from the assumptions $s<t$, $s,t\in[m/2,m]$, and a standard computation.
		
		Furthermore, since $s\geq m/2$, we have
		\begin{equation}\label{eq:tilde_lambda_4}
			\log\frac{t}{s}	= \log\Bigl(1+\frac{t-s}{s}\Bigr) 	\leq\frac{t-s}{s}\leq 2\, 	\frac{t-s}{m}.
		\end{equation}
		Gathering \eqref{eq:tilde_lambda_2}, \eqref{eq:tilde_lambda_3}, and \eqref{eq:tilde_lambda_4} finally implies \eqref{eq:regularity_estimate_lambda_tilde} with $\mathbf C\coloneqq 4C_\lambda/\log 2$ and concludes the proof of the lemma.
	\end{proof}
	\vv

	To prove the bounds for suppressed kernels, we also need to introduce a dominating function $\Lambda$ such that $\Lambda (\cdot, r)$ has local Lipschitz bounds at a scale $r>0$. This can be done as in \cite[Proposition 1.3]{HYY12}.
	
	\begin{lemma}\label{lemma:aux_dominating_function_3}
		Let $(X,d,\mu)$ be a geometrically doubling metric measure space, and let $\lambda\colon X\times(0,\infty)\to(0,\infty)$ be a dominating function satisfying properties {\rm (i) - (ii)} above.
		
		For
		\[
		\widetilde\lambda(x,r)\coloneqq	\frac1{\log 2}\int_r^{2r}\lambda(x,\rho)\frac{d\rho}{\rho},	\qquad x\in X,\ r>0,
		\]
		we define
		\[
		\Lambda(x,r)\coloneqq\inf_{z\in X}\widetilde\lambda(z,r+d(x,z)),	\qquad x\in X,\ r>0.
		\]
		
		We have that
		\begin{equation}\label{eq:comparability_lambda_Lambda}
			\Lambda(x,r)\leq C_\lambda\, \lambda(x,r),	\qquad \text{ for every } x\in X,\ r>0.
		\end{equation}
		Furthermore, $\Lambda$ satisfies the following properties:
		\begin{enumerate}
			\item[(a)] For every $x\in X$, the function $\Lambda(x,\cdot )$ is non-decreasing. Moreover, $\Lambda(x,r)>0$ and $\mu(B(x,r))\leq \Lambda(x,r)$ for every $x\in X$, $r>0$.		
			\item[(b)] For every $x\in X$ and $r>0$, we have
			\begin{equation}\label{eq:doubling_aux_Lambda}
				\Lambda(x,2r)\leq C^2_\lambda\, \Lambda(x,r).
			\end{equation}
			\item[(c)]  For every $x\in X$, every $m>0$, and every $s,t\in[m/2,m]$, we have
			\begin{equation}\label{eq:regularity_estimate_Lambda_r}
				|\Lambda (x,t)- \Lambda (x,s)|\leq 	2\mathbf C\,C_\lambda\,  \frac{|t-s|}{m} \, \Lambda(x,m),
			\end{equation}
			where $\mathbf C>0$ is as in Lemma \ref{lemma:aux_lambda_2}.		
			\item [(d)] For every $r>0$ and every $x,y \in X$ with $d(x,y)\leq r$, we have
			\[
			\Lambda(x,r)\leq C_\lambda \Lambda (y,r).
			\]
			\item[(e)] If $r>0$ and $x,x'\in X$ are such that \(d(x,x')\leq r\), then
			\begin{equation}\label{eq:bound_in_x_aux_function_Lambda}
				|\Lambda(x,r) - \Lambda(x',r)|\leq \mathbf C\,C_\lambda\,  \frac{d(x,x')}{r}\,\Lambda(x,r).
			\end{equation}
		\end{enumerate}
	\end{lemma}
	{
		\begin{proof}
			Let $x\in X$ and $r>0$. For every $z\in X$, we have that $B(x,r)\subset B(z,r+d(x,z)).$	
			Therefore,
			\[
			\mu(B(x,r))	\leq\lambda(z,r+d(x,z))\overset{\eqref{eq:comparability_lambda_lambda_tilde}}{\leq}\widetilde\lambda(z,r+d(x,z))
			\]
			and, taking the infimum over $z\in X$, we obtain $\mu(B(x,r))\leq \Lambda(x,r)$.
			
			To prove \eqref{eq:comparability_lambda_Lambda}, we choose $z=x$ in the definition of $\Lambda$ and obtain
			\[
			\Lambda(x,r)\leq\widetilde\lambda(x,r)\overset{\eqref{eq:comparability_lambda_lambda_tilde}}{\leq} C_\lambda\, \lambda(x,r).
			\]
			
			We now prove properties (a)-(e).
			For every $x\in X$, the monotonicity of $\Lambda(x,\cdot)$ follows from the monotonicity of $\widetilde\lambda(z,d(x,z) + \cdot )$, for every fixed \(z\in X\).
			
			To prove \eqref{eq:doubling_aux_Lambda}, it suffices to note that, for every $x,z\in X$ and $r>0$, we have
			\[
			\widetilde\lambda\bigl(z,2r+d(x,z)\bigr)\leq\widetilde\lambda\bigl(z,2(r+d(x,z))\bigr)	\leq	C^2_\lambda\,  \widetilde\lambda(z,r+d(x,z)),
			\]
			where in the last inequality we used \eqref{eq:UD2_property} and \eqref{eq:comparability_lambda_lambda_tilde}.
			Hence, we take the infimum over $z\in X$ and obtain \eqref{eq:doubling_aux_Lambda}. Since $\mu\neq0$, every $x$ has a sufficiently large ball of positive measure. The domination in and a finite iteration of (b) therefore imply $\Lambda(x,r)>0$ for every $r>0$.
			
			Let us prove \eqref{eq:regularity_estimate_Lambda_r}. By symmetry, we assume that $s<t$ and, for $z\in X$, we denote
			\[
			\tilde s\coloneqq s+d(x,z),\qquad \tilde t\coloneqq t+d(x,z),\qquad  \widetilde m\coloneqq\tilde t.
			\]
			Then $\tilde s,\tilde t\in[\widetilde m/2,\widetilde m]$,
			$\widetilde m\geq m/2$, and $\tilde t\leq2\tilde s$.
			The preceding lemma and the doubling property of $\widetilde\lambda$ give
			\begin{equation}\label{eq:lip_aux_tilde_lambda}
				\begin{split}
					0\leq\widetilde\lambda(z,\tilde t)-\widetilde\lambda(z,\tilde s)
					&\leq\mathbf C\,\frac{t-s}{\widetilde m}\,
					\widetilde\lambda(z,\widetilde m)\leq2\mathbf C C_\lambda\,\frac{t-s}{m}\,
					\widetilde\lambda(z,\tilde s).
				\end{split}
			\end{equation}
			Thus $\widetilde\lambda(z,\tilde t)\leq
			(1+2\mathbf C C_\lambda(t-s)/m)\widetilde\lambda(z,\tilde s)$.
			Taking the infimum over $z$ and using monotonicity yields
			\[
			0\leq\Lambda(x,t)-\Lambda(x,s)
			\leq2\mathbf C C_\lambda\,\frac{t-s}{m}\,\Lambda(x,s)
			\leq2\mathbf C C_\lambda\,\frac{t-s}{m}\,\Lambda(x,m),
			\]
			as required.
			
			To prove (d), if $d(x,y)\leq r$, then $r+d(x,z)\leq2(r+d(y,z))$. The doubling property of $\widetilde\lambda$ gives $\widetilde\lambda(z,r+d(x,z))\leq C_\lambda\widetilde\lambda(z,r+d(y,z))$; taking the infimum over $z$ proves the claim.
			
			We are left with proving property (e). Let $x,x'\in X$ and assume	$d(x,x')\leq r.$ For \(z\in X\), we define
			\[
			s \coloneqq r+d(x,z)\qquad \text{ and  }\qquad	t \coloneqq r+d(x',z),
			\]
			so that, by the triangle inequality,
			\[
			|s-t|\leq d(x,x')\leq r\leq \min\{s,t\}.
			\]
			
			Moreover, if we denote $m\coloneqq \max\{s,t\}$, then $s,t\in[m/2,m]$, so \eqref{eq:regularity_estimate_lambda_tilde} and the doubling property of $\widetilde \lambda$ imply
			\begin{equation*}
				\begin{split}
					|\widetilde\lambda(z,s)-\widetilde\lambda(z,t)|	&\leq \mathbf C\, \frac{|s-t|}{m}\, \widetilde\lambda(z,m)\leq \mathbf C\, \frac{d(x,x')}{r}\, \widetilde\lambda(z,m)\leq \mathbf C\,C_\lambda\, \frac{d(x,x')}{r}\, \widetilde\lambda(z,t),
				\end{split}
			\end{equation*}
			so we have
			\[
			\widetilde\lambda(z,s)\leq \Bigl(1+\mathbf C\,C_\lambda\, \frac{d(x,x')}{r}\Bigr)\widetilde\lambda(z,t).
			\]
			
			Hence, we finally get
			\[
			\Lambda(x,r)\leq \Bigl(1+\mathbf C\,C_\lambda\, \frac{d(x,x')}{r}\Bigr) \Lambda (x',r)\quad \text{ and }\quad 	\Lambda(x',r)\leq \Bigl(1+\mathbf C\,C_\lambda\, \frac{d(x,x')}{r}\Bigr) \Lambda (x,r),
			\]
			which proves \eqref{eq:bound_in_x_aux_function_Lambda}.	
		\end{proof}
	}
	
	\vv
	
	\subsection{CZ-type bounds for suppressed kernels}
	In this subsection, we define suppressions of Calder\'on-Zygmund kernels in a metric measure space $(X,d,\mu)$ with a dominating function $\lambda$. As in Remark \ref{remark:upper_doubling}, we say that $K\colon(X\times X)\setminus\{(x,x):x\in X\}\to\mathbb C$ is a \emph{Calder\'on-Zygmund kernel} if there exist $\gamma\in(0,1]$ and $C>0$ such that
	\begin{equation}\label{eq:size}
		|K(x,y)|\leq\frac{C}{\lambda(x,d(x,y))},\qquad x,y\in X,\quad x\neq y,
	\end{equation}
	and, for all $x,x',y\in X$ with $m\coloneqq\max\{d(y,x),d(y,x')\}$ and $2d(x,x')\leq m$,
	\begin{equation}\label{eq:holder_alpha}
		|K(x,y)-K(x',y)|+|K(y,x)-K(y,x')|\leq\frac{C}{\lambda(y,m)}\Bigl(\frac{d(x,x')}{m}\Bigr)^\gamma.
	\end{equation}
	
	By Lemma \ref{lemma:aux_dominating_function_3}-(a) and (b), $\Lambda$ is a dominating function and, by \eqref{eq:comparability_lambda_Lambda}, we have $1/\lambda\leq C_\lambda/\Lambda$, so \eqref{eq:size} and \eqref{eq:holder_alpha} also hold with $\Lambda$ in place of $\lambda$, after increasing $C$. Since, in general, only $\Lambda\leq C_\lambda\lambda$ holds, growth conditions with respect to $\lambda$ for measures other than $\mu$ do not transfer to $\Lambda$, and such conditions, such as \eqref{eq:B-large-scale-growth} and \eqref{eq:B-hypotheses-lemma-312} below, are understood with respect to $\Lambda$. From now on, we denote $\Lambda$ by $\lambda$, and by $C_\lambda$ the largest constant in Lemma \ref{lemma:aux_dominating_function_3}, so that the following properties hold:
	\begin{enumerate}
		\item[(i)] The function $\lambda (x, \cdot)$ is non-decreasing and $\mu(B(x,r))\le \lambda(x,r)$, for all $x\in X$, $r>0$.
		\item[(ii)] We have that
		\begin{equation}\label{eq:UD2_property_new}
			\lambda(x,2r)\leq C_\lambda\,\lambda(x,r), \qquad \text{ for every }x\in X, r>0.
		\end{equation}
		\item[(iii)]  For every $x\in X$, every $m>0$, and every $s,t\in[m/2,m]$, we have
		\begin{equation}\label{eq:regularity_estimate_Lambda_r_new}
			|\lambda (x,t)- \lambda (x,s)|\leq 	C_\lambda\,  \frac{|t-s|}{m} \, \lambda(x,m).
		\end{equation}
		\item [(iv)] For every $r>0$ and every $x,y \in X$ with $d(x,y)\leq r$, we have
		\begin{equation}\label{eq:comparability_lambda_new}
			\lambda(x,r)\leq C_\lambda \,  \lambda (y,r).
		\end{equation}
		\item[(v)] If $r>0$ and $x,x'\in X$ are such that $d(x,x')\leq r$, then
		\begin{equation}\label{eq:bound_in_x_aux_function_Lambda_new}
			|\lambda(x,r) - \lambda(x',r)|\leq C_\lambda  \, \frac{d(x,x')}{r}\,\lambda(x,r).
		\end{equation}
	\end{enumerate}
	Let $d_\lambda\geq0$ be such that $\lambda(x,2r)\leq2^{d_\lambda}\lambda(x,r)$ for every $x\in X$ and $r>0$. By \eqref{eq:UD2_property_new}, we can take $d_\lambda=\log_2C_\lambda$. Iterating this inequality, we obtain
	\[
	\lambda(z,t)\leq2^{d_\lambda}(t/s)^{d_\lambda}\lambda(z,s),\qquad z\in X,\quad 0<s\leq t.
	\]
	We fix $M\geq\max\{1,d_\lambda/2\}$, and we remark that the proofs below use only these two lower bounds. In Section \ref{sec:aux_operator}, where $\lambda(x,r)=r^n$ and we can take $d_\lambda=n$, the suppressed kernel \eqref{eq:suppressed_kernel} corresponds to $M=n$, which is admissible since $n\geq1$.

	\begin{definition}[Product-suppressed kernel]\label{def:B-product-suppression}
		Let $\Xi\colon X\to[0,\infty)$ be $1$-Lipschitz, that is,
		\begin{equation}\label{eq:B-Xi-Lipschitz}
			|\Xi(x)-\Xi(x')|\leq d(x,x'),\qquad x,x'\in X.
		\end{equation}
		For $x\neq y$, write $r=d(x,y)$ and set
		\begin{equation}\label{eq:B-product-factor}
			A_\Xi(x,y)\coloneqq\frac{r^{2M}}{r^{2M}+\Xi(x)^M\Xi(y)^M}.
		\end{equation}
		The product-suppressed kernel associated with $K$ is
		\begin{equation}\label{eq:B-product-kernel}
			K_\Xi(x,y)\coloneqq K(x,y)A_\Xi(x,y),	\qquad x,y\in X,\quad x\neq y.
		\end{equation}
	\end{definition}
	
	The factor in \eqref{eq:B-product-factor} is symmetric and takes values in $(0,1]$. Moreover, if $\Xi(x)=0$ or $\Xi(y)=0$, then $A_\Xi(x,y)=1$.
	In particular, $K_\Xi(x,y)=K(x,y)$ whenever $\Xi(x)=0$ or $\Xi(y)=0$.
	
	\begin{lemma}[Upper-doubling product suppression]\label{lem:B-product-suppression}
		Let $K$ satisfy \eqref{eq:size} and \eqref{eq:holder_alpha}, and let $\Xi$ satisfy \eqref{eq:B-Xi-Lipschitz}.  Then $K_\Xi$ is an	upper-doubling Calder\'on-Zygmund kernel with the same exponent
		$\gamma$.  More precisely, the following assertions hold, with constants depending only on $C$, $C_\lambda$, $M$, and $\gamma$.
		
		\smallskip
		\noindent\textup{(i)} For every $x\neq y$,
		\begin{equation}\label{eq:B-usual-size}
			|K_\Xi(x,y)|\leq \frac{C}{\lambda(x,d(x,y))}.
		\end{equation}
		
		\smallskip
		\noindent\textup{(ii)} With the convention	$1/\lambda(z,0)\coloneqq +\infty$, we have
		\begin{equation}\label{eq:B-suppressed-size}
			|K_\Xi(x,y)|\lesssim
			\min\Bigl\{
			\frac1{\lambda(x,\Xi(x))},
			\frac1{\lambda(y,\Xi(y))}
			\Bigr\}.
		\end{equation}
		
		\smallskip
		\noindent\textup{(iii)} If $x,x',y\in X$ satisfy $2d(x,x')\leq m\coloneqq\max\{d(y,x),d(y,x')\}$, then
		\begin{equation}\label{eq:B-product-regularity}
			\begin{split}
				&|K_\Xi(x,y)-K_\Xi(x',y)| +|K_\Xi(y,x)-K_\Xi(y,x')|\lesssim	\frac{1}{\lambda(y,m)}	\Bigl(\frac{d(x,x')}{m}\Bigr)^\gamma .
			\end{split}
		\end{equation}
		
		\smallskip
		\noindent\textup{(iv)} The identity
		\begin{equation}\label{eq:B-kernel-one-zero}
			K_\Xi(x,y)=K(x,y)
		\end{equation}
		holds whenever $\Xi(x)=0$ or $\Xi(y)=0$.
	\end{lemma}
	\begin{proof}
		Let $x\neq y$ and $r\coloneqq d(x,y)$.
		Since $0<A_\Xi\leq1$, \eqref{eq:B-usual-size} follows immediately from \eqref{eq:size}.  We prove the estimate in	\eqref{eq:B-suppressed-size} involving $x$, since the estimate involving $y$ follows by the same argument.
		
		We set $a\coloneqq\Xi(x)$ and $b\coloneqq\Xi(y)$.  Without loss of generality, we can assume that $a>0$.
		If
		$r\geq a/4$, then monotonicity and two applications of
		\eqref{eq:UD2_property_new} give
		\(
		\lambda(x,a)\leq C_\lambda^2\, \lambda(x,r),
		\)
		and hence \eqref{eq:B-usual-size} yields
		$|K_\Xi(x,y)|\leq CC_\lambda^2/\lambda(x,a)$.  Suppose next that
		$r<a/4$.  By \eqref{eq:B-Xi-Lipschitz},
		\(
		b\geq a-r\geq a/2,
		\)
		and therefore
		\[
		A_\Xi(x,y)\leq 2^M\Bigl(\frac ra\Bigr)^{2M}.
		\]
		Hence, by \eqref{eq:size},
		\[
		|K_\Xi(x,y)|\leq\frac{C\,2^M}{\lambda(x,r)}\, \Bigl(\frac ra\Bigr)^{2M}\leq\frac{C\,2^{M+d_\lambda}}{\lambda(x,a)}\, \Bigl(\frac ra\Bigr)^{2M-d_\lambda}\leq\frac{C\,2^{M+d_\lambda}}{\lambda(x,a)},
		\]
		where the last inequality follows from $2M\geq d_\lambda$ and $r<a$. This
		proves \eqref{eq:B-suppressed-size}.
		\vv
		
		To prove \eqref{eq:B-product-regularity}, we fix $x,x',y\in X$ with $2d(x,x')\leq m$, and we write $\delta\coloneqq d(x,x')$, $r\coloneqq d(x,y)$, $r'\coloneqq d(x',y)$, $a\coloneqq\Xi(x)$, $a'\coloneqq\Xi(x')$, and $b\coloneqq\Xi(y)$. Then $r,r'\in[m/2,m]$ and, by \eqref{eq:B-Xi-Lipschitz}, $|a-a'|\leq\delta$, $|a-b|\leq r$, and $|a'-b|\leq r'$. The function $A_\Xi$ is symmetric, thus we obtain
		\[
		\begin{split}
			&|K_\Xi(x,y)-K_\Xi(x',y)|+|K_\Xi(y,x)-K_\Xi(y,x')|\\
			&\qquad\qquad\leq\bigl(|K(x,y)-K(x',y)|+|K(y,x)-K(y,x')|\bigr)A_\Xi(x,y)\\
			&\qquad \qquad\qquad\quad+\bigl(|K(x',y)|+|K(y,x')|\bigr)\,|A_\Xi(x,y)-A_\Xi(x',y)|
			\eqqcolon\mathfrak I_1+\mathfrak I_2.
		\end{split}
		\]
		Note that $A_\Xi\leq1$, thus \eqref{eq:holder_alpha} yields $\mathfrak I_1\leq C\lambda(y,m)^{-1}(\delta/m)^\gamma$.
		
		To estimate $\mathfrak I_2$, we first note that, since $r'\in[m/2,m]$, \eqref{eq:size}, \eqref{eq:comparability_lambda_new}, and \eqref{eq:UD2_property_new} give $|K(x',y)|+|K(y,x')|\lesssim\lambda(y,m)^{-1}$. Moreover, $A_\Xi(x,y)=F(r,a,b)$ and $A_\Xi(x',y)=F(r',a',b)$, where $F(\rho,s,t)\coloneqq\rho^{2M}/(\rho^{2M}+s^Mt^M)$, and a direct computation shows that
		\[
		|\partial_\rho F(\rho,s,t)|+|\partial_sF(\rho,s,t)|\lesssim_M\frac1\rho\qquad\text{if }\rho>0\text{ and }|s-t|\leq\rho.
		\]
		
		For $\partial_sF$, we may assume that $\rho=1$, by homogeneity, and we use that $t\geq s/2$ if $s>4$. Since the set $\{(\rho,s):\rho\geq m/2,\ |s-b|\leq\rho\}$ is convex and contains $(r,a)$ and $(r',a')$, the mean value theorem implies
		\[
		|A_\Xi(x,y)-A_\Xi(x',y)|\lesssim_M\frac{|r-r'|+|a-a'|}{m}\leq 2\, \frac{\delta}{m}.
		\]
		
		Hence, since $\delta\leq m$ and $\gamma\leq1$, we have $\mathfrak I_2\lesssim\lambda(y,m)^{-1}(\delta/m)^\gamma$, which, together with the bound for $\mathfrak I_1$, proves \eqref{eq:B-product-regularity}.
	\end{proof}
	
	\vv
	
	For a finite complex Radon measure $\nu$ and $t>0$, set
	\begin{align}
		T_t\nu(x)&\coloneqq\int_{d(x,y)>t}K(x,y)\,d\nu(y),
		& T_{\Xi,t}\nu(x)&\coloneqq\int_{d(x,y)>t}K_\Xi(x,y)\,d\nu(y),
		\label{eq:B-direct-truncations}\\
		T_t^*\nu(x)&\coloneqq\int_{d(x,y)>t}K(y,x)\,d\nu(y),
		& T_{\Xi,t}^*\nu(x)&\coloneqq\int_{d(x,y)>t}K_\Xi(y,x)\,d\nu(y),
		\label{eq:B-adjoint-truncations}
	\end{align}
	and define the truncated upper-doubling maximal function by
	\begin{equation}\label{eq:B-lambda-maximal}
		\mathcal M_{\lambda,t}\nu(x)\coloneqq \sup_{s\geq t}\frac{|\nu|(B(x,s))}{\lambda(x,s)}.
	\end{equation}
	
	\begin{lemma}\label{lem:B-hard-comparison}
		Let $\nu$ be a finite complex Radon measure, let $t>0$, and let $x\in X$ be such that $\Xi(x)\leq t$.  Then
		\begin{equation}\label{eq:B-hard-comparison}
			\begin{split}
				&|T_{\Xi,t}\nu(x)-T_t\nu(x)|+|T_{\Xi,t}^*\nu(x)-T_t^*\nu(x)|\lesssim\mathcal M_{\lambda,t}\nu(x).
			\end{split}
		\end{equation}
	\end{lemma}
	\begin{proof}
		We write $a\coloneqq \Xi(x)$, and we fix $y$ with $r\coloneqq d(x,y)>t$.  By \eqref{eq:B-Xi-Lipschitz}, $\Xi(y)\leq a+r\leq2r$, and hence
		\begin{equation}\label{eq:B-factor-error}
			0\leq1-A_\Xi(x,y)
			\leq\frac{a^M\Xi(y)^M}{r^{2M}}
			\leq2^M\left(\frac ar\right)^M.
		\end{equation}
		On the annulus $2^kt<d(x,y)\leq2^{k+1}t$, the contribution to the first difference in \eqref{eq:B-hard-comparison} is bounded by a constant times
		\[
		2^{-kM}\, \frac{|\nu|(B(x,2^{k+2}t))}{\lambda(x,2^kt)}\lesssim C_\lambda^2\,  2^{-kM}\mathcal M_{\lambda,t}\nu(x).
		\]
		Summation over $k\geq0$ proves the direct estimate.  The adjoint estimate follows in the same way because \eqref{eq:comparability_lambda_new} gives
		$\lambda(y,d(x,y))\approx_{C_\lambda}\lambda(x,d(x,y))$.
	\end{proof}
	\vv
	
	{
		We define
		\begin{equation}\label{eq:B-hard-at-suppression-scale}
			H_{\Xi,\varepsilon}\nu(x)\coloneqq T_{\max\{\varepsilon,\Xi(x)\}}\nu(x)
			\qquad\text{ and }\qquad 
			H_{\Xi,\varepsilon}^*\nu(x)	\coloneqq T_{\max\{\varepsilon,\Xi(x)\}}^*\nu(x).
		\end{equation}
		Lemmas \ref{lem:B-product-suppression} and \ref{lem:B-hard-comparison} imply, for every $\varepsilon>0$,
		\begin{equation}\label{eq:B-smooth-hard-comparison}
			\begin{split}
				&|T_{\Xi,\varepsilon}\nu(x)-H_{\Xi,\varepsilon}\nu(x)|
				+|T_{\Xi,\varepsilon}^*\nu(x)-H_{\Xi,\varepsilon}^*\nu(x)|\lesssim
				\mathcal M_{\lambda,\max\{\varepsilon,\Xi(x)\}}\nu(x).
			\end{split}
		\end{equation}
		Indeed, if $\varepsilon\geq\Xi(x)$, this is	Lemma \ref{lem:B-hard-comparison}.  If $\varepsilon<\Xi(x)$, the part of $T_{\Xi,\varepsilon}\nu(x)$ over $\varepsilon<d(x,y)\leq\Xi(x)$ is bounded, by
		\eqref{eq:B-suppressed-size}, by a constant times
		\[
		\frac{|\nu|(B(x,2\Xi(x)))}{\lambda(x,\Xi(x))}
		\leq C_\lambda\mathcal M_{\lambda,\Xi(x)}\nu(x),
		\]
		and the remaining part is treated by Lemma \ref{lem:B-hard-comparison} with $t=\Xi(x)$.  The estimate for the adjoint is identical.
	}
	
	For a non-negative Borel function $\Xi$, we define the suppressed maximal function
	\begin{equation}\label{eq:B-suppressed-lambda-maximal}
		\mathcal M^\lambda_\Xi\nu(x)\coloneqq\sup_{s\geq\Xi(x),\ s>0}
		\frac{|\nu|(B(x,s))}{\lambda(x,s)},
	\end{equation}
	and prove its boundedness properties in the following lemma.
	\begin{lemma}
		\label{lem:B-maximal-mapping}
		Let $\Xi\colon X\to[0,\infty)$ be Borel measurable, let $m_1,m_2$ be non-negative Radon measures, and suppose that, for some $c_1,c_2>0$,
		\begin{equation}\label{eq:B-large-scale-growth}
			m_j(B(x,r))\leq c_j\,\lambda(x,r), \qquad \text{ for }x\in X, r>0,\quad r\geq\Xi(x),\quad \text{ and } j=1,2.
		\end{equation}
		
		Then $\mathcal M^\lambda_\Xi$ maps $L^1(m_2)$ to $L^{1,\infty}(m_1)$ and maps $L^p(m_2)$ to $L^p(m_1)$ for $1<p<\infty$, with constants depending only on
		$p$, the geometric-doubling constant, $C_\lambda$, $c_1$, and $c_2$.
	\end{lemma}
	\begin{proof}
		If $\mathcal M^\lambda_\Xi(fm_2)(x)>t$, we choose $r_x>0$ with $r_x\geq\Xi(x)$ such that
		\[
		\int_{B(x,r_x)}|f|\,dm_2>t\lambda(x,r_x).
		\]
		
		Since we obtain estimates uniform on the radii, without loss of generality we can consider them bounded. Then, we apply the $5r$-covering lemma, and obtain a pairwise disjoint subfamily $\{B_j=B(x_j,r_j)\}_j$ such that $\{x\in X: \mathcal M^\lambda_\Xi(fm_2)(x)>t\}$ is covered by $\bigcup_j5B_j$.  By \eqref{eq:B-large-scale-growth} and  \eqref{eq:UD2_property_new}, we get
		\[
		\begin{split}
			m_1\bigl(\{\mathcal M^\lambda_\Xi(fm_2)>t\}\bigr)
			&\leq\sum_jm_1(5B_j)\lesssim c_1\sum_j\lambda(x_j,r_j)\\
			&\leq\frac{c_1}{t}\, \sum_j\int_{B_j}|f|\,dm_2	\leq\frac{c_1}{t}\|f\|_{L^1(m_2)},
		\end{split}
		\]
		which proves that $\mathcal M^\lambda_\Xi$ is bounded from $L^1(m_2)$ to $L^{1,\infty}(m_1)$.
		
		Moreover, \eqref{eq:B-large-scale-growth} for $m_2$ gives the endpoint estimate	$\|\mathcal M^\lambda_\Xi(fm_2)\|_{L^\infty(m_1)}\leq
		c_2\|f\|_{L^\infty(m_2)}$, so Marcinkiewicz interpolation completes the proof.
	\end{proof}
	\vv
	
	We conclude with the analogue of Lemma \ref{lemma:lemm_bd_suppressed_1}. Let $\nu$ be a finite complex Radon measure, and let $x\in X$, $R_0\geq0$, $a_0>0$, and $c_0>0$ be such that
	\begin{equation}\label{eq:B-hypotheses-lemma-312}
		|\nu|(B(x,r))\leq c_0\lambda(x,r)
		\qquad\text{ and }\qquad 
		|T_t\nu(x)|\leq a_0	,
	\end{equation}
	for $t,r>0$, $t,r\geq R_0$.
	
	If $\Xi(x)\geq R_0$, then, for $\varepsilon>0$ and $t\coloneqq\max\{\varepsilon,\Xi(x)\}\geq R_0$, \eqref{eq:B-smooth-hard-comparison} and \eqref{eq:B-hypotheses-lemma-312} give
	\begin{equation}\label{eq:B-upper-lemma-312}
		|T_{\Xi,\varepsilon}\nu(x)|\leq|T_t\nu(x)|+C\mathcal M_{\lambda,t}\nu(x)\leq a_0+Cc_0,
	\end{equation}
	and the same bound holds for $T^*_{\Xi,\varepsilon}\nu(x)$ if $|T^*_t\nu(x)|\leq a_0$ for $t>0$ with $t\geq R_0$.
	
	In summary, Lemma \ref{lemma:kernel_suppressed} is the case $\lambda(x,r)=r^n$ and $M=n$ of Lemmas \ref{lem:B-product-suppression} and \ref{lem:B-hard-comparison}, Lemma \ref{lem:B-maximal-mapping} is the analogue of Lemma \ref{lem:suppressed maximal}, and \eqref{eq:B-upper-lemma-312} is the analogue of Lemma \ref{lemma:lemm_bd_suppressed_1}.
	
	On the other hand, we remark that the mixed Cotlar inequality of Lemma \ref{lemma:mixed_cotlar_transfer}, Sections \ref{sec:corona}-\ref{sec:main_theorem_antisym}, and Theorem \ref{theorem:main_theorem_hofmann_problem} are stated and proved only under polynomial-growth assumptions.
	
	\section{Proof of Theorem \ref{theorem:matrix_Tb}}\label{appendix:matrix_Tb}
	
	We follow the proof of Theorem \ref{theorem:main_theorem_hofmann_problem} in Sections \ref{section:exceptional_set}-\ref{sec:main_theorem_antisym}, including the reduction in Lemma \ref{lemma:reduction_cutoffs}, with $V$-valued functions and with the norm of $V$ and the operator norm in place of absolute values.
	For the sake of brevity, we further assume $p,q<\infty$; otherwise, the proof can be obtained by minor modifications.
	 The bilinear pairings $\langle f,g\rangle_{\mathfrak m}$ are replaced by $\int\langle f,g\rangle_V\,d\mathfrak m$, and the transposed kernels, such as $K(y,x)$ and $L_\varepsilon(y,x)$ in Sections \ref{sec:corona}-\ref{section:suppressed_mixed_good}, by their adjoints $K(y,x)^*$ and $L_\varepsilon(y,x)^*$. The constants may also depend on $\dim V$, $p$, $q$, and $C_b$. We list the changes, section by section.
	
	\begin{itemize}
		\item \textit{Section \ref{subsection:bounded_testing_reduction}.} Let $B\in\mathscr B$ and $\mathfrak a_B\coloneqq\avint_B\mathfrak b_B\,d\mu$. By \eqref{eq:matrix_Tb_accretive} and \eqref{eq:matrix_Tb_size}, $\langle\mathfrak a_Bv,v\rangle_V\geq|v|^2$, so that $\|\mathfrak a_B^{-1}\|\leq1$, and $\|\mathfrak a_B\|\leq C_b$. We consider the $\mathcal L(V)$-valued measure $\widehat\nu_B\coloneqq\mathfrak b_B\mathfrak a_B^{-1}\mu$, its total variation $|\widehat\nu_B|\coloneqq\|\mathfrak b_B\mathfrak a_B^{-1}\|\,\mu$, and its polar decomposition $\widehat\nu_B=\widehat b_B\,|\widehat\nu_B|$, where $\widehat b_B\coloneqq\mathfrak b_B\mathfrak a_B^{-1}/\|\mathfrak b_B\mathfrak a_B^{-1}\|$ is well defined on $B$, since \eqref{eq:matrix_Tb_accretive} gives 
		\[
			|\mathfrak b_B(x)\mathfrak a_B^{-1}v|\geq|\mathfrak a_B^{-1}v|\geq C_b^{-1}|v|, \qquad \text{ for }v\in V\text{ and }\mu\text{-a.e. }x\in B.
		\]
		Then $\widehat\nu_B(B)=\mu(B)I$ and, by H\"older's inequality,
		\begin{equation}\label{eq:matrix_normalized}
			\mu(B)\leq|\widehat\nu_B|(B)\leq C_b\,\mu(B)\qquad\text{ and }\qquad |\widehat\nu_B|(F)\leq C_b\,\mu(B)^{1/p}\mu(F)^{1/p'},
		\end{equation}
		for every Borel set $F\subset B$. The same construction for $\mathfrak b^*_B$, with $\mathfrak a^*_B\coloneqq\avint_B\mathfrak b^*_B\,d\mu$, gives a measure $\widehat\nu^*_B\coloneqq\mathfrak b^*_B(\mathfrak a^*_B)^{-1}\mu$ which satisfies \eqref{eq:matrix_normalized} with $q$ in place of $p$. Thus, \eqref{eq:measure_testing_size}, with $\widehat\nu_B(B)=\mu(B)I$ and $\widehat\nu^*_B(B)=\mu(B)I$, and \eqref{eq:measure_testing_small_sets} hold for both systems, for $\mathfrak C_1\geq2C_b$ such that $C_b^s(2\mathfrak C_1/C_b)^{1-s}\leq1/2$ for $s\in\{p,q\}$, and $\varepsilon_0\coloneqq(200\,C_b\mathfrak C_1)^{-\max\{p',q'\}}$. Moreover, $T_{\mu,\varepsilon}(\mathfrak b_B\mathfrak a_B^{-1})=(T_{\mu,\varepsilon}\mathfrak b_B)\mathfrak a_B^{-1}$, so that $T_\sharp\widehat\nu_B\leq T_{\mu,\sharp}\mathfrak b_B$, and, for $0<\tau<\tau_0$ and $0<\zeta<\zeta_0$, \eqref{eq:matrix_Tb_testing} and Chebyshev's inequality give \eqref{eq:measure_testing_L0}, for every $\eta$, with $\mathfrak C_2\coloneqq(2C_\tau/\zeta)^{1/\min\{p,q\}}$. The same holds for $\mathfrak b^*_B$ and $T^*$, with $q$ in place of $p$.
		
		\item \textit{Sections \ref{section:exceptional_set}, \ref{sec:corona}, and \ref{section:L2bd}.} We set $\nu_B\coloneqq\widehat\nu_B+\mu|_{X\setminus B}\,I$ and $|\nu_B|\coloneqq|\widehat\nu_B|+\mu|_{X\setminus B}$, and we extend $\widehat b_B$ by $I$ on $X\setminus B$. In the definition of the family $\mathcal A(w)$ of non-accretive cubes, the condition $|\nu(S)|<\eta|\nu|(S)$ is replaced by $|\nu(S)v|<\eta\,|\nu|(S)\,|v|$ for some $v\in V\setminus\{0\}$, and the first inequality in Lemma \ref{lemma:lemma38} reads $|\nu(S^*)v|\geq\eta|\nu|(S^*)|v|$, for $v\in V$, and its proof is unchanged. For Lemma \ref{lemma:estim_T_i}, if $S$ is a cube, then
		\[
		\nu(S)=\Bigl[\int_{S\cap B}\mathfrak b_B\,d\mu+\mu(S\setminus B)\,\mathfrak a_B\Bigr]\mathfrak a_B^{-1},
		\]
		where, by \eqref{eq:matrix_Tb_accretive}, the operator in brackets satisfies $\langle[\,\cdot\,]v,v\rangle_V\geq\mu(S)|v|^2$. Hence, $|\nu(S)v|\geq C_b^{-1}\mu(S)|v|$ and, since $\eta=(2\mathfrak C_1)^{-1}$, every $S\in\mathcal A(w)$ satisfies $|\nu|(S)>L_*\mu(S)$, with $L_*\coloneqq2\mathfrak C_1/C_b\geq4$. If $Z$ is the union of a family of pairwise disjoint cubes of $\mathcal A(w)$, then $|\widehat\nu_B|(Z)>L_*\mu(Z\cap B)$ and, by \eqref{eq:matrix_normalized}, 
		\[
			|\widehat\nu_B|(Z)<C_b^pL_*^{1-p}\, \mu(B)\leq\mu(B)/2\leq(1-\eta)\, |\widehat\nu_B|(B),
		\]
		and the same holds for the system of $T^*$, with $q$ in place of $p$. Hence, Lemma \ref{lemma:estim_T_i} follows by monotone convergence, as in its proof. The other constructions and estimates of these sections, including Lemma \ref{lemma:mixed_cotlar_transfer}, only involve $|\widehat\nu_B|$, $|\nu_B|$, the truncations $T_\varepsilon\widehat\nu_B$, and stopping conditions, and they hold without changes.
		
		\item \textit{Section \ref{section:suppressed_mixed_good}.}
		\begin{itemize}
			\item[$\circ$] \textit{Martingale differences.} In \eqref{eq:top_averages_def}, \eqref{eq:def_martingales}, and in the definitions of $\Delta^*_{i,Q}$, $D^j_{I_j}$, and $c_{R,Q}$, we replace $f\,\widehat b_i$ by $\widehat b_i^*f$, we read the quotients $\langle f\rangle_{i,P}/\langle\widehat b_i\rangle_{i,P}$ and $\langle\widehat b_i^*f\rangle_{i,P}/\langle\widehat b_i\rangle_{i,P}$ as $\langle\widehat b_i\rangle_{i,P}^{-1}\langle f\rangle_{i,P}$ and $(\langle\widehat b_i\rangle_{i,P}^*)^{-1}\langle\widehat b_i^*f\rangle_{i,P}$, and the multiplication by $\widehat b_i$ is on the left. For instance, if $Q\in\bbD_i^{\tr}$ and $P\in\Ch(Q)\cap\bbD_i^{\tr}$, then
			\[
			\Delta_{i,Q}f=\widehat b_i\bigl(A_P^{-1}f_P-A_Q^{-1}f_Q\bigr)\quad\text{ on }P,\qquad\text{ where }\quad A_P\coloneqq\langle\widehat b_i\rangle_{i,P}\quad\text{ and }\quad f_P\coloneqq\langle f\rangle_{i,P},
			\]
			and, by the stopping conditions, $|A_Pv|\geq\eta|v|$ for every transit cube $P$ and $v\in V$. Then $\int\Delta_{i,Q}f\,d|\nu_i|=0$ and $\Delta_{i,Q}(\widehat b_iv)=0$, for $v\in V$, which give \eqref{eq:property_id_mart_diff}, and $\Delta^*_{i,Q}$ is the adjoint of $\Delta_{i,Q}$ with respect to $\langle\cdot,\cdot\rangle_{|\nu_i|}$.
			
			\item[$\circ$] \textit{Martingale decompositions.} We write 
			\[
				A_P^{-1}f_P-A_Q^{-1}f_Q=A_P^{-1}(f_P-f_Q)+A_P^{-1}(A_Q-A_P)A_Q^{-1}f_Q,
			\]
			and similarly for the adjoints. The first term is controlled by the martingale differences of $f$ with respect to $|\nu_i|$. For the second one, let $\alpha_Q\coloneqq\sum_P|\nu_i|(P)\,\|A_P-A_Q\|_{\mathrm{HS}}^2$, where $\|\cdot\|_{\mathrm{HS}}$ is the Hilbert-Schmidt norm and $P$ runs over the transit children of $Q$. The orthogonality of the martingale differences of the entries of $\widehat b_i$ gives $$\sum_{Q\in\bbD_i^{\tr},\,Q\subset R}\alpha_Q\leq(\dim V)\,|\nu_i|(R),\qquad  \text{ for }R\in\bbD(w_i),$$ so that Carleson's embedding theorem gives $\sum_Q\alpha_Q|f_Q|^2\lesssim\|f\|^2_{L^2(|\nu_i|)}$. On a terminal child $P$ of $Q$, we have $\Delta_{i,Q}f=f-\widehat b_iA_Q^{-1}f_Q$, these children are pairwise disjoint, and $|f_Q|$ is bounded by the dyadic maximal function of $|f|$ with respect to $|\nu_i|$. Moreover, the partial sums of \eqref{eq:dec_f_mart} over the first $k$ generations below $Q_2(w_i)$ equal $\widehat b_iA_P^{-1}f_P$ on the transit cubes $P$ of the next generation, and $f$ on the rest of $Q_2(w_i)$, and they converge to $f$, with domination by $\eta^{-1}$ times this maximal function, since $A_P\to\widehat b_i(x)$ and $f_P\to f(x)$ along the transit cubes $P\ni x$, for $|\nu_i|$-a.e. $x$. This gives \eqref{eq:dec_f_mart}, \eqref{eq:dec_f_mart_adj}, and the upper bounds in \eqref{eq:dec_f_mart_2} and in its adjoint version, and the lower bounds and the unconditional convergence follow by duality, since $\Delta_{i,Q}^2=\Delta_{i,Q}$.
			
			\item[$\circ$] \textit{Coefficients.} The identities that use $|\widehat b_i|=1$ are replaced by the following consequence of Jensen's inequality: for every transit cube $P$ and $v\in V$,
			\begin{equation}\label{eq:matrix_coefficients}
				\eta\,|v|\,|\nu_i|(P)^{1/2}\leq\|\chi_P\widehat b_iv\|_{L^2(|\nu_i|)}\leq|v|\,|\nu_i|(P)^{1/2}.
			\end{equation}
			Hence, where now $D^1_If_{\mathtt w},c_{R,Q}(g_{\mathtt w})\in V$, we have $$|D^1_If_{\mathtt w}|^2|\nu_1|(I)\leq\eta^{-2}\|\chi_I\Delta_{1,I^*}f_{\mathtt w}\|^2_{L^2(|\nu_1|)}$$ and $$|c_{R,Q}(g_{\mathtt w})|^2|\nu_2|(P)\leq\eta^{-2}\|\chi_P\Delta_{2,R}g_{\mathtt w}\|^2_{L^2(|\nu_2|)},$$ for $P\in\{R_Q,S\}$, in the proof of Lemma \ref{lemma:S_delta_1}, in \eqref{eq:property_coef_c_Q_R}, and in the proof of Lemma \ref{lemma:s_41_tr_s_2}.
			
			\item[$\circ$] \textit{Weak boundedness.} Lemmas \ref{lemma:main_technical_lemma_wbp}, \ref{lemma:key_supp_WBP_2}, and \ref{lemma:key_supp_WBP_in_text} hold with $\widehat b_1v$ and $\widehat b_2v'$ in place of $\widehat b_1$ and $\widehat b_2$, for $v,v'\in V$, and with an additional factor $|v||v'|$ on the right-hand sides, since their proofs only use $\|\widehat b_i\|=1$ and the bounds \eqref{eq:global_testing_before_localization_sec5} and \eqref{eq:localized_maximal_testing_sec5}, now for the operator norms of $T_{\Xi_{\mathtt w},\varepsilon}\widehat\nu_1$, $T^*_{\Xi_{\mathtt w},\varepsilon}\widehat\nu_2$, $S_{|\nu_1|,\varepsilon,\rho}\widehat b_1$, and $S^*_{|\nu_2|,\varepsilon,\rho}\widehat b_2$. Thus, for a nice ball $B$ in the proof of Lemma \ref{lemma:S_delta_1}, the pairing in \eqref{eq:S_delta_1_2} equals $\langle S_{|\nu_1|,\varepsilon}(\chi_B\widehat b_1D^1_If_{\mathtt w}),\chi_B\widehat b_2D^2_Jg_{\mathtt w}\rangle_*$, and we estimate it by Lemma \ref{lemma:key_supp_WBP_in_text} with $v=D^1_If_{\mathtt w}$ and $v'=D^2_Jg_{\mathtt w}$.
			
			\item[$\circ$] \textit{Paraproduct.} In the proof of Lemma \ref{lemma:s_41_tr_s_1}, by \eqref{eq:identity_coefficients_c_Q_R}, the term $\widehat S_1$ is the sum of $\langle\Delta_{1,Q}f_{\mathtt w},(S^*_{|\nu_2|,\varepsilon}\widehat b_2)\,c_{R,Q}(g_{\mathtt w})\rangle$, the coefficients $c_{R,Q}(g_{\mathtt w})$ telescope as in the scalar case, and we obtain $\widehat S_1=\langle f_{\mathtt w},\Pi_{\Xi_{\mathtt w}}(\widetilde g)\rangle$, with
			\[
			\Pi_{\Xi_{\mathtt w}}(\phi)\coloneqq\sum_{Q\in\mathcal F_\Pi}\bigl(\Delta^*_{1,Q}S^*_{|\nu_2|,\varepsilon}\widehat b_2\bigr)\,\langle\widehat b_2\rangle_{2,R(Q)}^{-1}\langle\phi\rangle_{2,R(Q)},
			\]
			where $\Delta^*_{1,Q}$ acts on the columns of the $\mathcal L(V)$-valued function $S^*_{|\nu_2|,\varepsilon}\widehat b_2$. Then the proof applies with $a_Q\coloneqq\|\Delta^*_{1,Q}S^*_{|\nu_2|,\varepsilon}\widehat b_2\|^2_{L^2(|\nu_1|;\mathrm{HS})}$, up to the factors $\dim V$ and $\eta^{-2}$.
		\end{itemize}
		The other estimates of Section \ref{section:suppressed_mixed_good} only involve the norms of the functions, the cancellation of the martingale differences, and the kernel estimates, and they hold without changes.
		
		\item \textit{Section \ref{sec:main_theorem_antisym}.} The arguments of Section \ref{section:probabilisti_argument_subsection} apply, Lemma \ref{lemma:final_boundedness_suppressed_operators} holds, with the parameters chosen as in its proof, and the first part of the proof of Proposition \ref{proposition:bounded_testing_normalized} gives \eqref{eq:corollary49_maximal_sec8} for $T$ and $T^*$, and for $V$-valued functions $g$. Let $\{v_j\}_j$ be an orthonormal basis of $V$, and let $T_{jk}$ be the SIO associated with the scalar kernel $K_{jk}(x,y)\coloneqq\langle K(x,y)v_k,v_j\rangle_V$. Then $(T_{jk})_{\mu,\sharp}g\leq T_{\mu,\sharp}(g\,v_k)$ and, since $K_{jk}(y,x)=\langle v_k,K(y,x)^*v_j\rangle_V$, the transpose of $T_{jk}$ satisfies $(T^*_{jk})_{\mu,\sharp}g\leq T^*_{\mu,\sharp}(g\,v_j)$. Thus, \eqref{eq:corollary49_maximal_sec8} holds for $T_{jk}$, and the remainder of the proof of Proposition \ref{proposition:bounded_testing_normalized}, which only uses this estimate and the structural assumptions, shows that $T_{jk}$ is bounded on $L^2(\mu)$. Hence, $T_\mu$ is bounded on $L^2(\mu;V)$, and the cut-offs are removed as in Lemma \ref{lemma:reduction_cutoffs}.
	\end{itemize}

	\color{black}


\begin{thebibliography}{AHMMMTV16}
		
		\bibitem[AAAHK11]{AAAHK11}
		M.~A.~Alfonseca, P.~Auscher, A.~Axelsson, S.~Hofmann, and  S.~Kim.
		\newblock{Analyticity of layer potentials and $L^2$ solvability of boundary value problems for divergence form elliptic equations with complex $L^\infty$ coefficients.}
		\newblock{\em Adv. Math.}{, 226 (5): 4533--4606, 2011.}
		
		\bibitem[AA11]{AA11}
		P.~Auscher and A.~Axelsson.
		\newblock{Weighted maximal regularity estimates and solvability of non-smooth elliptic systems I.}
		\newblock{\em Invent. Math.}{, 184 (1): 47--115, 2011.}
		
		\bibitem[AAH08]{AAH08}
		P.~Auscher, A.~Axelsson, and S.~Hofmann.
		\newblock{Functional calculus of Dirac operators and complex perturbations of Neumann and Dirichlet problems.}
		\newblock{\em J. Funct. Anal.}{, 255 (2): 374--448, 2008.}
		
		\bibitem[AAM10]{AAM10}
		P.~Auscher, A.~Axelsson, and A.~McIntosh.
		\newblock{Solvability of elliptic systems with square integrable boundary data.}
		\newblock{\em Ark. Mat.}{, 48 (2): 253--287, 2010.}
		
		\bibitem[AHLMT02]{AHLMT02}
		P.~Auscher, S.~Hofmann, M.~Lacey, A.~McIntosh, and P.~Tchamitchian.
		\newblock{The solution of the Kato square root problem for second order elliptic operators on $\mathbb R^n$.}
		\newblock{\em Ann. of Math. (2)}{, 156 (2): 633--654, 2002.}
		
		\bibitem[AHMTT02]{AHMTT02}
		P.~Auscher, S.~Hofmann, C.~Muscalu, T.~Tao, and C.~Thiele.
		\newblock{Carleson measures, trees, extrapolation, and $T(b)$ theorems.}
		\newblock{\em Publ. Mat.}{, 46: 257--325, 2002.}
		
		\bibitem[AMM13]{AMM13}
		P.~Auscher, A.~McIntosh, and M.~Mourgoglou.
		\newblock{On $L^2$ solvability of BVPs for elliptic systems.}
		\newblock{\em J. Fourier Anal. Appl.}{, 19 (3): 478--494, 2013.}
		
		\bibitem[AR13]{AR13}
		P.~Auscher and E.~Routin.
		\newblock{Local $Tb$ theorems and Hardy inequalities.}
		\newblock{\em J. Geom. Anal.}{, 23: 303--374, 2013.}
		
		\bibitem[AY09]{AY09}
		P.~Auscher and Q.~X.~Yang.
		\newblock{BCR algorithm and the $T(b)$ theorem.}
		\newblock{\em Publ. Mat.}{, 53: 179--196, 2009.}
		
		\bibitem[AKM06]{AKM06}
		A.~Axelsson, S.~Keith, and A.~McIntosh.
		\newblock{Quadratic estimates and functional calculi of perturbed Dirac operators.}
		\newblock{\em Invent. Math.}{, 163 (3): 455--497, 2006.}
		
		\bibitem[AHM$^3$TV16]{AHMMMTV16}
		J.~Azzam, S.~Hofmann, J.~M.~Martell, S.~Mayboroda, M.~Mourgoglou, X.~Tolsa, and A.~Volberg.
		\newblock{Rectifiability of harmonic measure. }
		\newblock{\em Geom. Funct. Anal.}{, 26 (3): 703-728, 2016.}
		
		\bibitem[AHMMT20]{AHMMT20}
		J.~Azzam, S.~Hofmann, J.~M.~Martell, M.~Mourgoglou, and X.~Tolsa.
		\newblock{Harmonic measure and quantitative connectivity: geometric characterization of the $L^p$-solvability of the Dirichlet problem.}
		\newblock{\em Invent. Math.}{, 222 (3): 881--993, 2020.}
		
		\bibitem[AMT17]{AMT17a}
		J.~Azzam, M.~Mourgoglou, and X.~Tolsa.
		\newblock{Mutual absolute continuity of interior and exterior harmonic measure implies rectifiability.}
		\newblock{\em Commun. Pure Appl. Math.}{, 70 (11): 2121-2163, 2017.}
		
		\bibitem[AMT20]{AMT20}
		J.~Azzam, M.~Mourgoglou, and X.~Tolsa.
		\newblock{A two-phase free boundary problem for harmonic measure and uniform rectifiability.}
		\newblock{\em Trans. Am. Math. Soc.}{, 373 (6): 4359-4388, 2020.}
		
		\bibitem[AMTV19]{AMTV19}
		J.~Azzam, M.~Mourgoglou, X.~Tolsa, and A.~Volberg.
		\newblock{On a two-phase problem for harmonic measure in general domains.}
		\newblock{\em Am. J. Math.}{, 141 (5): 1259-1279, 2019.}
		
		
		\bibitem[Be11]{Be11}
		L.~Berkovits.
		\newblock{Parabolic Muckenhoupt weights in the Euclidean space.}
		\newblock{\em J. Math. Anal. Appl.}, {379 (2): 524--537, 2011.}
		
		\bibitem[BKW21]{BKW21}
		K.~Bickel, T.~Korhonen, and B.~Wick.
		\newblock{Two-weight $Tb$ theorems for well-localized operators.}
		\newblock{\em Math. Nachr.}{, 294 (7): 1277--1294, 2021.}
		
		\bibitem[BHHLN22]{BHHLN22}
		S.~Bortz, J.~Hoffman, S.~Hofmann, J.L.~Luna-Garcia, and K.~Nyström.
		\newblock{Coronizations and big pieces in metric spaces. (Construction en couronnes et grands morceaux dans des espaces métriques.)}
		\newblock{\em Ann. Inst. Fourier}{, 72 (5): 2037--2078, 2022.}
		
		\bibitem[BHHLN26]{BHHLN26}
		S.~Bortz, J.~Hoffman, S.~Hofmann, J.L.~Luna-Garcia, and K.~Nyström.
		\newblock{Parabolic singular integrals with nonhomogeneous kernels.}
		\newblock{\em Bull. Lond. Math. Soc.}{, 58, No. 7, Article ID e70425, 16 p., 2026.}
		
		\bibitem[BHMN25a]{BHMN23}
		S.~Bortz, S.~Hofmann, J.M.~Martell, and K.~Nyström.
		\newblock{Solvability of the $L^p$ Dirichlet problem for the heat equation is equivalent to parabolic uniform rectifiability in the case of a parabolic Lipschitz graph.}
		\newblock{\em Invent. Math.}{, 239 (1): 165--217, 2025.}
		
		\bibitem[BHMN25b]{BHMN25}
		S.~Bortz, S.~Hofmann, J.~M.~Martell, and K.~Nyström.
		\newblock{Solvability of the $L^p$ Dirichlet problem for the heat equation implies parabolic uniform rectifiability.}
		\newblock{\em arXiv preprint,}{ 2510.22047, 2025.}
		
		\bibitem[BDS82]{BDS82}
		F.~Brackx, R.~Delanghe, and F.~Sommen.
		\newblock{\em Clifford analysis.}
		\newblock{Research Notes in Mathematics 76, Pitman, Boston, MA, 1982.}
		
		\bibitem[Br95]{Br95}
		R.~M.~Brown.
		\newblock{The Neumann problem on Lipschitz domains in Hardy spaces of order less than one.}
		\newblock{\em Pacific J. Math.}{, 171 (2): 389--407, 1995.}
		
		\bibitem[Ch90]{Ch90}
		M.~Christ.
		\newblock{A $T(b)$ theorem with remarks on analytic capacity and the Cauchy integral.}
		\newblock{\em Colloq. Math.}{, 60/61: 601--628, 1990.}
		
		\bibitem[CMM82]{CMM82}
		R.~R.~Coifman, A.~McIntosh, and Y.~Meyer.
		\newblock{L'int\'egrale de Cauchy d\'efinit un op\'erateur born\'e sur $L^2$ pour les courbes lipschitziennes.}
		\newblock{\em Ann. of Math. (2)}{, 116 (2): 361--387, 1982.}
		
		\bibitem[Co12]{Con12}
		J.~M.~Conde-Alonso.
		\newblock{A note on dyadic coverings and nondoubling Calder\'on-Zygmund theory.}
		\newblock{arXiv:1201.3513, 2012.}
		
		\bibitem[DK87]{DK87}
		B.~E.~J.~Dahlberg and C.~E.~Kenig.
		\newblock{Hardy spaces and the Neumann problem in $L^p$ for Laplace's equation in Lipschitz domains.}
		\newblock{\em Ann. of Math. (2)}{, 125 (3): 437--465, 1987.}
		
		\bibitem[Da98]{Da98}
		G.~David.
		\newblock{Unrectifiable 1-sets have vanishing analytic capacity.}
		\newblock{\em Rev. Mat. Iberoam.}{, 14 (2): 369--479, 1998.}
		
		\bibitem[DJ90]{DJ90}
		G.~David and D.~Jerison.
		\newblock{Lipschitz approximation to hypersurfaces, harmonic measure, and singular integrals.}
		\newblock{\em Indiana Univ. Math. J.}{, 39 (3): 831--845, 1990.}
		
		\bibitem[DJ84]{DJ84}
		G.~David and J.-L.~Journé.
		\newblock{A boundedness criterion for generalized Calder\'on-Zygmund operators.}
		\newblock{\em Ann. of Math.}{, 120: 371--397, 1984.}
		
		\bibitem[DJS85]{DJS85}
		G.~David, J.-L.~Journé, and S.~Semmes.
		\newblock{Opérateurs de Calderón-Zygmund, fonctions para-accrétives et interpolation. (Calderón-Zygmund operators, para-accretive functions and interpolation).}
		\newblock{\em Rev. Mat. Iberoam.}{, 1 (4): 1--56, 1985.}
		
		\bibitem[DM00]{DM00}
		G.~David and P.~Mattila.
		\newblock{Removable sets for Lipschitz harmonic functions in the plane.}
		\newblock{\em Rev. Mat. Iberoamericana}{, 16 (1): 137--215, 2000.}
		
		\bibitem[DS91]{DS91}
		G.~David and S.~Semmes.
		\newblock{Singular integrals and rectifiable sets in $\mathbb R^n$: Au-del\`a des graphes lipschitziens.}
		\newblock{\em Ast\'erisque}{, 193, 1991.}
		
		\bibitem[DHP23]{DHP23}
		M.~Dindo\v{s}, S.~Hofmann, and J.~Pipher.
		\newblock{Regularity and Neumann problems for operators with real coefficients satisfying Carleson conditions.}
		\newblock{\em J. Funct. Anal.}{, 285 (6), Paper No. 110024, 2023.}
		
		\bibitem[DPR17]{DPR17}
		M.~Dindo\v{s}, J.~Pipher, and D.~Rule.
		\newblock{Boundary value problems for second-order elliptic operators satisfying a Carleson condition.}
		\newblock{\em Comm. Pure Appl. Math.}{, 70 (7): 1316--1365, 2017.}
		
		\bibitem[FS83]{FS83}
		E.~Fabes and S.~Salsa.
		\newblock{Estimates of caloric measure and the initial-Dirichlet problem for the heat equation in Lipschitz cylinders.}
		\newblock{\em Trans. Amer. Math. Soc.}, {279 (2): 635--650, 1983.}
		
		\bibitem[FL24]{FL24}
		J.~Feneuil and L.~Li.
		\newblock{The $L^p$ Poisson-Neumann problem and its relation to the Neumann problem.}
		\newblock{\em arXiv preprint,}{ https://arxiv.org/abs/2406.16735, 2024.}
		
		\bibitem[Fe26]{Fe26}
		M.~Fernàndez-Vilaseca.
		\newblock{A $Tb$ type theorem for suppressed kernels.}
		\newblock{\em arXiv preprint,}{ https://arxiv.org/abs/2606.26401, 2026.}
		
		\bibitem[Ga25]{Ga25}
		J.~M.~Gallegos.
		\newblock{One-sided Rellich inequalities, regularity problem and uniform rectifiability.}
		\newblock{\em arXiv preprint,}{ https://arxiv.org/abs/2506.03431, 2025.}
		
		\bibitem[GMT25]{GMT25}
		J.~M.~Gallegos, M.~Mourgoglou, and X.~Tolsa.
		\newblock{Extrapolation of solvability of the regularity and the Poisson regularity problems in rough domains.}
		\newblock{\em J. Funct. Anal.}{, 288, Paper No. 110672, 2025.}
		
		\bibitem[GMT18]{GMT18}
		J.~Garnett, M.~Mourgoglou, and X.~Tolsa.
		\newblock{Uniform rectifiability from Carleson measure estimates and $\epsilon$-approximability of bounded harmonic functions.}
		\newblock{\em Duke Math. J.}{, 167 (8): 1473--1524, 2018.}
		
		\bibitem[GLQ93]{GLQ93}
		G.~I.~Gaudry, R.~Long, and T.~Qian.
		\newblock{A martingale proof of $L^2$ boundedness of Clifford-valued singular integrals.}
		\newblock{\em Ann. Mat. Pura Appl. (4)}{, 165: 369--394, 1993.}
		
		\bibitem[GeH20]{GeH20}
		A.~Genschaw and S.~Hofmann.
		\newblock{A weak reverse Hölder inequality for caloric measure.}
		\newblock{\em J. Geom. Anal.}{, 30 (2): 1530--1564, 2020.}
		
		\bibitem[GH17]{GH17}
		A.~Grau de la Herrán and S.~Hofmann.
		\newblock{Generalized local $Tb$ theorems for square functions.}
		\newblock{\em Mathematika}{, 63 (1): 1--28, 2017.}
		
		
		\bibitem[Gu25]{Gu25}
		I.~Guill\'en-Mola.
		\newblock{$L^p$-solvability of boundary value problems for the Laplacian in locally flat unbounded domains.}
		\newblock{\em arXiv preprint,}{ https://arxiv.org/abs/2503.23865, 2025.}
		
		\bibitem[He24]{He24}
		J.~Hernández.
		\newblock{Removable singularities for Lipschitz fractional caloric functions in time varying domains.}
		\newblock{Preprint, arXiv:2412.18402, 2024.}
		
		\bibitem[HMP24]{HMP24}
		J.~Hernández, J.~Mateu, and L.~Prat.
		\newblock{On the semi-additivity of the $1/2$-symmetric caloric capacity.}
		\newblock{Preprint, arXiv:2405.01195, 2024.}
		
		
		\bibitem[HJ25]{HJ25}
		J.~Hoffman and B.~Jaye.
		\newblock{On Singular Integrals and Quantitative Rectifiability in Parabolic Space and the Heisenberg Group.}
		\newblock{Preprint, arXiv:2510.26934, 2025.}
		
		\bibitem[Ho07]{Ho07}
		S.~Hofmann.
		\newblock{A proof of the local $Tb$ theorem for standard Calder\'on-Zygmund operators.}
		\newblock{\em unpublished manuscript, arXiv:0705.0840, 2007.}
		
		\bibitem[Ho10]{Ho10}
		S.~Hofmann.
		\newblock{Local $T(b)$ theorems and applications in PDE.}
		\newblock{\em Harmonic Analysis and Partial Differential Equations}{, Contemp. Math. 505, 29--52, Amer. Math. Soc., Providence, RI, 2010.}
		
		\bibitem[HLMN17]{HLMN17}
		S.~Hofmann, P.~Le, J.~M.~Martell, and K.~Nyström.
		\newblock{The weak-$A_\infty$ property of harmonic and $p$-harmonic measures implies uniform rectifiability.}
		\newblock{\em Anal. PDE}, {10 (3): 513--558, 2017.}
		
		\bibitem[HL96]{HL96}
		S.~Hofmann and J.~L.~Lewis.
		\newblock{$L^2$ solvability and representation by caloric layer potentials in time-varying domains.}
		\newblock{\em Ann. of Math. (2)}, {144 (2): 349--420, 1996.}
		
		\bibitem[HLN03]{HLN03}
		S.~Hofmann, J.~L.~Lewis, and K.~Nyström.
		\newblock{Existence of big pieces of graphs for parabolic problems.}
		\newblock{\em Ann. Acad. Sci. Fenn. Math.}, {28 (2): 355--384, 2003.}
		
		\bibitem[HLN04]{HLN04}
		S.~Hofmann, J.~L.~Lewis, and K.~Nyström.
		\newblock{Caloric measure in parabolic flat domains.}
		\newblock{\em Duke Math. J.}, {122 (2): 281--346, 2004.}
		
		\bibitem[HMU14]{HMU14}
		S.~Hofmann, J.~M.~Martell, and I.~Uriarte-Tuero.
		\newblock{Uniform rectifiability and harmonic measure II: Poisson kernels in $L^p$ imply uniform rectifiability.}
		\newblock{\em Duke Math. J.}{, 163 (8): 1601--1654, 2014.}
		
		\bibitem[HoMc11]{HoMc11}
		S.~Hofmann and A.~McIntosh.
		\newblock{Boundedness and applications of singular integrals and square functions: a survey.}
		\newblock{\em Bull. Math. Sci.}{, 1 (2): 201--244, 2011.}
		
		\bibitem[HMM15]{HMM15}
		S.~Hofmann, M.~Mitrea, and A.~J.~Morris.
		\newblock{The method of layer potentials in $L^p$ and endpoint spaces for elliptic operators with $L^\infty$ coefficients.}
		\newblock{\em Proc. Lond. Math. Soc. (3)}{, 111 (3): 681--716, 2015.}
		
		\bibitem[HMT10]{HMT10}
		S.~Hofmann, M.~Mitrea, and M.~Taylor.
		\newblock{Singular integrals and elliptic boundary problems on regular Semmes-Kenig-Toro domains.}
		\newblock{\em Int. Math. Res. Not. IMRN}{, (14): 2567--2865, 2010.}
		
		\bibitem[HS25]{HS25}
		S.~Hofmann and D.~Sparrius.
		\newblock{The Neumann function and the $L^p$ Neumann problem in chord-arc domains.}
		\newblock{\em Adv. Nonlinear Stud.}{, 25 (2): 313--349, 2025.}
		
		\bibitem[Hy06]{Hy06}
		T.~Hyt\"onen.
		\newblock{An operator-valued $Tb$ theorem.}
		\newblock{\em J. Funct. Anal.}{, 234 (2): 420--463, 2006.}
		
		\bibitem[Hy10]{Hy10}
		T.~Hytönen.
		\newblock{A framework for non-homogeneous analysis on metric spaces, and the RBMO space of Tolsa.}
		\newblock{\em Publ. Mat., Barc.}{, 54 (2): 485--504, 2010.}
		
		\bibitem[Hy12]{Hy12}
		T.~Hytönen.
		\newblock{The sharp weighted bound for general Calder\'on-Zygmund operators.}
		\newblock{\em Ann. of Math.}{, 175: 1473--1506, 2012.}
		
		\bibitem[HK12]{HK12}
		T.~Hytönen and A.~Kairema.
		\newblock{Systems of dyadic cubes in a doubling metric space.}
		\newblock{\em Colloq. Math.}{, 126 (1): 1--33, 2012.}
		
		\bibitem[HLYY12]{HLYY12}
		T.~Hyt\"onen, S.~Liu, D.~Yang, and D.~Yang.
		\newblock{Boundedness of Calder\'on-Zygmund operators on non-homogeneous metric measure spaces.}
		\newblock{\em Canad. J. Math.}{, 64 (4): 892--923, 2012.}
		
		\bibitem[HM12]{HM12}
		T.~Hyt\"onen and H.~Martikainen.
		\newblock{Non-homogeneous $Tb$ theorem and random dyadic cubes on metric measure spaces.}
		\newblock{\em J. Geom. Anal.}{, 22 (4): 1071--1107, 2012}.
		
		\bibitem[HM14]{HM14}
		T.~Hyt\"onen and H.~Martikainen.
		\newblock{Non-homogeneous $T1$ theorem for bi-parameter singular integrals.}
		\newblock{\em Adv. Math.}{, 261: 220--273, 2014.}
		
		\bibitem[HN20]{HN20}
		T.~Hytönen and F.~Nazarov.
		\newblock{The local $Tb$ theorem with rough test functions.}
		\newblock{\em Adv. Math.}{, 372, Article no. 107306, 2020.}
		
		\bibitem[HV15]{HV15}
		T.~Hytönen and A.~Vähäkangas.
		\newblock{The local non-homogeneous $Tb$ theorem for vector-valued functions. }
		\newblock{\em Glasg. Math. J.}{, 57 (1): 17--82, 2015.}
		
		\bibitem[HYY12]{HYY12}
		T.~Hytönen, D.~Yang, and D.~Yang.
		\newblock{The Hardy space $H^1$ on non-homogeneous metric spaces. }
		\newblock{\em Math. Proc. Camb. Philos. Soc.}{, 153 (1): 9--31, 2012.}
		
		\bibitem[JK81]{JK81}
		D.~Jerison and C.~E.~Kenig.
		\newblock{The Neumann problem on Lipschitz domains.}
		\newblock{\em Bull. Amer. Math. Soc. (N.S.)}{, 4 (2): 203--207, 1981.}
		
		\bibitem[JK82]{JK82}
		D.~Jerison and C.~E.~Kenig.
		\newblock{Hardy spaces, $A_\infty$, and singular integrals on chord-arc domains.}
		\newblock{\em Math. Scand.}{, 50: 221--247, 1982.}
		
		\bibitem[KW88]{KW88}
		R.~Kaufman and J.-M.~Wu.
		\newblock{Parabolic measure on domains of class $\mathrm{Lip}\,\frac12$.}
		\newblock{\em Compositio Math.}, {65 (2): 201--207, 1988.}
		
		\bibitem[Ke94]{Ke94}
		C.~E.~Kenig.
		\newblock{\em Harmonic analysis techniques for second order elliptic boundary value problems.}
		\newblock{CBMS Regional Conference Series in Mathematics 83, American Mathematical Society, Providence, RI, 1994.}
		
		\bibitem[KP93]{KP93}
		C.~E.~Kenig and J.~Pipher.
		\newblock{The Neumann problem for elliptic equations with non-smooth coefficients.}
		\newblock{\em Invent. Math.}{, 113 (3): 447--509, 1993.}
		
		\bibitem[KP95]{KP95}
		C.~E.~Kenig and J.~Pipher.
		\newblock{The Neumann problem for elliptic equations with nonsmooth coefficients. II.}
		\newblock{\em Duke Math. J.}{, 81 (1): 227--250, 1995.}
		
		\bibitem[KR09]{KR09}
		C.~E.~Kenig and D.~J.~Rule.
		\newblock{The regularity and Neumann problem for non-symmetric elliptic operators.}
		\newblock{\em Trans. Amer. Math. Soc.}{, 361 (1): 125--160, 2009.}
		
		\bibitem[KS08]{KS08}
		A.~S.~Kim and Z.~Shen.
		\newblock{The Neumann problem in $L^p$ on Lipschitz and convex domains.}
		\newblock{\em J. Funct. Anal.}{, 255 (7): 1817--1830, 2008.}
		
		\bibitem[KM24]{KM24}
		J.~Kinnunen and K.~Myyryläinen.
		\newblock{Characterizations of parabolic Muckenhoupt classes.}
		\newblock{\em Adv. Math.}, {444, Paper No. 109612, 2024.}
		
		\bibitem[KM25]{KM25}
		J.~Kinnunen and K.~Myyryläinen.
		\newblock{Characterizations of parabolic reverse H\"older classes.}
		\newblock{\em J. Anal. Math.}, {155 (2): 609--656, 2025.}
		
		\bibitem[KMYZ24]{KMYZ24}
		J.~Kinnunen, K.~Myyryläinen, D.~Yang, and C.~Zhu.
		\newblock{Parabolic Muckenhoupt weights with time lag on spaces of homogeneous type with monotone geodesic property.}
		\newblock{\em Potential Anal.}, {60 (4): 1513--1569, 2024.}
		
		\bibitem[KS16a]{KS16a}
		J.~Kinnunen and O.~Saari.
		\newblock{Parabolic weighted norm inequalities and partial differential equations.}
		\newblock{\em Anal. PDE}, {9 (7): 1711--1736, 2016.}
		
		\bibitem[KS16b]{KS16b}
		J.~Kinnunen and O.~Saari.
		\newblock{On weights satisfying parabolic Muckenhoupt conditions.}
		\newblock{\em Nonlinear Anal.}, {131: 289--299, 2016.}
		
		\bibitem[LM16]{LM16}
		M.~Lacey and H.~Martikainen.
		\newblock{Local $Tb$ theorem with $L^2$ testing conditions and general measures: Calder\'on-Zygmund operators.}
		\newblock{\em Ann. Sci. Éc. Norm. Supér.}{, 49: 57--86, 2016.}
		
		\bibitem[LM17]{LM17}
		M.~Lacey and H.~Martikainen.
		\newblock{Local $Tb$ theorem with $L^2$ testing conditions and general measures: Square functions.}
		\newblock{\em J. Anal. Math.}{, 133: 71--89, 2017.}
		
		\bibitem[LV16]{LV16}
		M.~Lacey and A.~Vähäkangas.
		\newblock{On the local $Tb$ theorem: a direct proof under the duality assumption.}
		\newblock{\em Proc. Edinb. Math. Soc.}{, II. Ser. 59, No. 1: 193--222, 2016.}
		
		\bibitem[LM24]{LM24}
		M.~Laurel and M.~Mitrea.
		\newblock{\em Weighted Morrey spaces: Calder\'on-Zygmund theory and boundary problems.}
		\newblock{De Gruyter Studies in Mathematics 99, De Gruyter, Berlin, 2024.}

		\bibitem[LM95]{LM95}
		J.~L.~Lewis and M.~A.~M.~Murray.
		\newblock{The method of layer potentials for the heat equation in time-varying domains.}
		\newblock{\em Mem. Amer. Math. Soc.}, {114 (545), 1995.}
		
		\bibitem[LMS92]{LMS92}
		C.~Li, A.~McIntosh, and S.~Semmes.
		\newblock{Convolution singular integrals on Lipschitz surfaces.}
		\newblock{\em J. Amer. Math. Soc.}{, 5 (3): 455--481, 1992.}
		
		\bibitem[LN21]{LN21}
		M.~Litsgård and K.~Nyström.
		\newblock{On the fine properties of parabolic measures associated to strongly degenerate parabolic operators of Kolmogorov type.}
		\newblock{\em Adv. Math.}{, 387, Article ID 107833, 39 p., 2021.}
		
		\bibitem[M$^5$22]{MMMMM22}
		J.~J.~Mar\'in, J.~M.~Martell, D.~Mitrea, I.~Mitrea, and M.~Mitrea.
		\newblock{\em Singular integral operators, quantitative flatness, and boundary problems.}
		\newblock{Progress in Mathematics 344, Birkh\"auser/Springer, Cham, 2022.}
		
		\bibitem[Ma12a]{Ma12a}
		H.~Martikainen.
		\newblock{Representation of bi-parameter singular integrals by dyadic operators.}
		\newblock{\em Adv. Math.}{, 229 (3): 1734--1761, 2012.}
		
		\bibitem[Ma12b]{Ma12b}
		H.~Martikainen.
		\newblock{Vector-valued non-homogeneous $Tb$ theorem on metric measure spaces.}
		\newblock{\em Rev. Mat. Iberoam.}{, 28 (4): 961--998, 2012.}
		
		\bibitem[MM15]{MM15}
		H.~Martikainen and M.~Mourgoglou.
		\newblock{Boundedness of non-homogeneous square functions and $L^q$ type testing conditions with $q\in (1,2)$.}
		\newblock{\em Math. Res. Lett.}{, 22: 1417--1457, 2015.}
		
		\bibitem[MMT18]{MMT18}
		H.~Martikainen, M.~Mourgoglou, and X.~Tolsa.
		\newblock{Improved Cotlar’s inequality in the context of local $Tb$ theorems.}
		\newblock{\em J. Funct. Anal.}{, 274  (5): 1255--1275, 2018.}
		
		\bibitem[MMV17]{MMV17}
		H.~Martikainen, M.~Mourgoglou, and E.~Vuorinen.
		\newblock{Non-homogeneous square functions on general sets: suppression and big pieces methods.}
		\newblock{\em J. Geom. Anal.}{, 27: 3176--3227, 2017.}
		
		\bibitem[MMV21]{MMV21}
		H.~Martikainen, M.~Mourgoglou, and E.~Vuorinen.
		\newblock{A new approach to non-homogeneous local $Tb$ theorems.}
		\newblock{\em J. Anal. Math.}{, 143 (1): 95--121, 2021.}
		
		\bibitem[MV21]{MV21}
		H.~Martikainen and E.~Vuorinen.
		\newblock{Dyadic-probabilistic methods in bilinear analysis.}
		\newblock{\em Memoirs of the American Mathematical Society}{, 1344, 82 p., 2021.}
		
		\bibitem[MPT22]{MPT22}
		J.~Mateu, L.~Prat, and X.~Tolsa.
		\newblock{Removable singularities for Lipschitz caloric functions in time varying domains.}
		\newblock{\em Rev. Mat. Iberoam.}{, 38 (2): 547--588, 2022.}
		
		
		\bibitem[Mc89]{Mc89}
		A.~McIntosh.
		\newblock{Clifford algebras and the higher dimensional Cauchy integral.}
		\newblock{\em Approximation and function spaces (Warsaw, 1986)}{, Banach Center Publ. 22, 253--267, PWN, Warsaw, 1989.}
		
		\bibitem[Me26]{Me26}
		A.~Merlo.
		\newblock{On the WALA conjecture, Alberti representations and applications.}
		\newblock{\em arXiv preprint,}{ https://arxiv.org/abs/2605.29804, 2026.}
		
		\bibitem[MMP25]{MMP25}
		A.~Merlo, M.~Mourgoglou, and C.~Puliatti.
		\newblock{Layer potentials for elliptic operators with DMO-type coefficients: big pieces $Tb$
			theorem, quantitative rectifiability, and free boundary problems.}
		\newblock{\em arXiv preprint,}{ https://arxiv.org/abs/2505.23478, 2025.}
		
		\bibitem[MMM23a]{MMM23a}
		D.~Mitrea, I.~Mitrea, and M.~Mitrea.
		\newblock{\em Geometric harmonic analysis IV. Boundary layer potentials in uniformly rectifiable domains, and applications to complex analysis.}
		\newblock{Developments in Mathematics 75, Springer, Cham, 2023.}
		
		\bibitem[MMM23b]{MMM23}
		D.~Mitrea, I.~Mitrea, and M.~Mitrea.
		\newblock{\em Geometric harmonic analysis V. Fredholm theory and finer estimates for integral operators, with applications to boundary problems.}
		\newblock{Developments in Mathematics 76, Springer, Cham, 2023.}
		
		\bibitem[Mi94]{Mi94}
		M.~Mitrea.
		\newblock{\em Clifford wavelets, singular integrals, and Hardy spaces.}
		\newblock{Lecture Notes in Mathematics 1575, Springer-Verlag, Berlin, 1994.}
		
		\bibitem[MiTa26]{MiTa26}
		M.~Mitrea and P.~Takemura.
		\newblock{\em Singular integrals, Herz-type function spaces, and boundary problems.}
		\newblock{Progress in Mathematics 362, Birkh\"auser/Springer, Cham, 2026.}
		
		\bibitem[MPT25]{MPT25}
		M.~Mourgoglou, B.~Poggi, and X.~Tolsa.
		\newblock{Solvability of the Poisson-Dirichlet problem with interior data in $L^{p'}$-Carleson spaces and its applications to the $L^p$-regularity problem.}
		\newblock{\em J. Eur. Math. Soc. (JEMS)}{, published online first, 2025, https://doi.org/10.4171/JEMS/1660.}
		
		\bibitem[MP21]{MP21}
		M.~Mourgoglou and C.~Puliatti.
		\newblock{Blow-ups of caloric measure in time varying domains and applications to two-phase problems. }
		\newblock{\em J. Math. Pures Appl. (9)}{, 152: 1--68, 2021.}
		
		\bibitem[MT20]{MT20}
		M.~Mourgoglou and X.~Tolsa.
		\newblock{Harmonic measure and Riesz transform in uniform and general domains.}
		\newblock{\em J. Reine Angew. Math.}{, 758: 183--221, 2020.}
		
		\bibitem[MT24]{MT24}
		M.~Mourgoglou and X.~Tolsa.
		\emph{The regularity problem for the Laplace equation in rough domains}.
		Duke Math. J. \textbf{173} (2024), no.~9, 1731--1837.
		
		\bibitem[MT24b]{MT24b}
		M.~Mourgoglou and X.~Tolsa.
		\newblock{Solvability of the Neumann problem for elliptic equations in chord-arc domains with very big pieces of good superdomains.}
		\newblock{\em arXiv preprint,}{ https://arxiv.org/abs/2407.20385, 2024.}
		
		\bibitem[Mur85]{Mur85}
		M.~A.~M.~Murray.
		\newblock{The Cauchy integral, Calder\'on commutators, and conjugations of singular integrals in $\mathbb R^n$.}
		\newblock{\em Trans. Amer. Math. Soc.}{, 289 (2): 497--518, 1985.}
		
		\bibitem[NToV14]{NToV14}
		F.~Nazarov, X.~Tolsa, and A.~Volberg.
		\newblock{On the uniform rectifiability of AD-regular measures with bounded Riesz transform operator: the case of codimension 1.}
		\newblock{\em Acta Math.}, {213 (2): 237--321, 2014.}
		
		\bibitem[NTrV99]{NTrV99}
		F.~Nazarov, S.~Treil, and A.~Volberg.
		\newblock{The $Tb$-theorem on non-homogeneous spaces that proves a conjecture of Vitushkin.}
		\newblock{\em arXiv preprint arXiv:1401.2479}{, 2014 (original preprint: 1999).}
		
		\bibitem[NTrV02]{NTrV02}
		F.~Nazarov, S.~Treil and A.~Volberg.
		\newblock{Accretive system $Tb$-theorems on nonhomogeneous spaces.}
		\newblock{\em Duke Math. J.}{, 113 (2): 259--312, 2002.}
		
		\bibitem[NTrV03]{NTrV03}
		F.~Nazarov, S.~Treil, and A.~Volberg,
		\newblock{The $Tb$-theorem on non-homogeneous spaces.}
		\newblock{\em Acta Math.}{, 190 (2): 151--239, 2003.}
		
		\bibitem[PPT21]{PPT21}
		L.~Prat, C.~Puliatti, and X.~Tolsa.
		\newblock{$L^2$-boundedness of gradients of single-layer potentials and uniform rectifiability.}
		\newblock{\em Anal. PDE,}{ 14 (3): 717-791, 2021.}
		
		\bibitem[Pu22]{Pu19}
		C.~Puliatti.
		\newblock{Gradient of the single layer potential and quantitative rectifiability for general Radon measures.}
		\newblock{\em J. Funct. Anal.,}{ 282 (6), Article 109376, 2022.}
		
		\bibitem[Saw86]{Saw86}
		E.~Sawyer.
		\newblock{Weighted inequalities for the one-sided Hardy-Littlewood maximal functions.}
		\newblock{\em Trans. Amer. Math. Soc.}, {297 (1): 53--61, 1986.}
		
		\bibitem[Se89]{Se89}
		S.~Semmes.
		\newblock{A criterion for the boundedness of singular integrals on hypersurfaces.}
		\newblock{\em Trans. Amer. Math. Soc.}{, 311 (2): 501--513, 1989.}
		
		\bibitem[Se90]{Se90}
		S.~Semmes.
		\newblock{Analysis vs. geometry on a class of rectifiable hypersurfaces in $\mathbb R^n$.}
		\newblock{\em Indiana Univ. Math. J.}{, 39: 1005--1035, 1990.}
		
		\bibitem[Sh07]{Sh07}
		Z.~Shen.
		\newblock{The $L^p$ boundary value problems on Lipschitz domains.}
		\newblock{\em Adv. Math.}{, 216 (1): 212--254, 2007.}
		
		\bibitem[To03]{To03}
		X.~Tolsa.
		\newblock{Painlevé’s problem and the semiadditivity of analytic capacity.}
		\newblock{\em Acta Math.}{, 190 (1): 105--149, 2003.}
		
		
		\bibitem[To14]{To14}
		X.~Tolsa.
		\newblock{\em Analytic capacity, the Cauchy transform, and non-homogeneous Calder\'on-Zygmund theory.}
		\newblock{Progress in Mathematics 307. Birkhäuser/Springer. xiii, 396 p., 2014.}
		
		\bibitem[To17]{To16}
		X.~Tolsa.
		\newblock{Rectifiable measures, square functions involving densities, and the Cauchy transform.}
		\newblock{\em Mem. Am. Math. Soc.}{, 245, Issue 1158, 2017.}
		
		\bibitem[To26]{To26}
		X.~Tolsa.
		\newblock{Interactions between quantitative rectifiability, singular integrals, and boundary value problems for harmonic functions.}
		\newblock{\em arXiv preprint,}{ https://arxiv.org/abs/2607.16457, 2026.}
		
		\bibitem[Tor10]{Tor10}
		T.~Toro.
		\newblock{Potential analysis meets geometric measure theory.}
		\newblock{\em Proceedings of the International Congress of Mathematicians (Hyderabad, 2010), Vol. III}{, 1485--1500, Hindustan Book Agency, New Delhi, 2010.}
		
		\bibitem[Ve84]{Ve84}
		G.~C.~Verchota.
		\newblock{Layer potentials and regularity for the Dirichlet problem for Laplace's equation in Lipschitz domains.}
		\newblock{\em J. Funct. Anal.}{, 59 (3): 572--611, 1984.}
		
		\bibitem[Vi24]{Vi24}
		P.~Villarroya.
		\newblock{New local $T1$ theorems on non-homogeneous spaces.}
		\newblock{\em Publ. Mat., Barc.}{, 68 (2): 445--506, 2024.}
		
		\bibitem[Vo03]{Vo03}
		{A.~Volberg}.
		\newblock{\em Calderón-Zygmund capacities and operators on nonhomogeneous spaces.}
		\newblock{CBMS Regional Conference Series in Mathematics 100. Providence, RI: American Mathematical Society (AMS), 2003.}
		
		\bibitem[VZ18]{VZ18}
		A.~Volberg and P.~Zorin-Kranich.
		\newblock{Sparse domination on non-homogeneous spaces with an application to $A_p$ weights.}
		\newblock{\em Rev. Mat. Iberoam.}{, 34 (3): 1401--1414, 2018.}
		
		\bibitem[Wa12]{Wa12}
		N.A.~Watson.
		\newblock{\em Introduction to Heat Potential Theory.}
		\newblock{Mathematical Surveys and Monographs, Vol 182, American Mathematical Society, 2012.}
		
	\end{thebibliography}
\end{document}